\documentclass[11pt]{article}
\usepackage{amssymb,epsfig,amsmath,mathrsfs,bbm} 
\usepackage[T1]{fontenc}
\usepackage[dvips]{color}
\usepackage{comment}
\usepackage[utf8]{inputenc}
\usepackage[colorlinks]{hyperref}
\usepackage[acronym,toc]{glossaries}
\newglossary[slg]{symbols}{sym}{sbl}{List of Symbols}

\newtheorem{defi}{Definition}[section]
\newtheorem{prop}[defi]{Proposition}
\newtheorem{theo}[defi]{Theorem}
\newtheorem{conj}[defi]{Conjecture}
\newtheorem{lemm}[defi]{Lemma}
\newtheorem{coro}[defi]{Corollary}
\newtheorem{rema}[defi]{Remark}
\newtheorem{exem}[defi]{Example}
\newtheorem{exems}[defi]{Examples}

\newcommand{\bdefi}{\begin{defi}}
\newcommand{\edefi}{\end{defi}}
\newcommand{\bprop}{\begin{prop}}
\newcommand{\eprop}{\end{prop}}
\newcommand{\btheo}{\begin{theo}}
\newcommand{\etheo}{\end{theo}}
\newcommand{\blemm}{\begin{lemm}}
\newcommand{\brema}{\begin{rema}}
\newcommand{\erema}{\end{rema}}
\newcommand{\bexer}{\begin{exem}}
\newcommand{\eexer}{\end{exem}}
\newcommand{\bexems}{\begin{exems}}
\newcommand{\eexems}{\end{exems}}
\newcommand{\bconj}{\begin{conj}}
\newcommand{\econj}{\end{conj}}
\newcommand{\elemm}{\end{lemm}}
\newcommand{\bcoro}{\begin{coro}}
\newcommand{\ecoro}{\end{coro}}
\newcommand{\dem}{\noindent{\bf Proof. }}
\newcommand{\rem}{\noindent{\bf Remark. }}

\renewcommand{\cal}{\mathscr}

\newcommand{\A}{{\cal A}}

\newcommand{\T}{{\cal T}}
\newcommand{\B}{{\cal B}}
\renewcommand{\L}{{\cal L}}
\newcommand{\M}{{\cal M}}
\newcommand{\W}{{\cal W}}
\newcommand{\N}{{\cal N}}
\newcommand{\Q}{{\cal Q}}
\newcommand{\G}{{\cal G}}
\newcommand{\D}{{\cal D}}
\newcommand{\E}{{\cal E}}
\newcommand{\F}{{\cal F}}

\renewcommand{\H}{{\cal H}}
\newcommand{\OOO}{{\cal O}}
\newcommand{\C}{{\cal C}}

\newcommand{\PPP}{{\cal P}}

\newcommand{\maths}[1]{{\mathbb #1}}  

\newcommand{\RR}{\maths{R}}
\newcommand{\MM}{\maths{M}}
\newcommand{\NN}{\maths{N}}
\newcommand{\CC}{\maths{C}}

\renewcommand{\SS}{\maths{S}}
\newcommand{\DD}{\maths{D}}
\newcommand{\HH}{\maths{H}}

\newcommand{\ZZ}{\maths{Z}}
\newcommand{\PP}{\maths{P}}

\newcommand{\LL}{\maths{L}}

\newcommand{\UU}{\maths{U}}

\newcommand{\mmm}{{\mathfrak m}}

\newcommand{\ra}{\rightarrow}
\newcommand{\ov}[1]{{\overline #1}} 
\newcommand{\wt}[1]{{\widetilde{#1}}}
\newcommand{\wh}[1]{{\widehat{#1}}}
\newcommand{\ga}{\gamma}
\newcommand{\Ga}{\Gamma}
\newcommand{\bs}{\backslash}
\newcommand{\cinf}{{\operatorname{C}^\infty}}
\renewcommand{\ln}{\log}

\newcommand{\weakstar}{\overset{*}\rightharpoonup}
\newcommand{\cqfd}{\hfill$\Box$}
\newcommand{\card}{\operatorname{Card}}
\newcommand{\vol}{\operatorname{vol}}
\newcommand{\prT}{{T\wt p}}
\newcommand{\Par}{\operatorname{Par}}
\newcommand{\Per}{\operatorname{\PPP er}}
\newcommand{\Long}{\operatorname{length}}
\newcommand{\Isom}{\operatorname{Isom}}
\newcommand{\Hom}{\operatorname{Hom}}
\newcommand{\Ker}{\operatorname{Ker}}
\newcommand{\id}{\operatorname{id}}
\newcommand{\Vol}{\operatorname{Vol}}
\newcommand{\Axe}{\operatorname{Axe}}

\newglossaryentry{valeurabsolue}{type=symbols,
  name={\ensuremath{\lfloor x\rfloor}}, sort=aaaa,
  description={~~~~integral part of $x\in\RR$}}

\newglossaryentry{crossratio}{type=symbols,
  name={\ensuremath{[a,b,c,d]_F}}, sort=aaaaa,
  description={~~~~crossratio of pairwise distinct $a,b,c,d\in\wt M$
    with respect to $F$}}

\newglossaryentry{characteristicfunction}{type=symbols,
  name={\ensuremath{\mathbbm{1}_A}}, sort=aaaaaa,
  description={~~~~characteristic function of a subset $A$}}

\newglossaryentry{axetranslation}{type=symbols,
  name={\ensuremath{\Axe_\ga}}, sort=busemann,
  description={~~~~translation axis of an isometry $\ga$ of $\wt M$}}

\newglossaryentry{busemann}{type=symbols,
  name={\ensuremath{\beta_\xi(x,y)}}, sort=busemann,
  description={~~~~Busemann cocycle}}

\newglossaryentry{ssballriem}{type=symbols,
  name={\ensuremath{\B^{\rm ss}(v,r)}}, sort=bzssriem,
  description={~~~~Ball of center $v$ and radius $r$ for the 
induced Riemannian metric on $W^{\rm ss}(v)$}}

\newglossaryentry{suballriem}{type=symbols,
  name={\ensuremath{\B^{\rm su}(v,r)}}, sort=bzsuriem,
  description={~~~~Ball of center $v$ and radius $r$ for the 
induced Riemannian metric on $W^{\rm su}(v)$}}

\newglossaryentry{ssballHam}{type=symbols,
  name={\ensuremath{B^{\rm ss}(v,r)}}, sort=bzssHam,
  description={~~~~Ball of center $v$ and radius $r$ for the 
Hamenst\"adt distance on $W^{\rm ss}(v)$}}

\newglossaryentry{suballHam}{type=symbols,
  name={\ensuremath{B^{\rm su}(v,r)}}, sort=bzsuHam,
  description={~~~~Ball of center $v$ and radius $r$ for the 
Hamenst\"adt distance on $W^{\rm su}(v)$}}

\newglossaryentry{Gibbscocycle}{type=symbols,
  name={\ensuremath{C_{F,\,\xi}(x,y)}}, sort=cocycle,
  description={~~~~Gibbs cocycle for the potential $F$}}

\newglossaryentry{complementary}{type=symbols,
  name={\ensuremath{^c A}}, sort=complementary,
  description={~~~~complementary set of a subset $A$}}

\newglossaryentry{cone}{type=symbols, name={\ensuremath{\C_xB}},
  sort=cone, description={~~~~cone on $B\subset \partial_\infty\wt M$
    with vertex $x\in\wt M\cup \partial_\infty
    \wt M$}}

\newglossaryentry{conethick}{type=symbols,
  name={\ensuremath{\C^+_r(z,Z)}}, sort=conea,
  description={~~~~$r$-thickened cone over $Z\subset\partial_\infty
    \wt M$ with vertex $z\in\wt M$}}

\newglossaryentry{conethin}{type=symbols,
  name={\ensuremath{\C^-_r(z,Z)}}, sort=coneb,
  description={~~~~$r$-thinned cone over $Z\subset\partial_\infty \wt
    M$ with vertex $z\in\wt M$}}

\newglossaryentry{convexhull}{type=symbols,
  name={\ensuremath{\C\Lambda\Ga}}, sort=convex, description={~~~~convex
    hull of the limit set of $\Ga$}}

\newglossaryentry{criticalexponentusuel}{type=symbols,
  name={\ensuremath{\delta_{\Ga}}}, sort=delta,
  description={~~~~ critical exponent of $\Ga$}}

\newglossaryentry{criticalexponent}{type=symbols,
  name={\ensuremath{\delta_{\Ga,\,F}}}, sort=deltaa,
  description={~~~~ critical exponent of $(\Ga,F)$}}

\newglossaryentry{criticalexponenttwi}{type=symbols,
  name={\ensuremath{\delta_{\Ga,\,F,\,\chi}}}, sort=deltaaa,
  description={~~~~ (twisted) critical exponent of $(\Ga,F,\chi)$}}

\newglossaryentry{boundary}{type=symbols,
  name={\ensuremath{\partial_\infty\wt M}}, sort=deltap,
  description={~~~~boundary at infinity of $\wt M$}}

\newglossaryentry{boundarydeux}{type=symbols,
  name={\ensuremath{\partial^2_\infty\wt M}}, sort=deltapdeux,
  description={~~~~space of ordered pairs of distinct points of
    $\partial_\infty\wt M$}}

\newglossaryentry{dMandTunwtM}{type=symbols, name={\ensuremath{d}},
  sort=distance, description={~~~~distance on $\wt M$ and
    $T^1\wt M$, and on any metric space}}

\newglossaryentry{dMandTunwtMprime}{type=symbols, name={\ensuremath{d'}},
  sort=distancea, description={~~~~another distance on  $T^1\wt M$}}

\newglossaryentry{visualdist}{type=symbols, name={\ensuremath{d_x}},
  sort=distanceb, description={~~~~visual distance on 
$\partial_\infty\wt M$ seen from $x\in M$}}

\newglossaryentry{Hamdiss}{type=symbols, name={\ensuremath{d_{W^{\rm
          ss}(w)}}}, sort=distanceHam, description={~~~~Hamenst\"adt's
    distance on the strong stable leaf of $w\in T^1\wt M$}}

\newglossaryentry{Hamdisu}{type=symbols, name={\ensuremath{d_{W^{\rm
          su}(w)}}}, sort=distanceHam, description={~~~~Hamenst\"adt's
    distance on the strong unstable leaf of $w\in T^1\wt M$}}

\newglossaryentry{dirac}{type=symbols, name={\ensuremath{\D_z}},
  sort=ditrac, description={~~~~unit Dirac mass at a point $z$}}

\newglossaryentry{potentialgap}{type=symbols,
  name={\ensuremath{D_{F,\,x}(\xi,\eta)}}, sort=dzgap,
  description={~~~~potential gap seen from $x\in \wt M$ between
    $\xi,\eta\in\partial_\infty \wt M$}}

\newglossaryentry{stronunstabtang}{type=symbols,
  name={\ensuremath{E^{\,\rm su}(v)}}, sort=ee,
  description={~~~~tangent space at $v\in T^1\wt M$ or $v\in T^1 M$
    to $W^{su}(v)$}}

\newglossaryentry{stronstabtang}{type=symbols,
  name={\ensuremath{E^{\,\rm ss}(v)}}, sort=eee,
  description={~~~~tangent space at $v\in T^1\wt M$ or $v\in T^1 M$
    to $W^{\rm ss}(v)$}}

\newglossaryentry{unstabtang}{type=symbols,
  name={\ensuremath{E^{\,\rm u}(v)}}, sort=eeee,
  description={~~~~tangent space at $v\in T^1\wt M$ or $v\in T^1 M$
    to $W^{\rm u}(v)$}}

\newglossaryentry{stabtang}{type=symbols,
  name={\ensuremath{E^{\,\rm s}(v)}}, sort=eeeee,
  description={~~~~tangent space at $v\in T^1\wt M$ or $v\in T^1 M$
    to $W^{\rm s}(v)$}}

\newglossaryentry{neuttang}{type=symbols,
  name={\ensuremath{E^{0}(v)}}, sort=eeeeee,
  description={~~~~tangent space at $v\in T^1\wt M$ or $v\in T^1 M$
    to $t\mapsto \phi_tv$}}

\newglossaryentry{potentialtilde}{type=symbols, name={\ensuremath{\wt
      F}}, sort=Fonction, description={~~~~potential on $T^1\wt M$}}

\newglossaryentry{potential}{type=symbols, name={\ensuremath{F}},
  sort=Fonctionbas, description={~~~~potential on $T^1M$}}

\newglossaryentry{potentialsutilde}{type=symbols,
  name={\ensuremath{\wt F^{\,\rm su}}}, sort=Fonctionsu,
  description={~~~~strong unstable potential on $T^1\wt M$}}

\newglossaryentry{potentialsu}{type=symbols,
  name={\ensuremath{F^{\,\rm su}}}, sort=Fonctionsubas,
  description={~~~~strong unstable potential on $T^1 M$}}

\newglossaryentry{bisectorbcountfunc}{type=symbols,
  name={\ensuremath{G_{\Ga,\,F,\,x,\,y,\,U,\,V}}}, sort=gga,
  description={~~~~bisectorial orbital counting function}}

\newglossaryentry{twibisectorbcountfunc}{type=symbols,
  name={\ensuremath{G_{\Ga,\,F,\,\chi,\,x,\,y,\,U,\,V}}}, sort=ggaa,
  description={~~~~twisted bisectorial orbital counting function}}

\newglossaryentry{annbisectorbcountfunc}{type=symbols,
  name={\ensuremath{G_{\Ga,\,F,\,x,\,y,\,U,\,V,\,c}}}, sort=ggab,
  description={~~~~annular bisectorial orbital counting function}}

\newglossaryentry{anntwibisectorbcountfunc}{type=symbols,
  name={\ensuremath{G_{\Ga,\,F,\,\chi,\,x,\,y,\,U,\,V,\,c}}}, sort=ggaba,
  description={~~~~twisted annular bisectorial orbital counting function}}

\newglossaryentry{sectorbcountfunc}{type=symbols,
  name={\ensuremath{G_{\Ga,\,F,\,x,\,y,\,U}}}, sort=ggb,
  description={~~~~sectorial orbital counting function}}

\newglossaryentry{annsectorbcountfunc}{type=symbols,
  name={\ensuremath{G_{\Ga,\,F,\,x,\,y,\,U,\,c}}}, sort=ggba,
  description={~~~~annular sectorial orbital counting function}}

\newglossaryentry{twisectorbcountfunc}{type=symbols,
  name={\ensuremath{G_{\Ga,\,F,\,\chi,\,x,\,y,\,U}}}, sort=ggbb,
  description={~~~~twisted sectorial orbital counting function}}

\newglossaryentry{anntwisectorbcountfunc}{type=symbols,
  name={\ensuremath{G_{\Ga,\,F,\,\chi,\,x,\,y,\,U,\,c}}}, sort=ggbba,
  description={~~~~twisted annular sectorial orbital counting function}}

\newglossaryentry{orbcountfunc}{type=symbols,
  name={\ensuremath{G_{\Ga,\,F,\,x,\,y}}}, sort=ggc,
  description={~~~~orbital counting function}}

\newglossaryentry{annorbcountfunc}{type=symbols,
  name={\ensuremath{G_{\Ga,\,F,\,x,\,y,\,c}}}, sort=ggc,
  description={~~~~annular orbital counting function}}

\newglossaryentry{twiorbcountfunc}{type=symbols,
  name={\ensuremath{G_{\Ga,\,F,\,\chi,\,x,\,y}}}, sort=ggcc,
  description={~~~~twisted orbital counting function}}

\newglossaryentry{anntwiorbcountfunc}{type=symbols,
  name={\ensuremath{G_{\Ga,\,F,\,\chi,\,x,\,y,\,c}}}, sort=ggcca,
  description={~~~~twisted annular orbital counting function}}

\newglossaryentry{Gamma}{type=symbols, name={\ensuremath{\Ga}},
  sort=Gamma, description={~~~~nonelementary discrete group of
    isometries of $\wt M$}}

\newglossaryentry{iota}{type=symbols, name={\ensuremath{\iota}},
  sort=iota, description={~~~~flip map $v\mapsto -v$ on $T^1\wt M$ and
    $T^1M$}}

\newglossaryentry{iotasubx}{type=symbols, name={\ensuremath{\iota_x}},
  sort=iotasubx, description={~~~~antipodal map with respect to
    $x\in\wt M$}}

\newglossaryentry{isometrygroup}{type=symbols,
  name={\ensuremath{\Isom(\wt M)}}, sort=isometrygroup,
  description={~~~~isometry group of $\wt M$}}

\newglossaryentry{Kplus}{type=symbols, name={\ensuremath{K^+(z,r,Z)}},
  sort=Kplus, description={~~~~set of $v\in T^1\wt M$ such that
    $v_+\in Z$ and $v$ is $r$-close to $z$}}

\newglossaryentry{Kmoins}{type=symbols, name={\ensuremath{K^-(z,r,Z)}},
  sort=Kmoins, description={~~~~set of $v\in T^1\wt M$ such that
    $v_-\in Z$ and $v$ is $r$-close to $z$}}

\newglossaryentry{lambdaga}{type=symbols,
  name={\ensuremath{\Lambda\Ga}}, sort=lambda, description={~~~~limit
    set of $\Ga$}}

\newglossaryentry{lambdaconga}{type=symbols,
  name={\ensuremath{\Lambda_c\Ga}}, sort=lambdac,
  description={~~~~conical (or radial) limit set of $\Ga$}}

\newglossaryentry{lambdacongar}{type=symbols,
  name={\ensuremath{\Lambda_{c,\,r}\Ga}}, sort=lambdacr,
  description={~~~~filtration of the conical limit set of $\Ga$}}

\newglossaryentry{lambdaMyrga}{type=symbols,
  name={\ensuremath{\Lambda_{\rm Myr}\Ga}}, sort=lambdam,
  description={~~~~Myrberg limit set of $\Ga$}}

\newglossaryentry{translationlength}{type=symbols,
  name={\ensuremath{\ell(\ga)}}, sort=length,
  description={~~~~translation length of $\ga\in\Isom(\wt M)$}}

\newglossaryentry{lebesguemeasure}{type=symbols,
  name={\ensuremath{\L_g}}, sort=leo, description={~~~~Lebesgue
    measure along a periodic orbit $g$}}

\newglossaryentry{Lset}{type=symbols, name={\ensuremath{L_r(z,w)}},
  sort=lep, description={~~~~set of $(\xi,\eta)\in(\partial_\infty
    \wt M)^2$ such that $]\xi,\eta[$ meets first $B(z,r)$ then
    $B(w,r)$}}

\newglossaryentry{ln}{type=symbols, name={\ensuremath{\ln}}, sort=log,
  description={~~~~natural logarithm (with $\ln(e)=1$)}}

\newglossaryentry{Mtilde}{type=symbols, name={\ensuremath{\wt M}},
  sort=manifold, description={~~~~pinched negatively curved complete
    simply connected Riemannian manifold}}

\newglossaryentry{M}{type=symbols,
  name={\ensuremath{M}}, sort=manifoldbas,
  description={~~~~Riemannian orbifold $M=\Ga\bs \wt M$}}

\newglossaryentry{measuregibbs}{type=symbols,
  name={\ensuremath{\wt m_F}}, sort=measuregibbs,
  description={~~~~Gibbs measure on $T^1\wt M$ with
  potential $F$}}

\newglossaryentry{measuregibbsbas}{type=symbols,
  name={\ensuremath{m_F}}, sort=measuregibbsbas,
  description={~~~~Gibbs measure on $T^1 M$ with
  potential $F$}}

\newglossaryentry{neighbourhood}{type=symbols,
  name={\ensuremath{\N_rA}}, sort=neighbourhood,
  description={~~~~open $r$-neighbourhood of $A$}}

\newglossaryentry{rinterior}{type=symbols,
  name={\ensuremath{\N^-_rA}}, sort=neighbourhoodi,
  description={~~~~open $r$-interior of $A$}}

\newglossaryentry{normaliser}{type=symbols,
  name={\ensuremath{N(\Ga)}}, sort=normaliser,
  description={~~~~normaliser of $\Ga$ in $\Isom(\wt M)$}}

\newglossaryentry{omegatilde}{type=symbols,
  name={\ensuremath{\wt\Omega\Ga}}, sort=omega, description={~~~~set
    of $v\in T^1\wt M$ such that $v_-,v_+\in\Lambda\Ga$}}

\newglossaryentry{omegatildec}{type=symbols,
  name={\ensuremath{\wt\Omega_c\Ga}}, sort=omegaa, description={~~~~set
    of $v\in T^1\wt M$ such that $v_-,v_+\in\Lambda_c\Ga$}}

\newglossaryentry{omega}{type=symbols, name={\ensuremath{\Omega\Ga}},
  sort=omegabas, description={~~~~nonwandering set of the geodesic
    flow in $T^1M$}}

\newglossaryentry{omegac}{type=symbols, name={\ensuremath{\Omega_c\Ga}},
  sort=omegabasc, description={~~~~two-sided
  recurrent set of the geodesic
    flow in $T^1M$}}

\newglossaryentry{shadow}{type=symbols, name={\ensuremath{\OOO_xA}},
  sort=ooo, description={~~~~shadow of $A\subset \wt M$ seen form
    $x\in\wt M\cup\partial_\infty\wt M$}}

\newglossaryentry{geodesicflow}{type=symbols,
  name={\ensuremath{\phi_t}}, sort=phi, description={~~~~geodesic flow
    at time $t\in\RR$ on $T^1\wt M$ and on $T^1M$}}

\newglossaryentry{basepointproj}{type=symbols,
  name={\ensuremath{\pi}}, sort=pi, description={~~~~base point
    projections $T^1\wt M \ra \wt M$ and $T^1M \ra M$}}

\newglossaryentry{Gurevichpressure}{type=symbols,
  name={\ensuremath{P_{Gur}(\Ga,\,F)}}, sort=Pkgur,
  description={~~~~Gurevich pressure of $(\Ga,F)$}}

\newglossaryentry{Gurevichpressuretwi}{type=symbols,
  name={\ensuremath{P_{Gur}(\Ga,\,F,\,\chi)}}, sort=Pkgura,
  description={~~~~(twisted) Gurevich pressure of $(\Ga,F,\chi)$}}

\newglossaryentry{metricpressure}{type=symbols,
  name={\ensuremath{P_{\Ga,\,F}(m)}}, sort=Plressmes,
  description={~~~~metric pressure of a measure $m$}}

\newglossaryentry{topologicalpressure}{type=symbols,
  name={\ensuremath{P(\Ga,F)}}, sort=Plresstop,
  description={~~~~topological pressure of $(\Ga,F)$}}

\newglossaryentry{setperiod}{type=symbols,
  name={\ensuremath{\Per(s)}}, sort=pmeriod,
  description={~~~~Set of periodic orbits of $(\phi_t)_{t\in\RR}$
of length $\leq s$}}

\newglossaryentry{setperiodprim}{type=symbols,
  name={\ensuremath{\Per'(s)}}, sort=pmerioda,
  description={~~~~Set of primitive periodic orbits of $(\phi_t)_{t\in\RR}$
of length $\leq s$}}

\newglossaryentry{Poincarserieusuel}{type=symbols,
  name={\ensuremath{Q_{\Ga}=Q_{\Ga,\,x,\,y}}}, sort=poincare,
  description={~~~~Poincar\'e series of $\Ga$}}

\newglossaryentry{Poincarserie}{type=symbols,
  name={\ensuremath{Q_{\Ga,\,F}=Q_{\Ga,\,F,\,x,\,y}}}, sort=poincarea,
  description={~~~~Poincar\'e series of $(\Ga,F)$}}

\newglossaryentry{Poincarserietwi}{type=symbols,
  name={\ensuremath{Q_{\Ga,\,F,\,\chi,\,x,\,y}}}, sort=poincareb,
  description={~~~~(twisted) Poincar\'e series of $(\Ga,F,\chi)$}}

\newglossaryentry{tangentmap}{type=symbols,
  name={\ensuremath{Tf}}, sort=tangentmap,
  description={~~~~tangent map $TN\ra TN'$ of a smooth map $f:N\ra N'$}}

\newglossaryentry{tangentspace}{type=symbols,
  name={\ensuremath{T^1\wt M}}, sort=tangentspace,
  description={~~~~(total space of the) unit tangent bundle of $\wt M$}}

\newglossaryentry{tangentspacebas}{type=symbols,
  name={\ensuremath{T^1M}}, sort=tangentspacebas,
  description={~~~~orbifold $\Ga\bs T^1\wt M$}}

\newglossaryentry{bitangentspacebas}{type=symbols,
  name={\ensuremath{TT^1M}}, sort=tangentspacebasbi,
  description={~~~~orbifold $\Ga\bs TT^1\wt M$}}

\newglossaryentry{vminus}{type=symbols, name={\ensuremath{v_-}},
  sort=vminus, description={~~~~point at $-\infty$ of the
    geodesic line defined by $v\in T^1\wt M$}}

\newglossaryentry{vplus}{type=symbols, name={\ensuremath{v_+}},
  sort=vplus, description={~~~~point at $+\infty$ of the
    geodesic line defined by $v\in T^1\wt M$}}

\newglossaryentry{Wss}{type=symbols, name={\ensuremath{W^{\rm ss}(v)}},
  sort=Wss, description={~~~~strong stable
  leaf of $v\in T^1\wt M$ or $v\in T^1M$}}

\newglossaryentry{Wsu}{type=symbols, name={\ensuremath{W^{\rm su}(v)}},
  sort=Wsu, description={~~~~strong unstable
  leaf of $v\in T^1\wt M$ or $v\in T^1M$}}

\newglossaryentry{Ws}{type=symbols, name={\ensuremath{W^{\rm s}(v)}},
  sort=Wssz, description={~~~~stable
  leaf of $v\in T^1\wt M$ or $v\in T^1M$}}

\newglossaryentry{Wu}{type=symbols, name={\ensuremath{W^{\rm u}(v)}},
  sort=Wsuz, description={~~~~unstable
  leaf of $v\in T^1\wt M$ or $v\in T^1M$}}

\newglossaryentry{Wssfol}{type=symbols,
  name={\ensuremath{\wt{\W}^{\rm ss}}}, sort=Wssfol,
  description={~~~~strong stable foliation of $T^1\wt M$}}

\newglossaryentry{Wsufol}{type=symbols,
  name={\ensuremath{\wt{\W}^{\rm su}}}, sort=Wsufol,
  description={~~~~strong unstable foliation of $T^1\wt M$}}

\newglossaryentry{Wsfol}{type=symbols,
  name={\ensuremath{\wt{\W}^{\rm s}}}, sort=Wsszfol,
  description={~~~~stable foliation of $T^1\wt M$}}

\newglossaryentry{Wufol}{type=symbols,
  name={\ensuremath{\wt{\W}^{\rm u}}}, sort=Wsuzfol,
  description={~~~~unstable foliation of $T^1\wt M$}}

\newglossaryentry{Wssfolbas}{type=symbols,
  name={\ensuremath{{\W}^{\rm ss}}}, sort=Wssfolbas,
  description={~~~~strong stable foliation of $T^1 M$}}

\newglossaryentry{Wsufolbas}{type=symbols,
  name={\ensuremath{{\W}^{\rm su}}}, sort=Wsufolbas,
  description={~~~~strong unstable foliation of $T^1 M$}}

\newglossaryentry{Wsfolbas}{type=symbols,
  name={\ensuremath{{\W}^{\rm s}}}, sort=Wsszfolbas,
  description={~~~~stable foliation of $T^1 M$}}

\newglossaryentry{Wufolbas}{type=symbols,
  name={\ensuremath{{\W}^{\rm u}}}, sort=Wsuzfolbas,
  description={~~~~unstable foliation of $T^1 M$}}

\newglossaryentry{periodcountingfucnt}{type=symbols,
  name={\ensuremath{Z_{\Ga,\,F,\,W}}}, sort=Zperiod,
  description={~~~~period counting function of $(\Ga,F,W)$}}

\newglossaryentry{annperiodcountingfucnt}{type=symbols,
  name={\ensuremath{Z_{\Ga,\,F,\,W,\,c}}}, sort=Zperioda,
  description={~~~~annular period counting function of $(\Ga,F,W)$}}

\newglossaryentry{periodcountingfucnttwi}{type=symbols,
  name={\ensuremath{Z_{\Ga,\,F,\,\chi,\,W}}}, sort=Zperiodb,
  description={~~~~(twisted) period counting function of $(\Ga,F,\chi,W)$}}

\newglossaryentry{annperiodcountingfucnttwi}{type=symbols,
  name={\ensuremath{Z_{\Ga,\,F,\,\chi,\,W,\,c}}}, sort=Zperiodc,
  description={~~~~(twisted) annular period counting function of 
$(\Ga,F,\chi,W)$}}

\makeindex
\makeglossaries

\title{Equilibrium states in negative curvature}
\author{Fr\'ed\'eric Paulin \and Mark Pollicott \and Barbara Schapira}
\date{\today.} 

\begin{document} 
\bibliographystyle{../alphanum}
\maketitle
 
\begin{abstract} 
\noindent   
With their origin in thermodynamics and symbolic dynamics, Gibbs
measures are crucial tools to study the ergodic theory of the geodesic
flow on negatively curved manifolds. We develop a framework (through
Patterson-Sullivan densities) allowing us to get rid of compactness
assumptions on the manifold, and prove many existence, uniqueness and
finiteness results of Gibbs measures. We give many applications, to the
Variational Principle, the counting and equidistribution of orbit
points and periods, the unique ergodicity of the strong unstable
foliation and the classification of Gibbs densities on some Riemannian
covers.
\end{abstract}

\section{Introduction} 
\label{sec:intro} 

With their thermodynamic origin, Gibbs measures (or states) are very
useful in symbolic dynamics over finite alphabets (see for instance
\cite{Ruelle04,Zinsmeister96,Keller98}). Sinai, Bowen and Ruelle
introduced them in hyperbolic dynamics (which, via coding theory, has
strong links to symbolic dynamics), in particular in order to study
the weighted distribution of periodic orbits and the equilibrium
states for given potentials (see for instance \cite{Sinai72, Bowen75,
  BowRue75,ParPol90}). For instance, this allows the dynamical
analysis of the geodesic flow of negatively curved Riemannian
manifolds, provided its non-wandering set is compact. As a first step
to venture beyond the compact case, Gibbs measures in symbolic
dynamics over countable alphabets have been developed by Sarig
\cite{Sarig99,Sarig01a}. But since no coding theory which does not
lose geometric information is known for general non-compact manifolds,
this methodology is not well adapted to the non-compact case. Note that
a coding-free approach of transfer operators, in particular via
improved spectral methods, was put in place by Liverani, Baladi,
Gou\"ezel and others (see for instance, restricting the references to
the case of flows, \cite{Liverani04, ButLiv07, Tsuji10, Tsuji12,
  BalLiv12, GiuLivPol12}, as well as \cite{AviGou12} for an example in
a non locally homogeneous and non-compact situation). In this text, we
construct and study geometrically Gibbs measures for the geodesic flow
of negatively curved Riemannian manifolds, without compactness
assumption.

By considering the action on a universal Riemannian cover of its
covering group, this study is strongly related to the (weighted)
distribution of orbits of discrete groups of isometries of negatively
curved simply connected Riemannian manifolds. After work of Huber and
Selberg (in particular through his trace formula) in constant
curvature, Margulis (see for instance \cite{Margulis04}) made a
breakthrough, albeit in the unweighted lattice case, to give the
precise asymptotic growth of the orbits. Patterson \cite{Patterson76}
and Sullivan \cite{Sullivan79}, in the surface and constant curvature
case respectively, have introduced a nice approach using measures at
infinity, giving in particular a nice construction of the measure of
maximal entropy in the cocompact case. The Patterson-Sullivan theory
extends to variable curvature and non lattice case, and an optimal
reference (in the unweighted case) is due to Roblin \cite{Roblin03}.

After preliminary work of Hamenst\"adt \cite{Hamenstadt97} on Gibbs
cocycles, of Ledrappier \cite{Ledrappier95}, Coudène \cite{Coudene03}
and Schapira \cite{Schapira04a} using a slightly different approach,
and of Mohsen \cite{Mohsen07} (their contributions will be explained
as the story unfolds), the aim of this book is to develop an optimal
theory of Gibbs measures for the dynamical study of geodesic flows in
general negatively curved Riemannian manifolds, and of weighted
distribution of orbits of general discrete groups of isometries of
negatively curved simply connected Riemannian manifolds.

\medskip Let us give a glimpse of our results, starting by describing
the players. Let $M$ be a complete connected Riemannian manifold with
pinched negative curvature. Let $F:T^1M\ra \RR$ be a
H\"older-continuous map, called a {\it potential},\index{potential}
which is going to help us define the various weights. In order to
simplify the exposition of this introduction, we assume here that $F$
is invariant by the antipodal map $v\mapsto -v$ (see the main
body of the text for complete statements). Let $p:\wt M\ra M$ be a
universal Riemannian covering map, with covering group $\Ga$ and
sphere at infinity $\partial_\infty \wt M$, and let $\wt F=F\circ
p$. We make no compactness assumption on $M$: we only assume $\Ga$ to
be non-elementary (that is, non virtually nilpotent).  (In the body of
the text, we will also allow $M$ to be a good orbifold, hence $\Ga$ to
have torsion.)  Let $\phi= (\phi_t)_{t\in\RR}$ be the geodesic flow on
$T^1M$ and $\wt\phi=(\wt \phi_t)_{t\in\RR}$ the one on $T^1\wt M$. For
all $x,y\in\wt M$, let us define $\int_x^y\wt F=\int_0^{d(x,\,y)}\wt
F(\wt\phi_tv)\;dt$ where $v$ is the unit tangent vector at $x$ to a
geodesic from $x$ through $y$. For every periodic orbit $g$ of $\phi$,
let $\L_g$ be the Lebesgue measure along $g$ and $\int_gF=\L_g(F)$ the
{\it period}\index{period} of $g$ for the potential $F$.

We will study three numerical invariants of the weighted dynamics.

$\bullet$~ Let $x,y\in\wt M$, and $c>0$ large enough. The {\it critical
  exponent}\index{critical exponent} of $(\Ga,F)$ is an exponential
growth rate of the orbit points of $\Ga$ weighted by the potential
$F$:
$$
\delta_{\Ga,\,F}=\lim_{t\ra +\infty} \frac{1}{t}\ln
\;\sum_{\ga\in\Ga,\;t-c\leq d(x,\,\ga y)\leq t} 
\;e^{\int_x^{\ga y} \wt F}\;.
$$

$\bullet$~ For every $t\geq 0$, let $\Per(t)$ be the set of periodic
orbits of $\phi$ with length at most $t$, and $\Per'(t)$ its subset of
primitive ones. Let $W$ be a relatively compact open subset of $T^1M$
meeting the non-wandering set of $\phi$, and $c>0$ large enough. The
{\it Gurevich pressure}\index{Gurevich pressure} of $(M,F)$ is an
exponential growth rate of the closed geodesics of $M$ weighted by the
potential $F$:
$$
P_{Gur}(M,F)=\lim_{t\ra+\infty}\;\frac{1}{t}\ln
\sum_{g\in \Per(t)-\Per(t-c),\;g\cap W\neq \emptyset}\;e^{\int_g F}\;.
$$
The restriction for the periodic orbits under consideration to meet a
given compact set is necessary: for instance if $M$ is an infinite
cyclic cover of a compact manifold, and $F$ is invariant under the
cyclic covering group, then $T^1M$ has periodic orbits with the same
length and same period going out of every compact subset, and the left
hand side would otherwise be $+\infty$ for $t,c$ large enough.

$\bullet$ Let $\M$ be the set of $\phi$-invariant probability measures
on $T^1M$ and let $h_m(\phi)$ be the (metric, that is
measure-theoretic) entropy of the geodesic flow $\phi$ with respect to
$m\in\M$. The {\it topological pressure}\index{topological pressure}
of $(\phi,F)$ is the upper bound of the sum of the entropy and the
averaged potential of all states:
$$
P_{top}(\phi,F)=\sup_{m\in\M} \;\Big(h_m(\phi)+\int_{T^1M}F\,dm\Big)\;.
$$ 

When $F=0$, the quantity $\delta_{\Ga,\,F}$ is the usual critical
exponent $\delta_\Ga$ of $\Ga$, and $P_{top}(\phi,F)$ is the
topological entropy $h_{top}(\phi)$ of $\phi$. We prove in Subsection
\ref{subsec:loggrowth} and \ref{subsec:equalcritGur} that the limits
defining $\delta_{\Ga,\,F}$ and $P_{Gur}(M,F)$ exist (a result of
Roblin \cite{Roblin02} when $F=0$) and are independent of $x,y,W$ and
$c>0$ large enough (depending on $x,y,W$), and that, when
$\delta_{\Ga,\,F} >0$, the sums in their definition may be taken
respectively over $\{\ga\in\Ga,\; d(x,\ga y)\leq t\}$ and $\{g\in
\Per(t),\;g\cap W\neq \emptyset\}$.  

We have not assumed $\Ga$ to be finitely generated, so we prove in
Subsection \ref{subsec:semigroup} that $\delta_{\Ga,\,F}$ is the upper
bound of the critical exponents with potential $F$ of the free
finitely generated semi-groups in $\Ga$ (of Schottky type), a result
due to Mercat \cite{Mercat12} if $F=0$. For Markov shifts on countable
alphabets, the Gurevich pressure has been introduced by Gurevich
\cite{Gurevich69,Gurevich70} when the potential is equal to $0$, and
by Sarig \cite{Sarig99, Sarig01a, Sarig01b} in general. The Gurevich
pressure had not been studied much in a non symbolic non-compact
context before this work.

Our first result, generalising in particular results of Manning
\cite{Manning79} ($\delta_\Ga=h_{top}(\phi)$), Ruelle \cite{Ruelle81}
and Parry \cite{Parry88} when $M$ is compact, is the following one
(see Subsection \ref{subsec:equalcritGur} and Chapter
\ref{sec:variaprincip}).

\btheo \label{theo:triplequality}
We have
$$
P_{Gur}(M,F)=\delta_{\Ga,\,F}=P_{top}(\phi,F)\;.
$$
\etheo

\medskip Let us now define the Gibbs measures (see Chapter
\ref{sec:GPS}). We use Mohsen's extension (whose compactness
assumptions are not necessary) of the Patterson-Sullivan construction,
see \cite{Mohsen07}, which is more convenient than the one of
\cite{Ledrappier95,Coudene03,Schapira04a}.  The (normalised) {\it
  Gibbs cocycle}\index{Gibbs cocycle} of $(M,F)$ is the map $C$ from
$\partial_\infty\wt M\times\wt M\times \wt M$ to $\RR$ defined by
$$
(\xi,x,y)\mapsto C_\xi(x,y)= \lim_{t\ra+\infty}\int_y^{\xi_t}(\wt
F-\delta_{\Ga,\,F})-\int_x^{\xi_t}(\wt F-\delta_{\Ga,\,F})\;,
$$
where $t\mapsto \xi_t$ is any geodesic ray converging to a given
$\xi\in \partial_\infty\wt M$. When $F=0$, we have
$C_\xi(x,y)=\delta_\Ga \beta_\xi(x,y)$ where $\beta$ is the Busemann
cocycle. A (normalised) {\it Patterson density}%
\index{Patterson density} of $(\Ga,F)$ is a family of pairwise
absolutely continuous finite positive nonzero Borel measures
$(\mu_x)_{x\in\wt M}$ on the sphere at infinity $\partial_\infty\wt M$
such that, for all $\ga\in\Ga$, $x,y\in\wt M$ and for $\mu_y$-almost
every $\xi\in\partial_\infty \wt M$,
$$
\ga_*\mu_x=\mu_{\ga x}\;\;\;{\rm and}\;\;\;
\frac{d\mu_x}{d\mu_y}(\xi)=e^{-C_\xi(x,\,y)}\;.
$$
Mohsen has extended Patterson's construction of such a family, as well
as Sullivan's shadow lemma.

Fix $x_0\in\wt M$. For every $v\in T^1\wt M$, let $v_-,\pi(v),v_+$ be
respectively the point at $-\infty$, the origin, the point at
$+\infty$ of the geodesic line defined by $v$, and let $t$ be the
distance between $\pi(v)$ and the closest point to $x_0$ on this
line. Define a measure $\wt m_F$ on $T^1\wt M$ by
$$
d\wt m_F(v)= \frac{d\mu_{x_0}(v_-)}{e^{C_{v_-}(\pi(v),\,x_0)}}\;
\frac{d\mu_{x_0}(v_+)}{e^{C_{v_+}(\pi(v),\,x_0)}}\;dt\;.
$$
This measure on $T^1\wt M$ is invariant under $\Ga$ and $\wt\phi$,
hence induces a measure on $T^1M$ invariant under $\phi$, called the
{\it Gibbs measure}\index{Gibbs measure} (or {\it Gibbs
  state})\index{Gibbs state} of the potential $F$. The Gibbs measures
are Radon measures (that is, they are locally finite Borel positive
measures), but they are not always finite. 

In symbolic dynamics, the Gibbs measures are shift-invariant measures,
which give to any cylinder a weight defined by the Birkhoff sum of the
potential. Here, the analogs of the cylinders are the neighbourhoods
of a (pointed and oriented) geodesic line, defined by small
neighbourhoods of its two points at infinity, that we also weight
according to the Birkhoff integral of the potential in order to define
our Gibbs measures. In Subsection \ref{subsec:gibbsproperty}, we prove
that our Gibbs measures indeed satisfy the Gibbs property on the
dynamical balls.

\medskip We have a huge choice of potential functions. For instance,
we can define
$$
F^{\rm su}=-\,\frac{d}{dt}_{\mid t= 0}\log\operatorname{Jac}
\big({\phi_t}_{\mid W^{\rm su}(v)}\big)(v)\;,
$$
the negative of the pointwise exponential growth rate of the Jacobian
of the geodesic flow restricted to the strong unstable manifold (see
Chapter \ref{sec:Liouville}). When $M$ has dimension $n$ and constant
curvature $-1$, then $F^{\rm su}$ is constant, equal to $-(n-1)$. Some
references use the opposite sign convention on the potential to define
the topological pressure as $\sup_{m\in\M} \big(h_m(\phi)-\int_{T^1M}F\,
dm\big)$, hence the unstable Jacobian needs also to be defined with the
opposite sign as the one above.

When $M$ is compact, we recover, up to a scalar multiple, for $F=0$
(the Patterson-Sullivan construction of) the Bowen-Margulis measure,
and for $F=F^{\rm su}$ the Liouville measure. Using Gibbs measures is a way
of interpolating between these measures, and gives a wide range of
applications. 

We will prove (see Chapter \ref{sec:Liouville}) new results about when
the Liouville measure is a Gibbs measure. Note that there are many
examples of non-compact (with infinite volume) Riemannian covers of
compact Riemannian manifolds with constant sectional curvature $-1$
and with ergodic (hence conservative) geodesic flow for the Liouville
measure, see for instance \cite{Rees81}.

\btheo \label{theo:liougibbsintro} If $\wt M$ is a Riemannian cover of a
compact manifold, and if the geodesic flow $\phi$ of $M$ is
conservative with respect to the Liouville measure, then the
Liouville measure is, up to a scalar multiple, the Gibbs measure for
the potential $F^{\rm su}$ and 
$$
P_{Gur}(M,F^{\rm su})= \delta_{\Ga,\,F^{\rm su}} =
P_{top}(\phi, F^{\rm su})= 0\;.
$$ 
\etheo

\medskip Let us now state (a simplified version of) the main results
of our book. An {\it equilibrium state}\index{equilibrium state} for
the potential $F$ is a $\phi$-invariant probability measure $m$ on
$T^1M$ such that $h_m(\phi)+\int_{T^1M}F\,dm$ is equal to the
topological pressure $P_{top}(\phi,F)$ (hence realising the upper
bound defining it). The first result (see Chapter
\ref{sec:variaprincip}) says that the equilibrium states are the Gibbs
states (normalised to probability measures).

\btheo[Variational Principle]\label{theo:introvarprincip}%
\index{theorem@Theorem!Variational Principle} If $m_F$ is
finite, then $\frac{m_F}{\|m_F\|}$ is the unique equilibrium state.
Otherwise, there exists no equilibrium state.  
\etheo

When $\Ga$ is convex-cocompact, this result is due to Bowen and Ruelle
\cite{BowRue75}. Their approach, using symbolic dynamics, called the
{\it thermodynamic formalism},\index{thermodynamic formalism} is
described in many references, and has produced many extensions for
various situations (see for instance
\cite{Ruelle04,Zinsmeister96,Keller98}), under some compactness
assumption (except for the work of Sarig already mentioned). For
general $\Ga$ but when $F=0$, the result is due to Otal and Peigné
\cite{OtaPei04}, and we will follow their proof.

\bigskip Up to translating $F$ by a constant, which does not change
the Patterson density nor the Gibbs measure, we assume in the
following four results \ref{theo:introorbequidis} to
\ref{coro:intropercount}, that $\delta_{\Ga,\,F}>0$ (see the main body
of this book for the appropriate statements when $\delta_{\Ga,\,F}\leq
0$). The next result, of equidistribution of (ordered) pairs of Dirac
masses on an orbit weighted by the potential towards the product of
Patterson measures, is due to Roblin \cite{Roblin03} when $F=0$, and
we will follow his proof (as well as for the next three results). We
denote by $\gls{dirac}$ the unit Dirac mass at a point $z\in\wt M$.

\btheo [Two point orbital equidistribution] 
\label{theo:introorbequidis}
\index{theorem@Theorem!Two point orbital \\equidistribution}
If $m_F$ is finite and mixing, then we have, as $t$ goes to $+\infty$,
$$
\delta_{\Ga,\,F}\;\|m_F\|\;e^{-\delta_{\Ga,\,F}\; t}
\;\sum_{\ga\in\Ga\;:\;d(x,\ga y)\leq t} \;
e^{\int_{x}^{\ga y} \wt F}\;\D_{\ga^{-1}x}\otimes\D_{\ga y}\;\;\;
\stackrel{*}{\rightharpoonup}\;\;\;\mu_{y}\otimes\mu_{x}\;.
$$
\etheo

\bcoro[Orbital counting] \label{coro:introorbcount}%
\index{theorem@Theorem!Orbital counting} If $m_F$ is finite
and mixing, then, as $t$ goes to $+\infty$,
$$
\sum_{\ga\in\Ga\;:\;d(x,\ga y)\leq t} \; e^{\int_{x}^{\ga y} \wt
  F}\;\sim\;\frac{\|\mu_{y}\|\;\|\mu_{x}\|}
{\delta_{\Ga,\,F}\;\|m_F\|}\;e^{\delta_{\Ga,\,F}\;t} \;. 
$$ 
\ecoro

When $M$ is compact and $F=0$, this corollary is due to Margulis. When
$F=0$, it is due to Roblin in this generality. Theorem
\ref{theo:introorbequidis} also gives counting asymptotics of orbit
points $\ga x_0$ staying, as well as their reciprocal $\ga^{-1} x_0$,
in given cones, see Subsection \ref{subsec:growth} (improving, under
stronger hypotheses, the logarithmic growth result of Subsection
\ref{subsec:loggrowth}).

The next two results are due to Bowen \cite{Bowen72b} when $M$ is
compact (or even $\Ga$ convex-cocompact), see also \cite{ParPol83,
  Parry84}. But our proof is very different even under this
compactness assumption, in particular we avoid any coding theory,
hence any symbolic dynamics, and any zeta function or transfer
operator. They are due to Roblin \cite{Roblin03} when $F=0$ in this
generality. The assumptions in the second statement cannot be removed:
for instance if $M$ is an infinite cyclic cover of a closed manifold
(which is not geometrically finite), and $F$ is invariant under the
cyclic covering group, then the Gibbs measure is infinite, $T^1M$ has
infinitely many periodic orbits with the same length and same period,
and the left hand side is $+\infty$ for $t$ large enough. Due to the
rarity of non geometrically finite examples with finite and mixing
Bowen-Margulis measure, we do not know if the assumption that $\Ga$ is
geometrically finite can be removed.

\btheo [Equidistribution of periodic orbits] 
\label{theo:introperequidis}%
\index{theorem@Theorem!Equidistribution of periodic \\orbits}
If $m_F$ is finite and mixing, then, as $t$ goes to $+\infty$,
$$
\delta_{\Ga,\,F}\;e^{-\delta_{\Ga,\,F}\; t}
\;\sum_{g\in\Per'(t)} \;e^{\int_{g} F}\;\L_g\;\;\;
\stackrel{*}{\rightharpoonup}\;\;\;\frac{m_F}{\|m_F\|}\;.
$$
\etheo

\bcoro[Counting of periodic orbits] \label{coro:intropercount}%
\index{theorem@Theorem!Counting of periodic orbits} If $\Ga$ is
geometrically finite and if $m_F$ is finite and mixing, then, as $t$
goes to $+\infty$,
$$
\sum_{g\in\Per'(t)} \;e^{\int_{g} F}\;\sim\;
\frac{e^{\delta_{\Ga,\,F}\;t}}{\delta_{\Ga,\,F}\;t} \;. 
$$ 
\ecoro

Let us now describe additional results.  Consider the Poincar\'e series
of $(\Ga,F)$
$$
Q_{\Ga,\,F,\,x,\,y}(s)=\sum_{\ga\in\Ga} \;\; e^{\int_x^{\ga y} (\wt F-s)}\;,
$$
which, independently of $x,y\in\wt M$, converges if
$s>\delta_{\Ga,\,F}$ and diverges if $s<\delta_{\Ga,\,F}$. In Chapter
\ref{sec:GHTSappli}, extending works of Hopf, Tsuji, Sullivan, Roblin
when $F=0$ and following the proofs of \cite{Roblin03}, we give
criteria for the ergodicity and non ergodicity of the geodesic flow on
$T^1M$ endowed with a Gibbs measure. 

\btheo\label{theo:HTSRintro}
\index{theorem@Theorem!of Hopf-Tsuji-Sullivan-Roblin with potentials!ergodic}\index{Hopf-Tsuji-Sullivan-Roblin theorem with \\ potentials!ergodic}
The following conditions are equivalent
\begin{enumerate}
\item[(i)]  The Poincar\'e series of $(\Ga,F)$ diverges at 
  $s=\delta_{\Ga,\,F}$.
\item[(ii)] The conical limit set of $\Ga$ has positive measure with
  respect to $\mu_x$.
\item[(iii)] The dynamical system $(\partial_\infty\wt
  M\times\partial_\infty\wt M, \Ga, \mu_x\otimes \mu_x)$ is ergodic
  and conservative.
\item[(iv)] The dynamical system $(T^1M, \phi, m_F)$ is
  ergodic and conservative. 
\end{enumerate} 
\etheo

We discuss several applications in Subsection \ref{subsec:uniqueness}.
In particular, if the Poincar\'e series diverges at $s=
\delta_{\Ga,\,F}$, then the normalised Patterson density has no atom
and is unique up to a scalar multiple, and the Gibbs measure on $T^1M$
is unique up to a scalar multiple and ergodic when finite.

In view of (iii), it would be interesting to study if
$(\partial_\infty\wt M,\mu_x)$ is a (weak) $\Ga$-boundary in the sense
of Burger-Mozes \cite{BurMoz96} or Bader-Furman \cite{BadFur12,
  BadFurSha}.

In Subsection \ref{subsec:finiteness}, we give a criterion for the
finiteness of $m_F$ when $M$ is geometrically finite, extending the
one in \cite{DalOtaPei00} when $F=0$, as in Coudène
\cite{Coudene03}. The results \ref{theo:introorbequidis} to
\ref{coro:intropercount} require $m_F$ to be finite and mixing. But
given the finiteness assumption, the mixing one is mild. Indeed, by
Babillot's result \cite[Theo.~1]{Babillot02b}, since Gibbs measures
are quasi-product measures, if $m_F$ is finite, then $m_F$ is mixing
if and only if the geodesic flow $\phi$ is topologically mixing on
$T^1M$, a purely topological condition which is seemingly easier to
check and is conjecturally always satisfied.

It is indeed another interesting feature of Gibbs measures, that they
have strong ergodic properties. As soon as they are finite, they are
ergodic and even mixing as just said (when $\phi$ is topologically
mixing on $T^1M$). If in addition $P_{top}(\phi,F)>\sup F$ (for example
if $\sup F-\inf F < h_{top}(\phi)$), then the entropy of the Gibbs
measure $m_F$ is positive.  These properties are natural, since the
geodesic flow in negative curvature is a hyperbolic flow, and
therefore should have strong stochastic properties.  However,
generically (in Baire's sense), Coudène-Schapira \cite{CouSch13} have
announced that the invariant probability measures for $\phi$ are
ergodic, but not mixing and of entropy zero. Therefore, contrarily to
the generic case, the Gibbs measures provide an explicit family of
good measures reflecting the strong stochastic properties of the
geodesic flow.

\medskip
Consider the cocycle $c_{F}$ defined on the (ordered) pairs of vectors
$v,w\in T^1M$ in the same leaf of the strong unstable foliation
$\W^{\rm su}$ of the geodesic flow $\phi$ in $T^1M$ by
$$
c_{F}:(w,v)\mapsto \lim_{t\ra+\infty} \int_0^t
(F(\phi_{-s}w)-F(\phi_{-s}v)) \;ds\;.
$$
A nonzero family $(\nu_T)_T$ of locally finite (positive Borel)
measures on the transversals $T$ to $\W^{\rm su}$, stable by
restrictions, is {\it ${c}_{F}$-quasi-invariant}%
\index{quasi-invariant!transverse measure}%
\index{transverse measure!quasi-invariant} if for every holonomy map
$f:T\ra T'$ along the leaves of $\W^{\rm su}$ between two transversals
$T,T'$ to $\W^{\rm su}$, we have
$$
\frac{d\,f_*\nu_T}{d\,\nu_{T'}}(f(x))=e^{{c}_F(f(x),\,x)}\;.
$$
We will prove in Chapter \ref{sec:ergtheounistabfolia}, as in
\cite[Chap.~6]{Roblin03} when $F=0$ and in
\cite[Chap.~8.2.2]{Schapira03a} (unpublished), the following unique
ergodicity result.

\btheo \label{theo:uniergointro}%
\index{theorem@Theorem!Unique ergodicity of the strong unstable
  foliation} If $m_F$ is finite and mixing, then there exists, up to a
scalar multiple, one and only one ${c}_{F}$-quasi-invariant
family of transverse measures $(\mu_T)_T$ for $\W^{\rm su}$ such
that $\mu_T$ gives full measure to the set of points of $T$
which are negatively recurrent under the geodesic flow $\phi$ in
$T^1M$.  
\etheo

We give an explicit construction of $(\mu_T)_T$ in terms of the above
Patterson densities, and we deduce from Theorem
\ref{theo:uniergointro}, when $M$ is geometrically finite, the
classification of the ergodic ${c}_{F}$-quasi-invariant families of
transverse measures. This generalises works when $F=0$ of Furstenberg,
Dani, Dani-Smillie, Ratner for $M$ a surface, and of Roblin for
general $M$.

The main tool, which follows from Babillot's work, is the following
equidistribution result in $T^1M$ of the strong unstable leaves
towards the Gibbs measure $m_F$. For every $v\in T^1M$, let $W^{\rm
  su}(v)$ be its strong unstable leaf for $\phi$ and let $(\mu_{W^{\rm
    su}(v)})_{v\in T^1M}$ be the family of conditional measures on the
strong unstable leaves of the Gibbs measure $m_F$ associated to the
potential $F$. Note that conditional measures are in general defined
almost everywhere, and up to normalisation. Here, the nice structure
of the Gibbs measure allows to define them everywhere and not up to
normalisation (though they depend on the harmless normalisation of the
Patterson measures). Then for every relatively compact Borel subset
$K$ of $W^{\rm su}(v)$ (containing a positively recurrent vector in
its interior), for every $\psi\in\C_c(T^1M;\RR)$, we have
$$
\lim_{t\ra+\infty} \frac{1}{\int_K  e^{c_{F}(w,v)}
\;d\mu_{W^{\rm su}(v)}(w)}\;\int_K \psi\circ\phi_t(w)\;
e^{c_{F}(w,v)}\;d\mu_{W^{\rm su}(v)}(w)=
\frac{1}{\|m_F\|}\int_{T^1M}\psi \;dm_F\;.
$$
We also prove and use the sub-exponential growth of the mass of strong
unstable balls for the conditional measures on the strong unstable
leaves (we may no longer have polynomial growth as when $F=0$, though
we have no example in mind besides divergent leaves).

\medskip In the final chapter \ref{sec:gibbscover}, we extend the work
of Babillot-Ledrappier \cite{BabLed96}, Hamenst\"adt \cite{Hamenstadt02}
and Roblin \cite{Roblin05} when $F=0$ to study the Gibbs states on
Galois Riemannian covers of $M$.

Let $\chi:\Ga\ra\RR$ be a (real) character of $\Ga$.  For all $x,y\in
\wt M$, let us define the {\it (twisted) Poincar\'e series of
  $(\Ga,F,\chi)$}%
\index{Poincar\'e series!twisted by a character} by
$$
\gls{Poincarserietwi}(s)=\sum_{\ga\in\Ga} \;\;
e^{\chi(\ga)+\int_x^{\ga y} (\wt F-s)}\;.
$$
The {\it (twisted) critical exponent of $(\Ga,F,\chi)$}%
\index{critical exponent!twisted by a character} is the element
$\delta_{\Ga,\,F,\,\chi}$ in $\mathopen{[}-\infty,
+\infty\mathclose{]}$ defined by
$$
\gls{criticalexponenttwi}=\limsup_{n\ra +\infty}\;\frac{1}{n}\ln
\;\sum_{\ga\in\Ga,\;n-1< d(x,\,\ga y)\leq n} 
\;e^{\chi(\ga)+\int_x^{\ga y} \wt F}\;.
$$
When $F=0$ and $M$ is compact, since $H^1(M;\RR)\simeq
\operatorname{Hom}(\Ga,\RR)$, we recover the Poincar\'e series
associated with a cohomology class of $M$ introduced in
\cite{Babillot95}, and the critical exponent of $(\Ga,F,\chi)$ is then
called the {\it cohomological pressure}.%
\index{cohomological pressure}\index{pressure!cohomological}

A {\it (twisted) Patterson density}%
\index{Patterson density!twisted by a character} of dimension
$\sigma\in\RR$ for $(\Ga,F,\chi)$ is a family of finite nonzero
(positive Borel) measures $(\mu_{x})_{x\in\wt M}$ on $\partial_\infty
\wt M$ such that, for every $\ga\in\Ga$, for all $x,y\in\wt M$, for
every $\xi\in\partial_\infty \wt M$, we have
$$
\ga_*\mu_{x}=\;e^{-\chi(\ga)}\;\mu_{\ga x}\;,
$$
$$
\frac{d\mu_{x}}{d\mu_{y}}(\xi)=
e^{\lim_{t\ra+\infty}
\int_{x}^{\xi_t}(\wt F-\sigma)-\int_{y}^{\xi_t}(\wt F-\sigma)}\;,
$$
where $t\mapsto \xi_t$ is any geodesic ray converging to $\xi$.

Using the above-mentioned extension of the Hopf-Tsuji-Sullivan-Roblin
theorem (Theorem \ref{theo:HTSRintro}), as well as (an immediate
extension to the case with potential of) the Fatou-Roblin radial
convergence theorem for the ratio of the total masses of two Patterson
densities (see \cite[Théo.~1.2.2]{Roblin05}), we prove in Chapter
\ref{sec:gibbscover}, amongst other results, that the topological
pressure is unchanged by taking an amenable cover of $M$, and we give
a classification result of the Patterson densities on nilpotent covers
of $M$, for instance when $M$ is compact.

\btheo \label{theo:claspatdensnilpcovintro}
Let $\Ga'$ be a normal subgroup of $\Ga$, and $F':\Ga'
\bs T^1\wt M\ra \RR$ be the map induced by $\wt F$.

(1) If the quotient group $\Ga/\Ga'$ is amenable, then
$\delta_{\Ga,\,F}= \delta_{\Ga',\,F'}$.

(2) Let $\chi:\Ga\ra\RR$ be a character of $\Ga$. If
$\delta_{\Ga,\,F,\,\chi}<+\infty$ and $Q_{\Ga,\,F,\,\chi,\, x,\,y}(s)$
diverges at $s=\delta_{\Ga,\,F,\,\chi}$ (for instance if $\Ga$ is
convex-cocompact), then there exists a unique (up to a scalar
multiple) twisted Patterson density $\mu_{\Ga,\,F,\,\chi}=
(\mu_{\Ga,\,F,\,\chi, \,x}) _{x\in\wt M}$ of dimension
$\delta_{\Ga,\,F,\,\chi}$ for $(\Ga,F,\chi)$.

(3) If $\Ga$ is convex-cocompact and $\Ga/\Ga'$ is nilpotent, then the
set of ergodic Patterson densities for $(\Ga',F')$ is the set of
multiples of $\mu_{\Ga,\,F,\,\chi}$ for the characters $\chi$ of $\Ga$
vanishing on $\Ga'$.  
\etheo

\medskip We conclude this introduction by saying a few words
concerning again the motivations, the framework, and possible
developments.

\medskip As already alluded to, the main motivation comes from
statistical physics, including one-dimensional lattice gasses,
interactions, ensembles and transfer matrix, from the point of view of
equilibrium states for chaotic dynamics under general potentials,
pressure computations, weighted distribution of periods, etc. The
restriction to compact state spaces was starting to be
cumbersome. From a purely dynamical system point of view, using
potentials is a way to study the weighting of the dynamics: a subset
of the state space where the potential is very negative tends to have
a small measure with respect to the equilibrium state, by the
maximisation process involved in the Variational Principle. More
precisely, for all $v,w\in T^1M$, if $\phi^\RR w$ stays in a subset of
$T^1M$ where $F$ is small, and $\phi^\RR v$ stays in a subset of
$T^1M$ where $F$ is large, the ratio $m_F(B(v,T,r))/m_F(B(w,T,r))$ of
the masses for the equilibrium state of the potential $F$ of the
dynamical balls centred at $v,w$ should tend exponentially to
$0$. Furthermore, the family of Gibbs measures, when the potential is
varying, is a large family having excellent dynamical properties,
including a quasi-product structure with respect to the
stable/unstable manifolds, which provides a wide playground with many
possible applications, including multi-fractal analysis.

\medskip The framework of pinched negatively curved orbifolds, instead
of quotients $\Ga\bs X$ of $\operatorname{CAT}(-1)$ metric spaces $X$
as in \cite{Roblin03}, comes precisely from the necessity of
controlling the potential. Though the quotient space $\Ga\bs \G X$ of
the space $\G X$ of geodesic lines in $X$ is the natural replacement
for the state space which is unit tangent bundle, it is problematic to
define the potential $F$ on this replacement $\Ga\bs \G X$, if one
wants to be able to integrate $F$ along pieces of the orbits, as two
geodesic lines might share a nontrivial segment. The lower bound on
the curvature is also useful to control the H\"older structure on the
unit tangent bundle (see Lemma \ref{lem:dholdequivSasak}) and the
divergence on the potential on long pieces of orbits which are
initially close (see Lemma \ref{lem:exogeohyp} (ii)). We refer to
\cite{BroParPau13} where the authors introduce, in the case of metric
trees (which lie at the other end of the range of
$\operatorname{CAT}(-1)$ spaces), the new tools necessary to deal with
the potential in these singular state spaces.

\bigskip Here is, amongst many possible others, a list of research
themes connected to the above presentation, some of which the authors
plan to work on (see also \cite{ParPau13b} for applications to
counting geodesic arcs with weights, and \cite{PauPol13} for
Khintchine type results for the above equilibrium states).

\begin{enumerate}
\item It would be interesting to know when the Patterson densities
  associated to potentials (defined above and in Subsection
  \ref{subsec:gibbspattersondens}) are harmonic measures for random
  walks on an orbit of $\Ga$, see \cite{ConMun07} when $M$ is compact.
  It would also be interesting to extend to the case with nonzero
  potential the study of the spectral theory of Patterson densities,
  as done by Patterson and Sullivan in constant curvature (starting by
  studying the meromorphic extension of the Poincar\'e series associated
  to the potential, possibly using techniques of \cite{GiuLivPol12} or
  the recent microlocal approach of \cite{DyaZwo14}).
\item It would be interesting to study on which conditions on two
  potentials their associated Patterson densities belong to the same
  measure class (that is, are absolutely continuous with respect
  to each other). Random walk techniques often allow to construct
  singular measures (see for instance \cite{KaiLeP11, GadMahTio14}),
  hence this problem could be linked to the previous one.
\item
It would be interesting to study when our Gibbs measures (defined
above and in Subsection \ref{subsec:GibbsSulivanmeasure}) are harmonic
measures associated to elliptic second-order differential operators on
$M$ (see for instance \cite{Hamenstadt97b} when $M$ is compact).
\item 
  Recall (see for instance \cite[Chap.~16]{PolYur98}, \cite{AarNak00})
that a measure preserving flow $(\psi_t)_{t\in\RR}$ of a measured
space $(X,\mu)$ is {\it
  $2$-recurrent}\index{recurrent@$2$-recurrent}\index{multiple
  recurrence} if for every measurable subset $B$ in $X$ with
$\mu(B)>0$, there exists $n\in\NN-\{0\}$ such that
$$
\mu(\psi_{-n}B\cap B\cap \psi_nB)>0\;.
$$
(We do not know whether we may replace $n\in\NN-\{0\}$ by $t\in
[1,+\infty[$ and get an equivalent definition.) For every
$d\in\NN-\{0\}$, the notion of $d$-recurrence (replacing the above
centred formula by $\mu(B\cap \psi_{-n}B\cap\dots\cap
\psi_{-dn}B)>0$) has been introduced by Furstenberg in his proof of
Szemeredi's theorem, and the $1$-recurrence property is the
conservativity property (see Subsection \ref{subsec:GHTSR}).
Furstenberg proved that if $\mu$ is a probability measure, then
$(\psi_t)_{t\in\RR}$ is $d$-recurrent for every $d\in\NN-\{0\}$. But
when $\mu$ is infinite, there are examples of transformations which
are $1$-recurrent but not $2$-recurrent (see for instance
\cite{AarNak00}). It would be interesting to know when the geodesic
flow on $T^1M$ is multiply recurrent with respect to the Gibbs measure
associated to a given potential.  More precisely, let $N$ be a compact
connected negatively curved Riemannian manifold; let $N'=G'\bs \wt N$
where $G'$ is a (non necessarily normal) subgroup of the covering
group $G$ of a universal Riemannian cover $\wt N\ra N$; let $S$ be a
finite generating set of $G$; and let $X(G,S)$ be the Cayley graph of
$G$ with respect to $S$. It would be interesting to know when exactly
the geodesic flow of $N'$ is ergodic or $2$-recurrent, with respect to
the Liouville measure, or to any given Gibbs measure. We expect the
ergodicity to hold exactly when the appropriate (possibly simple for
the Liouville measure) random walk on the Schreier graph $G'\bs
X(G,S)$ of $G$ relative to $G'$ is recurrent. We expect the
$2$-recurrence to be more restrictive, as Aaronson and Nakada (private
communication) claim that the geodesic flow on the $\ZZ^d$-cover of a
compact connected hyperbolic manifold is $2$-recurrent for the
Liouville measure if and only if $d=0$ or $d=1$, and Rees proved in
\cite{Rees81} that it is ergodic (and hence conservative by
\cite[Prop.~1.2.1]{Aaronson97}) for the Liouville measure if and only
if $d=0,1,2$.
\item
It would be interesting to extend our work on equilibrium states and
Gibbs measures from the geodesic flow to the magnetic flows, seen as
a deformation family of the geodesic flow, as studied for instance by
Paternain and Grognet \cite{Grognet99}.
\item 
  It would be interesting to know under which conditions on $\Ga$
  we have $\delta_{\Ga,\,F^{\rm su}}\leq 0$ or, as when $M$ is
  compact, $\delta_{\Ga,\,F^{\rm su}}= 0$ (see Section
  \ref{sec:Liouville} for partial results).
\item As suggested by the referee, it would be interesting to
  generalise Lalley's result \cite{Lalley87}, using our techniques, to
  prove that a weighted closed geodesic chosen at random amongst the
  ones with lengths at most $T$ is close with high probability to the
  equilibrium state. This can be made precise by Large deviation
  results (see Y.~Kifer \cite{Kifer94}, as well as \cite{Pollicott95b}
  for some extension to Gibbs measures, both in the compact case), as
  well as by Central Limit theorems and multidimensional ones, that
  would be interesting to establish in our framework. The key point,
  if one wants to follow Kifer's steps, would be to study the
  upper-continuity of the pressure of the measures (in the non-compact
  case, which could be delicate, considering the countable state shift
  counterexamples).
\item It would be interesting to study the continuous dependence
  (using for instance Wasser\-stein-type distances, as in
  \cite{Villani09}) or higher regularity dependence (using the
  differential calculus of Ambrosio-Gigli-Savaré \cite{AmbGigSav08} or
  others) on the potential $F$ of its associated Gibbs measure
  $m_F$. In another direction, it would be interesting to know whether
  the integral against a fixed H\"older-continuous reference function of
  a Gibbs measure varies analytically with the potential, possibly
  using ideas of \cite{ButLiv07} and \cite{GiuLivPol12}. One could
  also deform the Riemannian metric, as was done in
  \cite{KatKniPolWei89, Contreras92}, both in the compact case.
\item It would be interesting to study the exponential decay of
  correlations for equilibrium states in our general framework (see
  the works of Dolgopyat, Liverani, Giulietti-Liverani-Pollicott
  \cite{GiuLivPol12}, and the first serious result in the (geometric)
  non-compact case of \cite{MohOh12}).
\item It would be interesting to study the ergodic theory of the Gibbs
  measures (starting with their existence, uniqueness and explicit
  construction !) for the Teichmüller flow on the quotient of the
  space $\Q^1(\Sigma)$ of unit-norm holomorphic quadratic
  differentials on a closed connected orientable surface $\Sigma$ of
  genus at least two by a non-elementary subgroup $\Ga$ of the mapping
  class group of $\Sigma$, for a potential $F$ on $\Q^1(\Sigma)$
  invariant under $\Ga$. Patterson-Sullivan measures (without
  potential) on Thurston's boundary of Teichmüller spaces have been
  constructed for instance by Hamenst\"adt \cite{Hamenstadt09} and
  Gekhtman \cite{Gekhtman14}. Note that, relating this last question
  to the first two ones, V.~Gadre \cite{Gadre13} has proved that the
  harmonic measures on Thurston's boundary of Teichmüller spaces
  (whose existence is due to Kaimanovich-Masur), coming from random
  walks with finite support on the mapping class groups, are singular
  with respect to the Lebesgue measure, see also \cite{GadMahTio14}.
\end{enumerate}

\bigskip \noindent {\small{\it Acknowledgments. } The origin of this
  work is a fruitful stay of the second author as invited professor at
  the Ecole Normale Supérieure in 2009 during the first author's term
  there. The first author thanks the University of Warwick for visits
  during which part of the work took place. We thank J.~Parkkonen for
  his corrections on a preliminary version of this text, J.~Buzzi for
  his comments, and M.~Brin, S.~Crovisier, Y.~Coudène, F.~Ledrappier,
  J.~Aaronson for their help with Chapter \ref{sec:Liouville}. The
  first author thanks the Mittag-Leffler Institute where Chapter
  \ref{sec:ergtheounistabfolia} was written. The third author
  benefited from the ANR grant JCJC10-0108 Geode. Most of all, we
  would like to thank the anonymous referee for her/his tremendous
  work on the first version of this book. By her/his extremely
  thorough reading (her/his report was 26 pages long), she/he has
  corrected and improved many results, allowing to remove many
  superfluous hypotheses, and providing complete proofs in several 
  instances. In particular, 

$\bullet$~ the statements
  \ref{lem:elemproppressure} (viii), \ref{theo:limsupenfaitlimbis},
  \ref{lem:appendixAreferee}, \ref{coro:loggrowthbis},
  \ref{prop:atomlessreferee}, \ref{lem:annulaireferee},
  \ref{lem:neglinonprimreferee}, \ref{rem:refereeComment19}, 

$\bullet$~ the remark (2) following Definition \ref{def:gibbsprop}, 
  and what it allows later on,

$\bullet$~ the
  deduction of \ref{maingrowthRoblinbis} from \ref{maingrowthRoblin};
  \ref{coro:corocinq} from \ref{coro:coroun}, \ref{coro:corodeux},
  \ref{coro:corotrois} and \ref{coro:coroquatre};
  \ref{theo:equidisperorbbis} from \ref{theo:equidisperorb};
  \ref{theo:equidisperorbetroitebis} from
  \ref{theo:equidisperorbetroite}, 

$\bullet$~ the removal of the boundedness
  assumption on the potential to obtain the equality between its
  pressure and its critical exponent (see Equation
  \eqref{eq:egalitevariationnelle}), and 

$\bullet$~ a major improvement of
  Theorem \ref{theo:LiouvilleisGibbs} (where the present
  conservativity assumption was originally replaced by Furstenberg's
  $2$-recurrence property), 

$\bullet$~ a complete reorganisation of
  Section \ref{sec:gibbscover},

  \noindent should be attributed to her/him.}

\newpage
\tableofcontents

\newpage
\section{Background on negatively curved 
manifolds} 
\label{sec:negacurvnot} 

In this text, the triple $(\wt M,\Ga,\wt F)$ will denote the following
data.

Let $\gls{Mtilde}$ be a complete simply connected Riemannian manifold,
with dimension at least $2$ and pinched sectional curvature $-b^2\leq
K\leq -1$, where $b\geq 1$. Let $\gls{Gamma}$ be a non-elementary (see
the definition of non-elementary below) discrete group of isometries of
$\wt M$.  Let $\gls{potentialtilde} :T^1\wt M\ra \RR$ be a
H\"older-continuous $\Ga$-invariant map, called a {\it
  potential}\index{potential} (see the definition of the
H\"older-continuity we will use in Subsection \ref{subsec:holdercont}).

We refer for instance to \cite{BalGroSch85,BriHae99} for
general background on negatively curved manifolds. We only recall
in this chapter the notation and results about them that we will use
in this text.

\bigskip \medskip
\noindent{\bf General notation. } Here is some general notation that
will be used in this text.

Let $A$ be a subset of a set $E$. We denote by
$\gls{characteristicfunction}:E\ra \{0,1\}$ the characteristic
(or indicator) function of $A$: $\mathbbm{1}_A(x)=1$ if $x\in A$, and
$\mathbbm{1}_A(x)=0$ otherwise. We denote by $\gls{complementary}=E-A$
the complementary subset of $A$ in $E$.

We denote by $\gls{ln}$ the natural logarithm (with $\ln (e)=1$). 

For every real number $x$, we denote by $\gls{valeurabsolue}$ the
integral part of $x$, that is the largest integer not greater that
$x$. 

Given two maps $f,g:[0,+\infty\mathclose{[}\ra \mathopen{]}0,
+\infty[\,$, we will write $f\asymp g$ if there exists $c>0$ such that
for every $t\geq 0$, we have $\frac{1}{c}\;f(t)\leq g(t)\leq c\;f(t)$.

For every smooth map $f:N\ra N'$ between smooth manifolds, we denote
by $\gls{tangentmap}:TN\ra TN'$ its tangent map. 

We denote by $\|\mu\|$ the total mass of a finite positive measure
$\mu$. 

If $(X,\A)$ and $(Y,\B)$ are measurable spaces, $f:X\ra Y$ a
measurable map, and $\mu$ is a measure on $X$, we denote by $f_*\mu$
the image measure of $\mu$ by $f$, with $f_*\mu(B)=\mu(f^{-1}(B))$ for
every $B\in\B$. 

For every integer $n\geq 2$, we denote by $\HH^n_\RR$ (any model of)
the real hyperbolic space of dimension $n$ (its sectional curvature is
constant with value $-1$).

Given a topological space $X$, we denote by $\C_c(X;\RR)$ the vector
space of continuous maps $f:X \ra \RR$ with compact support.

\subsection{Uniform local H\"older-continuity}
\label{subsec:holdercont} 

Recall that two distances $d$ and $d'$ on a set $E$ are (uniformly
locally) {\it H\"older-equivalent}\index{Holder@H\"older!equivalence}
if there exist $c,\epsilon>0$ and $\alpha\in \mathopen{]}0,
1\mathclose{]}$ such that for all $x,y\in E$ with $d(x,y)\leq
\epsilon$, we have
$$
\frac{1}{c}\;d(x,y)^{\frac{1}{\alpha}}\leq d'(x,y)\leq 
c\;d(x,y)^\alpha\;.
$$ This relation is an equivalence relation on the set of distances on
$E$, and a {\it H\"older structure}\index{Holder@H\"older!structure}
on $E$ is the choice of such an equivalence class. Note that two
H\"older-equivalent distances on a set $E$ induce the same topology on
$E$. A metric space will be endowed with the H\"older structure of its
distance. A map $f:E\ra E'$ between sets endowed with H\"older
structures is {\it
  H\"older-continuous}\index{Holder@H\"older!continuity} if for any
distances $d$ and $d'$ in the H\"older structures of $E$ and $E'$,
there exist $c,\epsilon>0$ and $\alpha\in\mathopen{]}0,1\mathclose{]}$
    such that for all $x,y\in E$ with $d(x,y)\leq \epsilon$, we have
$$
d'(f(x),f(y))\leq c\;d(x,y)^\alpha\;.
$$
We will call $c,\epsilon$ and $\alpha$ the {\it H\"older
  constants}\index{Holder@H\"older!constants} of $f$, even though they
are not unique and depend on the choices of the distances $d$ and $d'$
in their equivalence class.

A map $f:E\ra E'$ between sets endowed with H\"older structures is
{\it locally H\"older-continuous}\index{locally H\"older-continuous}
if for every $x\in E$, there exists a neighbourhood $U$ of $x$ such
that the restriction of $f$ to $U$ is H\"older-continuous.

A H\"older structure on a set $E$ is {\it
  natural}\index{Holder@H\"older!structure!natural} under a group of
bijections $G$ of $E$ if any element of $G$ is H\"older-continuous.

\medskip
\noindent{\bf Remarks. } Note that this notion of H\"older-continuity
is a global one, and that some other texts only require the above
centred inequality to hold, for some given $c,\alpha$ and for every
$x$, only if $y$ is close enough to $x$ (that is, in some ball of
centre $x$ whose radius might be smaller than $1$ and might depend on
$x$). When the topology on $E$ defined by $d$ is compact, these two
definitions coincide, but we are trying to avoid any compactness
assumption in this text. For instance in symbolic dynamics, the
potentials may be only required to have this second property (or even
less, see for instance \cite{Keller98}): since two sequences in the
same strong stable leaf are eventually equal, the Gibbs cocycle is
then well-defined. But the fact that the Gibbs cocycle in our case is
well defined will use this global H\"older-continuity (see Subsection
\ref{subsec:Gibbscocycle}).

\medskip
Let $f:X \ra Y$ be a map between metric spaces. If $X$ is geodesic,
then the definition of the H\"older-continuity of $f$ is equivalent to
its following strengthened version: $f$ is {\it
  H\"older-continuous}\index{Holder@H\"older!continuity} if and only
if there exist $c>0$ and $\alpha\in\mathopen{]}0,1\mathclose{]}$ such
    that for all $x,y\in E$ with $d(x,y)\leq 1$, we have
$$
d(f(x),f(y))\leq c\;d(x,y)^\alpha\;.
$$
We will use this definition throughout this text, whenever the source
space is geodesic.

Another consequence of the global H\"older-continuous property of $f$
is that $f$ then has at most linear growth: if $X$ is a geodesic
space, the above condition implies that for all $x,y$ in $X$ with
$d(x,y)\geq 1$, we have $d(f(x),f(y))\leq 3c\,d(x,y)$. This will turn
out to be useful in particular in Subsection
\ref{subsec:classifquasinvtransv}.

\subsection{Boundary at infinity, isometries and
 the Busemann cocycle}
\label{subsec:geombus}

We denote by $\gls{dMandTunwtM}$ or $d_{\wt M}$ the Riemannian
distance on $\wt M$. We denote by $\gls{boundary}$ the boundary at
infinity of $\wt M$.  We endow $\wt M\cup \partial_\infty \wt M$ with
the cone topology, homeomorphic to the closed unit ball of $\RR^n$,
where $n$ is the dimension of $\wt M$. For every $x\in\wt M$, recall
that the Gromov-Bourdon {\it visual distance}\index{visual distance}
$\gls{visualdist}$ on $\partial_\infty \wt M$ seen from $x$ (see
\cite{Bourdon95}) is defined by
\begin{equation}\label{eq:defidistvis}
d_x(\xi,\eta)=
\lim_{t\ra+\infty} e^{\frac{1}{2}(d(\xi_t,\,\eta_t)-d(x,\,\xi_t)-d(x,\,\eta_t))}\;,
\end{equation}
where $t\mapsto \xi_t,\eta_t$ are any geodesic rays ending at
$\xi,\eta$ respectively. By the triangle inequality, for all
$x,y\in\wt M$ and $\xi,\eta\in\partial_\infty \wt M$, we have
\begin{equation}\label{eq:lipequivdistvis}
e^{-d(x,\,y)}\leq \frac{d_x(\xi,\eta)}{d_y(\xi,\eta)}\leq e^{d(x,\,y)}\;.
\end{equation}
We endow $\partial_\infty \wt M$ with the H\"older structure which is
the equivalence class of any visual distance. It is natural under the
isometry group of $\wt M$ and its induced topology on $\partial_\infty
\wt M$ coincides with the cone topology.

For the following definitions, let $r>0$, $x\in\wt M \cup \partial
_\infty \wt M$, $A\subset \wt M$ and $B\subset \partial_\infty\wt
M$. We denote by $\gls{shadow}$ the {\it shadow of $A$ seen from
  $x$},\index{shadow} that is, the subset of $\partial_\infty \wt M$
consisting of the endpoints of the geodesic rays (if $x\in\wt M$) or
lines (if $x\in\partial_\infty \wt M$) starting from $x$ and meeting
$A$. Note that $\OOO$ stands for ``ombre''. We denote by $\gls{cone}$
the {\it cone on $B$ with vertex $x$}\index{cone} in $\wt M$, that is,
the union of the geodesic rays or lines starting from $x$ and ending
at $B$. We denote by $\gls{neighbourhood}$ the open $r$-neighbourhood
of $A$, that is, 
$$
\N_rA=\{x\in \wt M\;:\;d(x,A)<r\}\;.
$$  
We denote by $\gls{rinterior}$ the {\it open $r$-interior}\index{open
  $r$-interior} of $A$, that is,
$$
\N^-_rA= \{x\in \wt M\;:\;d(x,\,^c\!A)>r\}\;,
$$ 
which is contained in $A$.  These sets satisfy the following
equivariance properties: for every isometry $\ga$ of $\wt M$, we have
$$
\ga\OOO_xA=\OOO_{\ga x}(\ga A),\;\;\;\ga\C_xB=\C_{\ga x}(\ga B),
\;\;\;\ga\N_rA=\N_r(\ga A),\;\;\;\ga\N^-_rA=\N^-_r(\ga A)\;.
$$

For every isometry $\ga$ of $\wt M$, we denote by
$$
\gls{translationlength}= \inf_{x\in\wt M} d(x,\ga x)
$$ 
the {\it translation length}\index{translation!length} of $\ga$ on
$\wt M$. Recall that $\ga$ is {\it elliptic}\index{elliptic} if it has
a fixed point in $\wt M$, {\it parabolic}\index{parabolic} if non
elliptic and $\ell(\ga)= 0$, and {\it loxodromic}\index{loxodromic}
otherwise. In this last case, we denote by
$$
\gls{axetranslation}=\{x\in\wt M\;:\;d(x,\ga x)=\ell(\ga)\}
$$ 
the {\it translation axis}\index{translation!axis} of $\ga$, which is
isometric to $\RR$. 

The isometry group $\gls{isometrygroup}$ of $\wt M$ is endowed with
the compact-open topology. Let $\Ga'$ be a discrete group of
isometries of $\wt M$. Its {\it limit set},\index{limit set} which is
the set of accumulation points in $\partial_\infty\wt M$ of any orbit
of $\Ga'$ in $\wt M$, will be denoted by $\gls{lambdaga}'$, and the
convex hull in $\wt M$ of this limit set by $\gls{convexhull}'$.
Recall that $\Ga'$ is {\it non-elementary}\index{nonelementary} if
$\card(\Lambda\Ga')\geq 3$. Since $\wt M$ has pinched negative
curvature, the following assertions are equivalent (see for instance
\cite{Bowditch95}):

$\bullet$~ $\Ga'$ is  non-elementary,

$\bullet$~ $\Ga'$ contains a free subgroup of rank $2$,

$\bullet$~ $\Ga'$ is not virtually solvable,

$\bullet$~ $\Ga'$ is not virtually nilpotent.

\noindent 
Also recall that $\Ga'$ is {\it
  convex-cocompact}\index{convex-cocompact} if $\Ga'$ is non-elementary
and $\Ga'\bs\C\Lambda\Ga'$ is compact.

A loxodromic element $\ga\in\Ga'$ is {\it
  primitive}\index{primitive}\index{loxodromic!primitive} in $\Ga'$ if
there are no $\alpha\in\Ga'$ and $k\geq 2$ such that $\ga=\alpha^k$.
Every loxodromic element of $\Ga'$ is a power of a primitive element
of $\Ga'$.

If $\Ga'$ is non-elementary, then $\Lambda\Ga'$ is the closure of the
set of fixed points of loxodromic elements of $\Ga'$, and it is the
smallest nonempty, closed, $\Ga'$-invariant subset of
$\partial_\infty\wt M$. In particular, if $\Ga''$ is an infinite
normal subgroup of $\Ga'$, then 
$$
\Lambda\Ga''=\Lambda\Ga'\;.
$$

We will use the following dynamical properties of the action of $\Ga'$
on its limit set $\Lambda\Ga'$: if $\Ga'$ is non-elementary, then

$\bullet$~ the action of $\Ga'$ on $\Lambda\Ga'$ is minimal (that is,
every orbit is dense);

$\bullet$~ the diagonal action of $\Ga'$ on $\Lambda\Ga'\times
\Lambda\Ga'$ is transitive (that is, there exists a dense orbit);

$\bullet$~ the set of pairs of fixed points of the loxodromic elements
of $\Ga'$ is dense in $\Lambda\Ga'\times \Lambda\Ga'$.

\medskip
The {\it conical} (or {\it radial}) {\it limit set}%
\index{limit set!radial}\index{limit set!conical} $\gls{lambdaconga}'$
of $\Ga'$ is the set of points $\xi\in\partial_\infty \wt M$ such that
there exists a sequence of orbit points of $x$ under $\Ga'$ converging
to $\xi$ while staying at bounded distance from a geodesic ray ending
at $\xi$. We have $\Lambda_c\Ga'=\Lambda\Ga'$ if $\Ga'$ is
convex-cocompact. If $\Ga'$ is non-elementary, the conical limit set
$\Lambda_c\Ga'$ of $\Ga'$ is dense in its limit set $\Lambda\Ga'$, and
the point at infinity of a geodesic ray $\rho$ in $\wt M$ belongs to
$\Lambda_c\Ga'$ if and only if the image of $\rho$ in $M$ comes back
in some compact subset at times tending to $+\infty$.

A point $\xi$ in $\partial_\infty\wt M$ is a {\it Myrberg
  point}\index{Myrberg point} of $\Ga'$ (see for instance
\cite{Tukia94}) if for all $\eta\neq \eta'$ in $\Lambda\Ga$ and $x\in
\wt M$, there exists a sequence $(\ga_n)_{n\in\NN}$ in $\Ga'$ such
that $\lim_{n\ra\infty}\ga_nx=\eta$ and
$\lim_{n\ra\infty}\ga_n\xi=\eta'$. We will denote by
$\gls{lambdaMyrga}'$\index{limit set!Myrberg} the set of Myrberg
points of $\Ga'$. It is clearly a subset of $\Lambda_c\Ga'$, and even
a proper subset of $\Lambda_c\Ga'$, since the fixed points of the
loxodromic elements of $\Ga'$ are conical limit points but not Myrberg
points.

\medskip
In this text, we are not assuming $\Ga$ to be torsion free.

Note that the action of $\Ga$ on $\Lambda\Ga$ is not necessarily
faithful (viewing the real hyperbolic plane $\HH^2_\RR$ as the
intersection with the horizontal plane of the ball model of the real
hyperbolic space $\HH^3_\RR$, consider the subgroup generated by
a non-elementary discrete isometry group of $\HH^2_\RR$ and the
hyperbolic reflexion of $\HH^3_\RR$ with fixed point set $\HH^2_\RR$).
But the pointwise stabiliser $\operatorname{Fix}_\Ga (\Lambda\Ga)$ of
$\Lambda\Ga$ in $\Ga$ is a finite normal subgroup of $\Ga$.  

\blemm\label{lem:fixlimset}
The set of fixed points in $\Lambda\Ga$ of any element $\ga$ of $\Ga-
\operatorname{Fix}_\Ga (\Lambda\Ga)$ is a closed subset with empty
interior in $\Lambda\Ga$.
\elemm

\dem This is well known, but difficult to locate in a reference. We
may assume $\ga$ to be elliptic. If $\ga$ fixes a fixed point of a
loxodromic element $\alpha$, then $\ga$ pointwise fixes its
translation axis $\Axe_{\alpha}$: otherwise, $(\alpha^{-n}
\ga\alpha^n)_{n\in\NN}$ is a sequence of pairwise distinct elements of
$\Ga$ mapping any given point in $\Axe_{\alpha}$ at bounded distance from
itself, which contradicts the discreteness of $\Ga$.  Since the set of
pairs of endpoints of the translation axes of loxodromic elements of
$\Ga$ is dense in $\Lambda\Ga\times \Lambda\Gamma$, if $\ga$ pointwise
fixes an open subset of $\Lambda\Ga$, then $\ga$ is the identity on
$\Lambda\Ga$.  \cqfd

\medskip In particular, by Baire's theorem, the union of the sets of
fixed points in $\Lambda\Ga$ of the elements of $\Ga-
\operatorname{Fix}_\Ga (\Lambda\Ga)$ is a $\Ga$-invariant Borel
subset, whose complement is a dense $G_\delta$-set. We refer to Lemma
\ref{lem:fixomegaga} for the consequences of this lemma for  the fixed
points of the elements of $\Ga$ inside $T^1\wt M$.

\medskip The {\it Busemann cocycle}%
\index{Busemann cocycle}\index{cocycle!Busemann} is the map
$\beta: \partial_\infty\wt M\times \wt M\times\wt M\ra \RR$, defined
by
$$
(\xi,x,y)\mapsto 
\gls{busemann}=\lim_{t\ra+\infty} d(x,\xi_t)-d(y,\xi_t)\;,
$$
where $t\mapsto \xi_t$ is any geodesic ray ending at $\xi$.  For every
$\xi\in\partial_\infty \wt M$, the {\em horospheres centred at
  $\xi$}\index{horosphere}\index{centre!of a horosphere} are the level
sets of the map $y\mapsto \beta_\xi(y,x_0)$ from $\wt M$ to $\RR$, and
the (closed) {\em horoballs centred at
  $\xi$}\index{horoball}\index{centre!of a horoball} are its sub-level
sets, for any $x_0\in\wt M$.

\subsection{The geometry of the unit tangent bundle}
\label{subsec:unittanbun}

We denote by $\gls{basepointproj}:\gls{tangentspace} \ra \wt M$ the
unit tangent bundle of $\wt M$, and we endow $T^1\wt M$ with Sasaki's
Riemannian metric (see below).  We denote by $\gls{iota}:T^1\wt M \ra
T^1\wt M$ the flip map $v\mapsto -v$.

Consider the Riemannian orbifold $\gls{M}=\Ga\bs \wt M$ and let us
denote by $\gls{tangentspacebas}$ and $\gls{bitangentspacebas}$ the
quotient Riemannian orbifolds $\Ga\bs T^1\wt M$ and $\Ga\bs TT^1\wt
M$.  We denote by $\gls{potential}:T^1M\ra \RR$ and again by
$\gls{iota}:T^1 M \ra T^1 M$ and $\gls{basepointproj}:T^1 M\ra M$ the
quotient maps of $\wt F$, $\iota$ and $\pi$.

For every $v\in T^1\wt M$, we denote by $\gls{vminus}$ and
$\gls{vplus}$ the points at $-\infty$ and $+\infty$ of the geodesic
line $\ell_v:\RR\ra \wt M$ with $\dot{\ell_v}(0)=v$; we have $(\iota
v)_\pm=v_\mp$.

We denote by $\gls{omegatilde}$ (respectively $\gls{omegatildec}$) the
closed (respectively Borel) $\Ga$-invariant set of elements $v\in
T^1\wt M$ such that $v_-$ and $v_+$ both belong to $\Lambda\Ga$
(respectively $\Lambda_c\Ga$), and by $\gls{omega}$ and $\gls{omegac}$
their images by the canonical projection $T^1\wt M\ra T^1M$.

\blemm\label{lem:fixomegaga} The union $A$ of the sets of fixed points
in $\wt\Omega\Ga$ of the elements of $\Ga- \operatorname{Fix}_\Ga
(\Lambda\Ga)$ is $\Ga$-invariant, closed, and properly contained in
$\wt\Omega\Ga$.  
\elemm

\dem The invariance of $A$ is clear and its closeness follows by
discreteness. The set of fixed points in $\wt\Omega\Ga$ of a given
element of $\Ga- \operatorname{Fix}_\Ga (\Lambda\Ga)$ is closed, with
empty interior in $\wt\Omega\Ga$ (otherwise, its set of points at
$+\infty$ would have nonempty interior in $\Lambda\Ga$, which would
contradict Lemma \ref{lem:fixlimset}).  The set $A$ is hence
properly contained in $\wt\Omega\Ga$, since its complement is a dense
$G_\delta$-subset of $\wt\Omega\Ga$, by Baire's theorem.  
\cqfd

\medskip
For every $x\in\wt M$, let $\gls{iotasubx}:\partial_\infty \wt M\ra
\partial_\infty \wt M$ be the {\it antipodal map with respect to
  $x$},\index{antipodal map with respect to a point} that is the
involutive map $\xi\mapsto v_-$ where $v$ is the unique element of
$T^1\wt M$ such that $\pi(v)=x$ and $v_+=\xi$. Note that for every
isometry $\ga$ of $\wt M$ and for every $x$ in $\wt M$, we have 
$$
\iota_{\ga x}=\ga\circ\iota_x\circ \ga^{-1}\;.
$$

Let us now recall a parametrisation of $T^1\wt M$ in terms of the
boundary at infinity of $\wt M$. Let $\gls{boundarydeux}=
(\partial_\infty \wt M\times \partial_\infty \wt M) - \Delta$, where
$\Delta$ is the diagonal in $\partial_\infty \wt
M\times \partial_\infty \wt M$. For every base point $x_0$ in $\wt M$,
the space $T^1\wt M$ may be identified with $\partial_\infty^2\wt M
\times \RR$, by the map which maps a unit tangent vector $v$ to the
triple $(v_-, v_+,t)$ where $t$ is the algebraic distance on
$\ell_v(\RR)$ (oriented from $v_-$ to $v_+$) between $\ell_v(0)$ and
the closest point of $\ell_v(\RR)$ to $x_0$. This parametrisation
(called the {\it Hopf parametrisation} defined by $x_0$)\index{Hopf
  parametrisation} differs from the one defined by another base point
$x'_0$ only by an additive term on the third factor (independent of
the time $t$).

\medskip We end this subsection by giving some information on the
Sasaki metric on $T\wt M$ (see for instance
\cite[Chap.~IV]{Ballmann95}).  Recall that it has the following
defining properties: The vector bundle $TT\wt M\ra T\wt M$ has a
Whitney sum decomposition $TT\wt M=V\oplus H$ such that, for every
$v\in T\wt M$, the direct sum decomposition of the fibres
$$
T_vT\wt M=V_v\oplus H_v
$$
is orthogonal for Sasaki's scalar product, and if $\pi:T\wt M\ra \wt
M$ also denotes the full tangent bundle of $\wt M$,

$\bullet$~ $T\pi_{\mid H_v}:H_v\ra T_{\pi(v)}\wt  M$ is an isometric
linear isomorphism,

$\bullet$~ $V_v=\Ker T_v\pi=T_v(T_{\pi(v)}\wt M)=T_{\pi(v)}\wt M$
(equality as Euclidean spaces),

$\bullet$~ if $\nabla$ is the Levi-Civita connection of $\wt M$, if
$p_V:TT\wt M\ra V$ is the bundle map which is the linear projection
from $T_vT\wt M$ on $V_v$ parallel to $H_v$ for every $v\in T\wt M$,
then for every vector field $X:\wt M\ra T\wt M$ on $\wt M$, with
$TX:T\wt M\ra TT\wt M$ its tangent map, we have
$$
\nabla_vX=p_V\circ TX(v)\;.
$$

The map $\pi:T^1\wt M \ra \wt M$ is a Riemannian submersion with
totally geodesic fibres, and if $\ga$ is an isometry of $\wt M$, then
its tangent map, still denoted by $\ga$, is an isometry of $T^1\wt M$.
The flip map $\iota$ is also an isometry of the Sasaki metric. The
Riemannian distance $d_S$ induced by the Sasaki metric satisfies (see
\cite{Solorzano10} for generalisations)
\begin{equation}\label{eq:distSasaki}
d_S(v,w)=\inf_{\alpha}\sqrt{\ell(\alpha)^2+
\|\;(v\,-\parallel_\alpha \!w)\;\|^2}\;,
\end{equation}
where the lower bound is taken on the smooth paths $\alpha:
\mathopen{[}0,1\mathclose{]}\ra\wt M$ of length $\ell(\alpha)$ from
$\pi(w)$ to $\pi (v)$ and $\parallel_\alpha$ is the parallel transport
along $\alpha$. Note that for all $v,w\in T^1\wt M$

\begin{equation}\label{eq:comparpidistSasaki}
d_S(v,w)\geq d(\pi(v),\pi(w))\;.
\end{equation}

\medskip We will endow $T^1\wt M$ with its quite usual distance
$\gls{dMandTunwtM}=d_{T^1\wt M}$, defined as follows: for all $v,v'\in
T^1\wt M$,
\begin{equation}\label{eq:formdistT1wtM}
d(v,v')=\frac{1}{\sqrt{\pi}}\;
\int_\RR d(\pi(\phi_{t}v),\pi(\phi_{t}v'))\,e^{-t^2}\;dt\;.
\end{equation}
Note that for all $v\in T^1\wt M$ and $t\in\RR$, we have
$$
d(v,\phi_tv)=|t|\;,
$$ 
and for all $v,v'\in T^1\wt M$ and $\ga\in\Isom(\wt M)$, we have 
\begin{equation}\label{eq:invardantipod}
d(\ga \,v,\ga \,v')= d(v,v')\;\;\;{\rm and}\;\;\;
d(\iota \,v,\iota \,v')=d(v,v')\;.
\end{equation}

We will sometimes consider another distance $\gls{dMandTunwtMprime}
=d'_{T^1\wt M}$ on $T^1\wt M$, defined by (using the convexity of the
distance in $\wt M$ to obtain the second equality)
\begin{align}
\forall\;v,w\in T^1\wt M,\;\;\; d'_{T^1\wt M}(v,w) & =
\max_{s\in\mathopen{[}-1,0\mathclose{]}} 
d_{\wt M}\big(\pi(\phi_sv),\pi(\phi_sw)\big)
\label{eq:defidprime}
\\ & =\max\big\{d_{\wt M}\big(\pi(v),\pi(w)\big),
d_{\wt M}\big(\pi(\phi_{-1}v),\pi(\phi_{-1}w)\big)\big\}\;.
\nonumber
\end{align}
The distance $d'$ is (Lipschitz-)equivalent to the Riemannian distance
$d_S$ induced by the Sasaki metric on $T^1\wt M$, since $\wt M$ has
pinched negative sectional curvature (see for instance \cite[page
70]{Ballmann95}).  Note that the flip map $\iota$, which is an
isometry of the Sasaki metric, is hence Lipschitz for this metric, with
Lipschitz constant depending only on the bounds on the sectional
curvature.

\medskip The following result should be well known, but we provide a
proof in the absence of a precise reference.

\blemm \label{lem:dholdequivSasak} The distance $d$ on $T^1\wt M$ is
H\"older-equivalent to the Riemannian distance $d_S$ induced by
Sasaki's Riemannian metric.  
\elemm

In particular, when we will consider H\"older-continuous functions
$\wt F:T^1\wt M\ra \RR$, we may use any one of the three distances
$d_S$, $d=d_{T^1\wt M}$ and $d'=d'_{T^1\wt M}$, remembering that the
H\"older constants depend on this distance.

\medskip 
\dem 
It is sufficient to prove that the distances $d$ and
$d'$ on $T^1\wt M$ are H\"older-equivalent.

\smallskip\noindent
\begin{minipage}{8cm} ~~~ Let $v,w\in T^1\wt M$. Define
  $x_t=\pi(\phi_t v)$ and $y_t=\pi(\phi_t w)$ for every $t\in\RR$, as
  well as $\epsilon_0=d(x_0,y_0)$ and $\epsilon_{-1}=
  d(x_{-1},y_{-1})$.  Let $t\mapsto z_t$ be the geodesic ray
  parametrised by $\mathopen{[}-1,+\infty\mathclose{[}$ from $x_{-1}$
  to the point at infinity $w_+$.
\end{minipage}
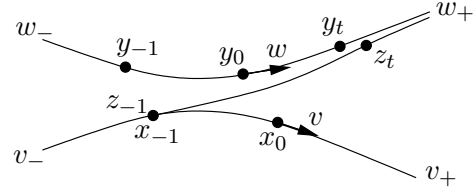
\begin{minipage}{6.9cm}
\begin{center}
\input{fig_comparsasaki.pstex_t}
\end{center}
\end{minipage}

\medskip By the definition of $d'$, for every $t\in\mathopen{[}-1,
0\mathclose{]}$, we have
$$
d(x_t,y_t)\leq d'(v,w)=\max\{\epsilon_0, \epsilon_{-1}\}\;.
$$
By convexity, for every $t\geq -1$, we have $d(z_t,y_t)\leq
\epsilon_{-1}$.  Since $\wt M$ has a finite lower bound on its
curvature, there exists a constant $c_1>0$ such that for every $t\geq
0$, we have $d(x_t,z_t)\leq d(x_0,z_0)\;e^{c_1\,t}$. Hence, by the
triangle inequality, for every $t\geq 0$, we have
$$
d(x_t,y_t)\leq \epsilon_{-1}+\epsilon_0\;e^{c_1\,t}\;.
$$
Similarly, considering the geodesic ray from $x_0$ to $w_-$, for every
$t\leq -1$, we have $d(x_t,y_t)\leq \epsilon_0+\epsilon_{-1}\;
e^{c_1\,(1-t)}$.  Therefore $\sqrt{\pi}\;d(v,w)$ is bounded from above
by
$$
\int_{-\infty}^{-1}(\epsilon_0+\epsilon_{-1}\; 
e^{c_1\,(1-t)})\;e^{-t^2}\;dt+
\int_{-1}^0\max\{\epsilon_0, \epsilon_{-1}\}\;e^{-t^2}\;dt+
\int_{0}^{+\infty}(\epsilon_{-1}+\epsilon_0\;e^{c_1\,t})\;e^{-t^2}\;dt\;,
$$
and there exists a constant $c_2>0$ such that $d(v,w)\leq
c_2\;(\epsilon_{-1}+\epsilon_0)\leq 2\,c_2\;d'(v,w)$.

Conversely, assume that $d'(v,w)\leq 1$. Then
$\epsilon_0,\epsilon_1\leq 1$ and
\begin{align*}
\sqrt{\pi}\;d(v,w)&\geq 
\int_{0}^{\epsilon_0/4}d(x_t,y_t)\;e^{-t^2}\;dt+
\int_{-1}^{-1+\epsilon_{-1}/4}d(x_t,y_t)\;e^{-t^2}\;dt\\ & \geq
\int_{0}^{\epsilon_0/4}(\epsilon_0-2t)\;e^{-t^2}\;dt+
\int_{-1}^{-1+\epsilon_{-1}/4}(\epsilon_{-1}-2(t+1))\;e^{-t^2}\;dt\\ & 
\geq 
\frac{\epsilon_0}{4}\times \frac{\epsilon_0}{2}\times \frac{1}{e}
\;+\;
\frac{\epsilon_{-1}}{4}\times \frac{\epsilon_{-1}}{2}\times \frac{1}{e}
\geq \frac{1}{16\,e}(\epsilon_0+\epsilon_{-1})^2
\geq\frac{1}{16\,e}(d'(v,w))^2\;.
\end{align*}
This proves the result.
\cqfd

\subsection{Geodesic flow, (un)stable foliations and 
the Hamenst\"adt distances}
\label{subsec:geodflowback}

We denote by $\big(\gls{geodesicflow}:T^1\wt M\ra T^1\wt M)_{t\in\RR}$
the geodesic flow of $\wt M$, which satisfies $\iota\circ\phi_t =
\phi_{-t}\circ\iota$. In the Hopf parametrisation, the geodesic flow is
the action of $\RR$ by translation on the third factor. We denote
again by $\gls{geodesicflow}:T^1 M \ra T^1 M$ the quotient map of
$\phi_t$.

Given $x\neq y$ in $\wt M$, the unit tangent vector at $x$ {\it
  pointing towards}\index{pointing towards} $y$ is the unique element
$v\in T^1_x\wt M$ such that there exists $t>0$ with $\pi(\phi_tv)=y$.

By \cite[Coro.~3.8]{Eberlein72}, the closed subset $\Omega\Ga$ defined
in the previous subsection is the (topological) {\it non-wandering
  set}\index{non-wandering set!(topological)} of the geodesic flow in
$T^1M$, that is, the set of points $v\in T^1M$ such that, for every
neighbourhood $U$ of $v$, there exists a sequence $(t_n)_{n\in\NN}$ in
$\RR$ converging to $+\infty$ such that, for every $n\in\NN$, the
intersection $U\cap \phi_{t_n}U$ is nonempty. The support of any Borel
probability measure on $T^1M$ which is invariant under the geodesic
flow is contained in $\Omega\Ga$. Note that $\Omega_c\Ga$ is dense in
$\Omega\Ga$, and is the {\it two-sided recurrent set}\index{two-sided
  recurrent set} of the geodesic flow in $T^1M$, that is, the set of
unit tangent vectors whose geodesic line comes back in some compact
subset of $M$ at times tending to $+\infty$ and to $-\infty$. In this
book, we will say that an element $v\in T^1M$ is positively
(respectively negatively) {\it recurrent}\index{recurrent} (under the
geodesic flow) if there exists a compact subset $K$ of $T^1M$ and a
sequence $(t_n)_{n\in\NN}$ in $\RR$ converging to $+\infty$ such that
$\phi_{t_n}v$ (respectively $\phi_{-t_n}v$) belongs to $K$.

\medskip
We now recall the definitions and properties of the dynamical
foliations associated to the geodesic flow on $T^1\wt M$.  For every
$v\in T^1\wt M$, we define the {\it strong stable
  leaf}\index{leaf!strong stable}\index{strong!stable leaf} of $v$,
{\it strong unstable leaf}\index{leaf!strong unstable}%
\index{strong!unstable leaf} of $v$, {\it (weak or central) stable
  leaf}\index{leaf!(weak or central) stable}%
\index{stable leaf} of $v$ and {\it (weak or central) unstable
  leaf}\index{leaf!(weak or central) unstable}\index{unstable!leaf} of
$v$ respectively by
$$
\gls{Wss}=\{w\in T^1\wt M\;:\; 
\lim_{t\ra+\infty}\;d(\phi_tv,\phi_tw)=0\}\;,
$$
$$
\gls{Wsu}=\{w\in T^1\wt M\;:\; 
\lim_{t\ra-\infty}\;d(\phi_tv,\phi_tw)=0\}\;,
$$
$$
\gls{Ws}=\{w\in T^1\wt M\;:\;\; \exists s\in\RR,\; 
\lim_{t\ra+\infty}\;d(\phi_{t+s}v,\phi_tw)=0\}\;,
$$
$$
\gls{Wu}=\{w\in T^1\wt M\;:\;\; \exists s\in\RR,\;
\lim_{t\ra-\infty}\;d(\phi_{t+s}v,\phi_tw)=0\}\;.
$$
Another set of notation is sometimes used, replacing $W^{\rm ss},
W^{\rm su}, W^{\rm s}$ and $W^{\rm u}$ by respectively $W^{\rm s},
W^{\rm s},W^{\rm cs}$ and $W^{\rm cu}$, which is easier to relate to direct
sums of the appropriate tangent spaces, but we will keep the above
one, which is quite traditional.

These four families define continuous foliations of $T^1\wt M$ with
smooth leaves, denoted by $\gls{Wssfol}$, $\gls{Wsufol}$,
$\gls{Wsfol}$, $\gls{Wufol}$ and called the {\it strong stable, strong
  unstable, stable, unstable foliations}\index{foliation!strong
  stable}\index{foliation!strong
  unstable}\index{foliation!stable}\index{foliation!unstable} of the
geodesic flow $(\phi_t)_{t\in\RR}$ on $T^1\wt M$. 

The smooth submanifolds $\pi(W^{\rm ss}(v))$ and $\pi(W^{\rm su}(v))$
of $\wt M$ are the horospheres passing through $\pi(v)$ with centres
respectively $v_+$ and $v_-$, and conversely $W^{\rm ss}(v)$ and
$W^{\rm su}(v)$ are the subsets of $T^1\wt M$ consisting respectively
of the inner and outer unit normal vectors to these two horospheres.
Note that $W^{\rm s}(v)$ and $W^{\rm u}(v)$ are invariant under the
geodesic flow $(\phi_t)_{t\in\RR}$.  Furthermore, for every $v\in
T^1\wt M$, the map $w\mapsto w_+$ from $W^{\rm su}(v)$ to
$\partial_\infty \wt M-\{v_-\}$ is a homeomorphism, and the map
$w\mapsto w_-$ from $W^{\rm ss}(v)$ to $\partial_\infty \wt M-\{v_+\}$
is also a homeomorphism.

For every $v\in T^1M$, denote again by $\gls{Wss}$ (respectively
$\gls{Wsu}$, $\gls{Ws}$, $\gls{Wu}$) the image by the canonical
projection $T^1\wt M\ra T^1M$ of $W^{\rm ss}(\wt v)$ (respectively
$W^{\rm su}(\wt v)$, $W^{\rm s}(\wt v)$, $W^{\rm u}(\wt v)$), where
$\wt v$ is any lift of $v$ to $T^1\wt M$. Note that $W^{\rm ss}(\iota
v) =\iota (W^{\rm su}(v))$, $W^{\rm su}(\iota v)=\iota (W^{\rm ss}
(v))$, $W^{\rm s}(\iota v)=\iota (W^{\rm u}(v))$, $W^{\rm u}(\iota v)
=\iota (W^{\rm s}(v))$ for every $v\in T^1\wt M$ or $v\in T^1 M$.  The
sets $W^{\rm ss}(v)$ (respectively $W^{\rm su}(v)$, $W^{\rm s}(v)$,
$W^{\rm u}(v)$) define a partition of $T^1M$, denoted by
$\gls{Wssfolbas}$ (respectively $\gls{Wsufolbas}$, $\gls{Wsfolbas}$,
$\gls{Wufolbas}$), which is a foliation outside the image in $T^1M$ of
the set of fixed points of elements of $\Ga-\{\id\}$ in $T^1\wt M$.

\medskip Let us also recall a natural family of distances on the
strong stable leaves of $T^1\wt M$ (see for instance
\cite[Appendix]{HerPau97}, slightly modifying the definition of
\cite{Hamenstadt89}, as well as \cite[\S 2.2]{HerPau10} for a
generalisation when a horosphere is replaced by the boundary of
any nonempty closed convex subset).  For every $w\in T^1\wt M$, let
$\gls{Hamdisu}$ be the {\it Hamenst\"adt distance}%
\index{Hamenst\"adt distance} on the strong unstable leaf of $w$,
defined as follows: for all $v,v'\in W^{\rm su}(w)$,
$$
d_{W^{\rm su}(w)}(v,v') = \lim_{t\ra+\infty} 
e^{\frac{1}{2}d(\pi(\phi_{t}v),\;\pi(\phi_{t}v'))-t}\;.
$$
This limit exists, and the Hamenst\"adt distance is a distance inducing
the original topology on $W^{\rm su}(w)$. For all $v,v'\in W^{\rm su}
(w)$ and for every isometry $\ga$ of $\wt M$, we have 
$$
d_{W^{\rm su}(\ga w)}(\ga v,\ga v')= d_{W^{\rm su}(w)}(v,v')\;.
$$  
For all $w\in T^1\wt M$, $s\in\RR$ and $v,v'\in W^{\rm su} (w)$, we
have
\begin{equation}\label{eq:dilatHamdist}
d_{W^{\rm su} (\phi_s w)}(\phi_s v,\phi_s v')=e^{s}\;d_{W^{\rm su}(w)}(v,v')\;.
\end{equation}
The following lemma compares the Hamenst\"adt distance $d_{W^{\rm su}
  (w)}$ with the distance $d$ on $\wt M$ and the distance $d$ on
$T^1\wt M$.

\blemm There exists $c>0$ such that, for all $w\in T^1\wt M$ and
$v,v'\in W^{\rm su} (w)$, we have
\begin{equation}\label{eq:compardistHamen}
\max\{\;\frac{1}{c}\; d(v,v'),\;d(\pi(v),\pi(v'))\;\}
\leq d_{W^{\rm su} (w)}(v,v')\leq 
e^{\frac{1}{2}d(\pi(v),\,\pi(v'))}\;.
\end{equation}
\elemm

\dem We give proofs only for the sake of completeness. The inequality
on the right hand side follows from the triangle inequality.
 
The inequality $d(v,v')\leq c\; d_{W^{\rm su} (w)}(v,v')$ for some
universal constant $c>0$ is due to \cite[Lem.~3]{ParPau13a}, as
follows.

\smallskip\noindent
\begin{minipage}{8.7cm} ~~~ We may assume that $v\neq v'$. Let
  $x_t=\pi(\phi_tv)$ and $x'_t= \pi(\phi_tv')$. By the convexity
  properties of the distance in $\wt M$, the map from $\RR$ to $\RR$
  defined by $t\mapsto d(x_t,x'_t)$ is increasing, with image
  $\mathopen{]}0,+\infty\mathclose{[}$.  Let $S\in\RR$ be such that
  $d(x_S, x'_S)=1$.  For every $t\ge S$, let $p$ and $p'$ be the
  closest point projections of $x_S$ and $x'_S$ respectively on the
  geodesic segment $\mathopen{[}x_t,x'_t\mathclose{]}$.
\end{minipage}
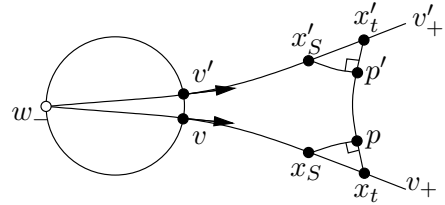
\begin{minipage}{6.2cm}
\begin{center}
\input{fig_comparHamd.pstex_t}
\end{center}
\end{minipage}

\medskip 
In the hyperbolic upper half-plane $\HH^2_\RR$, the
points $i$, $-1+2i$ and $1+2i$ are pairwise at distance
$2\log\frac{1+\sqrt{5}}{2} <1$. Hence by convexity, if $x,x'$ are the
points of $\HH^2_\RR$ with horizontal coordinates $1, -1$ respectively,
with same vertical coordinate and at hyperbolic distance $1$, then the
distance from $x,x'$ to the geodesic line with endpoints $1,-1$ is at
most $1$.  By comparison, we therefore have $d(p,x_S),\;d(p',x'_S)\leq
1$. Hence, by convexity and the triangle inequality,
\begin{align*}
d(x_t,x'_t)&\geq d(x_t,p)+d(p',x'_t)\\& 
\geq d(x_t,x_S)-1+ d(x'_t,x'_S)-1=2(t-S-1)\;.
\end{align*}
Thus by the definition of the Hamenst\"adt distance $d_{W^{\rm su}(w)}$,
we have
\begin{equation}\label{eq:minodistsu}
d_{W^{\rm su}(w)}(v,v')\geq e^{-S-1}\;.
\end{equation}

By the triangle inequality, if $t\geq S$, then
$$
d(x_t,x'_t)\leq
d(x_t,x_S)+d(x_S,x'_S)+ d(x'_S,x'_t)=2(t-S)+1\;.
$$
Since $\wt M$ is $\operatorname{CAT}(-1)$, if $t\leq S$, we have by
comparison
$$
d(x_t,x'_t)\leq e^{t-S}\,d(x_S,x'_S)=e^{t-S}\;.
$$
Therefore, by the definition of the distance $d$ on $T^1\wt M$ (see
Equation \eqref{eq:formdistT1wtM}), we have
$$
d(v,v')\leq \int_{-\infty}^{S}e^{t-S}\,e^{-t^2}\;dt+
\int_{S}^{+\infty}(2(t-S)+1)\,e^{-t^2}\;dt=\operatorname{O}(e^{-S})\;.
$$
The result hence follows from Equation \eqref{eq:minodistsu}.

\medskip
The last inequality $d(\pi(v),\pi(v')) \leq d_{W^{\rm su} (w)}(v,v')$ is
due to \cite[\S 2]{ParPau13b}, as follows.

\smallskip\noindent
\begin{minipage}{8cm} ~~~ Let $x=\pi(v)$, $x'=\pi(v')$ and $\rho=
  d_{W^{\rm su}(w)}(v,v')$.  Consider the ideal triangle $\Delta$ with
  vertices $v_+,v'_+$ and $v_-=v'_-$.  Let $p\in\mathopen{]}v_-,
  v_+\mathclose{[}\,$, $p'\in\mathopen{]}v'_-,v'_+\mathclose{[}$ and
  $q\in \mathopen{]}v_+,v'_+\mathclose{[}$ be the pairwise tangency
  points of horospheres centred at the vertices of $\Delta$:
  $\beta_{v_-}(p,p')=0$, $\beta_{v_+}(p,q) =0$ and $\beta_{v'_+}(p',q)
  =0$.  By definition of the Hamenst\"adt distance, the algebraic
  distance from $x$ to $p$ on the geodesic line
  $\mathopen{]}v_-,v_+\mathclose{[}$ (oriented from $v_-$ to $v_+$) is
  $-\ln \rho$.
\end{minipage} 
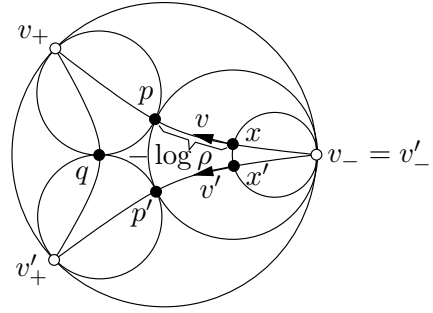
\begin{minipage}{6.9cm}
\begin{center}
\input{fig_horotang.pstex_t}
\end{center}
\end{minipage}

\medskip \noindent
\begin{minipage}{9cm}~~~ Consider the ideal triangle
  $\overline{\Delta}$ in the hyperbolic upper half-plane $\HH^2_\RR$,
  with vertices $-\frac{1}{2}$, $\frac{1}{2}$ and $\infty$. Let
  $\overline{p}=(-\frac 12,1)$, $\overline{p}'=(\frac 12,1)$ and
  $\overline{q}=(0,\frac{1}{2})$ be the pairwise tangency points of
  horospheres centred at the vertices of $\overline{\Delta}$. Let
  $\overline{x}$ and $\overline{x}'$ be the point at algebraic
  (hyperbolic) distance $-\ln \rho$ from $\overline{p}$ and
  $\overline{p}'$, respectively, on the upwards oriented vertical line
  through them. By comparison, we have
  $d(x,x')\leq d(\overline{x},\overline{x}')\leq 1/e^{-\ln
    \rho}=\rho$. \cqfd
\end{minipage} 
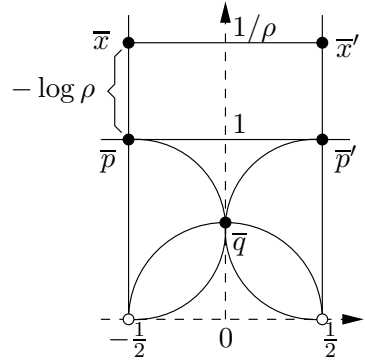
\begin{minipage}{5.9cm}
\begin{center}
\input{fig_triangideal.pstex_t}
\end{center}
\end{minipage}

\bigskip We define analogously the {\it Hamenst\"adt
  distance}\index{Hamenst\"adt distance} $\gls{Hamdiss}$ on the
strong stable leaf of $w\in T^1\wt M$ by
$$
d_{W^{\rm ss}(w)}(v,v') = \lim_{t\ra+\infty} 
e^{\frac{1}{2}d(\pi(\phi_{-t}v),\;\pi(\phi_{-t}v'))-t}\;.
$$
It satisfies analogous properties, in particular that for all $w\in
T^1\wt M$, $s\in\RR$ and $v,v'\in W^{\rm ss} (w)$, we have
\begin{equation}\label{eq:contractHamdist}
  d_{W^{\rm ss} (\phi_s w)}(\phi_s v,\phi_s v')=
e^{-s}\;d_{W^{\rm ss}(w)}(v,v')\;.
\end{equation}
The Hamenst\"adt distances on the strong stable and strong unstable
leaves are related as follows: for all $w\in T^1\wt M$ and $v,v'\in
W^{\rm su}(w)$, we have
\begin{equation}\label{eq:disthamantipod}
d_{W^{\rm ss}(\iota \,w)}(\iota \,v,\iota \,v')=d_{W^{\rm su}(w)}(v,v')\;.
\end{equation}

\medskip We denote by $\gls{suballHam}$ the open ball of center $v$
and radius $r$ for the Hamenst\"adt distance $d_{W^{\rm su}(v)}$ on
$W^{\rm su}(v)$, and similarly for $\gls{ssballHam}$. They are in
general different from the open balls of centre $v$ and radius $r$ for
the induced Riemannian distances, that we denote by respectively
$\gls{suballriem}$ and $\gls{ssballriem}$. The reader is advised not
to mix the two notations: in general variable curvature, the
Riemannian balls do not have the exact scaling property that the
Hamenst\"adt balls satisfy (see Equation \eqref {eq:dilatboulhamen} and
\eqref{eq:contractboulhamen}), hence using the second ones is
preferable in most of the cases. We will only use the Riemannian balls
in Subsection \ref{subsec:proofvaraprincip}, for the purpose of
staying close to some reference.  The distances $d_{W^{\rm su} (v)}$
and balls $B^{\rm su}(v,r)$ were denoted by $d_{H_+(v)}$ and
$B^+(v,r)$ in \cite{Roblin03} and \cite{Schapira04a}. For all $v\in
T^1\wt M$, $r>0$, $\ga\in\Ga$ and $t\in\RR$, by Equation
\eqref{eq:dilatHamdist}, we have
$$
\ga B^{\rm su}(v,r)=B^{\rm su}(\ga v,r)\;,
$$
and
\begin{equation}\label{eq:dilatboulhamen}
\phi_t\big(B^{\rm su}(v,r)\big)=B^{\rm su}(\phi_tv,e^tr)\;.
\end{equation}
Similarly, for all $\ga\in\Ga$ and $t\in\RR$, the open balls
$B^{\rm ss}(v,r)$ of centre $v\in T^1\wt M$ and radius $r>0$ for the
Hamenst\"adt distance $d_{W^{\rm ss}(v)}$ on the strong stable leaves
$W^{\rm ss}(v)$ satisfy $\ga B^{\rm ss}(v,r)=B^{\rm ss}(\ga v,r)$ and
\begin{equation}\label{eq:contractboulhamen}
\phi_t\big(B^{\rm ss}(v,r)\big)=B^{\rm ss}(\phi_tv,e^{-t}r)\;.
\end{equation}
By Equation \eqref{eq:disthamantipod}, we have
\begin{equation}\label{eq:ballhamantipod}
\iota\, B^{\rm su}(v,r)=B^{\rm ss}(- v,r)\;\;{\rm and}\;\;
\iota\, B^{\rm ss}(v,r)=B^{\rm su}(- v,r)\;.
\end{equation}

\subsection{Some exercises in hyperbolic geometry}
\label{exohypgeom}

We end this chapter by the following (well known, see for
instance \cite{ParPau10GT}) exercises in hyperbolic geometry.

\blemm \label{lem:exogeohyp} (i) For all $x,y,z$ in $\wt M$ such that
$d(x,z)=d(y,z)$, for every $t\in\mathopen{[}0, d(x,z)\mathclose{]}$,
if $x_t$ (respectively $y_t$) is the point on $\mathopen{[}x,
z\mathclose{]}$ (respectively $\mathopen{[}y,z\mathclose{]}$) at
distance $t$ from $x$ (respectively $y$), then
$$
d(x_t,y_t)\leq e^{-t}\sinh d(x,y)\;.
$$
(ii) For every $\alpha>0$, there exists $\beta>0$ (depending only on
$\alpha$ and the bounds on the sectional curvature of $\wt M$) such
that for all $x,y,z$ in $\wt M$ such that $d(x,z)=d(y,z)\geq 2$ and
$d(x,y)\leq \alpha$, for every $t\in\mathopen{[}0,
d(x,z)\mathclose{]}$, with $x_t$ and $y_t$ defined as above,
$$
d(x_t,y_t)\leq \beta\,d(x_1,y_1)\,e^{-t}\;.
$$
\elemm

\dem (i) We may assume that $x\neq y$, otherwise $x_t=y_t$ and the
result is trivially true. Since the sectional curvature of $\wt M$ is
at most $-1$, and by comparison, we may assume that $\wt M$ is the
hyperbolic upper half-space $\HH^2_\RR$. Let $m$ be the orthogonal
projection of $z$ on the geodesic line through $x$ and $y$, which is
the midpoint of $\mathopen{[}x,y\mathclose{]}$ by symmetry.

Replacing $z$ by the point at infinity of the geodesic ray starting
from $m$ and passing through $z$ (any one perpendicular to
$\mathopen{[}x,y\mathclose{]}$ if $m=z$) increases $d(x_t,y_t)$ and
does not change $d(x,y)$. Hence we may assume that $z$ is the point at
infinity in $\HH^2_\RR$, and that $x$ and $y$ lie on the (Euclidean)
circle of centre $0$ and radius $1$, have same (Euclidean) height,
with the real part of $x$ positive, the real part of $y$ negative.

\medskip\noindent
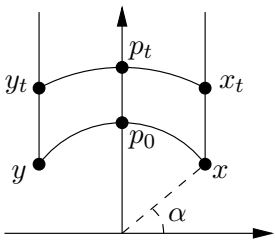
\begin{minipage}{6cm}
\begin{center}
\input{fig_expopinch.pstex_t}
\end{center}
\end{minipage}
\begin{minipage}{8.9cm} ~~~ Let $p_t$ be the intersection point of
  $\mathopen{[}x_t, y_t\mathclose{]}$ with the vertical axis, which is
  by symmetry the midpoint of $\mathopen{[}x_t,y_t\mathclose{]}$.  If
  $\alpha$ is the (Euclidean) angle at $0$ between the ho\-ri\-zontal
  axis and the (Euclidean) line from $0$ passing through $x$, then an
  easy computation in hyperbolic geometry (see also \cite{Beardon83},
  page 145) gives $\sinh d(x,p_0) =\cos \alpha/\sin\alpha$.
\end{minipage}

\smallskip\noindent 
Similarly, $\sinh d(x_t,p_t) =\cos \alpha/(e^t\sin\alpha)$. So that
$$
d(x_t,y_t)=2\,d(x_t,p_t)\leq 2\sinh d(x_t,p_t)=2\,e^{-t}\sinh d(x,p_0)
=2\,e^{-t}\sinh\frac{d(x,y)}{2}\;,
$$
which proves the first result, as $2\sinh\frac{t}{2}\leq \sinh t$ if
$t\geq 0$.
 
\medskip (ii) We start by the observation that for all $a>0$ and
$s\in\mathopen{[}0, a\mathclose{]}$, we have $s\leq\sinh s\leq
s\,\frac{\sinh a}{a}$.

By (i), we hence only have to prove that there exists a constant $c>0$
such that $d(x_1,y_1)\geq c\,d(x,y)$.  By comparison, we may assume
that $\wt M$ is the real hyperbolic plane with constant curvature
$-b^2$. By the hyperbolic sine formula in this plane, we have
$$
\frac{\sinh\big(b\;\frac{d(x_1,\,y_1)}{2}\big)}
{\sinh\big(b\;d(x_1,z)\big)}=\frac{\sinh\big(b\;\frac{d(x,\,y)}{2}\big)}
{\sinh\big(b\;d(x,z)\big)}\;.
$$ Since the map $s\mapsto \frac{\sinh t}{\sinh (t+b)}$ is non
decreasing as $b>0$, since $d(x,z)=d(x_1,z)+1\geq 2$, and since
$d(x_1,y_1)\leq d(x,y)\leq \alpha$ by convexity, we hence have, by
the preliminary observation, that
$$
d(x_1,y_1)\geq \frac{\sinh b}{\sinh (2\,b)}\;
\frac{\frac{b\,\alpha}{2}}{\sinh \frac{b\,\alpha}{2}}\;d(x,y)\;,
$$
as required. \cqfd

\bcoro\label{coro:controlexpopinchunit} (i) For every
$T\in\mathopen{[}0,+\infty\mathclose{[}\,$, for all $v$ and $w$ in
$T^1\wt M$ such that $\pi(\phi_Tv)=\pi(\phi_Tw)$, for every
$t\in\mathopen{[}0, T\mathclose{]}$,
$$
d'_{T^1\wt M}(\phi_tv,\phi_tw)\leq 
e^{-t}\sinh\big(2+d(\pi(v),\pi(w))\big)\;.
$$

(ii) For every $\alpha'>0$, there exists $\beta'>0$ (depending only on
$\alpha'$ and the bounds on the sectional curvature of $\wt M$) such
that for every $T\in\mathopen{[}1,+\infty\mathclose{[}\,$, for all $v$
and $w$ in $T^1\wt M$ such that $\pi(\phi_Tv)=\pi(\phi_Tw)$ and
$d(\pi(v),\pi(w))\leq \alpha'$, for every
$t\in\mathopen{[}0,T\mathclose{]}$,
$$
d'_{T^1\wt M}(\phi_tv,\phi_tw)\leq 
\beta'\,d(\pi(v),\pi(w))\,e^{-t}\;.
$$
\ecoro

\dem (i) By the definition of the distance $d'_{T^1\wt M}$ (see
Equation \eqref{eq:defidprime}), this follows from Lemma
\ref{lem:exogeohyp} (i), applied to $x=\pi(\phi_{-1}v)$,
$y=\pi(\phi_{-1}w)$ and $z=\pi(\phi_Tv)$, since $d(x,y)\leq
d(\pi(v),\pi(w))+2$.

(ii) This follows similarly from Lemma \ref{lem:exogeohyp} (ii), with
$\alpha=\alpha'+2$, since with $x,y,z$ as above and $x_1,y_1$ as in
this lemma, we have $x_1=\pi(v), y_1=\pi(w)$ and
$d(x,z)=d(y,z)=T+1\geq 2$.  \cqfd

\medskip Also recall without proof the following well known results in
complete simply connected Riemannian manifolds with sectional
curvature at most $-1$ (variants of Anosov's closing lemma).

\blemm \label{lem:quasigeodpresquegeod} For all
$\ell,\epsilon,\theta>0$, there exists $\theta_0= \theta_0(\ell,
\epsilon, \theta)>0$ such that for every piecewise geodesic path
$\omega= \bigcup_{0\leq i<N}\mathopen{[}x_i,x_{i+1}\mathclose{]}$ in
$\wt M$ with $N\in\NN\cup \{+\infty\}$, if its exterior angle at $x_i$
is at most $\theta_0$ for $0<i<N-1$ and if the length of each segment
$\mathopen{[}x_i,x_{i+1}\mathclose{]}$ is at least $\ell$, then

$\bullet$~ $\omega$ is contained in the $\epsilon$-neighbourhood of
the geodesic segment $\omega'=\mathopen{[}x_0,x_{N}\mathclose{]}$, if
$N$ is finite, or of the geodesic ray $\omega'=\mathopen{[}x_0,
\xi\mathclose{]}$ where $\xi\in\partial_\infty \wt M$ is the point at
infinity to which the sequence $(x_i)_{i\in\NN}$ converges, if $N$ is
infinite;

$\bullet$~ the angles between $\omega$ and $\omega'$
at the finite endpoints of $\omega'$ are at most $\theta$. 
\elemm

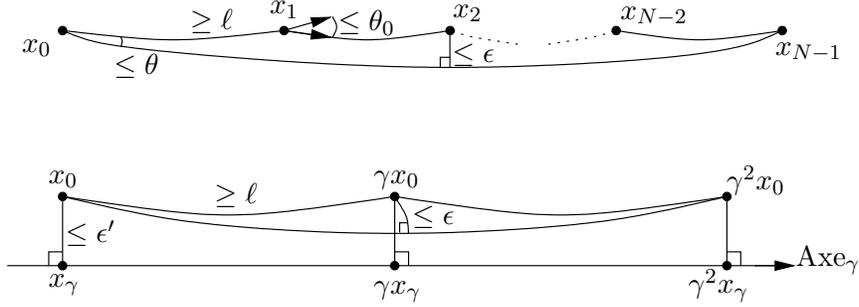
\begin{center}
\input{fig_closinglemma.pstex_t}
\end{center}

\blemm \label{lem:crithyperboelem} For all $\ell,\epsilon'>0$, there
exists $\epsilon=\epsilon(\ell,\epsilon')\in\mathopen{]}0,
1\mathclose{]}$ with $\lim_{\epsilon'\ra 0}\epsilon=0$ such that for
every isometry $\ga$ of any proper geodesic CAT$(-1)$-space $X$, for
every $x_0$ in $X$, if $d(x_0,\ga x_0)\geq \ell$ and $d(\ga
x_0,\mathopen{[}x_0,\ga^2x_0\mathclose{]})\leq \epsilon$, then $\ga$
is loxodromic and $d(x_0,\Axe_\ga)\leq\epsilon'$ where $\Axe_\ga$ is
the translation axis of $\ga$ in $X$.  
\elemm

\subsection{Pushing measures by branched covers}
\label{subsec:pushmeas}

In this text, we will often construct measures on $T^1M$ starting from
$\Ga$-invariant measures on $T^1\wt M$. Given a Galois covering
$f:X\ra Y$ of topological spaces, and a (positive Borel) measure $m$
on $X$ invariant under the covering group, it is well known that there
exists a unique measure on $Y$ (which is not the push-forward measure
$f_*m$) such that the map $f$ locally preserves the measure.  In the
case of ramified covers, the analogous construction is not so well
known (see for instance \cite[\S 2.4]{ParPau13b}). The fact that we
are not assuming $\Ga$ to be torsion free requires it.

\medskip
Let $\wt X$ be a locally compact metrisable space, endowed with a
proper (but not necessarily free) action of a discrete group $G$. Let
$p:\wt X\ra X= G\bs\wt X$ be the canonical projection. Let $\wt \mu$
be a locally finite $G$-invariant measure on $\wt X$.

Note that the map $N$ from $\wt X$ to $\NN-\{0\}$ sending a point
$x\in X$ to the order of its stabiliser in $G$ is upper
semi-continuous. In particular, for every $n\geq 1$, the $G$-invariant
subset $\wt X_n=N^{-1}(\{n\})$ is locally closed, hence locally
compact metrisable. With $X_n=p(\wt X_n)$, the restriction $p_{\mid
  \wt X_n}: \wt X_n\ra X_n$ is a local homeomorphism. Since $\wt\mu$
is $G$-invariant, there exists a unique measure $\mu_n$ on $X_n$ such
that the map $p_{\mid \wt X_n}$ locally preserves the measure. Now,
considering a measure on $X_n$ as a measure on $X$ with support in
$X_n$, let us define
$$
\mu=\sum_{n\geq 1}\frac{1}{n}\;\mu_n\;,
$$
which is a locally finite measure on $X$, called the measure {\it
  induced by}\index{measure!induced by a branched covering} $\wt \mu$
on $X$. Note that the factor $1/n$, besides being natural, is
necessary to get the continuity of see $\wt \mu\mapsto \mu$, see
below.

Note that if $\wt \mu$ gives measure $0$ to the set
$N^{-1}(\mathopen{[}2,+\infty\mathclose{[})$ of fixed points of
nontrivial elements of $G$, then $\mu=\mu_1$, and the above
construction is not really needed.

If $\wt X'$ is another locally compact metrisable space, endowed with
a proper action of a discrete group $G'$, if $\wt \mu'$ is a locally
finite $G'$-invariant measure on $\wt X'$, with induced measure $\mu'$
on $G'\backslash X'$, then the measure on the product of the quotient
spaces $G\bs X\;\times\; G'\bs X'$, induced by the product measure
$\wt\mu\otimes\wt \mu'$, is the product $\mu\otimes\mu'$ of the induced
measures.

The following result is easy to check.

\blemm \label{lem:pushmeasure} The map $\wt \mu\mapsto \mu$ from the
space of locally finite measures on $\wt X$ to the space of locally
finite measures on $X$, both endowed with their weak-star topologies,
is continuous.  
\elemm

We conclude this subsection with an example which indicates why the
factor $1/n$ before $\mu_n$ in the formula defining $\mu$ is
necessary. Let $G$ be the finite group of rotations of order $n$ in
the real plane $\wt X=\CC$, so that the quotient map $p:\CC\ra
X=G\bs\CC$ is a branched cover of order $n$. Denote by $\D_x$ the
unit Dirac mass at a point $x$. For all $\rho\geq 0$, let
$\wt\mu_\rho=\sum_{k=0}^{n-1} \D_{\rho\,e^{ik\theta}}$, which is a
locally finite $G$-invariant measure on $\CC$, which depends
continuously on $\rho$. If $\rho>0$, then, by the regular cover
situation, the measure induced by $\wt\mu_\rho$ is the unit Dirac mass
at the image of $\rho$ in $X$, which converges, as $\rho$ tends to
$0$, to the unit Dirac mass at the image $[0]$ of $0$ in $X$. Since
$\wt\mu_0$ is $n$ times the unit Dirac mass at $0$, in order to have
the required continuity property, we do have to define the measure
induced by $\wt\mu_0$ as $1/n$ times the push-forward of $\wt\mu_0$ by
$0\mapsto [0]$.

\section{A Patterson-Sullivan theory for Gibbs states} 
\label{sec:GPS}

Let $(\wt M,\Ga,\wt F)$ be as in the beginning of Chapter
\ref{sec:negacurvnot}: $\wt M$ is a complete simply connected
Riemannian manifold, with dimension at least $2$ and pinched sectional
curvature at most $-1$; $\Ga$ is a non-elementary discrete group of
isometries of $\wt M$; and $\wt F :T^1\wt M\ra \RR$ is a
H\"older-continuous $\Ga$-invariant map. We use the notation
$M=\Ga\bs\wt M$, $T^1M=\Ga\bs T^1\wt M$ and $F:T^1M\ra \RR$ (the map
induced by $\wt F$) introduced in Chapter \ref{sec:negacurvnot}.

\medskip We recall in this chapter the construction of the Gibbs
measure on the unit tangent bundle of a negatively curved manifold
associated with a given potential, due to O.~Mohsen \cite{Mohsen07}.
More precisely, the subsections \ref{subsec:gap},
\ref{subsec:crossratio}, \ref{subsec:GibbsSulivanmeasure},
\ref{subsec:gibbsproperty}, \ref{subsec:condmesgibbs} are new (except
for the strong influence of \cite{Hamenstadt97}), the others are
extracted (up to an adaptation to the non-cocompact case) from
\cite{Mohsen07}. We refer to \cite{Hamenstadt97, Ledrappier95,
  Coudene03, Schapira04a} for different approaches. The Gibbs cocycles
used in these four references are the same up to signs as the one we
will define in Subsection \ref{subsec:Gibbscocycle}. Note that the last
three references use $\sum_{\ga\in\Ga} e^{s\int_{x}^{\gamma y} \wt F}$
for their Poincar\'e series, whereas we use $\sum_{\ga\in\Ga}
e^{\int_{x}^{\gamma y} (\wt F-s)}$.  This additive instead of
multiplicative approach greatly simplifies the techniques.

\subsection{Potential functions and their periods} 
\label{subsec:potentials}

For all $x,y$ in $\wt M$, let us define
$$
\int_x^y \wt F=\int_0^{d(x,\,y)}\wt F(\phi_t(v))\;dt\;,
$$
where $v$ is a unit tangent vector such that $\pi(v)=x$ and
$\pi(\phi_{d(x,y)}v)=y$ (unique if $y\neq x$).

\medskip
\noindent{\bf Remark. } A general extension of the techniques of this
work when $\wt M$ is any complete locally compact CAT$(-1)$ metric
space as in \cite{Roblin03} (even assuming that its geodesic segments
are extendible to geodesic lines) seems difficult at this point. The
correct analog of the unit tangent bundle for such an extension is
probably not the usual one of the space of all geodesic lines. Indeed,
there could be many geodesic lines containing a given segment
$\mathopen{[}x,y\mathclose{]}$. Hence defining $\int_x^y \wt F$ would
require a choice of one of these geodesic lines, or restrictions on
$\wt F$ or some averaging process, hence appropriate additional
information in order not to lose the invariance under the group
$\Ga$. A better analog would be the space of germs of geodesic
segments, but the geodesic flow then does not extend. See the case of
trees in \cite{BroParPau13}.

\medskip
Note that
\begin{equation}\label{eq:timereversal}
\forall \;\ga\in\Ga,\;\;\;\int_{\ga x}^{\ga y} \wt F=\int_x^y\wt F
\;\;\;\;\;\;{\rm and}\;\;\;\;\;\;
\int_y^x \wt F=\int_x^y\wt F\circ\iota\;.
\end{equation}
For every loxodromic element $\ga$ in $\Ga$, the {\it 
period of $\ga$}\index{period} for the potential $F$ is 
$$
\operatorname{Per}_F(\ga)=\int_x^{\ga x} \wt F
$$ 
for any point $x$ on the translation axis of $\ga$ in $\wt M$.  Note
that, for all $n\in\NN-\{0\}$ and $\alpha\in\Ga$, we have
\begin{equation}\label{eq:periodinvconj}
\operatorname{Per}_F(\alpha\ga\alpha^{-1})=\operatorname{Per}_F(\ga),
\;\;\;\operatorname{Per}_F(\ga^n)=n\;\operatorname{Per}_F(\ga)
\;\;\;{\rm and}\;\;\; 
\operatorname{Per}_F(\ga^{-1})= \operatorname{Per}_{F\circ\iota}(\ga)\;.
\end{equation}

Let $\wt F^* :T^1\wt M\ra \RR$ be another H\"older-continuous
$\Ga$-invariant map. Say that $\wt F^*$ is {\it
  cohomologous}\index{cohomologous}\index{potential!cohomologous} to
$\wt F$ (see for instance \cite{Livsic72}) if there exists a
H\"older-continuous $\Ga$-invariant map $\wt G :T^1\wt M\ra \RR$,
differentiable along every flow line, such that
$$
\wt F^*(v)-\wt F(v)=\frac{d}{dt}_{\mid t=0}\wt G(\phi_tv)\;.
$$
We will say that $\wt F^*$ is {\it cohomologous to $\wt F$ via $\wt
  G$}\index{cohomologous!via a function $\wt G$}%
\index{potential!cohomologous!via a function $\wt G$} when we want to
emphasise $\wt G$.  Two cohomologous potentials have the same periods:
if $\wt F^*$ is cohomologous to $\wt F$, then, denoting by $F^*$ the map
induced on $T^1M$ by $\wt F^*$, for every loxodromic element $\ga\in
\Ga$, we have
$$
\operatorname{Per}_F(\ga)=\operatorname{Per}_{F^*}(\ga)\;.
$$

\brema\label{rem:livsi} Let $\wt F,\wt F^* :T^1\wt M\ra \RR$ be
H\"older-continuous $\Ga$-invariant maps such that
$\operatorname{Per}_F(\ga) = \operatorname{Per}_{F^*}(\ga)$ for every
loxodromic element $\ga\in\Ga$. Then $\wt F^*$ is cohomologous to $\wt
F$ in restriction to $\wt \Omega\Ga$.  
\erema

\dem Since the diagonal action of $\Ga$ on $\Lambda\Ga\times
\Lambda\Ga$ is topologically transitive, the action of the geodesic
flow on the (topological) non-wandering set $\Omega\Ga$ is
topologically transitive. The validity of Anosov's closing lemma (see
Lemma \ref{lem:quasigeodpresquegeod}) does not require any compactness
assumption on $M$. The proof of Liv\v{s}ic's theorem in
\cite[Theo.~19.2.4]{KatHas95} then extends. \cqfd

\medskip
Note that if $\wt F^*$ and $\wt F$ are cohomologous via the map $\wt
G$, then $\wt F^*\circ \iota$ and $\wt F\circ \iota$ are cohomologous
via the map $-\,\wt G\circ \iota$ (beware the sign). We will say that
$\wt F$ (or $F$) is {\it reversible}\index{reversible}%
\index{potential!reversible} if $\wt F$ is cohomologous to $\wt F\circ
\iota$.

\medskip We end this subsection by the following technical lemma on
potential functions (see also \cite[Lem.~3]{Coudene03} with the
multiplicative approach).

\blemm\label{lem:technicholder} For every $r_0>0$, there exist two
constants $c_1>0$ and $c_2\in\mathopen{]}0,1\mathclose{]}$ (depending,
besides the $r_0$ dependence of $c_1$, only on the H\"older constants
of $\wt F$ and the bounds on the sectional curvature of $\wt M$) such
that for all $x,y,z$ in $\wt M$, we have
$$
\Big|\;\int_x^{z}\wt F\;-\int_{y}^{z}\wt F\;\Big|\leq
c_1\,e^{d(x,\,y)}+ d(x,y)\max_{\pi^{-1}(B(x,\,d(x,\,y)))}|\wt F|\;,
$$
and if furthermore $d(x,y)\leq r_0$, then 
$$
\Big|\;\int_x^{z}\wt F\;-\int_{y}^{z}\wt F\;\Big|\leq
c_1\,d(x,y)^{c_2}+ d(x,y)\max_{\pi^{-1}(B(x,\,d(x,\,y)))}|\wt F|\;.
$$
\elemm

Note that in the first c entered equation (as well as the one in Lemma
\ref{lem:holderconseq}), the constant $c_1$ actually does not depend
on any $r_0$, hence could be given some other name. But the
formulation of the statement helps to keep at a minimum the number of
constants, and facilitates future simultaneous references to both c
entered equation.

\medskip
\dem By symmetry, we may assume that $d(x,z)\geq d(y,z)$. The result
is true if $y=z$, hence we assume that $y\neq z$.

\medskip\noindent\begin{minipage}{8.9cm} ~~~ Let $x'$ be the point on
  $\mathopen{[}x,z\mathclose{]}$ at distance $d(y,z)$ from $z$. Let
  $v$ (respectively $w$) be the unit tangent vector at $x'$
  (respectively $y$) pointing towards $z$.
\end{minipage}
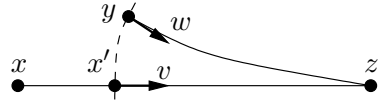
\begin{minipage}{6cm}
\begin{center}
\input{fig_xyzholder.pstex_t}
\end{center}
\end{minipage}

\medskip The closest point $p$ of $y$ on $\mathopen{[}x,
z\mathclose{]}$ lies in $\mathopen{[}x',z\mathclose{]}$ by
convexity. Hence $d(x,x')\leq d(x,p)\leq d(x,y)$, since closest point
maps do not increase distances. Since the distance function from a
given point to a point varying on a geodesic line is convex, we have
$d(x',y)\leq d(x,y)$.

Let $c>0$ and $\alpha\in\mathopen{]}0,1\mathclose{]}$ be the H\"older
constants of $\wt F$ for the distance $d'=d'_{T^1\wt M}$ on $T^1\wt M$
defined in Equation \eqref{eq:defidprime}, that is $|\wt F(u)-\wt
F(u')|\leq c\,d'(u,u')^\alpha$ for all $u,u'$ in $T^1\wt M$ with
$d'(u,u')\leq 1$.  By Corollary \ref{coro:controlexpopinchunit} (i),
there exists a universal constant $c'>0$ (for instance
$c'=\frac{e^2}{2}$) such that for every $t\in\mathopen{[}0,
d(y,z)\mathclose{]}$, we have $d'(\phi_tv,\phi_tw)\leq
c'e^{d(x',\,y)-t}$.

Let $t\in\mathopen{[}0,d(y,z)\mathclose{]}$. If
$d'(\phi_tv,\phi_tw)\leq 1$, then
$$
|\wt F(\phi_tv)- \wt F(\phi_tw)| \leq c \,d'(\phi_tv,\phi_tw)^\alpha
\leq c \,\big(c'\,e^{d(x',\,y)-t}\big)^\alpha\;.
$$
If $d'(\phi_tv,\phi_tw)\geq 1$, then by the remark at the end of 
Subsection \ref{subsec:holdercont}, we have
$$
|\wt F(\phi_tv)- \wt F(\phi_tw)| \leq 3\,c \,d'(\phi_tv,\phi_tw)
\leq 3\,c \,c'\,e^{d(x',\,y)-t}\;.
$$
Hence,
\begin{align*}
  \Big|\;\int_x^{z}\wt F-\int_{y}^{z}\wt F\;\Big| & =
  \Big|\;\int_{-d(x,\,x')}^{d(y,\,z)}\wt F(\phi_tv)\,dt-
  \int_{0}^{d(y,\,z)}\wt F(\phi_tw)\,dt\;\Big|\\ & \leq
  \int_{-d(x,\,x')}^{0}|\wt F(\phi_tv)|\,dt+ \int_{0}^{d(y,\,z)}|\wt
  F(\phi_tv)- \wt F(\phi_tw)|\,dt\\ & \leq
  d(x,x')\max_{\pi^{-1}(B(x,\,d(x,\,x')))}|\wt F|+ \int_{0}^{+\infty}
  c\,\big(c'\,e^{d(x',\,y)-t}\big)^\alpha+3\,c\,c'\,e^{d(x',\,y)-t}\,dt 
\\ & \leq
  d(x,y)\max_{\pi^{-1}(B(x,\,d(x,\,y)))}|\wt
  F|+\big(3\,c\,c'+\frac{c(c')^\alpha}{\alpha}\big)\,e^{d(x,\,y)}\;.
\end{align*}
The first result follows. The second may be proved similarly, using
Corollary \ref{coro:controlexpopinchunit} (ii).  
\cqfd

\medskip
\noindent {\bf Remarks. } (1) When $x,y,z\in\C\Lambda\Ga$, we may
replace, in each assertion of the above lemma,
$\max_{\pi^{-1}(B(x,\,d(x,\,y)))}|\wt F|$ by
$\max_{\pi^{-1}(B(x,\,d(x,\,y))\,\cap\,\C\Lambda\Ga)}|\wt F|$.

\medskip (2) Using the equality $\int_x^{x'}\wt F-\int_{y}^{y'}\wt F=
\big(\int_x^{x'}\wt F-\int_{y}^{x'}\wt F\big)+ \big(\int_{x'}^{y}\wt
F\circ\iota-\int_{y'}^{y}\wt F\circ\iota\big)$ and the fact that
$\iota$ is Lipschitz with constants depending only on the bounds on
the sectional curvature of $\wt M$, we have a similar control on $
\big|\;\int_x^{x'}\wt F-\int_{y}^{y'}\wt F\;\big|$ for all $x,x',y,y'$
in $\wt M$, in terms of $\max\{d(x,y),d(x',y')\}$.

\subsection{The Poincar\'e series and the critical exponent of
  $(\Ga,F)$}
\label{subsec:GibbsPoincareseries}

Let us fix $x,y\in\wt M$. The {\it Poincar\'e series}%
\index{Poincar\'e series} of $(\Ga,F)$ is the map $\gls{Poincarserie}:
\RR\ra\mathopen{[}0,+\infty\mathclose{]}$ defined by
$$
Q_{\Ga,\,F,\,x,\,y}(s)=\sum_{\ga\in\Ga} \;\; 
e^{\int_x^{\ga y} (\wt F-s)}\;.
$$
The {\it critical exponent}\index{critical exponent} of $(\Ga,F)$ is
the element $\gls{criticalexponent}$ in $\mathopen{[}-\infty,
+\infty\mathclose{]}$ defined by
$$
\delta_{\Ga,\,F}=\limsup_{n\ra +\infty}\;\frac{1}{n}\ln
\;\sum_{\ga\in\Ga,\;n-1< d(x,\,\ga y)\leq n} \;e^{\int_x^{\ga y} \wt F}\;.
$$
When $F=0$, the Poincar\'e series $\gls{Poincarserieusuel}$ of
$(\Ga,0)$ is the usual {\it Poincar\'e series}%
\index{Poincar\'e series} of $\Ga$, and the critical exponent
$\delta_{\Ga,\,0}$ is the {\it critical exponent}%
\index{critical exponent} $\gls{criticalexponentusuel}$ of $\Gamma$
(which belongs to $\mathopen{]}0,+\infty\mathclose{[}$ since $\Ga$ is
non-elementary and $\wt M$ has pinched negative curvature, see for
instance \cite{Roblin03}).

We will prove in the following Lemma \ref{lem:elemproppressure} (v)
that $\delta_{\Ga,\,F}>-\infty$. If $\delta_{\Ga,\,F}<+\infty$, we say
that $(\Ga,F)$ is {\it of divergence type}\index{divergence type} if
the series $Q_{\Ga,\,F}(\delta_{\Ga,\,F})$ diverges, and {\it of
  convergence type}\index{convergence type} otherwise.

If $\delta_{\Ga,\,F}<+\infty$, the {\it normalised
  potential}\index{normalised potential}\index{potential!normalised}
is the H\"older-continuous map $\wt F-\delta_{\Ga,\,F}$ on $T^1\wt M$
(or its induced map $F-\delta_{\Ga,\,F}$ on $T^1 M$). Some references
use the normalised potential with the opposite sign.

\medskip Since $\wt F$ is H\"older-continuous (see Remark (2) at the
end of the previous Subsection \ref{subsec:potentials}), the critical
exponent $\delta_{\Ga,\,F}$ of $(\Ga,F)$ does not depend on $x,y$, and
satisfies the following elementary properties.

\blemm\label{lem:elemproppressure} 
(i) The Poincar\'e series of $(\Ga,F)$ converges if
$s>\delta_{\Ga,\,F}$ and diverges if $s<\delta_{\Ga,\, F}$. The
Poincar\'e series of $(\Ga,F)$ diverges at $\delta_{\Ga,\,F}$ if
and only if the Poincar\'e series of $(\Ga,F+\kappa)$ diverges
at $\delta_{\Ga,\,F+\kappa}$, for every $\kappa\in\RR$, and in
particular
\begin{equation}\label{eq:invarpressuretranslat}
\forall\; \kappa\in\RR,\;\;\;\delta_{\Ga,\,F+\kappa}=
\delta_{\Ga,\,F}+\kappa\;.
\end{equation} 
(ii) We have
\begin{equation}\label{eq:pressuretimereversalinva}
\forall\;s\in\RR,\; Q_{\Ga,\,F\circ\iota,\,x,\,y}(s)=Q_{\Ga,\,F,\,y,\,x}(s)
\;\;\;{\rm and}\;\;\;
\delta_{\Ga,\,F\circ \iota}=\delta_{\Ga,\,F}\;.
\end{equation}
(iii) If $\Ga'$ is a non-elementary subgroup of $\Ga$, denoting by
$F':\Ga'\backslash T^1\wt M\ra\RR$ the map induced by $\wt F$, we have
$$
\delta_{\Ga',\,F'}\leq \delta_{\Ga,\,F}\;.
$$
(iv) We have the upper and lower bounds
$$
\delta_\Gamma + \inf_{\pi^{-1}(\C\Lambda\Ga)}\wt F\;\leq\; 
\delta_{\Ga,\,F}\;\leq \;
\delta_\Gamma + \sup_{\pi^{-1}(\C\Lambda\Ga)}\wt F\;.
$$
(v) We have  $\delta_{\Ga,\,F}>-\infty$.

\noindent (vi) The map $F\mapsto \delta_{\Ga,\,F}$ is convex,
sub-additive and $1$-Lipschitz for the uniform norm on the vector space
of real continuous maps on $\pi^{-1}(\C\Lambda\Ga)$, that is, if
$\wt{F^*}:T^1\wt M\ra \RR$ is another H\"older-continuous $\Ga$-invariant
map, inducing $F^*:\Ga\backslash T^1\wt M\ra\RR$, if
$\delta_{\Ga,\,F},\delta_{\Ga,\,F^*}<+\infty$, then
$$
|\;\delta_{\Ga,\,F^*}-\delta_{\Ga,\,F}\;|\leq 
\sup_{v\in \pi^{-1}(\C\Lambda\Ga)} |\;\wt {F^*}(v)-\wt F(v)\;|\;,
$$
$$
\delta_{\Ga,\,F+F^*}\leq \delta_{\Ga,\,F}+\delta_{\Ga,\,F^*}\;,
$$
and, for every $t\in\mathopen{[}0,1\mathclose{]}$, 
$$
\delta_{\Ga,\,tF+(1-t)F^*}\leq t\,\delta_{\Ga,\,F}+(1-t)\,\delta_{\Ga,\,F^*}\;.
$$

\noindent (vii) For every $c>0$, we have 
\begin{equation}\label{eq:pressuresumleqnbis}
\delta_{\Ga,\,F}=\limsup_{n\ra +\infty}\;\frac{1}{n}\ln
\;\sum_{\ga\in\Ga,\;n-c<d(x,\,\ga y)\leq n} \;e^{\int_x^{\ga y} \wt F}\;,
\end{equation} 
and if $\delta_{\Ga,\,F}\geq 0$, then 
\begin{equation}\label{eq:pressuresumleqn}
\delta_{\Ga,\,F}=\limsup_{n\ra +\infty}\;\frac{1}{n}\ln
\;\sum_{\ga\in\Ga,\;d(x,\,\ga y)\leq n } \;e^{\int_x^{\ga y} \wt F}\;.
\end{equation} 

\noindent (viii) If $\Ga''$ is a discrete cocompact group of isometries
of $\wt M$ such that $F$ is $\Ga''$-invariant, denoting by
$F'':\Ga''\backslash T^1\wt M\ra\RR$ the map induced by $\wt F$,
we have
$$
\delta_{\Ga,\,F}\leq \delta_{\Ga'',\,F''}\;.
$$

\noindent (ix) If $\Ga$ is the infinite cyclic group\footnote{Note
  that $\Ga$ is elementary, but we can still define the Poincaré
  series of $(\Ga,F)$ and its critical exponent.} generated by a
loxodromic isometry $\ga$ of $\wt M$, then $(\Ga,F)$ is of divergence
type and
$$
\delta_{\Ga,\,F}=\max\big\{\;\frac{\operatorname{Per}_F(\ga)}{\ell(\ga)},\;
\frac{\operatorname{Per}_{F\circ\iota}(\ga)}{\ell(\ga)}\;\big\}\,.
$$
Furthermore, $\sum_{\ga\in\Ga,\;d(x,\,\ga y)\leq t } \;e^{\int_x^{\ga y} \wt F}
\asymp\left\{\begin{array}{ll}
e^{\delta_{\Ga,\,F}\;t}&{\rm if}\;\delta_{\Ga,\,F}>0\\
t&{\rm if}\;\delta_{\Ga,\,F}=0 \\ 1&{\rm if}\;\delta_{\Ga,\,F}<0 
\end{array}\right.$.
\elemm

We will prove in Subsection \ref{subsec:loggrowth} that the upper
limits in Equation \eqref{eq:pressuresumleqn} if $c$ is large enough, and
in Equation \eqref{eq:pressuresumleqn} if $\delta_{\Ga,\,F}>0$ are in
fact limits.

An interesting problem is to study whether or not the critical
exponent $\delta_{\Ga,\,F}$ of $(\Ga,F)$ is equal to the upper bound of
the critical exponents $\delta_{\Ga_0,\,F_0}$ where $\Ga_0$ ranges over
the convex-cocompact subgroups of $\Ga$ and $F_0:\Ga_0\backslash
T^1\wt M\ra\RR$ is the map induced by $\wt F$. Replacing subgroups by
subsemigroups, we answer this positively in Subsection
\ref{subsec:semigroup}.

\medskip
\dem The verifications of (i), (ii) and (iii) are elementary. For
instance, Equation \eqref{eq:pressuretimereversalinva} follows from
both parts of Equation \eqref{eq:timereversal} and the change of
variable $\ga\mapsto \ga^{-1}$ in the Poincar\'e series.

\medskip
To prove Assertion (iv), note that if $x$ is a point in the convex
hull of the limit set $\C\Lambda\Gamma$, then, for every $\ga\in\Ga$,
the geodesic segment between $x$ and $\ga x$ is contained in
$\C\Lambda\Gamma$. Hence
$$
d(x,\ga x)\;\Big(\inf_{\pi^{-1}(\C\Lambda\Ga)}\wt F-s\Big)\leq
{\int_x^{\ga x} (\wt F-s)}\leq d(x,\ga x)\;
\Big(\sup_{\pi^{-1}(\C\Lambda\Ga)}\wt F -s\Big)\;.
$$
This proves Assertion (iv) for instance by taking the exponential,
summing over $\ga\in\Ga$ and using the first assertion of (i).

\medskip
To prove Assertion (v), let $\Ga'$ be a non-elementary convex-cocompact
subgroup of $\Ga$ (for instance a Schottky subgroup of $\Ga$). Denote
by $F':\Ga'\backslash T^1\wt M\ra\RR$ the map induced by $\wt F$. Since
$|\wt F|$ is $\Ga$-invariant and bounded on compact subsets of $\wt
M$, by Assertion (iv), we have $\delta_{\Ga',\,F'}>-\infty$. Assertion
(v) then follows from Assertion (iii).

To prove Assertion (vi), let $c=\sup_{\pi^{-1}(\C\Lambda\Ga)}
|\wt F-\wt{F^*}|$ and $x\in\C\Lambda\Ga$. We have
$$
\int_x^{\ga x}(\wt{F}-(s+c))\leq\int_x^{\ga x}(\wt {F^*}-s)\leq 
\int_x^{\ga x}(\wt{F}-(s-c))\;,
$$ 
and the first claim follows.  To prove the two remaining
claims of (vi), for all $s>\delta_{\Ga,\,F}$ and
$s^*>\delta_{\Ga,\,F^*}$, we have
$$
Q_{\Ga,\,F,\,x,\,y}(s)Q_{\Ga,\,F^*,\,x,\,y}(s^*)\geq Q_{\Ga,\,F+F^*,\,x,\,y}(s+s^*)
$$
and, by the convexity of the exponential, 
$$
Q_{\Ga,\,tF+(1-t)F^*,\,x,\,y}(ts+(1-t)s^*)\leq 
tQ_{\Ga,\,F,\,x,\,y}(s)+(1-t)Q_{\Ga,\,F^*,\,x,\,y}(s^*)\;.
$$
Hence $s+s^*\geq \delta_{\Ga,\,F+F^*}$ and $\delta_{\Ga,\,tF+(1-t)F^*}
\leq ts+(1-t)s^*$ by the first assertion of (i). The result follows by
letting $s$ and $s^*$ converge respectively to $\delta_{\Ga,\,F}$ and
$\delta_{\Ga,\,F^*}$.

\medskip
Assertion (vii) is also elementary. For its first claim, one
cuts the interval $\mathopen{]}t-c,t\mathclose{]}$ at the points $t-k$
for $k=0,\cdots, \lfloor c\rfloor$ if $c\geq 1$, or the interval
$\mathopen{]}t-1,t\mathclose{]}$ at the points $t-ck$ for $k=0,\cdots,
\lfloor \frac{1}{c}\rfloor$ if $c\leq 1$, and one uses the inequalities
$\max_{1\leq i\leq N}a_i\leq \sum_{1\leq i\leq N}a_i\leq N \max_{1\leq
  i\leq N}a_i$.

For the second claim, the right hand side of Equation
\eqref{eq:pressuresumleqn} is clearly not smaller than
$\delta_{\Ga,\,F}$.  To prove the other inequality, for every
$\epsilon >0$, let $n_0\in\NN$ be such that $\sum_{n-1 < d(x,\,\ga y)
  \leq n} e^{\int_x^y F} \leq e^{(\delta_{\Ga,\,F}+\epsilon)n}$ for
all $n\geq n_0$. Then for all $n\geq n_0$, we have 
$$
\sum_{d(x,\,\ga y)\leq n} e^{\int_x^y F} \leq \sum_{k=0}^{n} 
e^{(\delta_{\Ga,\,F}+\epsilon)k}+\operatorname{O}(1) 
\leq c'\;e^{(\delta_{\Ga,\,F}+\epsilon)n}
$$ 
for some $c'>0$, since $\delta_{\Ga,\,F} +\epsilon>0$. The result
follows.  

\medskip 
Let us prove Assertion (viii). Since $\Ga''\bs \wt M$ is
compact, there exists $c\geq 0$ such that for every $\ga\in\Ga$, there
exists $\ga''\in\Ga''$ be such that $d(\ga x,\ga'' x)\leq c$.  Hence
for every $s>\delta_{\Ga'',\,F''}$, we have, with $r_0=c$ and $c_1$ as
in Lemma \ref{lem:technicholder} and by $\Ga$-equivariance of $\wt F$,
\begin{align*}
Q_{\Ga,\,F,\,x,\,x}(s)&\leq \sum_{\ga\in\Ga} \;
e^{\int_x^{\ga'' x} (\wt F-s) +c_1\,e^{d(\ga x,\,\ga'' x)}+d(\ga x,\,\ga'' x)
\max_{\pi^{-1}(B(\ga x,\,d(\ga x,\,\ga ''y)))}|\wt F|}\\ & \leq
e^{c_1\,e^c+c\,\max_{\pi^{-1}(B(x,\,c))}|\wt F|}\;
Q_{\Ga'',\,F'',\,x,\,x}(s)<+\infty\;.
\end{align*}
Hence $\delta_{\Ga,\,F}\leq s$, and the result follows by letting $s$
tend to $\delta_{\Ga'',\,F''}$. 

\medskip 
Finally, in order to prove Assertion (ix), if $x$ belongs to
the translation axis of $\ga$, we have
\begin{align*}
\sum_{\ga\in\Ga} \;
e^{\int_x^{\ga x} (\wt F-s)}& =\sum_{n\in\NN} \;
e^{\int_x^{\ga^n x} (\wt F-s)}+\sum_{n\in\NN-\{0\}} \;
e^{\int_x^{\ga^{-n} x} (\wt F-s)}\\ &=\sum_{n\in\NN} \;
e^{n(\operatorname{Per}_F(\ga)-s\,\ell(\ga))}+\sum_{n\in\NN-\{0\}} \;
e^{n(\operatorname{Per}_{F\circ\iota}(\ga)-s\,\ell(\ga))}\,.
\end{align*}
Hence $Q_{\Ga,\,F,\,x,\,x}(s)$ converges if and only if
$\operatorname{Per}_F(\ga)-s\,\ell(\ga)<0$ and
$\operatorname{Per}_{F\circ\iota}(\ga)-s\,\ell(\ga)<0$, which proves the
first two claims of Assertion (xi). We have, again if $x$ belongs to
the translation axis of $\ga$,
$$
\sum_{\ga\in\Ga,\;d(x,\,\ga x)\leq t} \;e^{\int_x^{\ga x} \wt F}=
\sum_{0\leq k\leq t/\ell(\ga)} \;e^{k\,\operatorname{Per}_F(\ga)}+
\sum_{0< k\leq t/\ell(\ga)} \;e^{k\,\operatorname{Per}_{F\circ\iota}(\ga)}\,.
$$
The last claim follows, by a standard geometric series argument.
\cqfd

\medskip Here are a few immediate consequences. The convergence or
divergence of the Poincar\'e series $Q_{\Ga,\,F }(s)$ is independent of
the points $x$ and $y$. The Poincar\'e series $Q_{\Ga,\,F\circ\iota}(s)$
converges if and only if $Q_{\Ga,\,F }(s)$ does. 

The critical exponent $\delta_{\Ga,\,F}$ of $(\Ga,F)$ is positive if
$\inf_{T^1\wt M}\wt F> -\delta_\Ga$. For instance, $\delta_{\Ga,\,F}
>0 $ if $F$ is a small perturbation of zero (more precisely if
$\|F\|_\infty=\sup_{T^1\wt M}|\wt F|<\delta_\Ga$), or if $F\geq 0$,
since $\Ga$ is non-elementary and therefore $\delta_\Ga>0$.  Recall
that $\delta_\Ga< +\infty$ since $\wt M$ has pinched curvature. Hence
if $F$ is bounded from above, then $\delta_{\Ga,\,F}<+\infty$.

When $\delta_{\Ga,\,F}<+\infty$, the additive convention for the
Poincar\'e series and the fact that its interesting behaviour occurs at
$s=\delta_{\Ga,\,F}$ suggest it is meaningful to consider the
normalised potential $F-\delta_{\Ga,\,F}$. By Equation
\eqref{eq:invarpressuretranslat}, the critical exponent of the
normalised potential is $0$. Moreover, adding a large enough constant
to $F$ does not change the normalised potential and allows the
critical exponent of $(\Ga,F)$ to be positive (this reduction
will be used repeatedly in many results of the chapters
\ref{sec:growth} and \ref{sec:ergtheounistabfolia}).

The relationship between $\delta_{\Ga,\,F}$ and $\delta_{\Ga,\,sF}$
for $s>0$ seems unclear in general.

\medskip \rem Let $\wt F^*:T^1\wt M\ra \RR$ be another
H\"older-continuous $\Ga$-invariant map, which is cohomologous to $\wt
F$ via $\wt G$. Let $F^*$ be the map induced on $T^1M=\Ga\bs T^1\wt M$
by $\wt F^*$. Note that for all $x,y\in \wt M$ and $\ga\in\Ga$, by
$\Ga$-invariance of $\wt G$,
$$
\Big|\int_x^{\ga y} \wt F -\int_x^{\ga y} \wt F^*\;\Big|=
\big|\wt G(\phi_{d(x,\ga y)}v)-\wt G(v)\big|\leq
2\sup_{w\,\in\, T^1_x\wt M\,\cup \,T^1_y\wt M}|\wt G(w)|\;,
$$
where $v$ is a unit tangent vector such that $\pi(v)=x$ and
$\pi(\phi_{d(x,\ga y)}v)=\ga y$.  Hence $F$ and $F^*$ have the same
critical exponent:
$$
\delta_{\Ga,\,F^*}= \delta_{\Ga,\,F}\;.
$$
Furthermore, $(\Ga,F^*)$ is of divergence type if and only if
$(\Ga,F)$ is of divergence type.

\subsection{The Gibbs cocycle of $(\Ga,F)$}
\label{subsec:Gibbscocycle}

By the H\"older-continuity of $F$ and the properties of asymptotic
geodesic rays in $\wt M$, there exists a well-defined map
$C_{F}: \partial_\infty\wt M\times \wt M\times\wt M\ra \RR$, called
the {\it Gibbs cocycle}\index{Gibbs cocycle}\index{cocycle!Gibbs} for
the potential $F$, defined by
$$
(\xi,x,y)\mapsto \gls{Gibbscocycle}=\lim_{t\ra+\infty} 
\int_y^{\xi_t}\wt F-\int_x^{\xi_t}\wt F
$$ 
for $t\mapsto \xi_t$ any geodesic ray ending at $\xi$. Do note the
apparent order reversal of $x$ and $y$ in this formula, which allows us
to have the required sign in front of the Busemann cocycle in Equation
\eqref{eq:translatgibbscocycle}.  We will compare our cocycle with the
one introduced by Hamenst\"adt in \cite{Hamenstadt97} at the end of
Subsection \ref{subsec:crossratio}.  

If $x$ belongs to the geodesic ray from $y$ to $\xi$, then
\begin{equation}\label{eq:cocyclealongray}
  C_{F,\,\xi}(x,y)= \int_y^{x} \wt F\;.
\end{equation}
In particular, for every $w\in T^1\wt M$, for all $x$ and $y$ on the
image of the geodesic line defined by $w$, with $w_-,x,y,w_+$ in this
order, we have
\begin{equation}\label{eq:changemoinsplus}
C_{F\circ\iota,\,w_-}(x,y)=\int_y^x \wt F\circ\iota=\int_x^y \wt F
=C_{F,\,w_+}(y,x)=-C_{F,\,w_+}(x,y)\;.
\end{equation}
Note that when $F=-1$, the Gibbs cocycle equals the Busemann cocycle
(defined at the end of Subsection \ref{subsec:geombus})
$$
\beta_\xi(x,y)=\lim_{t\ra+\infty} d(x,\xi_t)-d(y,\xi_t)\,
$$
with $t\mapsto \xi_t$ as above. Hence, for every $s\in\RR$,
\begin{equation}\label{eq:translatgibbscocycle}
C_{F-s,\,\xi}(x,y)=C_{F,\,\,\xi}(x,y)+s\;\beta_\xi(x,y)\;.
\end{equation}

The Gibbs cocycle satisfies the following cocycle property: for all
$\xi$ in $\partial_\infty\wt M$ and $x,y,z$ in $\wt M$,
\begin{equation}\label{eq:cocycleprop}
C_{F,\,\xi}(x,z)=C_{F,\,\xi}(x,y)+C_{F,\,\xi}(y,z)\;\;{\rm and}
\;\;C_{F,\,\xi}(x,y)=-\,C_{F,\,\xi}(y,x)\;,
\end{equation}
and the following invariance property: for all $\ga$ in $\Ga$,
$\xi$ in $\partial_\infty\wt M$ and $x,y$ in $\wt M$,
\begin{equation}\label{eq:equivcocprop}
C_{F,\,\ga\xi}(\ga x,\ga y)=C_{F,\,\xi}(x,y)\;.
\end{equation}

\noindent{\bf Remarks. } (1) If $\wt F^* :T^1\wt M\ra \RR$ is a
H\"older-continuous $\Ga$-invariant map, which is cohomologous to $\wt
F$, then the cocycle $C_{F^*}-C_{F}$ is a coboundary. Indeed, and more
precisely, for all $x,y\in \wt M$, for every $\xi\in\partial_\infty
\wt M$, if $\wt F^*$ is cohomologous to $\wt F$ via the map $\wt G
:T^1\wt M\ra \RR$, if $v_{z\xi}$ is the tangent vector at the origin
of the geodesic ray from a point $z\in\wt M$ to $\xi$, then
\begin{equation}\label{eq:cohomologuecocycle}
C_{F^*,\,\xi}(x,y)-C_{F,\,\xi}(x,y)=\wt G(v_{x\xi})-\wt G(v_{y\xi})\;.
\end{equation}
In particular, $C_{F^*}-C_{F}$ is bounded on $\partial_\infty\wt
M\times \wt M\times\wt M$ if $\wt G$ is bounded.

\medskip (2) Note that the Gibbs cocycle $C_{F}$ entirely determines
the potential function $\wt F$, since for every $v\in T^1\wt M$, we have
$$
\wt F(v)=\lim_{t\ra 0^+}\frac{1}{t}\,C_{F,\,v_+}(\pi(v),\pi(\phi_tv))\;.
$$

\medskip (3) For every $x_0\in\wt M$, the map $(\xi,\ga)\mapsto
C_{F,\,\xi}(x_0,\ga x_0)$ from $\Lambda\Ga\times \Ga$ to $\RR$
determines, up to a cohomologous potential, the restriction of $\wt F$
to $\wt \Omega\Ga$.

Indeed, let $y_0$ be the closest point to $x_0$ on the translation
axis of a loxodromic element $\ga\in\Ga$. Let $\ga_+$ be the
attractive fixed point of $\ga$. By the equations
\eqref{eq:cocycleprop}, \eqref{eq:equivcocprop} and
\eqref{eq:cocyclealongray}, we have
\begin{align*}
C_{F,\,\ga^+}(x_0,\ga
x_0) =C_{F,\,\ga^+}(x_0,y_0)+C_{F,\,\ga^+}(y_0,\ga y_0)+
C_{F,\,\ga^+}(\ga y_0,\ga x_0)
=\int_{y_0}^{\ga y_0} \wt F\;.
\end{align*}
Hence 
$$
\operatorname{Per}_F(\ga)= C_{F,\,\ga^+}(x_0,\ga x_0)\;.
$$
The result now follows from Remark \ref{rem:livsi}.

\blemm \label{lem:holderconseq} For every $r_0>0$, there exist
$c_1,c_2,c_3,c_4>0$ with $c_2,c_4\leq 1$ (depending only on $r_0$, the
H\"older constants of $\wt F$ and the bounds on the sectional
curvature of $\wt M$) such that the following assertions hold.

(1) For all $x,y\in \wt M$ and $\xi\in\partial_\infty \wt M$,
$$
|\;C_{F,\,\xi}(x,y)\;|\leq 
\;c_1\,e^{d(x,\,y)} \;+\; d(x,y)\max_{\pi^{-1}(B(x,\,d(x,\,y)))}|\wt F|\;,
$$
and if furthermore $d(x,y)\leq r_0$, then 
$$
|\;C_{F,\,\xi}(x,y)\;|\leq
c_1\,d(x,y)^{c_2}+ d(x,y)\max_{\pi^{-1}(B(x,\,d(x,\,y)))}|\wt F|\;.
$$

(2) For every $r\in\mathopen{[}0,r_0\mathclose{]}$, for all $x,y'$ in
$\wt M$, for every $\xi$ in the shadow $\OOO_xB(y',r)$ of the ball
$B(y',r)$ seen from $x$, we have
$$
\Big|\;C_{F,\,\xi}(x,y')+\int_x^{y'} \wt F\;\Big|\leq 
c_3 \;r^{c_4}+2r\max_{\pi^{-1}(B(y',\,r))}|\wt F|\;. 
$$
\elemm

\dem 
(1) This follows from Lemma \ref{lem:technicholder} by
letting $z$ tend to $\xi$ along a geodesic ray.

(2) Let $t\mapsto\xi_t$ be the geodesic ray starting from $x$ and
ending at $\xi$, and let $p$ be the closest point on it to $y'$, which
satisfies $d(y',p)\leq r$.  Since $\iota$ is Lipschitz for $d'_{T^1M}$
(see the paragraph before Lemma \ref{lem:dholdequivSasak}) with
constant depending only on the bounds on the sectional curvature, the
assertion (2) follows from the second assertion of Lemma
\ref{lem:technicholder}, using
$$
\int_{x}^{\xi_t}\wt F-\int_{y'}^{\xi_t}\wt F-\int_{x}^{y'}\wt F=
\Big(\int_{p}^{x}\wt F\circ\iota-\int_{y'}^{x}\wt F\circ\iota\Big)+
\Big(\int_{p}^{\xi_t}\wt F-\int_{y'}^{\xi_t}\wt F\Big)\;.
\;\;\;\Box
$$

\medskip \rem By Remark (1) at the end of Subsection
\ref{subsec:potentials}, when $x,y\in\C\Lambda\Ga$ and
$\xi\in\Lambda\Ga$, we may replace $\max_{\pi^{-1}(B(x,\,d(x,y)))}|\wt
F|$ by $\max_{\pi^{-1}(B(x,\,d(x,y))\cap\,\C\Lambda\Ga)}|\wt F|$ in
Assertion (1) of the above lemma.

\bprop \label{prop:gibbscocyholder} (1) The map
$C_F: \partial_\infty\wt M\times \wt M\times\wt M\ra \RR$ is
continuous.

(2) Assume that $\wt F$ is bounded. Then $C_F$ is locally
H\"older-continuous, and the maps $\xi\mapsto C_{F,\,\xi}(x,y)$ for all
$x,y\in\wt M$ and $(x,y)\mapsto C_{F,\,\xi}(x,y)$ for every
$\xi\in\partial_\infty\wt M$ are H\"older-continuous.  
\eprop

In particular, the Busemann cocycle (obtained by taking $F=-1$) is
locally H\"older-continuous.

If we assume only that $\wt F$ is bounded on $\C\Lambda\Ga$, then the
same proof gives that the restriction of $C_F$ to $\Lambda\Ga\times
\wt M\times\wt M\ra \RR$ is locally H\"older-continuous.

\medskip
\dem Let $x,y,x',y'\in\wt M$ and $\xi,\xi'\in\partial_\infty\wt M$.

\medskip
(1) For every $t\geq 0$, let $x_t$ be the point at distance $t$ on the
geodesic ray $\mathopen{[}x,\xi\mathclose{[}$. Let $y_t$, $x'_t$ and
$y'_t$ be the closest points to $x_t$ on respectively
$\mathopen{[}y,\xi\mathclose{[}$, $\mathopen{[}x',\xi'\mathclose{[}$
and $\mathopen{[}y',\xi'\mathclose{[}$. By (a version with
$z\in \partial_\infty\wt M$ of) Lemma \ref{lem:exogeohyp} (i), there
exists $c'_1>0$ (depending only on $d(x,y)$) such that $d(x_t,y_t)\leq
c'_1\,e^{-t}$.

By the linear growth property of H\"older-continuous functions (see the
end of Subsection \ref{subsec:holdercont}), there exists $c'_2>0$
(depending only on $\wt F(x)$ and the H\"older constants of $\wt F$)
such that $\max_{\pi^{-1}(B(x_t,\,2))}|\wt F|\leq c'_2(t+1)$ for every
$t\geq 0$.

Let $\epsilon>0$. Let $c_1,c_2$ be the positive constants so that
Lemma \ref{lem:technicholder} and Lemma \ref{lem:holderconseq} (1)
hold true for $r_0=1$.  Let us fix $t$ large enough (which depends only on
$\epsilon$, $d(x,y)$, $\wt F(x)$ and the H\"older constants of $\wt F$),
so that $2\,c'_1\,e^{-t}\leq 1$ and $c_1(2\,c'_1\,e^{-t})^{c_2}+
2\,c'_1\,e^{-t}\,c'_2(t+1)\leq \epsilon$.

By the continuity property of geodesic rays and closest point maps, if
$\xi'$, $x'$ and $y'$ are close enough to $\xi$, $x$ and $y$
respectively, we have $d(x,x'),d(y,y'),d(x_t,x'_t),d(y_t,y'_t)\leq
c'_1\,e^{-t}$ and $d(x'_t,y'_t)\leq 2\,c'_1\,e^{-t}$. Hence, using the
cocycle property \eqref{eq:cocycleprop} and Equation
\eqref{eq:cocyclealongray} for the following equality,
\begin{align*}
&|\;C_{F,\,\xi'}(x',y')-C_{F,\,\xi}(x,y)\;| \\=\; &
\Big|\;\int_{y'}^{y'_t}\wt F+C_{F,\,\xi'}(x'_t,y'_t)-\int_{x'}^{x'_t}\wt F-
\int_{y}^{y_t}\wt F-C_{F,\,\xi}(x_t,y_t)+\int_{x}^{x_t}\wt F\;\Big|
\\ \leq\; &
\Big|\int_{y'}^{y'_t}\wt F-\int_{y}^{y_t}\wt F\;\Big|+
\Big|\int_{x'}^{x'_t}\wt F-\int_{x}^{x_t}\wt F\;\Big|+
|\;C_{F,\,\xi'}(x'_t,y'_t)\;|+|\;C_{F,\,\xi}(x_t,y_t)\;|\;.
\end{align*}
Then by the definition of $t$ and by Lemma \ref{lem:technicholder} and
Lemma \ref{lem:holderconseq} (1), we have $|\;C_{F,\,\xi'}(x',y')-
C_{F,\,\xi}(x,y)\;|\leq 6\,\epsilon$, and Assertion (1) follows.

\medskip (2) To prove Assertion (2), let us first prove that there
exist constants $c>0$ and $\alpha\in\mathopen{]}0,1\mathclose{[}$
(depending only on $d(x,y)$, on the H\"older constants of $\wt F$, on
the upper bound of $|\wt F|$ and on the bounds on the sectional
curvature of $\wt M$) such that if $d_x(\xi,\xi')\leq
e^{-d(x,\,y)-2}$, then
$$
|C_{F,\,\xi'}(x,y)-C_{F,\,\xi}(x,y)|\leq 
c\; d_x(\xi,\xi')^{\alpha}\;.
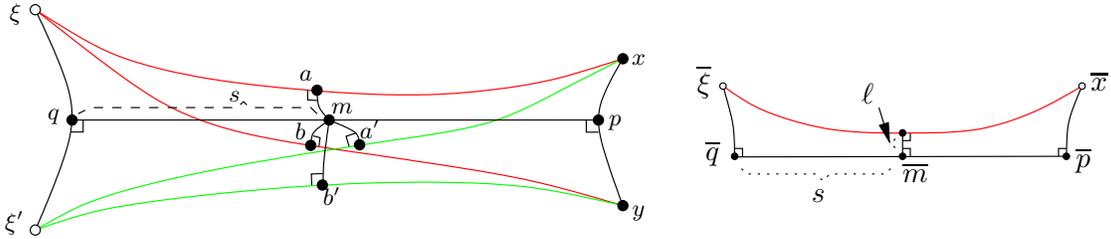
$$
\begin{center}
\input{fig_regholdcocyc.pstex_t}
\end{center}

Let $\mathopen{[}p,q\mathclose{]}$ be the shortest arc between
$\mathopen{[}x,y\mathclose{]}$ and $\mathopen{]}\xi,
\xi'\mathclose{[}\;$, with $p\in \mathopen{[}x,y\mathclose{]}$. Let
$m$ be the midpoint of $\mathopen{[}p,q\mathclose{]}$ and $s=d(p,m) =
d(m,q)=\frac{1}{2}d(p,q)$. Let $a,b,a',b'$ be the closest points to $m$
on respectively $\mathopen{[}x, \xi\mathclose{[}\,$, $\mathopen{[}y,
\xi\mathclose{[}\,$, $\mathopen{[}x, \xi'\mathclose{[}\,$,
$\mathopen{[}y, \xi'\mathclose{[}\,$.

Since $q$ is the closest point to $p$ on $\mathopen{]}\xi,
\xi'\mathclose{[}\,$, and by comparison, the distance from $q$ to
$\mathopen{[}p,\xi\mathclose{[}$ and to $\mathopen{[}p,
\xi'\mathclose{[}$ is at most $\log(1+\sqrt{2})$.  Hence by the
triangle inequality and the definition of the visual distance (see
Equation \eqref{eq:defidistvis}), we have
$$
e^{-2s}\leq d_p(\xi,\xi')\leq e^{-2s +2\log(1+\sqrt{2})}\;.
$$
By the properties of the visual distances (see Equation
\eqref{eq:lipequivdistvis}), we have
$$
e^{-d(x,\,y)}\leq e^{-d(x,\,p)}\leq \frac{d_x(\xi,\xi')}{d_p(\xi,\xi')}
\leq e^{d(x,\,p)}\leq e^{d(x,\,y)}\;.
$$
In particular, 
$$
s\geq -\frac{1}{2}\log d_p(\xi,\xi')\geq
\frac{1}{2}(-\log d_x(\xi,\xi')-d(x,y))\geq 1\;.
$$  
By comparison, since the angle at $q$ between $p$ and $\xi$ is
$\frac{\pi}{2}$, and since, if $p\neq x$, the angle at $p$ between $q$
and $x$ is at least $\frac{\pi}{2}$ (and exactly $\frac{\pi}{2}$ if
$p\neq y$), the distance from $m$ to $a$ is at most the distance
$\ell$ in the real hyperbolic plane $\HH^2_\RR$ between the midpoint
$\ov m$ of a segment $\mathopen{[}\ov p,\ov q\mathclose{]}$ of length
$2s$ to a geodesic line between $\ov\xi,\ov
x\in\partial_\infty\HH^2_\RR$ where the angle at $\ov q$ between $\ov
p$ and $\ov \xi$ and the angle at $\ov p$ between $\ov q$ and $\ov x$
are exactly $\frac{\pi}{2}$. By a well known formula in hyperbolic
geometry (see for instance \cite[page 157]{Beardon83}), we have $\sinh
\ell =\frac{1}{\sinh s}$.  Note that for all $s\geq 1$, we have $\sinh
s\geq \frac{e^s}{4}$.  Hence
$$
d(a,m)\leq \ell\leq \sinh\ell=\frac{1}{\sinh s}\leq 4\;e^{-s}\leq 
4\;d_p(\xi,\xi')^{\frac{1}{2}}\leq 
4\;e^{\frac{d(x,\,y)}{2}}\,d_x(\xi,\xi')^{\frac{1}{2}}\;.
$$
Since the computations are unchanged by permuting $x$ and $y$, as
well as $\xi$ and $\xi'$, we have
$$
\max\{d(a,m),\,d(b,m),\,d(a',m),\,d(b',m)\}\leq 
4\;e^{\frac{d(x,\,y)}{2}}\,d_x(\xi,\xi')^{\frac{1}{2}}\;.
$$
Hence by the triangle inequality,
\begin{equation}\label{eq:majoabcd}
\max\{d(a,b),\,d(a',b'),\,d(a,a'),\,d(b,b')\}\leq 
8\;e^{\frac{d(x,\,y)}{2}}\,d_x(\xi,\xi')^{\frac{1}{2}}\leq 
8\;e^{\frac{d(x,\,y)}{2}}\;.
\end{equation}
Now, using respectively 

$\bullet$~ the cocycle property \eqref{eq:cocycleprop},

$\bullet$~ Equation \eqref{eq:cocyclealongray},

$\bullet$~ the triangle inequality,

$\bullet$~ Lemma \ref{lem:technicholder} and Lemma
\ref{lem:holderconseq} with $r_0=8\;e^{\frac{d(x,\,y)}{2}}$, for some
constants $c_1>0$ and $c_2\in\mathopen{]}0,1\mathclose{[}$ (depending
only on $r_0$, the H\"older constants of $\wt F$, and the bounds on
the sectional curvature of $\wt M$) and $\kappa=\sup_{\wt M}|\wt F|$,

$\bullet$~ Equation \eqref{eq:majoabcd},

\noindent we have
\begin{align*}
&|C_{F,\,\xi'}(x,y)-C_{F,\,\xi}(x,y)|\\
=\;&|-C_{F,\,\xi'}(a',x)+C_{F,\,\xi'}(a',b')+C_{F,\,\xi'}(b',y)+
C_{F,\,\xi}(a,x)-C_{F,\,\xi}(a,b)-C_{F,\,\xi}(b,y)|\\ =\;&
\Big|-\int_x^{a'}\wt F+C_{F,\,\xi'}(a',b')+\int_y^{b'}\wt F
+\int_x^a\wt F-C_{F,\,\xi}(a,b)-\int_y^b\wt F\Big|\\ \leq\;&
 \Big|\int_x^a\wt F-\int_x^{a'}\wt F\;\Big|+|C_{F,\,\xi'}(a',b')|+
|C_{F,\,\xi}(a,b)|+
\Big|\int_y^b\wt F-\int_y^{b'}\wt F\;\Big|
\\ \leq\;& c_1\,d(a,a')^{c_2}+\kappa \,d(a,a')+
c_1\,d(a',b')^{c_2}+\kappa \,d(a',b')
\\ &\;+ c_1\,d(a,b)^{c_2}+\kappa \,d(a,b)+ 
c_1\,d(b,b')^{c_2}+\kappa \,d(b,b')
\\ \leq\;&  8^{1+c_2}\max\{c_1,\kappa\}\; e^{\frac{c_2\,d(x,y)}{2}}\,
d_x(\xi,\xi')^{\frac{c_2}{2}}\;.
\end{align*}

This proves the claim made at the beginning of this proof, as well as
the second claim in the statement of Proposition
\ref{prop:gibbscocyholder} (2).

Now, by Equation \eqref{eq:cocycleprop}, we have 
$$
|C_{F,\,\xi'}(x',y')-C_{F,\,\xi}(x,y)|\leq
|C_{F,\,\xi'}(x,y)-C_{F,\,\xi}(x,y)|+
|C_{F,\,\xi'}(x,x')|+
|C_{F,\,\xi'}(y,y')|\;.
$$
Using the initial claim and again Lemma \ref{lem:holderconseq} (1),
the result follows.  \cqfd

\subsection{The potential gap of $(\Ga,F)$}
\label{subsec:gap}

For all $x\in\wt M$ and $(\xi, \eta)\in \partial_\infty^2 \wt M$, let us
define the {\it $F$-gap seen from $x$}\index{potential!gap seen from a
  point}\index{gap!of potential seen from a point} between $\xi$ and
$\eta$ by

\medskip\noindent
\begin{minipage}{10.2cm}
$$
\gls{potentialgap}=\exp\frac{1}{2}\Big(\lim_{t\ra+\infty} 
\int_x^{\eta_t}\wt F-\int_{\xi_t}^{\eta_t}\wt F+
\int_{\xi_t}^{x}\wt F\Big)\;,
$$ 
where $t\mapsto \xi_t$, $t\mapsto \eta_t$ are any geodesic rays ending
at $\xi$ and $\eta$ respectively (see the picture on the right). By
the H\"older-continuity of $F$ and the properties of asymptotic
geodesic rays in $\wt M$, the limit does exist and is independent of
the choices.
\end{minipage}
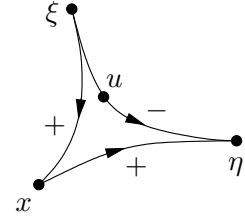
\begin{minipage}{4.7cm}\begin{center}
\input{fig_gapleft.pstex_t}
\end{center}
\end{minipage}

\medskip
We denote by
$$
D_F: \wt M\times \partial_\infty^2 \wt M\ra \RR
$$ 
the map defined by $(x,\xi,\eta)\mapsto D_{F,\,x}(\xi,\eta)$, which we
will call the {\it gap map}\index{gap!map} of the potential $F$.  We
will compare our map $D_F$ with the analogous one introduced by
Hamenst\"adt in \cite{Hamenstadt97} at the end of Subsection
\ref{subsec:crossratio}.  When $F$ is the constant function with value
$-1$, then $D_{F,\,x}(\xi,\eta)={d_x}(\xi,\eta)$, where $d_x$ is the
visual distance on $\partial_\infty \wt M$ seen from $x$ (see
Subsection \ref{subsec:geombus}). So that for every $s\in\RR$, we have
$$
D_{F-s,\,x}(\xi,\eta)=D_{F,\,x}(\xi,\eta)\;{d_x}(\xi,\eta)^{s}\;.
$$

The gap map $D_F$ satisfies the following elementary properties.

\blemm \label{lem:justepourdeffigurecidess}
(1) For any point $u$ on the geodesic line between $\xi$ and
$\eta$, we have
\begin{equation}\label{eq:ecartvisuel}
D_{F,\,x}(\xi,\eta)=
e^{-\frac{1}{2}(C_{F,\,\eta}(x,\,u)+C_{F\circ\iota,\,\xi}(x,\,u))}\;.
\end{equation}
For every  $\ga\in\Ga$,
\begin{equation}\label{eq:equivgap}
D_{F,\,\ga x}(\ga \xi,\ga\eta)=D_{F,\,x}(\xi,\eta) \;\;\;{\rm and}\;\;\; 
D_{F,\,x}(\eta,\xi)=D_{F\circ \iota,\,x}(\xi,\eta)\;.
\end{equation}

(2) If $\wt F$ is bounded, the gap map $D_F: \wt M
\times \partial_\infty^2 \wt M\ra \RR$ is locally H\"older-continuous.  
\elemm

\dem (1) This is immediate, using Equation \eqref{eq:timereversal}.

(2) Let us first prove that there exists $c>0$ such that for all
$(\xi,\eta),(\xi',\eta') \in \partial_\infty^2 \wt M$ and for every
point $u$ on the geodesic line $(\xi\eta)$ between $\xi$ and $\eta$,
if $d_u(\xi,\xi'), d_u(\eta,\eta') \leq e^{-2}$, then there exists
$u'$ on the geodesic line $(\xi'\eta')$ such that
$$
d(u,u')\leq c\;\max\{d_u(\xi,\xi'),d_u(\eta,\eta')\}\;.
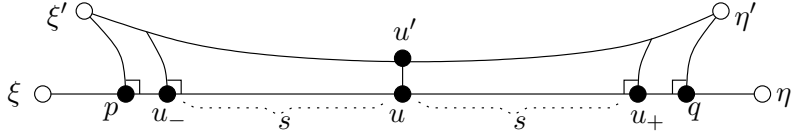
$$
\begin{center}
\input{fig_holdcontgap.pstex_t}
\end{center}
Define 
$$
s=\min\{-\log d_u(\xi,\xi'),-\log d_u(\eta,\eta')\}-
\log(1+\sqrt{2})\geq 1\;.
$$
Let $p$ be the closest point to $\xi'$ on the geodesic ray
$\mathopen{[}u,\xi\mathclose{[}\,$, which by comparison is at distance
at most $\log(1+\sqrt{2})$ from the geodesic line $(\xi\xi')$. By the
triangle inequality and the definition of the visual distance (see
Equation \eqref{eq:defidistvis}), we have $d(u,p)\geq -\log
d_u(\xi,\xi')- \log(1+\sqrt{2})$.  Similarly, the closest point $q$ to
$\eta'$ on $\mathopen{[}u,\eta\mathclose{[}\,$ satisfies $d(u,q)\geq
-\log d_u(\eta,\eta')- \log(1+\sqrt{2})$.  By convexity, the point
$u_-$ at distance $s$ from $u$ on $\mathopen{[}u,\xi\mathclose{[}$ and
the point $u_+$ at distance $s$ from $u$ on $\mathopen{[}u,
\eta\mathclose{[}$ are the orthogonal projections of points on the
extended geodesic line $\mathopen{[}\xi', \eta'\mathclose{]}$.  By
comparison (as in the proof of Proposition \ref{prop:gibbscocyholder}
(2) using the formulae of \cite[page 157]{Beardon83}), the point $u$
is at distance at most $4\,e^{-s}$ from the geodesic line
$(\xi'\eta')$. The initial claim follows.

Now, for all $x,x'\in \wt M$ and $(\xi,\eta),(\xi',\eta')
\in \partial_\infty^2 \wt M$, fix $u$ on the geodesic line
$(\xi\eta)$, and let $u'$ be the point defined above, which exists if
$\xi'$ and $\eta'$ are close enough to $\xi$ and $\eta$ respectively.
By Equation \eqref{eq:ecartvisuel}, we have
$$
|D_{F,\,x}(\xi,\eta)-D_{F,\,x'}(\xi',\eta')|=
|e^{-\frac{1}{2}(C_{F,\,\eta}(x,u)+C_{F\circ\iota,\,\xi}(x,u))}
-e^{-\frac{1}{2}(C_{F,\,\eta'}(x',u')+C_{F\circ\iota,\,\xi'}(x,u'))}|\;.
$$
Assertion (2) of Lemma \ref{lem:justepourdeffigurecidess} then follows
from Proposition \ref{prop:gibbscocyholder} (2).  
\cqfd

\medskip\noindent{\bf Remarks.}  (1) Let $\wt F' :T^1\wt M\ra \RR$ be
a H\"older-continuous $\Ga$-invariant map, which is cohomologous to
$\wt F$ via the map $\wt G :T^1\wt M\ra \RR$. Then for every $x\in \wt
M$, for all $\xi,\eta\in\partial_\infty \wt M$, it follows from the
equations \eqref{eq:ecartvisuel} and \eqref{eq:cohomologuecocycle},
and from the fact that the potentials $\wt F'\circ\iota$ and $\wt
F\circ\iota$ are cohomologous via the map $-\wt G\circ\iota$, that
\begin{equation}\label{eq:gapcoho}
D_{F',\,x}(\xi,\eta)=D_{F,\,x}(\xi,\eta)\;
e^{\frac{1}{2}(\wt G(v_{x\xi})-\wt G\circ\iota (v_{x\eta}))}\;.
\end{equation} 
In particular, if $F$ and $F\circ \iota$ are cohomologous, then by
Equation \eqref{eq:equivgap}, the ratio $\frac{D_{F,\,x}(\eta,\xi)}
{D_{F,\,x}(\xi,\eta)}$ is uniformly bounded in $(\eta,\xi)$, and also
in $(x,\eta,\xi)$ if the function $\wt G$ such that $\wt F$ and $\wt
F\circ\iota$ are cohomologous via  $\wt G$ is bounded.

\medskip (2) By Equation \eqref{eq:ecartvisuel}, by the cocycle
property of $C_F$ and by Lemma \ref{lem:holderconseq} (1), for all
$x,y$ in $\wt M$, there exists a constant $c_{x,\,y}>0$ such that, for
all $\xi,\eta\in\partial_\infty\wt M$, we have
$$
\frac{1}{c_{x,\,y}}\leq\frac{D_{F,\,x}(\xi,\eta)}{D_{F,\,y}(\xi,\eta)}\leq
c_{x,\,y}\;.
$$

\medskip (3) By the thin triangles property of CAT$(-1)$ spaces, if
$p$ is the closest point to $x$ on the geodesic line between $\xi$ and
$\eta$ and $r=\ln(1+\sqrt{2})$, we have $\eta,\xi\in\OOO_x
B(p,r)$. Hence by Lemma \ref{lem:holderconseq} (2), if $\wt F$ is
bounded, there exists a constant $c_5\geq 1$ (depending only on $\|\wt
F\|_\infty$ and on the H\"older constants of $\wt F$) such that,
for all $x\in \wt M$ and $\xi,\eta\in \partial_\infty\wt M$, we have
\begin{equation}\label{eq:compargapdist}
\frac{1}{c_5}\leq \frac{D_{F,\,x}(\xi,\eta)}
{e^{\frac{1}{2}(\int_x^p\wt F+\int_x^p\wt F\circ\iota)}}\leq c_5\;.
\end{equation}
These inequalities are still satisfied, when $\wt F$ is only assumed
to be bounded on $\pi^{-1}(\C\Lambda\Ga)$, if $x\in \C\Lambda\Ga$ and
$\xi,\eta\in \Lambda\Ga$.

\medskip (4) Let us prove in this remark that, under some assumptions
on the potential, a minor modification of the gap map $D_F$ yields a
distance on $\Lambda\Ga$ (see also \cite[Section 2.6]{Schapira03a}).

\blemm\label{lem:gapdistun} Let $x\in\C\Lambda\Ga$. Assume that $\wt
F$ is bounded on $\pi^{-1}(\C\Lambda\Ga)$. Assume that $\wt F$ is
reversible or that there exists a constant $c>0$ such that, for all $y,z\in
\C\Lambda\Ga$ such that $y\in\mathopen{[}x,z\mathclose{]}$, we have
\begin{equation}\label{eq:hypinegtrianultra}
\int_x^z\wt F\leq c+\int_x^y\wt F\;.
\end{equation} 
Then for every $\epsilon>0$ small enough, there exist a distance
$d_{F,\,x,\,\epsilon}$ on $\Lambda\Ga$ and  $c_\epsilon>0$ such
that for all $\xi,\eta\in \Lambda\Ga$, we have
\begin{equation}\label{eq:compardistgap}
\frac{1}{c_\epsilon} \;D_{F,\,x}(\xi,\eta)^\epsilon\leq
d_{F,\,x,\,\epsilon}(\xi,\eta)\leq c_\epsilon\;
D_{F,\,x}(\xi,\eta)^\epsilon\;.
\end{equation} 
\elemm

We will give, in the proof of Proposition \ref{prop:gapdistdeux},
conditions on $(\Ga, F)$ implying that the hypothesis
\eqref{eq:hypinegtrianultra} is satisfied when $F$ is a normalised
potential. Note that the hypothesis \eqref{eq:hypinegtrianultra} does
hold when $\wt F\leq 0$ (for any $\Ga$), as for instance when
$\wt F=\wt F^{\rm su}$ (see Chapter \ref{sec:Liouville}).

\medskip \dem By the techniques of approximation by trees in
$\operatorname{CAT} (-1)$-spaces (see for instance \cite[\S
2.2]{GhyHar90}), there exists a universal constant $c'>0$ such that
for all $x\in \wt M$ and $\xi_1,\xi_2,\xi_3 \in \partial_\infty\wt M$
such that $\xi_1\neq\xi_2$, if $p$ is the closest point to $x$ on the
geodesic line between $\xi_1$ and $\xi_2$, then at least one of the
two following cases holds:

\medskip
\noindent\begin{minipage}{9.4cm}
  $\bullet$~ the distance between $p$ and the closest point to
  $x$ on the geodesic line between $\xi_3$ and one of $\xi_1,\xi_2$ is
  at most $c'$,
\end{minipage}
\begin{minipage}{5.5cm}
\begin{center}
\input{fig_approxtree1.pstex_t}
\end{center}
\end{minipage}

\medskip
\noindent\begin{minipage}{9.4cm}
  $\bullet$~ or the distance between the closest point $q$ to $\xi_3$
  on the segment $\mathopen{[}x,p\mathclose{]}$ and the closest point
  to $x$ on the geodesic line between $\xi_3$ and $\xi_1$ is at most
  $c'$.
\end{minipage}
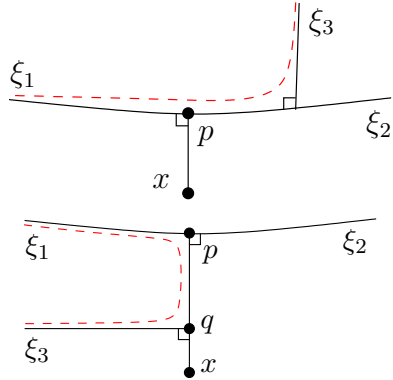
\begin{minipage}{5.5cm}
\begin{center}
\input{fig_approxtree2.pstex_t}
\end{center}
\end{minipage}

\medskip
Note that $p$, as well as $q$ in the second case, is contained in
$\C\Lambda\Ga$ if $\xi_1,\xi_2,\xi_3 \in \Lambda\Ga$.

Note that in the second case, if $\wt F$ and $\wt F\circ \iota$ are
cohomologous via $\wt G$, with $v_{p,\,q}$ a unit tangent vector at
$p$ to the geodesic segment $\mathopen{[}p,q\mathclose{]}$ (unique if
$p\neq q$, with $v_{q,\,p}=-v_{p,\,q}$ if $p=q$), we have, for every
$\ga\in\Ga$, and by the $\Ga$-invariance of $\wt G$,
\begin{equation}\label{eq:majovaetvient}
\Big|\int_{x}^{\ga x}\wt F -\int_{\ga x}^{x}\wt F\;\Big|=
\Big|\int_{x}^{\ga x}(\wt F -\wt F\circ
\iota)\Big|= |\,\wt G(-v_{\ga x,\,x})-\wt G(v_{x,\,\ga x})|
\leq 2\,\max_{T^1_x\wt M}\;|\,\wt G\,|\;.
\end{equation}

Note that Formula \eqref{eq:compargapdist} is also valid, up to
enlarging $c_5$, if $p$ is replaced by any point $p'$ at distance less
than $c'$ from $p$. Indeed, with $F'=F$ or $F'=F\circ\iota$, by Lemma
\ref{lem:technicholder}, the quantity $\big|\int_x^{p'}\wt F'
-\int_x^p\wt F'\big|$ is bounded by a constant depending only on $c'$
and $\kappa=\sup_{\pi^{-1}(\N_{c'}(\C\Lambda\Ga))}\,|\,\wt
F\,|$. Furthermore, $\kappa$ is finite since $|\,\wt F\,|$ is bounded
on $\pi^{-1}(\C\Lambda\Ga)$ and has at most linear growth by the
remark at the end of Subsection \ref{subsec:holdercont}.

Hence, by Equation \eqref{eq:compargapdist}, considering the two above
cases (the first one being easier) and using in the second case either
Formula \eqref{eq:majovaetvient} or Formula
\eqref{eq:hypinegtrianultra}, we see that there exists a constant
$c_6$ such that for all $\xi_1,\xi_2,\xi_3 \in \Lambda\Ga$,
$$
D_{F,\,x}(\xi_1,\xi_2)\leq 
c_6\max\{D_{F,\,x}(\xi_1,\xi_3),D_{F,\,x}(\xi_2,\xi_3)\}\;.
$$
The constant $c_6$ depends on $\sup_{v\in\pi^{-1}(\C\Lambda\Ga)} |\wt
F(v)|$, on the H\"older constants of $\wt F$ and either on
$\max_{v\,\in \,T^1_x\wt M}\;|\,\wt G\,|$ if $\wt F$ and $\wt F\circ
\iota$ are assumed to be cohomologous via $\wt G$ or on the constant
$c$ appearing in Equation \eqref{eq:hypinegtrianultra}.  Note that,
also by Formula \eqref{eq:compargapdist},
$$
D_{F,\,x}(\xi_1,\xi_2)\leq {c_5}^2 D_{F,\,x}(\xi_2,\xi_1)\;.
$$
These weak ultrametric triangle inequality and weak symmetry imply
(see for instance \cite[\S 7.3]{GhyHar90}) that there exists
$c_7\in\mathopen{]}0,1\mathclose{]}$ (depending only on $c_5,c_6$)
such that for every $\epsilon\in \mathopen{]}0,c_7\mathclose{]}$,
there exist a distance $d_{F,\,x,\,\epsilon}$ on $\Lambda\Ga$ and a
constant $c_\epsilon>0$ (depending only on $\epsilon,c_5,c_6$) such
that for all $\xi,\eta\in \Lambda\Ga$, we have
$$
\frac{1}{c_\epsilon} \;D_{F,\,x}(\xi,\eta)^\epsilon\leq
d_{F,\,x,\,\epsilon}(\xi,\eta)\leq c_\epsilon\;
D_{F,\,x}(\xi,\eta)^\epsilon\;.
$$
This proves Lemma \ref{lem:gapdistun}. \cqfd

\medskip
By the above Remark (2), for all $x,y\in\C\Lambda\Ga$, the distances
$d_{F,\,x,\,\epsilon}$ and $d_{F,\,y,\,\epsilon}$ are equivalent. Note that if
$F$ is constant, equal to $-1$, then we may take $c_7=1$ and
$d_{F,\,x,\,\epsilon}={d_x}^{\epsilon}$, where $d_x$ is the visual
distance on $\partial_\infty \wt M$ (see for instance
\cite{Bourdon95}).

Note that if in the assumptions of Lemma \ref{lem:gapdistun}, we
replace $\C\Lambda\Ga$ by $\wt M$, then we get a distance
$d_{F,\,x,\,\epsilon}$ on the full boundary $\partial_\infty\wt M$ instead
of just $\Lambda\Ga$, that satisfies Equation \eqref{eq:compardistgap}.

\subsection{The crossratio of $(\Ga,F)$}
\label{subsec:crossratio}

Let $\partial_\infty^4\wt M$ be the space of pairwise distinct
quadruples of points in $\partial_\infty\wt M$. The {\it
  crossratio}\index{crossratio} of the potential $F$ is the map from
$\partial_\infty^4\wt M$ to $\RR$ defined by
$$
\gls{crossratio}=\exp\;\frac{1}{2}\lim_{t\ra+\infty} 
\Big(\int_{a_t}^{d_t}\wt F
-\int_{b_t}^{d_t}\wt F+\int_{b_t}^{c_t}\wt F-\int_{a_t}^{c_t}\wt F \Big)\;,
$$

\medskip\noindent
\begin{minipage}{9.2cm}
  where $t\mapsto a_t,t\mapsto b_t,t\mapsto c_t,t\mapsto d_t$ are
  geodesic rays converging to $a,b,c,d$ respectively (see the figure
  on the right). The limit does exist and is independent of the
  choices of these geodesic rays, again by the H\"ol\-der-\-continuity
  of $F$ and the properties of asymptotic geodesic rays in $\wt M$.
  An easy cancellation argument shows that, for every $x\in\wt M$,
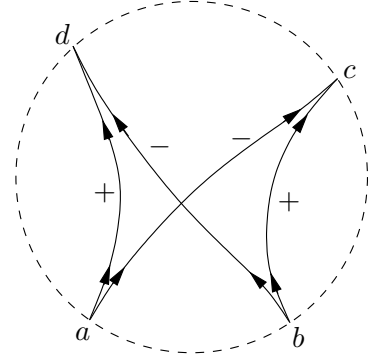
\begin{equation}\label{eq:frombirbir}
\mathopen{[}a,b,c,d\mathclose{]}_F=
\frac{D_{F,\,x}(a,c)D_{F,\,x}(b,d)}{D_{F,\,x}(a,d)D_{F,\,x}(b,c)}\;.
\end{equation}
In particular, the right hand side of this equation does not depend
on $x$. 
\end{minipage}
\begin{minipage}{5.7cm}\begin{center}
\input{fig_gapright.pstex_t}
\end{center}
\end{minipage}

\medskip
We summarise the elementary properties of the  
crossratio of $F$ in the next statement.

\blemm \label{lem:propricrossratio}
(i) The  crossratio of $F$ is a positive map from
$\partial_\infty^4\wt M$ to $\RR$ such that, for every $(a,b,c,d)
\in \partial_\infty^4\wt M$,

(1) {\rm ($\Ga$-invariance)} for every $\ga\in\Ga$, we have 
$\mathopen{[}\ga a,\ga b,\ga c,\ga d\mathclose{]}_F=
\mathopen{[}a,b,c,d\mathclose{]}_F$;

(2) $\mathopen{[}a,b,c,d\mathclose{]}_F=
\mathopen{[}b,a,d,c\mathclose{]}_F$;

(3) $\mathopen{[}a,b,c,d\mathclose{]}_F=
{\mathopen{[}a,b,d,c\mathclose{]}_F}^{-1}$;

(4) $\mathopen{[}a,b,c,d\mathclose{]}_{F\circ\iota}=
\mathopen{[}c,d,a,b\mathclose{]}_F$;

(5) $\mathopen{[}a,b,c,d\mathclose{]}_F=
\mathopen{[}a,b,c,\xi\mathclose{]}_F\;
\mathopen{[}a,b,\xi,d\mathclose{]}_F$ for every
$\xi\in\partial_\infty \wt M-\{a,b,c,d\}$. 

\medskip
\noindent (ii) If the potential $F$ is bounded, then its crossratio is
locally H\"older-continuous.

\medskip
\noindent(iii) If $\wt F':T^1\wt M\ra \RR$ is a H\"older-continuous
$\Ga$-invariant map, cohomologous to $\wt F$, then their associated
crossratios coincide: 
$$
\mathopen{[}a,b,c,d\mathclose{]}_F=\mathopen{[}a,b,c,d\mathclose{]}_{F'}\;.
$$

\medskip
\noindent
(iv) For every loxodromic element $\ga$ in $\Ga$, with attractive and
repulsive fixed points $\ga_+$ and $\ga_-$ respectively on
$\partial_\infty \wt M$, we have
$$
\lim_{n\ra+\infty}\;\frac{1}{n}\ln \,\mathopen{[}\ga_-,
\ga^n\xi,\xi,\ga_+\mathclose{]}_F=
\frac{1}{2}\big(\operatorname{Per}_F\ga + 
\operatorname{Per}_F(\ga^{-1})\big)\;,
$$
uniformly for $\xi$ in a compact subset of $\partial_\infty \wt
M-\{\ga_-,\ga_+\}$.  
\elemm

Assertion (5) is a multiplicative cocycle property in the last
two variables of the crossratio of $F$. Using (4) and (5), it follows
that the  crossratio of $F$ is also a multiplicative cocycle in the
first two variables, that is 
$$
\mathopen{[}a,b,c,d\mathclose{]}_F=\mathopen{[}a,\xi,c,d\mathclose{]}_F
\;\mathopen{[}\xi,b,c,d\mathclose{]}_F
$$
for every $\xi\in\partial_\infty \wt M-\{a,b,c,d\}$.

\medskip \dem The assertions (2), (3) and (4) follow easily from the
definition of the crossratio of $F$. The equations
\eqref{eq:frombirbir} and \eqref{eq:equivgap} imply the assertions (1)
and (5). Assertion (ii) follows from Equation \eqref{eq:frombirbir}
and Lemma \ref{lem:justepourdeffigurecidess} (2). The equations
\eqref{eq:frombirbir} and \eqref{eq:gapcoho} imply Assertion (iii).

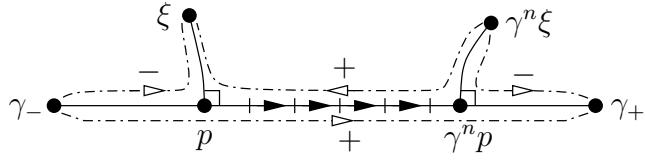
\begin{center}
\input{fig_periocrossratio.pstex_t}
\end{center}

To prove the last assertion, let us define $p$ to be the closest point to
$\xi$ on the translation axis $\Axe_\ga$ of $\ga$. Note that $\ga^np$
is then the closest point to $\ga^n\xi$ on $\Axe_\ga$. By
hyperbolicity, the four geodesic lines along which one integrates $\wt
F$ to define the crossratio are at bounded distance from the union of
$\Axe_\ga$ and of the geodesic rays $\mathopen{]}\xi,p\mathclose{]}$
and $\mathopen{]}\ga^n\xi,\ga^np\mathclose{]}$.  Furthermore, the
segment between $p$ and $\ga^np$ stays very close to the geodesic line
between $\xi$ and $\ga^n\xi$ except for a segment of bounded length
near each of its endpoints. By definition of the periods,
$$
\int_p^{\ga^np}\wt F =n\operatorname{Per}_F\ga\;\;\;
{\rm and}\;\;\; \int_{\ga^n p}^{p}\wt F
=\int_{p}^{\ga^{-n} p}\wt F =n \operatorname{Per}_F (\ga^{-1})
$$
After dividing by $n$, the other contributions vanish as $n$ goes to
$+\infty$, by Lemma \ref{lem:holderconseq} (1) and since $\wt F$ is
bounded on the $\log(1+\sqrt{2})$-neighbourhood of $\Axe_\ga$. The
result follows.  \cqfd

\medskip \noindent{\bf Remarks. }  (1) Let us compare our maps
$C_{F,\,\xi},D_{F,\,x},\mathopen{[}\cdot,\cdot,\cdot,\cdot\mathclose{]}_F$
with those introduced in Hamenst\"adt's paper \cite{Hamenstadt97}
(which assumes that $\Ga$ is torsion free, $M=\Ga\backslash \wt M$ is
compact and $F$ and $F\circ \iota$ are cohomologous, requirements we
do not assume in this book).  Let $\zeta:T^1\wt M\times\RR\ra \RR$ be
the map defined by
$$
\zeta(v,t)= \int_{0}^{t}\wt F(\phi_s v)\,ds\;.
$$ 
This map is locally H\"older-continuous, invariant under $\Ga$ (acting
trivially on $\RR$) and satisfies $\zeta(v,s+t)=
\zeta(v,s)+\zeta(\phi_s(v),t)$. Then Hamenst\"adt's map $k_\zeta:\wt
M\times \wt M\times \partial_\infty\wt M\ra \RR$ defined in Lemma 1.1
in loc.~cit.~has the opposite sign to our Gibbs cocycle: for all
$x,y\in\wt M$ and $\xi\in\partial_\infty\wt M$, we have
\begin{equation}\label{eq:comparhamencocycl}
k_{\zeta}(x,y,\xi)=-C_{F,\,\xi}(x,y)\;.
\end{equation}
Let 
$$
DT\wt M=\{(v,w)\in T^1\wt M\times T^1\wt M\;:\;\pi(v)=\pi(w)\;\;
{\rm and}\;\;v\neq w\}\;.
$$ 

\smallskip\noindent
\begin{minipage}{9.7cm}
  Hamenst\"adt constructs (just above Lemma 1.2 in loc.~cit.) a
  locally H\"older-continuous $\Ga$-invariant map $\alpha_\zeta:DT\wt
  M\ra \RR$, and it may be proved using the equations
  \eqref{eq:ecartvisuel} and \eqref{eq:comparhamencocycl} that when
  $(v,w)\in DT\wt M$ (see the picture on the right),
$$ 
\alpha_\zeta(v,w)-\ln D_{F,\,\pi(v)}(v_+,w_+)
$$
is uniformly bounded (and even equal to $0$ if $F\circ\iota=F$).
\end{minipage}
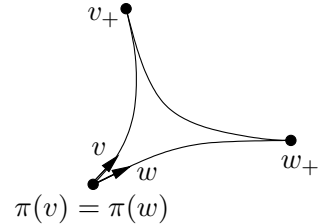
\begin{minipage}{4.7cm}\begin{center}
\input{fig_gapmiddle.pstex_t}
\end{center}
\end{minipage}

\bigskip (2) Assume in this remark that $F$ and $F\circ\iota$ are
cohomologous. It follows from the assertions (4), (iii) and (iv) of
the above Lemma \ref{lem:propricrossratio} and from the last equality
in Equation \eqref{eq:periodinvconj} that the crossratio of $F$ then
satisfies the extra symmetry
$$
\mathopen{[}a,b,c,d\mathclose{]}_{F}=\mathopen{[}c,d,a,b\mathclose{]}_F
\;\;{\rm for~every}\;\; (a,b,c,d) \in \partial_\infty^4\wt M\;,
$$ 
and that
$$
\lim_{n\ra+\infty}\;\frac{1}{n}\ln\; \mathopen{[}\ga_-,\ga^n\xi,\xi,
\ga_+\mathclose{]}_F=\operatorname{Per}_F\ga\;.
$$ 
In particular, the crossratio of $F$ then determines the periods of
$F$. By Remark \ref{rem:livsi}, the crossratio of the potential
function $F$ determines the cohomology class of the restriction of $F$
to the (topological) non-wandering set $\Omega\Ga$. When $\Ga$ is
torsion free and cocompact, this reproves one inclusion in
\cite[Theo.~A]{ Hamenstadt97}.

\subsection{The Patterson densities of $(\Ga,F)$}
\label{subsec:gibbspattersondens}

Let $\sigma\in\RR$. A {\it Patterson density}\index{Patterson density}
of dimension $\sigma$ for $(\Ga,F)$ is a family of finite nonzero
(positive Borel) measures $(\mu_{x})_{x\in\wt M}$ on $\partial_\infty
\wt M$, such that, for every $\ga\in\Ga$, for all $x,y\in\wt M$, for
every $\xi\in\partial_\infty \wt M$, we have
\begin{equation}\label{eq:equivdensity}
\ga_*\mu_{x}=\mu_{\ga x}\;,
\end{equation}
\begin{equation}\label{eq:radonykodensity}
d\mu_{x}(\xi)=e^{-C_{F-\sigma,\,\xi}(x,\,y)}\;d\mu_{y}(\xi)\;.
\end{equation}
In particular, the measures $\mu_{x}$ for $x\in\wt M$ are in the same
measure class. The Radon-Nikodym derivative $\frac{d\mu_{x}}
{d\mu_{y}}$ is only defined almost everywhere, but by the continuity
of the map $\xi\mapsto C_{F-\sigma,\,\xi}(x,\,y)$ (see Proposition
\ref{prop:gibbscocyholder} (1)), we may and we will take
$\frac{d\mu_{x}}{d\mu_{y}} (\xi)=e^{-C_{F-\sigma,\,\xi}(x,\,y)}$ for
every $\xi\in\partial_\infty \wt M$.

Note that a Patterson density of dimension $\sigma$ for $(\Ga,F)$ is
also a Patterson density of dimension $\sigma+s$ for $(\Ga,F+s)$,
for every $s\in\RR$. If $F=0$, then we get the usual notion of a
Patterson density for the group $\Ga$, see for instance 
\cite{Patterson76,Sullivan79,Nicholls89,Bourdon95,Roblin03}.

Proposition \ref{prop:gibbscocyholder} (1) (or Lemma
\ref{lem:holderconseq} (1)) implies that the map $x\mapsto \mu_x$ from
$\wt M$ to the space of finite measures on $\partial_\infty \wt M$,
endowed with the weak-star topology, is continuous, as for every
$x\in\wt M$ and $\epsilon>0$, if $y$ is close enough to $x$, then
$\frac{d\mu_x}{d\mu_y}\in\mathopen{[}e^{-\epsilon},
e^{\epsilon}\mathclose{]}$ by Equation \eqref{eq:radonykodensity}.

The support of the measure $\mu_x$ is independent of $x\in \wt M$ by
Equation \eqref{eq:radonykodensity}. It is a closed nonempty
$\Ga$-invariant subset of $\partial_\infty \wt M$ by Equation
\eqref{eq:equivdensity}. Hence it contains the limit set $\Lambda\Ga$
of $\Ga$. We will be mostly interested in the case when the support of
$\mu_x$ is equal to $\Lambda\Ga$, as this is the case when $\sigma=
\delta_{\Ga,\,F}<+\infty$ and $(\Ga,F)$ is of divergence type, by
Corollary \ref{coro:uniqpatdens}. We will prove in Corollary
\ref{coro:descripdensconfsuppdifflimset} that there exists a Patterson
density of dimension $\sigma$ for $(\Ga,F)$ whose support contains
strictly $\Lambda\Ga$ if and only if $\Lambda\Ga\neq\partial_\infty
\wt M$ and the Poincar\'e series $\sum_{\ga\in\Ga}e^{\int_x^{\ga x}
  (\wt F-\sigma)}$ converges. The basic existence result is the
following one.

\bprop[Patterson] \label{prop:existPattdens} If
$\delta_{\Ga,\,F}<\infty$, then there exists at least one Patterson
density of dimension $\delta_{\Ga,\,F}$ for $(\Ga,F)$, with support
exactly equal to $\Lambda\Gamma$.  
\eprop

\dem This is proved in \cite{Mohsen07} (see also \cite{Hamenstadt97}),
following Patterson's construction.  Since we will need the
construction, we recall it here.  

Fix a point $y$ in $\wt M$. For every $z\in\wt M$, let $\D_z$ be the
unit Dirac mass at $z$. Let $h:\mathopen{[}0,+\infty\mathclose{[}\;\ra
\mathopen{]}0,+\infty\mathclose{[}$ be a nondecreasing map such that

$\bullet$~ for every $\epsilon>0$, there exists $r_\epsilon\geq 0$ such
that $h(t+r)\leq e^{\epsilon t}h(r)$ for all $t\geq 0$ and $r\geq
r_\epsilon$;

$\bullet$~ if $\overline{Q}_{x,\,y}(s)=\sum_{\ga\in\Ga} \;
e^{\int_x^{\ga y} (\wt F-s)}\;h(d(x,\ga y))$, then $\overline{Q}_{x,\,y}(s)$
diverges if and only if the inequality $s\leq \delta_{\Ga,\,F}$ holds.

If $(\Ga,F)$ is of divergence type, we may take $h=1$ constant. Define
the measure
$$
\mu_{x,\,s}= \frac{1}{\overline{Q}_{y,\,y}(s)}\sum_{\ga\in\Ga} \;
e^{\int_x^{\ga y} (\wt F-s)}\;h(d(x,\ga y))\;\D_{\ga y}
$$ 
on $\wt M$. It is proved in \cite{Mohsen07} that such a map $h$
exists, and that there exists a sequence $(s_k)_{k\in\NN}$ in
$\mathopen{]}\delta_{\Ga,\,F},+\infty\mathclose{[}$ converging to
$\delta_{\Ga,\,F}$ such that for every $x\in\wt M$, the sequence of
measures $(\mu_{x,\,s_k})_{k\in\NN}$ weak-star converges to a measure
$\mu_x$ on the compact space $\wt M\cup\partial_\infty \wt M$, such
that $(\mu_{x})_{x\in\wt M}$ is a Patterson density of dimension
$\delta_{\Ga,\,F}$ for $(\Ga,F)$. Since $\overline{Q}_{y,\,y}
(\delta_{\Ga,\,F})=+\infty$ and since the support of $\mu_{x,\,s}$ in
$\wt M\cup \partial_\infty\wt M$ is equal to $\Ga y\,\cup
\,\Lambda\Ga$, the support of $\mu_{x}$ is contained in $\Lambda\Ga$,
hence equal to $\Lambda\Ga$. 
\cqfd

\medskip We will often denote by $(\mu_{F,\,x})_{x\in\wt M}$ a density
as in Proposition \ref{prop:existPattdens}. We will prove in Subsection
\ref{subsec:uniqueness} that if $(\Ga,F)$ is of divergence type, then
$(\mu_{F,\,x})_{x\in\wt M}$ is unique up to a scalar multiple. Note
that the definition of this density involves only the normalised
potential $F-\delta_{\Ga,\,F}$.  Hence the assumption that
$\delta_{\Ga,\,F}$ is positive is harmless, by Equation
\eqref{eq:invarpressuretranslat}, and we will often make it in Chapter
\ref{sec:growth} and Chapter \ref{sec:ergtheounistabfolia}.

\medskip The main basic tool concerning Patterson densities, which
will be very useful, is the following lemma, proved in \cite{Mohsen07}
(see also \cite[Lem.~4]{Coudene03} with the multiplicative rather than
additive convention) along the lines of Sullivan's shadow lemma (see
\cite[p.~10]{Roblin03}).  We will prove a more general result in
Proposition \ref{prop:principombres} in Subsection
\ref{subsec:shadowprinciple}.

\blemm [Mohsen's shadow lemma] \label{lem:shadowlemma} Let
$\sigma\in\RR$, let $(\mu_{x})_{x\in\wt M}$ be a Patterson density
of dimension $\sigma$ for $(\Ga,F)$, and let $K$ be a compact subset of
$\wt M$.  If $R$ is large enough, there exists
$C>0$ such that for all $\ga\in\Ga$ and  $x,y\in K$,
$$
\frac{1}{C}\;e^{\int_x^{\ga y} (\wt F-\sigma)}\leq 
\mu_{x}\big(\OOO_xB(\ga
y,R)\big)\leq C\;e^{\int_x^{\ga y} (\wt F-\sigma)}\;.
$$
\elemm

\dem Let $\sigma$, $(\mu_{x})_{x\in\wt M}$ and $K$ be as in the
statement. Let us prove that if $R'$ is large enough, there exists
$C'>0$ such that for all $\ga\in\Ga$ and $x,y\in K$, we have
\begin{equation}\label{eq:ombrerreduc}
\frac{1}{C'}
\leq \mu_{\ga y}(\OOO_xB(\ga y,R'))
\leq C'\;.
\end{equation} 
Assuming this, let us prove Lemma \ref{lem:shadowlemma}.  By Equation
\eqref{eq:radonykodensity}, we have
$$
\mu_{x}(\OOO_xB(\ga y,R'))=\int_{\xi\in\OOO_xB(\ga y,\,R')} 
e^{-C_{F-\sigma,\,\xi}(x,\,\ga y)}\;d\mu_{\ga y}(\xi)\;.
$$ 
By Lemma \ref{lem:holderconseq} (2) applied with $r=r_0=R'$, since
$\wt F$ is continuous and $\Ga$-invariant, hence is uniformly bounded
on $\Ga\pi^{-1}(\N_{R'}K))$, which contains $\pi^{-1}(B(\ga y,R'))$
for all $\ga\in\Ga$ and $y\in K$, there exists $C''>0$ such that for
all $x,y\in K$ and $\xi\in\OOO_xB(\ga y,R')$,
$$
\big|\,C_{F-\sigma,\,\xi}(x,\ga y)+
\int_x^{\ga y} (\wt F-\sigma)\;\big|\leq C''\;.
$$
Hence 
\begin{multline*}
e^{-C''}\,e^{\int_x^{\ga y} (\wt F-\sigma)}\;\mu_{\ga y}(\OOO_xB(\ga y,R'))
\\
\leq \mu_{x}(\OOO_xB(\ga y,R'))
\\ \leq 
e^{C''}\,e^{\int_x^{\ga y} (\wt F-\sigma)}\;\mu_{\ga y}(\OOO_xB(\ga y,R'))\;.
\end{multline*}
Therefore the result follows with $R=R'$ and $C=C'e^{C''}$.

\medskip
Let us now prove the upper bound in Equation \eqref{eq:ombrerreduc}.
Fix $y_0\in K$, and let
$$
C'''=\sup_{x,\,y\in K,\;\xi\in \partial_\infty\wt M}\;|C_{F-\sigma,\,\xi}(x,y)|\;,
$$ 
which is finite by Lemma \ref{lem:holderconseq} (1), since $K$ is
compact. Then, using Equation \eqref{eq:equivdensity} for the equality
and Equation \eqref{eq:radonykodensity} for the last inequality, we
have
$$
\mu_{\ga y}(\OOO_xB(\ga y,R'))
\leq \|\mu_{\ga y}\|=\|\mu_{y}\| \leq e^{C'''}\;\|\mu_{y_0}\|\;,
$$ 
and the upper bound holds if $C'\geq e^{C'''}\;\|\mu_{y_0}\|$.

\medskip 
Finally, to prove the lower bound in Equation \eqref{eq:ombrerreduc},
we assume for a contradiction that there exist sequences
$(x_i)_{i\in\NN}$ and $(y_i)_{i\in\NN}$ in $K$, $(\ga_i)_{i\in\NN}$ in
$\Ga$, and $(R_i)_{i\in\NN}$ in $\mathopen{]}0,+\infty\mathclose{[}$
converging to $+\infty$, such that
\begin{equation}\label{eq:hypoabsurombr}
\lim_{i\ra+\infty} \;
\mu_{\ga_i y_i}(\OOO_{x_i}B(\ga_i y_i,R_i))\;=0\;.
\end{equation}
Up to a subsequence, we have $\lim x_i=x\in K$, $\lim y_i= y\in K$ and
$\lim\ga_i^{-1} x_i =\xi\in\wt M\cup\partial_\infty\wt M$. By Equation
\eqref{eq:hypoabsurombr}, the shadow $\OOO_{x_i}B(\ga_i y_i,R_i)$ is
different from $\partial_\infty\wt M$, hence $d(x_i,\ga_i y_i)\geq
R_i$, which tends to $+\infty$ as $i\ra+\infty$. Therefore
$\xi\in\partial_\infty\wt M$.

Let $V$ be a relatively compact open subset of $\partial_\infty\wt
M-\{\xi\}$.  Then $V$ is contained in $\OOO_{\ga_i ^{-1}x_i}
B(y_i,R_i)$ for $i$ large enough, since $\lim R_i=+\infty$ and
$(y_i)_{i\in\NN}$ stays in the compact subset $K$ (this is the
original key remark on shadows by Sullivan). For every $i\in \NN$, we
have, using Equation \eqref{eq:equivdensity} for the first equality
and the equivariance property of shadows for the second one,
\begin{align*}
\mu_y(V)&=(\ga_i^{-1})_*\mu_{\ga_i y}(V)\leq 
(\ga_i^{-1})_*\mu_{\ga_i y}(\OOO_{\ga_i ^{-1}x_i} B(y_i,R_i))
=\mu_{\ga_i y}(\OOO_{x_i} B(\ga_i y_i,R_i))\\ &\leq
e^{C'''}\mu_{\ga_i y_i}(\OOO_{x_i} B(\ga_i y_i,R_i))\;,
\end{align*} 
this last inequality by the invariance property
\eqref{eq:equivcocprop} of the Gibbs cocycle and by Equation
\eqref{eq:radonykodensity}.  By Equation \eqref{eq:hypoabsurombr},
letting $i\ra+\infty$, we hence have $\mu_y(V)=0$.  Therefore the
measure $\mu_y$ is supported in $\{\xi\}$. Since the support of a
Patterson density of $(\Ga,F)$ is $\Ga$-invariant, this implies that
$\xi$ is fixed by $\Ga$. This is a contradiction since $\Ga$ is
non-elementary.  
\cqfd

\medskip A first corollary of Mohsen's shadow lemma shows the
inexistence of a Patterson density of dimension less than the critical
exponent of $(\Ga,F)$. Its proof is similar to the one given by
Sullivan \cite{Sullivan79} in the case $F=0$ and constant curvature
(see also \cite[page 147]{Roblin02}).

\bcoro \label{coro:dimgeqpress} Let $\sigma\in\RR$ and let
$(\mu_{x})_{x\in\wt M}$ be a Patterson density of dimension $\sigma$
for $(\Ga,F)$.

(1) For all $x,y\in \wt M$, there exists $c>0$ such that for every
$n\in\NN$, we have
$$
\sum_{\ga\in\Ga\;:\;n-1< d(x,\,\ga y)\leq n}\;
e^{\int_x^{\ga y} (\wt F-\sigma)}\leq c\;.
$$

(2) We have $\sigma\geq \delta_{\Ga,\,F}$.
\ecoro

In particular, if $\delta_{\Ga,\,F}=+\infty$, there exists no
Patterson density for $(\Ga,F)$.

\medskip
We will describe in Proposition \ref{prop:classifdimdens} the exact
set of elements $\sigma\in\RR$ for which there exists a Patterson
density of dimension $\sigma$ for $(\Ga,F)$ with support equal to
$\Lambda\Ga$.

\medskip \dem Let $x,y\in \wt M$. Let $R$ and $C$ be as in Mohsen's
shadow lemma \ref{lem:shadowlemma} for $(\mu_{x})_{x\in\wt M}$ and
$K=\{x,y\}$. Let $\kappa=\card\{\ga\in\Ga\;:\;d(y,\ga y)\leq 1+3R\}$,
which is finite by discreteness. For every $n\in\NN$, let
$\Ga_n=\{\ga\in\Ga\;:\;n-1< d(x,\ga y)\leq n\}$.

If $\ga\in \Ga_n$ and $\xi\in\OOO_xB(\ga y,R)$, if $p_\ga$ is the
closest point to $\ga y$ on the geodesic ray from $x$ to $\xi$, then
$n-R-1\leq d(x,p_\ga)\leq n$ by the triangle inequality. Hence, for
all $\ga,\ga'\in\Ga_n$, if $\xi\in \OOO_xB(\ga y,R)\cap \OOO_xB(\ga'
y,R)$, then, again by the triangle inequality,
$$
d(\ga' y,\ga y)\leq 
d(\ga' y, p_{\ga'})+\big|d(x,p_{\ga'})-d(x, p_{\ga})\big|+
d(p_{\ga},\ga y)\leq 1+3R\;.
$$
Thus, for every $n\in\NN$, a point $\xi\in \partial_\infty\wt M$
belongs to at most $\kappa$ subsets $\OOO_xB(\ga y,R)$ as $\ga$ ranges
over $\Ga_n$. Therefore
$$
\sum_{\ga\in\Ga_n}\mu_{x}\big(\OOO_xB(\ga y,R)\big)\leq 
\kappa\;\mu_{x}\Big(\bigcup_{\ga\in\Ga_n}\OOO_xB(\ga y,R)\Big)\;.
$$

Now, by the lower bound in Lemma \ref{lem:shadowlemma}, and since
$\mu_x$ is a finite measure,
\begin{align}\nonumber
\sum_{\ga\in\Ga_n}\;e^{\int_x^{\ga y} \wt F}&\leq
e^{\sigma\,n}\sum_{\ga\in\Ga_n}e^{\int_x^{\ga y} (\wt F-\sigma)}\leq
C\;e^{\sigma\,n}\sum_{\ga\in\Ga_n}\mu_{x}\big(\OOO_xB(\ga y,R)\big)
\\&\label{eq:grandocontage}
\leq C\kappa\;\|\mu_{x}\|\;e^{\sigma\,n}\;.
\end{align}

The first assertion follows. By taking the logarithm, dividing by $n$
and taking the upper limit as $n$ goes to $+\infty$, the second
assertion then follows, by the definition of the critical exponent of
$(\Ga,F)$.  \cqfd

\medskip Let us give another consequence of Mohsen's shadow lemma
\ref{lem:shadowlemma}, concerning the doubling properties of the
Patterson densities.  Recall (see for instance \cite[page
3]{Heinonen01}) that a measured metric space $(X,d,\mu)$, that is a
metric space $(X,d)$ endowed with a Borel positive measure $\mu$, is
{\it doubling}\index{doubling measure}\index{measure!doubling} if
there exists $c\geq 1$ such that for all $x\in X$ and $r>0$,
$$
\mu(B(x,2\,r))\leq c\,\mu(B(x,r))
$$
where $B(x,r)$ is the ball of centre $x$ and radius $r$ in $X$. 
Note that, up to changing $c$, the number $2$ may be replaced by 
any constant larger than $1$.

\bprop\label{prop:doubling} Let $\sigma\in\RR$ and let
$(\mu_{x})_{x\in\wt M}$ be a Patterson density of dimension $\sigma$
for $(\Ga,F)$. 

(1) For every compact subset $K$ of $\wt M$, for every $R>0$ large
enough, there exists $C=C(R)>0$ such that for all $\ga \in \Ga$ and
$x,y\in K$, we have
$$
\mu_x(\OOO_{x}B(\ga y,5R))\le C\,\mu_x(\OOO_xB(\ga y,R))\;.
$$

(2) For every $x\in\wt M$, if $\Ga$ is convex-cocompact, then the
measured metric space $(\Lambda\Ga, d_x, \mu_x)$ is doubling.  
\eprop

\dem (1) This is immediate from Mohsen's shadow lemma
\ref{lem:shadowlemma}

(2) Let $x\in\wt M$. Recall (see the beginning of Subsection
\ref{subsec:geombus}) that $d_x$ is the visual distance on
$\partial_\infty \wt M$ seen from $x$, and we again denote by $d_x$
its restriction to $\Lambda\Ga$. Denote by $B_x(\xi,r)$ the ball of
centre $\xi$ and radius $r$ for $d_x$.  For every $\xi
\in \partial_\infty \wt M$, let $\rho_\xi:\mathopen{[}0, +\infty
\mathclose{[}\ra \wt M$ be the geodesic ray from $x$ to $\xi$.

By for instance \cite[Lem.~2.1]{HerPau10}, for every $s>0$, there
exists $a(s)\geq 1$ such that for every $t$ large enough, for every
$\xi\in\partial_\infty \wt M$, we have
$$
B(\xi, s\, e^{-t})  \subset \OOO_x B(\rho_\xi(t), s) 
\subset B(\xi,  a(s) e^{-t})\;.
$$

If $\Ga$ is convex-cocompact, there exists $\Delta>0$ such that for
all $\xi\in \Lambda\Ga$ and $t\geq 0$, there exists $\ga_t \in \Ga$
with $d(\rho_\xi(t),\ga_t x)\leq \Delta$. Let $m\in \NN-\{0\}$ be such
that $5^m \geq 2 \,a(1+2\Delta)$. Let $C'=\prod_{i=0}^{m-1} C(5^i)$. For
every $r>0$ small enough, let $t=-\log r +\log a(1+2\Delta)$, which is
large enough. Then 
$$
B_x(\xi,2\,r)\subset B_x(\xi, 5^m e^{-t}) \subset
\OOO_x B(\rho_\xi(t), 5^m) \subset \OOO_x B(\ga_t x , 5^m(1+\Delta))\;,
$$
and 
$$
\OOO_x B(\ga_t x, 1+\Delta) \subset \OOO_x B(\rho_\xi(t), 1+2\Delta)
\subset B(\xi,  a(1+2\Delta) e^{-t}) = B(\xi,r)\;.
$$
Using $m$ times Assertion (1), we have
$$
\mu_x\big(\OOO_x B(\ga_t x , 5^m(1+\Delta))\big) \leq
C' \mu_x\big(\OOO_x B(\ga_t x , 1+\Delta)\big)\;.
$$
Hence if $r>0$ is small enough, for every $\xi\in \Lambda\Ga$, we have
$$
\mu_x (B(\xi,2\,r)) \leq C'  \mu_x(B(\xi,r))\;.
$$ 
By compactness and since the support of $\mu_x$ is $\Lambda\Gamma$,
this proves the result.  \cqfd

\medskip
Here is yet another consequence of Mohsen's shadow lemma
\ref{lem:shadowlemma}, concerning the closeness of the gap map defined
in Subsection \ref{subsec:gap} to an actual distance.

\bprop\label{prop:gapdistdeux}  
Assume that $\delta=\delta_{\Ga,\,F}<+\infty$.  Let $x\in\C\Lambda
\Ga$. Assume that $\Ga$ is convex-cocompact or that $\wt F$ is
reversible and $\wt F$ is bounded on $\pi^{-1}(\C\Lambda\Ga)$. Then
for every $\epsilon>0$ small enough, there exist a distance
$d_{F-\delta,\,x,\,\epsilon}$ on $\Lambda\Ga$ and $c_\epsilon>0$ such
that for all $\xi,\eta\in \Lambda\Ga$, we have
$$
\frac{1}{c_\epsilon} \;D_{F-\delta,\,x}(\xi,\eta)^\epsilon\leq
d_{F-\delta,\,x,\,\epsilon}(\xi,\eta)\leq c_\epsilon\;
D_{F-\delta,\,x}(\xi,\eta)^\epsilon\;.
$$
Furthermore, if $\sup_{T^1\wt M}\;\wt F <\delta$, then
$d_{F-\delta,\,x,\,\epsilon}$ induces the original topology on
$\Lambda\Ga$.  
\eprop

We do not know if the assumption $\sup_{T^1\wt M}\;\wt F <\delta$
is necessary for the last claim.

\medskip
\dem The case $\wt F$ is reversible and $\wt F$ is bounded on
$\pi^{-1}(\C\Lambda\Ga)$ follows from Lemma \ref{lem:gapdistun}.
Assume that $\Ga$ is convex-cocompact.
Then $\wt F$ is bounded on $\pi^{-1}(\C\Lambda\Ga)$, say by $\kappa
\geq 0$. Let $\kappa'$ be the diameter of $\Ga\bs \C\Lambda\Ga$.  For
all $x,y,z\in \C\Lambda\Ga$ such that $y\in \mathopen{[}x,
z\mathclose{]}$, let $\alpha,\beta \in \Ga$ be such that $d(y,\alpha
x),d(z,\beta x)\leq \kappa'$.  Then by Lemma \ref{lem:technicholder}
with $r_0=\kappa'$, there exists $c>0$ (depending only on $\kappa$,
$\kappa'$, the H\"older constants of $\wt F$ and the bounds on the
sectional curvature) such that, with $\wt F'=\wt F-\delta$,
\begin{align*}
\int_x^z \wt F'&=\int_x^y \wt F'+\int_y^z \wt F'\\ &
\leq \int_x^y \wt F'+ \int_{\alpha x}^{\beta x}\wt F'+
\Big|\int_y^z \wt F'-\int_{\alpha x}^z \wt F'\Big|+
\Big|\int_{\alpha x}^z \wt F'-\int_{\alpha x}^{\beta x} \wt F'\Big|
\\ & \leq \int_x^y \wt F'+\int_{\alpha x}^{\beta x}\wt F'+c\;
\leq \int_x^y \wt F'+\;
\sup_{\ga\in\Ga} \;\int_{x}^{\ga x}\wt F'+c\;.
\end{align*}
The existence of a Patterson density of dimension $\sigma=
\delta_{\Ga,\,F} <+\infty$, its finiteness, and the lower bound in
Mohsen's shadow lemma imply that for all $x,y\in\wt M$,
$$
\sup_{\ga\in\Ga} \int_x^{\ga y} (\wt F-\delta)<+\infty\;.
$$
This proves that Equation \eqref{eq:hypinegtrianultra} holds for the
normalised potential $F-\delta$.  Hence, the first claim of
Proposition \ref{prop:gapdistdeux} follows by Lemma
\ref{lem:gapdistun}.

\medskip Let us now prove the final claim of Proposition
\ref{prop:gapdistdeux}. The function $\wt F-\delta$ is bounded from
below and from above by a negative constant. Hence by Equation
\eqref{eq:compargapdist}, there exist $c,c'>0$ such that ${d_x}^c\leq
D_{F-\delta,\,x}\leq {d_x}^{c'}$ on $\Lambda\Ga\times \Lambda\Ga$. The
claim that the distance $d_{F-\delta,\,x,\,\epsilon}$ on $\Lambda\Ga$
induces its topology therefore follows from the fact that the visual
distance $d_x$ induces the topology on $\partial_\infty\wt M$. \cqfd


\bigskip
\noindent{\bf Remark. }  Let $\wt F^* :T^1\wt M\ra \RR$ be a
H\"older-continuous $\Ga$-invariant map, which is cohomologous to $\wt
F$ via the map $\wt G :T^1\wt M\ra \RR$. Let $\sigma\in\RR$ and let
$(\mu_{x})_{x\in\wt M}$ be a Patterson density of dimension $\sigma$
for $(\Ga,F)$. Then the family of measures $(\mu^*_{x})_{x\in\wt M}$
defined by setting, for all $x\in\wt M$ and  $\xi\in\partial_\infty\wt M$,
$$
d\mu^*_x(\xi)=e^{-\wt G(v_{x\xi})}\;d\mu_x(\xi)\;,
$$
where $v_{x\xi}$ is the tangent vector at $x$ of the geodesic ray from
$x$ to $\xi$, is also a Patterson density of dimension $\sigma$ for
$(\Ga,F^*)$. Indeed, the invariance property \eqref{eq:equivdensity}
for $(\mu^*_{x})_{x\in\wt M}$ follows from the one for
$(\mu_{x})_{x\in\wt M}$ and the invariance of $\wt G$. And the
absolutely continuous property \eqref{eq:radonykodensity} for
$(\mu^*_{x})_{x\in\wt M}$ follows from the one for $(\mu_{x})_{x\in\wt
  M}$ and Equation \eqref{eq:cohomologuecocycle}.

\subsection{The Gibbs states of $(\Ga,F)$}
\label{subsec:GibbsSulivanmeasure}

Let $\sigma\in\RR$, let $(\mu^\iota_{x})_{x\in\wt M}$ be a Patterson
density of dimension $\sigma$ for $(\Ga,F\circ \iota)$, and let
$(\mu_{x})_{x\in\wt M}$ be a Patterson density of the same dimension
$\sigma$ for $(\Ga,F)$. Using the Hopf parametrisation with
respect to any base point $x_0$ of $\wt M$, we define the {\it Gibbs
  measure on $T^1\wt M$ associated with the pair of Patterson
  densities $\big((\mu^\iota_{x})_{x\in\wt M}, (\mu_{x})_{x\in\wt
    M}\big)$}\index{Gibbs measure!on $T^1\wt M$!associated with a pair
  of densities} (of dimension $\sigma$) as the measure $\wt m$ on
$T^1\wt M$ given by
$$
d\,\wt m(v)= \frac{d\mu^\iota_{x_0}(v_-)\,d\mu_{x_0}(v_+)\,dt}
{D_{F-\sigma,\,x_0}(v_-,v_+)^2} \;.
$$
Using Equation \eqref{eq:ecartvisuel}, the following expression of
this measure, again in the Hopf parametrisation $v\mapsto (v_-,v_+,t)$
with respect to the base point $x_0$, will be useful
\begin{equation}\label{eq:formulemesgibbs}
d\,\wt m(v)= 
e^{C_{F\circ \iota-\sigma,\,v_-}(x_0,\,\pi(v))+C_{F-\sigma,\,v_+}(x_0,\,\pi(v))}
\;d\mu^\iota_{x_0}(v_-)\,d\mu_{x_0}(v_+)\,dt\;.
\end{equation}

The measure $\wt m$ on $T^1\wt M$ is
independent of $x_0$ by the equations \eqref{eq:formulemesgibbs},
\eqref{eq:radonykodensity} and \eqref{eq:cocycleprop}, hence is
invariant under the action of $\Ga$ by the equations
\eqref{eq:equivdensity} and \eqref{eq:equivgap}. It is invariant under
the geodesic flow. Hence (see Subsection \ref{subsec:pushmeas}) it
defines a measure $m$ on $T^1 M=\Ga\bs T^1\wt M$ which is invariant
under the quotient geodesic flow. We call $m$ the {\it Gibbs measure
  on $T^1M$ associated with the pair of Patterson densities
  $\big((\mu^\iota_{x})_{x\in\wt M}, (\mu_{x})_{x\in\wt
    M}\big)$}.\index{Gibbs measure!on $T^1 M$!associated with a pair
  of densities} We will justify the terminology in Subsection
\ref{subsec:gibbsproperty}.

\medskip 
The (Borel positive) measure 
\begin{equation}\label{eq:geodesicurrent}
d\lambda(\xi,\eta)= \frac{d\mu^\iota_{x_0}(\xi)\,d\mu_{x_0}(\eta)}
{D_{F-\sigma,\,x_0}(\xi,\eta)^2}
\end{equation}
on $\partial^2_\infty\wt M$, which is locally finite and invariant
under the diagonal action of $\Ga$ on $\partial^2_\infty\wt M$, is a
{\it geodesic current}\index{geodesic current} for the action of $\Ga$
on the Gromov-hyperbolic proper metric space $\wt M$ in the sense of
Ruelle-Sullivan-Bonahon (see for instance \cite{Bonahon91} and
references therein).

\medskip When $\delta_{\Ga,\,F}< \infty$, we will denote by
$\gls{measuregibbs}$, and call the {\it Gibbs measure on $T^1\wt M$
  with potential $F$},\index{Gibbs measure!on $T^1\wt M$ with
  potential $F$} the Gibbs measure on $T^1\wt M$ associated with any
given pair of densities $\big((\mu_{F\circ\iota,\,x})_{x\in\wt M},
(\mu_{F,\,x})_{x\in\wt M}\big)$ (which have the same dimension
$\delta_{\Ga,\,F}$ by Equation \eqref{eq:pressuretimereversalinva}).
Its induced measure on $T^1 M$ will be denoted by
$\gls{measuregibbsbas}$ and called the {\it Gibbs measure on $T^1 M$
  with potential $F$}%
\index{Gibbs measure!on $T^1 M$ with potential $F$}. By the uniqueness
statement in Subsection \ref{subsec:uniqueness}, if $(\Ga,F)$ is of
divergence type (we will prove in Corollary \ref{coro:uniqGibstate}
that this is the case for instance if $m_F$ is finite), then $\wt m_F$
and $m_F$ are uniquely defined, up to a scalar multiple.

Since the supports of the Patterson densities
$(\mu_{F\circ\iota,\,x})_{x\in\wt M}$ and $(\mu_{F,\,x})_{x\in\wt M}$ are
equal to $\Lambda\Ga$, the support of $\wt m_F$ is the set
$\wt\Omega\Ga$ of tangent vectors to the geodesic lines in $\wt M$
whose endpoints both lie in $\Lambda\Ga$. Using the Hopf
parametrisation, this support corresponds to $\big((\partial_\infty^2
\wt M) \cap(\Lambda\Ga\times\Lambda\Ga)\big)\times \RR$. Hence the
support of $m_F$ is the (topological) non-wandering set $\Omega\Ga$ of
the geodesic flow in $T^1M$.

Since only the normalised potential $F-\delta_{\Ga,\,F}$ is involved
in the definition of $\wt m_F$ and $m_F$, the assumption that
$\delta_{\Ga,\,F}$ is positive is harmless, by Equation
\eqref{eq:invarpressuretranslat}. We will often make this assumption
in Chapter \ref{sec:growth} and Chapter \ref{sec:ergtheounistabfolia}.

\bigskip
\noindent {\bf Remark (1) }  
Since $d\,\wt m(v)= \frac{d\mu^\iota_{x_0}(v_-)\,d\mu_{x_0}(v_+)\,dt}
{D_{F\circ \iota-\sigma,\,x_0}(v_+,\,v_-)^2}$ by Equation
\eqref{eq:equivgap}, the measure $\iota_*\wt m$ is the Gibbs measure
on $T^1\wt M$ associated with the switched pair of Patterson densities
$\big((\mu_{x})_{x\in\wt M},(\mu^\iota_{x})_{x\in\wt M}\big)$ for
$(\Ga,F)$ and $(\Ga,F\circ \iota)$ respectively (and similarly on $T^1
M$).

\medskip
\noindent {\bf Remark (2) }  Another parametrisation of $T^1\wt M$
(used for instance in Subsection \ref{subsec:condmesgibbs})
also depending on the choice of a base point $x_0$ in $\wt M$, is the
map from $T^1\wt M$ to $\partial_\infty^2\wt M \times \RR$ sending $v$
to $(v_-, v_+,s)$ where $v_-$ and $v_+$ are as above and $s=
\beta_{v_-}(\pi(v),x_0)$ (we may also use the different time parameter
$s= \beta_{v_+}(x_0,\pi(v))$).

\medskip\noindent\begin{minipage}{6cm} ~~~ 
For every $(\eta,\xi)$ in $\partial_\infty^2 \wt M$, let
$p_{\eta,\,\xi}$ be the closest point to $x_0$ on the geodesic line
between $\eta$ and $\xi$.
\end{minipage}
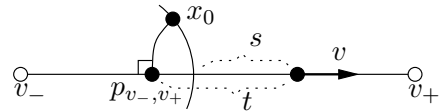
\begin{minipage}{8.9cm}
\begin{center}
\input{fig_autparahopf.pstex_t}
\end{center}
\end{minipage}

\medskip For every $v\in T^1\wt M$, with $(v_-, v_+,t)$ the original
Hopf parametrisation, since
$$
s-t=\beta_{v_-}(\pi(v),x_0)-\beta_{v_-}(\pi(v),p_{v_-,v_+})=
\beta_{v_-}(p_{v_-,v_+},x_0)
$$
depends only on $v_-$ and $v_+$, the measures $d\mu^\iota_{x_0}(v_-)
\,d\mu_{x_0}(v_+)\,dt$ and $d\mu^\iota_{x_0}(v_-)\, d\mu_{x_0}(v_+) \,
ds$ are equal. We hence may (and will in the proof of Proposition
\ref{prop:disintegrGibbs}) use the second one to define the Gibbs
measure associated with the above pair of densities.

\medskip
\noindent {\bf Remark (3) } It is sometimes useful (see for instance
the proof of Theorem \ref{theo:uniergodtransv}) to consider the
following more general construction. Let $x_0, x'_0\in \wt M$ and
$\sigma^\iota,\sigma\in\RR$ with possibly $\sigma^\iota\neq\sigma$.
Let $(\mu^\iota_{x})_{x\in\wt M}$ be a Patterson density of dimension
$\sigma^\iota$ for $(\Ga,F\circ \iota)$, and let $(\mu_{x})_{x\in\wt
  M}$ be a Patterson density of dimension $\sigma$ for
$(\Ga,F)$. Using the Hopf parametrisation with respect to $x_0$ (the
original one or the one defined in the previous Remark (2)), we define
the {\it Gibbs measure}%
\index{Gibbs measure!on $T^1\wt M$!associated with a pair of densities
  of different dimensions} associated with the pair of densities
$\big((\mu^\iota_{x})_{x\in\wt M}, (\mu_{x})_{x\in\wt M}\big)$ as the
measure $\wt m$ on $T^1\wt M$ given by
\begin{equation}\label{eq:formulemesgibbsgene}
d\,\wt m(v)= 
e^{C_{F\circ \iota-\sigma^\iota,\,v_-}(x'_0,\,\pi(v))+C_{F-\sigma,\,v_+}(x_0,\,\pi(v))}
\;d\mu^\iota_{x'_0}(v_-)\,d\mu_{x_0}(v_+)\,dt\;.
\end{equation}
It is again independent of $x_0, x'_0$ by the equations
\eqref{eq:radonykodensity} and \eqref{eq:cocycleprop}, hence is
invariant under the action of $\Ga$ by the equations
\eqref{eq:equivdensity} and \eqref{eq:equivcocprop}, and defines a
measure $m$ on $T^1M$, by Subsection \ref{subsec:pushmeas}.%
\index{Gibbs measure!on $T^1 M$!associated with a pair of densities of
  different dimensions} But it is no longer invariant under the
geodesic flow, as it now satisfies, for all $t\in\RR$ and $v\in T^1\wt
M$, using Equation \eqref{eq:cocyclealongray} to get the second equality,
\begin{align*}
d(\phi_t)_*\wt m(v)&= e^{C_{F\circ\iota-\sigma^\iota,\,v_-}
(\pi(v),\,\pi(\phi_{-t}v))+
C_{F-\sigma,\,v_+}(\pi(v),\,\pi(\phi_{-t}v))} \;d\,\wt m(v)\;.
\\ &=
e^{(\sigma^\iota-\sigma)t} \;d\,\wt m(v)\;.
\end{align*}
When $F=0$, these measures have been considered for instance by Burger
\cite{Burger90} and Knieper \cite{Knieper98} in particular cases, and
by Roblin \cite[page 12]{Roblin03} in general, and are sometimes
called {\it Burger-Roblin measures}\index{Burger-Roblin
  measures}\index{measure!of Burger-Roblin}, see \cite{OhShaCounting,
  Kim13}. When $F$ is non-constant, see for instance
\cite{Schapira03a}.

\medskip
\noindent {\bf Remark (4) } The Gibbs measure $m_F$ is not changed if we
replace $F$ by a cohomologous potential. Indeed, let $\wt F^* :T^1\wt
M\ra \RR$ be a H\"older-continuous $\Ga$-invariant map cohomologous to
$\wt F$ via the map $\wt G :T^1\wt M\ra \RR$. We have
$\delta_{\Ga,\,F^*}= \delta_{\Ga,\,F}$ by the remark at the end of
Subsection \ref{subsec:GibbsPoincareseries}. Choose the Patterson
densities of the remark at the end of Subsection
\ref{subsec:gibbspattersondens} to construct $\wt m_{F^*}$ (such a
choice is not important when $(\Ga,F)$ is of divergence type, see
Subsection \ref{subsec:uniqueness}). By Equation
\eqref{eq:formulemesgibbs}, by the remark at the end of Subsection
\ref{subsec:gibbspattersondens} and by Equation
\eqref{eq:cohomologuecocycle}, we then have
$$
\wt m_{F^*}=\wt m_F\;.
$$

\subsection{The Gibbs property of Gibbs states}
\label{subsec:gibbsproperty}

Let us indicate why our terminology of Gibbs measures (introduced in
Subsection \ref{subsec:GibbsSulivanmeasure}) is appropriate, by
explicitly making an analogy with the symbolic dynamics case and by
proving they satisfy the Gibbs property on the dynamical balls of the
geodesic flow (see Remark (3) in the beginning of Chapter
\ref{sec:ergtheounistabfolia} for other explanations).

Gibbs measures were first introduced in statistical mechanics,
and are naturally associated, via the thermodynamic formalism%
\index{thermodynamic formalism} (see for instance
\cite{Ruelle04,Keller98,Zinsmeister96}), to symbolic dynamics. Let us
recall the definitions of the pressure of a potential in a general
dynamical system, and its equilibrium states, adapting them to the
non-compact case, as in \cite{Sarig99} for the countable Markov shift
case.

If $\phi'=(\phi'_t)_{t\in\ZZ}$ or $\phi'=(\phi'_t)_{t\in\RR}$ is a
discrete-time or continuous-time one-parameter group of homeomorphisms
of a metric space $X'$ and if $F':X'\ra\RR$ is a continuous map, let us
define the {\it pressure $P(F')$ of the potential}%
\index{pressure!of a potential} $F'$ as the upper bound of
$h_{m'}(\phi')+\int_{T^1M}F'\,dm'$, on all $\phi'$-invariant
probability measures $m'$ on $X'$ such that the negative part
$\max\{0,-F'\}$ of $F'$ is $m'$-integrable, where
$h_{m'}(\phi')=h_{m'}(\phi'_1)$ is the (metric) entropy of the
one-parameter group of homeomorphisms $\phi'$ with respect to
$m'$. Define an {\it equilibrium state}\index{equilibrium state} of
$F'$ as any such $m'$ realising this upper bound. When $\phi'=\phi$ is
the geodesic flow on $X'=T^1M$, we will come back to these definitions
in Chapter \ref{sec:variaprincip}.

\medskip Let us now recall what is a Gibbs measure in symbolic
dynamics.  Let $A$ be a countable set, called an {\it
  alphabet}\index{alphabet} (see \cite{Sarig99,Sarig01a,Sarig01b} for
the infinite case), endowed with the discrete distance $d(a,a')=1$ if
$a\neq a'$, and $d(a,a')=0$ otherwise. Let $\Sigma$ be the product
topological space $A^\ZZ$, which is compact if $A$ is finite. Let us
endow it with the distance (inducing its topology)
$$
d((x_i)_{i\in\ZZ},(x'_i)_{i\in\ZZ})=\sum_{i\in\ZZ}
e^{-i^2}d(x_i,x'_i)\;.
$$ 
Let $\sigma:\Sigma\ra\Sigma$ be the (full) {\it shift
  map},\index{shift map} which is the homeomorphism of $\Sigma$
defined by $(x_i)_{i\in\ZZ}\mapsto (y_i)_{i\in\ZZ}$ where
$y_i=x_{i+1}$ for every $i\in\ZZ$. Note that if $\pi':\Sigma\ra A$ is
the map $(x_i)_{i\in\ZZ}\mapsto x_0$, then for all $x,x'\in\Sigma$, we
have
$$
d(x,x')=\sum_{i\in\ZZ} \;d(\pi'(\sigma^ix),\pi'(\sigma^ix'))\;e^{-i^2}
$$
(compare with Equation \eqref{eq:formdistT1wtM}).
For all $n,n'\in\NN$ and $a_{-n'},\dots, a_n$ in $A$, let
$$
\mathopen{[}a_{-n'},\dots, a_n\mathclose{]}=\{(x_i)_{i\in\ZZ}\in \Sigma\;:\; 
x_i=a_i\;\;{\rm for}\;-n'\leq i\leq n\}\;,
$$
which is called a {\it cylinder}\index{cylinder} in $\Sigma$. Recall
that given $(x_i)_{i\in\ZZ}\in\Sigma$, the set of cylinders
$\{\mathopen{[}x_{-n'},\dots, x_n\mathclose{]}\;:\; n,n'\in\NN\}$ is a
(natural) neighbourhood basis of $(x_i)_{i\in\ZZ}$ (even if we take
$n'=n$). Finally, let $F:\Sigma\ra \RR$ be a H\"older-continuous
map. For instance when $A=\{-1,+1\}$, the energy map $F:\Sigma\ra \RR$
defined by $(x_i)_{i\in\ZZ}\mapsto x_0(x_1+ x_{-1})$ describes the
one-dimensional Ising model\index{Ising model} (see for instance
\cite[page 96]{Keller98}).

The Gibbs measures are then shift-invariant measures, whose mass of a
cylinder is a weight defined by integrating the normalised potential
along the corresponding piece of orbit.  More precisely (compare for
instance with \cite[page 100]{Keller98}), a {\it Gibbs measure of
  potential $F$}\index{Gibbs measure!in symbolic dynamics} is a
$\sigma$-invariant measure $m'$ on $\Sigma$ such that there exists
$c'(F)\in\RR$ such that, for every finite subset $A'$ of $A$, there
exists $c\geq 1$ such that for all $x=(x_i)_{i\in\ZZ}\in\Sigma$ and
$n,n'\in\NN$ with $x_{-n'},x_n\in A'$, we have
\begin{equation}\label{eq:propGibbssymboldyn}
\frac{1}{c}\leq \frac{m'(\mathopen{[}x_{-n'},\dots, x_n\mathclose{]})}
{e^{\sum_{i=-n'}^n F(\sigma^i x)-(n+n'+1)c'(F)}}\leq c\;.
\end{equation}

In this symbolic framework, when the alphabet $A$ is finite, with
$\phi'= (\sigma^i)_{i\in\ZZ}$, $X'=\Sigma$ and $F'=F$, it is
well-known (see for instance \cite{Bowen75, Ruelle04, Keller98,
  Zinsmeister96} as well as \cite{HayRue92}) that a probability Gibbs
measure of $F$ exists (with constant $c'(F)$ equal to the pressure
of $F$), is unique, and is the unique equilibrium state of $F$ (see
\cite{Sarig99,Sarig01a,Sarig01b} for what remains true in the
countable alphabet case).

\medskip The analogy between the geodesic flow $\phi=
(\phi_t)_{t\in\RR}$ and the full shift $(\sigma^i)_{i\in\ZZ}$ proceeds
as follows.
For all $v\in T^1\wt M$ and $r>0$, $T,T'\geq 0$, let us define
$$
B(v;T,T',r)=\{w\in T^1\wt M\;:\;\; \sup_{t\in\mathopen{[}-T',T\mathclose{]}}\;
d(\pi(\phi_tv),\pi(\phi_tw))<r\}\;,
$$
called a {\it dynamical (or Bowen) ball}%
\index{dynamical!ball}\index{Bowen ball} around $v$. They satisfy the
following invariance properties: for all $s\in[-T',T]$ and $\ga\in\Ga$,
$$
\phi_s B(v;T,T',r)= B(\phi_s v;T-s,T'+s,r)\;\;\;{\rm and}\;\;\;
\ga B(v;T,T',r)= B(\ga v;T,T',r)\;.
$$
Note the following inclusion properties of the dynamical balls: if
$r\leq s$, $T\geq S$, $T'\geq S'$, then $B(v;T,T',r)$ is contained in
$B(v;S,S',s)$. 

\blemm\label{lem:containdynaboul} (1) For all $r'\geq r>0$, there
exists $T_{r,\,r'}\geq 0$ such that for all $v\in T^1\wt M$ and
$T,T'\geq 0$, the dynamical ball $B(v;T+T_{r,\,r'},T'+T_{r,\,r'},r')$
is contained in $B(v;T,T',r)$.

(2) If $B(v;T,T',r)$ meets $B(u;S,S',r)$ and $T\geq S$, $T'\geq S'$,
then $v$ belongs to the dynamical ball $B(u;S,S',2r)$.  
\elemm

\dem 
(1) This follows by the properties of long geodesic segments with
endpoints at bounded distance in a $\operatorname{CAT}(-1)$ space (see
for instance \cite{BriHae99}).

\medskip (2) Let $z\in B(v;T,T',r)\cap B(u;S,S',r)$. By the triangle
inequality and since $T\geq S$, $T'\geq S'$, we have
$$
\max_{t\in\mathopen{[}-S',S\mathclose{]}}\;
d(\pi(\phi_tv),\pi(\phi_tu))\leq \max_{t\in\mathopen{[}-T',T\mathclose{]}}\;
d(\pi(\phi_tv),\pi(\phi_tz))+\max_{t\in\mathopen{[}-S',S\mathclose{]}}\;
d(\pi(\phi_tz),\pi(\phi_tu))<2r\;.
$$
Hence $v\in B(u;S,S',2r)$.
\cqfd

\medskip It is easy to see that for every fixed $r>0$, the set
$\{B(v;T,T',r) \;:\; T,T'\geq 0\}$ is a (natural) neighbourhood basis
of $v\in T^1\wt M$ (even if we take $T'=T$), analogous to the set of
cylindrical neighbourhoods of a sequence in $\Sigma$. We will see in
Lemma \ref{lem:dynamballHopfparamet} below that, as it was mentioned
in the introduction, these (small) neighbourhoods of $v$ are defined
by (small) neighbourhoods of the two points at infinity $v_-,v_+$ of
the geodesic line defined by $v$. For every $v\in T^1M$, let us define
$B(v;T,T',r')$ as the image by the canonical projection $T^1\wt M\ra
T^1M=\Ga\bs T^1\wt M$ of $B(\wt v;T,T',r')$, for any preimage $\wt v$
of $v$ in $T^1\wt M$.
 
Adapting the definition of \cite{KatHas95} to the non-compact case, we
want to consider the measures on $T^1M$ invariant under the geodesic
flow, whose masses of appropriate dynamical balls are weights defined
by integrating the normalised potential along the corresponding pieces
of orbits, as follows. 

\bdefi \label{def:gibbsprop} A $\phi$-invariant measure $m'$ on $T^1M$
{\rm satisfies the Gibbs property}%
\index{Gibbs property}\index{measure!satisfying the Gibbs property} for
the potential $F$ with constant $c(F)\in\RR$ if for every compact
subset $K$ of $T^1 M$, there exist $r>0$ and $c_{K,\,r}\geq 1$ such
that for all large enough $T,T'\geq 0$, for every $v\in T^1 M$ with
$\phi_{-T'}v,\phi_{T}v\in K$, we have
$$
\frac{1}{c_{K,\,r}}\leq \frac{m'(B(v;T,T',r))}
{e^{\int_{-T'}^T (F(\phi_t v)-c(F))\,dt}}\leq c_{K,\,r}\;.
$$
\edefi

\noindent{\bf Remarks. } (1) The constant $c(F)$ is then uniquely
defined.

\medskip (2) Since 
$$
\frac{m'(B(v;T,T',r))} {e^{\int_{-T'}^T (F(\phi_t
    v)-c(F))\,dt}}=\frac{m'(B(\phi_{-T'}v;T+T',0,r))}
{e^{\int_{0}^{T+T'} (F(\phi_t (\phi_{-T'}v))-c(F))\,dt}}
$$ 
by invariance of the measures, and by a compactness argument, and
again by this invariance property, a $\phi$-invariant measure $m'$ on
$T^1M$ satisfies the Gibbs property for the potential $F$ with
constant $c(F)\in\RR$ if and only if for every compact subset $K$ of
$T^1 M$, there exist $r>0$ and $c_{K,\,r}\geq 1$ such that one of the
following three assertions holds

$\bullet$~ for all large enough $T\geq 0$, for every $v\in T^1 M$ with
$v,\phi_{T}v\in K$, we have
$$
\frac{1}{c_{K,\,r}}\leq \frac{m'(B(v;T,0,r))}
{e^{\int_{0}^T (F(\phi_t v)-c(F))\,dt}}\leq c_{K,\,r}\;;
$$

$\bullet$~ for all $T\geq 0$ and $v\in T^1 M$ with
$v,\phi_{T}v\in K$, we have
$$
\frac{1}{c_{K,\,r}}\leq \frac{m'(B(v;T,0,r))}
{e^{\int_{0}^T (F(\phi_t v)-c(F))\,dt}}\leq c_{K,\,r}\;;
$$

$\bullet$~ for all $T,T'\geq 0$ and $v\in T^1 M$ with
$\phi_{-T'}v,\phi_{T}v\in K$, we have
$$
\frac{1}{c_{K,\,r}}\leq \frac{m'(B(v;T,T',r))}
{e^{\int_{-T'}^T (F(\phi_t v)-c(F))\,dt}}\leq c_{K,\,r}\;.
$$ 

Furthermore, fixing $T_0\geq 0$, a $\phi$-invariant measure $m'$ on
$T^1M$ satisfies the Gibbs property for the potential $F$ with
constant $c(F)\in\RR$ if and only if for every compact subset $K$ of
$T^1 M$, there exist $r>0$ and $c_{K,\,r,\,T_0}\geq 1$ such that for
all large enough $T\geq 0$, for every $v\in T^1 M$ with
$\phi_{-T_0}v,\phi_{T}v\in K$, we have
$$
\frac{1}{c_{K,\,r,\,T_0}}\leq \frac{m'(B(v;T,T_0,r))}
{e^{\int_{-T_0}^T (F(\phi_t v)-c(F))\,dt}}\leq c_{K,\,r,\,T_0}\;;
$$

\medskip (3) Using Lemma \ref{lem:containdynaboul} and this
equivalence of definitions, it is easy to check that, up to changing
the constant $c_{K,\,r}$, the Gibbs property does not depend on the
constant $r>0$. 

Indeed let $r'\geq r>0$ and let $T_{r,\,r'}$ be as in Lemma
\ref{lem:containdynaboul}. Since $|F|$ is bounded, say by $\kappa\geq
0$, on the compact subset $\bigcup_{t\in[-T_{r,\,r'},T_{r,\,r'}]}
\phi_tK$ for every compact subset $K$ of $T^1 M$, we have, for all
$T,T'\geq 0$ and for all $v\in T^1 M$ with $\phi_{-T'}v,\phi_{T}v\in K$,
\begin{multline*}
e^{-2\,T_{r,\,r'}(\kappa+|c(F)|)}\;
\frac{m'(B(v;T+T_{r,\,r'},T'+T_{r,\,r'},r'))}
{e^{\int_{-T'-T_{r,\,r'}}^{T+T_{r,\,r'}} (F(\phi_t v)-c(F))\,dt}}\\
\leq \frac{m'(B(v;T,T',r))}
{e^{\int_{-T'}^T (F(\phi_t v)-c(F))\,dt}}\leq\frac{m'(B(v;T,T',r'))}
{e^{\int_{-T'}^T (F(\phi_t v)-c(F))\,dt}}\;.
\end{multline*}
Hence the Gibbs property for $r'$ implies the Gibbs property for
$r$. The converse implication is proved by noticing that, for all
$T,T'\geq 0$ large enough (and in particular at least $T_{r,r'}$) and
$v\in T^1 M$ with $\phi_{-T'}v,\phi_{T}v\in K$,
\begin{multline*}
\frac{m'(B(v;T,T',r))}
{e^{\int_{-T'}^T (F(\phi_t v)-c(F))\,dt}}\leq\frac{m'(B(v;T,T',r'))}
{e^{\int_{-T'}^T (F(\phi_t v)-c(F))\,dt}}\;\\ \leq e^{2\,T_{r,\,r'}(\kappa+|c(F)|)}\;
\frac{m'(B(v;T-T_{r,\,r'},T'-T_{r,\,r'},r))}
{e^{\int_{-T'+T_{r,\,r'}}^{T-T_{r,\,r'}} (F(\phi_t v)-c(F))\,dt}}\;.
\end{multline*}

(4) We will prove in Proposition \ref{prop:uniqGibbpro} that, under
some recurrence properties (as a very particular case when $M$ is
compact), there exists at most one $\phi$-invariant measure on $T^1M$
which satisfies the Gibbs property for the potential $F$.

\medskip
The following result shows that our terminology of Gibbs measures
is indeed appropriate.

\bprop \label{prop:gibbsgibbs} Assume that $\delta_{\Ga,\,F}<+\infty$.
Let $m$ be the Gibbs measure on $T^1M$ associated with a pair of
Patterson densities $\big((\mu^\iota_{x})_{x\in\wt M},
(\mu_{x})_{x\in\wt M} \big)$ of dimension $\delta_{\Ga,\,F}$ for
$(\Ga,F\circ \iota)$ and $(\Ga,F)$ respectively. Then $m$ satisfies
the Gibbs property for the potential $F$, with constant $c(F)=
\delta_{\Ga,\,F}$. 
\eprop

In Chapter \ref{sec:variaprincip}, we will prove that the critical
exponent $\delta_{\Ga,\,F}$ is equal to the pressure $P(\Ga,F)$, so
that $c(F)=P(\Ga,F)$ in accordance with the case of symbolic dynamics
over a finite alphabet.

\medskip
\dem We start with the following lemma, which describes the dynamical
balls in terms of points at infinity.

\blemm \label{lem:dynamballHopfparamet} For all $r>0$ and $T,T'\geq
0$, for every $v\in T^1\wt M$, using the Hopf parametrisation with
base point $x=\pi(v)$ and setting $x_T=\pi(\phi_Tv)$ and
$x_{-T'}=\pi(\phi_{-T'}v)$, we have
$$
B(v;T,T',r)\subset \OOO_xB(x_{-T'},2r)\times \OOO_xB(x_{T},2r)
\times\mathopen{]}-r, r\mathclose{[}\;.
$$ 
Conversely, for every $r>0$, there exists $T_r>0$ such that for all
$T,T'\geq T_r$ and $v\in T^1\wt M$, using the Hopf parametrisation
with base point $x=\pi(v)$ and taking $x_{-T'}$ and $x_T$ as above, we
have
$$
\OOO_xB(x_{-T'},r)\times \OOO_xB(x_{T},r)\times\mathopen{]}-1,
1\mathclose{[}\;\subset B(v;T,T',2r+2) \;.
$$
\elemm

\dem To prove the first claim, for every $w\in B(v;T,T',r)$, let $p$
be the closest point to $x$ on the geodesic line defined by $w$ and let
$y_T$ be the point at distance $T$ from $x$ on the geodesic ray
$\mathopen{[}x,w_+\mathclose{[}\;$.  

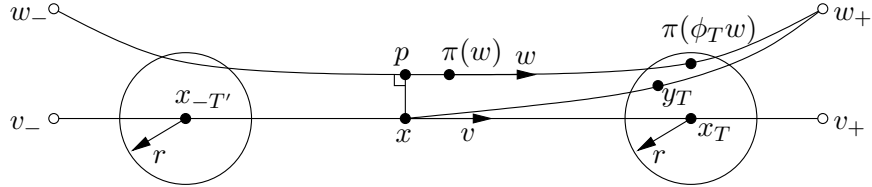
\begin{center}
\input{fig_bouledyn.pstex_t}
\end{center}

Since the closest point maps do not increase the distances, we have
$$
d(p,\pi(w))\leq d(\pi(v),\pi(w))< r\;.
$$
By the triangle inequality and by convexity, we have
$$
d(y_T,x_{T})\leq d(y_T,\pi(\phi_Tw))+d(\pi(\phi_Tw),\pi(\phi_Tv))\leq 
d(x,\pi(w))+r<2r\;,
$$
hence $w_+$ belongs to $\OOO_xB(x_{T},2r)$.  With a similar argument
for $w_-$, this proves the first claim.

\medskip To prove the second claim, let $\eta\in \OOO_xB(x_{-T'},r)$
and let $\xi\in\OOO_xB(x_{T},r)$, which is different from $\eta$ if
$T$ and $T'$ are larger than a constant depending only on $r$. Let $p$
be the closest point to $x$ on the geodesic line
$\mathopen{]}\eta,\xi\mathclose{[}\,$. Since the height of a vertex of
a geodesic triangle in a $\operatorname{CAT}(-1)$ space whose
comparison angle is close to $\pi$ is small, if $T$ and $T'$ are
larger than a constant depending only on $r$, we have $d(x,p)\leq
1$. Let $w\in T^1\wt M$ be such that $d(p,\pi(w))< 1$ and $w_-=\eta$,
$w_+=\xi$.

\begin{center}
\input{fig_bouledynbis.pstex_t}
\end{center}

Let us prove that $w\in B(v;T,T',2r+2)$, which yields the second
claim. By convexity, we only have to prove that $d(\pi(\phi_{-T'}w),
x_{-T'})<2r+2$ and $d(\pi(\phi_{T}w), x_{T}) <2r+2$. Since the argument
for the first inequality is similar, we only prove the second
inequality.

Let $z$ be the closest point to $x_T$ on the geodesic ray
$\mathopen{[}x,\xi\mathclose{[}\,$, which satisfies $d(z,x_T)\leq
r$. Note that $d(x,z)\geq d(x,x_T)-d(x_T,z)\geq T-r$ by the inverse
triangle inequality. Let $y_T$ be the point at distance $T$ from $x$
on the geodesic ray $\mathopen{[}x,w_+\mathclose{[}\,$, which
satisfies
$$
d(\pi(\phi_{T}w),y_T)\leq d(\pi(w),x)\leq d(\pi(w),p)+d(p,x)< 2
$$ 
by convexity. By the convexity of balls, we have $z\in
\mathopen{[}x,y_T\mathclose{]}$.  Then
\begin{align*}
d(\pi(\phi_{T}w), x_{T}) & \leq d(\pi(\phi_{T}w),y_T) +d(y_T,x_{T})
<2  +d(y_T,z)+d(z,x_{T})\\ & 
\leq 2+\big(d(y_T,x)-d(x,z)\big)+r\leq 2+r+T-(T-r)
=2r+2\;.\;\;\;\Box
\end{align*}

\medskip Now, let $\wt m$ be the Gibbs measure on $T^1\wt M$
associated with a pair of Patterson densities
$\big((\mu^\iota_{x})_{x\in\wt M}, (\mu_{x})_{x\in\wt M}\big)$ of
dimension $\delta=\delta_{\Ga,\,F}$ for $(\Ga,F\circ \iota)$ and
$(\Ga,F)$ respectively. 

Let $\wt K$ be a compact subset of $T^1\wt M$. Let $r>2$ and $C>0$ be
such that the conclusion of Mohsen's shadow lemma
\ref{lem:shadowlemma} holds true for $R=\frac{r-2}{2}$ and $R=2r$, for
the Patterson densities $(\mu^\iota_{x})_{x\in\wt M}$ and
$(\mu_{x})_{x\in\wt M}$, and for all $x,y$ in the compact subset
$\pi(\wt K)$ of $\wt M$. Let $T_*=T_{\frac{r-2}{2}}$ be given by Lemma
\ref{lem:dynamballHopfparamet}. Let $x= \pi(v)$, which stays in the
compact set $\N_{T_*}(\pi(\wt K))$. 

Let us prove that there exists $c'>0$ such that for all $T\geq 0$
large enough (in particular $T\geq T_*$) and for every $v\in T^1\wt M$
such that $\phi_{-T_*}v\in \wt K$ and $\phi_Tv\in \Ga \wt K$, we have
\begin{equation}\label{eq:GibbsGibbshaut}
\frac{1}{c'}\leq \frac{\wt m(B(v;T,T_*,r))}
{e^{\int_{-T_*}^T (\wt F(\phi_t v)-\delta)\,dt}}\leq c'\;.
\end{equation}
Note that $\phi_{-T_*}B(v;T,T_*,r)$ is contained in a compact subset of
$T^1\wt M$ (depending only on $\wt K$ and $r$). Hence the multiplicity
of the restriction to $\phi_{-T_*}B(v;T,T_*,r)$ of the canonical
projection $T^1\wt M\ra T^1 M$ is bounded, by the discreteness of
$\Ga$. Since $\wt m$ is $\phi$-invariant, Equation
\eqref{eq:GibbsGibbshaut} hence implies Proposition
\ref{prop:gibbsgibbs}, using Remark (2) following Definition
\ref{def:gibbsprop}.

By Lemma \ref{lem:holderconseq} (1) with $r_0=r$, there exists a
constant $c_1>0$ such that for all $v\in T^1\wt M$, $T\geq 0$ and
$w\in B(v;T,T_*,r)$, since $d(x,\pi(w))\leq r$, we have
\begin{align*}
&\max\big\{\,|\;C_{F\circ\iota-\delta,\,w_-}(x,\pi(w))\;|,\;
|\;C_{F-\delta,\,w_+}(x,\pi(w))\;|\,\big\}\\&\leq\; c'_1=
\;c_1\,e^{r} \;+\; r\,\max_{\pi^{-1}(\N_{r+T_*}(\pi(\wt K)))}\;|\wt F-\delta|\;.
\end{align*}
Note that $c'_1$ is finite since $\wt F$ is bounded on compact
subsets. We use the Hopf parametrisation with respect to the base
point $x$.  Let $x_T=\pi(\phi_Tv)$ and
$x_{-T_*}=\pi(\phi_{-T_*}v)$, which belong to $\Ga\pi(\wt K)$.

By Equation \eqref{eq:formulemesgibbs} and the second claim in Lemma
\ref{lem:dynamballHopfparamet}, we hence have
$$
\wt m(B(v;T,T_*,r))\geq 2\,e^{-2\,c'_1}\;
\mu^\iota_x(\OOO_xB(x_{-T_*},\frac{r-2}{2}))\; 
\mu_x(\OOO_xB(x_{T},\frac{r-2}{2})).
$$ 
Mohsen's shadow lemma \ref{lem:shadowlemma} and Equation
\eqref{eq:timereversal} imply that
\begin{align*}
\wt m(B(v;T,T_*,r))&\geq 2\,e^{-2\,c'_1}\,C^{-2}\;
e^{\int_x^{x_{-T_*}}(\wt F\circ\iota-\delta)}\;e^{\int_x^{x_{T}}(\wt F-\delta)}\\ & 
=2\,e^{-2\,c'_1}\,C^{-2}\;
e^{\int_{-T_*}^T(\wt F(\phi_tv)-\delta)\,dt}\;.
\end{align*}
This proves the lower bound in Equation \eqref{eq:GibbsGibbshaut}.  

The proof of the upper bound is similar. By Equation
\eqref{eq:formulemesgibbs} and the first claim in Lemma
\ref{lem:dynamballHopfparamet}, we have
$$
\wt m(B(v;T,T_*,r))\leq 2\,r\,e^{-2\,c'_1}\;
\mu^\iota_x(\OOO_xB(x_{-T_*},2\,r))\; 
\mu_x(\OOO_xB(x_{T},2\,r)).
$$ 
Mohsen's shadow lemma \ref{lem:shadowlemma} and Equation
\eqref{eq:timereversal} imply that
\begin{align*}
\wt m(B(v;T,T_*,r))&\leq 2\,r\,e^{2\,c'_1}\,C^2\;
e^{\int_x^{x_{-T_*}}(\wt F\circ\iota-\delta)}\;e^{\int_x^{x_{T}}(\wt F-\delta)}\\ & 
=2\,r\,e^{2\,c'_1}\,C^2\;
e^{\int_{-T_*}^T(\wt F(\phi_tv)-\delta)\,dt}\;.
\end{align*}
This proves the upper bound in Equation \eqref{eq:GibbsGibbshaut}, and
Proposition \ref{prop:gibbsgibbs} follows. \cqfd

\medskip 
\noindent \rem A more usual definition of the dynamical balls
is the following one. For every $T\geq 0$, consider the
$\Ga$-invariant distance $d'_{\phi,T}$ on $T^1\wt M$ defined by
$$
\forall\;v,w\in T^1\wt M,\;\;\;\; d'_{\phi,\,T}(v,w)=
\max_{t\in \mathopen{[}0,T\mathclose{]}} d'(\phi_tv,\phi_tw)\;,
$$
where $d'=d'_{T^1\wt M}$ is the distance on $T^1\wt M$ defined in
Equation \eqref{eq:defidprime}. 
For all $\epsilon>0$ and $T\geq 0$, let
$$
B'(v;\epsilon,T)=\{w\in T^1\wt  M\;:\; 
\max_{t\in \mathopen{[}0,T\mathclose{]}} d'(\phi_tv,\phi_tw)\;<\epsilon\}
$$
be the open ball with centre $v\in T^1\wt M$ and radius $\epsilon$ for
this distance $d'_{\phi,\,T}$. By the definition of $d'$, we have
$$
B'(v;\epsilon,T)=B(v;T,-1,\epsilon)\;.
$$
In particular, by the remarks (2) and (3) following Definition
\ref{def:gibbsprop}, a $\phi$-invariant measure $m'$ on $T^1M$
satisfies the Gibbs property for the potential $F$ with constant
$c(F)$ if and only if for every compact subset $K$ of $T^1 M$,
there exists (or equivalently for every) $\epsilon>0$, there exists
$c_{K,\,\epsilon}\geq 1$ such that for every (or equivalently for
every large enough) $T\geq 0$, for every $v\in T^1 M$ with
$v,\phi_{T}v\in K$, we have
$$
\frac{1}{c_{K,\,\epsilon}}\leq \frac{m'(B'(v;\epsilon,T))}
{e^{\int_{0}^T (F(\phi_t v)-c(F))\,dt}}\leq c_{K,\,\epsilon}\;.
$$

Note that if $d'=d'_{T^1\wt M}$ is replaced by any H\"older-equivalent
distance $d''$ (for instance the Riemannian distance $d_S$ induced by
the Sasaki metric, or the distance $d=d_{T^1\wt M}$ defined in
Equation \eqref {eq:formdistT1wtM}), and if $B''(v;\epsilon,T)$ are
the corresponding dynamical balls, then there exist $\epsilon_0, c>0$
and $\alpha\in\mathopen{]}0,1\mathclose{]}$ such that, for all
$\epsilon \in \mathopen{]}0,\epsilon_0\mathclose{]}$ and $T\geq 0$,
$$
B'(v;\frac{1}{c}\,\epsilon^{\frac{1}{\alpha}},T) \subset B''(v;\epsilon,T)
\subset B'(v;c\,\epsilon^\alpha,T)\;.
$$
Since $\wt M$ has pinched negative curvature and by the divergence
properties of the geodesic lines, for all $r>0$, there exist
$\epsilon_r, a_r,b_r,c_r,T_r>0$ such that for all
$\epsilon\in\mathopen{]}0,\epsilon_r\mathclose{]}$ and $T\geq T_r$,
\begin{multline*}
B(v;T-a_r\log\epsilon+c_r,-a_r\log\epsilon+c_r,r)\\ 
\subset B'(v;\epsilon,T) \subset \\ 
B(v;T-b_r\log\epsilon-c_r,-b_r\log\epsilon-c_r,r)\;.
\end{multline*}
When $\wt M$ has constant curvature $-1$, we may take $a_r=b_r=1$. But
in general variable curvature, as $\epsilon$ tends to $0$, the
relationship between the dynamical balls $B'(v;\epsilon,T)$ and our
ones $B'(v;T,T',r)$ becomes less precise.

\subsection{Conditional measures of Gibbs states 
on strong (un)stable leaves}
\label{subsec:condmesgibbs}

The aim of this subsection is to describe the conditional measures of
the Gibbs measures constructed in Subsection
\ref{subsec:GibbsSulivanmeasure} on the leaves of the strong unstable
and strong stable foliations. We refer for instance to \cite[Chap.~3,
\S 70]{DelMey75} for general information on the dis\-integration of
measures.

\medskip We start by describing the family of measures on the strong
unstable and strong stable leaves which will allow us to determine
the conditional measures of the Gibbs measures.

Let $\sigma\in\RR$. A Patterson density $(\mu_{x})_{x\in\wt
  M}$ of dimension $\sigma$ for $(\Ga,F)$ defines a family of measures
$(\mu_{W^{\rm su}(v)})_{v\in T^1\wt M}$ on the strong unstable leaves, in
the following way. Recall that for every $v\in T^1\wt M$, the map
$w\mapsto w_+$ from $W^{\rm su}(v)$ to $\partial_\infty \wt M-\{v_-\}$ is
a homeomorphism.  This allows us to define a measure $\mu_{W^{\rm su}(v)}$
on $W^{\rm su}(v)$ by
\begin{equation}\label{eq:defmeasstrongunstable}
d\mu_{W^{\rm su}(v)}(w)=
e^{C_{F-\sigma,\,w_{+}}(x_0,\,\pi(w))}\,d\mu_{x_0}(w_+)\;,
\end{equation}
where $x_0$ is any point of $\wt M$. This measure is nonzero, since
$\Ga$ being non-elementary, the support of $\mu_{x_0}$ is not reduced
to $\{v_-\}$. By the absolute continuity property
\eqref{eq:radonykodensity} of $(\mu_x)_{x\in\RR}$ and by the cocycle
property \eqref{eq:cocycleprop} of the Gibbs cocycle $C_{F-\sigma}$,
the measure $\mu_{W^{\rm su}(v)}$ does not depend on the choice of
$x_0$. By the invariance properties \eqref{eq:equivdensity} and
\eqref{eq:equivcocprop}, the family of measures
$(\mu_{W^{\rm su}(v)})_{v\in T^1\wt M}$ is $\Ga$-equivariant: for all
$v\in T^1\wt M$ and $\ga\in\Ga$,
\begin{equation}\label{eq:propunmesfortinstabl}
\ga_*\mu_{W^{\rm su}(v)}=\mu_{W^{\rm su}(\ga v)}\;.
\end{equation}
Note that $\mu_{W^{\rm su}(v)}=\mu_{W^{\rm su}(v')}$ if $v$ and $v'$ are in
the same strong unstable leaf. If the support of $\mu_x$ is
$\Lambda\Ga$, the support of $\mu_{W^{\rm su}(v)}$ is
\begin{equation}\label{eq:supportmesfortinstabl}
\operatorname{Supp}(\mu_{W^{\rm su}(v)})=
\{w\in W^{\rm su}(v)\;:\; w_+\in\Lambda\Ga\}\;.
\end{equation}
We will also consider the measure $\mu_{W^{\rm su}(v)}$ on $W^{\rm
  su}(v)$ as a measure on $T^1\wt M$ with support contained in $W^{\rm
  su}(v)$.  The map $v\mapsto\mu_{W^{\rm su}(v)}$ is continuous for
the weak-star topology on the space of measures on $T^1\wt M$. The
family of measures $(\mu_{W^{\rm su}(v)})_{v\in T^1\wt M}$ uniquely
determines the Patterson density $(\mu_{x})_{x\in\wt M}$, by Equation
\eqref{eq:defmeasstrongunstable}.  The geodesic flow
$(\phi_t)_{t\in\RR}$ preserves the measure class of these measures:
for all $t\in\RR$ and $w\in W^{\rm su}(v)$, we have
\begin{equation}\label{eq:propdemesfortinstabl}
\frac{d\,(\phi_t)_*\mu_{W^{\rm su}(v)}}
{d\,\mu_{W^{\rm su}(\phi_t v)}}(\phi_t w)=
e^{C_{F-\sigma,\,w_+}(\pi(\phi_tw),\,\pi(w))}= 
e^{\int_0^t(\wt F(\phi_sw)-\sigma)\,ds}\;,
\end{equation}
(which is equal to $e^{\int_{\pi(w)}^{\pi (\phi_tw)}(\wt F-\sigma)}$ if
$t\geq 0$ and to $e^{-\int_{\pi (\phi_tw)}^{\pi(w)}(\wt F-\sigma)}$
otherwise, by Equation \eqref{eq:cocyclealongray}).

\medskip\noindent\begin{minipage}{8.9cm} ~~~ For all $v,v'\in T^1\wt
  M$ such that $v_-\neq v'_-$, let $v_{v'}$ be the unique element of
  $W^{\rm su}(v)$ such that $(v_{v'})_+ = v'_-$. The map $\theta^{su} =
  \theta^{su}_{v',v}: W^{\rm su}(v) -\{v_{v'}\} \ra W^{\rm su}(v')-
  \{v'_{v}\}$ sending $w$ to the unique element $\theta^{su}(w)$ such
  that $\theta^{su}(w)_+=w_+$, is a homeomorphism. An easy computation
  using Equation \eqref{eq:defmeasstrongunstable} and Equation
  \eqref{eq:cocycleprop} shows that for every $w\in
  W^{\rm su}(v)-\{v_{v'}\}$, we have
\end{minipage}
\begin{minipage}{6cm}
\begin{center}
\input{fig_familmesure.pstex_t}
\end{center}
\end{minipage}
\medskip

\begin{equation}\label{eq:proptrmesfortinstabl}
  \frac{d\,(\theta^{su})_*\mu_{W^{\rm su}(v)}}{d\,\mu_{W^{\rm su}(v')}}
(\theta^{su}(w))=  e^{C_{F-\sigma,\,w_+}(\pi(\theta^{su}(w)),\,\pi(w))}\;.
\end{equation}

Conversely, a family of nonzero measures $(\mu_{W^{\rm su}(v)})_{v\in
  T^1\wt M}$ on $T^1\wt M$, constant on the strong unstable leaves,
with $\mu_{W^{\rm su}(v)}$ supported in $W^{\rm su}(v)$, satisfying the
properties \eqref{eq:propunmesfortinstabl},
\eqref{eq:propdemesfortinstabl} and \eqref{eq:proptrmesfortinstabl}
defines a unique Patterson density $(\mu_{x})_{x\in\wt M}$ of
dimension $\sigma$ for $(\Ga,F)$ satisfying Equation
\eqref{eq:defmeasstrongunstable}.

By the equivariance property \eqref{eq:propunmesfortinstabl}, the family
$(\mu_{W^{\rm su}(v)})_{v\in T^1\wt M}$ induces on $T^1 M$ a weak-star
continuous family of measures, which we will denote by
$(\mu_{W^{\rm su}(v)})_{v\in T^1M}$, constant on the strong unstable
leaves, with $\mu_{W^{\rm su}(v)}$ supported in $W^{\rm su}(v)$ for every
$v\in T^1M$.

\bigskip A Patterson density $(\mu_{x})_{x\in\wt M}$ of dimension
$\sigma$ for $(\Ga,F)$ also defines a family of measures
$(\mu_{W^{\rm u}(v)})_{v\in T^1\wt M}$ on the unstable leaves, in the
following way. 

Let $x_0\in \wt M$ and $v\in T^1\wt M$. We have a
homeomorphism from $W^{\rm su} (v) \times\RR$ to $W^{\rm u}(v)$,
defined by $(w,t)\mapsto w'=\phi_tw$, whose inverse is the map
$$
w'\mapsto\big(w=\phi_{\beta_{v_-}(\pi(v),\,\pi(w'))}w', 
t=\beta_{v_-}(\pi(w'),\pi(v))\big)\;.
$$ 
For future use, note that these homeomorphisms, as $v$ ranges over
$T^1\wt M$, satisfy the following equivariance property: For every
$\ga\in\Ga$, the following diagram commutes:
$$
\begin{array}{ccc} W^{\rm su} (v) \times\RR& \ra & W^{\rm u}(v)\\
^{\ga\times\id}\downarrow\;& & \;\downarrow{}^\ga\\
W^{\rm su} (\ga v) \times\RR& \ra & W^{\rm u}(\ga v)\;.
\end{array}
$$ 

Using this homeomorphism from $W^{\rm su} (v) \times\RR$ to $W^{\rm
  u}(v)$, as well as the homeomorphism from the unstable leaf $W^{\rm
  u}(v)$ to $(\partial_\infty \wt M-\{v_-\})\times \RR$ defined by
$$
w'\mapsto \big(w'_+,t=\beta_{v_-} (\pi(w'),\pi(v))\big)\;,
$$ 
we may define a nonzero (positive Borel) measure $\mu_{W^{\rm u}(v)}$ on
$W^{\rm u}(v)$, by
\begin{align}
d\mu_{W^{\rm u}(v)}(w')&=
e^{C_{F-\sigma,\,w_+}(\pi(w),\,\pi(\phi_t w))}\;
d\mu_{W^{\rm su}(v)}(w)dt\nonumber\\ &
=e^{C_{F-\sigma,\,w'_+}(x_0,\,\pi(w'))}\;d\mu_{x_0}(w'_+)dt\;.\label{eq:defimus}
\end{align}
In particular, the second formula is independent of the choice of
$x_0$.  Note that when $\wt F=0$, the first equality simplifies as
$$
d\mu_{W^{\rm u}(v)}(w')= e^{-\sigma\,t}\;d\mu_{W^{\rm su}(v)}(w)dt\;.
$$

By the invariance property \eqref{eq:equivcocprop} of the Gibbs cocycle
and by Equation \eqref{eq:propunmesfortinstabl}, the family of
measures $(\mu_{W^{\rm u}(v)})_{v\in T^1\wt M}$ is equivariant under
$\Ga$: for all $v\in T^1\wt M$ and $\ga\in\Ga$,
\begin{equation}\label{eq:propunmesinstabl}
\ga_*\mu_{W^{\rm u}(v)}=\mu_{W^{\rm u}(\ga v)}\;.
\end{equation}
Note that $\mu_{W^{\rm u}(v)}=\mu_{W^{\rm u}(v')}$ if $v$ and $v'$ are in the
same unstable leaf, by Equation \eqref{eq:defimus} and since the
parametrisation of $W^{\rm u}(v)$ by $(\partial_\infty\wt M-\{v_-\})
\times\RR$ is only changed, when passing from $v$ to $v'$, by adding a
constant to the parameter $t$, as the Lebesgue measure $dt$ is
invariant under translations.  If the support of $\mu_x$ is
$\Lambda\Ga$, the support of $\mu_{W^{\rm u}(v)}$ is
\begin{equation}\label{eq:supportmesfaibinstabl}
\operatorname{Supp}(\mu_{W^{\rm u}(v)})=
\{w\in W^{\rm u}(v)\;:\; w_+\in\Lambda\Ga\}\;.
\end{equation}
The map $v\mapsto\mu_{W^{\rm u}(v)}$ is continuous for the weak-star
topology on the space of measures on $T^1\wt M$. The family of
measures $(\mu_{W^{\rm u}(v)})_{v\in T^1\wt M}$ uniquely determines
the family of measures $(\mu_{W^{\rm su}(v)})_{v\in T^1\wt M}$, hence
the Patterson density $(\mu_{x})_{x\in\wt M}$.  The geodesic flow
$(\phi_t)_{t\in\RR}$ preserves the measure class of these measures:
for all $t\in\RR$ and $w\in W^{\rm u}(v)$, we have
\begin{equation}\label{eq:propdemesinstabl}
\frac{d\,(\phi_t)_*\mu_{W^{\rm u}(v)}}{d\,\mu_{W^{\rm u}(v)}}(\phi_t w)=
e^{\int_0^t(\wt F(\phi_sw)-\sigma)\,ds}\;.
\end{equation}

\medskip\noindent\begin{minipage}{8.9cm} ~~~ For all $v,v'\in T^1\wt
  M$, the map $\theta^u = \theta^u_{v',v}: W^{\rm u}(v) -\{w\in W^{\rm
    u}(v)\;:\; w_+= v'_-\} \ra W^{\rm u}(v')- \{w'\in W^{\rm u} (v')
  \;:\;w'_+= v_-\}$ sending $w$ to the unique element $w'=\theta^u(w)$
  in $W^{\rm ss}(w)$ such that $w'_-=v'_-$ is a homeomorphism. An easy
  computation using Equation \eqref{eq:proptrmesfortinstabl} shows
  that for every $w\in W^{\rm u}(v)$ with $w_+\neq v'_-$, we have
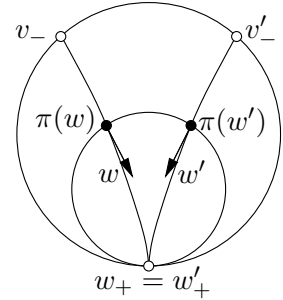
\begin{equation}\label{eq:proptrmesinstabl}
  \frac{d\,(\theta^u)_*\mu_{W^{\rm u}(v)}}{d\,\mu_{W^{\rm u}(v')}}(\theta^u(w))=
  e^{C_{F-\sigma,\,w_+}(\pi(\theta^u(w)),\,\pi(w))}\;.
\end{equation}
\end{minipage}
\begin{minipage}{6cm}
\begin{center}
\input{fig_changtransv.pstex_t}
\end{center}
\end{minipage}

\bigskip Similarly, given $\sigma^\iota\in\RR$ and a Patterson density
$(\mu^\iota_{x})_{x\in\wt M}$ of dimension $\sigma^\iota$ for
$(\Ga,F\circ\iota)$, we define weak-star continuous families
$(\mu^\iota_{W^{\rm ss}(v)})_{v\in T^1\wt M}$ and
$(\mu^\iota_{W^{\rm s}(v)})_{v\in T^1\wt M}$ of nonzero measures on
$T^1\wt M$, constant on respectively the strong stable leaves and
stable leaves, equivariant under $\Ga$, as follows:

$\bullet$~ using the homeomorphism from $W^{\rm ss}(v)$ to
$\partial_\infty\wt M-\{v_+\}$ defined by $w\mapsto w_-$, we set
\begin{equation}\label{eq:defmeasstrongstable}
d\mu^\iota_{W^{\rm ss}(v)}(w)=
e^{C_{F\circ \iota -\sigma^\iota,\,w_{-}}(x_0,\,\pi(w))}\,d\mu^\iota_{x_0}(w_-)\;,
\end{equation}
which is independent of $x_0\in\wt M$, 

$\bullet$~ using the homeomorphism from $W^{\rm ss} (v) \times\RR$ to
$W^{\rm s}(v)$, defined by $(w,t)\mapsto w'=\phi_tw$, and the
homeomorphism from $W^{\rm s}(v)$ to $(\partial_\infty \wt M-
\{v_+\})\times \RR$ defined by $w'\mapsto \big(w'_-,t=\beta_{v_+}
(\pi(v), \pi(w'))\big)$, we set
\begin{align}
d\mu^\iota_{W^{\rm s}(v)}(w')&=
e^{C_{F\circ \iota-\sigma^\iota,\,w_-}(\pi(w),\,\pi(\phi_t w))}
\;d\mu^\iota_{W^{\rm ss}(v)}(w)dt\label{eq:defimspourabscont}\\ &
=e^{C_{F\circ \iota-\sigma^\iota,\,w'_-}(x_0,\,\pi(w'))}
\;d\mu^\iota_{x_0}(w'_-)dt\;.\label{eq:defims}
\end{align}
As above, for all $v\in T^1\wt M$ and $\ga\in\Ga$, we have
\begin{equation}\label{eq:propunmesstabl}
\ga_*\mu^\iota_{W^{\rm ss}(v)}=\mu^\iota_{\ga W^{\rm ss}(v)}
\;\;\;{\rm and}\;\;\;
\ga_*\mu^\iota_{W^{\rm s}(v)}=\mu^\iota_{\ga W^{\rm s}(v)}\;.
\end{equation}
As above, if the support of $\mu^\iota_x$ is $\Lambda\Ga$, then the
support of $\mu^\iota_{W^{\rm ss}(v)}$ is $\{w\in W^{\rm ss}(v)\;: \;w_-\in
\Lambda\Ga \}$, and, for every $t\in\RR$,
$$
\forall\;w\in W^{\rm ss}(v),\;\;\; \frac{d\,(\phi_t)_*
\mu^\iota_{W^{\rm ss}(v)}}{d\,\mu^\iota_{W^{\rm ss}(\phi_t v)}}
(\phi_t w)=e^{-\int_0^t(\wt F(\phi_sw)-\sigma^\iota)\,ds}\;,
$$
$$
\forall\;w\in W^{\rm s}(v),\;\;\; \frac{d\,(\phi_t)_*
\mu^\iota_{W^{\rm s}(v)}}{d\,\mu^\iota_{W^{\rm s}(v)}}
(\phi_t w)=e^{-\int_0^t(\wt F(\phi_sw)-\sigma^\iota)\,ds}\;.
$$
For all $v,v'\in T^1\wt M$, the map $\theta^s : \{w\in W^{\rm s}(v)\;:\;
w_-\neq v'_+\} \ra \{w'\in W^{\rm s}(v')\;:\;w'_-\neq v_+\}$ sending $w$
to the unique element $w'$ in $W^{\rm su}(w)\cap W^{\rm s}(v')$, is a
homeomorphism. For every $w$ in the domain of $\theta^s$, we
have
\begin{equation}\label{eq:proptrmesstabl}
  \frac{d\,(\theta^s)_*\mu^\iota_{W^{\rm s}(v)}}{d\,\mu^\iota_{W^{\rm s}(v')}}
(\theta^s(w))=  e^{C_{F\circ\iota-\sigma^\iota,\,w_-}(\pi(\theta^s(w)),\,\pi(w))}\;.
\end{equation}

\medskip Now, before determining the conditional measures of the Gibbs
measures, we describe the local product structures of $T^1\wt M$
defined by the strong unstable and strong stable foliations, giving the 
local fibrations allowing to disintegrate the Gibbs measures.

Fix $w\in T^1\wt M$, and let 
$$
U_w=T^1\wt M-(W^{\rm s}(-w)\cup
W^{\rm u}(-w))= \{v\in T^1\wt M\;:v_+\neq w_-, v_-\neq w_+\}
$$ 
be the open neighbourhood of $w$ in $T^1\wt M$ consisting of the unit
tangent vectors that belong neither to the stable nor to the unstable
leaf of $-w$.  Note that $U_w$ is dense in $T^1\wt M$. For every $v\in
U_w$,

$\bullet$~ let $v_{su}$ be the intersection point of $W^{\rm su}(w)$
and $W^{\rm s}(v)$, which is the unique point $v_{su}$ of $W^{\rm su}
(w)$ such that $(v_{su})_+=v_+$,

$\bullet$~ let $v_{ss}$ be the intersection point of $W^{\rm ss}(w)$ and
$W^{\rm u}(v)$, which is the unique point $v_{ss}$ of $W^{\rm ss}(w)$ such
that $(v_{ss})_-=v_-$,

$\bullet$~ and let   $t_{su}\in\RR$ be the unique real number
such that $\phi_{t_{su}}(v_{su})\in W^{\rm ss}(v)$. 

\noindent Then the map $v\mapsto (v_{ss},v_{su},t_{su})$ from $U_w$ to
$W^{\rm ss}(w)\times W^{\rm su}(w)\times \RR$ is a homeomorphism. Note that,
when $v_-$ and $v_+$ are fixed, so is $(v_{su})_+$, and the difference
of the time parameters $t-t_{su}$ between the Hopf parametrisation and
this one is constant.

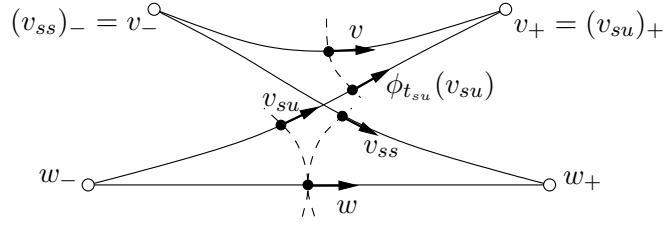
\begin{center}
\input{fig_gibbstabinstab.pstex_t}
\end{center}

The image by this homeomorphism of the Gibbs measures is absolutely
continuous with respect to product measures. Indeed, let $\sigma
\in\RR$ and let $\wt m$ be the Gibbs measure on $T^1\wt M$ associated
with a pair of Patterson densities $(\mu^\iota_{x}) _{x\in\wt M}$ and
$(\mu_{x})_{x\in\wt M}$ of dimension $\sigma$ for $(\Ga,F\circ \iota)$
and $(\Ga,F)$ respectively. Fix $x\in \wt M$. By definition of the
measures on the strong unstable and stable leaf, we have
$$
d\mu_{W^{\rm su}(w)}(v_{su})=
e^{C_{F-\sigma,\,(v_{su})_+}(x,\,\pi(v_{su}))}\;d\mu_x((v_{su})_+)
$$
and
$$
d\mu^\iota_{W^{\rm ss}(w)}(v_{ss})=
e^{C_{F\circ\iota-\sigma,\,(v_{ss})_-}(x,\,\pi(v_{ss}))}
\;d\mu^\iota_x((v_{ss})_-)\;.
$$
By Equation \eqref{eq:formulemesgibbs} and the cocycle property, for
all $v\in U_w$, we hence have
\begin{align}\label{eq:mesGibbmesstabinstab}
d\wt m(v)&
=e^{C_{F\circ\iota-\sigma,\,v_-}(\pi(v_{ss}),\,\pi(v))+
C_{F-\sigma,\,v_+}(\pi(v_{su}),\,\pi(v))}\;
d\mu^\iota_{W^{\rm ss}(w)}(v_{ss})\,d\mu_{W^{\rm su}(w)}(v_{su})\,dt_{su}\;.
\end{align}
This quasi-product property, either in the form of Equation
\eqref{eq:formulemesgibbs} or of Equation
\eqref{eq:mesGibbmesstabinstab}, is crucial in our book.  In
particular, if a measurable subset $A$ of $W^{\rm su}(w)$ has measure
$0$ for $\mu_{W^{\rm su}(w)}$, then the measurable set $\bigcup_{v\in
  A} W^{\rm s}(v)$ has measure $0$ for $\wt m$.

\medskip A slightly different way of understanding this property is as
follows. Let us fix $v\in T^1\wt M$. Let us define $U_{v_-}=\{w\in
T^1\wt M\;:\;w_+\neq v_-\}$, which is an open and dense subset of
$T^1\wt M$, invariant under the geodesic flow. If $\mu_x(\{v_-\})=0$
(this does not depend on $x\in \wt M$), for instance if the measure
$\mu_x$ has no atom, then $U_{v_-}$ has full measure with respect to
$\wt m$.  We have $U_{\ga v_-}=\ga U_{v_-}$ for every isometry $\ga$
of $\wt M$, and in particular $U_{v_-}$ is invariant under the
isometries of $\wt M$ fixing $v_-$.

The map $\psi_{W^{\rm u}(v)}$ from $U_{v_-}$ to $W^{\rm u}(v)$,
sending $w\in U_{v_-}$ to the unique element $w'$ in $W^{\rm ss}
(w)\cap W^{\rm u}(v)$, is a continuous fibration over the unstable
leaf $W^{\rm u}(v)$, whose fibre over $w'\in W^{\rm u}(v)$ is
precisely the strong stable leaf $W^{\rm ss}(w')$ of $w'$ (see the
left hand picture below). This map depends only on the unstable leaf
of $v$. For every isometry $\ga$ of $\wt M$, we have $\psi_{W^{\rm u}
  (\ga v)} \circ\ga=\ga\circ \psi_{W^{\rm u}(v)}$. For every
$t\in\RR$, we have $\psi_{W^{\rm u}(v)}\circ \phi_t=\phi_t\circ
\psi_{W^{\rm u}(v)}$: the fibration $\psi_{W^{\rm u}(v)}$ commutes
with the geodesic flow.

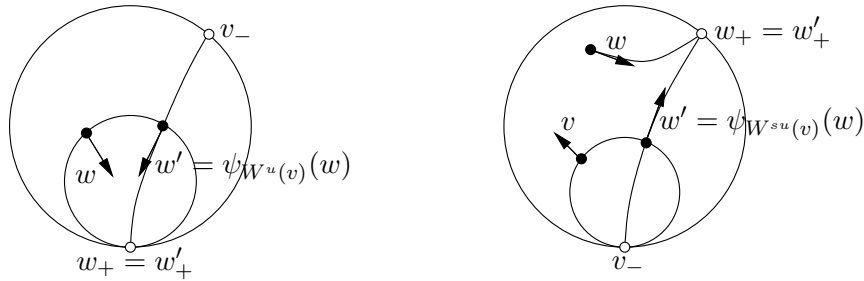
\begin{center}
\input{fig_fibstabunstab.pstex_t}
\end{center}

Similarly, for every $v\in T^1\wt M$, the map $\psi_{W^{\rm su}(v)}$ from
$U_{v_-}$ to $W^{\rm su}(v)$, sending $w\in U_{v_-}$ to the unique element
$w'$ in $W^{\rm s}(w)\cap W^{\rm su}(v)$, is a continuous fibration over the
strong unstable leaf $W^{\rm su}(v)$, whose fibre over $w'\in W^{\rm su}(v)$
is precisely the stable leaf $W^{\rm s}(w')$ of $w'$ (see the right hand
picture above). For every isometry $\ga$ of $\wt M$, we have
$\psi_{W^{\rm su}(\ga v)} \circ\ga=\ga\circ \psi_{W^{\rm su}(v)}$. For every
$t\in\RR$, we have $\psi_{W^{\rm su}(v)}\circ \phi_t= \psi_{W^{\rm su}(v)}$:
the fibration $\psi_{W^{\rm su}(v)}$ is invariant under the geodesic flow.

\medskip The aforementioned disintegration result, with respect to the
strong unstable foliation, of the Gibbs measures, in the general
setting of Remark (3) of Subsection \ref{subsec:GibbsSulivanmeasure},
is the following one.

\bprop \label{prop:disintegrGibbs} Let $\wt m$ be the Gibbs measure on
$T^1\wt M$ associated with a pair of Patterson densities
$(\mu^\iota_{x}) _{x\in\wt M}$ and $(\mu_{x})_{x\in\wt M}$ of
dimensions $\sigma^\iota$ and $\sigma$ for $(\Ga,F\circ \iota)$ and
$(\Ga,F)$, respectively. Let $v\in T^1\wt M$.

(1) The restriction to $U_{v_-}$ of the measure $\wt m$ disintegrates
by the fibration $\psi_{W^{\rm u}(v)}$ over the measure $\mu_{W^{\rm u} 
  (v)}$, with conditional measure on the fibre $W^{\rm ss}(w')$ of
$w'$ the measure $e^{C_{F,\,w'_+}(\pi(w'),\, \pi(w))}
\;d\mu^\iota_{W^{\rm ss}(w')}(w)$: for every $w\in T^1\wt M$ such that
$w_+\neq v_-$, we have
$$
d\wt m(w)=\int_{w'\in W^{\rm u}(v)}\;e^{C_{F,\,w'_+}(\pi(w'),\,
  \pi(w))} \;d\mu^\iota_{W^{\rm ss}(w')}(w)\;d\mu_{W^{\rm u}(v)}(w')\;.
$$

(2) The restriction to $U_{v_-}$ of the measure $\wt m$ disintegrates
by the fibration $\psi_{W^{\rm su}(v)}$ over the measure
$\mu_{W^{\rm su}(v)}$, with conditional measure on the fibre $W^{\rm s}(w')$
of $w'$ the measure $e^{C_{F,\,w'_+}(\pi(w'),\, \pi(w))}
\;d\mu^\iota_{W^{\rm s}(w')} (w)$: for every $w\in T^1\wt M$ such that
$w_+\neq v_-$, we have
$$
d\wt m(w)=\int_{w'\in W^{\rm su}(v)}\;e^{C_{F,\,w'_+}(\pi(w'),\,
  \pi(w))} \;d\mu^\iota_{W^{\rm s}(w')}(w)\;d\mu_{W^{\rm su}(v)}(w')\;.
$$
\eprop

We leave to the reader the analogous statements obtained by exchanging
the stable and unstable foliations. When $F=0$ and
$(\mu^\iota_{x})_{x\in\wt M}=(\mu_{x})_{x\in\wt M}$, we recover the
well-known fact, due to Margulis when $M$ is compact, that
$\mu_{W^{\rm su}(v)}$ and $\mu_{W^{\rm ss}(v)}$ are the conditional measures
of $\wt m$ along the leaves of the foliations $\wt \W^{\rm su}$ and
$\wt\W^{\rm ss}$, respectively.

\medskip \dem Fix $v\in T^1\wt M$. For every continuous map $\varphi
\in\C_c(U_{v_-};\RR)$ with compact support in the domain $U_{v_-}$ of
$\psi_{W^{\rm u}(v)}$, let
$$
I_\varphi=\int_{w\in T^1\wt M}\;\varphi(w)\;d\wt m(w)= 
\int_{w\in U_{v_-}}\varphi(w)\;d\wt m(w)\;.
$$ 
Using the Hopf parametrisation $w\mapsto \big(w_-,w_+,s= \beta_{w_+}
(x_0,\pi(w))\big)$ as in Remark (2) of Subsection
\ref{subsec:GibbsSulivanmeasure}, and Equation
\eqref{eq:formulemesgibbsgene} with $x'_0=x_0$ fixed in $\wt M$, we
have
\begin{align*}
I_\varphi=
\int_{w_+\in\,\partial_\infty\wt M-\{v_-\}}
\int_{w_-\in\,\partial_\infty\wt M-\{w_+\}}\int_{s\in\RR}\;&\varphi(w)\;
e^{C_{F\circ \iota-\sigma^\iota,\,w_-}(x_0,\,\pi(w))+C_{F-\sigma,\,w_+}(x_0,\,\pi(w))}\\&
d\mu^\iota_{x_0}(w_{-})\, d\mu_{x_0}(w_{+})\,ds\;.
\end{align*}
For every $w\in U_{v_-}$, let $w'=\psi_{W^{\rm u}(v)}(w)$ be the
unique point in $W^{\rm ss}(w)\cap W^{\rm u}(v)$. By the cocycle
property of the Gibbs cocycle, and since $\pi(w)$ and $\pi(w')$ are in
the same horosphere centred at $w_+=w'_+$, we have
\begin{align*}
C_{F-\sigma,\,w_+}(x_0,\pi(w))-C_{F-\sigma,\,w'_+}(x_0,\pi(w')) & =
C_{F-\sigma,\,w'_+}(\pi(w'),\pi(w))\\ & = C_{F,\,w'_+}(\pi(w'),\pi(w))\;.
\end{align*}

\smallskip\noindent
\begin{minipage}{8cm}~~~ When $w_+=w'_+$ is fixed, the time parameter
  $s=\beta_{w_+}(x_0,\pi(w))=\beta_{w_+}(x_0,\pi(w'))$ differs from
  the time parameter $t=\beta_{v_-} (\pi(w'),\pi(v))$ by a constant
  (equal to $\beta_{v_-}(\pi(v),x'_0)$ where $x'_0$ is the
  intersection point of the geodesic line $\mathopen{]}v_-,
  w_+\mathclose{[}$ and the horosphere centred at $w_+$ through
  $x_0$, by an easy computation), hence $ds=dt$.
\end{minipage}
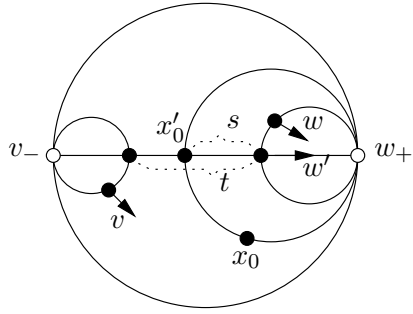
\begin{minipage}{6.9cm}
\begin{center}
\input{fig_stdesintgib.pstex_t}
\end{center}
\end{minipage}

\medskip 
By Equation \eqref{eq:defmeasstrongstable} and Equation
\eqref{eq:defimus}, we therefore have
$$
I_\varphi=\int_{w'\in\,W^{\rm u}(v)}
\int_{w\in\,W^{\rm ss}(w')}
\varphi(w)\;e^{C_{F,\,w'_+}(\pi(w'),\,\pi(w))}\;
d\mu^\iota_{W^{\rm ss}(w')}(w)\;d\mu_{W^{\rm u}(v)}(w')\;.
$$
This is the first claim of Proposition \ref{prop:disintegrGibbs}.
The second one is proved similarly.
\cqfd

\section{Critical exponent and Gurevich pressure} 
\label{sec:critgur}

Let $(\wt M,\Ga,\wt F)$ be as in the beginning of Chapter
\ref{sec:negacurvnot}: $\wt M$ is a complete simply connected
Riemannian manifold, with dimension at least $2$ and pinched sectional
curvature at most $-1$; $\Ga$ is a non-elementary discrete group of
isometries of $\wt M$; and $\wt F :T^1\wt M\ra \RR$ is a
H\"older-continuous $\Ga$-invariant map. The aim of this chapter is to
prove results on the critical exponent $\delta_{\Ga,\,F}$ which are
valid under these general hypotheses (in particular, $\Ga$ is only
assumed to be non-elementary).

\subsection{Counting orbit points and periodic geodesics} 
\label{subsec:countfunct} 

One way to study the repartition of an orbit of $\Ga$ in $\wt M$ is to
count the number of its elements in relatively compact subsets of $\wt
M$. We will count them with weights given by the potential $\wt F$.
We refer to Subsection \ref{subsec:geombus} for the definition of the
cone $\C_zB$ with vertex $z\in \wt M$ on $B\subset \partial_\infty
\wt M$. 

For all $s\geq 0$ and $c>0$, for all $x,y\in\wt M$ and for all
open subsets $U$ and $V$ of $\partial_\infty\wt M$, let
$$
\gls{bisectorbcountfunc}(t)= 
\sum_{\ga\in\Ga\;:\;d(x,\ga y)\leq t\,,\; 
\ga y\in \C_xU\,,\; \ga^{-1} x\in \C_yV} 
\;e^{\int_x^{\ga y} \wt F}\;,
$$
and
$$
\gls{annbisectorbcountfunc}(t)= 
\sum_{\ga\in\Ga\;:\;t-c<d(x,\ga y)\leq t\,,\; 
\ga y\in \C_xU\,,\; \ga^{-1} x\in \C_yV} 
\;e^{\int_x^{\ga y} \wt F}\;,
$$
The map $s\mapsto G_{\Ga,\,F,\,x,\,y,\,U,\,V}(s)$ will be called the
{\it bisectorial orbital counting function}%
\index{counting function!bisectorial orbital}\index{bisectorial
  orbital counting function} of $(\Ga, F,U,V)$, and $s\mapsto
G_{\Ga,\,F,\,x,\,y,\,U,\,V,\,c}(s)$ the {\it annular bisectorial
  orbital counting function}%
\index{counting function!annular!bisectorial orbital}%
\index{annular!bisectorial orbital counting function} of $(\Ga,
F,U,V)$. When $V=\partial_\infty \wt M$, we denote them by $s\mapsto
\gls{sectorbcountfunc}(s)$ and $s\mapsto \gls{annsectorbcountfunc}(s)$
and call them the {\it sectorial orbital counting function}%
\index{counting function!sectorial orbital}\index{sectorial orbital
  counting function} and {\it annular sectorial orbital counting
  function}%
\index{counting function!annular!sectorial orbital}%
\index{annular!sectorial orbital counting function} of
$(\Ga,F,U)$. When $U=V=\partial_\infty \wt M$, we denote them by
$s\mapsto \gls{orbcountfunc}(s)$ and $s\mapsto
\gls{annorbcountfunc}(s)$, and call them the {\it orbital counting
  function}%
\index{counting function!orbital}\index{orbital counting function} and
{\it annular orbital counting function}%
\index{counting function!annular!orbital}%
\index{annular!orbital counting function} of $(\Ga,F)$.  When $F=0$,
we recover the usual orbital counting functions of $\Ga$, see for
instance \cite{Margulis04,Roblin03} and \cite{Babillot02a} as well as
its references.

\medskip Other interesting counting functions are the ones counting
periodic orbits of the geode\-sic flow on $T^1M$, again with weights
given by the potential.

For every periodic (not necessarily primitive) orbit $g$ of length
$\ell(g)$ of the geodesic flow $(\phi_t)_{t\in\RR}$ on $T^1M$, let
$\gls{lebesguemeasure}$ be the {\it Lebesgue measure along $g$},%
\index{Lebesgue's measure along a periodic orbit}%
\index{measure!of Lebesgue along a periodic orbit} that is, the
measure on $T^1M$ with support $g$, such that, for every continuous
map $f:T^1M\ra \RR$ and for any $v\in g$,
$$
\L_g(f)=\int_{0}^{\ell(g)} f(\phi_tv)\;dt\;.
$$ 
The {\it period of $g$ for the potential $F$}\index{period!of a
  periodic orbit for $F$} is (for any $v\in g$)
$$
\int_gF=\L_g(F)=\int_{0}^{\ell(g)} F(\phi_tv)\;dt\;.
$$
One needs to be careful with the terminology and not confuse the
{length} $\ell(g)=\L_g(1)$ of $g$ and its {period} $\int_gF=\L_g(F)$
(which depends on $F$). 

For every $s\geq 0$, we denote by $\gls{setperiod}=\Per_{T^1M}(s)$ the
set of periodic orbits of the geodesic flow in $T^1M$ with length at
most $s$; we allow non-primitive ones, that is, periodic orbits winding
 $n$ times around a primitive one, in which case its length and
period are $n$ times the length and period of the primitive one, for
any $n\in\NN$.
We denote by $\gls{setperiodprim}=\Per'_{T^1M}(s)$ the subset of
primitive ones.

For every relatively compact open subset $W$ of $T^1M$, let us define the
{\it period counting function}%
\index{period!counting series}\index{counting function!period} of
$(\Ga,F,W)$ by
$$
s\mapsto \gls{periodcountingfucnt}(s)=
\sum_{g\in \Per(s),\;g\cap W\neq \emptyset}\;e^{\int_g F}\;,
$$
and, for every $c>0$, the {\it annular period counting function}%
\index{annular!period counting series}%
\index{counting function!annular!period} of $(\Ga,F,W)$ by
$$
s\mapsto \gls{annperiodcountingfucnt}(s)=
\sum_{g\in \Per(s)-\Per(s-c),\;g\cap W\neq \emptyset}\;e^{\int_g F}\;.
$$
Since $M$ is not assumed to be compact, it is important to assume that
$W$ is relatively compact in $T^1M$, since for instance when $M$ is an
infinite Galois Riemannian cover of a compact negatively curved
Riemannian manifold, and if $F$ is invariant under the covering group,
then $T^1M$ contains infinitely many periodic orbits of the same
period and of the same length at most $s$, for every $s$ large enough.
Conversely, if $\Ga$ is geometrically finite (see the definition of
geometrically finite discrete groups of isometries in Subsection
\ref{subsec:finiteness}), there exists a relatively compact open
subset of $T^1M$ meeting every periodic orbit. Note that
$Z_{\Ga,\,F,\,W}(s)$ is nonzero for $s$ large enough if and only if
$W$ meets the (topological) non-wandering set $\Omega\Ga$ of the
geodesic flow in $T^1M$. The maps $c\mapsto G_{\Ga,\,F,\,x,\,y,
  \,U,\,V,\,c} (t)$, $c\mapsto G_{\Ga,\,F,\,x,\,y, \,U,\,c} (t)$,
$c\mapsto G_{\Ga,\,F,\,x,\,y, \,c} (t)$ and $c\mapsto
Z_{\Ga,\,F,\,W,\,c}(s)$ are nondecreasing.

We define the {\it Gurevich pressure}\index{Gurevich pressure} of
$(\Ga,F)$ by
\begin{equation}\label{eq:defiGurepress}
\gls{Gurevichpressure}=
\limsup_{s\ra+\infty}\;\frac{1}{s} \ln Z_{\Ga,\,F,\,W,\,c}(s)\;,
\end{equation}
where $W$ is any relatively compact open subset of $T^1M$ meeting
$\Omega\Ga$ and $c$ any positive real number. The fact that the
Gurevich pressure does not depend on $c>0$ is proved as the first
claim of Assertion (vii) of Lemma \ref{lem:elemproppressure}. We will
prove in Subsection \ref{subsec:equalcritGur} that the Gurevich
pressure does not depend on $W$ too, and that the above upper limit is
a limit if $c$ is large enough (see also Remark \ref{rem:Gureprimit}
for the problem of counting primitive versus non-primitive periodic
orbits). Note that $\int_g(F+\kappa)=\int_gF +\kappa\ell(g)$ for every
periodic orbit $g$ of the geodesic flow on $T^1M$, and that
$P_{Gur}(\Ga,F+\kappa)=P_{Gur}(\Ga,F)+\kappa$ for every
$\kappa\in\RR$.

The Gurevich pressure is a (exponential) asymptotic growth rate of the
periodic orbits of the geodesic flow, weighted by the potential. For
Markov shifts on countable alphabets, it has been introduced by
Gurevich \cite{Gurevich69,Gurevich70} when the potential vanishes, and
by Sarig \cite{Sarig99,Sarig01a,Sarig01b} in general. These last two
works have motivated our study. To our knowledge, the Gurevich
pressure had not been studied in a non symbolic non-compact context.

\medskip We will give in Chapter \ref{sec:growth} precise asymptotic
results as $t$ goes to $+\infty$ of these counting functions, under
stronger hypotheses on $(\wt M,\Ga,\wt F)$. In the next subsection, we
will only give weaker (logarithmic) asymptotic results, valid in
general.

\medskip For instance, the following (very weak) result is an easy
consequence of Mohsen's shadow lemma \ref{lem:shadowlemma}.

\bcoro \label{coro:bigOloggrowth} 
For every $c>0$, we have $G_{\Ga,\,F,\,x,\,y, \,U,\,V,\,c} (t)=
O(e^{\delta_{\Ga,\,F}\;t})$ as $t$ goes to $+\infty$.  If
$\delta_{\Ga,\,F}>0$, we have $G_{\Ga,\,F,\,x,\,y, \,U,\,V,\,c}
(t)\leq G_{\Ga,\,F,\,x,\,y, \,U,\,V} (t)= O(e^{\delta_{\Ga,\,F}\;t})$
as $t$ goes to $+\infty$.  
\ecoro

\dem Since $G_{\Ga,\,F,\,x,\,y,\, U,\,V}\leq G_{\Ga,\,F,\,x,\,y}$
and $G_{\Ga,\,F,\,x,\,y,\, U,\,V,\,c}\leq G_{\Ga,\,F,\,x,\,y,\,c}$, we
may assume that $U=V=\partial_\infty \wt M$. We may also assume that
$\delta_{\Ga,\,F}<+\infty$. Let $(\mu_{x})_{x\in\wt M}$ be a Patterson
density of dimension $\sigma=\delta_{\Ga,\,F}$ (see Proposition
\ref{prop:existPattdens}) for $(\Ga,F)$. Then the result follows from
Corollary \ref{coro:dimgeqpress} (1), using a geometric series
argument when $\delta_{\Ga,\,F}>0$.  \cqfd

\subsection{Logarithmic growth of the orbital counting 
functions} 
\label{subsec:loggrowth}

In this subsection, we give logarithmic counting results which, though
less precise than the results in Subsection \ref{subsec:growth},
are valid in a much greater generality (they only require $\Ga$ to be
non-elementary).

\medskip The first result shows that the upper limit defining the
critical exponent $\delta_{\Ga,\,F}$ of $(\Ga,F)$ is in fact a
limit. Its proof follows closely the proof of the main result of
\cite{Roblin02} (corresponding to the case $F=0$).

\btheo \label{theo:limsupenfaitlim} Let $\wt M$ be a complete simply
connected Riemannian manifold, with dimension at least $2$ and pinched
sectional curvature at most $-1$, and $x,y\in\wt M$. Let $\Ga$ be a
non-elementary discrete group of isometries of $\wt M$.  Let $\wt F :
T^1\wt M\ra \RR$ be a H\"older-continuous $\Ga$-invariant map. If 
$\delta_{\Ga,\,F}>0$, then
$$
\delta_{\Ga,\,F}=\lim_{n\ra +\infty}\;
\frac{1}{n} \ln \;\sum_{\ga\in\Ga,\;d(x,\,\ga y)\leq n}
\;e^{\int_x^{\ga y} \wt F}\;.
$$  
\etheo

With the previous notation, this can be written as 
$$
\lim_{s\ra+\infty}\; \frac{1}{s} \ln \;G_{\Ga,\,F,\,x,\,y} (s)
=\delta_{\Ga,\,F}\;,
$$
that is, the orbital counting function grows logarithmically as
$s\mapsto e^{\delta_{\Ga,\,F}\,s}$. We will prove in Corollary
\ref{coro:loggrowth} that the sectorial and bisectorial counting
functions grow similarly (under obvious conditions).

\medskip \dem Let $\delta=\delta_{\Ga,\,F}>0$. For every $x'\in\wt M$,
with $a_n= \sum_{\ga\in\Ga,\;d(x',\,\ga y)\leq n} \;e^{\int_{x'}^{\ga y}
  \wt F}$, we have seen in Equation \eqref{eq:pressuresumleqn} that
$\delta= \limsup_{n\ra +\infty}\; \frac{1}{n} \ln a_n$.

For a contradiction, assume that $\liminf_{n\ra +\infty}\; \frac{1}{n}
\ln a_n <\delta$. Note that this lower limit does not depend on
$x'$. Hence there exist a sequence $(n_k)_{k\in\NN}$ of positive
integers and $\sigma \in\mathopen{]}-\infty,\delta\mathclose{[}$ such
that for every $x'\in\wt M$, for every $k\in\NN$ large enough,
\begin{equation}\label{eq:majoank}
a_{n_k}\leq e^{\sigma n_k}\;.
\end{equation}
Let us construct a Patterson density of dimension $\sigma$ for
$(\Ga,F)$, which contradicts Corollary \ref{coro:dimgeqpress} (2) since
$\sigma<\delta$.

Let $\D_z$ be the (unit) Dirac mass at a point $z\in\wt M$. For all
$t\in \mathopen{[}0,+\infty\mathclose{[}$ and $x'\in\wt M$, let
$$
\nu_{x',\,t} = \frac{\sum_{\ga \in \Ga, \;d(x',\,\ga y) \leq t}
\;e^{\int_{x'}^{\ga y}(\wt F-\sigma)}\;\D_{\ga y}}
{\sum_{\ga\in\Ga,\;d(y,\,\ga y)\leq t}\;
e^{\int_{y}^{\ga y}(\wt F-\sigma)}}\;.
$$
By compactness, for every $r\in \NN$, there exists a sequence
$\big(k_i(r)\big)_{i\in\NN}$ of positive integers such that the
probability measures $\nu_{y, \,n_{k_i(r)}-r}$ weak-star converge to a
probability measure $\mu_{y,\,r}$. Since the Poincar\'e series
$\sum_{\ga \in \Ga}\;e^{\int_{y}^{\ga y}(\wt F-\sigma)}$ diverges (as
$\sigma<\delta$), the support of the measure $\mu_{y,\,r}$ is
$\Lambda\Ga$. For every $x'\in \wt M$, let us define
$$
d\mu_{x',\,r}(\xi)=e^{C_{F-\sigma,\,\xi}(x',y)}\;d\mu_{y,\,r}(\xi)\;,
$$ 
and let us prove that $\nu_{x',\,n_{k_i(r)}-r}$ weak-star converges to
$\mu_{x',\,r}$ as $i\ra+\infty$ if $r\geq d(x',y)$.

Note that if $r\geq d(x',y)$, for every $k\in\NN$ large enough, by the
triangle inequality and Equation \eqref{eq:majoank}, we have
\begin{align*}
  &\Big\|\sum_{\ga\in\Ga,\;d(x',\,\ga y)\leq n_k-r}\;e^{\int_{x'}^{\ga
      y}(\wt F-\sigma)}\;\D_{\ga y} -\sum_{\ga\in\Ga,\;d(y,\,\ga y)\leq
    n_k-r}\;e^{\int_{x'}^{\ga y}(\wt F-\sigma)}\;\D_{\ga y}\Big\| \\ &\leq
  \sum_{\ga\in\Ga,\;n_k-r -d(x',\,y)\leq d(x',\,\ga y) \leq n_k-r
    +d(x',\,y)}\; e^{\int_{x'}^{\ga y}(\wt F-\sigma)} \leq \;e^{(2r-n_k)
    \sigma}\;a_{n_k}\leq e^{2r\sigma}\;.
\end{align*}
Since the denominator of $\nu_{x',\,n_{k_i(r)}-r}$ tends to $+\infty$
as $i\ra+\infty$ and since
$$
\Big|\int_{x'}^{\ga y}(\wt F-\sigma)-\int_{y}^{\ga y}(\wt F-\sigma)
-C_{F-\sigma,\,\xi}(x',y)\Big|
$$
is arbitrarily small if $\ga y$ is close enough to
$\xi\in \partial_\infty \wt M$, we have the convergence we looked for.

Now, for all $\ga\in\Ga$, $t\geq 0$ and $x'\in \wt M$, we clearly have
$\ga_*\nu_{x',\,t}=\nu_{\ga x',\,t}$.  Hence if $r\geq
\max\{d(x',y),d(\ga x',y)\}$, then $\ga_*\mu_{x',\,r}=\mu_{\ga
  x',\,r}$ by taking limits.  By compactness, there exists a sequence
$(r_j)_{j\in\NN}$ in $\mathopen{[}0,+\infty\mathclose{[}$ converging
to $+\infty$ such that the probability measures $\mu_{y,\,r_j}$
weak-star converge to a probability measure $\mu_{y}$. If
$(\mu_{x'})_{x'\in\wt M}$ is the family of measures on
$\partial_\infty \wt M$ defined by
$$
d\mu_{x'}(\xi)=e^{C_{F-\sigma,\,\xi}(x',y)}\;d\mu_{y}(\xi)\;,
$$
the sequence $(\mu_{x',\,r_j})_{j\in\NN}$ weak-star converges to
$\mu_{x'}$, and hence $(\mu_{x'})_{x'\in\wt M}$ is a Patterson density
of dimension $\sigma$ for $(\Ga,F)$, a contradiction.  \cqfd

\medskip \rem Here is a short proof of this result, using the
sub-additivity (or rather supermultiplicativity) ideas of
\cite{DalPeiSam10}. Let us fix $x,y\in\wt M$. This paper proves that

$\bullet$~ there exist a finite subset $P$ of $\Ga$ and $c\geq 0$ such
that for all $\alpha,\beta\in \Ga$, there exists $\ga\in P$ such that
the Hausdorff distance between the geodesic segment $\mathopen{[}x,
\alpha\ga\beta y\mathclose{]}$ and $\mathopen{[}x, \alpha
y\mathclose{]} \cup\mathopen{[}\alpha y,\alpha\ga x\mathclose{]}\cup
\mathopen{[}\alpha\ga x,\alpha\ga \beta y\mathclose{]}$ is at most
$c$;

$\bullet$ let $(b_n)_{n\in\NN}$ be a sequence of nonnegative real
numbers such that there exist $C>0$ and $N\in\NN$ such that for all
$n,m\in\NN$, we have 
$$
b_n b_m\leq C\sum_{i=-N}^{N} b_{n+m+i}\;,
$$ 
then, with $a_n=\sum_{k=0}^{n-1}b_n$, the limit of ${a_n}^{\frac{1}{n}}$
as $n\ra +\infty$ exists (and hence is equal to its upper limit).

The first point implies that there exists $c'\geq 0$ (depending on
$x,y$) such that for all $\alpha,\beta\in \Ga$, there exists $\ga\in
P$ satisfying
$$
\big|\,d(x,\alpha\ga\beta y)-d(x,\alpha y)- d(x, \beta y)\,\big|\leq c'
$$
and, since $\wt F$ is H\"older-continuous and bounded on compact
subsets of $\wt M$,
$$
\Big|\int_x^{\alpha\ga\beta y} \wt F -\int_x^{\alpha y} \wt F- 
\int_x^{\beta y}\wt F \Big|\leq c'\;.
$$
Hence the sequence $\big(b_n= \sum_{\ga\in\Ga,\;n-1< d(x,\,\ga y)\leq n}
\;e^{\int_{x}^{\ga y} \wt F}\big)_{n\in\NN}$ satisfies the assumptions
of the second point, and Theorem \ref{theo:limsupenfaitlim} follows.

\bigskip As mentioned by the referee (who provided its proof), we have
the following version of Theorem \ref{theo:limsupenfaitlim} when the
critical exponent is possibly nonpositive. It may also be deduced from
Theorem \ref{theo:limsupenfaitlim} by adding a big enough constant to
$F$.

\btheo \label{theo:limsupenfaitlimbis} Let $\wt M$ be a complete simply
connected Riemannian manifold, with dimension at least $2$ and pinched
sectional curvature at most $-1$, and $x,y\in\wt M$. Let $\Ga$ be a
non-elementary discrete group of isometries of $\wt M$.  Let $\wt F :
T^1\wt M\ra \RR$ be a H\"older-continuous $\Ga$-invariant map. For
every $c>$ large enough, we have
$$
\delta_{\Ga,\,F}=\lim_{n\ra +\infty}\;
\frac{1}{n} \ln \;\sum_{\ga\in\Ga,\;n-c<d(x,\,\ga y)\leq n}
\;e^{\int_x^{\ga y} \wt F}\;.
$$  
\etheo

Note that $c$ depends on $x,y$. With the previous notation, this can
be written as
$$
\lim_{s\ra+\infty}\; \frac{1}{s} \ln \;G_{\Ga,\,F,\,x,\,y,\,c} (s)
=\delta_{\Ga,\,F}\;,
$$
that is, the annular orbital counting function grows logarithmically
as $s\mapsto e^{\delta_{\Ga,\,F}\,s}$. We will prove in Corollary
\ref{coro:loggrowthbis} that the annular sectorial and bisectorial
counting functions grow similarly (under obvious conditions).

\medskip 
\dem 
Since $\Ga$ is non elementary, let $\ga_0$ be a loxodromic element of
$\Ga$, with translation axis $\Axe_{\ga_0}$ and translation length
$\ell(\ga_0)$. Let $c= \ell(\ga_0) + 2\,d(x,y) +
2\,d(y,\Axe_{\ga_0})$.  For every $n\in\NN$ with $n\geq d(x,y)$, let
us first prove that there exists $\ga\in\Ga$ such that $n-c\leq
d(x,\ga y)\leq n$.  The set $\{k'\in\NN \;:\; d(x,\ga_0^{k'}y)\leq
n\}$ being non empty and finite, let $k$ be its maximum and let us
prove that $\ga= \ga_0^k$ satisfies the required inequalities. By the
triangle inequality, for every $k'\in\NN$, we have
$$
-d(x,y)+k'\,\ell(\ga_0) \leq 
d(x,\ga_0^{k'}y)\leq d(x,y)+k'\,\ell(\ga_0)+2\,d(y,\Axe_{\ga_0})\,.
$$
Since $d(x,\ga_0^{k+1}y)> n$ by the maximality of $k$, we have
\begin{align*}
d(x,\ga_0^{k}y)&\geq -d(x,y)+k\,\ell(\ga_0) =
d(x,y)+(k+1)\ell(\ga_0)+2\,d(y,\Axe_{\ga_0})-c\\ & 
\geq d(x,\ga_0^{k+1}y)-c > n-c\,,
\end{align*}
as required.

Now, consider the sequence $\big(b_n= \sum_{\ga\in\Ga,\;n-c< d(x,\,\ga
  y) \leq n} \;e^{\int_{x}^{\ga y}\wt F} \big)_{n\in\NN}$, and note
that $b_n>0$ for all $n\geq d(x,y)$. As in the above short proof of
Theorem \ref{theo:limsupenfaitlim}, the result follows from the
following version of its second point: Let $(b_n)_{n\in\NN}$ be a
sequence of nonnegative real numbers such that there exist $C>0$ and
$N<N'$ in $\NN$  such that $b_{n}>0$ for every $n\geq N'$ and for all
$n,m\in\NN$, we have
$$
b_n b_m\leq C\sum_{i=-N}^{N} b_{n+m+i}\;,
$$ 
then for every $N''\geq N+N'-1$, with $a_n=\sum_{k=n-N''}^{n}b_n$,
the limit of $\frac{1}{n}\ln \,a_n$ as $n\ra +\infty$ exists.

By the standard comparison of a sum of a bounded number of nonnegative
real numbers and their maximum, this follows from the following lemma,
by taking $u_n=\frac{1}{C(2N+1)}\,b_n$.

\blemm\label{lem:appendixAreferee} 
Let $(u_n)_{n\in\NN}$ be a sequence of nonnegative real numbers
such that there exist $N<N'$ in $\NN$ such that $u_{n}>0$ for every
$n\geq N'$ and for all $n,m\in\NN$, we have
$$
u_n u_m\leq \max_{-N\leq i\leq N} u_{n+m+i}\;.
$$ 
Then, with $N''\geq N+N'$ and $v_n=\max_{n-N''< i\leq n} u_i$,
the limit of $\frac{1}{n}\ln \,v_n$ as $n\ra +\infty$ exists.
\elemm

\dem Let $\delta = \limsup_{n\ra+\infty} \frac{1}{n}\ln \,v_n$. Let us
prove that $\liminf_{n\ra+\infty} \frac{1}{n}\ln \,v_n\geq \delta$,
which implies the result. We may assume that $\delta>-\infty$. Let
$\sigma<\delta$ and $\epsilon\in\mathopen{]}0,1\mathclose{]}$ with
$2\epsilon\neq \sigma$.

Let $n'>N''+N$ be large enough so that $\max \{\frac{N|\sigma|}
{n''+N}, \;\frac{N''|\sigma|}{n''}\}\leq \epsilon$ for every
$n''\in\NN$ with $n'-N''\leq n''\leq n'$, and $\frac{1}{n'}\ln
\,v_{n'}\geq \sigma$. There exists $n''\in\NN$ with $n'-N''\leq
n''\leq n'$ and $v_{n'}=u_{n''}$, so that $n''> N$ and
$$
\frac{1}{n''}\ln\,u_{n''}=  \frac{n'}{n''}\,\frac{1}{n'}\ln \,v_{n'}
\geq \frac{n'}{n''}\,\sigma\geq\sigma-\frac{N''|\sigma|}{n''}\geq 
\sigma-\epsilon\,.
$$
There exists by induction an increasing sequence $(n_k)_{k\in\NN}$ in
$\NN$ starting with $n_0=n''$, such that $|n_{k+1}-(n_k+n_0)|\leq N$
and $u_{n_k}\,u_{n_0}\leq u_{n_{k+1}}$ for every $k\in\NN$. We have by
induction $n_k\leq (k+1)(n_0+N)$ and
$$
u_{n_k} \geq u_{n_0}^{k+1}=u_{n''}^{k+1}\geq e^{(\sigma-\epsilon)\,(k+1)\,n''} 
\geq e^{(\sigma-\epsilon-\,\frac{N\,|\sigma|}{n''+N})\, n_k}
\geq e^{(\sigma-2\,\epsilon) n_k}\,.
$$

Now, for every 
$$
n\geq \max\big\{\frac{|\sigma-2\,\epsilon|(2N+N'+n'')}
{\epsilon},\;\frac{(n''+N+1)|\ln u_{N'}|}{\epsilon}\big\}\,,
$$
let us prove that $\frac{1}{n}\ln \,v_n\geq \sigma-4\epsilon$, which
implies the result.  Let $k\in\NN$ be maximal so that $n_k\leq
n$. Again by induction, let $(m_i)_{i\in\NN}$ be an increasing
sequence such that $|m_{i+1}-(m_i+N')| \leq N$ and
$u_{n_k}\,u_{N'}^{i+1}\leq u_{n_{k}+m_i}$ for every $i\in\NN$. There
exists $\ell\in \NN$ such that $n -(N+N')< n_{k}+m_\ell\leq n$, and
$\ell\leq m_\ell\leq n_{k+1}-n_{k}\leq n'' +N$ by the maximality of
$k$. Then
$$
n\geq n_{k}\geq n -(N+N') -m_\ell\geq n-(2N+N'+n'')\geq
n\,\big(1-\frac{\epsilon}{|\sigma-2\,\epsilon|}\big)
$$ 
and, since $N''\geq N+N'$,
$$
v_n\geq u_{n_{k}+m_\ell}\geq u_{n_k}\,u_{N'}^{\ell+1}\geq 
e^{(\sigma-2\,\epsilon) n_k}\,e^{-(\ell+1)|\ln u_{N'}|}\geq 
e^{(\sigma-4\epsilon)\,n}\,,
$$
as required. \cqfd

\bcoro\label{coro:loggrowth} With the assumptions of Theorem
\ref{theo:limsupenfaitlim}, 

(1) if $U$ is an open subset of $\partial_\infty \wt M$ meeting
$\Lambda\Gamma$, then $\lim_{t\ra+\infty} \frac{1}{t}\ln
G_{\Ga,\,F,\,x,\,y, \,U}(t) = \delta_{\Ga,\,F}$;

(2) if $U$ and $V$ are any two open subsets of $\partial_\infty \wt M$
meeting the limit set $\Lambda\Gamma$, then $\lim_{t\ra+\infty}
\frac{1}{t}\ln G_{\Ga,\,F,\,x,\,y, \,U,\,V}(t) = \delta_{\Ga,\,F}$.  
\ecoro

\dem 
We adapt the proof of \cite[Coro.~1]{Roblin02} which corresponds
to the case $F=0$.  The validity of these statements for all open sets
$U$ and $V$ is independent of $x$ and $y$, since for all $z,z'\in\wt
M$ and all open subsets $W,W'$ of $\partial_\infty \wt M$ such that
$\overline{W}\subset W'$, there exists a compact subset $K$ of $\wt M$
such that $\C_zW\subset K\cup\C_{z'}W'$, and since the intersection of
$K$ with any orbit of $\Ga$ is a finite.  Recall that the convex hull
$\C\Lambda\Ga$ of $\Lambda\Ga$ is $\Ga$-invariant. We may assume that
$y\in \C\Lambda\Ga$, hence that the orbit of $y$ is contained in
$\C\Lambda\Ga$.

\medskip (1) For every $t\in\mathopen{[}0,+\infty\mathclose{[}\,$, for
every open subset $U'$ of $\wt M\cup \partial_\infty \wt M$ and for
every $z\in\wt M$, let
$$
a_{t,\,U',\,z}= \sum_{\ga\in\Ga\;:\;d(z,\ga y)\leq t\,,\; \ga y\in U'}
e^{\int_{z}^{\ga y}\wt F}\;.
$$
Note that for every $\ga\in\Ga$, by Equation \eqref{eq:timereversal}
and a change of variable in this sum, we have
\begin{equation}\label{eq:equivmescontcon}
a_{t,\,U',\,z}=a_{t,\,\ga U',\,\ga z}\;.
\end{equation}
When $U'$ is the open set $\C_xU-\{x\}$, the difference
$G_{\Ga,\,F,\,x,\,y,\,U}(t)- a_{t,\,U',\,x}$ is equal to $1$ if $x$
belongs to the $\Ga$-orbit of $y$, and $0$ otherwise. Since
$G_{\Ga,\,F,\,x,\,y,\,U}(t)\leq G_{\Ga,\,F,\,x,\,y}(t)$ and by Theorem
\ref{theo:limsupenfaitlim}, we have
$$
\limsup_{t\ra+\infty} \frac{1}{t}\ln G_{\Ga,\,F,\,x,\,y,\,U}(t) \leq
\delta_{\Ga,\,F}\;.
$$ 
Hence to prove the first assertion of Corollary
\ref{coro:loggrowth}, we only have to prove that
$$
\liminf_{t\ra+\infty} \frac{1}{t}\ln a_{t,\,U',\,x} \geq
\delta_{\Ga,\,F}
$$ 
for every open subset $U'$ of $\wt M\cup \partial_\infty \wt M$
meeting $\Lambda\Gamma$.

Let $U'$ be such an open subset. Since $\Ga$ is non-elementary, every
orbit in $\Lambda \Ga$ is dense. Hence by compactness of $\Lambda
\Ga$, there exist $\ga_1,\dots,\ga_k$ in $\Ga$ such that $\Lambda \Ga
\subset \bigcup_{i=1}^k\ga_i U'$. In particular $\C\Lambda\Ga -
\bigcup_{i=1}^k\ga_i U'$ is compact. Hence there exists $c'_1\geq 0$
such that for every $t\in\mathopen{[}0,+\infty\mathclose{[}\,$, we
have
\begin{equation}\label{eq:techmincontcon}
G_{\Ga,\,F,\,x,\,y}(t)\leq c'_1+\sum_{i=1}^k \;a_{t,\,\ga_i U',\,x}\;.
\end{equation}

Let $r=\max_{1\leq i\leq k} \;d(x,\ga_i x)$ and $t\in\mathopen{[}r,
+\infty\mathclose{[}\,$.  Note that for all $i\in\{1,\dots, k\}$ and
$\ga\in\Ga$, by the triangle inequality, if $d(\ga_i^{-1}x,\ga y)\leq
t-r$, then $d(x,\ga y)\leq t$. Furthermore, by Lemma
\ref{lem:technicholder}, there exists $c'_2>0$ such that for all
$i\in\{1,\dots, k\}$ and $\ga\in\Ga$, we have
$$
\Big|\int_x^{\ga y}\wt F-\int_{\ga_i^{-1} x}^{\ga y}\wt F\Big|
\leq c'_2\;.
$$
Hence, for all $i\in\{1,\dots, k\}$ and $t\in\mathopen{[}r,
+\infty\mathclose{[}\,$, we have, using Equation
\eqref{eq:equivmescontcon},
$$
a_{t-r,\,\ga_i U',\,x}=
a_{t-r,\,U',\,\ga_i^{-1}x}\leq e^{c'_2}\;a_{t,\,U',\,x}\;.
$$
Therefore, by Equation \eqref{eq:techmincontcon}, we have
\begin{equation}\label{eq:findemloggrowthsect}
a_{t,\,U',\,x}\geq 
\frac{1}{k\;e^{c'_2}}\big(G_{\Ga,\,F,\,x,\,y}(t-r) -c'_1\big)\;.
\end{equation}
By taking the logarithm, dividing by $t$ and taking the lower limit
as $t\ra+\infty$, the result then follows from Theorem
\ref{theo:limsupenfaitlim}.

\medskip (2) The proof of the second assertion reduces to the first
one, using a similar approach as the reduction from the first one to
Theorem \ref{theo:limsupenfaitlim}. For every
$t\in\mathopen{[}0,+\infty\mathclose{[}\,$, for all open subsets
$U',V'$ of $\wt M\cup \partial_\infty \wt M$ and for all $z,w\in\wt
M$, we now introduce
$$
b_{t,\,U',\,V',\,z,\,w}= 
\sum_{\ga\in\Ga\;:\;d(z,\,\ga w)\leq 
t\,,\; \ga w\in U',\;\ga^{-1} z\in V'}
e^{\int_{z}^{\ga w}\wt F}\;,
$$
which satisfies $b_{t,\,U',\,V',\,z,\,w}=b_{t,\,U',\,\alpha V',\,z, \,
  \alpha w}$ for every $\alpha\in\Ga$. As above, we only have to prove
that 
$$
\liminf_{t\ra+\infty} \frac{1}{t}\ln b_{t,\,U',\,V',\,x,\,y}
\geq \delta_{\Ga,\,F}
$$
for all open subsets $U',V'$ of $\wt M \cup \partial_\infty \wt M$
meeting $\Lambda\Gamma$.

Let $U',V'$ be two such open subsets. As above, there exist
$\alpha_1,\dots,\alpha_\ell$ in $\Ga$ such that $\Lambda \Ga \subset
\bigcup_{i=1}^\ell\alpha_i V'$, and hence there exists $c'_3\geq 0$
such that for every $t\in\mathopen{[}0,+\infty\mathclose{[}\,$, we
have
$$
a_{t,\,U',\,x}\leq 
c'_3+\sum_{i=1}^\ell \;b_{t,\,U',\,\alpha_i V',\,x,\,y}\;.
$$
Now let $r=\max_{1\leq i\leq \ell} \;d(y,\alpha_i y)$.  For all
$i\in\{1,\dots, \ell\}$, $t\in\mathopen{[}r,+\infty\mathclose{[}$ and
$\ga\in\Ga$, by the triangle inequality, if $d(x,\ga \alpha_i^{-1}y)
\leq t-r$, then $d(x,\ga y)\leq t$. Furthermore, by Equation
\eqref{eq:timereversal} and by Lemma \ref{lem:technicholder}, there
exists $c'_4>0$ such that for all $i\in\{1,\dots, \ell\}$ and
$\ga\in\Ga$, we have
$$
\Big|\int_x^{\ga y}\wt F-\int_{x}^{\ga \alpha_i^{-1} y}\wt F\Big|
=\Big|\int_{y}^{\ga^{-1} x}\wt F\circ\iota-
\int_{\alpha_i^{-1}y}^{\ga^{-1} x}\wt F\circ\iota\Big| \leq c'_4\;.
$$
Hence, for all $i\in\{1,\dots, k\}$ and $t\in\mathopen{[}r,
+\infty\mathclose{[}\,$, we have
$$
b_{t-r,\,U',\,\alpha_i V',\,x,\,y}=
b_{t-r,\,U',\,V',\,x,\,\alpha_i^{-1}y}\leq 
e^{c'_4}\;b_{t,\,U',\,V',\,x,\,y}\;.
$$
Therefore
$$
b_{t,\,U',\,V',\, x,\,y}\geq 
\frac{1}{\ell\; e^{c'_4}}\big(a_{t-r,\,U',\,x} -c'_3\big)\;,
$$
and we conclude, using Assertion (1), as in the end of the proof of
this assertion. \cqfd

\medskip We have the following version of Corollary
\ref{coro:loggrowth} when the critical exponent is possibly
nonpositive. It may also be deduced from Corollary
\ref{coro:loggrowth} by adding a big enough constant to $F$.

\bcoro\label{coro:loggrowthbis} With the assumptions of Theorem
\ref{theo:limsupenfaitlimbis}, if $c>0$ is large enough,

(1) if $U$ is an open subset of $\partial_\infty \wt M$ meeting
$\Lambda\Gamma$, then $\lim_{t\ra+\infty} \frac{1}{t}\ln
G_{\Ga,\,F,\,x,\,y, \,U,\,c}(t) = \delta_{\Ga,\,F}$;

(2) if $U$ and $V$ are any two open subsets of $\partial_\infty \wt M$
meeting the limit set $\Lambda\Gamma$, then $\lim_{t\ra+\infty}
\frac{1}{t}\ln G_{\Ga,\,F,\,x,\,y, \,U,\,V,\,c}(t) = \delta_{\Ga,\,F}$.  
\ecoro

\dem (1) The proof is similar to that of Corollary
\ref{coro:loggrowth} (1). We define, for $t,U',z$ as previously and
$c>0$ large enough,
$$
a_{t,\,U',\,z,\,c}= \sum_{\ga\in\Ga\;:\;t-c<d(z,\ga y)\leq t\,,\; \ga y\in U'}
e^{\int_{z}^{\ga y}\wt F}\;.
$$
We have $a_{t,\,U',\,z,\,c}=a_{t,\,\ga U',\,\ga z,\,c}$ for every $\ga\in\Ga$,
and we also only have to prove that 
$$
\liminf_{t\ra+\infty} \frac{1}{t}\ln a_{t,\,U',\,x,\,c} \geq
\delta_{\Ga,\,F}\,,
$$
since the converse inequality is immediate, as
$G_{\Ga,\,F,\,x,\,y,\,U,\,c}\leq G_{\Ga,\,F,\,x,\,y,\,c}$.  With
$\ga_1,\dots,\ga_k,r$ and $c'_2$ as previously, since $t-r-(c-2r)<
d(\ga_i^{-1}x,\ga y)\leq t-r$ implies that $t-c<d(x,\ga y) \leq t$, we
have $a_{t-r,\,\ga_iU',\,z,\,c-2r} \leq e^{c'_2}\;a_{t,\,U',\,
  z,\,c}$. Hence
$$
a_{t,\,U',\,x,\,c}\geq \frac{1}{k\;e^{c'_2}}
\big(G_{\Ga,\,F,\,x,\,y,\,c-2r}(t-r) -c'_1\big)\,,
$$
and the result follows from Theorem \ref{theo:limsupenfaitlimbis} as
previously.
  
\medskip (2) The proof is similar to that of Corollary
\ref{coro:loggrowth} (2) by considering, for $t$, $U'$, $V'$, $z$,
$r$, $\ell$, $c'_3$, $c'_4 $ as in its proof and $c>0$ large enough,
\begin{equation}\label{eq:bisecmap}
b_{t,\,U',\,V',\,z,\,w,\,c}= 
\sum_{\ga\in\Ga\;:\;t-c<d(z,\,\ga w)\leq t\,,\; \ga w\in U',\;
\ga^{-1} z\in V'}
e^{\int_{z}^{\ga w}\wt F}\,,
\end{equation}
and by proving as previously that
$$
b_{t,\,U',\,V',\, x,\,y,\,c}\geq 
\frac{1}{\ell\; e^{c'_4}}\big(a_{t-r,\,U',\,x,\,c-2r} -c'_3\big)$$
so that
\begin{equation}\label{eq:pourgure}
\liminf_{t\ra+\infty} \frac{1}{t}\ln b_{t,\,U',\,V',\,x,\,y,\,c}
\geq \delta_{\Ga,\,F}\;.\;\;\;\Box
\end{equation}

\subsection{Equality between critical exponent and 
Gurevich pressure}
\label{subsec:equalcritGur}

The aim of this subsection is stated in its title: we now prove that
the logarithmic growth rate of the periodic geodesics (weighted by the
potential) is equal to the logarithmic growth rate of the orbit points
(weighted by the potential).

\btheo \label{theo:equalcritgur}
Let $\wt M$ be a complete simply connected Riemannian manifold,
with dimension at least $2$ and pinched sectional curvature at most
$-1$, and $x,y\in\wt M$. Let $\Ga$ be a non-elementary discrete group
of isometries of $\wt M$.  Let $\wt F : T^1\wt M\ra \RR$ be a
H\"older-continuous $\Ga$-invariant map. Let $W$ be a relatively
compact open subset of $T^1M$ meeting the (topological) non-wandering set
$\Omega\Ga$. Then
$$
P_{Gur}(\Ga,F)=\delta_{\Ga,\,F}\,.
$$
If $c>0$ is large enough, then
$$
P_{Gur}(\Ga,F)=\lim_{s\ra+\infty}\;\frac{1}{s} \ln Z_{\Ga,\,F,\,W,\,c}(s)\,,
$$
and if $P_{Gur}(\Ga,F)>0$, then
$$
P_{Gur}(\Ga,F)=\lim_{s\ra+\infty}\;\frac{1}{s} \ln Z_{\Ga,\,F,\,W}(s)\,.
$$
\etheo

This proves in particular that the Gurevich pressure does not depend
on $W$ and that the upper limit defining it (see Equation
\eqref{eq:defiGurepress}) is a limit if $c>0$ is large enough. 

\medskip
\dem Let $\prT: T^1\wt M\ra T^1M=\Ga\bs T^1\wt M$ be the canonical
projection.

Let us first prove that 
$$
\limsup_{s\ra+\infty}\;\frac{1}{s} \ln
Z_{\Ga,\,F,\,W,\,c}(s)\leq \delta_{\Ga,\,F}
$$ 
if $c>0$ is large enough, and that
$$
\limsup_{s\ra+\infty}\;\frac{1}{s} \ln
Z_{\Ga,\,F,\,W}(s)\leq \delta_{\Ga,\,F}
$$  
if $\delta_{\Ga,\,F}>0$. Since $W$ is relatively compact, there
exists a compact subset $K$ of $\wt M$ such that $\prT(\pi^{-1}(K))$
contains $W$. Let $r$ be the diameter of $K$, and fix $x\in K$.

For all $s\geq 0$ and $g$ in $\Per(s)$ such that $g\cap W\neq
\emptyset$, let $\ga_g$ be one of the loxodromic elements of $\Ga$
whose translation axis $\Axe_{\ga_g}$ meets $K$, whose translation
length is the length $\ell(g)$ of $g$ and such that for every $y\in
\Axe_{\ga_g}$, the image by $\prT$ of the unit tangent vector at $y$
pointing towards $\ga_g y$ belongs to $g$. Note that the number of
these elements $\ga_g$ is at least equal to the cardinality of the
pointwise stabiliser of $\Axe_{\ga_g}$, that is to the multiplicity of
$g$.

Let $x_g$ be the closest point to $x$ on $\Axe_{\ga_g}$, which
satisfies $d(x,x_g)\leq r$ since $x\in K$ and $\Axe_{\ga_g}$ meets
$K$. We have by the triangle inequality
$$
\ell(g)\leq d(x,\ga_g x)\leq d(x,x_g)+d(x_g,\ga_g x_g)+d(\ga_gx,\ga_gx_g)
\leq \ell(g)+2r\leq s+2r\;.
$$
Furthermore, by (two applications of) Lemma \ref{lem:technicholder}
with $r_0=r$, there exists a constant $c'_5\geq 0$ (depending only on
$r$, the H\"older constants of $\wt F$, the bounds on the sectional
curvature and $\max_{\pi^{-1}(B(x,\,r))}|\wt F|$) such that
$$
\Big|\;\int_x^{\ga_g x}\wt F\;-\int_{g}\wt F\;\Big|
=\Big|\;\int_x^{\ga_g x}\wt F\;-\int_{x_g}^{\ga_g x_g}\wt F\;\Big|
\leq c'_5\;.
$$
Hence $Z_{\Ga,\,F,\,W}(s)\leq e^{c'_5}\;G_{\Ga,\,F,\,x,\,x}(s+2r)$ and
$Z_{\Ga,\,F,\,W,\,c}(s)\leq e^{c'_5}\;G_{\Ga,\,F,\,x,\,x,\,c+2r}
(s+2r)$, which proves our first claim, by Theorems
\ref{theo:limsupenfaitlim} and \ref{theo:limsupenfaitlimbis}.

\medskip Since $Z_{\Ga,\,F,\,W,\,c}\leq Z_{\Ga,\,F,\,W}$ for every
$c>0$, in order to prove Theorem \ref{theo:equalcritgur}, we now only
have to prove that, if $c>0$ is large enough,
$$
\liminf_{s\ra+\infty}\;\frac{1}{s} \ln Z_{\Ga,\,F,\,W,\,c}(s)\geq
\delta_{\Ga,\,F}\;.
$$
Let $v\in T^1\wt M$ such that $\prT(v)\in W\cap \Omega\Ga$, and now let
$x=\pi(v)$. Note that the points $v_-$ and $v_+$ both belong to
$\Lambda\Ga$. 

By standard arguments (see for instance Lemma
\ref{lem:crithyperboelem} and \cite[page 150-151]{GhyHar90}), there
exist $U'$ and $V'$ small enough neighbourhoods in $\wt
M\cup\partial_\infty \wt M$ of $v_+$ and $v_-$ respectively, such that
for every $\ga\in\Ga$ such that $\ga x\in U'$ and $\ga^{-1} x\in V'$,
then $\ga$ is a loxodromic element and $v$ is close to its
translation axis $\Axe_\ga$, in the sense that the point $x$ is at
distance at most $1$ from some point $x_\ga$ in $\Axe_\ga$, and that if
$v_\ga$ is the unit tangent vector at $x_\ga$ pointing towards $\ga
x_\ga$, then $p(v_\ga)\in W$ (recall that $W$ is open). Note that
$U'$ and $V'$ meet $\Lambda\Ga$.

Let $s\geq 0$ and let $\ga\in\Ga$ be such that $s-(c-2)<d(x,\ga x)\leq
s$, $\ga x\in U'$ and $\ga^{-1} x\in V'$. Note that the orbit $g_\ga$
under the geodesic flow of $\prT(v_\ga)$ is periodic, and that its
length satisfies
$$
s-c<d(x,\ga x)-2\leq \ell(g_\ga)\leq d(x,\ga x)\leq s\;,
$$
and as above, there exists $c'_6\geq 0$ such that
$$
\Big|\;\int_{g_\ga}\wt F\;-\int_x^{\ga x}\wt F\;\Big|
=\Big|\;\int_{x_\ga}^{\ga x_\ga}\wt F\;-\int_{x}^{\ga x}\wt F\;\Big|
\leq c'_6\;.
$$
Hence $Z_{\Ga,\,F,\,W,\,c}(s)\geq e^{-c'_6}\;
b_{s,\,U',\,V',\,x,\,x,\,c-2}$ with the notation of Equation
\eqref{eq:bisecmap}, which proves our second claim, by Equation
\eqref{eq:pourgure}.  
\cqfd

\medskip The next result is an obvious corollary of Theorem
\ref{theo:equalcritgur}, also following from Remark \ref{rem:livsi}
and the remark at the end of Subsection
\ref{subsec:GibbsPoincareseries}.

\bcoro Under the assumption of the previous theorem, if $\wt F^* :
T^1\wt M\ra \RR$ is another H\"older-continuous $\Ga$-invariant map
with induced map $F^*:T^1M\ra\RR$, such that the periods of every $\ga
\in \Ga$ for $\wt F$ and $\wt F^*$ are the same, then the critical
exponents of $(\Ga,F)$ and $(\Ga,F^*)$ are the same:
$$
\delta_{\Ga,\,F} = \delta_{\Ga,\,F^*} \;.
$$
\ecoro

\subsection{Critical exponent of Schottky semigroups}
\label{subsec:semigroup}

The main aim of this subsection is to prove (see Theorem
\ref{theo:critexpsubsemigroup}) that the critical exponent
$\delta_{\Ga,\,F}$ is the upper bound of the critical exponents
$\delta_{G,\,F}$ (defined as for groups) of finitely generated
subsemigroups $G$ of $\Ga$, which have strong geometric properties (of
Schottky type). The constructions of this subsection will also be
useful for the proof of the Variational Principle in Subsection
\ref{subsec:proofvaraprincip}. 

Though the similar problem for free subgroups instead of free
subsemigroups is open in our case, it is in general a much harder one,
when dealing with growth problems. For instance, there exist finitely
generated soluble groups with exponential growth, and though the most
frequent way to prove they have this type of growth is by proving that
they contain free subsemigroups on (at least) two generators, they
simply do not contain any free subgroup on two generators. A famous
result due to Doyle is that there is a gap away from $2$ for the
critical exponents of the classical Schottky subgroups (which are
free) of the non-elementary Kleinian groups (that is, of
non-elementary discrete subgroups of $\operatorname{PSL}_2(\CC)$, seen
as the orientation preserving group of isometries of the real
hyperbolic space $\HH^3_\RR$ of dimension $3$). Hence free subgroups
have rigidity properties that free semigroups do not have (see for
instance \cite{Mercat12}).

\medskip
Let $\wt M$ be a complete simply connected Riemannian manifold, with
dimension at least $2$ and pinched negative sectional curvature at
most $-1$.  Let $G$ be a discrete semigroup of isometries of $\wt M$.
We refer to \cite{Mercat12} for general information about discrete
semigroups of isometries of CAT$(-1)$-spaces.

Let $\wt F :T^1\wt M\ra \RR$ be a H\"older-continuous $G$-invariant
map.  We define the {\it limit set}\index{limit set!of a semigroup}
$\Lambda G$ of $G$ (as the set of accumulation points in
$\partial_\infty\wt M$ of any orbit $G x_0$), the {\it Poincar\'e
  series}\index{Poincar\'e series!of a semigroup}
$$
Q_{G,\,F,\,x,\,y}(s)=\sum_{\ga\in G} \;\; e^{\int_x^{\ga y} (\wt F-s)}
$$ 
of $(G,F)$, and the {\it critical exponent}\index{critical exponent!of
  a semigroup}
$$
\delta_{G,\,F}=\limsup_{n\ra
  +\infty}\;\frac{1}{n}\ln \;\sum_{\ga\in G,\;n-1 <d(x,\,\ga y)\leq n}
\;e^{\int_x^{\ga y} \wt F}
$$ 
of $(G,F)$, exactly as for the case of groups.

Recall that a map $f:X\ra Y$ between metric spaces is {\it
  quasi-isometric}\index{quasi-isometric} if there exist $\lambda\geq
1$ and $c\geq 0$ such that for all $x,y\in X$, we have
$$
-c+\frac{1}{\lambda}\;d(x,y)\leq d(f(x),f(y))\leq \lambda\;d(x,y) +c\;.
$$
If $G$ is generated (as a semigroup) by a finite set $S$, and is
endowed with the (semigroup) word metric defined by $S$ (see below
when $G$ is free on $S$), we say that $G$ is {\it
  convex-cocompact}\index{convex-cocompact!semigroup} if the map from
$G$ to $\wt M$ defined by $g\mapsto g x_0$ is quasi-isometric (for any
$x_0\in \wt M$).

\medskip
We first give a construction of convex-cocompact free semigroups of
isometries of $\wt M$ (of Schottky type). We have not tried to get
optimal constants.

\bprop \label{prop:schottkysemigroups} For every $\epsilon\in
\mathopen{]}0,1\mathclose{]}$, there exists $\theta>0$ such that for
every set $S$ of isometries of $\wt M$, for all $N\geq 5$ and $C\geq
6$, for all $x_0\in \wt M$ and $v_0\in T^1_{x_0}\wt M$ such that the
following conditions are satisfied:
 \begin{enumerate}
\item[(1)] 
  $N\leq d(x_0,\alpha x_0)< N+1$ for every $\alpha\in S$,
\item[(2)] 
  $d(\alpha x_0,\beta x_0)\geq C$ for all distinct $\alpha$
  and $\beta$ in $S$,
\item[(3)] 
  for every $\alpha\in S$, the angle at $x_0$ between $v_0$
  and the unit tangent vector at $x_0$ pointing towards $\alpha x_0$,
  and the angle at $\alpha x_0$ between $-\alpha v_0$ and the unit
  tangent vector at $\alpha x_0$ pointing towards $x_0$, are at most
  $\theta$ (see the picture below),
\end{enumerate}
then the semigroup $G$ generated by $S$ is free on $S$ and
convex-cocompact, and for all $\ga,\ga'$ in $G$, the piecewise
geodesic segment $\mathopen{[}x_0,\ga x_0\mathclose{]}\cup
\mathopen{[}\ga x_0,\ga\ga'x_0\mathclose{]}$ is contained in the
$\epsilon$-neighbourhood of $\mathopen{[}x_0,\ga\ga'x_0\mathclose{]}$.
Furthermore
$$
k(N-2\epsilon)\leq d(x_0,\ga x_0)\leq k(N+1)
$$ 
for every word $\ga$ in the elements of $S$ of length $k$.  Lastly,
for all words $\ga$ and $\ga'$ in the elements of $S$, none of them
being an initial subword of the other, there exist $x\in
\mathopen{[}x_0,\ga x_0\mathclose{]}$ and $x'\in \mathopen{[}x_0,\ga'
x_0\mathclose{]}$ such that $d(x_0,x)=d(x_0,x')$ and $d(x,x')\geq
1-8\epsilon$.  \eprop

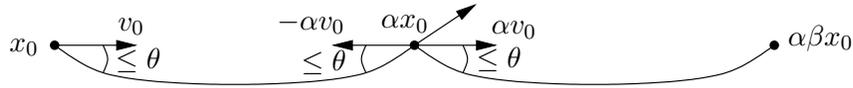
\begin{center}
\input{fig_construcschotky.pstex_t}
\end{center}

\dem Fix $\epsilon\in\mathopen{]}0,1\mathclose{]}$, and let
$\theta=\min\{\frac{\pi}{4}, \frac{1}{2}
\theta_0(5,\epsilon,\frac{\pi}{4})\}$, where
$\theta_0(\cdot,\cdot,\cdot)$ has been defined in Lemma
\ref{lem:quasigeodpresquegeod}. Let $S,N,C,x_0,v_0$ be as in the
statement. Note that the assumptions (1) and (2) imply that $S$ is
finite. Before proving that $G$ is free and convex-cocompact, we start
with a few preliminary remarks.

For every word $m=\alpha_1\dots\alpha_k$ in elements of $S$, we denote
by $\ell(m)=k\in\NN$ its length (the only word with zero length being
the empty word) and, for $0\leq j\leq k$, by $m_j=\alpha_1\dots
\alpha_j$ its initial subword of length $j$ (in particular $m_0$ is
the empty word).  If $m=\alpha_1\dots\alpha_k$ and
$m'=\beta_1\dots\beta_{k'}$ are two words in elements of $S$, let
$$
j(m,m')=\max\{i\in\NN\;:\;\forall\; i'\leq i,\;\;\alpha_{i'}=
\beta_{i'}\}
$$
be the length of the maximal common initial subword of $m$ and
$m'$. It is $0$ if and only if $m$ or $m'$ is the empty word or
$\alpha_1\neq \beta_1$. The distance between the words $m$
and $m'$ is
$$
d(m,m')=\ell(m)+\ell(m')-2j(m,m')\;.
$$

Consider the piecewise geodesic path
$$
\omega_m=\mathopen{[}m_0x_0,m_1 x_0\mathclose{]}\cup
\mathopen{[}m_1x_0,m_2 x_0\mathclose{]}
\cup\dots\cup \mathopen{[}m_{k-1}x_0,m_kx_0\mathclose{]}
$$
(with $\omega_m=\{x_0\}$ by convention if $m$ is the empty word).  By
Assumption (3), the exterior angles of $\omega_m$ at the points $m_i
x_0$ for $1\leq i\leq k-1$ are at most $2\theta$ (see the picture
above). By Assumption (1), the length of each segment of $\omega_m$ is
at least $N\geq 5$. Hence, by Lemma \ref{lem:quasigeodpresquegeod} and
the definition of $\theta$, the piecewise geodesic path $\omega_m$ is
contained in the $\epsilon$-neighbourhood of
$\mathopen{[}x_0,mx_0\mathclose{]}$. In particular, this proves the
ante-penultimate assertion of Proposition \ref{prop:schottkysemigroups}, by
convexity.

Since $N\geq 5$, $\epsilon\leq 1$ and $2\theta\leq \frac{\pi}{2}$, the
closest points to $m_ix_0$ on $\mathopen{[}x_0,mx_0\mathclose{]}$ for
$0\leq i\leq k$ are in this order on this segment.

\smallskip\noindent
\begin{minipage}{12cm}~~~ Indeed, if $x,y,z$ are points in $\wt M$
  such that $N\leq d(x,y)\leq N+1$, $N\leq d(y,z)\leq N+1$,
  $\measuredangle_y(x,z)\geq \frac{\pi}{2}$, having closest points
  respectively $p,q,r$ at distance at most $\epsilon$ on a geodesic
  segment, with the absurd hypothesis that $r\in \mathopen{[}p,
  q\mathclose{]}$, then by the triangle inequality and since closest
  point maps do not increase distances, we have
\begin{align*}
d(x,z)&\leq d(p,r)+2\epsilon=d(p,q)-d(r,q)+2\epsilon\leq 
d(x,y)-d(z,y)+4\epsilon
\\ &\leq N+1-N+4\epsilon=1+4\epsilon\leq 5\leq N \;.
\end{align*}
By an angle comparison, the angle $\measuredangle_y(x,z)$ is less than
$\frac{\pi}{2}$, a contradiction. [Another argument is that since
$\measuredangle_y(x,z)\geq\frac{\pi}{2}$, we have by a distance
comparison that $d(x,z)\geq d(x,y)+d(y,z)-2\log (1+\sqrt{2})\geq
2N-2\log (1+\sqrt{2})>N$, a contradiction.]
\end{minipage}
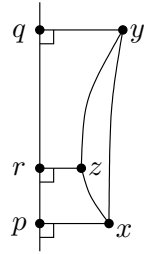
\begin{minipage}{2.9cm}
\begin{center}
\input{fig_nonretour.pstex_t}
\end{center}
\end{minipage}

\bigskip We now give our last preliminary remark, which yields the
penultimate assertion of Proposition \ref{prop:schottkysemigroups}. By
the triangle inequality and Assumption (1), we have
$$
d(mx_0,x_0)\leq \sum_{i=0}^{k-1}d(m_ix_0,m_{i+1}x_0)\leq
k(N+1)=(N+1)\,\ell(m)\;,
$$
on one hand, and on the other hand, with $p_i$ the closest point to
$m_ix_0$ on $\mathopen{[}x_0,mx_0\mathclose{]}$ for $0\leq i\leq k$,
using the above ordering remark,
\begin{align*}
d(mx_0,x_0)&=\sum_{i=0}^{k-1}d(p_i,p_{i+1})\geq 
\sum_{i=0}^{k-1}(d(m_ix_0,m_{i+1}x_0)-2\epsilon)
\\ &\geq k(N-2\epsilon)= (N-2\epsilon)\,\ell(m)\;.
\end{align*}
Note that $N-2\epsilon>0$ by the assumptions on $N$ and $\epsilon$.
In particular, if $m$ is not the empty word, then $mx_0\neq x_0$.

\medskip After these preliminary remarks, let us prove that $G$ is
free on $S$. Let $m=\alpha_1\dots\alpha_k$ and $m'=\beta_1\dots
\beta_{k'}$ be two distinct words in elements of $S$, and let us prove
that $mx_0\neq m'x_0$, which yields the result.

Assume for a contradiction that $mx_0= m'x_0$.  Let $j=j(m,m')$. Up to
multiplying $m$ and $m'$ by the inverse of their maximal common
initial subword, we may assume that $j=0$. We have $\min\{k,k'\}>0$,
otherwise, up to exchanging $m$ and $m'$, we have $k\neq 0$ and
$k'=0$, so that $mx_0=x_0$ and $m$ is not the empty word, a
contradiction to the last preliminary remark. Hence $\alpha_1$ and
$\beta_1$ exist and are distinct.

\medskip\noindent
\begin{minipage}{7.9cm} ~~~ Let $p$ and $q$ be the closest points to
  $\alpha_1 x_0$ and $\beta_1 x_0$ on $\mathopen{[}x_0, m x_0
  \mathclose{]}$.  Up to permuting $m$ and $m'$, we may assume that
  $x_0,p,q,m x_0$ are in this order on $\mathopen{[}x_0,
  mx_0\mathclose{]}$.
\end{minipage}
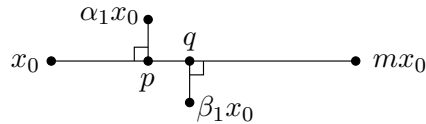
\begin{minipage}{7cm}
\begin{center}
\input{fig_alignmot.pstex_t}
\end{center}
\end{minipage}

\medskip By the triangle inequality, since closest point
maps do not increase distances, by Assumption (1) and by the
properties of $\epsilon$ and $C$, we have
\begin{align*}
d(\alpha_1 x_0,\beta_1x_0)&\leq d(p,q)+2\epsilon=
d(q,x_0)-d(p,x_0)+2\epsilon\\ &
\leq d(\beta_1x_0,x_0)-d(\alpha_1 x_0,x_0)+3\epsilon
\leq N+1 -N+3\epsilon=3\epsilon+1< C\;.
\end{align*}
This contradicts Assumption (2). Hence $G$ is indeed free on $S$.

\medskip Let us now prove that the map from $G$ to $\wt M$
defined by $g\mapsto g x_0$ is quasi-isometric, which says that $G$ is
convex-cocompact. Let $m=\alpha_1\dots\alpha_k$ and $m'= \beta_1\dots
\beta_{k'}$ be two words in elements of $S$.  Let us prove that
\begin{equation}\label{eq:semiquasiisometri}
(N-2\epsilon)d(m,m')-2(\ln(1+\sqrt{2})+N-\epsilon)
\leq d(mx_0,m'x_0)\leq (N+1)\,d(m,m')\;,
\end{equation}
which implies the result. 

Let $j=j(m,m')$. Up to multiplying $m$ and $m'$ by the inverse of
their maximal common initial subword, we may assume that $j=0$. By the
last preliminary remark, we may assume that neither $m$ nor $m'$ is
the empty word (we would not have to do this in the group case, but
$G$ is here only a semigroup).  We may also assume that $m\neq
m'$. Hence $\alpha_1$ and $\beta_1$ exist and are distinct.

\medskip\noindent
\begin{minipage}{6.4cm} ~~~ Let $p$ be the closest point to $\alpha_1
  x_0$ on $\mathopen{[}x_0,m x_0\mathclose{]}$.  Let $r$ and $q$ be
  the closest points to $p$ and $\beta_1 x_0$ respectively on
  $\mathopen{[}x_0,m' x_0\mathclose{]}$. Let $\overline{p}$ and
  $\overline{q}$ be the closest points to $p$ and $q$ respectively on
  $\mathopen{[}m x_0,m' x_0\mathclose{]}$.
\end{minipage}
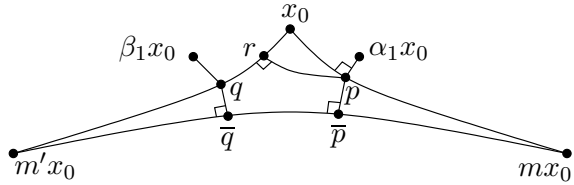
\begin{minipage}{8.5cm}
\begin{center}
\input{fig_semiconvcomp.pstex_t}
\end{center}
\end{minipage}

\medskip Consider the geodesic triangle in the CAT$(-1)$-space $\wt M$
with vertices $mx_0,x_0,m'x_0$. Each of its sides is contained in the
$\ln(1+\sqrt{2})$-neighbourhood of the union of the two other sides.

Assume for a contradiction that $d(p,r)\leq \ln(1+\sqrt{2})$. As when we
proved that $G$ is free, we have
$$
d(\alpha_1 x_0,\beta_1x_0)\leq d(r,q)+2\epsilon+\ln(1+\sqrt{2})=
|d(q,x_0)-d(r,x_0)|+2\epsilon+\ln(1+\sqrt{2})\;.
$$
But by Assumption (1), we have $N-\epsilon\leq d(q,x_0)\leq N+1$,
$$ 
d(r,x_0)\leq d(p,x_0)\leq d(\alpha_1 x_0,x_0)\leq N+1
$$ 
and
\begin{align*}
d(r,x_0)& \geq d(p,x_0)-\ln(1+\sqrt{2})\geq
d(\alpha_1 x_0,x_0)-\epsilon-\ln(1+\sqrt{2}) \\ &\geq
N-\epsilon-\ln(1+\sqrt{2})
\;.
\end{align*}
In particular, $|d(q,x_0)-d(r,x_0)|\leq 1+\epsilon+ \ln(1+\sqrt{2})$
and, since $\epsilon\leq 1$,
$$
d(\alpha_1 x_0,\beta_1x_0)\leq 1+3\epsilon+2\ln(1+\sqrt{2})<6\;.
$$
Since $C\geq 6$, this contradicts Assumption (2).

By the above property of the geodesic triangles, we hence have
$d(p,\overline p)\leq \ln(1+\sqrt{2})$. Similarly, $d(q,\overline
q)\leq \ln(1+\sqrt{2})$.

Let us prove the upper bound in Equation \eqref{eq:semiquasiisometri}.
Indeed, by the triangle inequality and the last preliminary remark,
since $j(m,m')=0$, we have
$$
d(mx_0,m'x_0)\leq d(mx_0,x_0)+d(x_0,m'x_0)\leq 
(N+1)k+(N+1) k'=(N+1)d(m,m')\;.
$$

Let us prove the lower bound in Equation
\eqref{eq:semiquasiisometri}. By convexity, the points $m
x_0,\overline{p},\overline{q},m'x_0$ are in this order on
$\mathopen{[}m x_0,m' x_0\mathclose{]}$. Hence by the triangle
inequality, since $\alpha_1$ is the first letter of $m$ and $\beta_1$
is the first letter of $m'$, by the last preliminary remark, and since
$j(m,m')=0$, we have
\begin{align*}
d(mx_0,m'x_0)&\geq d(mx_0,\overline{p})+d(\overline{q},m'x_0)\geq
d(mx_0,{p})+d({q},m'x_0)-2\ln(1+\sqrt{2})\\ &
\geq
d(m x_0,\alpha_1 x_0)+d(\beta_1x_0,m'x_0)-2\ln(1+\sqrt{2})-2\epsilon
\\ &\geq (N-2\epsilon)(k-1)+(N-2\epsilon)(k'-1)
-2\ln(1+\sqrt{2})-2\epsilon
\\ &=(N-2\epsilon)d(m,m')-2(\ln(1+\sqrt{2})+N-\epsilon)\;,
\end{align*}
as required.

\medskip Let us finally prove the last assertion of Proposition
\ref{prop:schottkysemigroups}. Let $m$ and $m'$ be two words in the
elements of $S$, none of them being an initial subword of the
other. They are nontrivial (since the empty word is an
initial subword of any word), and we may write $m= m_0\alpha_1m_1$ and
$m'=m_0\alpha'_1m'_1$ where $m_0,m_1,m'_1$ are (possibly empty) words
in the elements of $S$ with $\ell(m_0)=j(m,m')$ and
$\alpha_1,\alpha'_1$ are distinct elements of $S$. Let
$$a=d(x_0,m_0x_0)\;\;\;{\rm and}\;\;\;
b=\min\{d(m_0x_0,m x_0),\;d(m_0x_0,m'x_0)\}\;.
$$ 
Let $y$ and $y'$ be the points of $\mathopen{[}m_0x_0,m
x_0\mathclose{]}$ and $\mathopen{[}m_0x_0, m' x_0\mathclose{]}$
respectively at distance $b$ from $m_0x_0$. Let $z$ and $z'$ be the
closest points to $y$ and $y'$ on $\mathopen{[}x_{0},m
x_{0}\mathclose{]}$ and $\mathopen{[}x_{0},m' x_{0}\mathclose{]}$
respectively. Let $w$ and $w'$ be the closest points to $m_0x_0$ on
$\mathopen{[}x_{0},m x_{0}\mathclose{]}$ and $\mathopen{[}x_{0},m'
x_{0}\mathclose{]}$ respectively. Let $x$ and $x'$ be the points on
$\mathopen{[}x_{0},m x_{0}\mathclose{]}$ and $\mathopen{[}x_{0},m'
x_{0}\mathclose{]}$ respectively at the same distance $a+b-3\epsilon$
from $x_{0}$. Let us prove that $x$ and $x'$ are well defined and
satisfy the requirements of the last assertion of Proposition
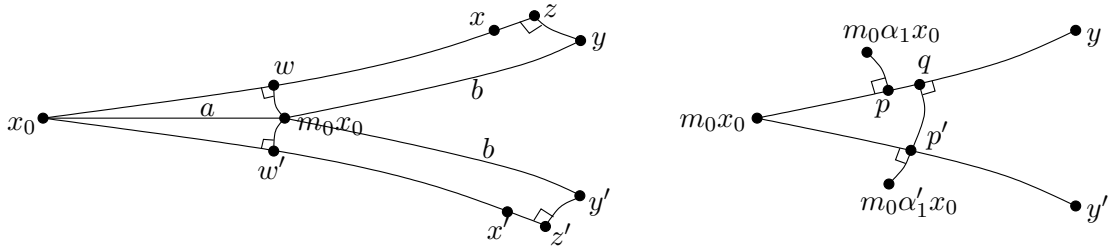
\ref{prop:schottkysemigroups}.

\begin{center}
\input{fig_lotsepardyna.pstex_t}
\end{center}

We have
$$
\max\{d(m_0x_0,w),d(m_0x_0,w'),d(y,z),d(y',z')\}\leq \epsilon
$$
since $\mathopen{[}x_0,m_0 x_0\mathclose{]}\cup \mathopen{[}m_0x_0, m
x_0\mathclose{]}$ is contained in the $\epsilon$-neighbourhood of
$\mathopen{[}x_0,m x_0\mathclose{]}$ and similarly upon replacing $m$
by $m'$.

By the last preliminary remark, we have
$$
b=\min\{d(x_0,\alpha_1m_1 x_0),\;d(x_0,\alpha'_1m'_1x_0)\}\geq N-2\epsilon
\geq 3\epsilon\;.
$$
Hence $a+b-3\epsilon\geq 0$, and the points $x_0,w,z$ and $x_0,w',z'$
are in this order on $\mathopen{[}x_{0},m x_{0}\mathclose{]}$ and
$\mathopen{[}x_{0},m' x_{0}\mathclose{]}$ respectively.  By the
triangle inequality, and since closest point maps do not increase the
distances, we have
\begin{align*}
d(x_0,x)&=a+b-3\epsilon= d(x_0,m_0x_0)+d(m_0x_0,y)-3\epsilon\\ &\leq 
d(x_0,w)+d(w,z)+2d(m_0x_0,w)+d(y,z)-3\epsilon\\ &
\leq d(x_0,z)=d(x_0,w)+d(w,z)\leq a+b\;.
\end{align*}
In particular, $a+b-3\epsilon\leq d(x_0,z)\leq d(x_0,m x_0)$ and $x$
is well defined.  Similarly, $x'$ is well defined.

By the triangle inequality,
\begin{align*}
d(x,y)&\leq d(x,z)+d(z,y)=(d(x_0,z)-d(x_0,x))+d(z,x)\\ &\leq 
a+b-(a+b-3\epsilon)+\epsilon=4\epsilon\;.
\end{align*}
Similarly, $d(x',y')\leq 4\epsilon$, hence again by the triangle inequality
$$
d(x,x')\geq d(y,y')-d(x,y)-d(x',y')\geq d(y,y')-8\epsilon\;. 
$$
Let us prove that $d(y,y')\geq 1$, which implies the last
assertion of Proposition \ref{prop:schottkysemigroups}.

Assume for a contradiction that $d(y,y')< 1$. Let $p$ and $p'$ be the
closest points to $m_0\alpha_1 x_0$ and $m_0\alpha'_1 x_0$ on
$\mathopen{[}m_0x_0, m x_0\mathclose{]}$ and $\mathopen{[}m_0x_0, m'
x_0\mathclose{]}$ respectively. We have
$$
\max\{d(m_0\alpha_1x_0,p),\;d(m_0\alpha'_1x_0,p')\}\leq \epsilon\;.
$$
Up to permuting $m$ and $m'$, since $y=m x_0$ or $y'=m'x_0$ by the
definition of $y$ and $y'$, we may assume that $p'\in
\mathopen{[}m_0x_0, y'\mathclose{]}$. Let $q$ be the closest point to
$p'$ on $\mathopen{[}m_0x_0, m x_0\mathclose{]}$. By convexity, we
have $d(p',q)\leq d(y',y)$.  By the triangle inequality, since closest
point maps do not increase distances, and by the hypothesis (1), we
have
\begin{align*}
d(p,q)&=|d(m_0x_0,p)-d(m_0x_0,q)|
\leq |d(m_0x_0,p)-d(m_0x_0,p')|+d(p',q)
 \\ &
\leq |d(m_0x_0,m_0\alpha_1x_0)-d(m_0x_0,m_0\alpha'_1x_0)|+ \epsilon + d(p',q)
 \\ &< ((N+1)-N)+ \epsilon + d(y',y)= 1+ \epsilon + d(y',y)\;.
\end{align*}
Therefore
\begin{align*}
d(\alpha_1x_0,\alpha'_1x_0)&=d(m_0\alpha_1x_0,m_0\alpha'_1x_0)
\leq d(m_0\alpha_1x_0,p)+d(p,q)+d(q,p')+d(p',m_0\alpha'_1x_0) \\& \leq 
1+ 3\epsilon + 2 d(y',y)< 6\leq C\;.
\end{align*}
This contradicts the hypothesis (2).
This ends the proof of Proposition \ref{prop:schottkysemigroups}.
\cqfd

\bigskip We now construct (in Proposition
\ref{prop:constructsemigroup}) finite subsets $S$ of $\Ga$ satisfying
the assumptions (1)-(3) of Proposition \ref{prop:schottkysemigroups}.
These sets $S$ satisfy another property (4), which will allow us to
prove in Theorem \ref{theo:critexpsubsemigroup} that the critical
exponent $\delta_{G,\,F}$, where $G$ is the subsemigroup generated by
$S$, is close to the critical exponent $\delta_{\Ga,\,F}$.

We will follow arguments contained in the proof of
\cite[Theo.~1]{OtaPei04}, though the necessity of controlling the
integral of the potential $\wt F$ along pieces of orbits of the
geodesic flow, and the new application Theorem
\ref{theo:critexpsubsemigroup}, require a much more precise
construction. The origin of the techniques goes back to the paper
\cite{BisJon97}, as extended independently by U.~Hamenst\"adt
(unpublished) and \cite{Paulin97d} to prove that the Hausdorff
dimension of the radial limit set $\Lambda_c\Ga$ of $\Ga$ is the
(standard) critical exponent of $\Ga$.

\bprop\label{prop:constructsemigroup}  
For all $\delta'\in\mathopen{]}-\infty, \delta_{\Ga,\,F}\mathclose{[}$
and $\theta,C>0$, there exist $x_0\in \wt M$ and $v_0\in T^1_{x_0}\wt
M$ such that for every $L>0$, there exist $N>L$ and a finite subset
$S$ of $\Ga$ such that
\begin{enumerate}
\item[(1)] 
  $N\leq d(x_0,\alpha x_0)< N+1$ for every $\alpha\in S$;
\item[(2)] 
  $d(\alpha x_0,\beta x_0)\geq C$ for all distinct $\alpha$
  and $\beta$ in $S$;
\item[(3)] 
  for every $\alpha\in S$, the angle at $x_0$ between $v_0$
  and the unit tangent vector at $x_0$ pointing towards $\alpha x_0$,
  and the angle at $\alpha x_0$ between $-\alpha v_0$ and the unit
  tangent vector at $\alpha x_0$ pointing towards $ x_0$, are at most
  $\theta$;
\item[(4)] 
  $\displaystyle \sum_{\alpha\in S} e^{\int_{x_0}^{\alpha x_0}  \wt F}\; 
  \geq e^{\delta' N}$.
\end{enumerate}  
\eprop

\dem Let us fix a real constant $\delta'<\delta= \delta_{\Ga,\,F}$ and
positive ones $\theta,C$. First assume that $\delta'>0$ (see the end
of the proof to remove this assumption).  We may assume that $\theta<
\frac{\pi}{2}$. Let $\epsilon,\theta'>0$ be small enough parameters,
and let $L'$ be a large enough parameter.

\medskip\noindent
\begin{minipage}{9.9cm} ~~~ Since $\wt M$ is CAT$(-1)$, there exist
  $\rho_0=\rho_0(\theta)>0$ and
  $\rho=\rho(\theta,\theta')\geq\rho_0+1$ such that for every $x'\in
  \wt M$, for all $v',w'\in T^1_{x'}\wt M$ whose angle
  $\measuredangle_{x'}(v',w')$ at $x'$ is a least $\frac{\theta}{4}$,
  if $y'=\pi(\phi_\rho v')$ and $z'=\pi(\phi_\rho w')$, then the
  angles at $y'$ and $z'$ of the geodesic triangle with vertices
  $x',y',z'$ are at most $\theta'$, and the distance from $x'$ to the
  opposite side $\mathopen{[}y',z'\mathclose{]}$ is at most $\rho_0$.
\end{minipage}
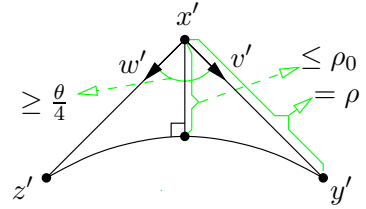
\begin{minipage}{5cm}
\begin{center}
\input{fig_trianglefin.pstex_t}
\end{center}
\end{minipage}

\medskip
Since $\Ga$ is non-elementary, there exist two distinct points $\xi$
and $\eta$ in $\Lambda\Ga$ such that the orbit of $(\xi,\eta)$ in
$\Lambda\Ga\times\Lambda\Ga$ under the diagonal action of $\Ga$ is
dense in $\Lambda\Ga\times\Lambda\Ga$.  Let $u$ be a unit tangent
vector to the geodesic line from $\xi$ to $\eta$, and let $x=\pi(u)$.

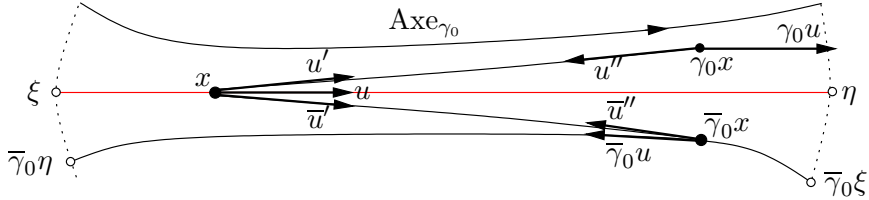
\begin{center}
\input{fig_udensebis.pstex_t}
\end{center}

\medskip Again since $\Ga$ is non-elementary, the pairs of endpoints of
translation axes of loxodromic elements of $\Ga$ are dense in
$\Lambda\Ga\times\Lambda\Ga$. In particular there exists $\ga_0\in\Ga$
whose translation axis $\Axe_{\ga_0}$ has endpoints close to $\xi$ and
$\eta$.  Up to replacing $\ga_0$ by a positive or negative power, we
may assume that $\ga_0$ translates from the endpoint close to $\xi$
towards the endpoint close to $\eta$, and that the geodesic segment
$\mathopen{[}x,\ga_0 x\mathclose{]}$ has length at least $L'$ and is
contained in the $\epsilon$-neighbourhood of $\Axe_{\ga_0}$.  Denote
by $u'$ the unit tangent vector at $x$ pointing towards $\ga_0x$ and
by $u''$ the unit tangent vector at $\ga_0x$ pointing towards $x$. We
may also assume that the angles $\measuredangle_x(u,u')$ and
$\measuredangle_{\ga_0 x} (\ga_0u, -u'')$ are at most
$\frac{\theta}{8}$. In particular, by the spherical triangle
inequality, we have
\begin{equation}\label{eq:angleuprimuprimprim}
\angle_{\ga_0 x}(u'', -\ga_0 u')\leq\angle_{\ga_0 x}(u'', -\ga_0 u)+
\angle_{\ga_0 x}(-\ga_0 u, -\ga_0 u')\leq 2\;\frac{\theta}{8}
\leq\frac{\theta}{2}\;.
\end{equation}

By density and up to taking powers, there exists $\overline{\ga}_0\in
\Ga$ such that $\overline \ga_0\xi$ is close to $\eta$,
$\overline\ga_0 \eta$ is close to $\xi$, the geodesic segment
$\mathopen{[}x,\overline\ga_0 x\mathclose{]}$ has length at least
$L'$, and, with $\overline u'$ the unit tangent vector at $x$ pointing
towards $\overline\ga_0x$ and $\overline u''$ the unit tangent vector
at $\overline\ga_0x$ pointing towards $x$, the angles
$\measuredangle_x(u,\overline u')$ and $\measuredangle_{\overline\ga_0
  x} (\overline\ga_0u,\overline u'')$ are at most $\frac{\theta}{8}$.

\medskip Since $\wt M$ has pinched negative curvature and $\rho$
depends only on $\theta,\theta'$, there exists $\theta''= \theta''
(\epsilon, \theta,\theta') \in \mathopen{]}0, \frac{\theta}{8}
\mathclose{]}$ such that for all $v',v''$ in $T^1_x\wt M$ with angle
$\measuredangle_x (v',v'')$ at most $\theta''$, for every $t\in
\mathopen{[}0,\rho\mathclose{]}$ (or equivalently for $t=\rho$ by
convexity), we have
$$
d(\pi(\phi_tv'),\pi(\phi_t v''))\leq \epsilon\;.
$$

For every $n\in\NN$, let
$$
A_n=\{y\in \wt M\;:\; n\leq d(x,y)<n+1\}\;.
$$
We have
$$
\limsup_{n\ra+\infty}\sum_{\ga\in\Ga\;:\;\ga x\in A_n}
\;e^{\int_{x}^{\ga x} (\wt F-\delta')}=+\infty\;.
$$ Indeed, by the definition of the critical exponent
$\delta=\delta_{\Ga,\,F}$ of $(\Ga,F)$ (at the beginning of Subsection
\ref{subsec:GibbsPoincareseries}), if $\delta''\in
\mathopen{]}\delta',\delta\mathclose{[}$, there exists an increasing
    sequence $(n_k)_{k\in\NN}$ in $\NN$ such that $\frac{1}{n_k} \ln
    \sum_{\ga\in\Ga,\;\ga x\in A_{n_k}} \;e^{\int_x^{\ga y} \wt F}
    \geq \delta''$, and hence $\sum_{\ga\in\Ga\;:\;\ga x\in A_{n_k}}
    \;e^{\int_{x}^{\ga x} (\wt F-\delta')} \geq
    e^{(\delta''-\delta')n_k-\delta'}$, which tends to $+\infty$ as
    $k\ra+\infty$.

Since the isometric action of $\Ga$ on $\wt M$ is proper, the group
$\Ga$ is the union of finitely many subsets $E$ such that
\begin{equation}\label{eq:defiCprim}
d(\alpha x,\beta x)\geq C'=
C+\max\{\,2\,\rho_0,\;2\,d(x,\ga_0 x),\;2\,d(x,\overline{\ga}_0 x)\,\}
\end{equation}
for all $\alpha\neq \beta$ in $E$.
There exists such a subset $E_0$ satisfying
$$
\limsup_{n\ra+\infty}\sum_{\ga\in E_0\;:\;\ga x\in A_n}
\;e^{\int_{x}^{\ga x} (\wt F-\delta')}=+\infty\;.
$$

Consider the sequence of the finite subsets $Z_n$ of the unit tangent
vectors at $x$ pointing towards the elements of $E_0x\,\cap A_n$, for
$n\in\NN$. By looking at the accumulation points of these subsets and
using the compactness of $T^1_x\wt M$, there exist $v,w$ in $T^1_x\wt
M$, an increasing sequence $(n_k)_{k\in\NN}$ in $\mathopen{[}2\,
\rho_0, +\infty\mathclose{[} \cap \,\NN$ and, for every $k\in\NN$, a
finite subset $S_k$ of $\Ga$ such that
\begin{enumerate}
\item[(i)] 
  for all $\alpha$ and $\beta$ in $S_k$, we have $n_k\leq
  d(x,\alpha x)\leq n_k+1$ and $d(\alpha x,\beta x)\geq C'$ if
  $\alpha\neq \beta$;
\item[(ii)] 
  for every $\alpha\in S_k$, the angle at $x$ between $v$
  and the unit tangent vector $v_\alpha$ at $x$ pointing towards
  $\alpha x$, and the angle at $\alpha x$ between $\alpha w$ and the
  unit tangent vector $w_\alpha$ at $\alpha x$ pointing towards $x$,
  are at most $\theta''$;
\item[(iii)] 
    $\displaystyle \sum_{\alpha\in S_k} e^{\int_{x}^{\alpha x}\wt F}\; 
    \geq k\;e^{\delta' n_k}$.
\end{enumerate}

The end of the proof of Proposition \ref{prop:constructsemigroup},
which would be easy if $v=-w$, is subdivided in a few cases, the first
of them being the easiest case when $v=-w$ almost holds, and the other
ones dealing with the situation when the equality $v=-w$ is far from
holding, and where we try to go back to the easiest case by rounding up
finitely many geodesic segments in $M$ with the same origin and
endpoint, and making them composable in an almost geodesic way.

\medskip 
\noindent{\bf Case 1. } First assume that $\measuredangle_x(v,w)\geq
\frac{\theta}{4}$.

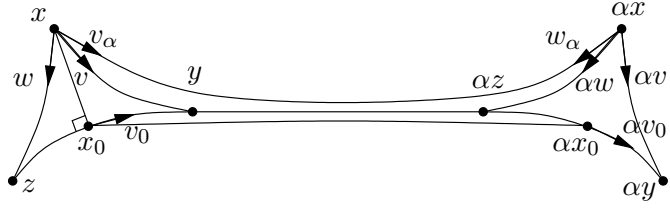
\begin{center}
\input{fig_HamPaucase1.pstex_t}
\end{center}

Let $y=\pi(\phi_\rho v)$ and $z=\pi(\phi_\rho w)$. Define $x_0$ to be the
closest point to $x$ on $\mathopen{[}y,z\mathclose{]}$ and $v_0$ as
the unit tangent vector at $x_0$ pointing towards $y$. In particular,
by the definition of $\rho_0$, we have $d(x,x_0)\leq \rho_0$.

Let $k\in\NN$ and $\alpha,\beta\in S_k$. By the triangle inequality
and Assertion (i) above, we have
\begin{equation}\label{eq:majodeplacecasun}
0<n_k-2\rho_0\leq d(x_0,\alpha x_0)\leq n_k+1+2\rho_0\;,
\end{equation}
and, if $\alpha\neq \beta$, by the definition of $C'$ in Equation
\eqref{eq:defiCprim},
$$
d(\alpha x_0,\beta x_0)\geq C'-2\rho_0\geq C\;.
$$

If $k$ is large enough, then $d(x,\alpha x)$ is large enough by
Assertion (i) and, since $\measuredangle_x(v,v_\alpha)\leq \theta''$
and $\measuredangle_{\alpha x}(\alpha w,w_\alpha)\leq \theta''$ by
Assertion (ii) above, the points $y$ and $\alpha z$ are, by the
definition of $\theta''$, at distance at most $\epsilon$ to the
geodesic segment $\mathopen{[}x,\alpha x\mathclose{]}$.

If two geodesic segments in a non-positively curved manifold have
lengths at least $1$ and have sufficiently close endpoints, then their
tangent vectors (parametrising them proportionally to arc-length by
$\mathopen{[}0,1\mathclose{]}$) are uniformly close. In particular (or
by using Lemma \ref{lem:technicholder}), since $\rho\geq 1$, if
$\epsilon$ is small enough and if $k$ is large enough, there exists a
constant $c_2>0$ (depending only on $\epsilon$, the H\"older constants
of $\wt F$ and $\max_{\pi^{-1}(B(y,\,\epsilon)\,\cup
  \,B(z,\,\epsilon))}\,|\wt F|$) such that
\begin{equation}\label{eq:controlintegralpotentun}
  \Big|\int_{x}^{\alpha x} \wt F-\int_{x}^{y} \wt F-
\int_{y}^{\alpha z} \wt F-\int_{\alpha z}^{\alpha x} \wt F\;\Big|\leq c_2\;.
\end{equation}

By the triangle inequality and the definition of $\rho$, we have
$d(x_0,y)\geq \rho-\rho_0\geq 1$, and similarly $d(\alpha z,\alpha
x_0) \geq 1$. If $k$ is large enough, we have, for every $\alpha \in
S_k$,
$$
d(y,\alpha z) \geq d(x,\alpha x) -2\rho\geq n_k-2\rho\geq 1\;.
$$
Recall that $d(y,\mathopen{[}x,\alpha x\mathclose{]})\leq \epsilon$,
$d(\alpha z,\mathopen{[}x,\alpha x\mathclose{]}) \leq \epsilon$, and
that $\measuredangle_y(\phi_\rho v, \phi_{d(x_0,y)} v_0) \leq \theta'$
by the definition of $\rho$.  Let $\epsilon'>0$ be small enough, and
assume that $\theta'$ and $\epsilon$ are small enough compared to
$\theta$ and $\epsilon'$.  Then, if $k$ is large enough, the piecewise
geodesic path
$$
\omega=\mathopen{[}x_0,y\mathclose{]}\cup \mathopen{[}y,
\alpha z\mathclose{]}\cup \mathopen{[}\alpha z,\alpha x_0\mathclose{]}
$$ 
has its exteriors angles at $y$ and $\alpha z$ at most $\theta_0(1,
\epsilon', \theta)$, where $\theta_0(\cdot,\cdot,\cdot)$ has been
defined in Lemma \ref{lem:quasigeodpresquegeod}.  Therefore by this
lemma, for every $\alpha \in S_k$, the angles at $x_0$ between $v_0$
and the unit tangent vector at $x_0$ pointing towards $\alpha x_0$,
and the angle at $\alpha x_0$ between $-\alpha v_0$ and the unit
tangent vector at $\alpha x_0$ pointing towards $ x_0$, are at most
$\theta$; furthermore the piecewise geodesic path $\omega$ is at
distance at most $\epsilon'$ from $\mathopen{[}x_0,\alpha
x_0\mathclose{]}$.

As in the proof of Equation \eqref{eq:controlintegralpotentun}, there
exists a constant $c_3>0$ (depending only on $\epsilon'$, the H\"older
constants of $\wt F$ and $\max_{\pi^{-1}(B(y,\,\epsilon')\,\cup
  \,B(z,\,\epsilon'))}\,|\wt F|$) such that
\begin{equation}\label{eq:controlintegralpotentdeu}
  \Big|\int_{x_0}^{\alpha x_0} \wt F-\int_{x_0}^{y} \wt F-
\int_{y}^{\alpha z} \wt F-\int_{\alpha z}^{\alpha x_0} \wt F\;\Big|\leq c_3\;.
\end{equation}
With $c_4=\rho\max_{\pi^{-1}(B(x,\,\rho))}|\wt F|$, since $d(x_0,y)\leq
d(x,y)=\rho$ and similarly with $z$ instead of $y$, we have
\begin{equation}\label{eq:controlintegralpotenttro}
 \Big|\int_{x_0}^{y} \wt F\;\Big|,
\;\;\Big|\int_{x}^{y} \wt F\;\Big|,
\;\;\Big|\int_{\alpha z}^{\alpha x} \wt F\;\Big|,\;\;
\Big|\int_{\alpha z}^{\alpha x_0} \wt F\;\Big|\leq c_4\;.
\end{equation}

By the formulas \eqref{eq:controlintegralpotentun},
\eqref{eq:controlintegralpotentdeu},
\eqref{eq:controlintegralpotenttro}, and by Assertion (iii) above,
we have, if $k$ is large enough,
$$
\sum_{\alpha\in S_k} e^{\int_{x_0}^{\alpha x_0}\wt F}\; \geq
e^{-c_2-c_3-4c_4}\sum_{\alpha\in S_k} e^{\int_{x}^{\alpha x}\wt F}\;
\geq k\;e^{-c_2-c_3-4c_4}\;e^{\delta' n_k}\;.
$$
Hence, using Equation \eqref{eq:majodeplacecasun}, there exists $N\in
\,\NN\,\cap\mathopen{]}n_k-2\rho_0-1,n_k+2\rho_0+2\mathclose{[}$
such that
\begin{align*}
 \sum_{\tiny\begin{array}{c}\alpha\in S_k\\ N\leq d(x_0,\alpha x_0)< N+1\end{array}} 
e^{\int_{x_0}^{\alpha x_0}\wt F}\; 
    \geq  \frac{k\;e^{-c_2-c_3-4c_4}}{4\rho_0+3}\;e^{\delta' n_k}\geq  
\frac{k\;e^{-c_2-c_3-4c_4-\delta'(2\rho_0+2)}}{4\rho_0+3}\;e^{\delta' N}\;.
\end{align*}
Given $L>0$, if $k$ is large enough, then
$\frac{k\;e^{-c_2-c_3-4c_4-\delta'(2\rho_0+2)}}{4\rho_0+3}\geq 1$
and $N\geq n_k-2\rho_0-1> L$, and defining
$$
S=\{\alpha \in S_k\;:\; N\leq d(x_0,\alpha x_0)< N+1\}\;,
$$
the objects $x_0,v_0,N,S$ defined above satisfy the properties (1)-(4)
required in Proposition \ref{prop:constructsemigroup}.

\medskip 
\noindent{\bf Case 2. } Assume that $\measuredangle_x(u,-v)\geq
\frac{\theta}{2}$ and that $\measuredangle_x(u,w)\geq
\frac{\theta}{2}$.

Consider the loxodromic element $\ga_0$ of $\Ga$ and the unit tangent
vectors $u'$ and $u''$ defined above. We will define in due course the
notation in the following figure.

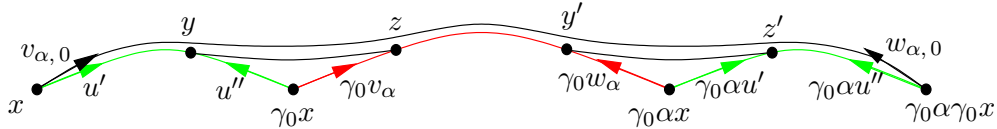
\begin{center}
\input{fig_HamPaucase2.pstex_t}
\end{center}

Let $k\in\NN$ and $\alpha\in S_k$. In order to get the following
inequalities, we use: 

$\bullet$~ the spherical triangle inequality for the first one;

$\bullet$~ the assumption of Case 2, the properties of $u'$ and $u''$,
and Assertion (ii) above for the second one;

$\bullet$~ and the assumption on $\theta''$ for the last one:
\begin{align*}
\measuredangle_{\ga_0x}(u'',\ga_0v_\alpha)&\geq
\measuredangle_{\ga_0x}(-\ga_0u,\ga_0v)
-\measuredangle_{\ga_0x}(-\ga_0u,u'')
-\measuredangle_{\ga_0x}(\ga_0v_\alpha,\ga_0v)\\ & \geq
\frac{\theta}{2}-\frac{\theta}{8}-\theta''\geq \frac{\theta}{4}\;.
\end{align*}
Similarly, $\measuredangle_{\ga_0\alpha x}(\ga_0w_\alpha,\ga_0\alpha u')\geq 
\frac{\theta}{4}$.

Hence, by the definition of $\rho$, the angles at $y=\pi(\phi_\rho
u'')$ and $z=\pi(\phi_\rho (\ga_0 v_\alpha))$ of the geodesic triangle
with vertices $\ga_0 x,y,z$ are at most $\theta'$, and similarly the
angles at $y'=\pi(\phi_\rho (\ga_0 w_\alpha))$ and $z'=\pi(\phi_\rho
(\ga_0\alpha u'))$ of the geodesic triangle with vertices $\ga_0
\alpha x,y',z'$ are at most $\theta'$.

Assume that $L'\geq \rho+1$, so that by the construction of $\ga_0$
$$
d(x,y)=d(x,\ga_0x)-\rho\geq L'-\rho\geq 1\;,
$$ 
and similarly $d(z',\ga_0 \alpha\ga_0 x)\geq 1$. By the triangle
inequality and the definition of $\rho$, we have $d(y,z)\geq
2(\rho-\rho_0)\geq 2$, and similarly $d(y',z')\geq 2$. Assume that $k$
is large enough, in particular so that we have
$$
d(z,y')= d(\ga_0 x,\ga_0\alpha x)-d(\ga_0 x,z)-d(y',\ga_0\alpha x) 
\geq n_k-2\rho\geq 1\;.
$$
Let $\epsilon'>0$ be small enough, and assume that $\theta'$ is small
enough, so that in particular $\theta'\leq \theta_0(1,\epsilon',
\frac{\theta}{2})$.  The piecewise geodesic path
$$
\omega=\mathopen{[}x,y\mathclose{]}\cup\mathopen{[}y,z\mathclose{]}
\cup \mathopen{[}z,y'\mathclose{]} \cup \mathopen{[}y',z'\mathclose{]}
\cup \mathopen{[}z',\ga_0\alpha\ga_0x\mathclose{]}
$$ 
has exterior angles at $y,z,y',z'$ at most $\theta'$.  Hence, by Lemma
\ref{lem:quasigeodpresquegeod}, the angles at $x$ between $u'$ and the
unit tangent vector $v_{\alpha,\,0}$ at $x$ pointing towards $\ga_0
\alpha \ga_0 x$, and the angle at $\ga_0\alpha\ga_0x$ between
$\ga_0\alpha u''$ and the unit tangent vector $w_{\alpha,\,0}$ at
$\ga_0\alpha\ga_0 x$ pointing towards $x$, are at most
$\frac{\theta}{2}$:
\begin{equation}\label{eq:proprianglecasdeu}
  \measuredangle_{x}(u',v_{\alpha,\,0})\leq \frac{\theta}{2},\;\;\;
\measuredangle_{\ga_0\alpha\ga_0 x}(\ga_0\alpha u'',w_{\alpha,\,0})
\leq \frac{\theta}{2}\;.
\end{equation}
Furthermore the piecewise geodesic path $\omega$ is at distance at
most $\epsilon'$ from $\mathopen{[}x,\ga_0\alpha \ga_0 x\mathclose{]}$.

Hence, as in Case 1, there exists a constant $c_5>0$ (depending only
on $\epsilon'$, the H\"older constants of $\wt F$ and
$\max_{\pi^{-1}(B(y,\,\epsilon')\,\cup \,B(z,\,\epsilon')\,\cup
  \,B(y',\,\epsilon')\,\cup \,B(z',\,\epsilon'))}\,|\wt F|$) such that
\begin{equation}\label{eq:controlintegralpotentuncasdeu}
  \Big|\int_{x}^{\ga_0\alpha\ga_0 x} \wt F\;-\int_{x}^{y} \wt F\;-
\int_{y}^{z} \wt F\;-\int_{z}^{y'} \wt F\;-\int_{y'}^{z'} \wt F\;-
\int_{z'}^{\ga_0\alpha\ga_0 x} \wt F\;\Big|\leq c_5\;.
\end{equation}
By invariance and additivity, we have
\begin{equation}\label{eq:controlintegralpotentdeucasdeu}
\int_{x}^{\alpha x} \wt F\;=\int_{\ga_0 x}^{\ga_0 \alpha x} \wt F\;=
\int_{\ga_0 x}^{z} \wt F\;+\int_{z}^{y'} \wt F\;+
\int_{y'}^{\ga_0 \alpha x} \wt F\;.
\end{equation}
Note that $d(x,y)\leq d(x,\ga_0x)$, $d(y,z)\leq 2\rho$ and
$d(y',z')\leq 2\rho$. Therefore, with $c_6=(d(x,\ga_0x)+2\rho)
\max_{\pi^{-1}(B(x,\,d(x,\ga_0x)+2\rho))} |\wt F|$, we have
$$
\Big|\int_{\ga_0 x}^{z} \wt F\;\Big|,\;
\Big|\int_{y'}^{\ga_0 \alpha x} \wt F\;\Big|,\;
\Big|\int_{x}^{y} \wt F\;\Big|,\;\Big|\int_{y}^{z} \wt F\;\Big|,\;
\Big|\int_{y'}^{z'} \wt F\;\Big|,\;
\Big|\int_{z'}^{\ga_0\alpha\ga_0 x} \wt F\;\Big|\leq c_6\;.
$$
Hence, by the formulas \eqref{eq:controlintegralpotentuncasdeu}
and \eqref{eq:controlintegralpotentdeucasdeu}, we have 
$$
\int_{x}^{\ga_0\alpha\ga_0 x} \wt F\;\geq 
\Big(\int_{x}^{\alpha x} \wt F \;\Big)-c_5-6c_6\;.
$$
Assertion (iii) above implies that
$$
\sum_{\alpha\in S_k} e^{\int_{x}^{\ga_0\alpha \ga_0x}\wt F}\; \geq
 k\;e^{-c_5-6c_6}\;e^{\delta' n_k}\;.
$$

Define 
$$
x_0=x\;\;\;{\rm and}\;\;\;v_0=u'\;.
$$ 
By Equation \eqref{eq:proprianglecasdeu}, for $k$ large enough and for
every $\alpha\in S_k$, we have 
$$
\measuredangle_{x_0}(v_0,v_{\alpha,\,0})
\leq \frac{\theta}{2}\leq\theta\;.
$$
By the spherical triangle inequality and since the isometry
$\ga_0\alpha$ preserves the angles for the first inequality, and by
the formulas \eqref{eq:angleuprimuprimprim} and
\eqref{eq:proprianglecasdeu} for the second one, we have
$$
\measuredangle_{\ga_0\alpha\ga_0 x_0}(-\ga_0\alpha\ga_0 u', w_{\alpha,\,0}) 
\leq \measuredangle_{\ga_0 x}(-\ga_0 u',u'')  
+\measuredangle_{\ga_0\alpha\ga_0 x}(\ga_0\alpha u'',w_{\alpha,\,0}) 
\leq \frac{\theta}{2}+\frac{\theta}{2}=\theta\;.
$$

For $k$ large enough and for all $\alpha,\beta\in S_k$, we have by
the triangle inequality and Assertion (i) above
$$
  n_k-2\,d(x,\ga_0 x)\leq d(x_0,\ga_0\alpha\ga_0 x_0)
\leq n_k+1+2\,d(x,\ga_0 x)\;,
$$
and, if $\alpha\neq \beta$, by the definition of $C'$ in Equation
\eqref{eq:defiCprim},
\begin{align*}
d(\ga_0\alpha\ga_0 x_0,\ga_0\beta\ga_0 x_0)&=
d(\alpha\ga_0 x,\beta\ga_0 x)\geq d(\alpha x,\beta x)
-d(\alpha x,\alpha \ga_0 x) - d(\beta x,\beta \ga_0 x) 
\\ &\geq C'-2\,d(x,\ga_0 x)\geq C\;.
\end{align*}
As in Case 1, there exists $N\in \NN\,\cap \mathopen{]}n_k-2\,d(x,\ga_0 x)-1,
n_k+2\,d(x,\ga_0 x)+2\mathclose{[}$ such that
\begin{align*}
 \sum_{\tiny\begin{array}{c}\alpha\in S_k\\ 
N\leq d(x_0,\ga_0\alpha\ga_0 x_0)< N+1\end{array}} 
e^{\int_{x_0}^{\ga_0\alpha\ga_0 x_0}\wt F}\; &\geq  
\frac{k\;e^{-c_5-6c_6}}{4\,d(x,\ga_0 x)+3}\;e^{\delta' n_k}
\\ &\geq  
\frac{k\;e^{-c_5-6c_6-\delta'(2\,d(x,\ga_0 x)+2)}}
{4\,d(x,\ga_0 x)+3}\;e^{\delta' N}\;.
\end{align*}
Given $L>0$, if $k$ is large enough, then 
$\frac{k\; e^{-c_5-6c_6-\delta'(2\,d(x,\ga_0 x)+2)}}
{4\,d(x,\ga_0 x)+3} \geq 1$ and $N\geq n_k-2\,d(x,\ga_0 x)-1> L$,
and defining
$$
S=\{\ga_0\alpha\ga_0\;:\; \alpha\in S_k\;\; {\rm and}\; \; 
N\leq d(x_0,\ga_0\alpha\ga_0 x_0)< N+1\}\;,
$$
the objects $x_0,v_0,N,S$ defined above satisfy the properties (1)-(4)
required in Proposition \ref{prop:constructsemigroup}.

\medskip 
\noindent{\bf Case 3. }  Assume  that $\measuredangle_x(v,w)<
\frac{\theta}{4}$ and that $\measuredangle_x(u,-v)<
\frac{\theta}{2}$.

Consider the element $\overline{\ga}_0$ of $\Ga$ and the unit tangent
vectors $\overline{u}'$ and $\overline{u}''$ defined above.

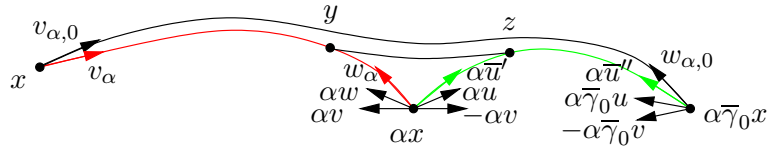
\begin{center}
\input{fig_HamPaucase3.pstex_t}
\end{center}

Let $k\in\NN$ and $\alpha\in S_k$. By the spherical triangle
inequality, by the assumptions of Case 3, by Assertion (ii) above, by
the properties of $\overline{u}'$, and since
$\theta''\leq \frac{\theta}{8}$ and $\theta\leq \frac{\pi}{2}$, we
have
\begin{align*}
\measuredangle_{\alpha x}(w_\alpha,\alpha\overline{u}')&\geq
\measuredangle_{\alpha x}(\alpha v,- \alpha v)\\ &\;\;\;\;
-\measuredangle_{\alpha x}(\alpha v,\alpha w)
-\measuredangle_{\alpha x}(\alpha w,w_\alpha)
-\measuredangle_{\alpha x}(- \alpha v,\alpha u)
-\measuredangle_{\alpha x}(\alpha u,\alpha\overline{u}')\\ & \geq
\pi-\frac{\theta}{4}-\theta''-\frac{\theta}{2}-\frac{\theta}{8}
\geq \pi-\theta\geq \frac{\pi}{2}\geq \frac{\theta}{4}\;.
\end{align*}

As in Case 2, if $L'$ and $k$ are large enough, for every $\epsilon'$
small enough, if $\theta'$ is small enough, so that in particular
$L'\geq \rho+1$, $n_k\geq \rho+1$, $\theta'\leq \theta_0(1,\epsilon',
\frac{\theta}{4})$, with $y=\pi(\phi_\rho w_\alpha)$ and $z=
\pi(\phi_\rho (\alpha \overline{u}'))$, considering the piecewise
geodesic path
$$
w=\mathopen{[}x,y\mathclose{]}\cup\mathopen{[}y,z\mathclose{]} \cup 
\mathopen{[}z,\alpha\overline{\ga}_0x\mathclose{]}\;,
$$
we have for every $\alpha\in S_k$:

$\bullet$~ the angles at $y$ and $z$ of the geodesic triangle with
vertices $\alpha x,y,z$ are at most $\theta'$, by the definition of
$\rho$, hence the exterior angles of $\omega$ at $y,z$ are at most
$\theta'\leq \theta_0(1,\epsilon', \frac{\theta}{4})$;

$\bullet$~ we have 
$$
d(x,y)\geq n_k-\rho\geq 1\;,\;\;d(y,z)\geq
2(\rho-\rho_0)\geq 2\;,$$
$$
d(z,\alpha\overline{\ga}_0x)= d(\alpha x,
\alpha\overline{\ga}_0x)-\rho\geq L'-\rho\geq 1\;.
$$ 
hence we may apply Lemma \ref{lem:quasigeodpresquegeod} to $\omega$
with $\ell=1$;

$\bullet$~ the angles at $x$ between $v_\alpha$ and the unit tangent
vector $v_{\alpha,\,0}$ at $x$ pointing towards $\alpha \overline{\ga}_0
x$, and the angle at $\alpha \overline{\ga}_0 x$ between
$\alpha\overline{u}''$ and the unit tangent vector $w_{\alpha,\,0}$ at
$\alpha \overline{\ga}_0 x$ pointing towards $x$, are at most
$\frac{\theta}{4}$;

$\bullet$~ the path $\omega$ is at distance at most $\epsilon'$ from
$\mathopen{[}x,\alpha \overline{\ga}_0 x\mathclose{]}$;

$\bullet$~ there exists a constant $c_7>0$ (depending only
on $\epsilon'$, the H\"older constants of $\wt F$ and
$\max_{\pi^{-1}(B(y,\,\epsilon')\,\cup \,B(z,\,\epsilon'))}\,|\wt F|$) such that
$$
  \Big|\int_{x}^{\alpha\overline{\ga}_0 x} \wt F-\int_{x}^{y} \wt F\;-
\int_{y}^{z} \wt F\;-\int_{z}^{\alpha\overline{\ga}_0 x} \wt F\;\Big|
\leq c_7\;;
$$

$\bullet$~ if $c_8=(d(x,\overline{\ga}_0x)+2\rho)
\max_{\pi^{-1}(B(x,\,d(x,\overline{\ga}_0x)+2\rho))} |\wt F|$, we have
$$
\Big|\int_{y}^{\alpha x} \wt F\;\Big|,\;
\Big|\int_{y}^{z} \wt F\;\Big|,\;
\Big|\int_{z}^{\alpha\overline{\ga}_0 x} \wt F\;\Big|\leq c_8\;.
$$

\medskip
Since $\int_{x}^{\alpha x} \wt F=\int_{x}^{y} \wt F+\int_{y}^{\alpha
  x} \wt F$, and by the last two points, we have 
$$
\int_{x}^{\alpha\overline{\ga}_0 x} \wt F\;\geq 
\Big(\int_{x}^{\alpha x} \wt F \;\Big)-c_7-3\,c_8\;.
$$
Assertion (iii) above implies that
$$
\sum_{\alpha\in S_k} e^{\int_{x}^{\alpha\overline{\ga}_0 x} \wt F}\; \geq
 k\;e^{-c_7-3\,c_8}\;e^{\delta' n_k}\;.
$$

Define 
$$
x_0=x\;\;\;{\rm and}\;\;\;v_0=v\;.
$$ 
By Assertion (ii) above and the third point above, we have
$$
\measuredangle_{x_0}(v_0,v_{\alpha,\,0})
\leq \measuredangle_{x}(v,v_{\alpha})+
\measuredangle_{x}(v_{\alpha},v_{\alpha,\,0})\leq
\theta''+\frac{\theta}{4}\leq\theta\;.
$$
By the assumptions of Case 3, by the properties of $\overline{u}''$,
and by the third point above, we have
\begin{align*}
\measuredangle_{\alpha\overline{\ga}_0 x}
(-\alpha\overline{\ga}_0 v_0, w_{\alpha,\,0}) &
\leq \measuredangle_{\alpha\overline{\ga}_0 x}
(-\alpha\overline{\ga}_0 v_0,\alpha\overline{\ga}_0 u)
+\measuredangle_{\alpha\overline{\ga}_0 x}
(\alpha\overline{\ga}_0 u,\alpha\overline{u}'')
+\measuredangle_{\alpha\overline{\ga}_0 x}
(\alpha\overline{u}'', w_{\alpha,\,0}) \\ &
\leq \frac{\theta}{2}+\frac{\theta}{8}+\frac{\theta}{4}\leq\theta\;.
\end{align*}

For all $\alpha,\beta\in S_k$, by the triangle inequality and
Assertion (i) above, we have
$$
n_k-d(x,\overline{\ga}_0 x)\leq d(x_0,\alpha\overline{\ga}_0 x_0)
\leq n_k+1+d(x,\overline{\ga}_0 x)\;,
$$
and, if $\alpha\neq \beta$, by the definition of $C'$ in Equation
\eqref{eq:defiCprim},
$$
d(\alpha\overline{\ga}_0 x_0,\beta\overline{\ga}_0 x_0) \geq 
d(\alpha x,\beta x)
-2\,d(x,\overline{\ga}_0 x) \geq C'-2\,d(x,\overline{\ga}_0 x)\geq C\;.
$$
As in Case 1, there exists $N\in \NN\,\cap 
\mathopen{]}n_k-d(x,\overline{\ga}_0 x)-1, 
n_k+d(x,\overline{\ga}_0 x)+2\mathclose{[}$
such that
$$
 \sum_{\tiny\begin{array}{c}\alpha\in S_k\\ 
N\leq d(x_0,\alpha\overline{\ga}_0 x_0)< N+1\end{array}} 
e^{\int_{x_0}^{\alpha\overline{\ga}_0 x_0}\wt F}\; 
\geq  
\frac{k\;e^{-c_7-3c_8-\delta'(d(x,\overline{\ga}_0 x)+2)}}
{2\,d(x,\overline{\ga}_0 x)+3}\;e^{\delta' N}\;.
$$
Given $L>0$, if $k$ is large enough, then $\frac{k\;e^{-c_7-3c_8- \delta'
    (d(x,\overline{\ga}_0 x)+2)}}{2\,d(x,\overline{\ga}_0 x)+3}
\geq 1$ and $N\geq n_k-d(x,\overline{\ga}_0 x)-1> L$, and defining
$$
S=\{\alpha\overline{\ga}_0\;:\; \alpha\in S_k\;\; {\rm and}\; \; 
N\leq d(x_0,\alpha\overline{\ga}_0 x_0)< N+1\}\;,
$$
the objects $x_0,v_0,N,S$ defined above satisfy the properties (1)-(4)
required in Proposition \ref{prop:constructsemigroup}.

\medskip 
\noindent{\bf Case 4. }  The remaining case, when
$\measuredangle_x(v,w) < \frac{\theta}{4}$ and $\measuredangle_x(u,w)<
\frac{\theta}{2}$, follows similarly to Case 3: we take $x_0=x, v_0=w$
and $S$ the appropriate subset of $\{\overline{\ga}_0\alpha\;:\;
\alpha\in S_k\}$.

\medskip
This ends the proof of Proposition \ref{prop:constructsemigroup} under
the assumption that $\delta'>0$. Up to adding the same big enough
constant to $\wt F$ and to $\delta'$, the conclusion of each of the
four steps above (as $k$ was allowed to be taken large enough) proves
that Proposition \ref{prop:constructsemigroup} remains valid when
$\delta'\leq 0$.  \cqfd

\bigskip
The main objective of this subsection is to show the following result.

\btheo\label{theo:critexpsubsemigroup} Let $\wt M$ be a complete
simply connected Riemannian manifold, with dimension at least $2$ and
pinched negative sectional curvature at most $-1$. Let $\Ga$ be a
non-elementary discrete group of isometries of $\wt M$. Let $\wt F
:T^1\wt M\ra \RR$ be a H\"older-continuous $\Ga$-invariant map.

Then the critical exponent $\delta_{\Ga,\,F}$ is equal to the upper bound
of the critical exponents $\delta_{G,\,F}$, where $G$ ranges over the
convex-cocompact free subsemigroups of $\Ga$.  
\etheo

When $F=0$, this result (with $\wt M$ any proper geodesic
CAT$(-1)$-space) is due to Mercat \cite{Mercat12}.

\medskip \dem It is immediate that $\delta_{\Ga,\,F}\geq
\delta_{G,\,F}$ for any subsemigroup $G$ of $\Ga$. Hence
$\delta_{\Ga,\,F}$ is at least the upper bound of the critical
exponents $\delta_{G,\,F}$ of the convex-cocompact free subsemigroups
$G$ of $\Ga$. We now prove the converse inequality. Up to adding a
constant to $F$, we may assume that $\delta_{\Ga,\,F}>0$.

Let $s\in\mathopen{]}0,\delta_{\Ga,\,F}\mathclose{[}$ and
  $\delta'\in\mathopen{]}s,\delta_{\Ga,\,F}\mathclose{[}\,$.  Let
  $\epsilon \in\mathopen{]}0,1\mathclose{]}$ be small enough. Let
$\theta>0$ (depending only on $\epsilon$) be given by Proposition
\ref{prop:schottkysemigroups}. Let $x_0\in \wt M$ and $v_0\in T^1\wt
M$ be given by Proposition \ref{prop:constructsemigroup} (depending on
$\delta',\theta,C$).

Since two geodesic segments of lengths at least $1$ which are close
enough have their unit tangent vectors close (see Lemma
\ref{lem:technicholder}), there exists $c>0$ (depending only on
$\epsilon$, $x_0$, the H\"older constants of $\wt F$ and on
$\max_{\pi^{-1}(B(x_0,1))} \,|\wt F-s|$) such that for all
$x,y,z\in\Ga x_0$ and for every piecewise geodesic segment
$\mathopen{[}x,y\mathclose{]}\cup \mathopen{[}y,z\mathclose{]}$, with
segments of length at least $1$, which is contained in the
$\epsilon$-neighbourhood of $\mathopen{[}x,z\mathclose{]}$, we have
$$
\big|\int_x^z(\wt F-s)\;-\int_x^y(\wt F-s)\;-\int_y^z(\wt
F-s)\;\big|\leq c\;.
$$

Let $L=\max\{5,\frac{c+s}{\delta'-s}\}$.  Let $N>L$ and $S$ a
finite subset of $\Ga$ be given by Proposition
\ref{prop:constructsemigroup} (depending on $\delta',\theta,L,C$). By
Proposition \ref{prop:schottkysemigroups}, whose assumptions are
satisfied by the first three assertions of Proposition
\ref{prop:constructsemigroup}, the semigroup $G$ generated by $S$ is
free on $S$ and convex-cocompact, and for all $\ga\in G$ and
$\alpha\in S$, the piecewise geodesic segment $\mathopen{[}x_0,\ga
x_0\mathclose{]}\cup \mathopen{[}\ga x_0,\ga\alpha x_0\mathclose{]}$,
whose segments have length at least $1$, is contained in the
$\epsilon$-neighbourhood of $\mathopen{[}x_0,\ga\alpha x_0
\mathclose{]}$. In particular, by the definition of $c$ and by
invariance, we have
\begin{equation}\label{eq:carmen}
\Big|\int_{x_0}^{\ga\alpha x_0}(\wt F-s)\;-\int_{x_0}^{\ga x_0}(\wt F-s)\;-
\int_{x_0}^{\alpha x_0}(\wt F-s)\;\Big|\leq c\;.
\end{equation}

For every $k\in\NN$, denote by $S_k$ the set of words of length $k$ in
elements of $S$. Let
$$
\Sigma_k=\sum_{\ga\in S_k} e^{\int_{x_0}^{\ga x_0}(\wt F-s)}\;,
$$
so that, since $G$ is the free semigroup generated by $S$, we have
$$
Q_{G,\,F,\,x_0,\,x_0}(s)=\sum_{k\in\NN}\Sigma_k\;.
$$
By Assertion (1) of Proposition \ref{prop:constructsemigroup}, for
every $\alpha\in S$, we have $-s d(x_0,\alpha x_0)\geq -s(N+1)$.  By
Equation \eqref{eq:carmen}, by Assertion (4) of Proposition
\ref{prop:constructsemigroup}, and since $N>L\geq
\frac{c+s}{\delta'-s}$, we hence have
\begin{align*}
  \Sigma_{k+1}=\sum_{\ga\in S_k,\;\alpha\in S}
  e^{\int_{x_0}^{\ga\alpha x_0}(\wt F-s)} &\geq e^{-c}\sum_{\ga\in
    S_k,\;\alpha\in S} e^{\int_{x_0}^{\ga x_0}(\wt F-s)+
    \int_{x_0}^{\alpha x_0}(\wt F-s)}\\
  &=e^{-c}\Sigma_k\sum_{\alpha\in S} e^{\int_{x_0}^{\alpha x_0}(\wt
    F-s)}\geq e^{-c-s}\Sigma_k\;e^{(\delta'-s)N} \geq \Sigma_k\;.
\end{align*}
Therefore the Poincar\'e series $Q_{G,\,F,\,x_0,\,x_0}(s)$ diverges. Hence
$\delta_{G,\,F}\geq s$.  As $s<\delta_{\Ga,\,F}$ was arbitrary,
this proves Theorem \ref{theo:critexpsubsemigroup}.
\cqfd

\section{A Hopf-Tsuji-Sullivan-Roblin theorem for 
Gibbs states and applications}
\label{sec:GHTSappli}

Let $\wt M,\Ga,\wt F$ be as in the beginning of Chapter
\ref{sec:negacurvnot}: $\wt M$ is a complete simply connected
Riemannian manifold, with dimension at least $2$ and pinched sectional
curvature at most $-1$; $\Ga$ is a non-elementary discrete group of
isometries of $\wt M$; and $\wt F :T^1\wt M\ra \RR$ is a
H\"older-continuous $\Ga$-invariant map. The aim of this chapter is to
give criteria for the ergodicity and non ergodicity of the geodesic
flow on the space $T^1M=\Ga\backslash T^1\wt M$ endowed with a Gibbs
measure. We also discuss several applications, in particular to the
uniqueness of a Gibbs measure on $T^1M$ as soon as there exists a
finite such measure.

\subsection{Some geometric notation}
\label{subsec:geomnotation}

We introduce in this subsection the notation for several geometric
objects that will be used in the proofs of the subsections
\ref{subsec:GHTSR} and \ref{subsec:productconv}. 

Though in this Section \ref{sec:GHTSappli}, this notation is not
needed for the statement of the main results \ref{theo:critnonergo}
and \ref{theo:critergo}, hence could logically be put at the
beginning of their proofs, we prefer to make a separate initial
subsection, precisely to facilitate the references to it in Section
\ref{sec:growth}. The reader may skip to the beginning of Subsection
\ref{subsec:GHTSR} and come back here when needed.

\medskip For every $z\in \wt M$, every subset $Z$ of $\partial_\infty
\wt M$ and every $r>0$, let us define the {\it
  $r$-thickened}\index{cone!$r$-thickened} and {\it
  $r$-thinned}\index{cone!$r$-thinned} cones over $Z$ with vertex $z$
as
$$
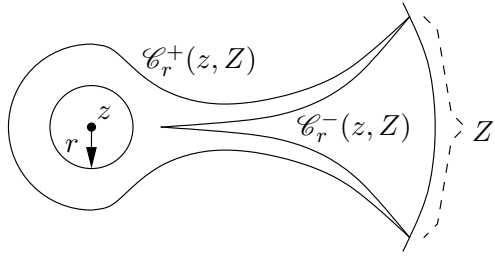
\gls{conethick}=\N_r\Big(\bigcup_{z'\in B(z,\,r)}\C_{z'}Z\Big)
\;\;\;{\rm and}\;\;\;
\gls{conethin}=\N^-_r\Big(\bigcap_{z'\in B(z,\,r)}\C_{z'}Z\Big)\;.
$$

\begin{center}
\input{fig_CrpmzZ.pstex_t}
\end{center}

For all $r>0$ and $z\in\wt M$, we will say that $v\in T^1\wt M$ is
{\it $r$-close to}\index{close@$r$-close to} $z$ if there exists
$s\in\mathopen{]}-\frac{r}{2}, \frac{r}{2}\mathclose{[}$ such that
$\pi(\phi_sv)$ is the closest point to $z$ on the geodesic line
defined by $v$, and $d(z,\pi(\phi_sv))<r$. For every $\ga\in\Ga$, the
vector $\ga v$ is $r$-close to $\ga z$ if and only if $v$ is $r$-close
to $z$.

For every subset $Z$ of $\partial_\infty \wt M$, let $\gls{Kplus}$ be
the set of $v\in T^1\wt M$ such that $v_+\in Z$ and $v$ is $r$-close
to $z$, and $\gls{Kmoins}=\iota K^+(z,r,Z)$. Note that $\ga
K^\pm(z,r,Z)= K^\pm(\ga z,r,\ga Z)$ for every $\ga\in\Ga$. Also denote
by
$$
K(z,r)=K^+(z,r,\partial_\infty\wt M)=
K^-(z,r,\partial_\infty\wt M)
$$
the relatively compact open set of vectors $v$ that are $r$-close to
$z$, which satisfies $\ga K(z,r)=K(\ga z,r)$ for every isometry $\ga$
of $\wt M$.

\medskip
\noindent\begin{minipage}{5.5cm}
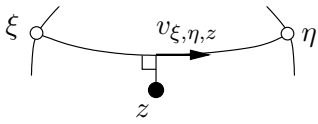
\begin{center}
\input{fig_vxietaz.pstex_t}
\end{center}
\end{minipage}
\begin{minipage}{9.4cm}
  ~~~ For distinct points $\xi$ and $\eta$ in $\partial_\infty
  \wt M$, let $v_{\xi,\,\eta,\,z}$ be the unit tangent vector based at the
  closest point to $z$ on the geodesic line between $\xi$ and $\eta$
  pointing towards $\eta$. We have $\ga v_{\xi,\,\eta,\,z}
  =v_{\ga \xi,\,\ga\eta,\,\ga z}$ for every isometry $\ga$ of $\wt M$.
\end{minipage}

\medskip\noindent\begin{minipage}{8.9cm} ~~~ For all $r>0$ and
  $z,w\in\wt M$ such that $d(z,w)>2r$, let us define $\gls{Lset}$ as the set
  of $(\xi,\eta)$ in $\partial_\infty \wt M\times \partial_\infty \wt
  M$ such that the geodesic line starting from $\xi$ and ending at
  $\eta$ meets first $B(z,r)$ and then $B(w,r)$. We have $\ga
  L_r(z,w)=L_r(\ga z,\ga w)$ for every isometry $\ga$ of $\wt M$.
\end{minipage}
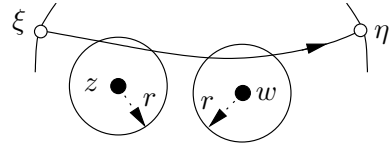
\begin{minipage}{6cm}
\begin{center}
\input{fig_ellrzw.pstex_t}
\end{center}
\end{minipage}

\medskip
For all $r>0$ and $z,w\in\wt M$, let us define
$$
\OOO^+_r(z,w)=\bigcup_{z'\in B(z,\,r)} \OOO_{z'}B(w,r)\;\;\;{\rm
  and}\;\;\; \OOO^-_r(z,w)=\bigcap_{z'\in B(z,\,r)} \OOO_{z'}B(w,r)\;.
$$
For all $r>0$, $w\in\wt M$ and $z\in \partial_\infty \wt M$, let us define 
$$
\OOO^+_r(z,w)=\OOO^-_r(z,w)=
\OOO_{z}B(w,r)\;.
$$  
Note that, for every fixed $w\in \wt M$, the closure of
$\OOO^\pm_r(z,w)$ converges to the closure of $\OOO^\pm_r(z_0,w)$ in
the space of closed subsets of $\wt M\cup\partial_\infty\wt M$, as $z$
tends to $z_0$ in $\wt M\cup\partial_\infty\wt M$, uniformly on $r$ in
a bounded interval. We have $\ga \OOO^\pm_r(z,w)=\OOO^\pm_r(\ga z,\ga
w)$ for every isometry $\ga$ of $\wt M$.

\blemm \label{lem:LrOOOpm} For all $r>0$ and $z,w\in\wt M$ such that
$d(z,w)>2r$, the set $L_r(z,w)$ contains $\OOO^-_{r}(w,z)\times
\OOO^-_{r}(z,w)$.  
\elemm

\dem Let $(\xi,\eta)\in \OOO^-_{r}(w,z)\times \OOO^-_{r}(z,w)$ and let
us prove that $(\xi,\eta)\in L_r(z,w)$.

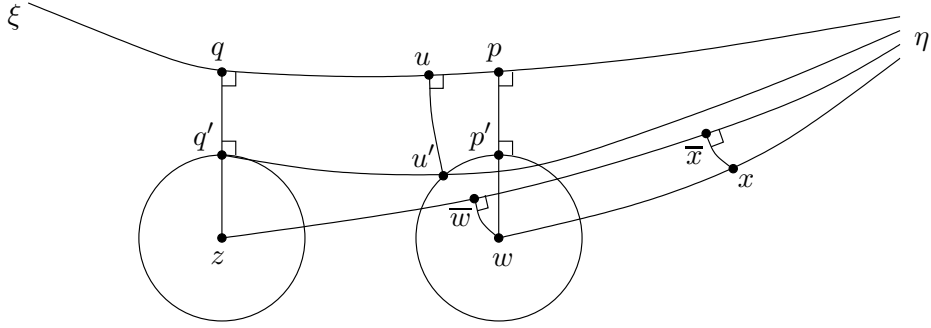
\begin{center}
\input{fig_LrOOOpm.pstex_t}
\end{center}

For every $x$ on the geodesic ray $\mathopen{[}w,\eta\mathclose{[}\,$,
denote by $\ov x$ its closest point on the geodesic ray
$\mathopen{[}z,\eta\mathclose{[}\,$.  Since $\eta\in
\OOO^-_{r}(z,w)\subset \OOO_{z}B(w,r)$, the point $\ov w$ is at
distance at most $r$ from $w$. Since closest point maps do not
increase the distances, by convexity, by the triangle inequality and
since $d(z,w)>2r$, we have
$$
d(z,x)\geq d(z,\ov x)\geq d(z,\ov w)\geq d(z,w)-d(w,\ov w)>2r-r=r\;.
$$ Hence $\eta\notin \OOO_{w}B(z,r)\supset \OOO^-_{r}(w,z)$. This
implies that $\xi\neq \eta$, hence the geodesic line $\mathopen{]}\xi,
\eta\mathclose{[}$ between $\xi$ and $\eta$ is well defined.

Let $p$ and $q$ be the closest points to $w$ and $z$ respectively on
the geodesic line $\mathopen{]}\xi, \eta\mathclose{[}$. We need to
  prove that $\xi, q,p,\eta$ are in this order on $\mathopen{]}\xi,
\eta\mathclose{[}$, and that $d(w,p),d(z,q)\leq r$.

Let $p'$ and $q'$ be the points on the geodesic segments
$\mathopen{[}w,p\mathclose{]}$ and $\mathopen{[}z, q\mathclose{]}$ at
distance $\min\{r,d(w,p)\}$ and $\min\{r,d(z,q)\}$ from $w$ and $z$,
respectively, and let us prove that $p'$ and $q'$ are actually equal
to $p$ and $q$, respectively (see the above picture, where for a
contradiction we assume that neither equality holds). Note that $p'$
and $q'$ are the closest points to $p$ and $q$ on the balls $B(w,r)$
and $B(z,r)$ respectively.

Up to permuting $p$ and $q$ as well as $\xi$ and $\eta$, we may assume
that $d(p,w)\geq d(q, z)$.  Since $\eta\in \OOO^-_{r}(z,w)$, the
geodesic ray $\mathopen{[}q', \eta\mathclose{[}$ meets $B(w,r)$ at at
least one point $u'$.  Let $u$ be the closest point to $u'$ on the
geodesic line between $\xi$ and $\eta$. Since $d(z,w)>2r$, we have
$u'\neq q'$. Hence, by convexity,
$$
d(q,q')\geq d(u,u')\geq d(p,p')\;
$$ (the second inequality holds since $d(p',p)$ is by convexity the
smallest distance between a point in the ball $B(w,r)$ and a point on
the geodesic line $\mathopen{]}\xi, \eta\mathclose{[}\,$). 

Assume for a contradiction that $q\neq q'$: the first inequality is
then strict, and this contradicts the fact that $d(p,w)=r+d(p,p')\geq
d(q, z)= r+ d(q,q')$. Hence $q=q'$, therefore $p=p'$ and the points
$z,w$ are at distance at most $r$ from the geodesic line
$\mathopen{]}\xi, \eta\mathclose{[}\,$.

Since $d(q,z)\leq r$ and $\eta\in \OOO^-_{r}(z,w)$, the geodesic ray
$\mathopen{[}q,\eta\mathclose{[}$ meets $B(w,r)$. Since $d(z,w)>2r$,
the point $q$ does not belong to the intersection of $B(w,r)$ and
$\mathopen{[}q,\eta\mathclose{[}\,$. Hence by convexity, $p$ belongs
to this intersection, and $\xi, q,p,\eta$ are in this order on the
geodesic line $\mathopen{]}\xi, \eta\mathclose{[}\,$. Therefore
$(\xi,\eta)\in L_r(z,w)$, as required.  \cqfd

\medskip
The following result that we state without proof is a simple extension
of Mohsen's shadow lemma \ref{lem:shadowlemma}, with a similar proof
(see \cite[p.~10]{Roblin03}).

\blemm \label{lem:thickshadowlemma} Let $\sigma\in\RR$, let
$(\mu_{x})_{x\in\wt M}$ be a Patterson density of dimension $\sigma$
for $(\Ga,F)$, and let $x,y\in\wt M$.  If $R$ is large enough, there
exists $C>0$ such that for every $\ga\in\Ga$,
$$
\frac{1}{C}\;e^{\int_x^{\ga y} (\wt F-\sigma)}\leq
\mu_{x}\big(\OOO^-_xB(\ga y,R)\big)\leq
\mu_{x}\big(\OOO^+_xB(\ga y,R)\big)\leq 
C\;e^{\int_x^{\ga y} (\wt F-\sigma)}\;.
$$
\elemm

\subsection{The Hopf-Tsuji-Sullivan-Roblin theorem for Gibbs 
states}
\label{subsec:GHTSR} 

Fix $\sigma\in\RR$. Let $(\mu^\iota_{x})_{x\in\wt M}$ and
$(\mu_{x})_{x\in\wt M}$ be Patterson densities of the same dimension
$\sigma$ for $(\Ga,F\circ \iota)$ and $(\Ga,F)$ respectively, and let
$\wt m$ (respectively $m$) be their associated Gibbs measure on
$T^1\wt M$ (respectively $T^1M$), see Subsection
\ref{subsec:GibbsSulivanmeasure}.

In this subsection, we give criteria for the ergodicity and non
ergodicity of the geodesic flow $(\phi_t)_{t\in\RR}$ on
$T^1M=\Ga\backslash T^1\wt M$ endowed with Gibbs measures.  These
criteria generalise the classical Hopf-Tsuji-Sullivan-Roblin results
for the Bowen-Margulis measures (that is, when $F=0$). Their proofs,
given at the end of this subsection, follow closely the general lines
of those presented in \cite[Chap.~1E]{Roblin03}.

\medskip Let us recall some definitions on a dynamical system
$(\Omega,G,\nu)$, where $\Omega$ is a Hausdorff locally compact
$\sigma$-compact topological space, $G$ is a locally compact
$\sigma$-compact topological group acting continuously on $\Omega$,
and $\nu$ is a positive regular (hence $\sigma$-finite) Borel measure
on $\Omega$, which is quasi-invariant under $G$.

Recall that $(\Omega,G,\nu)$ is {\it ergodic}%
\index{dynamical system (measurable)!ergodic}\index{ergodic} if for
any $G$-invariant measurable subset $A$ of $\Omega$, either $A$ or
$^cA$ has zero measure for $\nu$.

Recall that $\mathbbm{1}_A$ denotes the characteristic function of a
subset $A$. A (measurably) {\it wandering set}%
\index{wandering  set!(measurably)}%
\index{dynamical system (measurable)!wandering set of} of
$(\Omega,G,\nu)$ is a Borel subset $E$ of $\Omega$ such that for
$\nu$-almost every $\omega\in E$, the integral
$\int_G\mathbbm{1}_E(g\omega)\;dg$ is finite, where $dg$ is a fixed
left Haar measure on $G$. The {\it dissipative part}\index{dissipative
  part}\index{dynamical system (measurable)!dissipative part of} of
$(\Omega,G,\nu)$ (well defined up to a zero measure subset) is any
$G$-invariant Borel subset which contains every wandering set up to a
zero measure subset, and does not contain any non-wandering set of
positive measure. Note that any countable union of wandering sets is
contained in the dissipative part up to a zero measure subset. The
{\it conservative part}\index{conservative part}%
\index{dynamical system (measurable)!conservative part of} of
$(\Omega,G,\nu)$ is the complement of its dissipative part (and is
also defined up to a zero measure subset).  Recall that
$(\Omega,G,\nu)$ is {\it conservative}\index{conservative}%
\index{dynamical system (measurable)!conservative} (respectively {\it
  completely dissipative})\index{completely dissipative}%
\index{dynamical system (measurable)!completely dissipative} if its
dissipative part (respectively conservative part) has zero measure.

\medskip The first result gives criteria for the non-ergodicity
of the geodesic flow.

\btheo \label{theo:critnonergo}%
\index{theorem@Theorem!of Hopf-Tsuji-Sullivan-Roblin with potentials!non-ergodic}\index{Hopf-Tsuji-Sullivan-Roblin theorem with \\ potentials!non-ergodic} The
following assertions are equivalent.
\begin{enumerate}
\item[(1)] The Poincar\'e series of $(\Ga,F)$ at the point
  $\sigma$ converges: $Q_{\Ga,\, F}(\sigma) < + \infty$.
\item[(2)] There exists $x\in \wt M$ such that $\mu_x(\Lambda_c
  \Ga)=0$.
\item[(3)] There exists $x\in \wt M$ such that $\mu^\iota_x(\Lambda_c
  \Ga)=0$.
\item[(4)] There exists $x\in \wt M$ such that the dynamical system
  $(\partial_\infty^2 \wt M, \Ga, (\mu^\iota_x\otimes
  \mu_x)_{\mid \partial_\infty^2 \wt M})$ is non-ergodic and
  completely dissipative.
\item[(5)] The dynamical system $(T^1M, (\phi_t)_{t\in\RR}, m)$ is
  non-ergodic and completely dissipative.
\end{enumerate}
\etheo

We may replace Assertion (2) by

\smallskip
{\em (2') For every $x\in \wt M$, we have $\mu_x(\Lambda_c
  \Ga)=0$.}

\smallskip
\noindent and Assertion  (3) by

\smallskip
{\em (3') For every $x\in \wt M$, we have $\mu^\iota_x(\Lambda_c
  \Ga)=0$.}

\smallskip
\noindent since two elements in a Patterson density are pairwise
absolutely continuous.  By equivalence of measures, we may also replace
``There exists $x\in \wt M$ such that'' by ``For every $x\in \wt M$,''
in Assertion (4).

\medskip \rem This result indicates that the case when $Q_{\Ga,\,
  F}(\sigma) < + \infty$ (and in particular the case when $(\Ga,F)$ is
not of divergence type and $\sigma=\delta_{\Ga,\,F}$) is not
interesting from the dynamical system viewpoint. It is also not the
one that occurs the most frequently, as we will see later
on. Furthermore, when $Q_{\Ga,\, F}(\sigma) < + \infty$, uninteresting
Patterson densities of dimension $\sigma$ for $(\Ga,F)$ exist if
$\Lambda\Ga\neq\partial_\infty\wt M$, as follows. Fix $x_0$ in $\wt M$
and $\xi\in\partial_\infty\wt M-\Lambda\Ga$.  Denote by $\D_\eta$ the
(unit) Dirac mass at a point $\eta\in\partial_\infty\wt M$. Then for
every $x\in X$, the series
$$
\nu_x^\xi= \sum_{\ga\in\Ga} e^{-C_{F-\sigma,\,\ga\xi}(x,\ga x_0)}
  \;\D_{\ga\xi}
$$
converges. Indeed, its convergence is independent of $x$ and $x_0$ by
the cocycle property of $C_{F-\sigma}$ (Equation
\eqref{eq:cocycleprop}) and by Lemma \ref{lem:holderconseq}
(1). Hence, for every fixed $r>0$, we may assume that
$x\in\C\Lambda\Ga$ and that $x_0$ is close enough to $\xi$ so that
$\xi\in\OOO_{\ga^{-1}x}B(x_0,r)$ for every $\ga\in\Ga$. Then by the
invariance property of $C_{F-\sigma}$ (Equation
\eqref{eq:equivcocprop}) and by Lemma \ref{lem:holderconseq} (2),
there exists $c>0$ such that for every $\ga\in\Ga$, we have
$$
-C_{F-\sigma,\,\ga\xi}(x,\ga x_0)=-C_{F-\sigma,\,\xi}(\ga^{-1}x,x_0)\leq 
c+\int_{\ga^{-1}x}^{x_0}(\wt F-\sigma)=c+\int_{x}^{\ga x_0}(\wt F-\sigma)\;.
$$
Hence the convergence (in norm) of the series $\nu_x^\xi$ follows from
the convergence of the series $Q_{\Ga,\, F,\,x,\,x_0}(\sigma)$. It is
easy to check that $(\nu_x^\xi)_{x\in\wt M}$ is a Patterson density of
dimension $\sigma$ for $(\Ga,F)$. Taking linear combinations of them
as $\xi$ varies and weak-star limits gives other examples. The
measures $\nu_x^\xi$ are atomic and ergodic (but not properly ergodic,
indeed of type $I_\infty$ in the Hopf-Krieger classification, see for
instance \cite{FelMoo77}: they give full measure to an infinite orbit
of $\Ga$).

\medskip
The next result gives criteria for the ergodicity of the
geodesic flow.

\btheo\label{theo:critergo}%
\index{theorem@Theorem!of Hopf-Tsuji-Sullivan-Roblin with potentials!ergodic}\index{Hopf-Tsuji-Sullivan-Roblin theorem with \\ potentials!ergodic} The following assertions are equivalent.
\begin{enumerate}
\item[(i)]  The Poincar\'e series of $(\Ga,F)$ at the point
  $\sigma$ diverges: $Q_{\Ga,\, F}(\sigma) = + \infty$.
\item[(ii)] There exists $x\in \wt M$ such that
  $\mu_x(\,^c\Lambda_c \Ga)=0$.
\item[(iii)] There exists $x\in \wt M$ such that
  $\mu^\iota_x(\,^c\Lambda_c \Ga)=0$.
\item[(iv)] There exists $x\in \wt M$ such that the dynamical system
  $(\partial_\infty^2 \wt M, \Ga, (\mu^\iota_x\otimes
  \mu_x)_{\mid \partial_\infty^2 \wt M})$ is ergodic and 
  conservative.
\item[(v)] The dynamical system $(T^1M, (\phi_t)_{t\in\RR}, m)$ is
  ergodic and conservative. 
\end{enumerate}  
\etheo

We may also replace Assertion (ii) by

\smallskip
{\em (ii') For every $x\in \wt M$, we have $\mu_x(\,^c\Lambda_c
  \Ga)=0$.}

\smallskip
\noindent and Assertion (iii) by

\smallskip
{\em (iii') For every $x\in \wt M$, we have $\mu^\iota_x(\,^c\Lambda_c
  \Ga)=0$.}

\smallskip
\noindent By equivalence of measures, we may also replace
``There exists $x\in \wt M$ such that'' by ``For every $x\in \wt M$,''
in Assertion (iv).

\medskip
The proof of these two theorems will be based on the following list of
implications, as in the case $F=0$ in \cite[Chap.~1E]{Roblin03}. Let
us fix $x\in\wt M$.

\bprop\label{prop:listequivHTSRG} 
 \begin{enumerate}
 \item[(a)] If $Q_{\Gamma,\, F}(\sigma) < + \infty$ then $\mu_x(\Lambda_c
   \Ga)=0$ and $\mu^\iota_x(\Lambda_c \Ga)=0$.
 \item[(b)] If $\mu_x(\Lambda_c \Ga)=0$ and $\mu^\iota_x(\Lambda_c
   \Ga)=0$, then $(\partial_\infty^2 \wt M, \Ga, (\mu^\iota_x\otimes
   \mu_x)_{\mid \partial_\infty^2 \wt M})$ is completely dissipative.
 \item[(c)] If $(\partial_\infty^2 \wt M, \Ga, (\mu^\iota_x\otimes
   \mu_x)_{\mid \partial_\infty^2 \wt M})$ is ergodic, then
   $\mu^\iota_x\otimes\mu_x$ is atomless, the diagonal of
   $\partial_\infty \wt M \times \partial_\infty \wt M$ has zero
   measure for $\mu^\iota_x\otimes\mu_x$, and $(\partial_\infty \wt M
   \times \partial_\infty \wt M, \Gamma, \mu^\iota_x\otimes \mu_x)$ is
   ergodic and conservative.
 \item[(d)] For all $v \in T^1\wt M$, we have $v_+ \in \Lambda_c\Ga$
   if and only if there exists a relatively compact subset $K$ in
   $T^1\wt M$ such that $\int_0^\infty \mathbbm{1}_{\Ga K}\circ
   \phi_t(v) dt$ diverges. Therefore $(T^1M, (\phi_t)_{t\in\RR}, m)$
   is conservative (respectively completely dissipative) if
   and only if $\mu^\iota_x (\,^c\Lambda_c\Ga)=0$ and $\mu_x
   (\,^c\Lambda_c\Ga)=0$ (respectively $\mu^\iota_x(\Lambda_c \Ga)=0$
   and $\mu_x(\Lambda_c \Ga)=0$).
 \item[(e)] The dynamical system $(T^1 M, (\phi_t)_{t\in\RR}, m)$ is
   ergodic if and only if the dynamical system $(\partial_\infty^2 \wt
   M, \Ga, (\mu^\iota_x\otimes \mu_x)_{\mid \partial_\infty^2 \wt M})$ is
   ergodic.
 \item[(f)] If $Q_{\Gamma,\, F}(\sigma) = + \infty$ then
   $\mu_x(\Lambda_c\Ga) > 0$ and $\mu^\iota_x(\Lambda_c\Ga)>0$.
 \item[(g)] If $\mu_x(\Lambda_c\Ga) > 0$ then
   $\mu_x(\,^c\Lambda_c\Ga) = 0$.
 \item[(h)] If $(T^1M, (\phi_t)_{t\in\RR}, m)$ is
   conservative then it is ergodic.
\end{enumerate}
\eprop

\dem (a) Mohsen's shadow lemma \ref{lem:shadowlemma} applied to
$(\mu_x)_{x\in\wt M}$ with $K=\{x\}$ shows the existence of
$N_0\in\NN$ such that for all $N\geq N_0$, there exists $c_N>0$ such
that 
$$
\mu_x(\OOO_x(B(\ga x,N)))\leq C_N\;e^{\int_x^{\ga x}(\wt F-\sigma)}
\;.
$$  
By the definition of $\Lambda_c\Ga$, we have $\Lambda_c\Ga=
\bigcup_{N\in\NN,\;N\geq N_0}\Lambda(N)$ where
$$
\Lambda(N)=
\bigcap_{A \;{\rm finite~subset~of}\;\Ga}\;\;
\bigcup_{\ga\in \Ga-A} \OOO_x(B(\ga x,N))\;.
$$
The fact that $\mu_x(\Lambda(N))=0$ for every $N\geq N_0$, which implies that
$\mu_x(\Lambda_c \Ga)=0$, then follows by using the easy part of the
Borel-Cantelli lemma.

Since $Q_{\Gamma,\, F}(\sigma) < + \infty$ if and only if
$Q_{\Gamma,\, F\circ\iota} (\sigma) < + \infty$ (see Lemma
\ref{lem:elemproppressure} (ii)), the last claim follows similarly.

\medskip (b) By Fubini's theorem, the set $A = \{(\xi,
\eta)\in \partial_\infty^2 \widetilde M  \;:\;  \xi, \eta \not\in
\Lambda_c\Gamma \}$ has full measure with respect to
$(\mu^\iota_x\otimes \mu_x)_{\mid \partial_\infty^2 \wt M}$. For every
$(\xi, \eta)\in A$, let $\Ga_{\xi,\,\eta}$ be the finite nonempty set of
$\alpha\in\Ga$ such that the distance between $\alpha x$ and the
geodesic line between $\xi$ and $\eta$ is minimal; it satisfies $\ga
\Ga_{\xi,\,\eta} = \Ga_{\ga\xi,\,\ga\eta}$ for every $\ga\in\Ga$. For any
finite subset $F$ of $\Ga$, we denote $S_F = \{(\xi, \eta) \in A
\;:\; \Ga_{\xi,\, \eta} = F\}$.  We claim that the measurable subset
$S_F$ is wandering, since if $\gamma S_F$ meets $S_F$ then $\gamma F =
F$ and since the stabiliser (for the action by left translations of
$\Ga$ on itself) of a finite subset of $\Ga$ is finite. As a countable
union of wandering sets is contained in the dissipative part up to a
zero measure subset, and since $A=\bigcup_{F} S_F$ has full measure,
the result follows.

\medskip (c) Assume for a contradiction that $(\xi',\xi)
\in \partial_\infty \wt M\times\partial_\infty \wt M$ is an atom for
$\mu^\iota_x\otimes \mu_x$, that is, $\xi'$ is an atom for
$\mu^\iota_x$ and $\xi$ is an atom for $\mu_x$. This implies that
$\ga\xi$ is an atom for $\mu_x$ for every $\ga\in\Ga$ by
quasi-invariance. Let $\Ga_{\xi'} $ be the stabiliser of $\xi'$ in
$\Ga$, and note that its limit set $\Lambda\Ga_{\xi'}$ is finite (it
has $0$, $1$ or $2$ points).  Since $\Ga$ is non-elementary, there
exists $\ga_0 \in \Ga$ such that $\ga_0 \xi \neq \xi'$. By ergodicity,
the $\Ga$-orbit $\OOO$ of $(\xi', \ga_0\xi)$ in $\partial_\infty^2 \wt
M$ has full measure with respect to $(\mu^\iota_x\otimes
\mu_x)_{\mid \partial_\infty^2 \wt M}$. Since $\Gamma$ acts minimally
in $\Lambda\Ga$, there exists $\ga_1 \in \Ga$ such that $\ga_1\xi$
does not lie in the countable closed set $\Ga_{\xi'} (\ga_0 \xi) \cup
\Lambda\Ga_{\xi'}\cup\{\xi'\}$, whose complement in $\Lambda\Ga$ is
open and nonempty since $\Lambda\Ga$ is uncountable.  Note that
$(\xi',\ga_1\xi)$ does not belong to $\OOO$, since otherwise there
exists $\ga\in\Ga$ such that $\xi' = \ga\xi'$ and $\ga_1\xi =
\ga\ga_0\xi$, contradicting that $\ga_1\xi$ does not belong to
$\Ga_{\xi'} (\ga_0\xi)$. But $(\xi',\ga_1\xi)$ is an atom of
$(\mu^\iota_x\otimes \mu_x)_{\mid \partial_\infty^2 \wt M}$,
contradicting the fact that $\OOO$ has full measure with respect to
this measure.

Since the product measure $\mu^\iota_x\otimes\mu_x$ is atomless, the
diagonal of $\partial_\infty\wt M\times \partial_\infty\wt M$ has
measure zero for the product measure (by Fubini's theorem, since a
finite measure has at most countably many atoms). Hence the dynamical
system $(\partial_\infty \wt M \times \partial_\infty \wt M, \Gamma,
\mu^\iota_x\otimes \mu_x)$ is also ergodic.  Since $\Ga$ is countable
and since $\mu^\iota_x\otimes \mu_x$ is atomless, this ergodic
dynamical system is conservative (see for instance
\cite[page 51, Prop.~1.6.6]{Aaronson97}).

\medskip (d) The first claim that the conical limit points are the
endpoints of the geodesic rays whose image in $M=\Ga\backslash \wt M$
returns infinitely often in some relatively compact subset is immediate
from the definition of $\Lambda_c\Ga$. 

Note that the dynamical system $\D=(T^1M, (\phi_t)_{t\in\RR}, m)$ is
conservative if and only if every relatively compact subset $K$ of
$T^1M$ is contained in the conservative part of $\D$ up to $m$-measure
$0$. By the Halmos recurrence theorem (see for instance \cite[page 15,
Theo.~1.1.1]{Aaronson97} in the case of transformations), if $\D$ is
conservative, $m$-almost every element of any relatively compact
subset $K$ returns infinitely often (positively and negatively) in
$K$. Hence for $\wt m$-almost every $v\in T^1\wt M$, we have
$v_-,v_+\in\Lambda_c\Ga$. This holds if and only if $\mu^\iota_x
(\,^c\Lambda_c\Ga)=0$ and $\mu_x (\,^c\Lambda_c\Ga)=0$, since $d\wt
m(v)$ is a quasi-product measure of $d\mu^\iota_x(v_-)$, $d\mu_x(v_+)$
and $dt$. Conversely, if these two equalities hold, then by this
quasi-product property, for $m$-almost every $v\in T^1M$, for every
big enough relatively compact subset $K$ of $T^1M$, we have
$\int_\RR\mathbbm{1}_K(\phi_tv)\;dt=+\infty$. Hence any (measurably)
wandering set has $m$-measure $0$, and $\D$ is conservative.

\medskip (e) This is immediate by the construction of $m$.

\medskip (f) Assume that $Q_{\Gamma,\, F}(\sigma) = + \infty$. Let
$R>0$ be large enough so that Lemma \ref{lem:thickshadowlemma} (with
$y=x$) holds.  With $a\geq 0$ the maximum mass of an atom of
$\mu^\iota_x$, also assume that $R$ is large enough so that for every
$\ga\in\Ga$, we have $\mu^\iota_x(\OOO_R^-(\ga x,x))\geq
c_7=\frac{\|\mu^\iota_x\|-a}{2}$, with $\OOO_R^-(\ga x,x)$ defined in
Subsection \ref{subsec:geomnotation}. Let $K$ be the compact subset
$K(x,R)$ of $T^1\wt M$ defined in Subsection
\ref{subsec:geomnotation}. The proof of Assertion (f) relies on the
following two lemmas.

\blemm\label{lem:lemtechunpourf}
There exists $c_8> 0$ such that for every sufficiently large $T > 0$ 
$$
\int_0^T \int_0^T 
\Big(
\sum_{\alpha,\, \beta \in \Ga} 
\wt m(K \cap \phi_{-t} \alpha K \cap \phi_{-s-t} \beta K)
\Big)
\;ds\; dt \leq c_8
\Big(
\;\sum_{\ga \in \Ga,\, d(x,\ga x) \leq T} e^{\int_x^{\ga x}(\wt F -\sigma)}\;
\Big)^2\;.
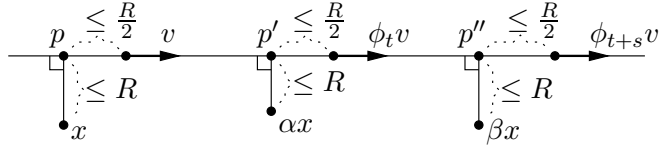
$$
\elemm

\begin{center}
\input{fig_Kalphbet.pstex_t}
\end{center}

\dem Let $s,t\geq 0$ and $\alpha,\beta\in\Ga$ with $d(x,\alpha
x),(\alpha x,\beta x), d(x,\beta x)> 2R$.  By the definition of
$K=K(x,R)$ in Subsection \ref{subsec:geomnotation}, by the triangle
inequality (and since closest point projections on geodesic lines do
not increase the distance), if $v\in K \cap \phi_{-t} \alpha K \cap
\phi_{-s-t} \beta K$, if $p,p',p''$ are the closest points to
respectively $x,\alpha x,\beta x$ on the geodesic line defined by $v$,
then $v_-,p,p',p'',v_+$ are in this order on this geodesic line (see
the above picture), and
$$
d(x,\alpha x) +
d(\alpha x,\beta x)\leq d(x,\beta x)+ 4R\;.
$$ 
By the definition of $L_R(x,\ga x)$ in Subsection
\ref{subsec:geomnotation}, we also have $(v_-,v_+)\in L_R(x,\ga x)$.
Furthermore,
$$
d(x,\alpha x) - 3R \leq t \leq d(x,\alpha x) + R\;\;\;{\rm and}\;\;\;
d(\alpha
x,\beta x) - 3R \leq s \leq d(\alpha x,\beta x) + R\;.
$$
Since $d(\alpha x,\mathopen{[}x,\beta x\mathclose{]})\leq 2R$, we also
observe by Lemma \ref{lem:technicholder} that there exists a constant
$c'\geq 0$ (depending only on $R$, on the H\"older constants of $\wt F$,
on the bounds on the sectional curvature and on
$\max_{\pi^{-1}(B(x,\,2R))}|\wt F-\sigma|$) such that
$$
\Big|\int_x^{\beta x} (\wt F-\sigma) - \int_x^{\alpha x} (\wt F-\sigma) -
\int_{\alpha x}^{\beta x} (\wt F-\sigma)\;\Big| \leq c'\;.
$$

By the definitions of $m$ in Subsection
\ref{subsec:GibbsSulivanmeasure}, since the time parameter of any unit
tangent vector $v\in K \cap \phi_{-t} \alpha K \cap \phi_{-s-t} \beta
K$ for Hopf's parametrisation with respect to the base point $x$
varies at most in an interval of length $R$, we have
$$
\wt m(K \cap \phi_{-t} \alpha K \cap \phi_{-s-t} \beta K)\leq 
R \int_{(\xi,\,\eta)\in L_R(x,\,\beta x)}
\frac{d\mu^\iota_x(\xi) \; d\mu_x(\eta)}{D_{F-\sigma,\,x}
(\xi,\eta)^2}\;.
$$
Let $(\xi,\eta)\in L_R(x,\beta x)$ and $p$ be the closest point to $x$
on the geodesic line between $\xi$ and $\eta$. By Equation
\eqref{eq:ecartvisuel} saying that 
$$
D_{F-\sigma,\,x}(\xi,\eta)=
e^{-\frac{1}{2}(C_{F-\sigma,\,\eta}(x,\,p)
  +C_{F\circ\iota-\sigma,\,\xi}(x,\,p))}\;,
$$ 
by Lemma \ref{lem:holderconseq} (1) implying that
$$
|\,C_{F-\sigma,\, \eta} (x,p)\,| \leq c_1 e^R + R
\sup_{\pi^{-1}(B(x,\, R))} |\wt F-\sigma|
$$ 
and similarly by replacing $F$ by $F\circ\iota$ and $\eta$ by $\xi$,
there exists a constant $c>0$ such that
\begin{equation}\label{eq:minomajogaps}
\frac{1}{c}\leq D_{F-\sigma,\,x}(\xi,\eta)\leq c\;.
\end{equation}
Since $L_R(x,\beta x)\subset \partial_\infty\wt M\times
\OOO^+_R(x,\beta x)$, and by Lemma \ref{lem:thickshadowlemma}, we
hence have, 
$$
\wt m(K \cap \phi_{-t} \alpha K \cap \phi_{-s-t} \beta K)\leq 
c^2 \;C\, R \;\|\mu^\iota_x\|\;  e^{\int_x^{\beta x} (\wt F - \sigma)}\;.
$$

Therefore, since $m(K \cap \phi_{-t} \alpha K \cap \phi_{-s-t} \beta
K)=0$ outside intervals of $s$ and $t$ of lengths $4R$, we have, up to
increasing the constant $C$ in order to take into account the
contribution of the elements $\alpha,\beta$ in $\Ga$ with $d(x,\alpha
x)$, $d(\alpha x,\beta x)$ or $d(x,\beta x)$ at most $2R$,
\begin{align*} {}&\int_0^T \int_0^T \Big( \sum_{\alpha,\, \beta \in \Ga}
  \wt m(K \cap \phi_{-t} \alpha K \cap \phi_{-s-t} \beta K) \Big)
  \;ds\; dt \\ &\leq (4R)^2
  \sum_{\mbox{\tiny
      $\begin{array}{c}\alpha\in\Ga,\;d(x,\,\alpha x)\leq T+3R\\
        \beta\in\Ga,\;d(\alpha x,\,\beta x)\leq T+ 3R \end{array}$}}
c^2\; C\, R \;\|\mu^\iota_x\|\;e^{c'}\;
e^{\int_x^{\alpha x} (\wt F-\sigma) + \int_{\alpha x}^{\beta x} (\wt  F-\sigma)}\\&
\leq C_3 \;\Big(\sum_{\ga\in\Ga,\;d(x,\,\ga x)\leq T+3R}
e^{\int_x^{\ga x} (\wt F-\sigma)}\;\Big)^2\;,
\end{align*}
for some constant $C_3>0$. By Corollary \ref{coro:dimgeqpress} (1), the
sum
$$
\sum_{\ga\in\Ga,\;T<d(x,\,\ga x) \leq T+3R} e^{\int_x^{\ga x} (\wt F-\sigma)}
$$ 
is uniformly bounded as $T$ tends to $+\infty$. Since the Poincar\'e
series of $(\Ga,F)$ diverges at the point $\sigma$, the result
follows. \cqfd

\blemm\label{lem:lemtechdeuxpourf}
There exists $c_9 > 0$ such that for every sufficiently large $T > 0$ 
$$
\int_0^T 
\Big(
\sum_{\ga \in \Ga} \wt m(K \cap \phi_{-t} \ga K)
\Big)
\; dt \geq c_9 \;\sum_{\ga \in \Ga,\, d(x,\ga x) \leq T} 
e^{\int_x^{\ga x}(\wt F -\sigma)} \;.
$$
\elemm

\dem The proof is similar to the previous one. Let $\ga\in\Ga$ be such
that $3R\leq d(x,\ga x)\leq T-3R$. Then by the definition of
$K=K(x,R)$, we have (see Subsection \ref{subsec:geomnotation} for the
definition of $v_{\xi,\,\eta,\,x}$)
\begin{align*}
  \int_0^T \wt m(K \cap \phi_{-t} \ga K)\;dt &= \int_{(\xi,\,\eta)\in
    L_R(x,\,\ga x)} \frac{d\mu^\iota_x(\xi) \; d\mu_x(\eta)}{D_{F-\sigma,\,x}
    (\xi,\eta)^2}\int_{-\frac{R}{2}}^{\frac{R}{2}} \int_0^T \mathbbm{1}_{\ga
    K}(\phi_{s+t}v_{\xi,\,\eta,\,x})\;dt\;ds\\ & = R^2\int_{(\xi,\eta)\in
    L_R(x,\ga x)} \frac{d\mu^\iota_x(\xi) \;
    d\mu_x(\eta)}{D_{F-\sigma,\,x}(\xi,\eta)^2}\\ & \geq \frac{1}{c^2}\;
  R^2\;\mu^\iota_x(\OOO^-_{R}(\ga x,x)) \; \mu_x(\OOO^-_{R}(x,\ga x))\;,
\end{align*}
using the fact that $L_R(x,\ga x)$ contains $\OOO^-_{R}(\ga x,x)\times
\OOO^-_{R}(x,\ga x)$ (see Lemma \ref{lem:LrOOOpm}) and the upper bound
in Equation \eqref{eq:minomajogaps}. By the assumption on $R$ at the
beginning of the proof of Assertion (f), and by the lower bound in
Lemma \ref{lem:thickshadowlemma}, we have
$$
\int_0^T \sum_{\ga\in\Ga} \wt m(K \cap \phi_{-t} \ga K)\;dt\geq 
\sum_{\ga \in \Ga,\, 3R\leq d(x,\,\ga x) \leq T-3R} 
\frac{1}{c^2}\;R^2\;c_7\;\frac{1}{C}\;e^{\int_x^{\ga x}(\wt F-\sigma)}\;.
$$ 
We conclude as in the end of the proof of the previous Lemma
\ref{lem:lemtechunpourf}.  \cqfd

\medskip Now, to prove Assertion (f) from these two lemmas, let $B$ be
the image of $K$ in $T^1M=\Ga\backslash T^1\wt M$, let $A_t = B \cap
\phi_{-t}B$, and let $\nu$ be the restriction of the measure $m$ to
the relatively compact set $B$.  Then $\int_0^\infty \nu(A_t)\; dt =
+\infty$ by Lemma \ref{lem:lemtechdeuxpourf} and the divergence of the
Poincar\'e series of $(\Ga,F)$ at the point $\sigma$. Furthermore, using
simultaneously Lemma \ref{lem:lemtechunpourf} and Lemma
\ref{lem:lemtechdeuxpourf},
$$
\int_0^T \int_0^T \nu(A_t \cap A_s)\;ds \;dt
\leq
2\int_0^T \int_0^T \nu(A_t \cap A_{s+t}) \;ds \;dt
\leq
2\frac{c_8}{{c_9}^2}
\Big(
\int_0^T \nu(A_t)\;dt
\Big)^2.
$$
By the appropriate version of the Borel-Cantelli Lemma (see for
instance \cite[page 20]{Roblin03}), we deduce that the set $\{v \in B
 \;:\; \int_0^\infty \mathbbm{1}_{A_t}(v)\; dt = +\infty\}$ has
positive measure with respect to $\nu$. By Assertion (d) and the fact
that $\wt m$ is a quasi-product measure, this implies that
$\mu_x(\Lambda_c\Ga)>0$.  Since $Q_{\Ga,\,F}(\sigma)$ diverges if and only
if $Q_{\Ga,\,F\circ\iota}(\sigma)$ diverges, we also have
$\mu^\iota_x(\Lambda_c\Ga)>0$.

\medskip (g) Observe first that, thanks to Proposition
\ref{prop:doubling} (1), the measure $\mu_x$ has the following
``doubling property of shadows'': for every $R>0$ large enough, there
exists $C=C(x,R)>0$ such that for every $\ga \in \Ga$, we have
$$
\mu_x(\OOO_{x}B(\ga x,5R))\le C\,\mu_x(\OOO_xB(\ga x,R))\;.
$$
Therefore, we deduce as in \cite[page 22-23]{Roblin03} the following
lemma. For every $R>0$, let $\Lambda(R)$ be the set of
$\xi\in\Lambda_c\Ga$ that are limits of points in $\Ga x$ at distance
at most $R$ of the geodesic ray starting from $x$ and ending at $\xi$.

\blemm \label{lem:mesurombreassg}
If $R>0$ is large enough, for every $\mu_x$-integrable map
$\varphi:\partial_\infty \wt M\ra \RR$, for $\mu_x$-almost every
$\xi\in \Lambda(R)$, when $d(x,\ga x)\ra +\infty$ and $\xi\in
\OOO_xB(\ga x,R)$, we have
$$
\frac{1}{\mu_x(\OOO_xB(\ga x,R))}
\int_{\OOO_xB(\ga x,R)} \varphi \,d\mu_x\to \varphi(\xi)\;. 
$$
\elemm

Since $\Lambda_c\Ga=\bigcup_{n\in\NN}\Lambda(n)$ and by the assumption
of assertion (g), choose $R>0$ large enough so that $\mu_x(\Lambda(R))
>0$, and Lemma \ref{lem:shadowlemma} with $K=\{x\}$ and Lemma
\ref{lem:mesurombreassg} are satisfied. Let
$B=\partial_\infty\widetilde{M}- \Lambda_c\Gamma$, let
$\varphi=\mathbbm{1}_B$ be the characteristic function of $B$ and let
$\xi\in \Lambda(R)$. Then $\varphi(\xi)=0$ and there exists a sequence
$(\gamma_n)_{n\in\NN}$ in $\Ga$ such that $d(x,\ga_n x)\ra+\infty$ and
$\xi\in \OOO_xB(\ga_nx,R)$.  By Lemma \ref{lem:mesurombreassg} (taking
$\xi$ for which it applies), we have
$$
\lim_{n\ra+\infty}\frac{\mu_x(B\cap \OOO_xB(\ga_n
  x,R))}{\mu_x(\OOO_xB(\ga_n x,R))}=0\;.
$$
Using Equation \eqref{eq:radonykodensity} and Lemma
\ref{lem:holderconseq} (2) to estimate the numerator, and Lemma
\ref{lem:shadowlemma} to estimate the denominator,  there
exists a constant $c>0$ (depending on $x,R,F$) such that
\begin{equation}\label{eq:poinclefdansg}
\frac{1}{c}\;\mu_{\ga_n x}(B\cap \OOO_xB(\ga_n x,R))
\leq \frac{\mu_x(B\cap \OOO_xB(\ga_n x,R))}{\mu_x(\OOO_xB(\ga_n x,R))}
\leq c\;\mu_{\ga_n x}(B\cap \OOO_xB(\ga_n x,R)) \;.
\end{equation}
By Equation \eqref{eq:equivdensity} and the invariance of $B$ under
$\Ga$, we hence have 
$$
\lim_{n\ra+\infty}\; \mu_x(B\cap \OOO_{\ga_n^{-1} x}B(x,R))=0\;.
$$
Choosing arbitrarily large values of $R$, we may construct an increasing
sequence $(U_k)_{k\in\NN}$ of nonempty open subsets of
$\partial_\infty X$, with $\lim_{k\ra+\infty}\mu_x(B\cap U_k)=
0$, such that the diameter of the complement of $\partial_\infty
X-U_k$ tends to $0$. Up to extracting a subsequence, the closed
subset $\partial_\infty X-U_k$ converges for the Hausdorff
distance to a singleton $\{\zeta\}$. We have
$\mu_x(B-\{\zeta\})=0$. Since $\Ga$ is non-elementary, there exists
$\ga\in\Ga$ such that $\ga\zeta\neq\zeta$. By quasi-invariance of the
measure $\mu_x$ and as $\Ga$ preserves $B$, we hence have $\mu_x(B)=0$,
as required.

\medskip (h) As in \cite{Sullivan79,Roblin03}, the proof relies on
well-known arguments of Hopf \cite{Hopf37}. The first part of it
(before Lemma \ref{lem:totconsmesprodnonatom}) could essentially be
replaced by a reference to the nice paper \cite{Coudene07}.

Fix a positive Lipschitz map $\rho: T^1M\ra \RR$ with
$\|\rho\|_{\LL^1(m)}=1$. A standard way of constructing such a $\rho$
is as follows. Let $d$ be the quotient distance on $T^1M$ of the
Riemannian distance on $T^1\wt M$, whose closed balls are compact. Fix
$v_0\in T^1M$ in the support of $m$. Let $f:\mathopen{[}0,+\infty
\mathclose{[} \ra\mathopen{]}0, +\infty \mathclose{[}$ be the
(positive, continuous, non-increasing) map, affine on each $[n,n+1]$
(hence Lipschitz since the slopes are bounded), such that $f(n)=
\frac{1}{2^nm(B(v_0,\,n+1))}$. Note that $v\mapsto f(d(v_0,v))$ is
integrable, since if $A_n=\{v\in T^1M\;:\; n\leq d(v_0,v)<n+1\}$,
then
$$
\int_{T^1M} f(d(v_0,v))\;dm(v)=
\sum_{n=0}^\infty \int_{A_n} f(d(v_0,v))\;dm(v)\leq
\sum_{n=0}^\infty \;\frac{m(A_n)}{2^n\,m(B(v_0,n+1))}\leq 2\;.
$$ 
Take for $\rho$ the integrable map $v\mapsto f(d(v_0,v))$
renormalised to have $\LL^1$-norm $1$.  Note that $\rho$ is Lipschitz,
as the composition of two Lipschitz maps.

As $(T^1M, (\phi_t)_{t\in\RR}, m)$ is conservative (by the
assumption of (h)), for $m$-almost every $v\in T^1M$, we have
$$
\int_0^{\infty}\rho\circ\phi_{-t}(v) dt=
\int_0^{\infty}\rho\circ\phi_t(v) dt=+\infty
$$ 
(see for instance \cite[Prop.~1.1.6, page 18]{Aaronson97} in the non
invertible case).
Let $\wt \rho\colon T^1\wt M\ra \RR$ be the lift of $\rho$ by the
canonical projection $T^1\wt M\ra \Ga\backslash T^1\wt M=T^1M$.

For every Lipschitz map $h\colon T^1M\ra \RR$ with compact support,
let $\wt h\colon T^1\wt M\ra \RR$ be its lift, and consider the
following maps, whose domains of definition in $T^1\wt M$ are denoted
by $\wt D^*$ and $\wt D_*$ respectively:
$$
\wt h^*\colon v\mapsto\lim_{T\ra +\infty}
\frac{\int_0^T \wt h\circ\phi_t (v)\,dt}
{\int_0^T\wt \rho\circ\phi_t (v)\,dt}
\quad \mbox{and}\quad 
\wt h_*\colon v\mapsto\lim_{T\ra +\infty}
\frac{\int_0^T \wt h\circ\phi_{-t} (v)\,dt}
{\int_0^T\wt \rho\circ\phi_{-t} (v)\,dt}\;.
$$ 
Note that $\wt D^*$, $\wt D_*$, $\wt h^*$ and $\wt h_*$ are
invariant under $\Ga$ and ($m$-almost everywhere) by the geodesic
flow, since $\int_0^{\infty}\rho\circ\phi_{\pm t}(v)\, dt=+\infty$ for
$m$-almost every $v$. Let $h^*$ and $h_*$ be the maps from the images
of $\wt D^*$ and $\wt D_*$ in $T^1M$ to $\RR$ induced by $\wt h^*$ and
$\wt h_*$ respectively. By the Lipschitz continuity of $h$ and $\rho$,
and by the fact that the distance $d(\phi_tv,\phi_tw)$ tends
exponentially to $0$ if $w\in W^{\rm ss}(v)$ as $t\ra +\infty$
(respectively $w\in W^{\rm su}(v)$ as $t\ra -\infty$), the domains $\wt
D^*$ and $\wt D_*$ are saturated respectively by the strong stable and
strong unstable leaves of the geodesic flow, and furthermore $\wt
h^*(v)$ depends only on $v_+$ and $\wt h_*(v)$ on $v_-$.

Let ${\cal I}$ be the $\sigma$-algebra of
$(\phi_t)_{t\in\RR}$-invariant Borel subsets of $T^1M$, let $\PP$ be the
probability measure $\rho m$ and let $E (\varphi\,|\,{\cal
  I})$ be the conditional expectation of $\varphi\in\LL^1(\PP)$
with respect to $\cal I$: the map $E (\varphi\,|\,{\cal I})$ is
$\PP$-integrable and almost everywhere invariant under
$(\phi_t)_{t\in\RR}$ and for every bounded measurable map
$\psi:T^1M\ra \RR$ which is almost everywhere invariant under
$(\phi_t)_{t\in\RR}$, we have
\begin{equation}\label{eq:defiespcond}
\int_{T^1M}\varphi\,\psi\;d\PP = 
\int_{T^1M} E (\varphi\,|\,{\cal I}) \,\psi\;d\PP \;.
\end{equation}
As $(T^1M, (\phi_t)_{t\in\RR}, m)$ is conservative (by the
assumption of (h)) and invertible, Hopf's ratio ergodic theorem (see
for instance \cite{Hopf37}, as well as \cite{Zweimuller04} for a short
proof) implies that $\wt D^*$ and $\wt D_*$ have full measure with
respect to $\wt m$, and that the maps $h^*$, $h_*$ and $E_{\rho m}
(\frac{h}{\rho}\,|\,{\cal I})$ are defined and coincide $m$-almost
everywhere. Let us denote by $\psi^*$ and $\psi_*$ two bounded
measurable maps on $\partial_\infty\wt M$ such that $\wt
h^*(v)=\psi^*(v_+)$ and $\wt h_*(v)=\psi^*(v_-)$ for $\wt m$-almost
every $v\in T^1\wt M$.
 
\blemm\label{lem:totconsmesprodnonatom}
If $(T^1M, (\phi_t)_{t\in\RR}, m)$ is conservative, then
the product measure $\mu^\iota_x\otimes \mu_x$ is atomless.
\elemm

\dem Assume for a contradiction that $(\xi',\xi)$ is an atom of
$\mu^\iota_x\otimes\mu_x$. Since $\Ga$ is non-elementary, there exists
$\ga\in\Ga$ such that $\ga\xi'$ is different from $\xi$ and from the
other fixed point of any loxodromic element of $\Ga$ (if there is one)
fixing $\xi$. Then by the quasi-invariance of $\mu^\iota_x$, the pair
$(\ga\xi',\xi)$ is another atom of $\mu^\iota_x\otimes\mu_x$. Let
$\epsilon>0$ be small enough. Let $v \in T^1\wt M$ be such that
$v_-=\ga\xi'$ and $v_+=\xi$.  Then the image of $\{\phi_tv \;:\;
t\in\mathopen{[}-\epsilon, \epsilon\mathclose{]}\}$ in $T^1M$ is a
measurable subset of positive measure with respect to $m$ which is a
wandering set for the geodesic flow, since $\{\ga\xi',\xi\}$ is not
the set of fixed points of a loxodromic element.  \cqfd

\medskip It follows from this lemma, as in the proof of Assertion (c),
that the diagonal of $\partial_\infty\wt M\times\partial_\infty\wt M$
has measure $0$ for the product measure $\mu^\iota_x\otimes \mu_x$.
Consider the probability measures $\ov\mu^\iota=
\frac{\mu^\iota_x}{\|\mu^\iota_x\|}$ and $\ov\mu=
\frac{\mu_x}{\|\mu_x\|}$. Since $\wt m$ is equivalent to the product
measure $d\mu^\iota_x\,d\mu_x\,dt$ on $\partial^2_\infty\wt M\times
\RR$, we have, using Fubini's theorem, the Hopf parametrisation of
Remark (2) of Subsection \ref{subsec:GibbsSulivanmeasure} (with
$x_0=x$), and the fact that $h_*$ and $h^*$ coincide almost everywhere
for the measure $d\ov\mu^\iota(v_-)\,d\ov\mu(v_+) \,dt$,
\begin{align*}
\int_{\partial_\infty\wt M} \psi^*(\eta)^2\,d\ov\mu(\eta)&=
\int_0^1\int_{\partial_\infty\wt M\times \partial_\infty\wt M} \;
 \psi^*(\eta)^2\,d\ov\mu^\iota(\xi)\,d\ov\mu(\eta)\,dt\\ &
=\int_{\partial_\infty^2\wt M\times\RR} \mathbbm{1}_{\mathopen{[}0,1\mathclose{]}}(t) \;
\psi^*(\eta)^2\,d\ov\mu^\iota(\xi)\,d\ov\mu(\eta) \,dt\\ &
=\int_{T^1\wt M} \mathbbm{1}_{\mathopen{[}0,1\mathclose{]}}(\beta_{v_-}(\pi(v),x))\;
h^*(v)^2\,d\ov\mu^\iota(v_-)\,d\ov\mu(v_+) \,dt\\ &
=\int_{T^1\wt M} \mathbbm{1}_{\mathopen{[}0,1\mathclose{]}}(\beta_{v_-}(\pi(v),x))\;
h_*(v)h^*(v)\,d\ov\mu^\iota(v_-)\,d\ov\mu(v_+) \,dt\\ &
=\int_0^1\int_{\partial_\infty\wt M\times \partial_\infty\wt M} 
\psi_*(\xi)\psi^*(\eta)\,d\ov\mu^\iota(\xi)\,d\ov\mu(\eta) \,dt
\\&=
\int_{\partial_\infty\wt M} \psi_*(\xi)\,d\ov\mu^\iota(\xi)
\int_{\partial_\infty\wt M} \psi^*(\eta)\,d\ov\mu(\eta)
\\ & =
\int_{\partial_\infty\wt M\times \partial_\infty\wt M\times\mathopen{[}0,1\mathclose{]}} 
\psi_*(\xi)\,d\ov\mu^\iota(\xi)\,d\ov\mu(\eta) \,dt
\int_{\partial_\infty\wt M} \psi^*(\eta)\,d\ov\mu(\eta)
\\ & =
\int_{\partial_\infty\wt M\times \partial_\infty\wt M\times\mathopen{[}0,1\mathclose{]}} 
\psi^*(\eta)\,d\ov\mu^\iota(\xi)\,d\ov\mu(\eta) \,dt
\int_{\partial_\infty\wt M} \psi^*(\eta)\,d\ov\mu(\eta)
\\ &=\Big(\int_{\partial_\infty\wt M} 
\psi^*(\eta)\,d\ov\mu(\eta)\Big)^2\;.
\end{align*}
By the equality case in the Cauchy-Schwarz inequality, we hence have
that $\psi^*$ is constant $\mu_x$-almost everywhere. By Fubini's
theorem and the fact that $\wt m$ is equivalent to the product measure
$d\mu^\iota_x\,d\mu_x\,dt$ on $\partial^2_\infty\wt M\times \RR$, the
preimage of a $\mu_x$-measure zero subset of $\partial_\infty \wt M$
by $v\mapsto v_+$ is a $\wt m$-measure zero subset of $T^1\wt M$. We
hence have that $\wt h^*$ is $\wt m$-almost everywhere constant.

Since the map $E (\frac{h}{\rho}\,|\,{\cal I})$ is $m$-almost
everywhere constant, it is $\PP$-almost everywhere equal to
$\int_{T^1M} \frac{h}{\rho}\;d\PP=\int_{T^1M} h\;dm$ by Equation
\eqref{eq:defiespcond}.  Let $\psi\in\LL^1(\PP)$ be invariant under
the geodesic flow. From the density in $\LL^1(\PP)$ of the subspace of
smooth maps with compact support (hence Lipschitz), we deduce that
$\int_{T^1M} \psi^2\;d\PP= (\int_{T^1M} \psi\;d\PP)^2$. By the
equality case in the Cauchy-Schwarz inequality, this implies that
$\psi$ is almost everywhere constant. Since $\PP$ and $m$ are
equivalent, this says that $(T^1M, (\phi_t)_{t\in\RR}, m)$ is
ergodic. \cqfd

\medskip
\noindent{\bf Proofs of Theorem \ref{theo:critnonergo} and Theorem
  \ref{theo:critergo}. }
We use the numbering of the assertions (1)-(5) of Theorem
\ref{theo:critnonergo}, (i)-(v) of Theorem \ref{theo:critergo} and
(a)-(g) of Proposition \ref{prop:listequivHTSRG}.  

The proof of Theorem \ref{theo:critnonergo} follows from the following
implications:

$\bullet$ (1), (2) and (3) are equivalent, by (a) and (f).

$\bullet$ (5) implies (2), by (d).

$\bullet$ (2) and (3) imply (4), by (b) and (c).

$\bullet$ (2) and (3) and (4) imply (5) by (e) and (d).

$\bullet$ Finally, assume that Assertion (4) holds, let us prove that
Assertion (1) holds. Since $(\partial_\infty^2 \wt M, \Ga,
(\mu^\iota_x\otimes \mu_x)_{\mid \partial_\infty^2 \wt M})$ is
non-ergodic, $(T^1M, (\phi_t)_{t\in\RR}, m)$ is non-ergodic by (e). By
the contrapositive of (h), $(T^1M, (\phi_t)_{t\in\RR}, m)$ is not
conservative. By the contrapositive of (d), we have
$\mu_x(\,^c\Lambda_c \Ga)>0$ or $\mu^\iota_x(\,^c\Lambda_c \Ga)>0$. By
the contrapositive of (g), which is also valid upon replacement of
$\mu_x$ by $\mu^\iota_x$, we have $\mu_x(\Lambda_c \Ga)=0$ or
$\mu^\iota_x(\Lambda_c \Ga)=0$. By the contrapositive of (f), this
proves (1).

\medskip
The proof of Theorem \ref{theo:critergo} follows from the following
implications:

$\bullet$ (ii) and (iii) are equivalent, by (g) which is also valid
upon replacement of $\mu_x$ by $\mu^\iota_x$, since (2) and (3) are.

$\bullet$ (i) implies (ii), by (f) and (g).

$\bullet$ (ii) and (iii) imply (v), by (d) and (h).

$\bullet$ (v) implies (iv), by (e) and (c).

$\bullet$ (iv) implies (i), by (c) and the contrapositives of (b) and
(a).  \cqfd

\subsection{Uniqueness of Patterson densities and 
of Gibbs states}
\label{subsec:uniqueness}

The first result in this subsection, generalising the case $F=0$ due
to Sullivan \cite{Sullivan84}, gives a useful criterion for being of
divergence type.

\bcoro \label{coro:critdivergence} Assume that there exists a Patterson
density $(\mu_x)_{x\in X}$ for $(\Ga,F)$ of dimension $\sigma\in\RR$
such that $\mu_x(\Lambda_c\Ga)>0$. Then $\sigma= \delta_{\Ga,\,F}$ and
$(\Ga,F)$ is of divergence type.  \ecoro

\dem We have $\sigma\geq \delta_{\Ga,\,F}$ by Corollary
\ref{coro:dimgeqpress} (2). If $\sigma> \delta_{\Ga,\,F}$, then the
Poincar\'e series $Q_{\Ga,\, F}$ converges at $\sigma$, and
$\mu_x(\Lambda_c\Ga)=0$ by Theorem \ref{theo:critnonergo}, a
contradiction. Hence $\sigma= \delta_{\Ga,\,F}$ and if $(\Ga,F)$ is of
convergence type, Theorem \ref{theo:critnonergo} gives the same
contradiction. \cqfd

\medskip
In particular, if $\Ga$ is convex-cocompact, then Patterson's
construction recalled in the proof of Proposition
\ref{prop:existPattdens} shows that $(\Ga,F)$ is of divergence type.

\bprop\label{prop:classifdimdens} Let $\sigma\in\RR$.

(1) If $\Ga$ is convex-cocompact, then there exists a Patterson
density for $(\Ga,F)$ of dimension $\sigma$ with support equal to
$\Lambda\Ga$ if and only if $\sigma=\delta_{\Ga,\,F}$.

(2) If $\Ga$ is not convex-cocompact, then there exists a Patterson
density for $(\Ga,F)$ of dimension $\sigma$ with support equal to
$\Lambda\Ga$ if and only if $\sigma\geq\delta_{\Ga,\,F}$.  
\eprop

\dem By Corollary \ref{coro:dimgeqpress}, we always have
$\sigma\geq\delta_{\Ga,\,F}$ if there exists a Patterson density for
$(\Ga,F)$ of dimension $\sigma$.  Since $\Lambda\Ga=\Lambda_c\Ga$ when
$\Ga$ is convex-cocompact, the first claim follows from Patterson's
construction (see Proposition \ref{prop:existPattdens}) for the
existence part when $\sigma=\delta_{\Ga,\,F}$, and from the
implication (1) implies (2) of Theorem \ref{theo:critnonergo} for the
converse statement.

To prove the second claim, by Corollary \ref{coro:dimgeqpress} and
Proposition \ref{prop:existPattdens}, we only have to prove that given
$\sigma>\delta_{\Ga,\,F}$, there exists a Patterson density for
$(\Ga,F)$ of dimension $\sigma$ with support equal to
$\Lambda\Ga$. Since $\Ga$ is not convex-cocompact, there exists a
sequence $(y_i)_{i\in\NN}$ in $\C\Lambda\Ga$ such that
$\lim_{i\ra+\infty} d(y_i,\Ga x_0)=+\infty$. For every $x\in\wt M$,
consider the nonzero (since $\sigma>\delta_{\Ga,\,F}$) measure
$$
\mu_{x,\,i}=\frac{\sum_{\ga\in\Ga}e^{\int_x^{\ga y_i}(\wt F-\sigma)} 
\D_{\ga y_i}}{\sum_{\ga\in\Ga} e^{\int_{x_0}^{\ga y_i}(\wt F-\sigma)}}\;,
$$
whose support is contained in $\C\Lambda\Ga$, and whose total mass is
bounded from above and below by a positive constant independent of $i$
by Lemma \ref{lem:technicholder}.  The family $(\mu_{x,\,i})_{x\in\wt
  M}$ has an accumulation point (for the product topology of the
weak-star topologies) $(\mu_{x})_{x\in\wt M}$ of finite nonzero
measures with support contained in (hence equal to) $\Lambda\Ga$,
which is easily proved to be a Patterson density for $(\Ga,F)$ of
dimension $\sigma$.  \cqfd

\medskip
The divergence type assumption has several nice consequences.

\bcoro\label{coro:uniqpatdens} Assume that $\delta_{\Ga,\,F}<+\infty$
and that $(\Ga,F)$ is of divergence type. Then there exists a
Patterson density $(\mu_{F,\,x})_{x\in\wt M}$ of dimension
$\delta_{\Ga,\,F}$ for $(\Ga,F)$, which is unique up to a scalar
multiple. If $\D_z$ is the unit Dirac mass at a point $z\in \wt M$,
then for every $x\in\wt M$, we have
\begin{equation}\label{eq:constructPattOK}
\frac{\mu_{F,\,x}}{\|\mu_{F,\,x}\|}=
\lim_{s\ra \,{\delta_{\Ga,\,F}}^+} \;\frac{1}{Q_{\Ga,\,F,\,x,\,x}(s)}\;
\sum_{\ga \in\Ga} e^{\int_x^{\ga x}(\wt F-s)}\;\D_{\ga x}\;.
\end{equation}
The support of the measures $\mu_{F,\,x}$ is $\Lambda\Ga$ and
furthermore $\mu_{F,\,x}(\Lambda\Ga-\Lambda_c \Ga) =0$.
We have $\mu_{F,\,x}(\Lambda\Ga-\Lambda_{\rm Myr} \Ga) =0$.  
\ecoro

Equation \eqref{eq:constructPattOK} says that Patterson's construction
of the Patterson densities (see the proof of Proposition
\ref{prop:existPattdens}) does not depend on choices of accumulation
points (and follows from the uniqueness property). This is due to
Patterson when $F=0$. In the context of Gibbs measures (with the
multiplicative convention, as explained in the beginning of Chapter
\ref{sec:GPS}), this has been first proved by Ledrappier
\cite{Ledrappier95}.

\medskip \dem The existence has been seen in Proposition
\ref{prop:existPattdens}. The uniqueness is classical (see
\cite{Sullivan79}): If $(\mu'_{F,\,x})_{x\in \wt M}$ is another
Patterson density of dimension $\delta_{\Ga,\,F}$ for $(\Ga,F)$, then
so is $(\mu_{F,\,x}+ \mu'_{F,\,x}) _{x\in \wt M}$. For any fixed $x\in
\wt M$, let $\mu=\mu_{F,\,x}$, $\mu'=\mu'_{F,\,x}$ and
$\nu=\mu+\mu'$. Note that $\mu$ is absolutely continuous with respect
to $\nu$, and that the Radon-Nikodym derivative $\frac{d\mu}{d\nu}$ is
$\Ga$-invariant and measurable, hence is $\nu$-almost everywhere equal
to a constant $\lambda>0$ by ergodicity (see Theorem
\ref{theo:critergo}, which implies that any Patterson density of
dimension $\delta_{\Ga,\,F}$ is ergodic under $\Ga$). Therefore
$\mu=\lambda \nu$ and $\mu'=(\lambda^{-1}-1)\mu$.

The fact that the support of $\mu_{F,\,x}$ is $\Lambda\Ga$ follows
from its construction (see Proposition \ref{prop:existPattdens}), and
the fact that $\mu_{F,\,x}$ gives full measure to the conical limit
set follows from Theorem \ref{theo:critergo}.

The proof that $\mu_x(\Lambda\Ga-\Lambda_{\rm Myr}\Ga)=0$ is the same
one as when $F=0$ in \cite{Roblin03} (generalising \cite{Tukia94}):
Let $\B$ be a countable basis of the topology of the support of the
Gibbs measure $\wt m_F$ associated to $(\mu_{F\circ \iota, \,x} )_{x
  \in \wt M}$ and $(\mu_{F,\,x})_{x\in\wt M}$, where $(\mu_{F\circ
  \iota,\,x} )_{x\in\wt M}$ is any Patterson density of dimension
$\delta_{\Ga,\,F}=\delta_{\Ga\circ \iota,\,F}$ for $(\Ga,F\circ
\iota)$, which exists by Proposition \ref{prop:existPattdens}.  For
every nonempty $U\in\B$, let $\M(U)$ be the set of $v\in T^1\wt M$ for
which there exist sequences $(\ga_i)_{i\in\NN}$ in $\Ga$ and
$(t_i)_{i\in\NN}$ in $\RR$ tending to $+\infty$ such that
$\phi_{t_n}(\ga_n v)\in U$.  Then $\M(U)$ is a $\Ga$-invariant
measurable subset of $T^1\wt M$. Since $(T^1M, (\phi_t)_{t\in\RR},
m_F)$ is conservative (see Theorem \ref{theo:critergo}),
the subset $U$ is contained in $\M(U)$. Since $\wt m_F(U)>0$ and by
ergodicity (see Theorem \ref{theo:critergo}), the set $\M(U)$ has full
measure. Since $\Lambda_{\rm Myr}\Ga=\{v_+\;:\;v\in \bigcap_{U\in\B}
\M(U)\}$ by definition (see Chapter \ref{sec:negacurvnot}), and since
by the quasi-product structure of $\wt m_F$, the preimage by $v\mapsto
v_+$ of a subset of positive measure has positive measure, we have
that $\mu_{F,\,x}(\Lambda\Ga-\Lambda_{\rm Myr}\Ga)=0$. \cqfd

\bprop \label{prop:atomlessreferee} 
Assume that $\delta_{\Ga,\,F}<+\infty$ and that $(\Ga,F)$ is of
divergence type. Let $(\mu_{F,\,x})_{x\in\wt M}$ be a Patterson
density of dimension $\delta_{\Ga,\,F}$ for $(\Ga,F)$, and let
$x,y\in\wt M$.

\smallskip
(1) The measure $\mu_{F,\,x}$ is atomless. 

\smallskip
(2) As $n\ra+\infty$, we have
$$
\max_{\ga\in\Ga,\;n-1<d(x,\ga y)\leq n}\;e^{\int_x^{\ga y}\wt F}\;=
\operatorname{o}\,(e^{\delta_{\Ga,\,F}\;n})\,.
$$
\eprop

\dem (1) Assume for a contradiction that $\xi\in\partial_\infty \wt M$
is an atom of $\mu_{F,\,x}$. Since $\mu_{F,\,x}$ gives full measure to
the conical limit set, we have $\xi\in\Lambda_c\Ga$, that is, there
exists a sequence of elements $(\ga_i)_{i\in\NN}$ in $\Ga$ such that
$(\ga_i \,x)_{i\in\NN}$ converges to $\xi$ while staying at bounded
distance from the geodesic ray $\mathopen{[}x,\xi\mathclose{[}\,$,
which implies in particular that $\xi\in \OOO_xB(\ga_i x,R)$ for any
$R>0$ large enough. Up to increasing $R$, by Mohsen's shadow lemma
\ref{lem:shadowlemma} with $K=\{x\}$, there exists $C>0$ such that
$\mu_{F,\,x}(\OOO_xB (\ga_i x,R))\leq C\,e^{\int_x^{\ga_i x} (\wt F-
  \delta_{\Ga,\,F})}$ for all $i\in\NN$.  Hence there exists
$\epsilon>0$ such that for all $i\in\NN$, using Equation
\eqref{eq:timereversal},
$$
e^{\int_x^{\ga_i^{-1} x} (\wt F\circ\iota- \delta_{\Ga,\,F})}=
e^{\int_x^{\ga_i x} (\wt F- \delta_{\Ga,\,F})}\geq\epsilon\,.
$$
Again by Mohsen's shadow lemma \ref{lem:shadowlemma} applied to a
Patterson density $(\mu_{F\circ\iota,\,x})_{x\in\wt M}$ of dimension
$\delta_{\Ga,\,F\circ\iota}=\delta_{\Ga,\,F}$ for $(\Ga,F\circ\iota)$,
there exists $\epsilon'>0$ such that $\mu_{F\circ\iota,\,x}(\OOO_xB
(\ga_i^{-1} x,R))\geq \epsilon'$ for all $i\in\NN$.  Up to a
subsequence, the sequence $(\ga_i^{-1} x)_{i\in\NN}$ converges to
$\eta \in \partial_\infty \wt M$, with
$\mu_{F\circ\iota,\,x}(\{\eta\}) >0$.  In particular, the measure
$\mu_{F\circ\iota,\,x}\otimes\mu_{F,\,x}$ has an atom at $(\eta,\xi)$,
a contradiction to Theorem \ref{theo:critergo} (iv) and Proposition
\ref{prop:listequivHTSRG} (c).

\medskip (2) Assume for a contradiction that there exists $c>0$ and a
sequence $(\ga_i) _{i\in\NN}$ in $\Ga$ such that
$\lim_{i\ra+\infty}d(x,\ga_i y)= + \infty$ and $e^{\int_x^{\ga_i y}\wt
  F}\geq c\;e^{\delta_{\Ga,\,F} \; d(x,\ga_i y)}$, hence
$e^{\int_x^{\ga_i y}(\wt F-\delta_{\Ga,\,F})} \geq c$.  Up to
extraction, the sequence $(\ga_i \,x)_{i\in\NN}$ converges to
$\xi\in\partial_\infty \wt M$. By Mohsen's shadow lemma
\ref{lem:shadowlemma}, $\xi$ is an atom of $\mu_{F,\,x}$, a
contradiction to (1). \cqfd

\bcoro\label{coro:descripdensconfsuppdifflimset}
Let $\sigma\in\RR$. There exists a Patterson
density of dimension $\sigma$ for $(\Ga,F)$ whose support
contains strictly $\Lambda\Ga$ if and only if
$\Lambda\Ga\neq\partial_\infty \wt M$ and the Poincar\'e series
$\sum_{\ga\in\Ga}e^{\int_x^{\ga x}(\wt F-\sigma)}$ converges.
\ecoro

\dem Assume first that $\Lambda\Ga\neq\partial_\infty \wt M$ and that
the Poincar\'e series $\sum_{\ga\in\Ga}e^{\int_x^{\ga x}(\wt
  F-\sigma)}$ converges. For every $\xi\in\partial_\infty \wt
M-\Lambda\Ga$, the construction of a Patterson density
$(\nu^{\xi}_x)_{x\in \wt M}$ of dimension $\sigma$ for $(\Ga,F)$ whose
support contains $\xi$, hence contains strictly $\Lambda\Ga$, has been
seen in the Remark following the statement of Theorem
\ref{theo:critnonergo}.

\medskip Conversely, assume that there exists a Patterson density
$(\mu_{x})_{x\in\wt M}$ of dimension $\sigma$ for $(\Ga,F)$ whose
support contains strictly $\Lambda\Ga$. In particular,
$\Lambda\Ga\neq\partial_\infty \wt M$. We have $\delta_{\Ga,\,F}\leq
\sigma$ by Corollary \ref{coro:dimgeqpress}, and in particular
$\delta_{\Ga,\,F}<+\infty$. If the Poincar\'e series
$\sum_{\ga\in\Ga}e^{\int_x^{\ga x}(\wt F-\sigma)}$ diverges, then
$\sigma=\delta_{\Ga,\,F}$, and $(\Ga,F)$ is of divergence type. Hence
by Corollary \ref{coro:uniqpatdens}, the support of $\mu_{x}$ is equal
to $\Lambda\Ga$, a contradiction.  
\cqfd

\bcoro \label{coro:uniqGibstate} Let $\sigma \in\RR$. If there exists
a Gibbs measure of dimension $\sigma$ for $(\Ga,F)$ which is finite on
$T^1M$, then $\sigma=\delta_{\Ga,\,F}$, the pair $(\Ga,F)$ is of
divergence type, there exists a unique (up to a scalar multiple) Gibbs
measure $m_F$ with potential $F$ on $T^1M$, its support is
$\Omega\Ga$, the geodesic flow $(\phi_t)_{t\in\RR}$ on $T^1M$ is
ergodic and conservative with respect to $m_F$, and the
union of the set of fixed points in $T^1\wt M$ of the elements of
$\Ga$ which do not pointwise fix $\Lambda\Ga$ has measure $0$ for $\wt
m_{F}$.  
\ecoro

It follows in particular from the first statement that if $\Ga$ is
convex-cocompact, then there exists no Patterson density of dimension
$\sigma>\delta_{\Ga,\,F}$ (see also Proposition
\ref{prop:classifdimdens}).

\medskip \dem Let $m$ be a finite Gibbs measure of dimension $\sigma$
for $(\Ga,F)$. Up to adding a constant to $F$, we may assume that
$\sigma>0$. Since $(\phi_t)_{t\in\RR}$ preserves $m$, Poincar\'e's
recurrence theorem (see for instance \cite[page 16]{Krengel85}) says
that $(T^1M, (\phi_t)_{t\in\RR}, m)$ is conservative. Hence
by Theorem \ref{theo:critnonergo}, the Poincar\'e series of $(\Ga,F)$
diverges at $\sigma$. Since there are no Patterson densities of
dimension strictly less than $\delta_{\Ga,\,F}$ (by Corollary
\ref{coro:dimgeqpress} (2)), and since $Q_{\Ga, \,F}(\sigma) < + \infty$ if
$\sigma>\delta_{\Ga,\,F}$, this implies that
$\sigma=\delta_{\Ga,\,F}$. Hence $(\Ga,F)$ is of divergence type. By
Theorem \ref{theo:critergo}, $(\phi_t)_{t\in\RR}$ is ergodic for
$m$. The uniqueness property of the Gibbs measure $m_F$ of potential
$F$ then follows by its construction from Corollary
\ref{coro:uniqpatdens}. Since the (topological) non-wandering set
$\Omega\Ga$ (see Subsection \ref{subsec:geodflowback}) is the image in
$T^1M$ of the set of the elements $v\in T^1\wt M$ such that
$v_-,v_+\in\Lambda\Ga$, the claim on the support of $m_F$ follows also
from Corollary \ref{coro:uniqpatdens}.

The union $A$ of the fixed points sets in $\wt\Omega\Ga$ of the
elements of $\Ga$ which do not pointwise fix $\Lambda\Ga$ is invariant
under the geodesic flow. It is $\Ga$-invariant, closed and properly
contained in $\wt\Omega\Ga$, by Lemma \ref{lem:fixomegaga}. Note that
$\wt\Omega\Ga$ is the support of $\wt m_F$ (see Subsection
\ref{subsec:GibbsSulivanmeasure}). Hence by ergodicity, the measure of
$A$ with respect to $\wt m_F$ is $0$. 
\cqfd

\medskip
\noindent{\bf Remark. } If $\delta_{\Ga,\,F}<+\infty$ and $(\Ga,F)$ is
of divergence type (which implies that $(\Ga,F\circ\iota)$ also is, by
Equation \eqref{eq:pressuretimereversalinva} in Lemma
\ref{lem:elemproppressure}), it follows from the uniqueness statement
in Corollary \ref{coro:uniqpatdens} that if we normalise the Patterson
density families for $(\Ga,F)$ and $(\Ga,F\circ\iota)$ given by
Corollary \ref{coro:uniqpatdens} to have total mass $1$ for a given
point $x_0\in\wt M$, and if we consider the associated Gibbs measures
$\wt m_F$ and $\wt m_{F\circ\iota}$, then

$\bullet$~ $(\mu_{F,\,x})_{x\in\wt M}=(\mu_{F+s,\,x})_{x\in\wt M}$
for every $s$ in $\RR$ (by Equation \eqref{eq:invarpressuretranslat});

\smallskip
$\bullet$~ $\iota_*\wt m_F=\wt m_{F\circ\iota}$ and $\iota_* m_F=
m_{F\circ\iota}$  (by Equation \eqref{eq:equivgap}); in particular, 
when $\wt F$ is reversible, we have  $(\mu_{F\circ\iota,\,x})_{x\in\wt M}
=(\mu_{F,\,x})_{x\in\wt M}$ and $\iota_*\wt m_F=\wt m_F$.

\smallskip $\bullet$~for every $s\in\RR$, if $F'$ is a potential
cohomologous to $F$, then $m_{F'+s}=m_F$ (by the remarks at the end of
Subsection \ref{subsec:gibbspattersondens} and Subsection
\ref{subsec:GibbsSulivanmeasure}).

\medskip In the next chapters, we will often assume that there exists
a finite Gibbs measure (of dimension $\delta_{\Ga,\,F}<+\infty$) on
$T^1M$ with potential $F$. It is then ergodic and unique up to
normalisation, by Corollary \ref{coro:uniqGibstate}. In particular,
$T^1M$, endowed with a given potential $F$, does not carry
simultaneously a finite Gibbs measure and an infinite one.

\section{Thermodynamic formalism and equilibrium 
states} 
\label{sec:variaprincip}

\medskip Let $\wt M,\Ga,\wt F$ be as in the beginning of Chapter
\ref{sec:negacurvnot}: $\wt M$ is a complete simply connected
Riemannian manifold, with dimension at least $2$ and pinched sectional
curvature at most $-1$; $\Ga$ is a non-elementary discrete group of
isometries of $\wt M$; and $\wt F :T^1\wt M\ra \RR$ is a
H\"older-continuous $\Ga$-invariant map. In this chapter, we denote by
$\phi=(\phi_t)_{t\in\RR}$ the geodesic flow, both on $T^1\wt M$ and on
$T^1M=\Ga\bs T^1\wt M$.

Given a measurable map $H:T^1M\ra\RR$, we denote by 
$$
H_-:v\mapsto \max\{0,-H(v)\}
$$ 
the {\it negative part}\index{negative part of a
  function} of $H$. Note that for every positive Borel measure $m$ on
$T^1M$ for which the negative part of $H$ is integrable, the integral
$\int_{T^1M} H\,dm$ is well defined in $\RR\cup\{+\infty\}$ as
$$
\int_{T^1M} H\,dm=\int_{T^1M}
\max\{0,H\}\,dm-\int_{T^1M} H_-\,dm\;.
$$ 
Note that if $m, m_1,m_2$ are positive Borel measures on $T^1M$ and if
$m=t m_1+(1-t)m_2$ for some $t\in\mathopen{]}0,1\mathclose{[}\,$, then
the negative part of $F$ is integrable with respect to $m$ if and only
if it is integrable with respect to both $m_1$ and $m_2$.

\medskip Let $\M_{F}(T^1M)$ be the space of $\phi$-invariant Borel
probability measures $m$ on $T^1M$ such that the negative part of $F$
is $m$-integrable, endowed with the weak-star topology. Note that the
support of any $\phi$-invariant probability measure is contained in
the (topological) non-wandering set $\Omega\Ga$. Note that
$\M_{F}(T^1M)$ is convex, but in general not closed in the convex
weak-star compact space of all $\phi$-invariant Borel probability
measures on $T^1M$.

For every $m\in \M_{F}(T^1M)$, the {\it (metric)
  pressure}\index{metric!pressure}\index{pressure!metric} of the
potential $F$ with respect to $m$ is the element of
$\RR\cup\{+\infty\}$ defined by
$$
\gls{metricpressure}=h_m(\phi)+\int_{T^1M}F\,dm\;,
$$
where $h_m(\phi)$ is the (metric) entropy of the geodesic flow $\phi$
with respect to $m$ (see \cite{KatHas95} or Subsection
\ref{subsec:measupartentrop} for the properties we will need on the
entropy).  The {\it (topological) pressure}%
\index{topological pressure}\index{pressure!topological}%
\index{pressure!of a potential} of the potential $F$ is the upper
bound of its metric (or measure-theoretic) pressures, that is, the
element of $\RR\cup\{+\infty\}$ defined by
$$
\gls{topologicalpressure}=\sup_{m\in\M_{F}(T^1M)} \;P_{\Ga,\,F}(m)\;.
$$ An element $m\in\M_{F}(T^1M)$ realising the upper bound, that is,
such that $P(\Ga,F)= \;P_{\Ga,\,F}(m)$, is called an {\it equilibrium
  state}\index{equilibrium state} for $(\Ga,F)$. Using the convexity
of $\M_{F}(T^1M)$, the affine property of the metric entropy (hence of
the metric pressure), and the ergodic decomposition of probability
measures invariant under the geodesic flow, the pressure $P(\Ga,F)$ of
$F$ is also the upper bound of the metric pressures of $F$ with
respect to the $\phi$-ergodic elements in $\M_{F}(T^1M)$.
Furthermore, by the Krein-Milman theorem, any equilibrium state is an
average of $\phi$-ergodic ones.

Note that for every $m\in\M_F(T^1M)$, we have $\iota_*m\in
\M_{F\circ\iota}(T^1M)$ and $P_{\Ga,\,F\circ\iota}(\iota_*m)=
P_{\Ga,\,F}(m)$.  Hence
$$
P(\Ga, F)=P(\Ga, F\circ\iota)\;.
$$

For all $\kappa\in\RR$ and $m\in\M_{F}(T^1M)$, we have
$P_{\Ga,\,F+\kappa} (m) = P_{\Ga,\,F}(m) +\kappa$ and $P(\Ga,F+\kappa)
= P(\Ga,F)+\kappa$; furthermore, $m$ is an equilibrium state for
$(\Ga, F+\kappa)$ if and only if $m$ is an equilibrium state for
$(\Ga,F)$.

\medskip The aim of this Chapter \ref{sec:variaprincip} is, following
the case $F\equiv 0$ in \cite{OtaPei04}, to prove the following result,
saying in particular that a finite Gibbs state, once renormalised to
be a probability measure, is the unique equilibrium state of a given
potential.

\btheo\label{theo:principevariationnel}%
\index{theorem@Theorem!Variational Principle}%
\index{variational principle} Let $\wt M$ be a
complete simply connected Riemannian manifold, with dimension at least
$2$ and pinched negative sectional curvature. Let $\Ga$ be a
non-elementary discrete group of isometries of $\wt M$. Let $\wt F
:T^1\wt M\ra \RR$ be a H\"older-continuous $\Ga$-invariant map such
that $\delta_{\Ga,\,F}<+\infty$.
\begin{enumerate}
\item[(1)] We have
\begin{equation}\label{eq:egalitevariationnelle}
P(\Ga,F)= \delta_{\Gamma,\,F}\;. 
\end{equation}
\item[(2)] If there exists a finite Gibbs measure $m_F$ for $(\Ga,F)$
  such that the negative part of $F$ is $m_F$-integrable, then $m^F=
  \frac{m_F}{\|m_F\|}$ is the unique equilibrium state for $(\Ga,F)$.
Otherwise, there exists no equilibrium state for $(\Ga,F)$. 
\end{enumerate}
\etheo

It follows in particular from Equation
\eqref{eq:egalitevariationnelle} that the topological pressure
$P(\Ga,F)$ is finite if the critical exponent $\delta_{\Gamma,\,F}$ is
finite.

\medskip \rem When $\Ga$ is torsion free and convex-cocompact, Theorem
\ref{theo:principevariationnel} has been proved in
\cite{Schapira04a}. In this case, the critical exponent $\delta_\Ga$
of $\Ga$ is the topological entropy $h$ of the geodesic flow of $M$
restricted to the (topological) non-wandering set $\Omega\Ga$ (see
\cite{Manning79}), and $\delta_{\Ga,\,F}$ is the topological pressure
$P(F)$ of the potential function $F:T^1M\ra\RR$ induced by $\wt F$
(see \cite{Ruelle81}), that is
$$
P(F)=\max_{m\in\M(T^1M)} \;\big(\;h_m(\phi)+\int_{T^1M} F\;dm\;\big)
$$ 
where $\M(T^1M)$ is the set of all Borel probability measures on
$T^1M$ invariant under the geodesic flow. Note that the reference
\cite{Schapira04a} uses the multiplicative convention for the Poincar\'e
series as mentioned in the beginning of Chapter \ref{sec:GPS}, but the
uniqueness property ensures that the two constructions of the Gibbs
measures coincide in this case. Note that any Gibbs measure is finite
in the convex-cocompact case, and uniqueness follows from general
argument of \cite{BowRue75}, hence the proof (using the Shadow lemma
and entropy computations of \cite{Kaimanovich90}) is much easier and
shorter in this particular case.

\medskip
An immediate consequence of the inequality $P(\Ga,F)\leq
\delta_{\Ga,\,F}$ is that the periods of the normalised potential are
nonpositive.

\bcoro
For every periodic orbit $g$ of the geodesic flow on $T^1M$, we have
$$
\int_g(F-\delta_{\Ga,\,F})\;\leq 0\;.
$$
\ecoro

\dem Let $\L_g$ be the Lebesgue measure along a periodic $g$ with
length $\ell(g)$ (see the beginning of Subsection
\ref{subsec:countfunct}).  The measure $\frac{\L_g}{\ell(g)}$ is a
probability measure (with compact support) on $T^1M$ invariant by the
geodesic flow, whose metric entropy is $0$ (since rotations of the
circle have entropy zero).  Hence
$$
\delta_{\Ga,\,F}\geq P(\Ga,F)\geq P_{\Ga,\,F}(m)=\int_{T^1M}F\;dm=
\frac{\L_g(F)}{\ell(g)}\;,
$$
which proves the result.
\cqfd

\medskip In particular, for every compact subset $K$ of $T^1M$, for
every measure $\mu$ on $T^1M$ which belongs to the weak-star closure
of the set of linear combinations (with nonpositive coefficients) of
the Lebesgue measures along the periodic orbits contained in $K$, we
have
$$
\int_{T^1M} (F-\delta_{\Ga,\,F})\;d\mu\;\leq 0\;.
$$

\subsection{Measurable partitions and entropy} 
\label{subsec:measupartentrop}

We refer to \cite{Rohlin67,Parry69,KatHas95} for the generalities (at
the bare minimum for what follows) on measurable partitions and
entropy recalled in this subsection, that the knowledgeable reader may
skip.

Recall that a {\it partition}\index{partition} $\zeta$ of a set $X$ is
a set of pairwise disjoint nonempty subsets of $X$ whose union is
$X$. We will denote by $\pi_\zeta:X\ra \zeta$ the canonical projection
map which maps every $x\in X$ to the (unique) element $z=\zeta(x)$ of
$\zeta$ containing $x$. Given two partitions $\zeta$ and $\zeta'$ of
$X$, we define their {\it refinement}\index{refinement} as the partition
$$
\zeta\vee\zeta'=
\{z\cap z'\;:\;z\in\zeta, z'\in\zeta', z\cap z'\neq \emptyset\}\;,
$$
and similarly for finitely many partitions.

Recall that two probability spaces with complete $\sigma$-algebras are
{\it isomorphic}%
\index{probability space!isomorphic}\index{isomorphism of probability
  spaces} if there exists a measure-preserving bimeasurable bijection
between full measure subsets of them. A probability space is {\it
  standard}\index{probability space!standard}\index{standard} if its
$\sigma$-algebra is complete and if it is isomorphic to a probability
space with underlying set the disjoint union of a compact interval
with a measurable countable set, where the restriction to the interval
is the Lebesgue $\sigma$-algebra and the Lebesgue measure, and where
the restriction to the countable set is the discrete $\sigma$-algebra.

Fix a standard probability space $(X,\A,m)$. The following result that
we state without proof is well known (see for instance \cite[Lemme
7]{OtaPei04}).

\blemm\label{lem:ergotflowtempstau}
If $(\phi_t)_{t\in\RR}$ is a $m$-ergodic measure preserving flow
on $X$, then there exists a countable subset $D$ of $\RR$ such that
$\phi_\tau$ is $m$-ergodic for every $\tau\in\RR-D$.
\elemm

For every partition $\zeta$ of the set $X$, denote by
$\widehat{\zeta}$ the sub-$\sigma$-algebra of elements of $\A$ which
are unions of elements of $\zeta$. Recall that $\zeta$ is {\it
  $m$-measurable}\index{partition!measurable}%
\index{measurable partition} (or {\it measurable} if $(\A,m)$ is
implicit) if there exist a full measure subset $Y$ in $X$ and a
sequence $(A_n)_{n\in\NN}$ in $\widehat{\zeta}$, satisfying the
following separation property: for all $z\neq z'$ in $\zeta$, there
exists $n\in\NN$ such that either we have $z\cap Y\subset A_n$ and
$z'\cap Y\subset X-A_n$ or we have $z'\cap Y\subset A_n$ and $z\cap
Y\subset X-A_n$. For instance, a finite partition of $X$ is measurable
if and only if each element of the partition is measurable.
Given two $m$-measurable partitions $\zeta,\zeta'$, we write
$\zeta\preceq \zeta'$ (the conventions differ amongst the references)
if $\zeta'(x)\subset \zeta(x)$ for $m$-almost every $x\in X$, and we
say that $\zeta$ and $\zeta'$ are {\it $m$-equivalent}%
\index{partition!measurable!equivalent}%
\index{equivalence!of measurable partitions}%
\index{measurable partition!equivalence of} if $\zeta\preceq \zeta'$
and $\zeta'\preceq \zeta$. Given a sequence of $m$-measurable
partitions $(\zeta_n)_{n\in\NN}$, there exists a $m$-measurable
partition, denoted by $\zeta=\bigvee_{n\in\NN}\zeta_n$ and unique up
to $m$-equivalence, such that $\zeta_n\preceq \zeta$ for every
$n\in\NN$, and if $\zeta'$ is a $m$-measurable partition such that
$\zeta_n\preceq \zeta'$ for every $n\in\NN$, then $\zeta\preceq
\zeta'$.

\medskip 
Fix a measurable partition  $\zeta$ of $(X,\A,m)$. The triple
$$
\big(\;\zeta,\;\;\A_\zeta=\pi_\zeta(\widehat{\zeta})=
\{\pi_\zeta(A)\;:\;A\in\widehat{\zeta}\},\;\; 
m_\zeta=(\pi_\zeta)_*m:A\mapsto m(\pi_\zeta^{-1}(A))\;\big)
$$
is a standard probability space, called the {\it factor space}%
\index{factor space}\index{probability space!factor space of} of
$(X,\A,m)$ by $\zeta$. In particular, saying that for $m_\zeta$-almost
every $z\in\zeta$, the set $z$ satisfies some property and for
$m$-almost every $x\in X$, the set $\zeta(x)$ satisfies this property are
equivalent.

There exists a family $(\A_z,m_z)_{z\in\zeta}$ (sometimes considered
as defined only for $m_\zeta$-almost every $z$ in $\zeta$), where
$(\A_z,m_z)$ is a standard probability measure on $z$, such that:

(1) for every $\A$-measurable map $f:X\ra \CC$, the map
$f_{\mid z}$ is $\A_z$-measurable for $m_\zeta$-almost every $z$ in
$\zeta$;

(2) the measure $m$ disintegrates with respect to the
projection $\pi_\zeta$, with family of {\it conditional
  measures}\index{conditional measures}\index{measure!conditional} on
the fibres the family $(m_z)_{z\in\zeta}$: for every $\A$-integrable
map $f:X\ra \CC$, the map $z\mapsto \int_z f_{\mid z}\;dm_z$ is
$m_\zeta$-integrable and
$$
\int_X f\;dm=\int_{z\in\zeta}\int_zf_{\mid z}\;d m_z\;dm_\zeta(z)\;.
$$
This last integral is equal to $\int_{x\in X}\int_{\zeta(x)}f_{\mid
  \zeta(x)}\; d m_{\zeta(x)}\;dm(x)$.  If two such families
$(\A_z,m_z)_{z\in\zeta}$ and $(\A'_z,m'_z)_{z\in\zeta}$ satisfy
properties (1) and (2), then $m_z=m'_z$ for $m_\zeta$-almost every
$z\in\zeta$. The properties (1) and (2) of such a family
$(\A_z,m_z)_{z\in\zeta}$ are satisfied if they are satisfied for $f$
the characteristic functions of elements of $\A$.

If $\zeta$ is finite or countable, then for every $z\in\zeta$ such
that $m(z)\neq 0$, we have $m_z=\frac{1}{m(z)}\,m_{\mid z}$.

\medskip Fix a measure-preserving transformation $\varphi:X\ra X$. We
denote by $\varphi^{-1}\zeta$ the $m$-measurable partition of $X$ by
the preimages of the elements of $\zeta$ by $\varphi$. We say that
$\zeta$ is {\it generating}\index{generating measurable partition}%
\index{measurable partition!generating}%
\index{partition!measurable!generating} under $\varphi$ if the
$m$-measurable partition $\bigvee_{k\in\NN}\varphi^{-k}\zeta$ is
$m$-equivalent to the partition by points.
%
%
For all measurable partitions $\zeta$ and $\zeta'$ of $X$, we define
the {\it entropy of $\zeta'$ relative to $\zeta$}\index{entropy!of a
  measurable partition relative to another one}%
\index{relative entropy} by
$$
H_m(\zeta'\,|\,\zeta)=\int_{x\in X}-\ln
m_{\zeta(x)}(\zeta'(x)\cap\zeta(x))\;dm(x)\;,
$$
with the convention that $-\ln 0=+\infty$. Note that if $\zeta$ and
$\zeta'$ are finite or countable, then 
$$
H_m(\zeta'\,|\,\zeta)=
\sum_{z\in\zeta,\,z'\in\zeta', \,m(z)\neq 0} - m(z'\cap z)
\ln\frac{m(z'\cap z)}{m(z)}\;,
$$ 
as more usual.  If $\zeta\preceq \varphi^{-1}\zeta$, we define the {\it
  entropy of $\varphi$ relative to $\zeta$}\index{entropy!of a
  measurable transformation relative to a measurable
  partition}\index{relative entropy} by
$$
h_m(\varphi,\zeta)=H_m(\varphi^{-1}\zeta\,|\,\zeta)\;.
$$
The {\it (metric) entropy}\index{metric!entropy}\index{entropy!metric}
$h_m(\varphi)$ of $\varphi$ with respect to $m$ is the upper bound of
$h_m(\varphi,\zeta)$ for $\zeta$ a measurable partition such that
$\zeta\preceq \varphi^{-1}\zeta$. The {\it (metric)
  entropy}\index{metric!entropy!of a flow}\index{entropy!metric!of a
  flow} $h_m(\phi)$ of a flow $\phi=(\phi_t)_{t\in\RR}$ of
measure-preserving transformations of $X$ is
$h_m(\phi)=h_m(\phi_1)$. We have $h_m(\phi_t)=|t|h_m(\phi)$ for every
$t\neq 0$.

\subsection{Proof of the Variational Principle} 
\label{subsec:proofvaraprincip}

The notation for this proof is the following one. Let $\delta
=\delta_{\Ga,\,F}$ and let $\wt m_F$ be a Gibbs measure on $T^1\wt M$
for $\Ga$ with potential $F$: there exist Patterson densities
$(\mu^\iota_{x})_{x\in\wt M}$ and $(\mu_{x})_{x\in\wt M}$ of dimension
$\delta$ for $(\Ga,F\circ \iota)$ and $(\Ga,F)$ such that, using the
Hopf parametrisation $T^1\wt M\ra\partial_\infty^2\wt M\times\RR$
identifying $v$ and $(v_-,v_+,t)$,
\begin{equation}\label{eq:rappeldefmesgibbs}
d\,\wt m(v)= e^{C_{F\circ \iota-\delta,\,v_-}(x,\,\pi(v))+
C_{F-\delta,\,v_+}(x,\,\pi(v))} d\mu^\iota_{x}(v_-)\,d\mu_{x}(v_+)
\,dt
\end{equation}
for any $x\in\wt M$. Let $m_F$ be the induced measure on $T^1 M$ (see
Subsection \ref{subsec:GibbsSulivanmeasure}).  Denote by
$(\mu_{W^{\rm su}(v)})_{v\in T^1\wt M}$ and $(\mu_{W^{\rm su}(v)})_{v\in
  T^1M}$ the families of measures on the strong unstable leaves of
$T^1\wt M$ and $T^1M$ associated with the Patterson density
$(\mu_{x})_{x\in\wt M}$ as in Subsection
\ref{subsec:condmesgibbs}.

\medskip Before beginning the proof of Theorem
\ref{theo:principevariationnel}, we introduce another cocycle.  For
all $v,w$ in the same strong unstable leaf in $T^1M$, we define
\begin{equation}\label{eq:defcocycleforunstab}
c_F(v,w)=\lim_{t\ra+\infty}
\int_0^t \big( F(\phi_{-s}v)-F(\phi_{-s}w)\big)\;ds\;,
\end{equation}
which exists since $F$ is H\"older-continuous and $d(\phi_{-t}v,
\phi_{-t}w)$ tends exponentially fast to $0$. Note that $c_F$ is
unchanged by adding a constant to $F$, and that it satisfies the
cocycle property
$$
c_F(v,v') + c_F(v',v'')= c_F(v,v'')
$$ 
for all $v,v',v''$ in the same strong unstable leaf in $T^1M$.  For
all $t\in\RR$ and $v,w$ in the same strong unstable leaf in $T^1M$,
we have
\begin{equation}\label{eq:cocyclebasvsflogeod}
c_F(\phi_tv,\phi_tw)=c_F(v,w)+
\int_0^t \big(F(\phi_{s}v)-F(\phi_{s}w)\big)\;ds\;.
\end{equation}
The relation between this cocycle and the Gibbs cocycle is the
following.  For all $v,w$ in the same strong unstable leaf in
$T^1M$, for all lifts $\wt v$ and $\wt w$ of $v$ and $w$,
respectively, in the same strong unstable leaf in $T^1\wt M$, we have,
for every $\sigma\in\RR$,
\begin{equation}\label{eq:relatpetitcgrandC}
c_F(v,w)= - \,C_{F\circ \iota -\sigma,\,\wt v_-}(\pi(\wt v),\pi(\wt w))\;.
\end{equation}
This follows from the definition of the Gibbs cocycle, since
$x=\pi(\wt v)$ and $y=\pi(\wt w)$ being in the same horosphere
centred at $\wt v_-$ implies that $C_{F\circ \iota -\sigma,\,\wt
  v_-}(x,y) =C_{F\circ \iota,\,\wt v_-}(x,y)$, and since the map
$t\mapsto \phi_t(-\wt v)=-\,\phi_{-t}(\wt v)$ from
$\mathopen{[}0,+\infty\mathclose{[}$ to $\wt M$ is the geodesic ray
from $x$ to $\wt v_-$.

\medskip 
\noindent{\bf Step 1. } The first step of the proof of Theorem
\ref{theo:principevariationnel} is a construction of measurable
partitions with nice geometric properties which allow entropy
computations.

Here is what we mean by ``nice geometric properties'', exactly as in
\cite{OtaPei04}. Given a probability measure $m$ on $T^1M$, a
partition $\zeta$ of $T^1M$ is called {\it
  $m$-subordinated}\index{partition!subordinated to a
  foliation}\index{subordinated partition} to the strong unstable
partition $\W^{\rm su}$ of $T^1M$ if for $m$-almost every $v$ in
$T^1M$, the set $\zeta(v)$ is a relatively compact neighbourhood of
$v$ in $W^{\rm su}(v)$.

The next result constructs a measurable partition, subordinated to the
strong unstable foliation and realising the upper bound in the
definition of the entropy, for every ergodic invariant probability
measure on $T^1M$. It is due to \cite{LedStr82} for Anosov
transformations of compact manifolds, and the adaptation we will use
to our geodesic flows is the content of the propositions 1 and 4 of
\cite{OtaPei04}.

\bprop[Otal-Peigné]\label{prop:partition} Let $m$ be a probability
measure on $T^1M$ invariant under the geodesic flow $(\phi_t)_{t\in\RR}$,
and let $\tau>0$ be such that $\phi_\tau$ is ergodic for $m$.  Then
there exists a $m$-measurable partition $\zeta$ of $T^1M$, which is
generating under $\phi_\tau$ and subordinated to $\W^{\rm su}$, such
that $(\phi_\tau^{-1}\zeta)(v)$ is contained in $\zeta(v)$ for
$m$-almost every $v$, and such that
$$
h_m(\phi_\tau)=h_m(\phi_\tau,\zeta)\;.
$$
\eprop
 
\rem The assumptions of \cite{OtaPei04} and ours are the same except
that we allow $\Ga$ to have fixed points on $\wt M$. Here are
a few comments. The proof of Proposition 1 in \cite{OtaPei04} starts
by fixing $u\in T^1M$ in the support of the measure $m$. For every
$r>0$, the authors consider the standard dynamical neighbourhood $U_r$
of $u$ defined by
$$
U_r=\bigcup_{|t|<r}\phi_t\Big(\bigcup_{v\in \B^{\rm ss}(u,\,r)}
\B^{\rm su}(v,r)\Big)\;,
$$
where $\B^{\rm ss} (u,r)$ is the open $r$-neighbourhood of $u$ in its
strong stable leaf for the induced Riemannian metric, and similarly
for $\B^{\rm su}(v,r)$. Noting that $U_r$ is the domain of a local
chart of the foliation $\W^{\rm su}$ if $r$ is small enough, they
define $\zeta$ as $\bigvee_{n=0}^{+\infty} \phi_{n\tau}
\widehat{\zeta'}$ where $\zeta'$ is the partition of $T^1M$ whose
elements are $T^1M-U_r$ and the local leaves of $\W^{\rm su}$ in
$U_r$. They then prove that for Lebesgue-almost every $r>0$ such that
$4r$ is strictly less than the injectivity radius of $M$ at the origin
$\pi(u)$ of $u$, the partition $\zeta$ is as required.

In our case, starting with the same $u$, we first consider a lift $\wt
u$ of $u$ in $T^1\wt M$. The neighbourhood $\wt{U}_r$ of $\wt u$
constructed as above is then invariant under the stabiliser $\Ga_{\wt
  u}$ of $\wt u$ in $\Ga$. We take $r$ small enough such that the
images of $\wt{U}_{4r}$ by the elements of the finite set
$\Ga_{\pi(\wt u)}-\Ga_{\wt u}$ are pairwise disjoint (where
$\Ga_{\pi(\wt u)}$ is the stabiliser of $\pi(\wt u)$ in $\Ga$) and
such that the restriction to $\wt{U}_{4r}$ of the canonical projection
from $T^1\wt M$ to $T^1M$ induces a homeomorphism from $\Ga_{\wt u}
\bs\wt{U}_{4r}$ onto an open neighbourhood of $u$. The construction in
\cite{OtaPei04} then provides the above result.

Furthermore, the assertion in Proposition \ref{prop:partition} that
$(\phi_\tau^{-1}\zeta)(v)$ is contained in $\zeta(v)$ for $m$-almost
every $v$ is not stated in \cite[Prop.~1]{OtaPei04} but is proved
there. It implies that $\zeta\preceq \phi_{-\tau}\zeta$ (as needed for
the computation of the relative entropy), but we will need an actual
containment and not a containment up to a set of $m$-measure $0$ in
the third step of the proof.

\medskip 
\noindent{\bf Step 2. } In order to apply Proposition
\ref{prop:partition} and the formula giving the entropy relative to a
measurable partition, we need to be able to compute the conditional
measures. In this second step of the proof of Theorem
\ref{theo:principevariationnel}, we show that a family of conditional
measures $(m_z)_{z\in\zeta}$ of a finite Gibbs measure (normalised to
be a probability measure) with respect to a subordinated measurable
partition is a family of measures absolutely continuous with respect
to (the restrictions to the elements of the partition of) the measures
$\mu_{W^{\rm su}(v)}$ on the strong unstable leaves, constructed using the
above cocycle.  The following lemma is a generalisation (up to
appropriate normalisations) of Proposition \ref{prop:disintegrGibbs},
when the partition was by the whole strong unstable leaves.

\blemm\label{lem:explicitconditionelles} Assume that $m_F$ is finite
and let $m^F=\frac{m_F}{\|m_F\|}$. Let $\zeta$ be a $m^F$-measurable
partition of $T^1M$ which is $m^F$-subordinated to $\W^{\rm su}$.
Then a family of conditional measures $(m^F_z)_{z\in\zeta}$ of $m^F$ with
respect to $\zeta$ is given by  
\begin{equation}\label{eq:condmesuputative}
dm^F_{z}(w)= \frac{\mathbbm{1}_{z}(w)}{\int_{z}e^{-c_F(w,u)}
\,d\mu_{W^{\rm su}(v)}(u)}\;d\mu_{W^{\rm su}(v)}(w) \;,
\end{equation}
where $w\in z\in\zeta$, for any $v\in T^1M$ such that $z\subset W^{\rm
  su}(v)$.  
\elemm

Note that when $F$ is constant (and in particular when $m_F$ is the
Bowen-Margulis measure), $m^F_{z}$ is just the normalised restriction
to $z$ of $\mu_{W^{\rm su}(v)}$.

\medskip \dem For every relatively compact Borel subset $z$ of a
strong unstable leaf $W^{\rm su}(v)$ in $T^1M$, such that the support
of $\mu_{W^{\rm su}(v)}$ meets the interior of $z$, let us define a
finite measure $m'_{z}$ on $z$ by
$$
dm'_{z}(w)= \frac{\mathbbm{1}_{z}(w)\;d\mu_{W^{\rm su}(v)}(w)}{\int_{z} 
   e^{-c_F(w,u)} \,d\mu_{W^{\rm su}(v)}(u)} \;,
$$
noting that the denominator is positive and finite and that the
numerator is a finite measure. Furthermore, $m'_{z}$ is a Borel
probability measure on $z$, since by the cocycle property of $c_F$,
$$
\|m'_{z}\|= \frac{\int_{z} 
   e^{-c_F(v,w)} \,d\mu_{W^{\rm su}(v)}(w)}{\int_{z} 
   e^{-c_F(v,u)} \,d\mu_{W^{\rm su}(v)}(u)}=1\;.
$$
Let us check the two properties for $(m'_{z})_{z\in\zeta}$ to be a family
of conditional measures for $m^F$ with respect to $\zeta$.

In order to check the measurability property, let $f:T^1M\ra \CC$ be a
measurable map, that we may assume to be bounded with compact
support. For every $n\in\NN$, let $X_n$ be the Borel subset of points
of $T^1\wt M$ whose stabiliser in $\Ga$ has cardinality $n$, which is
$\Ga$-invariant and $\phi$-invariant. Since $m_F$ is finite, hence
ergodic by Corollary \ref{coro:uniqGibstate}, we fix $n_0\in\NN$ such
that $X_{n_0}$ has full measure for $\wt m_F$. We may hence assume
that $f$ vanishes outside the image of $X_{n_0}$ in $T^1 M$. Up to
multiplying by $n_0$ the maps and measures when lifting them to
$T^1\wt M$ (see Subsection \ref{subsec:pushmeas}), we may assume that
$n_0=1$. For $m^F$-almost every $v\in T^1M$, since $z=\zeta(v)$ is a
measurable subset of $T^1M$, the restriction of $f$ to $z$ is
measurable, and
\begin{align*}
v\mapsto \int_{z} f_{\mid z}\;dm'_z& =
\int_{w\in T^1M} \frac{\mathbbm{1}_{\zeta(v)}(w)f(w)}{\int_{u\in
    T^1M}\mathbbm{1}_{\zeta(v)}(u) \,e^{-c_F(w,u)} \,d\mu_{W^{\rm su}(v)}(u)}
\,d\mu_{W^{\rm su}(v)}( w)
\end{align*}
is clearly measurable, since $v\mapsto \mu_{W^{\rm su}(v)}$ is weak-star
continuous and $\mathbbm{1}_{\zeta(v)}(w)=\mathbbm{1}_{\zeta(w)}(v)$.

In order to check the disintegration property, we give a sequence of
equalities, whose guide (to be read simultaneously with the
equalities below) is the following:

$\bullet$~ for the first equality, we use the definition of the
measures $(m'_{z})_{z\in\zeta}$ and the topological properties of the
elements of the partition $\zeta$;

$\bullet$~ for the second equality, we use a lifting by the canonical
projection $\prT:T^1\wt M\ra T^1M$; we denote by $\F_\Ga$ a measurable
(strict) fundamental domain for the action of $\Ga$ on $T^1\wt M$ and
by $\wt f=f\circ \prT$ the lift of $f$ to $T^1\wt M$; for every
$v\in T^1\wt M$, let 
$$
\wt\zeta(v)=W^{\rm su}(v)\cap (\prT)^{-1}\big(\zeta(\prT(v))\big)\;;
$$
note that for $\wt m_F$-almost every $v\in T^1\wt M$, the restriction
$\prT_{\mid\wt\zeta(v)}: \wt\zeta(v)\ra \zeta(\prT(v))$ is a
bi-measurable bijection: otherwise, there would exist $\ga\in \Ga$
sending an element in $W^{\rm su}(v)$ to a distinct element in $W^{\rm
  su}(v)$, thus (we assumed $\Ga$ to be torsion-free) $\ga$ is a
parabolic element fixing the parabolic fixed point $v_-$, and since
$\wt m_F$ is finite and the set of parabolic fixed points is
countable, the set of $v\in T^1\wt M$ such that $v_-$ is parabolic has
measure $0$ for $\wt m_F$ by Proposition \ref{prop:atomlessreferee}
(1) and Corollary \ref{coro:uniqGibstate};  we also use
Equation \eqref{eq:relatpetitcgrandC}; note that unless the function
which is integrated is $0$, we have $w_-=v_-$;

$\bullet$~ for the third  equality, we use the fact that
$\sum_{\ga\in\Ga}\mathbbm{1}_{\F_\Ga}(\ga w)=1$ by the definition of a
fundamental domain;

$\bullet$~ for the fourth equality, we use the change of variables
$v\mapsto \ga v$ and $w\mapsto \ga w$, the invariance under $\Ga$ of
$\wt f$, of the Gibbs cocycle and of the Gibbs measure $\wt m_F$, the
equivariance properties under $\Ga$ of the strong stable leaves and
their measures, as well as the equivariance property $\wt\zeta(\ga v)
=\ga\wt\zeta(v)$ for $\wt m_F$-almost every $v\in T^1\wt M$ which
follows from the construction of the lift $\wt\zeta(v)$ of $\zeta(v)$;

$\bullet$~ for the fifth equality, we first simplify by
$\sum_{\ga\in\Ga}\mathbbm{1}_{\F_\Ga}(\ga^{-1} v)=1$, and then, fixing
$x\in \wt M$, we apply the definition of the measures $\mu_{W^{\rm
    su}(v)}$ (see Equation \eqref{eq:defmeasstrongunstable}) and $\wt
m_F$ (see Equation \eqref{eq:rappeldefmesgibbs}); we use the Hopf
parametrisation $v=(v_-,v_+,t)$ with respect to $x$ and the
homeomorphism $w\mapsto w_+$ from $W^{\rm su}(v)$ to $\partial_\infty
\wt M-\{v_-\}$; we note that, since $m_F$ is finite, the diagonal has
measure $0$ for $\mu^\iota_x\otimes \mu_x$ (see just after Lemma
\ref{lem:totconsmesprodnonatom}) and the Patterson density
$(\mu_x)_{x\in\wt M}$ is atomless (see Proposition
\ref{prop:atomlessreferee} and Corollary \ref{coro:uniqGibstate});

$\bullet$~ for the sixth equality, we use the cocycle equality
$$
C_{F\circ \iota-\delta,\,v_-}(x,\pi(v))= C_{F\circ
  \iota-\delta,\,v_-}(x,\pi(w))+C_{F\circ
  \iota-\delta,\,v_-}(\pi(w),\pi(v))\;;
$$ 
we also use that since the relations $w\in W^{\rm su}(v)$ and
$\prT(w)\in\zeta(\prT(v))$ between elements $w,v\in T^1\wt M$ are
symmetric (as belonging to the same element of $\zeta$ is a symmetric
relation between two elements in $T^1M$), we have $w\in \wt\zeta(v)$
if and only if $v\in \wt\zeta(w)$, and this implies that
$\wt\zeta(v)=\wt\zeta(w)$; furthermore, unless the function which is
integrated is $0$, we have $w_-=v_-$ and $W^{\rm su}(v)= W^{\rm
  su}(w)$, and we denote $w=(v_-,w_+,t)$;

$\bullet$~ for the seventh equality, we apply Fubini's theorem
(exchanging the integrations on $v_+$ and $w_+$);

$\bullet$~ for the last ones, again since the map $v\mapsto v_+$ from
$W^{\rm su}(w)$ to $\partial_\infty \wt M-\{w_-\}$ is a homeomorphism, we
use the Hopf parametrisation $w=(w_-=v_-,w_+,t)$ with respect to $x$
and the fact that some ratio of two integrals is $1$:

\begin{align*}
  &\int_{v\in T^1M} \int_{w\in \zeta(v)} f(w)\,
  dm'_{\zeta(v)}(w)dm^F(v)\\ =& \int_{v\in T^1M} \int_{w\in W^{\rm su}(v)}
  \frac{\mathbbm{1}_{\zeta(v)}(w)f(w)\;d\mu_{W^{\rm su}(v)}(w)\;dm_F(v)}
  {\int_{u\in W^{\rm su}(v)} \mathbbm{1}_{\zeta(v)}(u)\,e^{-c_F(w,\,u)}
    \;d\mu_{W^{\rm su}(v)}(u)\;\|m_F\|}\\= & \int_{v\in T^1\wt M}
  \int_{w\in W^{\rm su}(v)} \mathbbm{1}_{\F_\Ga}(v)\,
    \frac{\mathbbm{1}_{\wt\zeta(v)}(w)\,\wt f(w)\;
    d\mu_{W^{\rm su}(v)}(w)\;d\wt m_F(v)} {\int_{u\in W^{\rm su}(v)}
    \mathbbm{1}_{\wt\zeta(v)}(u)\,e^{C_{F\circ\iota
  -\delta,\,w_-}(\pi(w),\,\pi(u))} \;d\mu_{W^{\rm su}(v)}(u)\;\|m_F\|}\\
= & \int_{v\in T^1\wt M}\int_{w\in W^{\rm su}(v)} 
    \sum_{\ga\in\Ga}\mathbbm{1}_{\F_\Ga}(v)\,
    \mathbbm{1}_{\F_\Ga}(\ga w)\\ &  \hspace*{4.5cm}
    \frac{\mathbbm{1}_{\wt\zeta(v)}(w)\,\wt f(w)\;
    d\mu_{W^{\rm su}(v)}(w)\;d\wt m_F(v)} {\int_{u\in W^{\rm su}(v)}
    \mathbbm{1}_{\wt\zeta(v)}(u)\,e^{C_{F\circ\iota
  -\delta,\,w_-}(\pi(w),\,\pi(u))} \;d\mu_{W^{\rm su}(v)}(u)\;\|m_F\|}\\
= & \int_{v\in T^1\wt M}\int_{w\in W^{\rm su}(v)} 
    \;\sum_{\ga\in\Ga}\mathbbm{1}_{\F_\Ga}(\ga^{-1}v)\,
    \mathbbm{1}_{\F_\Ga}(w)\\ &  \hspace*{4.5cm}
    \frac{\mathbbm{1}_{\wt\zeta(v)}(w)\,\wt f(w)\;
    d\mu_{W^{\rm su}(v)}(w)\;d\wt m_F(v)} {\int_{u\in W^{\rm su}(v)}
    \mathbbm{1}_{\wt\zeta(v)}(u)\,e^{C_{F\circ\iota
  -\delta,\,w_-}(\pi(w),\,\pi(u))} \;d\mu_{W^{\rm su}(v)}(u)\;\|m_F\|}
\\
=& \int_{\tiny\!\!
    \begin{array}{l}t\in \RR\\v_-\in \partial_\infty\wt M\\
      v_+\in \partial_\infty\wt M\\v_+\neq v_-\end{array}}\!\!
  \int_{\tiny\!\!
    \begin{array}{l}w_+\in \partial_\infty\wt M\\
w_+\neq v_-\end{array}}\frac{\mathbbm{1}_{\F_\Ga}(w)\,
\mathbbm{1}_{\wt\zeta(v)}(w)\,\wt f(w)}
  {\int_{u\in W^{\rm su}(v)}\mathbbm{1}_{\wt\zeta(v)}(u) \,
e^{C_{F\circ\iota -\delta,\,w_-}(\pi(w),\,\pi(u))+C_{F-\delta,\,u_+}(x,\,\pi(u))}
\;d\mu_{x}(u_+)}\\ &  \hspace*{1.8cm}
e^{C_{F-\delta,\,w_+}(x,\,\pi(w))+C_{F\circ \iota-\delta,\,v_-}(x,\,\pi(v))
+C_{F-\delta,\,v_+}(x,\,\pi(v))}\;
\frac{d\mu_{x}(w_+)d\mu^\iota_{x}(v_-)d\mu_{x}(v_+)dt}  {\|m_F\|}\\
  =& \int_{\tiny\!\!
    \begin{array}{l}t\in \RR\\v_-\in \partial_\infty\wt M\\
      v_+\in \partial_\infty\wt M\\v_+\neq v_-\end{array}}\!\!
  \int_{\tiny\!\!
    \begin{array}{l}w_+\in \partial_\infty\wt M\\
w_+\neq v_-\end{array}}\frac{\mathbbm{1}_{\F_\Ga}(w)\,
\mathbbm{1}_{\wt\zeta(w)}(v)\,\wt
    f(w)\,e^{C_{F\circ \iota-\delta,\,w_-}(\pi(w),\,\pi(v))
+C_{F-\delta,\,v_+}(x,\,\pi(v))}}
  {\int_{u\in W^{\rm su}(w)}\mathbbm{1}_{\wt\zeta(w)}(u) \,
e^{C_{F\circ\iota -\delta,\,w_-}(\pi(w),\,\pi(u))+C_{F-\delta,\,u_+}(x,\,\pi(u))}
\;d\mu_{x}(u_+)}\\ &  \hspace*{3.5cm}
\frac{e^{C_{F-\delta,\,w_+}(x,\,\pi(w))+C_{F\circ \iota-\delta,\,v_-}(x,\,\pi(w))}\,
d\mu_{x}(w_+)d\mu^\iota_{x}(v_-)d\mu_{x}(v_+)dt}  {\|m_F\|}\\
  =& \int_{\tiny\!\!
    \begin{array}{l}t\in \RR\\v_-\in \partial_\infty\wt M\\
      w_+\in \partial_\infty\wt M\\
w_+\neq v_-\end{array}}
  \frac{\int_{\tiny\!\!
    \begin{array}{l}v_+\in \partial_\infty\wt M\\
v_+\neq w_-\end{array}}\mathbbm{1}_{\wt\zeta(w)}(v)\,
e^{C_{F\circ \iota-\delta,\,w_-}(\pi(w),\,\pi(v))
+C_{F-\delta,\,v_+}(x,\,\pi(v))}\;d\mu_{x}(v_+)}
  {\int_{u\in W^{\rm su}(w)}\mathbbm{1}_{\wt\zeta(w)}(u) \,
e^{C_{F\circ\iota -\delta,\,w_-}(\pi(w),\,\pi(u))+C_{F-\delta,\,u_+}(x,\,\pi(u))}
\;d\mu_{x}(u_+)}\\ &  \hspace*{3.2cm}
\mathbbm{1}_{\F_\Ga}(w)\,\wt f(w)\, 
\frac{e^{C_{F-\delta,\,w_+}(x,\,\pi(w))+
C_{F\circ \iota-\delta,\,v_-}(x,\,\pi(w))}\,
d\mu_{x}(w_+)d\mu^\iota_{x}(v_-)dt}  {\|m_F\|}\\
  =&\int_{w\in T^1\wt M} \; \mathbbm{1}_{\F_\Ga}(w)\,\wt f(w)\,
 \frac{dm_F(w)}{\|m_F\|}=\int_{T^1M}  f\,dm^F\;.
\end{align*}

\noindent
This is exactly the required disintegration property, and ends the
proof of Lemma \ref{lem:explicitconditionelles}.  
\cqfd

\medskip 
\noindent{\bf Step 3. } In this penultimate step, we prove that for
every $\phi$-ergodic $\phi$-invariant probability measure $m$ on
$T^1M$ for which the negative part of $F$ is integrable, its metric
pressure $P_{\Ga,\,F}(m)$ is at most equal to the critical exponent
$\delta=\delta_{\Ga,\,F}$, with equality if and only if the given
Gibbs measure $m_F$ is finite, in which case $m$ is the normalisation
of $m_F$ to a probability measure. The interplay between these two
measures $m_F$ and $m$ is crucial in the following result.

\blemm\label{lem:steptrois} Let $m$ be a $\phi$-ergodic
$\phi$-invariant probability measure on $T^1M$ such that $\int_{T^1M}
F_-\,dm<+\infty$. Let $\tau>0$ be such that $\phi_\tau$ is ergodic
with respect to $m$, and let $\zeta$ be a partition associated with
$(m,\tau)$ by Proposition \ref{prop:partition}.

(1)  The real number
$$
G(v)=-\ln \int_{\zeta(v)}e^{-c_F(v,\,w)}d\mu_{W^{\rm su}(v)}(w)
$$ 
and the measure $m^F_{\zeta(v)}$ given by Equation
\eqref{eq:condmesuputative} are defined for $m$-almost every $v$ in
$T^1M$, and we have, for $m$-almost every $v$ in $T^1M$,
\begin{equation}\label{eq:formuleavecG}
-\ln m^F_{\zeta(v)}((\phi_\tau^{-1}\zeta)(v))=
\tau\delta-\int_0^{\tau } F(\phi_tv)\;dt \;+G(\phi_\tau v)-G(v)\;.
\end{equation}

(2) We have
\begin{equation}\label{eq:formulesansG}
\int_{T^1M}-\ln m^F_{\zeta(v)}((\phi_{\tau}^{-1}\zeta)(v))\;dm(v)=
\tau\delta-\tau \int_{T^1M} F\;dm\;.
\end{equation}

(3) If $\int_{T^1M} F_-\;dm_F< +\infty$ and if $m_F$ is
finite, then, with $m^F=\frac{m_F}{\|m_F\|}$, we have
$$
P_{\Ga,\,F}(m^F)=\delta\;.
$$

(4) We have $P_{\Ga,\,F}(m)\leq \delta$, hence
$$
P(\Ga,F)\leq\delta\;.
$$

(5) If $P_{\Ga,\,F}(m)=\delta$, then $m_F$ is finite and $m=m^F$, where
$m^F=\frac{m_F}{\|m_F\|}$.  
\elemm

\dem We will follow the same scheme of proof as the one for $F\equiv
0$ in \cite{OtaPei04}.

\medskip (1) Note that the measure $\mu_{W^{\rm su}(v)}$ is defined
for every $v\in T^1M$. By the disintegration properties of $\wt m$
proved in Subsection \ref{subsec:condmesgibbs}, for $m$-almost every
$v$, the element $v$ belongs to the support of $\mu_{W^{\rm su}(v)}$.
Furthermore, $\zeta(v)$ is a relatively compact neighbourhood of $v$
in $W^{\rm su}(v)$ for $m$-almost every $v$, since $\zeta$ is
$m$-subordinated to the strong unstable foliation.  Hence the measure
$$
dm^F_{\zeta(v)}(w)=
\frac{\mathbbm{1}_{\zeta(v)}(w)}{\int_{\zeta(v)} e^{-c_F(w,\,u)} \,
  d\mu_{W^{\rm su}(v)}(u)}\; d\mu_{W^{\rm su}(v)}(w)
$$ 
is well defined for $m$-almost every $v$. (Note that if $m_F$ is
finite and if $\zeta$ is also $m^F$-measurable and $m^F$-subordinated
to $\W^{\rm su}$, then $m^F_{\zeta(v)}$ is, by Lemma
\ref{lem:explicitconditionelles}, the conditional measure of $m^F$ on
$\zeta(v)$, but we are not assuming $m_F$ to be finite in Assertion
(1), (2) or (4).)

Recall that $(\phi_\tau^{-1}\zeta)(v)$ is contained in $\zeta(v)$ for
$m$-almost every $v$ by Proposition \ref{prop:partition}, and that by
definition 
$$
(\phi_\tau^{-1}\zeta)(v)=\phi_{-\tau}(\zeta(\phi_\tau v))\;.
$$
Let $w\in\phi_\tau W^{\rm su}(v)$ and $u\in W^{\rm su}(v)$. By
Equation \eqref{eq:cocyclebasvsflogeod}, we have
$$
c_F(\phi_{-\tau}w,v)=c_F(w,\phi_{\tau}v)+
\int_0^{-\tau}\big(F(\phi_tw) -F(\phi_{t+\tau} v)\big)\;dt\;.
$$
By Equation \eqref{eq:propdemesfortinstabl}, we have
$$
d\mu_{W^{\rm su}(v)}(\phi_{-\tau}w)=
e^{\int_{0}^{\tau} (F(\phi_{s-\tau}w)-\delta)\,dt}\;
d\mu_{W^{\rm su}(\phi_{\tau}v)}(w)=
e^{\int_{-\tau}^0 (F(\phi_tw)-\delta)\,dt}\;
d\mu_{W^{\rm su}(\phi_{\tau}v)}(w)\;.
$$
By the cocycle property, we have $c_F(\phi_{-\tau}w,u)=
c_F(\phi_{-\tau}w,v) +c_F(v,u)$ and $c_F(w,\phi_{\tau}v)=
-c_F(\phi_{\tau}v,w)$. Hence for $m$-almost every $v$, we have
\begin{align*}
&m^F_{\zeta(v)}((\phi_{\tau}^{-1})\zeta(v))\\ =\;&
\int_{\phi_{-\tau}(\zeta(\phi_\tau(v)))} \frac{d\mu_{W^{\rm su}(v)}(w')}
{\int_{\zeta(v)}e^{-c_F(w',\,u)}\,d\mu_{W^{\rm su}(v)}(u)}\\ =\;&
\int_{\zeta(\phi_\tau(v))} \frac{e^{c_F(\phi_{-\tau}w,\,v)}\;
d\mu_{W^{\rm su}(v)}(\phi_{-\tau}w)}
{\int_{\zeta(v)}e^{-c_F(v,u)}\,d\mu_{W^{\rm su}(v)}(u)} \\ =\;&
\frac{\int_{\zeta(\phi_\tau(v))} e^{c_F(w,\,\phi_{\tau}v)+
\int_0^{-\tau}(F(\phi_tw)-F(\phi_{t+\tau}v))\;dt+
\int_{-\tau}^0 (F(\phi_tw)-\delta)\,dt}
\;d\mu_{W^{\rm su}(\phi_{\tau}v)}(w)}
{\int_{\zeta(v)}e^{-c_F(v,\,u)}\,d\mu_{W^{\rm su}(v)}(u)}\\ =\;&
e^{-\tau\delta-\int_0^{-\tau } F(\phi_{t+\tau}v)\;dt \;-G(\phi_\tau v)+G(v)}\;.
\end{align*}
This proves Assertion (1).

\medskip (2) We are going to use the following quite classical result
that we state without proof (see for instance \cite[Lemme
8]{OtaPei04}).

\medskip\noindent{\bf Claim. } {\sl Let $T$ be a measure-preserving
  transformation of a probability space $(X,\A,\mu)$, let $H:X\ra\RR$
  be a measurable map such that $h=H\circ T-H$ has integrable negative
  part. Then $h$ is integrable and $\int_X h\;dm=0$.} 

\medskip By Assertion (1), for $m$-almost every $v\in T^1M$, since
$m^F_{\zeta(v)}$ is a probability measure, we have
$$
G\circ \phi_\tau(v) -G(v)\geq 
-\tau\delta+\int_0^{\tau } F(\phi_{t}v)\;dt\;.
$$ 
Since the negative part of $F$ is $m$-integrable, by Fubini's theorem
and the $\phi$-invariance of $m$, the map $g=G\circ \phi_\tau -G$
hence has $m$-integrable negative part.  Applying the above claim, we
may integrate Equation \eqref{eq:formuleavecG} over $v\in T^1M$ for
the measure $m$. The vanishing of $\int_{T^1M}g\;dm$ then yields
Equation \eqref{eq:formulesansG}.

\medskip (3) If $m_F$ is finite, then its normalised measure $m^F$ is
an ergodic probability measure on $T^1M$ (see Corollary
\ref{coro:uniqGibstate}). Hence there exists $\tau>0$ such that
$\phi_\tau$ is ergodic for $m^F$ (see Lemma
\ref{lem:ergotflowtempstau}). If $\int_{T^1M} F_-\;dm_F< +\infty$,
applying Equation \eqref{eq:formulesansG} to $m=m^F$ and
$\zeta=\zeta'$ a $m^F$-measurable partition associated with
$(m^F,\tau)$ by Proposition \ref{prop:partition}, we have, by the
definition of the entropy of $\phi_\tau$ relative to $\zeta'$ and
since $(\phi_{\tau}^{-1}\zeta')(v)$ is contained in $\zeta'(v)$ for
$m^F$-almost every $v$,
$$
h_{m^F}(\phi_\tau,\zeta')= \int_{v\in T^1M}-
\ln m^F_{\zeta'(v)}((\phi_{\tau}^{-1}\zeta')(v))\;dm^F(v)=
\tau\delta-\tau \int_{T^1M} F\;dm^F\;.
$$
By the definition of the entropy of $\phi$ with respect to $m^F$ and by
Proposition \ref{prop:partition}, we have
$$
h_{m^F}(\phi)=\frac{1}{\tau}h_{m^F}(\phi_\tau)=
\frac{1}{\tau}h_{m^F}(\phi_\tau,\zeta')=\delta-\int_{T^1M} F\;dm^F\;.
$$
This gives $P_{\Ga,\,F}(m^F)=\delta$.

\medskip (4) 
For $m$-almost every $v\in T^1M$, by Assertion (1),
let us define $\psi(v)\in\mathopen{[}0,+\infty\mathclose{]}$ by
$$
\psi(v)= \frac{m^F_{\zeta(v)}\big((\phi_\tau^{-1}\zeta)(v)\big)}
{m_{\zeta(v)}\big((\phi_\tau^{-1}\zeta)(v)\big)}
\quad \mbox{if} \quad 
m_{\zeta(v)}\big((\phi_\tau^{-1}\zeta)(v)\big)>0\,,\quad
\mbox{and}\quad \psi(v)=+\infty\quad \mbox{otherwise}\,.
$$

\medskip\noindent{\bf Key claim. (Otal-Peigné)} {\sl The map $\psi$ is
  $m$-measurable, the maps $\psi$ and $\ln \psi$ are $m$-integrable,
  and $\int_X \psi \;dm \leq 1$.} 

\medskip
\dem The proof when $F=0$ in the key result \cite[Fait 9]{OtaPei04} extends
immediately, by replacing $\mu_v^+$ in that reference by
$\mu_{W^{\rm su}(v)}$, and by noting, to prove the measurability of $\psi$
as indicated, that $v\mapsto \mu_{W^{\rm su}(v)}$ is weak-star continuous
as already mentioned, hence $v\mapsto \mu_{W^{\rm su}(v)}(A)$ is
measurable for every Borel subset $A$ of $T^1M$.  \cqfd

\medskip
By the definition of the (metric) entropy and by Proposition
\ref{prop:partition}, we have as above
$$
\tau h_m(\phi)=h_m(\phi_\tau)=h_m(\phi_\tau,\zeta)=
\int_{T^1M}-\ln m_{\zeta(v)}((\phi_\tau^{-1}\zeta)(v))\;dm(v)\;.
$$
Hence by Equation \eqref{eq:formulesansG},
$$
\int_{T^1M}\ln\psi\;dm=-\tau\delta+\tau\int_{T^1M}F\;dm+\tau h_m(\phi)\;.
$$
By Jensen's inequality and by the above Claim, we have
$$
\int_{T^1M}\ln\psi\;dm\leq\ln\Big(\int_{T^1M}\psi\;dm\Big)\leq 0\;.
$$
Hence
$$
P_{\Ga,\,F}(m)= h_m(\phi)+\int_{T^1M}F\;dm\leq \delta\;,
$$
as required.

\medskip The last claim of Assertion (4) follows by the definition of
$P(\Ga,F)$ (with the restriction to ergodic probability measures)
since $m$ is arbitrary and the existence of $\tau>0$ such that
$\phi_\tau$ is ergodic for $m$ is guaranteed by Lemma
\ref{lem:ergotflowtempstau}.

\medskip (5) The equality case in Jensen's inequality implies that
$\psi(v)=1$ for $m$-almost every $v\in T^1M$. Hence the measure
$m^F_{\zeta(v)}$ and the conditional measure $m_{\zeta(v)}$ coincide
on the $\sigma$-algebra generated by the restriction to $\zeta(v)$ of
$\phi_\tau^{-1}\zeta$ for $m$-almost every $v\in T^1M$.  Similarly by
replacing $\phi_\tau$ by $\phi_{k\tau}$, these measures coincide on
the $\sigma$-algebra generated by the restriction to $\zeta(v)$ of
$\phi_\tau^{-k}\zeta$ for $m$-almost every $v\in T^1M$. Since $\zeta$
is generating under $\phi_\tau$ (see Proposition
\ref{prop:partition}), these measures are equal, for $m$-almost every
$v\in T^1M$.

First assume that $m_F$ is finite. In particular the normalised
measure $m^F$ is ergodic by Corollary \ref{coro:uniqGibstate}. To
prove that $m$ and $m^F$ coincide, let us show, using Hopf's argument,
that $m(f)=m^F(f)$ for every $f:T^1M\ra\RR$ continuous with compact
support. Let $A_{f}$ be the set of $v\in T^1M$ such that the limit
$\lim_{t\ra+\infty} \frac{1}{t}\int_0^t f(\phi_sv)\;ds$ exists and is
equal to $m^F(f)$.  By Birkhoff's ergodic theorem, $A_{f}$ has full
measure with respect to $m^F$. By the uniform continuity of $f$, the
set $A_{f}$ is saturated by the stable foliation. By the quasi-product
structure (see the comment after Equation
\eqref{eq:mesGibbmesstabinstab}), for every $v\in T^1M$, the
intersection $A_{f}\cap W^{\rm su}(v)$ has full measure for
$\mu_{W^{\rm su}(v)}$. By the definition of $m^F_{\zeta(v)}$, this
implies that the intersection $A_{f}\cap \zeta(v)$ has full measure
with respect to $m^F_{\zeta(v)}$. By the above equality of measures,
we have that $A_{f}\cap \zeta(v)$ has full measure with respect to
$m_{\zeta(v)}$, for $m$-almost every $v$. By definition of the
conditional measures of $m$, this implies that $A_{f}$ has full
measure with respect to $m$. By Birkhoff's ergodic theorem applied
this time to $m$, we have $m(f)=m^F(f)$, as required.

Let us prove that $m_F$ is indeed finite. By the
Hopf-Tsuji-Sullivan-Roblin theorems \ref{theo:critnonergo} and
\ref{theo:critergo}, the dynamical system $(T^1M, (\phi_t)_{t\in\RR},
m_F)$ is either completely dissipative or conservative.

Assume that the first case holds. Let $f:T^1M\ra \RR$ be Lipschitz
continuous with compact support. Then for $m_F$-almost every $v\in
T^1M$, we have $\lim_{t\ra+\infty} f(\phi_sv)\;ds=0$. Considering the
set
$$
A'_{f}= \{v\in T^1M\;:\;\lim_{t\ra+\infty} \frac{1}{t}\int_0^t
f(\phi_sv)\;ds=0\}\;,
$$ 
which has full measure with respect to $m_F$, we have as above that
$m(f)=0$ for every such $f:T^1M\ra\RR$.  This contradicts the fact
that $m$ is a probability measure.

Assume that the second case holds. Then there exists a positive
Lipschitz continuous map $\rho: T^1M\ra \RR$ with
$\|\rho\|_{\LL^1(m_F)}=1$ such that, for $m_F$-almost every $v\in
T^1M$, we have 
$$
\int_0^{\infty} \rho\circ\phi_{-t}(v) dt=
\int_0^{\infty}\rho\circ\phi_t(v) dt= +\infty
$$ 
(see for instance the proof of Assertion (h) of Proposition
\ref{prop:listequivHTSRG}). The same argument as above, this time
using Hopf's ratio ergodic theorem, gives that
$\frac{m(f)}{m(\rho)}=\frac{m_F(f)}{m_F(\rho)}$. This implies that
$m_F$ is finite, since $m$ is a probability measure.  \cqfd

\bigskip 
\noindent{\bf Step 4. } In this last step of the proof of Theorem
\ref{theo:principevariationnel}, we prove the converse inequality
$\delta\leq P(\Ga,F)$.  Assertion (4) of Lemma
\ref{lem:steptrois} then implies that $\delta= P(\Ga,F)$. If $m_F$ is
finite and $\int_{T^1M} F_-\,dm_F< +\infty$, then the assertions (3)
and (5) of Lemma \ref{lem:steptrois} imply that the normalised Gibbs
measure $m^F$ is the unique equilibrium state. Otherwise, Assertion
(5) of Lemma \ref{lem:steptrois} also implies that there is no
equilibrium state (since $m=m^F$ implies that the negative part of the
potential $F$ is integrable with respect to $m^F$ if it is with
respect to $m$). Therefore Theorem \ref{theo:principevariationnel}
follows from the following result.

\blemm\label{lem:minopression} We have $\delta_{\Ga,\,F}\leq P(\Ga,F)$.  
\elemm

\dem We start the proof by recalling an alternative definition of
the topological pressure of a flow on a compact metric space.

Let $X$ be a compact metric space whose distance is denoted by $d$,
let $\varphi= (\varphi_{t})_{t\in\RR}$ be a continuous flow of
homeomorphisms of $X$, and let $G:X\ra \RR$ be a continuous
map.

For all $T\geq 0$ and $\epsilon>0$, a subset $E$ of $X$ is {\it
  $(\epsilon,T)$-separated}%
\index{separated@$(\epsilon,T)$-separated!subset of a dynamical metric
  space} if for all distinct points $x$ and $y$ in $E$, there exists
$t\in\mathopen{[}0,T\mathclose{]}$ such that $d(\varphi_tx,\varphi_ty)
\geq \epsilon$. Define $P_{d,\,\epsilon,\,T,\,G}$ as the upper bound
of
$$
\sum_{x\in E} e^{\int_{0}^T G(\varphi_tx)\;dt}
$$
on all $(\epsilon,T)$-separated subsets $E$ of $X$, and the {\it
  combinatorial pressure}\index{combinatorial
  pressure}\index{pressure!combinatorial } of $(G,\varphi)$ as
$$
P_X(G,\varphi)=
\sup_{\epsilon>0}\;\limsup_{T\ra+\infty}\;
\frac{1}{T}\;\ln P_{d,\,\epsilon,\,T,\,G}\;.
$$

The following Variational Principle says that the combinatorial
pressure and the topological pressure coincide. We denote by
$\M(X,\varphi)$ the weak-star compact convex set of probability
measures on $X$ invariant by the flow $\varphi$.

\btheo (see for instance \cite[Th 9.10]{Walters82}) 
\label{theo:dpression} We have
$$
P_X(G,\varphi)=
\sup_{m\in\M(X,\,\varphi)}\;\Big(h_m(\varphi)+\int_XG\;dm\Big)
\;.
$$
\etheo

Note that this result is stated for transformations in
\cite[Th 9.10]{Walters82}, but its proof extends to flows.

\medskip Let us resume the proof of Lemma \ref{lem:minopression}. We
may assume that $\delta=\delta_{\Ga,\,F}>0$. We are going to
construct, using Subsection \ref{subsec:semigroup}, for every
$\delta''\in\mathopen{]}0,\delta\mathclose{[}\,$, a compact subset $K$
of $T^1M$ invariant by the geodesic flow, such that the combinatorial
pressure of $(F_{\mid K},\phi_{\mid K})$ on $K$ satisfies
$$
P_K(F_{\mid K},\phi_{\mid K})\geq \delta''\;.
$$
(We endow $K$ with the distance $d$ which is the restriction to $K$ of
the quotient of the Riemannian distance of $T^1\wt M$ (for the Sasaki
metric).) Since $\M(K,\varphi)$ is contained in $\M_F(T^1M)$ (by
considering measures on $K$ as measures on $T^1M$ with support in
$K$), this implies that $P(\Ga,F)\geq \delta_{\Ga,\,F}$ as required in
Lemma \ref{lem:minopression}.

\medskip Let us now define the constants and objects we will use in
the proof of Lemma \ref{lem:minopression}. Let $\delta'\in\mathopen{]}
\delta'', \delta\mathclose{[}$ and $s\in\mathopen{]}\delta'',
\delta'\mathclose{[}\,$.  Let $L\geq 5$ (which will be chosen large
enough) and $\epsilon'\in\mathopen{]}0,\frac{1}{12}\mathclose{[}$
(which will be chosen small enough).  Let $\epsilon\in\mathopen{]}0,
\frac{1}{16}\mathclose{]}$ (depending only on $\epsilon'$ and tending
to $0$ as $\epsilon'\ra 0$) be given by Lemma
\ref{lem:crithyperboelem} with $\ell=3$. Let $\theta>0$ (depending
only on $\epsilon$) be defined by Proposition
\ref{prop:schottkysemigroups}.  Let $C=6$. Let $x_0\in\wt M$, $v_0\in
T^1_{x_0}\wt M$, $N>L$ and $S$ a finite subset of $\Ga$ (depending on
$\delta',\theta,C,L$) be given by Proposition
\ref{prop:constructsemigroup}. Let $G$ be the subsemigroup of $\Ga$
generated by $S$.

By Proposition \ref{prop:schottkysemigroups} (whose assumptions are
satisfied by the assertions of Proposition
\ref{prop:constructsemigroup}), the semigroup $G$ is free on $S$ and
for every $\ga\in G$ of nonzero length $\ell$ as a word in the
elements of $S$, we have 
$$
d(\ga x_0,\mathopen{[}x_0,\ga^2 x_0\mathclose{]})\leq
\epsilon\;\;\;{\rm and} \;\;\; d(x_0,\ga x_0)\geq
(N-2\epsilon)\ell\geq 3\;.
$$
Hence by Lemma \ref{lem:crithyperboelem}, every nontrivial element
$\ga$ of $G$ is loxodromic and if $x_\ga$ is the closest point to
$x_0$ on the translation axis $\Axe_\ga$ of $\ga$, then
$d(x_0,x_\ga)\leq \epsilon'$ (see the picture above Lemma
\ref{lem:crithyperboelem}). Let $v_\ga$ be the unit tangent vector at
$x_\ga$ pointing towards $\ga x_\ga$, which is tangent to the
translation axis $\Axe_\ga$ (see the picture below).

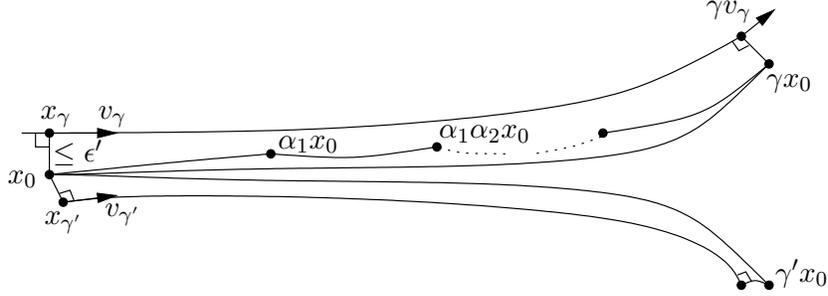
\begin{center}
\input{fig_minopression.pstex_t}
\end{center}

Let $\prT: T^1\wt M\ra T^1M=\Ga\bs T^1\wt M$ be the canonical
projection. For every $\eta>0$, denote by $\overline{\N}_\eta (A)$ the
closed $\eta$-neighbourhood of a subset $A$ of $\wt M$. Consider the
$\Ga$-invariant subset $\wt K$ of $\wt M$ defined by
$$
\wt K=\big\{v\in T^1\wt M\;:\;\forall\;t\in\RR,\;\;
\pi(\phi_tv)\in\Ga\bigcup_{\alpha\in S}
\overline{\N}_{\epsilon+\epsilon'} (\mathopen{[}x_0,\alpha x_0
\mathclose{]})\;\big\}\;.
$$
Its image $K=\prT(\wt K)$, which consists in the elements of $T^1M$
whose orbit under the geodesic flow stays at distance at most
$\epsilon+\epsilon'$ of $\bigcup_{\alpha\in S} \prT(\mathopen{[}x_0,
\alpha x_0\mathclose{]})$, is compact, since $S$ is finite. Note that
$K$ (which depends on $L$ and $\epsilon'$) is invariant under the
geodesic flow, by construction.

\medskip Now, to finish the proof, we are going to construct, for a
fixed $\epsilon''>0$ small enough, an increasing sequence
$(n_k)_{k\in\NN}$ in $\NN$ and for every $k\in\NN$, a $(\epsilon'',
n_k)$-separated subset of $K$ with a good control of the integral of
the potential on its flow lines.

For every nontrivial element $\ga\in G$, let $\ell\in\NN-\{0\}$ and
$\alpha_1,\dots,\alpha_\ell \in S$ be such that $\ga=\alpha_1\dots
\alpha_k$. Note that $\mathopen{[}x_0,\ga x_0\mathclose{]}$ is
contained in
$$
\overline{\N}_{\epsilon}(\mathopen{[}x_0,\alpha_1 x_0\mathclose{]}\cup
\mathopen{[}\alpha_1x_0,\alpha_1\alpha_2 x_0\mathclose{]} \cup \dots\cup
\mathopen{[}\alpha_1\dots\alpha_{k-1}x_0,\ga x_0\mathclose{]})
$$ 
by Proposition \ref{prop:schottkysemigroups}.

We hence have, for every $t\in\RR$,
$$
\pi(\phi_tv_\ga)\in \Axe_\ga=
\ga^\ZZ\mathopen{[}x_\ga,\ga x_\ga\mathclose{]}\subset
\Ga\,\overline{\N}_{\epsilon'}(\mathopen{[}x_0,\ga x_0\mathclose{]})
\subset \Ga\bigcup_{\alpha\in S} \overline{\N}_{\epsilon+\epsilon'} 
(\mathopen{[}x_0,\alpha x_0\mathclose{]})\;.
$$
This proves that $v_\ga$ belongs to $\wt K$ for every nontrivial
$\ga$ in $G$.

If $\epsilon'$ (hence $\epsilon$) is small enough and $L$ is large
enough, the proof of Theorem \ref{theo:critexpsubsemigroup} says that
the critical exponent $\delta_{G,\,F}$ of the subsemigroup $G$ is at
least $s$. In particular, since $\delta''<s$ and by the definition of
$\delta_{G,\,F}$, there exists an increasing sequence
$(n_k)_{k\in\NN}$ in $\NN$ such that for every $k\in\NN$, we have
\begin{equation}\label{eq:minoexpocritreferee}
\sum_{\ga\in G,\;n_k-1<d(x_0,\ga x_0)\leq n_k} e^{\int_{x_0}^{\ga x_0}\wt F}
\geq e^{\delta''\,n_k}\,.
\end{equation}
Let us define 
$$
A_k=\{\ga\in G,\;n_k-1<d(x_0,\ga x_0)\leq n_k\}\,.
$$

Furthermore, if $\epsilon'$ (hence $\epsilon$) is small enough, by
Lemma \ref{lem:technicholder}, there exists $c\geq 0$ (depending only
on $\epsilon$, on the H\"older constants of $\wt F$ and on
$\max_{\pi^{-1}(B(x_0,\,1))}\;|\wt F |$) such that for every
nontrivial element $\ga\in G$,
\begin{equation}\label{eq:deficonstantc}
\Big|\int_{x_\ga}^{\ga x_\ga}\wt F\;-\int_{x_0}^{\ga x_0}\wt F\;\Big|\leq c\;.
\end{equation}

We will need the following two lemmas.

\blemm \label{lem:sortirdescouronnes} For every $k\in\NN$, for all
distinct $\ga,\ga'$ in $A_k$, the element $\ga$ is not an initial
subword of $\ga'$.  
\elemm

\dem Assume for a contradiction that $\ga\in G$ is an initial subword
of $\ga'\in G$, different from $\ga'$, and that the distances from
$\ga x_0$ and $\ga' x_0$ to $x_0$ belong to $\mathopen{]}n_k-1,
n_k\mathclose{]}$. Let $\ga_1$ be a nontrivial word in $S$ such that
$\ga'=\ga\ga_1$. By Proposition \ref{prop:schottkysemigroups}, the
point $\ga x_0$ is at distance at most $\epsilon$ from the geodesic
segment $\mathopen{[}x_0,\ga' x_0\mathclose{]}$. Hence by the triangle
inequality
\begin{align*}
d(x_0,\ga_1x_0)& =d(\ga x_0,\ga\ga_1x_0) \leq 
d(x_0,\ga' x_0)-d(x_0,\ga x_0)+2\,\epsilon
\\ & \leq n_k-(n_k-1)+2\,\epsilon=1+2\,\epsilon\,.
\end{align*}
By Proposition \ref{prop:schottkysemigroups}, since the word-length
$\ell$ of $\ga_1$ is at least $1$, we have $d(x_0,\ga_1x_0)\geq
N-2\epsilon$.  Hence $L<N\leq 1+4\epsilon\leq 5$, a contradiction.
\cqfd

\blemm For every $\epsilon''\in\mathopen{]}0,\frac{1}{2} -6\epsilon'
\mathclose{[}\,$, for all $k\in\NN$ and $\ga\neq \ga'$ in $A_k$, there
exists $t\in \mathopen{[}0,n_k\mathclose{]}$ such that
$d(\pi(\phi_{t}v_\ga),\pi(\phi_{t}v_{\ga'})) \geq \epsilon''$.  
\elemm

\smallskip
\noindent
\begin{minipage}{8.5cm}
  \dem Since $\epsilon\leq \frac{1}{16}$, by the last claim of
  Proposition \ref{prop:schottkysemigroups} and by the previous Lemma
  \ref{lem:sortirdescouronnes}, there exists $x\in \mathopen{[}x_0,\ga
  x_0\mathclose{]}$ and $x'\in \mathopen{[}x_0,\ga' x_0\mathclose{]}$
  such that $d(x_0,x)=d(x_0,x')$ and $d(x,x')\geq \frac{1}{2}$. Let
  $y$ and $y'$ be the closest points to $x$ and $x'$ on
  $\mathopen{[}x_{\ga},\ga x_{\ga}\mathclose{]}$ and
  $\mathopen{[}x_{\ga'},\ga' x_{\ga'}\mathclose{]}$ respectively,
  which satisfy $d(x,y),d(x',y')\leq \epsilon'$ by convexity.  Let
  $t=d(x_0,x)-2\epsilon'$ and $z=\pi(\phi_{t}v_\ga)$ and
  $z'=\pi(\phi_{t}v_{\ga'})$.
\end{minipage}
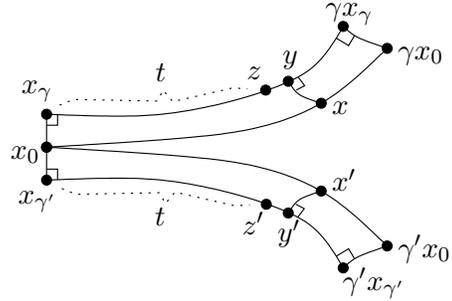
\begin{minipage}{6cm}
\begin{center}
\input{fig_separdynam.pstex_t}
\end{center}
\end{minipage}

\medskip We have $t\geq 0$ since otherwise $\frac{1}{2}\leq d(x,x')
\leq d(x,x_0)+ d(x_0,x')\leq 4\,\epsilon'$, which contradicts the fact
that $\epsilon'\leq \frac{1}{12}$. We have $t\leq d(x_0,x)\leq
d(x_0,\ga x_0)\leq n_k$.  

By the triangle inequality, and since closest point maps do not
increase the distances, we have
$$
d(x_\ga,z)=t=d(x_0,x)-2\epsilon'\leq d(x_\ga,y)\leq d(x_0,x)\;.
$$
Hence 
\begin{align*}
d(z,x)&\leq d(z,y)+d(y,x)=(d(x_\ga,y)-d(x_\ga,z))+d(y,x)\\ &\leq 
d(x_0,x)-(d(x_0,x)-2\epsilon')+\epsilon'=3\epsilon'\;.
\end{align*}
Similarly, $d(z',x')\leq 3\epsilon'$, hence by the
triangle inequality
$$
d(z,z')\geq d(x,x')-d(x,z)-d(x',z')\geq \frac{1}{2}-6\epsilon'
\geq \epsilon''\;.\;\;\; \mbox{$\Box$}
$$

\medskip Since $d_S(\phi_{t}v_\ga,\phi_{t}v_{\ga'})\geq
d(\pi(\phi_{t}v_\ga),\pi(\phi_{t}v_{\ga'}))$ by Equation
\eqref{eq:comparpidistSasaki}, this lemma implies that the set $\wt
E_k=\{v_\ga\;:\;\ga\in A_k\}$ is $(\epsilon'', n_k)$-separated.

\medskip
Let $C_K= \inf_{\wt K} \wt F$, which is finite since $F$ is continuous
on the compact subset $K$. With $\ell(\ga)=d(x_\ga,\ga x_\ga)$ the
translation length of $\ga\in A_k$, we have
$$
n_k-1-2\epsilon'\leq  d(x_0,\ga x_0)-2\epsilon'
\leq \ell(\ga)\leq d(x_0,\ga x_0)\leq n_k\;.
$$
Since $v_\ga\in \wt K$ and $\wt K$ is invariant under the geodesic
flow, we hence have, for all $k\in\NN$ and $\ga\in A_k$,
$$
\int_{\ell(\ga)}^{n_k}\wt F(\phi_tv_\ga) \; dt
\geq C_K(1+2\epsilon')\,.
$$

Since $K$ is compact and since the isometric action of $\Ga$ on $\wt
M$ is proper, there exists $\epsilon''\in\mathopen{]}0,\frac{1}{2}
-6\epsilon'\mathclose{[}$ (this interval is nonempty since $\epsilon'<
\frac{1}{12}$) and $N_K>0$ (when $\Ga$ is torsion free, we may take
$N_K=1$ and $\epsilon''>0$ less than the lower bound on $\pi(K)$ of
the injectivity radius of $M$) such that $\wt E_k$ may be subdivided
in at most $N_K$ subsets that the map $\prT:T^1\wt M\ra T^1M$ inject
into $(\epsilon'', n_k)$-separated subsets. For at least one of these
sets, let us call it $E_k$, we have, by Equations
\eqref{eq:deficonstantc} and \eqref{eq:minoexpocritreferee}
\begin{align*}
\sum_{v\in E_k} e^{\int_0^{n_k} F(\phi_tv)\;dt} &\geq \frac{1}{N_K}
\sum_{\ga\in A_k} e^{\int_0^{n_k}\wt F(\phi_tv_\ga)\;dt} \geq
\frac{e^{C_K(1+2\epsilon')}}{N_K}
\sum_{\ga\in A_k} e^{\int_{x_\ga}^{\ga x_\ga}\wt F}\\ & \geq
\frac{e^{C_K(1+2\epsilon') -c}}{N_K}
\sum_{\ga\in A_k} e^{\int_{x_0}^{\ga x_0}\wt F}\geq
\frac{e^{C_K(1+2\epsilon') -c}}{N_K}\;e^{\delta''\;n_k}\;.
\end{align*}
Taking the logarithm, dividing by $n_k$ and letting $k\ra+\infty$,
we have
$$
P_K(F_{\mid K},\phi_{\mid K})\geq 
\delta''\;.
$$
This is what we wanted to prove.  \cqfd

\section{The Liouville measure as a Gibbs measure}
\label{sec:Liouville}

Let $(\wt M,\Ga)$ be as in the beginning of Chapter
\ref{sec:negacurvnot}: $\wt M$ is a complete simply connected
Riemannian manifold, with dimension at least $2$ and pinched sectional
curvature at most $-1$ and $\Ga$ is a non-elementary discrete group of
isometries of $\wt M$.

The aim of this chapter is to investigate when the Liouville measure
of the geodesic flow on $T^1 M= \Ga\bs T^1\wt M$ coincides (up to a
multiplicative constant) with a Gibbs measure associated to some
H\"older potential.  This result is known when $M$ is compact, but we
provide a complete short proof in Subsection
\ref{subsec:liouvequilstate}.  Finding a necessary and sufficient
condition in the noncompact case is still open. The approach that we
follow in Subsection \ref{subsec:Gibbsapproach} actually allows us to
get more than what was previously known on the subject.  In
particular, we obtain that if $\wt M$ admits a discrete and cocompact
group of isometries, then the Liouville measure on $T^1M$ is a Gibbs
measure as soon as it is conservative (see Theorem
\ref{theo:LiouvilleisGibbs} and its following comment), but not
necessarily finite.

\medskip Recall that the {\it Liouville measure}%
\index{measure!of Liouville}\index{Liouville measure} $\vol_{T^1\wt
  M}$ on $T^1\wt M$ is the Riemannian volume of the Sasaki metric on
$T^1\wt M$ (see Subsection \ref{subsec:unittanbun} for a definition).
Equivalently, it disintegrates with respect to the projection $\pi:
T^1\wt M\ra \wt M$, with measure on the basis $\wt M$ the Riemannian
volume $\vol_{\wt M}$ of $\wt M$ and conditional measures on the
fibres $T^1_x\wt M$ the Riemannian spherical measure $\vol_{T^1_x\wt
  M}$ of the unit sphere of the Euclidean space $T_x\wt M$: for every
$f\in\C_c(T^1\wt M;\RR)$, we have
$$
\int_{v\in T^1\wt M}\;f(v)\;d\vol_{T^1\wt M}(v)=
\int_{x\in \wt M} \int_{v\in T^1_x\wt M}\;f(v)\;d\vol_{T^1_x\wt M}(v)\;
d\vol_{\wt M}(x)\;.
$$
It is invariant under the action of the geodesic flow (and the action
of $\ZZ/2\ZZ$ by time reversal). It induces (see Subsection
\ref{subsec:pushmeas}) a measure $\vol_{T^1M}$ on the Riemannian
orbifold $T^1M=\Ga\backslash T^1\wt M$, called the {\it Liouville
  measure}\index{measure!of Liouville}%
\index{Liouville measure} on $T^1M$, invariant under the geodesic
flow (and time reversal), which disintegrates similarly over the
Riemannian orbifold $M=\Ga\bs \wt M$ with fibres the Riemannian
orbifolds $T^1_xM=\Ga_{\wt x}\bs T^1_{\wt x}\wt M$ where $\wt x$ is a
lift of $x$, whose stabiliser in $\Ga$ is denoted by $\Ga_{\wt x}$:
$$
d\vol_{T^1M}(v)=
\int_{x\in M} d\vol_{T^1_xM}(v)\;d\vol_{M}(x)\;.
$$
In particular, the Liouville measure $\vol_{T^1M}$ is finite if and
only if $M$ has finite volume: since the fixed point set of a
Riemannian isometry different from $\id$ of a connected Riemannian
manifold has Riemannian measure $0$, if $M$ has dimension $m$, then
$$
\Vol(T^1M)=\Vol(\SS^{m-1})\;\Vol(M)\;.
$$
It can happen that $\vol_{T^1M}$ is infinite and ergodic, see for
example \cite{Rees81}, where she proves that the Liouville measure on
a Galois cover with covering group $\ZZ^d$ of a compact hyperbolic
manifold is ergodic with respect to the geodesic flow if and only if
$d\leq 2$.

Let us now define the potential on $T^1\widetilde{M}$ with respect to
which the Liouville measure is a Gibbs measure in all known cases. Let
$v\in T^1\wt M$ and $t\in\RR$. Denote by $\gls{stronunstabtang}$
(respectively $\gls{stronstabtang}$) the tangent space $T_vW^{\rm su}(v)$
at $v$ to $W^{\rm su}(v)$ (respectively $T_vW^{\rm ss}(v)$ at $v$ to
$W^{\rm ss}(v)$) and let $\gls{neuttang}=\RR\frac{d}{dt}_{\mid t=0}
(\phi_tv)$ be the tangent space at $v$ of its geodesic flow line. We
have
$$
T_vT^1\wt M=E^{\,\rm su}(v)\;\oplus\;E^{0}(v) \;\oplus\; E^{\,\rm ss}(v)\;.
$$
Let us also define 
$$
\gls{unstabtang}=E^{\,\rm su}(v)\;\oplus\;E^{0}(v)
\;\;\;{\rm and}\;\;\;
\gls{stabtang}= E^{0}(v)\;\oplus\;E^{\,\rm ss}(v) \;.
$$
For every $v\in T^1M$, we again denote by $E^{\,\rm su}(v)$, $E^{0}(v)$,
$E^{\,\rm ss}(v)$ the images of $E^{\,\rm su}(\wt v)$, $E^{0}(\wt v)$,
$E^{\,\rm ss}(\wt v)$, respectively, by the canonical map $TT^1\wt M\ra
TT^1M=\Ga\bs TT^1\wt M$, where $\wt v$ is any lift of $v$ by the
canonical projection $T^1\wt M\ra T^1M$.

Recall that the {\it Jacobian}\index{Jacobian} of a smooth map $f:N\ra
N'$ between two Riemannian manifolds is the map $\operatorname{Jac}f:
N\ra [0,+\infty[$ such that $\operatorname{Jac}f(x)$ is the absolute
value of the determinant (of the matrix in any orthonormal bases) of
the linear map $T_xf$ between the two Euclidean spaces $T_xN$ and
$T_{f(x)}N'$, for every $x\in N$.  For all $t\in\RR$ and $v\in T^1\wt
M$, the map $\phi_t$ on $T^1\wt M$ sends the smooth submanifold
$W^{\rm su}(v)$ to the smooth submanifold $W^{\rm su} (\phi_t v)$. We
endow each leaf of $\W^{\rm su}$ with the Riemannian metric induced by
the Sasaki metric on $T^1\wt M$.  We can therefore define the Jacobian
of the restriction of $\phi_{t}$ to the smooth submanifold $W^{\rm su}
(v)$ of $T^1\wt M$ at the point $v$, and we denote it by $\wt J^{\,\rm
  su} (v,t)>0$. It is smooth in $t$ and continuous in $(v,t)$ (and
studying its regularity in $v$ will be one of the main points
below). The flow property of the geodesic flow and the invariance of
the strong unstable foliation under the geodesic flow give the
following cocycle relation: for all $s,t\in\RR$ and $v\in T^1\wt M$,
\begin{equation}\label{eq:cocyclunstabjac}
\wt J^{\,\rm su}(v,0)=1\;\;\;{\rm and}\;\;\;
\wt J^{\,\rm su}(v,s+t)=\wt J^{\,\rm su}(\phi_tv,s)\;\wt J^{\,\rm su}(v,t)\;.
\end{equation}

Now, we define the following map
\begin{equation}\label{Liouville-potential}
\gls{potentialsutilde}(v)= 
-{\frac{d}{dt}}_{|t=0} \log \wt J^{\,\rm su}(v,t)=
{\frac{d}{dt}}_{|t=0} \log \wt J^{\,\rm su}(\phi_t v,-t)\,.
\end{equation}
The maps $v\mapsto \wt J^{\,\rm su}(v,t)$ and $v\mapsto \wt F^{\,\rm su} 
(v)$, being $\Ga$-invariant, define maps $v\mapsto J^{\rm su}(v,t)$
and $v\mapsto \gls{potentialsu}(v)$ on $T^1M$ (which, when $\Ga$ is
torsion free, may be defined directly on $T^1M$). We have $\wt
F^{\,\rm su}\leq 0$ by the dilation properties of the geodesic flow on
the unstable foliation. An estimation of the gap map $D_{F^{\rm
    su}} :\wt M \times \partial_\infty^2\wt M\ra \RR$ (see Subsection
\ref{subsec:gap}) of the potential $F^{\,\rm su}$ would be interesting.

Note that $\wt J^{\,\rm su}(v,t)=e^{(n-1)t}$ for all $v\in T^1\wt M$ and
$t\in\RR$ if $M$ has dimension $n$ and constant sectional curvature
$-1$, so that $\wt F^{\,\rm su}$ is then the constant map with value
$-(n-1)$.

We will prove in this chapter the following results. 

\btheo\label{theo:potentialholder} Let $\wt M$ be a complete simply
connected Riemannian manifold, with dimension at least $2$ and pinched
negative sectional curvature. Assume that the derivatives of the
sectional curvature are uniformly bounded. Then $\wt F^{\rm su}$ is
H\"older-continuous, bounded, reversible and invariant under the full
isometry group of $\wt M$ .  
\etheo

Note that the new hypothesis of having bounded derivatives of
sectional curvatures is satisfied if there exists a compact set in
$\wt M$ whose images by the group $\Isom(\wt{M})$ cover $\wt M$, for
instance if there exists a cocompact discrete group of isometries of
$\wt M$.  

Since $\wt F^{\,\rm su}$ is bounded under the hypotheses of Theorem
\ref{theo:potentialholder}, the critical exponent
$\delta_{\Ga,\,F^{\rm su}}$ is finite (see Lemma
\ref{lem:elemproppressure} (iv)). We will then denote by $\wt
m_{F^{\rm su}}$ and $m_{F^{\rm su}}$ the Gibbs measures on $T^1\wt M$
and $T^1M$ associated with a fixed pair of Patterson densities of
dimension $\delta_{\Ga,\,F^{\rm su}}$ for $(\Ga, F^{\,\rm su}
\circ\iota)$ and $(\Ga,F^{\,\rm su})$.

\btheo\label{theo:LiouvilleisGibbs} Let $\wt M$ be a complete simply
connected Riemannian manifold, with dimension at least $2$ and pinched
negative sectional curvature. Assume that the derivatives of the
sectional curvature are uniformly bounded.

If the geodesic flow of $M$ is conservative with respect to the
Liouville measure and if $\delta_{\Ga,\,F^{\rm su}}\leq 0$, then the
Liouville measure is proportional to the Gibbs measure $m_{F^{\rm
    su}}$ of the potential $F^{\rm su}$. Furthermore, $(\Gamma,
F^{\,\rm su})$ is of divergence type and
$$
P(\Ga, F^{\,\rm su})=\delta_{\Ga,\,F^{\rm su}}=P_{Gur}(\Ga,F^{\,\rm su})=0\;.
$$ 
\etheo

The results of the article of Bowen-Ruelle \cite{BowRue75} for Axiom A
flows imply, in the case of the geodesic flow on the unit tangent
bundle of a compact negatively curved Riemannian manifold, that its
Liouville measure coincides with a multiple of its Gibbs measure
$m_{F^{\rm su}}$.  See in particular Remark 5.2 in loc.~cit., and
Subsection \ref{subsec:liouvequilstate}.

\medskip \noindent {\bf Remarks. }  (1) The inequality
$\delta_{\Ga,\,F^{\rm su}}\leq 0$ is satisfied when $M$ is compact, by
Ruelle's inequality (see Equation \eqref{eq:negatifinequelle} below,
using the fact that $\delta_{\Ga,\,F^{\rm su}}=P(\Ga,F^{\rm su})$ by the
Variational Principle theorem \ref{theo:principevariationnel}
(1)). Hence if $\wt M$ admits a discrete and cocompact group of
isometries, then $\delta_{\Ga,\,F^{\rm su}}\leq 0$, by Lemma
\ref{lem:elemproppressure} (viii) and the fact that $\wt F^{\,\rm su}$ is
invariant under any isometry group of $\wt M$.

It would be interesting to know whether this inequality remains true
for any non-elementary discrete group of isometries $\Ga$ of $\wt M$.

\medskip (2) The assumptions that the derivatives of the sectional
curvature are uniformly bounded and that $\delta_{\Ga,\,F^{\rm su}}\leq 0$
are satisfied if $M$ is a Riemannian cover of a compact manifold, by
the previous remark. Hence Theorem \ref{theo:liougibbsintro} in the
introduction follows from Theorem \ref{theo:LiouvilleisGibbs}.

\subsection{The H\"older-continuity of the (un)stable Jacobian}
\label{subsec:holdcontunstabjac}

We will prove in this subsection the following H\"older regularity
result stated in Theorem \ref{theo:distributionisholder}.  As we make
no compactness assumption, the choice of distances is important.

As explained in Subsection \ref{subsec:geodflowback}, we endow $T\wt
M$ and then $TT\wt M$ with their Sasaki Riemannian metric, and
$T^1\wt M$ and then $TT^1\wt M$ with their induced Riemannian
metric. We denote by $\|\cdot\|$ their Riemannian norms. Note that the
distance on $T^1\wt M$ defined by Equation \eqref{eq:formdistT1wtM}
defines the same H\"older structure as the Riemannian distance, by Lemma
\ref{lem:dholdequivSasak}.

Given a Euclidean space $E$ with finite dimension $n$ and
$k\in\{1,\cdots, n-1\}$, we endow the {\it Grassmannian
  manifold}\index{Grassmannian!manifold} $\operatorname{Gr}_k(E)$ of
$k$-dimensional vector subspaces of $E$ with the distance which is the
Hausdorff distance between their unit spheres
$$
d_E(A,A')= \max_{w\in A, \;w'\in A',\;\|w\|=\|w'\|=1}\;
\Big\{\min_{v'\in A',\;\|v'\|=1}\|w-v'\|,
\min_{v\in A,\;\|v\|=1}\|w'-v\|\Big\}\;.
$$

Let $m$ be the dimension of $\wt M$. Let $E\ra T^1\wt M$ be a vector
bundle over $T^1\wt M$ of rank $r\geq 1$, endowed with a Riemannian
bundle metric and with a length distance inducing the Euclidean
distance on each fibre. We will be interested in two such Riemannian
vector bundles, the tangent bundle of $T^1\wt M$ endowed with the
Sasaki metric and its Riemannian distance, and the fibre bundle over
$T^1\wt M$ whose fibre over $v\in T^1\wt M$ is the product Euclidean
space $T_vT^1\wt M\times T_{\phi_1v}T^1\wt M$, with the distance
induced by the Riemannian distance on $TT^1\wt M\times TT^1\wt M$. For
$k\in\NN$, denote by $\operatorname{Gr}_k(E)\ra T^1\wt M$ the {\it
  Grassmannian bundle}\index{Grassmannian!bundle} of rank $k$ of $E\ra
T^1\wt M$, whose fibre over $v\in T^1\wt M$ is the Grassmannian
manifold $\operatorname{Gr}_k(E_v)$ of $k$-dimensional vector
subspaces of the Tiber $E_v$ of $E\ra T^1\wt M$ over $v$. For
instance, $v\mapsto E^{\rm ss}(v)$ and $v\mapsto E^{\rm su}(v)$ are
continuous sections of $\operatorname{Gr}_{m-1}(TT^1\wt M)\ra T^1\wt
M$. Upon identifying a linear map and its graph, the map $v\mapsto
T_v\phi_1$ is a continuous section of the Grassmannian bundle of rank
$2m-1$ of the second Riemannian vector bundle above.  Since the unit
spheres of the elements of $\operatorname{Gr}_k(E)$ are nonempty
compact subsets of the metric space $E$ if $1\leq k\leq r$, we endow
$\operatorname{Gr}_k(E)$ with the distance $d$, where $d(A,A')$ is the
Hausdorff distance between the unit spheres of the $k$-dimensional
Euclidean spaces $A,A'\in \operatorname{Gr}_k(E)$. The restriction of
$d$ to each fibre $\operatorname{Gr}_k(E_v)$ is the distance $d_{E_v}$
defined above.

\btheo\label{theo:distributionisholder} Let $\wt M$ be a complete
simply connected Riemannian manifold, with dimension at least $2$ and
pinched negative sectional curvature. Assume that the derivatives of
the sectional curvature are uniformly bounded.  Then the maps
$v\mapsto E^{\rm su}(v)$ and $v\mapsto E^{\rm ss}(v)$ defined on $T^1\wt M$
are (globally) H\"older-continuous.  
\etheo

The assumption that the derivatives of the sectional curvature are
uniformly bounded is necessary: Ballmann-Brin-Burns \cite{BalBriBur87}
have constructed a finite volume complete Riemannian surface with
pinched negative sectional curvature whose strong stable foliation is
not H\"older-continuous.

\medskip Since the distributions are Lipschitz (uniformly smooth) in
the direction of the flow lines, note that the maps $v\mapsto E^{\rm
  u} (v)$, $v\mapsto E^{\rm s}(v)$ are also (uniformly locally)
H\"older-continuous, and that the regularity of $\wt F^{\rm su}$
claimed in Theorem \ref{theo:potentialholder} holds. The reversibility
of the potential $F^{\rm su}$ follows from the fact that the Liouville
measure is invariant under the geodesic flow. The boundedness of $\wt
F^{\rm su}$ follows from the bounds on the sectional curvature and its
tangent map.

We will only prove the result for $v\mapsto E^{\rm ss}(v)$, giving
furthermore an explicit estimate on the H\"older exponent. The strong
unstable case is similar, by time reversal.

\medskip Let us give a few historical comments on Theorem
\ref{theo:distributionisholder}. When $\wt M$ has dimension $2$ and
has a cocompact discrete isometry group, the strong stable and
unstable foliations are $\operatorname{C}^1$, hence H\"older regular by
compactness (see \cite{Hopf40}). Furthermore, Hurder and Katok
\cite[Theo.~3.1, Coro.~3.5]{HurKat90} have proved that these
sub-bundles are $\operatorname{C}^{1,\alpha}$ for every $\alpha\in
\mathopen{]}0,1\mathclose{[}$ (see also \cite{HirPug75}), and that if
they are $\operatorname{C}^{1,1}$, then they are $\cinf$, in which
case by a theorem of Ghys \cite[p.~267]{Ghys87}, the geodesic flow is
$\cinf$-conjugated to the geodesic flow of a hyperbolic surface.

When $\wt M$ has dimension at least $3$ and has a cocompact discrete
isometry group, this result of H\"older regularity of the strong
stable and unstable foliations is due to Anosov \cite{Anosov69}.  We
follow Brin's and Brin-Stuck's proofs in \cite{Brin95}, which
modernise Anosov's original proof under the same assumptions, and
\cite[\S 6.1]{BriStu02} in the case of Anosov diffeomorphism (not
flows) on compact manifolds. Here, we provide details to explain why
this proof works under the assumptions of Theorem
\ref{theo:distributionisholder}.

\bigskip
\noindent{\bf Proof of Theorem \ref{theo:distributionisholder}. }
Recall that the geodesic flow of $\wt M$ is an {\it Anosov
  flow}\index{Anosov flow} for the Sasaki metric: there exist $c>0$
and $\lambda\in\mathopen{]}0,1\mathclose{[}$ (and one may take
$\lambda=e^{-1}$ since the sectional curvature of $\wt M$ is
normalised to have upper bound $-1$, by standard arguments of Jacobi
fields) such that for all $v\in T^1\wt M$, $V^{\rm ss}\in E^{\rm ss}(v)$,
$V^{\rm su}\in E^{\rm su}(v)$ and $t\ge 0$, we have
$$
\|T\phi_t(V^{\rm ss})\|\le c\,\lambda^t\, \|V^{\rm ss}\|
\;\;\;{\rm and}\;\;\;\|T\phi_{-t}(V^{\rm su})\|\le 
c\,\lambda^t\, \|V^{\rm su}\|\;,
$$
and for all $v\in T^1\wt M$, $V^{\rm s}\in E^{\rm s}(v)$, $V^{\rm u}
\in E^{\rm u}(v)$ and $t\ge 0$, we have
$$
\|T\phi_{-t}(V^{\rm s})\|\geq c^{-1}\|V^{\rm s}\|
\;\;\;{\rm and}\;\;\;
\|T\phi_t(V^{\rm u})\|\geq c^{-1}\|V^{\rm u}\|\;.
$$

We start by giving two lemmas, the second one, that we state without
proof, involving only linear algebra.

\blemm\label{lem:anglstabdistrib} Let $\wt M$ be a complete simply
connected Riemannian manifold, with dimension at least $2$ and pinched
negative sectional curvature. Then the angles between $E^{\rm ss}(v)$
and any one of $E^{\rm su}(v)$, $E^{0}(v)$ or $E^{\rm u}(v)$ are
uniformly bounded from below by a positive constant.  
\elemm

\dem We prove that the angle between $E^{\rm ss}(v)$ and $E^{\rm u}(v)$ is
bounded from below by a positive constant, the other cases are
analogous.

Using the operator norm on linear maps, let $C=\sup_{v\in T^1\wt M}
\|T_v\phi_1\|$, which is finite since $M$ has pinched negative
curvature (see for instance \cite[page 65]{Ballmann95}).  Let $c>0$
and $\lambda\in\mathopen{]}0,1\mathclose{[}$ be as above. Let
$N\in\NN$ be large enough so that $c^{-1}\geq 2c\lambda^{N}$. For all
$v\in T^1\wt M$, $V^{\rm ss}\in E^{\rm ss}(v)$ and $V^{\rm u}\in
E^{\rm u}(v)$ with $\|V^{\rm ss}\|=\|V^{\rm u}\|=1$, we have
\begin{align*}
C^N\|V^{\rm u}-V^{\rm ss}\|&\geq \|T\phi_N(V^{\rm u}-V^{\rm ss})\|\geq 
\|T\phi_N(V^{\rm u})\|-\|T\phi_N(V^{\rm ss})\|\\ &\geq 
c^{-1}\|V^{\rm u}\|-c\lambda^N \|V^{\rm ss}\|\geq c\lambda^N\;.
\end{align*}
The result follows. \cqfd

\blemm \label{lem:brinstuck} (Brin-Stuck \cite[Lem.~6.1.1]{BriStu02})
Let $c>0$, $\lambda\in\mathopen{]}0,1\mathclose{[}$, $\delta\in
\mathopen{[}0,1\mathclose{[}$ and $D\geq 1$, let $E$ be a finite
dimensional Euclidean vector space, let $A^1,A^2$ be two vector
subspaces of $E$ of the same dimension, and let $(L^1_n)_{n\in\NN},
(L^2_n)_{n\in\NN}$ be two sequences of linear endomorphisms of $E$
such that, for every $n\in\NN$, for $i=1,2$, for all $v\in A^i$ and
$w\in (A^i)^\perp$, we have
$$
\|L^1_n-L^2_n\|\leq \delta D^n,\;\;\; \|L^i_n(v)\|\leq 
c\,\lambda^n\, \|v\|,\;\;\; {\rm and}
\;\;\;\|L^i_n(w)\|\geq c^{-1} \|w\|\;.
$$
Then $d_E(A^1,A^2)\leq 3\;\frac{c^2}{\lambda}\;
\delta^{\;\frac{-\ln \lambda}{\ln D-\log\lambda}}$. 
\elemm

Up to replacing the Sasaki metric by the equivalent (by Lemma
\ref{lem:anglstabdistrib}) Riemannian metric whose norm is
\begin{equation}\label{eq:adaptreimmet}
\|V\|'=\sqrt{\|\,V^{\rm su}\,\|^2+\|\,V^0\,\|^2+\|\,V^{\rm ss}\,\|^2}
\end{equation}
for every $V=V^{\rm su}+V^0+V^{\rm ss}$ where $V^{\rm su}\in E^{\rm
  su}(v)$, $V^0\in E^{0}(v)$ and $V^{\rm ss}\in E^{\rm ss}(v)$, we
assume that the subspaces $E^{\rm su}(v), E^0(v)$ and $E^{\rm ss}(v)$
are pairwise orthogonal for every $v\in T^1\wt M$. For all $v,w\in
T^1\wt M$ and $W\in T_wT^1M$, we denote by $\parallel_w^v W$ the
parallel transport of $W$ along any fixed length-minimising geodesic
segment from $w$ to $v$. Note that the Sasaki metric on $T^1\wt M$ has
a positive lower bound on its injectivity radius at all points, and
this remains true with the modified metric. In particular, there
exists $\kappa>0$ such that for all $v,w\in T^1\wt M$, if
$d(v,w)<\kappa$, there exists a unique geodesic segment from $w$ to
$v$ (recall that $T^1\wt M$ is not nonpositively curved), and we will
use the parallel transports along these geodesic. Let $D\geq 4$ be a
constant such that, for all $v,w\in T^1\wt M$,
\begin{align*}
&d(\phi_1v,\phi_1w)\leq\sqrt{D}\; d(v,w),\;\;\; 
\|T_v\phi_1\|\leq\sqrt{D}
\\{\rm and}\;\;&
\big\|\;T_v\phi_1-\parallel_{\phi_1w}^{\phi_1v}\circ 
\;T_w\phi_1 \circ \parallel_v^w\big\| \leq \sqrt{D}\; d(v,w)\;,
\end{align*}
which exists since $v\mapsto T_v\phi_1$ is bounded and Lipschitz,
since the derivatives of the sectional curvature are uniformly bounded
(see for instance \cite[page 64]{Ballmann95}). For all $v,w\in T^1\wt
M$, let
$$
a_n(v)=\|T_v\phi_n\|\;\;\;{\rm and}\;\;\;b_n(v,w)=
\big\|\;T_v\phi_n-\parallel_{\phi_nw}^{\phi_nv}\circ 
\;T_w\phi_n \circ \parallel_v^w\big\|\;.
$$
Then $a_{n+1}(v)=\|T_{\phi_nv}\,\phi_1\circ T_v\phi_n\|\leq a_1(\phi_n
v) \;a_n(v)$ and
\begin{align*}
b_{n+1}(v,w)=&\|T_{\phi_nv}\,\phi_1\circ T_v\phi_n-
\parallel_{\phi_1\phi_nw}^{\phi_1\phi_nv}\circ T_{\phi_nw}\phi_1\circ 
\parallel_{\phi_nv}^{\phi_nw}\big(\parallel_{\phi_nw}^{\phi_nv}
\circ T_w\phi_n \circ \parallel_v^w\big)\|\\ & \leq a_1(\phi_n
v)\; b_n(v,w) + b_1(\phi_n v,\phi_n w)\; a_n(w)\;.
\end{align*}
Hence by induction and a geometric series argument, $a_n(v)\leq
D^{\frac{n}{2}}$ and
$$
b_n(v,w)\leq D^{\frac{n}{2}}\;\sum_{k=0}^{n-1}\,d(\phi_k v,\phi_k w)
\leq D^{\frac{n}{2}}\,\frac{D^{\frac{n}{2}}-1}{\sqrt{D}-1}\; d(v,w)
\leq D^n\,d(v,w)
\;.
$$

Now, let $v,w\in T^1\wt M$ be such that $d(v,w)< \kappa$. For every
$n\in\NN$, let $f_{v,\,n}: T_{\phi_nv}T^1\wt M \ra T_{v}T^1\wt M$ be an
isometric map. We are going to apply Lemma \ref{lem:brinstuck} with
$c,\lambda$ given by Anosov's flow property, with $D$ as above, with
$\delta= d(v,w)$, with $E$ the Euclidean space $T_vT^1\wt M$, and with
$$
A^1= E^{\rm ss}(v), \;\;\;A^2=\parallel_w^vE^{\rm ss}(w), 
\;\;\;L^1_n=f_{v,\,n}\circ T_v\phi_n, \;\;\;
L^2_n=f_{v,\,n}\circ \parallel_{\phi_nw}^{\phi_nv}\circ 
\;T_w\phi_n \;\circ\parallel_v^w\;,
$$
for every $n\in\NN$. The hypotheses of Lemma \ref{lem:brinstuck} are
easily checked (since $(A^1)^\perp= E^{\rm u}(v)$ and the parallel
transport is isometric), and we hence have
$$
d(E^{\rm ss}(v),\parallel_w^v\! E^{\rm ss}(w))\leq 3\;\frac{c^2}{\lambda}\;
d(v,w)^{\frac{-\ln \lambda}{\ln D-\log\lambda}}\;.
$$
Since $d(\,\parallel_w^v\! E^{\rm ss}(w),E^{\rm ss}(w))\leq d(v,w)$ by
Equation \eqref{eq:distSasaki}, the result hence follows by the
triangle inequality. \cqfd

\subsection{Absolute continuity of the strong unstable foliation}
\label{subsec:absolutecontinuity}

In order to prove Theorem \ref{theo:LiouvilleisGibbs}, it will be
useful to write the Liouville measure as the local product (up to a
density) of the Lebesgue measures on the unstable manifolds and the
Lebesgue measures on the strong stable manifolds.  This is due to
Anosov and Sinai in the compact case, but we follow closely Brin's
proof \cite{Brin95}, adapting it to the noncompact case.

Let us first recall some definitions. Let $N$ be a smooth Riemannian
manifold, and let $\F$ be a continuous foliation of $N$ with smooth
leaves. We denote by $\F(v)$ the leaf of $v\in N$. By a {\it
  transversal}\index{transversal} to $\F$, we mean in this subsection
a smooth submanifold $T$ of $N$ such that for every $v\in T$, we have
$T_vT\oplus T_v\F_v=T_vN$. We denote by $\vol_T$ the induced
Riemannian measure on $T$. For instance, every stable leaf is a
transversal to the strong unstable foliation $\F=\W^{\rm su}$ of $N=
T^1\wt M$ endowed with the Sasaki metric.

If $U$ is the domain of a foliated chart of $\F$, the {\it local
  leaf}\index{local leaf}\index{leaf!local} $\F_U(v)$ of $v\in U$ in
$U$ is the connected component of $v$ in $\F(v)\cap U$. We denote by
$\vol_{\F_U(v)}$ the induced Riemannian measure on $\F_U(v)$. We
denote by ${\cal U}\!\T$ the set of pairs $(U,T)$ where $U$ is the
domain of a foliated chart of $\F$ and $T$ is a transversal to $\F$
contained in $U$ such that $U=\bigcup_{v\in T} \F_U(v)$ and
$\F_U(v)\cap\F_U(w)=\emptyset$ for all $v\neq w$ in $T$. In
particular, we have a continuous fibration $U\ra T$, with fibres the
local leaves. For every $v\in N$, there exists $(U,T)\in {\cal U}\!\T$
such that $v\in T$. A {\it holonomy map}\index{holonomy map} of $\F$
is a homeomorphism $h:T\ra T'$ between two transversals to $\F$
contained in the domain $U$ of a foliated chart of $\F$, such that
$h(v)\in\F_U(v)$ for every $v\in T$ and $\F_U(v)\cap\F_U(w)=\emptyset$
for all $v\neq w$ in $T$. Note that these holonomy maps are in general
only continuous, and a priori not absolutely continuous with respect
to the Lebesgue measures. When $\F=\W^{\rm su}$ and $N= T^1\wt M$,
they are H\"older-continuous by Subsection
\ref{subsec:holdcontunstabjac} if the sectional curvature $M$ has
bounded derivatives. But this is not sufficient to imply that they are
absolutely continuous with respect to the Lebesgue measures.

The foliation $\F$ is {\it transversally absolutely continuous}%
\index{foliation!transversally absolutely continuous}%
\index{transversally absolutely continuous} if its holonomy maps
preserve the Lebesgue measure classes, that is, if for every holonomy
map $h:T\to T'$, the measure $h_*\vol_T$ is absolutely continuous with
respect to $\vol_{T'}$.  The foliation $\F$ is said to be {\it
  transversally absolutely continuous with locally uniformly bounded
  Jacobians} if furthermore, for every compact set $K$ in $N$ and for
every holonomy map $h:T\to T'$ of $\F$ with $T,T'\subset K$, the {\it
  Jacobian}\index{Jacobian} $\frac{d(h^{-1})_*\vol_T}{d\,\vol_{T'}}$
of the absolutely continuous map $h$ is (almost everywhere) bounded by
a constant depending only on $K$. By the change of variable formula,
if $\F$ is smooth, then $\F$ is transversally absolutely continuous,
and the above Jacobians of the holonomy maps are their Jacobians as
smooth maps.

The foliation $\F$ is {\it absolutely
  continuous}\index{foliation!absolutely continuous}\index{absolutely
  continuous} if for every $(U,T)\in {\cal U}\!\T$, the restriction to
$U$ of the Riemannian measure of $N$ disintegrates with respect to the
fibration $U\ra T$ over the induced Riemannian measure $\vol_T$ on
$T$, with conditional measures absolutely continuous with respect to
the induced Riemannian measures $\vol_{\F_U(v)}$ on the local leaves
$\F_U(v)$ for ($\vol_T$-almost) every $v\in T$: There exists a
measurable family $(\delta_v)_{v\in T}$ of positive measurable
functions $\delta_v:\F_U(v)\ra \RR$ such that for every measurable
subset $A$ of $U$, we have
$$
\vol_{T^1\wt M}(A)= \int_{v\in T} \int_{w\in \F_U(v)} 
\mathbbm{1}_A(w) \delta_v(w)\,d\vol_{\F_U(v)}(w)\,d\vol_T(v)\;.
$$
It is {\it absolutely continuous with locally bounded conditional
  densities} if furthermore, for every compact subset $K$ in $N$,
there exists $c_K\geq 0$ such that, for every $(U,T)\in {\cal U}\!\T$
with $U\subset K$, the above densities $\delta_v$ are bounded for
every $v\in T$ by the constant $c_K$.

By \cite[Prop.~3.5]{Brin95} and \cite[Rem.~3.9]{Brin95}, if $\F$ is
transversally absolutely continuous (respectively transversally
absolutely continuous with locally uniformly bounded Jacobians), then
$\F$ is absolutely continuous (respectively absolutely continuous with
locally bounded conditional densities), but the converse is not always
true.

\btheo\label{theo:transverseabsolutecontinuity} Let $\wt M$ be a
complete simply connected Riemannian manifold, with dimension at least
$2$ and pinched negative sectional curvature. Assume that the
derivatives of the sectional curvature are uniformly bounded.  Then
the strong stable, strong unstable, stable and unstable foliations are
transversally absolutely continuous with locally uniformly bounded
Jacobians, hence are absolutely continuous with locally bounded
conditional densities.  
\etheo

\dem We follow closely Brin's proof \cite[Theo.~5.1]{Brin95}, only
giving the arguments for the extension to the noncompact case. By the
smoothness of the geodesic foliation, and by a time reversal argument,
we only consider the case of the strong unstable foliation.

Let $K$ be a compact subset of $T^1\wt M$, let $h:A\ra A'$ be a
holonomy map contained in the domain $U$ of a foliated chart of
$\W^{\rm su}$, with $U\subset K$. For every $n\in\NN$, let $\F_n$ be
the smooth foliation of $T^1\wt M$, whose leaf $\F_n(v)$ through every
$v\in T^1\wt M$ is the set of outer unit normal vector to the sphere
$S(\pi(\phi_{-n}v),n)$ of radius $n$ centred at $\pi(\phi_{-n}v)$. As
$n\ra+\infty$, the foliation $\F_n$ converges to $\W^{\rm su}$
uniformly on compact sets: for every $\epsilon>0$, if $n$ is large
enough, for all $v,w\in K$ such that $w\in W^{\rm su}(v)$, if $w_n$ is
the outer unit normal vector to $S(\pi(\phi_{-n}v),n)$ at its closest
point to $\pi(w)$, we have $d(w,w_n)\leq \epsilon$. Hence, up to
shrinking $A,A'$ to neighbourhoods of given points $v\in A, v'\in A'$
respectively (and then using finite covering arguments), for $n$ large
enough, the foliation $\F_n$ is transversal to $A$ and $A'$, and there
exists a (smooth) holonomy map $h_n:A\ra A'$ for $\F_n$, which
converges uniformly to $h$ as $n\ra+\infty$.

By \cite[Lem.~2.7]{Brin95}, in order to prove Theorem
\ref{theo:transverseabsolutecontinuity}, it is sufficient to prove
that the (smooth) Jacobian $\operatorname{Jac}h_n$ of $h_n$ is
uniformly bounded on $U$, by a constant depending only on $K$. Let us
write
$$
h_n=\phi_{-n}\circ H_n\circ \phi_n\;,
$$ 
where $H_n:\phi_n(A)\ra\phi_n(A')$ is a holonomy map of the foliation
with leaves the fibres of the fibration $T^1\wt M\ra \wt M$.  For
$0\le k\le n-1$, we denote by 
$J_k,J'_k$ the Jacobian of $\phi_1$ at $\phi_k(v),\phi_k(v')$,
respectively. We have by the chain rule
$$
\operatorname{Jac}h_n(v)=\prod_{k=0}^{n-1} \frac{1}{J'_k} \;\;
\operatorname{Jac}H_n(\phi_nv)\;\;\prod_{k=0}^{n-1} J_k =
\operatorname{Jac}H_n(\phi_nv)\;\prod_{k=0}^{n-1}\frac{J_k}{J'_k}\;.
$$

Exactly as in the compact case (see \cite[page 94]{Brin95}), the
H\"older continuity of the strong unstable distribution proves that
there exist two constants $c,c'>0$ (depending only on $K$) such that,
for all $k\geq 0$ and $v\in A$, we have $|J_k-J'_k|\leq c\,
e^{-c'\,k}$.  The fact that the sectional curvature is pinched
between two negative constants implies that there exists a constant
$c''>0$ such that, for all $k\geq 0$ and $v\in A$, we have
$\frac{1}{c''}\leq J_k\leq c''$.  Since $1+x\leq e^x$ for all $x\geq
0$, we therefore have
$$
\prod_{k=0}^{n-1}\frac{J_k}{J'_k}\leq 
\prod_{k=0}^{n-1}(1+c\,c''\,e^{-c'\,k})
\leq e^{c\,c''\sum_{k=0}^{+\infty} e^{-c'\,k}}<+\infty\,.
$$

Exactly as in the compact case (see \cite[Lem.~4.2]{Brin95}), the
Anosov property of the geodesic flow in pinched negative sectionnal
curvature implies that at all of their points $v$, the transversals
$\phi_n(A)$ and $\phi_n(A')$ are uniformly (in $n$ for $n$ large
enough and in $v$) transverse to the unit tangent spheres
$T_{\pi(v)}\wt M$. Since $H_n$ is a holonomy map of a smooth
Riemannian foliation of a Riemannian manifold, between two
transversals whose angles with the leaves of the foliation are
uniformly bounded from below, the Jacobian $\operatorname{Jac}H_n$ is
uniformly bounded. Hence the Jacobian $\operatorname{Jac}h_n$ is
uniformly bounded, as required. \cqfd

\subsection{The Liouville measure as an equilibrium state}
\label{subsec:liouvequilstate}

In this subsection, we recall the definition of Lyapounov's exponents
and the main results about them (Oseledets's theorem, Ruelle's
inequality, the Pesin formula, see for instance \cite{Mane87,
  Ledrappier84, Mane81, KatMen95}), in order to prove that the Liouville
measure, once normalised, is the equilibrium state of an appropriate
potential, under compactness assumptions.

We denote here by $\phi=\phi_1$ the time one map of the geodesic flow
$(\phi_t)_{t\in\RR}$ on $T^1M$. We fix a probability measure $\mu$ on
$T^1M$ invariant under $(\phi_t)_{t\in\RR}$. Recall that $\log^+
t=\max\{\log t, 0\}$ for every $t>0$, and that $T^1\wt M$ is endowed
with Sasaki's Riemannian metric.

\medskip {\it Oseledets's
  theorem}\index{Oseledets!theorem}\index{theorem@Theorem!of
  Oseledets} (see for instance \cite{Oseledets68}, \cite[Theo.~S.2.9
page 665]{KatMen95}) asserts that if $\int_{T^1M}\log^+ \|T\phi^{\pm
  1} \|\,d\mu<\infty$ (in particular if $M$ is compact), then there
exists a Borel subset $\operatorname{Reg}$ of $T^1M$, invariant under
$(\phi_t)_{t\in\RR}$ and with full measure with respect to $\mu$, such
that for every $v\in\operatorname{Reg}$, there exists a unique direct
sum decomposition (of closed convex cones when $\Ga$ has torsion),
called the {\it Oseledets decomposition} at
$v$,\index{Oseledets!decomposition}
$$
T_vT^1M=\oplus_{i=1}^{s(v)}E_i(v)
$$
and unique real numbers $\chi_1(v)<\dots <\chi_{s(v)}(v)$, called the
{\it Lyapounov exponents}\index{Lyapounov's exponents} of $v$, such
that
\begin{enumerate}
\item[(1)] the map $v\mapsto s(v)$ from $\operatorname{Reg}$ to $\NN$
  is measurable; for all $i,r\in\NN$, the map $v\mapsto E_i(v)$ from
  $\{v\in\operatorname{Reg} \;:\;s(v)\geq i, \operatorname{dim} E_
  i(v)=r\}$ to the space of closed subsets of $TT^1M$ (endowed with
  Chabauty's topology) is measurable; and for all $v\in
  \operatorname{Reg}$ and $i\in\{1,\dots ,s(v)\}$, we have
  $T_v\phi(E_i(v)) =E_i(\phi v)$ and $\chi_i(\phi v)= \chi_i(v)$;
\item[(2)] for all $v\in \operatorname{Reg}$, $i\in\{1,\dots ,s(v)\}$
  and $V\in E_i(v)-\{0\}$,
$$
\lim_{n\to\pm \infty}\frac{1}{ n }\log \|T_v\phi^n(V)\|=\chi_i(v)\;;
$$
\item[(3)] for all $v\in \operatorname{Reg}$ and $i\in\{1,\dots ,
  s(v)\}$, if $\dim E_i(v)= k_i(v)$, then
$$
\lim_{n\to + \infty} \frac{1}{n}
\log J^{E_i}\phi^n(v)=k_i(v)\,\chi_i(v)\,,
$$ 
where $J^{E_i}\phi^n(v)$ denotes the Jacobian of the restriction of
$T_v\phi^n$ to $E_i(v)$;
\item[(4)] the angle between two subspaces of the Oseledets
  decomposition has at most subexponential decay: for all $v$ in
  $\operatorname{Reg}$, $i\neq j$ in $\{1, \dots , s(v)\}$, $V\in
  E_i(v)-\{0\}$ and $W\in E_j(v)-\{0\}$, we have
$$
\lim_{n\to\pm\infty}\frac{1}{|n|}\log 
\big|\,\sin\angle (T_v\phi^n(V), T_v\phi^n(W))\big|=0\;.
$$ 
\end{enumerate}

In particular, since $E^{\rm su}(v)=\oplus_{1\leq i\leq
  s(v),\;\chi_i(v)>0} \;E_i(v)$ for every $v\in\operatorname{Reg}$, we
deduce from the last two points that for $\mu$-almost every $v\in
T^1M$, we have
$$
\lim_{n\to +\infty}\frac{1}{n}\int_0^nF^{\rm su}(\phi_t v)\,dt=
-\lim_{n\to +\infty}\frac{1}{n}\log J^{\rm su}(v,n)=
-\sum_{1\leq i\leq s(v),\;\chi_i(v)>0}k_i(v)\,\chi_i(v)\;.
$$
Note that when $\mu$ is ergodic for $(\phi_t)_{t\in\RR}$, by
Birkhoff's ergodic theorem, for $\mu$-almost every $v\in T^1M$, this
equality becomes
\begin{equation}\label{eq:Lyapounovregular}
\int_{T^1M}F^{\rm su}\,d\mu=
-\sum_{1\leq i\leq s(v),\;\chi_i(v)>0} k_i(v)\,\chi_i(v)\;.
\end{equation} 

\medskip The {\it Ruelle inequality}\index{Ruelle inequality} (see
\cite{Ruelle78}) asserts that if $M$ is compact, then
\begin{equation}\label{eq:Ruelleinequality}
h_{\mu}(\phi)\leq \int_{v\in T^1M} 
\sum_{1\leq i\leq s(v),\;\chi_i(v)>0}k_i(v)\,\chi_i(v)\,d\mu(v)\;.
\end{equation}
It is unclear if it is possible to extend the proof of Ruelle's
inequality in the case where $M$ is noncompact, with perhaps
additional geometric assumptions besides the pinched negative
sectional curvature.

When $M$ is compact and $\mu$ is ergodic, by Equation
\eqref{eq:Lyapounovregular} and Equation \eqref{eq:Ruelleinequality}, we
have $h_{\mu}(\phi)+\int_{T^1M}F^{\rm su}\,d\mu\le 0$.  Since the
upper bound defining the topological pressure of a potential (see
Chapter \ref{sec:variaprincip}) may be taken on the ergodic
probability measures invariant under $(\phi_t)_{t\in\RR}$ by the convexity
properties of the metric entropy, we hence have
\begin{equation}\label{eq:negatifinequelle}
P(\Ga,F^{\rm su})\le 0\;.
\end{equation}

The {\it Pesin formula}\index{Pesin formula} (see for instance
\cite{Pesin77, Mane81}) asserts that if $\mu$ is absolutely continuous
with respect to the Lebesgue measure class, then we have equality in
Equation \eqref{eq:Ruelleinequality}:
\begin{equation}\label{eq:Pesinformula}
h_{\mu}(\phi)= \int_{v\in T^1M} 
\sum_{1\leq i\leq s(v),\;\chi_i(v)>0}k_i(v)\,\chi_i(v)\,d\mu(v)\;.
\end{equation}
If $M$ is noncompact, it is unclear when the Pesin formula remains
valid.

In particular, if $M$ is compact, then the Liouville measure
$\vol_{T^1M}$ is ergodic (see for instance \cite[page 95]{Brin95}),
and the probability measure $L$ proportional to the Liouville measure
satisfies
$$
P(\Ga,F^{\rm su})\le 0= h_{L}(\phi)+\int_{T^1M}F^{\rm su}\,dL\,.
$$
Hence $L$ is an equilibrium state for the potential $F^{\rm su}$. The
Variational Principle (Theorem \ref{theo:principevariationnel})
implies that $L=m_{F^{\rm su}}$, which proves the following expected
result.

\btheo\label{theo:liouvillegibbscascompact} If $M$ is compact, then
the normalised Liouville measure $\frac{\vol_{T^1M}}{\Vol(T^1M)}$
coincides with the normalised Gibbs measure $\frac{m_{F^{\rm su}}}
{\|m_{F^{\rm su}}\|}$ of the potential $F^{\rm su}$. 
\cqfd 
\etheo

\rem Assume in this remark that $M$ has finite volume with constant
curvature. Then $F^{\rm su}$ is constant, hence the Gibbs measure of
the potential $F^{\rm su}$ and the Bowen-Margulis measure (which is
the Gibbs measure of the potential $0$) coincide.  The
Patterson-Sullivan measure of $\Ga$ at the origin of the ball model of
the real hyperbolic $n$-space $\HH^n_\RR$ may be taken to be the
standard Riemannian measure of the unit sphere
$\SS^{n-1}= \partial_\infty \HH^n_\RR$.  The exact ratio between the
total masses of $T^1M$ for the Liouville measure and for the Gibbs
measure $m_{F^{\rm su}}$ is computed in \cite[Sect.~7]{ParPauRev},
yielding
$$
\frac{\|m_{F^{\rm su}}\|}{\Vol(T^1M)}=2^{n-1}\Vol(\SS^{n-1})\;.
$$

\subsection{The Liouville measure satisfies the Gibbs property}
\label{subsec:Gibbsapproach}

In this subsection, we prove that, under the previous bounded geometry
assumption, the Liouville measure satisfies a Gibbs property as
defined in Subsection \ref{subsec:gibbsproperty}, as well as a
leafwise version of it.

\medskip In the proof of this Gibbs property, as well as in the coming
Subsection \ref{subsec:proofLiouvilleisGibbs} in order to prove Theorem
\ref{theo:LiouvilleisGibbs}, we will use the following version of the
Gibbs property for the Liouville measure in restriction to the strong
unstable and strong stable foliations.  Recall that $B^{\rm su}(v,r)$
is the open ball of center $v$ and radius $r$ for the Hamenst\"adt
distance $d_{W^{\rm su}(v)}$ on $W^{\rm su}(v)$, and similarly for
$B^{\rm ss}(v,r)$ (see Subsection \ref{subsec:geodflowback}).

\blemm \label{lem:gibbsinstable} For every compact subset $K$ in $T^1
\wt M$ and for every $r>0$, there exists a constant $c'_{K,\,r}>0$
such that for all $v\in \Ga  K$ and $T\geq 0$ verifying $\phi_T v\in
\Ga   K$, we have
$$
\frac{1}{c'_{K,\,r}} \, e^{\int_{0}^T \wt F^{\,\rm su}(\phi_tv)\,dt} \le 
\vol_{W^{\rm su}(v)}(\phi_{-T}B^{\rm su}(\phi_T v,r))
\le c'_{K,\,r}\;  e^{\int_{0}^T \wt F^{\,\rm su}(\phi_tv)\,dt} \,,
$$
and for all $v\in \Ga K$ and $T\geq 0$ verifying $\phi_{-T} v\in \Ga
K$, we have
$$
\frac{1}{c'_{K,\,r}}\,  e^{\int_{-T}^0 \wt F^{\,\rm su}\circ\iota(\phi_tv)\,dt} \le 
\vol_{W^{\rm ss}(v)}(\phi_{T}B^{\rm ss}(\phi_{-T} v,r))
\le c'_{K,\,r}\;  e^{\int_{-T}^0   \wt F^{\,\rm su}\circ\iota(\phi_tv)\,dt} \,.
$$ 
\elemm

\dem Fix $K,r$ as in the statement. In order to prove the first claim,
observe that, for $v,T$ as in its statement,
\begin{align*}
\vol_{W^{\rm su}(v)}(\phi_{-T}B^{\rm su}(\phi_T v,r)) &
=\int_{B^{\rm su}(\phi_T v,\,r)} \, d\big((\phi_T)_*\vol_{W^{\rm su}(v)}\big)
\\ & =\int_{B^{\rm su}(\phi_T v,\,r)} 
\operatorname{Jac}\big((\phi_{-T})_{|W^{\rm su}(\phi_T v)}\big)\, 
d\vol_{W^{\rm su}(\phi_Tv)}\,.
\end{align*}

The H\"older-continuity of the strong unstable distribution and the
compactness of $K$ imply that there exists a constant $c$ depending
only on $K$ and $r$ such that, for $v,T$ as in the statement of the
first claim and for every $w\in B^{\rm su} (\phi_T v,\,r)$, we have
$$
\frac{1}{c}\, \operatorname{Jac}
\big((\phi_{-T})_{|W^{\rm su}(\phi_T v)}\big)(\phi_T v)\le 
\operatorname{Jac}\big((\phi_{-T})_{|W^{\rm su}(\phi_T v)}\big)(w)
\le c\;\operatorname{Jac}
\big((\phi_{-T})_{|W^{\rm su}(\phi_T v)}\big)(\phi_T v)\,.
$$

By the definition of the unstable Jacobian $\wt J^{\rm su}$, we have
$\operatorname{Jac}({\phi_{-T}}_{|W^{\rm su}(\phi_T v)})(\phi_T v)=\wt J^{\rm
  su} (\phi_T v, -T)$. Hence
\begin{multline}\label{eq:multlineBsuJsu}
\frac{1}{c}\;\wt J^{\rm su} (\phi_T v, -T) 
\vol_{W^{\rm su}(\phi_Tv)}(B^{\rm su}(\phi_T v,\,r))\\ \leq
\vol_{W^{\rm su}(v)}(\phi_{-T}B^{\rm su}(\phi_T v,r))\\ \leq 
c\;\wt J^{\rm  su} (\phi_T v, -T) 
\vol_{W^{\rm su}(\phi_Tv)}(B^{\rm su}(\phi_T v,\,r))  \;.
\end{multline}

Since $K$ is compact, the Hamenst\"adt ball of radius $r$ and center
$u\in K$ in its strong unstable leaf $W^{\rm su}(u)$ is contained in
the Riemannian ball of radius $r'$ and center $u$ in $W^{\rm su}(u)$,
and contains the Riemannian ball of radius $r''$ and center $u$ in
$W^{\rm su}(u)$, for some $r'$ and $r''$ depending only on $r$ and
$K$. For instance by the standard comparison theorems of volumes of
Riemannian balls and by the compactness of $K$, the volume
$\vol_{W^{\rm su}(w)}(B^{\rm su}(w,r))$ for $w\in \Ga K$ is hence
uniformly bounded from above and from below by positive constants
depending only on $K$ and $r$.

By the cocycle formula \eqref{eq:cocyclunstabjac}, for all
$s,t\in\RR$, we have $\log \wt J^{\,\rm su}(\phi_Tv,0)=0$ and
$$
\log \wt J^{\,\rm su}(\phi_Tv,s+t)=
\log \wt J^{\,\rm su}(\phi_t\phi_Tv,s) + 
\log \wt J^{\,\rm su}(\phi_Tv,t)\;,
$$
 hence, by the definition of the potential $\wt F^{\rm su}$,
\begin{align*}
\int_{0}^T \wt F^{\rm su}(\phi_tv)\,dt& =
\int_{-T}^0 \wt F^{\rm su}(\phi_t\phi_Tv)\,dt 
=\int_{-T}^0 -{\frac{d}{ds}}_{|s=0} 
\log \wt J^{\rm su}(\phi_t\phi_Tv,s)\,dt\\ 
& =\int_{-T}^0 -\frac{d}{dt} \log \wt J^{\rm su}(\phi_Tv,t)\,dt
=\log \wt J^{\rm su}(\phi_Tv,-T)\;.
\end{align*}
The first claim hence follows from Equation \eqref{eq:multlineBsuJsu}.

\medskip
In order to prove the second claim, observe that the map $\iota: v
\mapsto -v$ is an isometry of Sasaki's metric on $T^1\wt M$,
anti-commuting with the geodesic flow, exchanging the strong stable
and strong unstable balls (see Equation \eqref{eq:ballhamantipod}), so
that, for $v,T$ as in the statement of the second claim,
$$
\vol_{W^{\rm ss}(v)}(\phi_{T}B^{\rm ss}(\phi_{-T} v,r))=
\vol_{W^{\rm su} (-v)} (\phi_{-T}B^{\rm su}(\phi_T(-v),r))\;.
$$ 
Thus, since $-v$ and $\phi_T(-v)$ belong to $\Ga \iota(K)$ (and
$\iota(K)$ is compact), the first claim gives, for some constant
$c'>0$ depending only on $r$ and $K$,
$$ 
\frac{1}{c'}\, e^{\int_0^T \wt F^{\,\rm su}(\phi_s (-v))\,ds} 
\le \vol_{W^{\rm ss}(v)}(\phi_{T}B^{\rm ss}(\phi_{-T} v,r))
\le c'\; e^{\int_0^T \wt F^{\,\rm su}(\phi_s (-v))\,ds} \,.
$$
The equality 
$$
\int_0^T \wt F^{\,\rm su}(\phi_s(-v))\,ds=
\int_0^T \wt F^{\,\rm su}\circ \iota(\phi_{-s}v)\,ds=
\int_{-T}^0\wt F^{\,\rm su}\circ \iota(\phi_{ s}v)\,ds 
$$
then gives the second claim.  
\cqfd

\bprop\label{prop:LiouvilleJacobien} Assume that the derivatives of
the sectional curvature are uniformly bounded. Then the Liouville
measure on $T^1M$ satisfies the Gibbs property for the potential
$F^{\rm su}$ and the constant $c(F^{\rm su})=0$.  
\eprop

In the compact case, this was proved for hyperbolic diffeomorphisms in
Bowen-Ruelle in \cite[Lem.~4.2]{BowRue75}, and later in
Katok-Hasselblatt in \cite[Lem.~20.4.2]{KatHas95}.  It is possible to
adapt their arguments to flows in the noncompact case.  However, it is
shorter to use the absolute continuity properties of the strong
unstable foliation, proved in Section \ref{subsec:absolutecontinuity}.
We follow this approach here.

\medskip \dem Note that $F^{\rm su}$ is H\"older-continuous and
bounded by Theorem \ref{theo:potentialholder}. 

Note that by the definition of the dynamical balls in Subsection
\ref{subsec:gibbsproperty}, for every compact subset
$K$ of $T^1\wt M$, for all $r\in\mathopen{]}0,
1\mathclose{]}$, $v\in K$ and $T\geq 1$, the set $B(v;T, 0,r)$ is
contained in the compact subset $\pi^{-1}(\mathcal{N}_{1}(\pi( K)))$
of $T^1\wt M$. Hence the multiplicity of the restriction to $B(v;T,
0,r)$ of the map $T^1\wt M\ra T^1M$ is bounded by a constant depending
only on $K$, by the discreteness of the action of $\Ga$ on $T^1\wt M$.

Since the Liouville measure on $T^1\wt M$ is invariant under the
geodesic flow and by Remark (2) following the definition
\ref{def:gibbsprop} (in Subsection \ref{subsec:gibbsproperty}) of the
Gibbs property, we only have to prove that for every compact subset
$K$ of $T^1\wt M$, there exist $r\in\mathopen{]}0,1\mathclose{]}$,
$C_{K,\,r}\ge 1$ and $T_0\geq 0$ such that for all $v\in K$ and $T\geq
T_0$ with $\phi_{T}v\in \Ga K$, we have
$$
\frac{1}{C_{K,\,r}}\;e^{\int_{0}^T\wt F^{\rm su}(\phi_tv)\,dt}\leq 
\vol_{T^1\wt M}(B(v;T, 0,r))\leq 
C_{K,\,r}\; e^{\int_{0}^T \wt F^{\rm su}(\phi_tv)\,dt}\;.
$$

\medskip 
We first replace the dynamical balls by sets better adapted to the
local product structure. For all $r>0$, $v\in T^1\wt M$ and $T\ge 0$,
consider the sets
\begin{eqnarray*}
C(v;T,r)= \bigcup_{w\in B^{\rm ss}(v,\,r)}\;\bigcup_{|s|< r} 
B^{\rm su}(\phi_s w, re^{-T})\;.
\end{eqnarray*}
Note that $\ga C(v; T,r)= C(\ga v;T,r)$ for every $\ga\in\Ga$.

\blemm\label{lem:decompballfol} For every compact subset $K$ in
$T^1\widetilde{M}$ and every $r>0$, there exists $r'=r'_{K,\,r}>0$
such that for all $v\in \Ga K$ and $T\geq 3r$ satisfying $\phi_Tv \in
\Ga K$, we have
$$
C(v;T, \frac{r}{3} )\subset B(v;T,0,r)\subset C(v; T,r')\,.
$$ 
\elemm

\dem Let $v\in T^1\wt M$, $T\geq 0$, $r>0$ and $u\in
C(v;T,\frac{r}{3})$. By the definition of $C(v;T,\frac{r}{3})$, let
$s\in\mathopen{]}-\frac{r}{3}, \frac{r}{3} \mathclose{[}$ and $w\in
B^{\rm ss}(v,\frac{r}{3})$ be such that $u\in B^{\rm su}(\phi_s w,
\frac{r}{3}\,e^{-T})$. Then for every $t\in [0,T]$,  by the
triangle inequality, by Equation \eqref{eq:compardistHamen}, and by the
contraction properties \eqref{eq:dilatHamdist} and
\eqref{eq:contractHamdist} of the Hamenst\"adt distances, we have
\begin{align*}
d(\pi(\phi_tu),\pi(\phi_tv))& \leq d(\pi(\phi_tu),\pi(\phi_{t+s}w))
+d(\pi(\phi_{t+s}w),\pi(\phi_tw))+d(\pi(\phi_tw),\pi(\phi_tv))\\ & \leq
d_{W^{\rm su}(\phi_{t+s}w)}(\phi_tu,\phi_{t+s}w)
+|s|+d_{W^{\rm ss}(\phi_{t}w)}(\phi_tw,\phi_tv)
\\ & =e^{t}\;d_{W^{\rm su}(\phi_{s}w)}(u,\phi_{s}w)
+|s|+e^{-t}\;d_{W^{\rm ss}(w)}(w,v)\\ &
< e^{t}\,\frac{r}{3}\,e^{-T}+\frac{r}{3}+e^{-t}\,\frac{r}{3}\leq r\;.
\end{align*}
Hence $u\in B(v;T,0,r)$ by the definition of these dynamical
balls. This proves the inclusion on the left hand side.

\medskip
To prove the other inclusion, we start by the following remark.

\smallskip\noindent
\begin{minipage}{9.6cm} ~~~ 
For all $u,v\in T^1\wt M$, $T\geq 0$ and $r>0$ such that
$d(\pi(u),\pi(v))\leq r$ and $d(\pi(\phi_Tu),\pi(\phi_Tv))\leq r$, if
$T>2r$ then $u_+\neq v_-$. Otherwise, let $p$ be the closest point to
$\pi(\phi_Tu)$ on the geodesic line defined by $v$. Since $T>2r$ and
by convexity, we have $p\in \mathopen{]}v_-,\pi(v)]$,
$d(\pi(\phi_Tu),p)\leq r$ and $d(p,\pi(v))\geq T-2r$. Hence
\begin{align*}
d(\pi(\phi_Tv),\pi(\phi_Tu))&\geq d(\pi(\phi_Tv),p)-r\\ &=
d(\pi(\phi_Tv),\pi(v))+d(\pi(v),p)-r\\ &\geq T+(T-2r)-r>r\;,
\end{align*}
a contradiction. 
\end{minipage}\begin{minipage}{5.2cm}
\begin{center}
\input{fig_upludifvmin.pstex_t}
\end{center}
\end{minipage}

\medskip For every $\xi\in \partial_\infty\wt M$, let (see Section
\ref{subsec:condmesgibbs})
$$
U_{\xi}=\{u\in T^1\wt M\;:\; u_+\neq \xi\}\;.
$$
A similar proof shows that if furthermore $T\geq 3r$ then $u$ stays in
a compact subset of $U_{v_-}$.

Now, let $K,r$ be as in the statement. Let $v\in K$, $T\geq 3r$ and
$u\in T^1\wt M$ be such that $\phi_Tv \in \Ga K$ and $u\in
B(v;T,0,r)$, so that $d(\pi(u),\pi(v))\leq r$ and $d(\pi(\phi_Tu),
\pi(\phi_Tv))\leq r$. Then $d(u,v)\leq r+\pi$ and $d(\phi_Tu, \phi_Tv)
\leq r+\pi$, by the properties of Sasaki's metric. Recall that the map
from $U_{v_-}$ to $W^{\rm su}(v)$, sending $u\in U_{v_-}$ to the
unique $w\in W^{\rm su} (v)$ such that $u_+=w_+$, is a trivialisable
continuous fibration, with fiber over $w\in W^{\rm su} (v)$ the stable
leaf $W^{\rm s} (w)$. The Hamenstädt distances are proper on the
strong stable and strong unstable leaves, and vary continuously with
the leaves.  Every compact subset $K'$ of $U_{v'_-}$ is contained in
$U_{v''_-}$ if $v''$ is close enough to $v'$, and hence in $C(v'';0,r)$
for $r>0$ (depending only on $K'$ and $v'$) big enough.  By the
compactness of $K$ (with a finite covering argument), and by
equivariance, there hence exists $r'>0$ depending only on $K$ and $r$
such that $u\in C(v;0,r')$ and $\phi_Tu \in C(\phi_Tv;0,r')$.  By the
definition of the sets $C(v_0;0,r_0)$, let $s,s'\in\mathopen{]}-r', r'
\mathclose{[} \,$, and let $w\in B^{\rm ss} (v,r')$ and $w'\in B^{\rm
  ss} (\phi_Tv, r')$ be such that $u\in B^{\rm su} (\phi_s w, r')$ and
$\phi_Tu\in B^{\rm su} (\phi_{s'} w', r')$. 

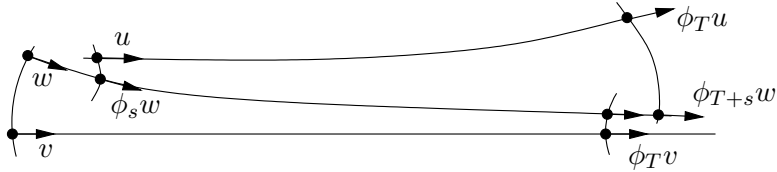
\begin{center}
\input{fig_celldynsimp.pstex_t}
\end{center}

By the uniqueness of the decomposition, we have $w'=\phi_Tw$ and
$\phi_{s'}w'=\phi_T\phi_{s}w$. Hence
$$
d_{W^{\rm su}(\phi_{s}w)}(u,\phi_{s}w)=
e^{-T} \;d_{W^{\rm su}(\phi_T\phi_{s}w)}(\phi_T u,\phi_T\phi_{s}w)
\leq e^{-T}\,r'\;.
$$
Hence $u\in C(v;T,r')$, which proves the inclusion on the right hand
side. \cqfd

\medskip By this lemma, we hence only have to prove that for every
compact subset $K$ in $T^1\wt M$ and every $r>0$, there exists
$C'_{K,\,r}\ge 1$ such that for all $v\in K$ and $T\geq 0$ with
$\phi_{T}v\in \Ga K$, we have
\begin{equation}\label{eq:dynballdecompfol}
\frac{1}{C'_{K,\,r}}\;e^{\int_{0}^T\wt F^{\rm su}(\phi_tv)\,dt}\leq 
\vol_{T^1\wt M}(C(v;T,r))\leq 
C'_{K,\,r}\; e^{\int_{0}^T \wt F^{\rm su}(\phi_tv)\,dt}\;.
\end{equation}

Let us fix $K,r$ as above. Since the foliations $\wt \W^{\rm su}$ and
$\wt \W^{\rm ss}$ are absolutely continuous with locally bounded
conditional densities (see Theorem
\ref{theo:transverseabsolutecontinuity}), there exists $c\geq 1$
(depending only on $K,r$) such that for all $v\in K$ and $T\geq 0$, we
have
\begin{multline} \label{eq:volabscont}
\frac{1}{c}\;\int_{B^{\rm ss}(v,\,r)}\int_{-r}^{r}
\vol_{W^{\rm su}(w)}(B^{\rm su}(\phi_{s}w, r e^{-T}))\;
ds\, d\vol_{W^{\rm ss}(v)}(w) \leq \vol_{T^1\wt M }(C(v;T,r)) \\
\leq  c \; \int_{B^{\rm ss}(v,\,r)}\int_{-r}^r
\vol_{W^{\rm su}(w)}(B^{\rm su}(\phi_{s}w, r e^{-T})) 
\;ds\, d\vol_{W^{\rm ss}(v)}(w)\;. 
\end{multline}

There exists a compact subset $K'$ in $T^1\wt M$ (depending only on
$K$ and $r$) such that for all $v\in K$,$T\geq 0$, $s\in[-r,r]$ and
$w\in B^{\rm ss} (v,\,r)$ with $\phi_Tv\in \Ga K$, we have $\phi_sw\in
K'$ and $\phi_T(\phi_sw)\in\Ga K'$. We have $B^{\rm su}(\phi_s w,r
e^{-T}) = \phi_{-T}B^{\rm su}(\phi_T(\phi_sw),r)$ by Equation
\eqref{eq:dilatboulhamen}.  By Lemma \ref{lem:gibbsinstable}, for all
$v,T,s,w$ as above, we hence have
$$
\frac{1}{c'_{K',\,r}}\;e^{\int_0^T \wt{F}^{\rm su}(\phi_{t+s}w)\,dt}\leq 
 \vol_{W^{\rm su}(w)}(B^{\rm su}(\phi_{s}w, r e^{-T}))\leq 
c'_{K',\,r} \;e^{\int_0^T \wt{F}^{\rm su}(\phi_{t+s}w)\,dt}\;.
$$
Furthermore, by Equation \eqref{eq:compardistHamen} for the strong
stable leaves, we have
$$
d(\pi(v),\pi(\phi_{s}w))\leq |s|+d(\pi(v),\pi(w))\leq 
r+d_{W^{\rm ss}(v)}(v,w)\leq 2r
$$ 
and
$$
d(\pi(\phi_Tv),\pi(\phi_{T+s}w))\leq |s|+d(\pi(\phi_Tv),\pi(\phi_Tw))
\leq r+e^{-T}d_{W^{\rm ss}(v)}(v,w)\leq 2r\;.
$$
Since the potential $\wt{F}^{\rm su}$ is H\"older-continuous and
bounded, by Lemma \ref{lem:technicholder} applied twice with $r_0=2r$,
there exists a constant $c'>0$ (depending only on $\|\wt{F}^{\rm su}
\|_\infty$ and $r$) such that for all $v\in K$, $s\in[-r,r]$ and $w\in
B^{\rm ss} (v,\,r)$, we have
$$
\Big|\int_0^T \wt{F}^{\rm su}(\phi_{t+s}w)\,dt-\int_0^T \wt{F}^{\rm
  su}(\phi_{t}v)\,dt\,\Big|=
\Big|\int_{\pi(\phi_{s}w)}^{\pi(\phi_{T+s}w)} 
\wt{F}^{\rm su}-\int_{\pi(v)}^{\pi(\phi_{T}v)} \wt{F}^{\rm
  su}\,\Big|\leq c' \;.
$$
Moreover, by arguments already seen in the proof of Lemma
\ref{lem:gibbsinstable} for the case of the strong unstable foliation,
there exists $c''\geq 1$ (depending only on $K$ and $r$) such that for
every $v\in K$,
$$
\frac{1}{c''}\leq \vol_{W^{\rm ss}(v)}(B^{\rm ss}(v,r))\leq c''\;.
$$
By Equation \eqref{eq:volabscont}, Equation
\eqref{eq:dynballdecompfol} therefore holds with $C'_{K,\,r}=c\, c''
e^{c'}\,c'_{K',\,r}\max\{2r,\frac{1}{2r}\}$. \cqfd

\subsection{Conservative Liouville measures are Gibbs 
measures}
\label{subsec:proofLiouvilleisGibbs}

The aim of this subsection is to give a proof of Theorem
\ref{theo:LiouvilleisGibbs}. We hence assume now that its hypotheses
are satisfied.

By Theorem \ref{theo:potentialholder} and the assumptions of Theorem
\ref{theo:LiouvilleisGibbs}, the potential $\wt F^{\,\rm su}$ is
H\"older-continuous and bounded.  Thus, $\delta=\delta_{\Ga,\,F^{\rm
    su}}$ is finite (and furthermore assumed to be nonpositive) and
there exists indeed a Gibbs measure $\wt m_{F^{\rm su}}$ on $T^1\wt M$
for the potential $\wt F^{\,\rm su}$ (see Subsection
\ref{subsec:GibbsSulivanmeasure}). We denote respectively by
$(\mu_{W^{\rm su} (v)})_{v\in T^1\wt M}$ and $(\mu^\iota_{W^{\rm ss} (v)} 
)_{v\in T^1\wt M}$ the conditional measures of $\wt m_{F^{\rm su}}$ 
on the strong unstable and strong stable foliations defined
in Subsection \ref{subsec:condmesgibbs}.

\medskip We have the following result for these conditional measures,
analogous to Lemma \ref{lem:gibbsinstable} for the leafwise Riemannian
measures. It may also be considered as a leafwise version of the Gibbs
property of the Gibbs measure $\wt m_{F^{\rm su}}$ that has been
proved in Proposition \ref{prop:gibbsgibbs}.

\blemm \label{lem:gibbsgibbsinstable} For every compact subset $K$ in
$T^1 \wt M$ and for every $r>0$, there exists a constant $c''_{K,\,r}
>0$ such that for all $v\in \Ga K$ and $T\geq 0$ such that $\phi_T
v\in \Ga K$, we have
$$
\frac{1}{c''_{K,\,r}} \, e^{\int_{0}^T (\wt F^{\,\rm su}(\phi_tv)-\delta)\,dt} \le 
\mu_{W^{\rm su}(v)}(\phi_{-T}B^{\rm su}(\phi_T v,r))
\le c''_{K,\,r}\;  e^{\int_{0}^T (\wt F^{\,\rm su}(\phi_tv)-\delta)\,dt} \,,
$$
and for all $v\in \Ga K$ and $T\geq 0$ such that $\phi_{-T} v\in \Ga
K$, we have
$$
\frac{1}{c''_{K,\,r} }\,  
e^{\int_{-T}^0 (\wt F^{\,\rm su}\circ\iota(\phi_tv)-\delta)\,dt} \le 
\mu^\iota_{W^{\rm ss}(v)}(\phi_{T}B^{\rm ss}(\phi_{-T} v,r))
\le c''_{K,\,r}\;  
e^{\int_{-T}^0 (\wt F^{\,\rm su}\circ\iota(\phi_tv)-\delta)\,dt} \,.
$$ 
\elemm

\dem Let us fix $K,r$ as in the statement. In order to prove the first
claim, observe that, by Equation \eqref{eq:propdemesfortinstabl}, for
all $v\in T^1\wt M$, $r>0$ and $T\geq 0$,
\begin{align}
\mu_{W^{\rm su}(v)}(\phi_{-T}B^{\rm su}(\phi_T v,r))& =
\int_{\phi_T w\in B^{\rm su}(\phi_T v,\,r)} 
\frac{d(\phi_T)_*\mu_{W^{\rm su}(v)}}
{d\mu_{W^{\rm su}(\phi_T v)}}(\phi_Tw)\, 
d\mu_{W^{\rm su}(\phi_T v)}(\phi_Tw)\nonumber\\ &=
\int_{\phi_T w\in B^{\rm su}(\phi_T v,\,r)} e^{\int_0^T(\wt F(\phi_t w)-\delta)\,dt}\, 
d\mu_{W^{\rm su}(\phi_T v)}(\phi_Tw)\;.\label{eq:controlcondmeasgibbsuball}
\end{align}
By convexity and by Equation \eqref{eq:compardistHamen}, we have
$$
d(\pi(v),\pi(w))\leq d(\pi(\phi_Tv),\pi(\phi_Tw))\leq 
d_{W^{\rm su}(\phi_Tv)}(\phi_Tv,\phi_Tw)\leq r
$$ 
for all $v\in T^1\wt M$, $r>0$ and $w\in \phi_{-T}B^{\rm su} (\phi_T
v,\,r)$. Hence, applying Lemma \ref{lem:technicholder} with $r_0=r$
and by the equivariance of $\wt F^{\,\rm su}$, there exists two
positive constants $c_1$ and $c_2$ such that for all $T\geq 0$ and
$v\in T^1\wt M$ such that $v,\phi_Tv\in \Ga K$, for every $w\in
\phi_{-T}B^{\rm su} (\phi_T v,\,r)$, we have
$$
\Big|\int_0^T(\wt F^{\rm su}(\phi_t w)-\delta)\,dt-
\int_0^T(\wt F^{\rm su}(\phi_t v)-\delta)\,dt\Big|\leq 
2c_1r^{c_2}+2r\max_{\pi^{-1}\big(\N_r(\pi(K))\big)}|\wt F^{\,\rm su}-\delta|\;.
$$
The right hand side of this inequality is a finite constant which
depends only $r$ and $K$. By the compactness of $K$ and the change of
leaf property stated in Equation \eqref{eq:proptrmesfortinstabl}, the
numbers $\mu_{W^{\rm su} (\phi_T v)}(B^{\rm su}(\phi_T v,\,r))$, for
all $T\geq 0$ and $v\in T^1\wt M$ such that $v,\phi_Tv\in \Ga K$, are
bounded from above and from below by positive constants depending only
on $K$ and $r$. The first claim of Lemma \ref{lem:gibbsgibbsinstable}
hence follows from Equation \eqref{eq:controlcondmeasgibbsuball}.

The second claim is proved analogously.  \cqfd

\medskip
Theorem \ref{theo:LiouvilleisGibbs} will easily follow from the
following result.

\blemm \label{lem:LiouabscontGibb} 
Under the assumptions of Theorem \ref{theo:LiouvilleisGibbs},
the Liouville measure $\vol_{T^1M}$ is absolutely continuous with
respect to $m_{F^{\rm su}}$.  
\elemm

\dem The idea of the proof is to use the local product structure of
both measures, and to prove that the conditional measures of the
Liouville measure are absolutely continuous with the ones of the Gibbs
measure, using the previous controls on leafwise balls and Vitali type
arguments.

Since $T^1\wt M$ is $\sigma$-compact and by $\Ga$-equivariance, we
only have to prove that for every big enough compact subset $K$ of
$T^1\wt M$, the restriction of $\vol_{T^1\wt M}$ to $K$ is absolutely
continuous with respect to the restriction of $\wt m_{F^{\rm su}}$ to
$K$. Let us fix a compact subset $K$ of $T^1\wt M$, containing the
domain $U$ of a foliated chart of the strong unstable foliation $\wt
\W^{\rm su}$. In particular, $K$ has positive Liouville measure. 
For every $v\in U$, we will denote by $W^{\rm su}_{\rm loc}(v)$ the
local strong unstable leaf through $v$, that is the connected
component of $v$ in the intersection $W^{\rm su}(v) \cap U$.  

Let $U'$ be the set of $v\in U$ such that the (measurable) set of
elements $w\in W^{\rm su}_{\rm loc}(v)$, for which there exists a
sequence $(t_k)_{k\in\NN}$ in $[0,+\infty[$ tending to $+\infty$ with
$\phi_{t_k}w\in \Ga K$ for every $k\in\NN$, has full Riemannian
measure in $W^{\rm su}_{\rm loc}(v)$. By the conservativity of the
Liouville measure of $T^1M$, and therefore its recurrence, since $K$
has positive Liouville measure, and by the local product structure of
the Liouville measure (see Theorem
\ref{theo:transverseabsolutecontinuity}), the set $U'$ has
full Liouville measure in $U$.

\medskip
Let us prove that for every $v\in U'$, the restriction of
$\vol_{W^{\rm su} (v)}$ to $W^{\rm su}_{\rm loc}(v)$ is bounded by a
constant times the restriction of $\mu_{W^{\rm su}(v)}$ to $W^{\rm
  su}_{\rm loc} (v)$. By a similar proof, when restricted to the local
strong stable leaf $W^{\rm ss}_{\rm loc}(v)$ through Liouville-almost
every $v$, the measure $\vol_{W^{\rm ss}(v)}$ is bounded by a constant
times $\mu^\iota_{W^{\rm ss}(v)}$. Since the Riemannian measure is
absolutely continuous with respect to the product of the Riemannian
measures along the strong stable leaves, strong unstable leaves and
flow lines by Theorem \ref{theo:transverseabsolutecontinuity}, and by
the disintegration property (see Proposition \ref{prop:disintegrGibbs}
(2)) of the Gibbs measure $\wt m_{F^{\rm su}}$ over the strong
unstable measure $\mu_{W^{\rm su}(v)}$, with conditional measures
$\mu^\iota_{W^{\rm s} (w)}$ which are by definition (see Equation
\eqref{eq:defimspourabscont}) absolutely continuous with respect to
the products of the strong stable measures $\mu^\iota_{W^{\rm ss}(w)}$
(which varies absolutely continuously on $w\in W^{\rm su}_{\rm loc}
(v)$ by Equation \eqref{eq:proptrmesstabl}) and the Lebesgue measure
along the flow line, the result follows.

\medskip Let us fix $r\in\mathopen{]}0,1]$. Let $A$ be a fixed
measurable subset of $T^1\wt M$ whose closure is contained in $W^{\rm
  su}_{\rm loc} (v)$. Let $A'=A\cap U'$, which has full Riemannian
measure in $A$, by the definition of $U'$. Let $V$ be any open
neighbourhood of $A$ in $W^{\rm su}_{\rm loc}(v)$.

For every $w\in A$, there exists $T_w\geq 0$ such that if $t\geq T_w$,
the set $\phi_{-t}B^{\rm su}(\phi_tw,r)$, which is equal to $B^{\rm
  su} (w,e^{-t}r)$ by Equation \eqref{eq:dilatboulhamen}, is contained
in $V$.  By the definition of $U'$, for every $w\in A'$, we may choose
$t_w\ge T_w$, such that $\phi_{t_w}w\in \Ga K$. In particular, $B^{\rm
  su}(w,e^{-t_w}r)$ is contained in $V$ for every $w\in A'$. Up to
enlarging slightly $K$, we can assume that $t_w\in\NN$ for all $w\in
A'$.

By a Vitali type of argument, let us now construct a family
$(w_i)_{i\in I}$ in $A'$ indexed by a finite or countable initial
segment $I$ in $\NN$, such that, with $t_i=t_{w_i}$ to simplify the
notation, the family $(B_i=B^{\rm su} (w_i, e^{-t_{i}}\,r))_{i\in I}$
has pairwise disjoint elements, and such that the family $(B'_i=
B^{\rm su} (w_i, 2\,e^{-t_{i}}\,r))_{i\in I}$ covers $A'$.

We proceed by induction. If $A'=\emptyset$, take $I= \emptyset$.
Otherwise, take $w_0\in A'$ such that the integer $t_0=t_{w_0}$ is
minimal. Assume that $w_0,\dots, w_k$ are constructed. Let $w_{k+1}\in
A'$ with $t_{k+1}=t_{w_{k+1}}$ minimal in the set of $w\in A'$ such
that $B^{\rm su} (w, e^{-t_w}\,r)$ does not meet $\bigcup_{0\leq i\leq
  k} B_i$, if this set is nonempty, otherwise the construction stops
at rank $k$, and we take $I=\{0,\dots, k\}$. Since $V$ is relatively
compact, the sequence $(t_{i})_{i\in I}$, if infinite, tends to
$+\infty$.  The construction then does give a finite or countable
family, since $T^1\wt M$ is separable.  By construction, for every
$w\in A'$, the ball $B^{\rm su} (w,e^{-t_w}\,r)$ meets $B_i$ for some
$i\in I$ and $t_w\geq t_{i}$.  Hence by the triangle inequality
$$
d_{W^{\rm su}(w)}(w,w_i)\leq e^{-t_w}\,r +e^{-t_i}\,r\leq 2\,e^{-t_i}\,r\;.
$$
Thus $w$ belongs to $B^{\rm su}(w_i,
2\,e^{-t_i}\,r)$.  This proves the existence of $I$ as required.

Now, since the elements of $(B_i)_{i\in I}$ are contained in $V$ and
are pairwise disjoint, by the $\sigma$-additivity of $\mu_{W^{\rm su}
  (v)}$, by Lemma \ref{lem:gibbsgibbsinstable}, since $\delta\leq 0$,
by Lemma \ref{lem:gibbsinstable}, since the family $(B'_i)_{i\in I}$
covers $A'$, and since $A'$ has full Riemannian measure in $A$, we
have
\begin{align*}
\mu_{W^{\rm su}(v)}(V) &\geq \mu_{W^{\rm su}(v)}\Big(\bigcup_{i\in I}B_i\Big)=
\sum_{i\in I}\mu_{W^{\rm su}(v)}(B_i) 
\geq \frac{1}{c''_{K,\,r}}\sum_{i\in I} 
e^{\int_0^{t_i} (\wt F^{\rm su}(\phi_s w_i)-\delta)\,ds}\\
&\geq \frac{1}{c''_{K,\,r}}\sum_{i\in I} 
e^{\int_0^{t_i} \wt F^{\rm su}(\phi_s w_i)\,ds}\geq
\frac{1}{c''_{K,\,r}\,c'_{K,\,2r}}\sum_{i\in I} \vol_{W^{\rm su}(v)}(B'_i)
\\& \geq \frac{1}{c''_{K,\,r}\,c'_{K,\,2r}} 
\vol_{W^{\rm su}(v)}\Big(\bigcup_{i\in I} B'_i\Big) \\
& \geq 
\frac{1}{c''_{K,\,r}\,c'_{K,\,2r}} \vol_{W^{\rm su}(v)} (A')=
\frac{1}{c''_{K,\,r}\,c'_{K,\,2r}} \vol_{W^{\rm su}(v)} (A)\;.
\end{align*}
Since $V$ is any open neighbourhood of $A$ in $W^{\rm su}_{loc} (v)$,
and as the constants $c''_{K,\,r}$ and $c'_{K,\,2r}$ do not depend on
$V$, by the regularity of $\mu_{W^{\rm su} (v)}$, we have
$$
\vol_{W^{\rm su}(v)}(A) 
\leq c''_{K,\,r}\, c'_{K,\,2r}\, \mu_{W^{\rm su}(v)} (A)
$$ 
for every Borel subset $A$ whose closure is contained in $W^{\rm
  su}_{loc} (v)$.  This proves the result.  
\cqfd

\medskip Let us finally conclude the proof of Theorem
\ref{theo:LiouvilleisGibbs}.  Recall that $\Omega_c\Gamma$ is the
subset of elements of $T^1M$ which are positively and negatively
recurrent under the geodesic flow. Since the Liouville measure on
$T^1M$ is conservative, the measurable set $\Omega_c\Ga$ has full
Liouville measure, the limit set of $\Ga$ is equal to
$\partial_\infty\wt M$ and the (topological) non-wandering set
$\Omega\Ga$ of the geodesic flow is equal to $T^1M$.

Lemma \ref{lem:LiouabscontGibb} implies that
$m_{F^{\rm su}}(\Omega_c\Ga)>0$, so that $m_ {F^{\rm su}}$ is not completely
dissipative. By the Hopf-Tsuji-Sullivan-Roblin theorems
\ref{theo:critnonergo} and \ref{theo:critergo}, this implies that
$(\Ga,F^{\rm su})$ is not of convergence type, hence is of divergence
type, and that $m_ {F^{\rm su}}$ is ergodic and conservative, hence also
gives full measure to $\Omega_c\Gamma$. A standard ergodicity argument
then implies that $\vol_{T^1M}$ and $m_{F^{\rm su}}$ are proportional, as
required. Their constants for the Gibbs property have to be equal,
which proves (by Theorem \ref{theo:equalcritgur} and by Theorem
\ref{theo:principevariationnel} (1)) the final equalities stated in
Theorem \ref{theo:LiouvilleisGibbs}. \cqfd

\medskip \rem The same proof shows that if the derivatives of the
sectional curvature of $\wt M$ are uniformly bounded, if
$\delta_{\Ga,\,F^{\rm su}}= 0$, and if the geodesic flow of $M$ is
recurrent with respect to the Gibbs measure $m_{F^{\rm su}}$ on
$T^1M$, then $m_{F^{\rm su}}$ is proportional to the conservative part of
the Liouville measure $\vol_{T^1M}$.

Note that the Ahlfors conjecture, proved by Calegari-Gabai
\cite{CalGab06} (after works of Bonahon and Canary) says that if $\wt
M=\HH^3_\RR$ is the real hyperbolic space of dimension $3$, if $\Ga$
is finitely generated, then either $\Lambda\Ga$ has Lebesgue measure
zero, or $\Lambda\Ga=\partial_\infty \HH^3_\RR$.  But when $\Gamma$ is
not finitely generated, it could happen that $\Lambda\Ga$ and
$\partial_\infty\wt{M}-\Lambda\Gamma$ are both of positive Lebesgue
measure in $\partial_\infty \HH^3_\RR$, so that the conservative and
the dissipative part of the Liouville measure $\vol_{T^1M}$ would both
be nontrivial. This explains how it could happen that $m_{F^{\rm su}}$
could be proportional to the conservative part of the Liouville
measure, but not to the Liouville measure.

\medskip A similar proof (using the Furstenberg $2$-recurrence
property defined in the introduction and the global Gibbs property as
in Proposition \ref{prop:LiouvilleJacobien} instead of their foliated
version as in Lemma \ref{lem:gibbsinstable}) also shows the following
fact, generalising known results when $M$ is compact (see for instance
\cite[Sect.~20.3]{KatHas95}).

\bprop \label{prop:uniqGibbpro}
Let $\wt M$ be a complete simply connected Riemannian manifold, with
dimension at least $2$ and pinched sectional curvature at most $-1$
and let $\Ga$ be a nonelementary discrete group of isometries of $\wt
M$. Let $\wt F:T^1\wt M\ra \RR$ be a $\Ga$-invariant
H\"older-continuous potential and let $c\in \RR$. Then there exists,
up to a multiplicative constant, at most one locally finite (Borel,
positive) measure on $T^1M$ invariant under the geodesic flow, ergodic
and $2$-recurrent, which satisfies the Gibbs property (see Definition
\ref{def:gibbsprop}) for the potential $F$ and the constant $c$.
\eprop

\section{Finiteness and mixing of Gibbs states}
\label{sec:finimixGibbs}

Let $(\wt M,\Ga,F)$ be as in the beginning of Chapter
\ref{sec:negacurvnot}: $\wt M$ is a complete simply connected
Riemannian manifold, with dimension at least $2$ and pinched sectional
curvature at most $-1$; $\Ga$ is a nonelementary discrete group of
isometries of $\wt M$; and $\wt F :T^1\wt M\ra \RR$ is a
H\"older-continuous $\Ga$-invariant map.  

Fix $x,y$ in $\wt M$. Assume that $\delta_{\Ga,\,F}<+\infty$ (see
Subsection \ref{subsec:GibbsPoincareseries} for comments on this
condition, which is not so important but without which we can not
even define Gibbs measures). Let $\wt m_F$ be the Gibbs measure on
$T^1\wt M$ associated with a pair of Patterson densities
$(\mu^\iota_{x}) _{x\in\wt M}$ and $(\mu_{x})_{x\in\wt M}$ for
respectively $(\Ga,F\circ\iota)$ and $(\Ga,F)$ of (common) dimension
$\delta_{\Ga,\,F\circ \iota} = \delta_{\Ga,\,F}$.  Denote by $\|m_F\|$
the total mass of the measure $m_F$ on $T^1M$ induced by $\wt m_F$.

Compared to the previous ones, the sections \ref{sec:growth},
\ref{sec:ergtheounistabfolia} and \ref{subsec:classnilpotcov} will
require stronger assumptions on the Gibbs measure $m_F$. We will
assume that it is finite and mixing under the geodesic flow. These
assumptions are already present in the case $F=0$ considered by
\cite{Roblin03}. But essentially, as explained in the coming
Subsection \ref{subsec:mixing}, only the finiteness of $m_F$ is
important (it is a nonempty assumption, even in the case $F=0$, see
for instance \cite{DalOtaPei00}).

Indeed, we have proved in Corollary \ref{coro:uniqGibstate} that, when
it is finite, the Gibbs measure $m_F$ is unique up to scaling and
is ergodic (and the Patterson densities $(\mu^\iota_{x}) _{x\in\wt M}$
and $(\mu_{x})_{x\in\wt M}$ are also unique up to scaling). We will
give one finiteness criterion in Subsection \ref{subsec:finiteness}.
But we prefer to start by stating the main dynamical tool to be used
in the coming chapters, saying that the mixing property is essentially
always true as soon as the Gibbs measure $m_F$ is finite.

\subsection{Babillot's mixing criterion for Gibbs states} 
\label{subsec:mixing}

Since any Gibbs measure is a quasi-product measure (see Subsection
\ref{subsec:GibbsSulivanmeasure} and Subsection
\ref{subsec:condmesgibbs}), we may apply Babillot's result
\cite[Theo.~1]{Babillot02b} to obtain the mixing property of the
geodesic flow. Recall that the {\it length spectrum}\index{length
  spectrum} of $M$ is the subset of $\RR$ consisting of the
translation lengths of the elements of $\Ga$, or, equivalently, of
the lengths of the periodic orbits of the geodesic flow on $T^1M$. A
continuous flow of homeomorphisms $(\varphi_t)_{t\in\RR}$ of a
topological space $X$ is {\it topologically
  mixing}\index{topologically mixing}\index{dynamical system
  (topological)!topologically mixing}\index{mixing (topologically)} if
for all nonempty open subsets $U$ and $V$ of $X$, there exists
$t_0\in\RR$ such that $\varphi_t(U)\cap V\neq\emptyset$ for every
$t\geq t_0$.

\btheo [Babillot]\label{theo:babmix}\index{Babillot!mixing theorem}%
\index{theorem@Theorem!of Babillot!of mixing of the geodesic flow} If
$\delta_{\Ga,\,F}<+\infty$ and if $m_F$ is finite, then the following
conditions are equivalent:
\begin{enumerate}
\item[(1)] the geodesic flow $(\phi_t)_{t\in\RR}$ on $T^1M$ is mixing
  for the measure $m_F$;
\item[(2)] the geodesic flow $(\phi_t)_{t\in\RR}$ on $T^1M$ is
  topologically mixing on its (topological) non-wandering set
  $\Omega\Ga$;
\item[(3)] 
  the length spectrum of $M$ is not contained in a discrete
  subgroup of $\RR$.
\end{enumerate}
\etheo

For instance, the last condition is satisfied if $\wt M$ is a
symmetric space, or if $\Ga$ contains a parabolic element, or if $M$
has dimension $2$, or if $\Lambda\Ga$ is not totally disconnected (see
for instance \cite{DalBo99,DalBo00} and their references).

Conjecturally, the non arithmeticity of the length spectrum (that is,
the validity of the third assertion above) should be always true,
hence the only assumption in this subsection should be the finiteness of
the Gibbs measure.

Also note that, at least since Margulis's thesis \cite{Margulis04},
the mixing hypothesis is standard in obtaining precise counting results
(see \cite{EskMcMul93} and the surveys \cite{Babillot02a,Oh10},
amongst many references).

\medskip In Chapter \ref{sec:ergtheounistabfolia}, to prove Theorem
\ref{theo:uniergodtransv}, we will use the following consequence of
the mixing property of the Gibbs measures.  The forthcoming result of
equidistribution of pieces of strong unstable leaves pushed by the
geodesic flow,  is due to
\cite[Coro.~3.2]{Roblin03} when $F=0$ and to
\cite[Theo.~3]{Babillot02b} under a more general quasi-product
hypothesis, satisfied by our Gibbs measures (see Subsection
\ref{subsec:condmesgibbs}).

For all $v,w$ in the same leaf of $\wt\W^{\rm su}$ in $T^1\wt M$, let us define 
\begin{equation}\label{eq:deficsubwtFbab}
c_{\wt F}(w,v)=C_{F\circ\iota-\delta_{\Ga,\,F},\,w_-}(\pi(v),\pi(w))\;,
\end{equation}
which is equal to $0$ if $\wt F=0$ (see Subsection
\ref{subsec:proofvaraprincip} for the cocycle $c_{F}$ associated to
$\W^{\rm su}$ in $T^1M$, and Subsection \ref{subsec:defiquasinvtransv}
for more information on $c_{\wt F}$).

\btheo [Babillot]\label{theo:mixingbab}%
\index{theorem@Theorem!of Babillot!of equidistribution of strong
  unstable leaves}\index{Babillot!equidistribution theorem} Let $\wt
M$ be a complete simply connected Riemannian manifold, with dimension
at least $2$ and pinched negative curvature at most $1$. Let $\Ga$ be
a nonelementary discrete group of isometries of $\wt M $.  Let $\wt
F:\wt M\ra \RR$ be a $\Gamma$-invariant H\"older-continuous map.
Assume that $\delta=\delta_{\Ga,\,F}$ is finite, and that the Gibbs
measure $m_F$ on $T^1M$ is finite and mixing under the geodesic flow.

Then for every $v\in T^1\wt M$ such that $v_-\in\Lambda\Ga$, for every
relatively compact Borel subset $B$ of $W^{\rm su}(v)$, for every
uniformly continuous and $m_F$-integrable $\psi:T^1M\ra\RR$, if $\wt
\psi=\psi\circ \prT$ where $\prT:T^1\wt M\ra T^1M$ is the canonical
projection, we have
$$
\lim_{t\ra+\infty}\int_B \wt \psi\circ\phi_t(w)\;
e^{c_{\wt F}(w,\,v)}\;d\mu_{W^{\rm su}(v)}(w)=
\int_B  e^{c_{\wt F}(w,\,v)}
\;d\mu_{W^{\rm su}(v)}(w)\;\frac{1}{\|m_F\|}\int_{T^1M}\psi \;dm_F\;.
$$
\etheo

As usual, we may replace continuous functions with compact support by
relatively compact Borel subsets with negligible boundary: Under the
hypotheses of the above theorem, with $v$ and $B$ as above, for every
relatively compact Borel subset $A$ of $T^1\wt M$ whose boundary has
measure $0$ for $\wt m_F$, we have
$$
\lim_{t\ra+\infty}\;\sum_{\ga\in \Ga}\int_{B\cap\phi_{-t}\ga A}
e^{c_{\wt F}(w,\,v)}\;d\mu_{W^{\rm su}(v)}(w)=
\int_B  e^{c_{\wt F}(w,\,v)}
\;d\mu_{W^{\rm su}(v)}(w)\;\frac{\wt m_F(A)}{\|m_F\|}\;.
$$

\dem We only give this proof for completeness, following very closely
Babillot's arguments.

Let $v$, $B$, $\psi$ and $\wt\psi$ be as in the statement. Since
$\lim_{t\ra +\infty} 0=0$, we may assume that $B$ meets the support of
$\mu_{W^{\rm su}(v)}$. By the linearity in $\psi$ of the formula to be
established in Theorem \ref{theo:mixingbab}, we may assume that
$\psi\geq 0$. We want to compute the limit as $t\ra+\infty$ of the
following ratio, whose denominator is different from $0$,
$$ 
R_t=\frac{\displaystyle \int_B \wt{\psi}\circ \phi_t(w)\; 
e^{c_{\wt{F}}(w,\,v)}\; d\mu_{W^{\rm su}(v)}(w)}{\displaystyle 
\int_B e^{c_{\wt{F}}(w,\,v)} \;d\mu_{W^{\rm su}(v)}(w)}\geq 0\;.
$$

We first thicken the subset $B$ of $W^{\rm su}(u)$ to a subset $\wt\B$
of $T^1\wt M$ in the following way. Let $A_+=\{w_+\;:\;w\in B\}$ be
the Borel subset of $\partial_\infty \wt{M}$ consisting of the
positive endpoints of the vectors of $B$. Let $A_-$ be a small enough
compact neighbourhood of $v_-$, disjoint from the closure of $A_+$
(which exists since $A_+$ is relatively compact in $\partial_\infty
\wt{M} -\{v_-\}$), and let $\eta>0$ be small enough. Note that
$\mu^\iota_{x_0}(A_-)>0$ since $v_-$ belongs to $\Lambda\Ga$, which is
the support of $\mu^\iota_{x_0}$ by Corollary \ref{coro:uniqpatdens}
(with the help of Corollary \ref{coro:uniqGibstate}). Using the Hopf
parametrisation $(w_-, w_+, t_w)\in\partial_\infty^2\wt{M}\times \RR$
of $w\in T^1\wt M$ given in Remark (2) in Subsection
\ref{subsec:GibbsSulivanmeasure} with respect to the base point
$x_0=\pi(v)$, let us define now
$$
\wt\B=\{w\in T^1\wt M\;:\; w_-\in A_-,\;\; w_+\in A_+,\;\; 
t_w-t_v\in [0,\eta]\}\;,
$$  
which is a relatively compact Borel subset of $T^1\wt M$. Up to
subdividing $B$ into finitely many Borel subsets, if $A_-$ and $\eta$
are small enough, by the ergodicity of $m_F$ and by the finite
additivity in $B$ of the formula to be established in Theorem
\ref{theo:mixingbab}, we may assume that for $\wt m_F$-almost every
$x\in\wt\B$, exactly $N$ elements of $\wt\B$ have the same image as
$x$ by the canonical projection $T^1\wt M\ra T^1M=\Ga\bs T^1\wt M$,
for some $N\in\NN-\{0\}$. We denote by $\B$ the image of $\wt\B$ by
this projection.

In these Hopf coordinates, the measure $\wt{m}_F$ can be written (see
Remark (2) in Subsection \ref{subsec:GibbsSulivanmeasure}) in the
following quasi-product form
$$
d\wt{m}_F(w)=f(w_-,w_+)
\;d\mu_{x_0}^\iota(w_-)\, d\mu_{x_0} (w_+)\, dt_w\;,
$$ 
where $f=(D_{F-\delta,\,x_0})^{-2}$ is a continuous and positive
function on $\partial_\infty^2\wt{M}$. By Equation
\eqref{eq:formulemesgibbs}, for every $w\in T^1\wt M$, we have
$$
f(w_-,w_+)=
e^{C_{F\circ \iota-\delta,\,w_-} (x_0,\,\pi(w)) +C_{F-\delta,\,w_+}
  (x_0,\,\pi(w))}\;.
$$ 
By Equation \eqref{eq:deficsubwtFbab}, since $d\mu_{W^{\rm
    su}(v)} (w)= e^{C_{F-\delta,\,w_{+}}(x_0,\,\pi(w))}\,
d\mu_{x_0}(w_+)$ by Equation \eqref{eq:defmeasstrongunstable}, since
$w_-=v_-$ and $t_w=t_v$ if $w\in W^{\rm su}(v)$, and since
$\pi(v)=x_0$, we have
$$
R_t = \frac{\displaystyle\int_{A_+} \wt{\psi}\circ 
\phi_t(v_-,w_+,t_v) \;f(v_-,w_+)\;d\mu_{x_0}(w_+)}{\displaystyle
\int_{A_+}  f(v_-,w_+)\;d\mu_{x_0}(w_+)}\;.
$$

Now, let us fix $\epsilon>0$. Note that $\wt{\psi}$ is uniformly
continuous, that $f$ is uniformly continuous and minorated by a
positive constant on compact sets, that every $w\in T^1\wt M$ such
that $w_+\neq v_-$ and the unit tangent vector with Hopf coordinates
$(v_-,w_+,t_v)$ are in the same stable leaf, on which $\phi_t$ for
$t\geq 0$ does not increase the distance $d'_{T^1\wt M}$, and that
$\wt\B$ is relatively compact. Hence, if $A_+$ and $\eta$ are small
enough, then for all $t\ge 0$, $r\in[0,\eta]$ and $w\in\B$, we have
\begin{enumerate}
\item[$\bullet$]~
$e^{-\epsilon}\le \frac{f(w_-,\,w_+)}{f(v_-,\,w_+)}\le 
e^\epsilon$, 
\item[$\bullet$]~ $\Big|\wt{\psi}\circ \phi_t(w_-,w_+,t_v+r)-
\wt{\psi}\circ \phi_t(v_-,w_+,t_v)\Big|\leq \epsilon\;$.
\end{enumerate}
For every $t\geq 0$, let 
$$
R'_t=\frac{\displaystyle \int_{A_-\times A_+ \times [0,\,\eta]} 
\wt{\psi}\circ \phi_t(w_-,w_+,t_v+r) \;f(w_-,w_+)\;
d\mu_{x_0}^\iota(w_-)\,d\mu_{x_0}(w_+)\,dr}{\displaystyle
\int_{A_-\times A_+  \times [0,\,\eta]}  f(w_-,w_+)\;
d\mu_{x_0}^\iota(w_-)\,d\mu_{x_0}(w_+)\,dr}\;.
$$
Using a simplification by $\eta\,\mu^\iota_{x_0}(A_-)>0$, the above
estimations show that 
$$
(R'_t-\epsilon)e^{-2\epsilon}\leq R_t\leq
(R'_t+\epsilon)e^{2\epsilon}
$$ 
for every $t\geq 0$.  Using a simplification by $N$ for the last
equality, we have
$$
R'_t=\frac{\displaystyle \int_{\wt\B}\wt{\psi}\circ\phi_t\;
d\wt{m}_F}{\displaystyle \wt{m}_F(\wt\B)}=
\frac{\displaystyle \int_{\B}\psi\circ\phi_t\;
d m_F}{\displaystyle m_F(\B)}\;.
$$
Since $\psi$ and $\mathbbm{1}_\B$ are $m_F$-integrable, the mixing
property of the geodesic flow with respect to the normalised measure
$\frac{1}{\|m_F\|}m_F$ implies that $R'_t$ converges, as $t\to
+\infty$, towards $ \frac{1}{\|m_F\|} \int_{T^1M}\psi \,dm_F$.  The
conclusion of the proof follows.  \cqfd

\bigskip
 The next subsection gives a finiteness criterion for the Gibbs
measures, hence allowing us to use Babillot's theorems \ref{theo:babmix}
and \ref{theo:mixingbab}.

\subsection{A finiteness criterion for Gibbs states} 
\label{subsec:finiteness}

A point $p\in \partial_\infty\wt M$ is a {\it bounded parabolic fixed
  point}\index{bounded parabolic fixed point}\index{parabolic!fixed
  point!bounded} of $\Ga$ if it is the fixed point of a parabolic
element of $\Ga$ and if its stabiliser $\Ga_p$ in $\Ga$ acts properly
with compact quotient on $\Lambda\Ga-\{p\}$. The discrete
nonelementary group of isometries $\Ga$ of $\wt M$ is said to be {\it
  geometrically finite}\index{geometrically finite} if every element
of $\Lambda\Ga$ is either a conical limit point or a bounded parabolic
fixed point of $\Ga$ (see for instance \cite{Bowditch95}).

\medskip The following finiteness criterion for Gibbs measures is due
to \cite{DalOtaPei00} when $F=0$. With the multiplicative approach
explained at the beginning of Chapter \ref{sec:GPS}, it is due to
Coudène \cite{Coudene03}.

\btheo \label{theo:DOPB}%
\index{theorem@Theorem!of Dal'Bo-Otal-Peigné with potentials}%
\index{Dal'Bo-Otal-Peigné theorem with \\ potentials} Assume that
$\delta_{\Ga,\,F}$ is finite, and that $\Ga$ is geometrically finite
with $(\Ga,F)$ of divergence type.  Then the Gibbs measure $m_F$ is
finite if and only if for every parabolic fixed point $p$ of $\Ga$,
the series
$$
\sum_{\alpha\in\Ga_p} \; 
d(x,\alpha y) \;e^{\int_x^{\alpha y} (\wt F-\delta_{\Ga,\,F})}
$$ 
converges, where $\Ga_p$ is the stabiliser of $p$ in $\Ga$.  
\etheo

\dem We will follow the scheme of proof of
\cite[Theo.~B]{DalOtaPei00}.  Note that the convergence of the above
series depends neither on $x$ nor on $y$, and that the result is
immediate if $\Ga$ has no parabolic element, since then the support of
the Gibbs measure, which is $\Omega\Ga$, is compact. Let
$\delta=\delta_{\Ga,\,F}$.

Let $\Par_\Ga$ be the set of parabolic fixed points of $\Ga$. Since
$\Ga$ is geometrically finite (see for instance \cite{Bowditch95}),
the action of $\Ga$ on $\Par_\Ga$ has only finitely many orbits, and
there exists a $\Ga$-equivariant family $(\H_p)_{p\in\Par_\Ga}$ of
pairwise disjoint closed horoballs, with $\H_p$ centred at $p$, such
that the quotient
$$
M_0=\Ga\bs\big(\C\Lambda\Ga-\bigcup_{p\in\Par_\Ga} \H_p\big)
$$ 
is compact. For every $p\in\Par_\Ga$, we denote by $\Ga_p$ the
stabiliser of $p$ in $\Ga$, and by $\F_p$ a relatively compact
measurable strict fundamental domain for the action of $\Ga_p$ on
$\Lambda\Ga-\{p\}$. The horoball $\H_p$ is {\it precisely
  invariant}\index{precisely invariant}%
\index{horoball!precisely invariant} under the stabiliser of $p$
in $\Ga$: The inclusion $\H_p\subset \wt M$ induces an injection
$(\Ga_p\bs\H_p)\ra(\Ga\bs \wt M)$ and we will identify $\Ga_p\bs\H_p$
with its image. Recall that $\pi:T^1 M\ra M$ is the canonical
projection.

Since the measure $m_F$ on $T^1 M$ is finite if and only if its
push-forward measure $\pi_*m_F$ on $M$ is finite, since the support of
$\pi_*m_F$ is contained in $\Ga\bs\C\Lambda\Ga$, since $M_0$ is
compact and since there are only finitely many orbits of parabolic
fixed points, the Gibbs measure $m_F$ is finite if and only if for
every parabolic fixed point $p$ of $\Ga$, we have
$$
\pi_*m_F(\Ga_p\bs\H_p)<+\infty\;.
$$

Fix such a point $p$. For every $(\xi,\eta)\in\partial^2_\infty\wt M$,
we will denote by $(\xi\eta)$ the geodesic line with endpoints $\xi$
and $\eta$ oriented from $\xi$ to $\eta$; if $\xi\neq p$ and
$(\xi\eta)$ meets $\H_p$, we denote by $x_{\xi,\,\eta}$ the first
intersection point of $(\xi\eta)$ with $\partial \H_p$.  Note
that since $\wt M$ is CAT$(-1)$, if $\xi\neq p$ and if the geodesic
lines $(\xi\eta)$ and $(\xi\eta')$ both meet $\H_p$, then
$d(x_{\xi,\,\eta},x_{\xi,\,\eta'})\leq 2\ln(1+\sqrt{2})$ (see for
instance \cite[Lemma 2.9]{ParPau10GT}).

\begin{center}
\input{fig_paraborne.pstex_t}
\end{center}

We define $x_0=x_{\xi_0,\,p}$ for any fixed $\xi_0\in\F_p$. In particular,
since $\F_p$ is relatively compact in $\partial_\infty\wt M-\{p\}$,
there exists $\kappa\geq 0$ such that for every $\xi\in\F_p$, for
every $\eta\in\partial_\infty\wt M-\{\xi\}$ such that $(\xi\eta)$
meets $\H_p$, we have
\begin{equation}\label{eq:proxixx}
d(x_0,x_{\xi,\,\eta})\leq \kappa\;.
\end{equation}

By the disjointness of the horoballs in the family
$(\H_p)_{p\in\Par_\Ga}$, a geodesic in the support of $\pi_*m_F$
meeting $\Ga_p\bs\H_p$ has a unique lift in $\wt M$ meeting $\H_p$
starting from the fundamental domain $\F_p$, unless it has a lift
which starts from $p$; its other endpoint is in $\alpha\F_p$ for some
$\alpha\in\Ga_p$, unless it is equal to $p$. Since $\Ga$ is of
divergence type, the Patterson densities $\mu_{x_0}$ and
$\mu^\iota_{x_0}$ are atomless (see Proposition
\ref{prop:atomlessreferee}) and the diagonal of $\partial_\infty\wt
M\times\partial_\infty\wt M$ has zero measure with respect to
$\mu^\iota_{x_0}\otimes\mu_{x_0}$ (see Assertion (c) of Proposition
\ref{prop:listequivHTSRG}).  Note that $\mu_{x_0}(\F_p)$ and
$\mu^\iota_{x_0}(\F_p)$ are positive, otherwise $\{p\}$ would be the
support of $\mu_{x_0}$ or $\mu^\iota_{x_0}$, hence would be
$\Ga$-invariant, contradicting the fact that $\Ga$ is nonelementary.

Hence, by the very definition of the Gibbs measure $\wt m_F$, we have
$$
\pi_*m_F(\Ga_p\bs\H_p)=
\sum_{\alpha\in\Ga_p}\int_{(\xi,\,\eta)\in\F_p\times\alpha \F_p}
\Long\big((\xi\eta)\cap\H_p\big)\;
\frac{d\mu^\iota_{x_0}(\xi)\,d\mu_{x_0}(\eta)}
{D_{F-\delta,\,x_0}(\xi,\eta)^2}\;.
$$
By Equation \eqref{eq:ecartvisuel}, by Equation
\eqref{eq:proxixx} and by Lemma \ref{lem:holderconseq} (1), there
exists $c>0$ (which depends only on $\kappa$, the constants in Lemma
\ref{lem:holderconseq} and $\max_{\pi^{-1}(B(x_0,\,\kappa))}|\wt F|$) such
that, for every $(\xi,\eta)\in\F_p\times\partial_\infty\wt M$ such
that $\xi\neq \eta$ and $(\xi\eta)$ meets $\H_p$, we have
$$
\frac{1}{c}\leq \frac{1}{D_{F-\delta,\,x_0}(\xi,\eta)^2}=
e^{C_{F\circ\iota-\delta,\,\xi}(x_0,\,x_{\xi,\,\eta})+
C_{F-\delta,\,\eta}(x_0,\,x_{\xi,\,\eta})}\leq c\;.
$$
Let $(\xi,\eta)\in\F_p\times\alpha \F_p$ be such that $\xi\neq \eta$
and $(\xi\eta)$ meets $\H_p$. The geodesic line from
$\alpha^{-1}\eta\in\F_p$ to $\alpha^{-1}\xi$ also meets $\H_p$. The
exiting point $y_{\xi,\,\eta}$ of $(\xi\eta)$ out of $\H_p$ is equal
to $\alpha\, x_{\alpha^{-1}\eta, \,\alpha^{-1}\xi}$.  Hence $d(\alpha
x_0, y_{\xi,\,\eta}) = d(x_0,x_{\alpha^{-1}\eta,\,\alpha^{-1}\xi})\leq
\kappa$.  Therefore, by the triangle inequality (see the picture
above), we have
$$
d(x_0,\alpha x_0)-2\kappa\leq 
\Long\big((\xi\eta)\cap\H_p\big)\leq d(x_0,\alpha x_0)+2\kappa\;.
$$
Hence 
\begin{align}
\frac{1}{c}\;\mu^\iota_{x_0}(\F_p)&\sum_{\alpha\in\Ga_p}\big(d(x_0,\alpha x_0)-
2\kappa\big)\mu_{x_0}(\alpha\F_p)\nonumber\\ \label{eq:dunepart}&
\leq\pi_*m_F(\Ga_p\bs\H_p) \leq
c\;\mu^\iota_{x_0}(\F_p)\sum_{\alpha\in\Ga_p}\big(d(x_0,\alpha x_0)+
2\kappa\big)\mu_{x_0}(\alpha\F_p)\;.
\end{align}
By the equations \eqref{eq:equivdensity} and
\eqref{eq:radonykodensity}, for every $\alpha\in\Ga_p$, we have
\begin{equation}\label{eq:ombre}
\mu_{x_0}(\alpha\F_p)=\mu_{\alpha^{-1} x_0}(\F_p)=
\int_{\xi\in \F_p} e^{-C_{F-\delta,\,\xi}(\alpha^{-1} x_0,\,x_0)}\;d\mu_{x_0}(\xi)\;.
\end{equation}
When $\alpha$ goes to infinity in the discrete group $\Ga_p$, the point
$\alpha^{-1} x_0$ converges to $p$. Hence, for every $\xi\in\F_p$,
except for finitely many $\alpha$, the set $\F_p$ is contained in the
shadow seen from $\alpha^{-1} x_0$ of $B(x_0,\kappa+1)$. Therefore, by
Lemma \ref{lem:holderconseq} (2), there exists a constant $c'>0$ such
that for all $\xi\in\F_p$ and  $\alpha\in\Ga_p$, we have
\begin{equation}\label{eq:reserviraDOPMoh}
\Big|\;C_{F-\delta,\,\xi}(\alpha^{-1} x_0,x_0)+\int_{\alpha^{-1} x_0}^{x_0} (\wt
F-\delta)\;\Big|\leq c'\;.
\end{equation}
Hence by Equation \eqref{eq:ombre} and by invariance, we have 
\begin{equation}\label{eq:delautre}
e^{-c'}\mu_{x_0}(\F_p)\;e^{\int_{x_0}^{\alpha x_0} (\wt F-\delta)}\leq
\mu_{x_0}(\alpha\F_p)\leq 
e^{c'}\mu_{x_0}(\F_p)\;e^{\int_{x_0}^{\alpha x_0} (\wt F-\delta)}\;.
\end{equation}
Putting together the equations \eqref{eq:dunepart} and
\eqref{eq:delautre}, there exists $C>0$ such that
\begin{align*}
\frac{1}{C}\;\sum_{\alpha\in\Ga_p}\big(d(x_0,\alpha x_0)-
C\big)e^{\int_{x_0}^{\alpha x_0} (\wt F-\delta)}&\leq\pi_*m_F(\Ga_p\bs\H_p)
\\ 
& \leq C\;\sum_{\alpha\in\Ga_p}\big(d(x_0,\alpha x_0)+
C\big)e^{\int_{x_0}^{\alpha x_0} (\wt F-\delta)}\;.
\end{align*}
This proves the result. \cqfd

\medskip A criterion for $(\Ga,F)$ to be of divergence type, when
$\Ga$ is geometrically finite, is the following one, again due to
\cite{DalOtaPei00} when $F=0$.

\btheo \label{theo:DOPA} Assume that $\delta_{\Ga,\,F}$ is finite and
that $\Ga$ is geometrically finite. If for every parabolic fixed point
$p$ of $\Ga$ with stabiliser $\Ga_p$, we have
$\delta_{\Ga_p,\,F}<\delta_{\Ga,\,F}$, then $(\Ga,F)$ is of divergence
type.  
\etheo

\dem Let $(\mu_x)_{x\in \wt M}$ be the Patterson density for $(\Ga,F)$
of dimension $\delta_{\Ga,\,F}$ constructed in Proposition
\ref{prop:existPattdens}: we use the notation $h,s_k,
\overline{Q}_{z,\,y}(s), \mu_{z,\,s}$ for $k\in\NN,
s>\delta_{\Ga,\,F},\,z\in\wt M$ introduced for that purpose. We fix an
arbitrary $y\in\wt M$. We start with an independent lemma (which is
necessary, by Proposition \ref{prop:atomlessreferee}).

\blemm \label{lem:nonatompara}
Let $p$ be a bounded parabolic fixed point of $\Ga$ such that
$\delta_{\Ga_p,\,F}<\delta_{\Ga,\,F}$.  Then $\mu_x(\{p\})=0$.  
\elemm
 
\dem Let $\H_p,\F_p,\kappa, \xi_0,x_0$ be as in the above proof of
Theorem \ref{theo:DOPB}. We may assume that $\xi_0\notin \Ga p$. Up to
replacing $\H_p$ by a smaller horoball centred at $p$, we may assume
that the closest point to $y$ on the geodesic line $\mathopen{]}\xi_0,
p\mathclose{[}$ between $\xi_0$ and $p$ does not belong to the
geodesic ray $\mathopen{[}x_0,p\mathclose{[}\,$. In particular, by the
convexity of the horoballs and since $\ga\H_p\cap \H_p$ is empty if
$\ga\notin\Ga_p$, the orbit $\Ga y$ does not meet $\H_p$.

For every $\ga\in\Ga$, choose a representative $\overline{\ga}$ of the
left coset of $\ga$ in $\Ga_p\bs\Ga$ such that $\overline{\ga}\xi_0
\in\F_p$ (which is possible since $\ga\xi_0\in \Lambda\Ga-\{p\}$ and
$\F_p$ is a fundamental domain for the action of $\Ga_p$ on
$\Lambda\Ga-\{p\}$).

Let $\epsilon>0$ such that
$\delta_{\Ga,\,F}>\delta_{\Ga_p,\,F}+\epsilon$.  Let $r_\epsilon\geq
0$ such that $h(t+r)\leq e^{\epsilon t}h(r)$ for all $t\geq 0$ and
$r\geq r_\epsilon$. There exists $\kappa'\geq \kappa$ large enough so
that if $z$ is the point of $\mathopen{[}x_0,p\mathclose{[}x$ at
distance $\max\{r_\epsilon, \kappa'+2\}$ of $x_0$, then for every
$\ga\in\Ga$, there exists a representative $\overline{\ga}$ of the
left coset of $\ga$ in $\Ga_p\bs\Ga$ such that $\ov\ga \,y\in
\OOO_zB(x_0,\kappa')$. In particular, $d(z,\Ga y)\geq r_\epsilon$. Let
$\epsilon'\in \mathopen{[}0, 1\mathclose{]}$ be such that
$$
\mu_z\big(\partial(\OOO_zB(x_0, \kappa'+\epsilon'))\big)=0\;.
$$

\begin{center}
\input{fig_paraborne2.pstex_t}
\end{center}

Let $X$ be the compact subset $\Ga y\cup\Lambda\Ga$ of $\wt
M\cup\partial_\infty \wt M$, which contains the support of the
measures $\mu_{z,\,s}$ and $\mu_z$ for all $z\in \wt M$ and
$s>\delta_{\Ga,\,F}$. Let 
$$
\C=X\cap
\bigcup_{\xi\in \OOO_zB(x_0,\,\kappa'+\epsilon')}\mathopen{[}z,\xi\mathclose{]}
$$ 
be the set of points of $X$ on a geodesic ray, including its point at
infinity, from $z$ through the closed ball
$B(x_0,\kappa'+\epsilon')$.

By Equation \eqref{eq:proxixx} and by convexity, the set $\F_p$ is
contained in the shadow $\OOO_zB(x_0,\kappa)$, hence in $\C$. Since
the point $\overline{\ga}y$ belongs to $\C$ for every $\ga\in\Ga$, we
have $X=\{p\}\cup\bigcup_{\alpha\in\Ga_p} \alpha\C$.

Let $(\alpha_i)_{i\in\NN}$ be an enumeration of $\Ga_p$, and for every
$n\in\NN$, let $U_n=X-\bigcup_{0\leq i\leq n} \alpha_i\C$. Then
$(U_n)_{n\in\NN}$ is a non-increasing sequence of neighbourhoods of $p$
in $X$, with intersection $\{p\}$.  Hence
$$
\mu_z(\{p\})=\lim_{n\ra +\infty} \mu_z(U_n)\;.
$$

By the convexity of horoballs and a standard argument of
quasi-geodesics, there exists a constant $c>0$ such that for $i$ large
enough and for every $\ga\in\Ga$, the point $\alpha_iz$ is at distance
at most $c$ from the geodesic segment $\mathopen{[}z,\alpha_i\ov\ga
y\mathclose{]}$ (see the picture above). Hence, by (two applications
of) Lemma \ref{lem:technicholder}, there exists a constant $c'>0$ such
that, for every large enough $i$, for every $\ga\in\Ga$, and for every
$s\in\mathopen{[}\delta_{\Ga,\,F},\delta_{\Ga,\,F}+1\mathclose{]}$, we
have
$$
\Big|\;\int_z^{\alpha_i\overline{\ga}y}(\wt F-s)-
\int_{z}^{\alpha_i z}(\wt F-s)\;-\int_z^{\overline{\ga}y}(\wt F-s)\;\Big|
\leq c'\;.
$$
Since $h$ is nondecreasing, by the triangle inequality and since
$d(z,\Ga y)\geq r_\epsilon$, we have
$$
h(d(z,\alpha_i\ov\ga y))\leq h(d(z,\alpha_iz)+d(z,\ov\ga y))\leq
e^{\epsilon\,d(z,\,\alpha_i z)}\;h(d(z,\ov\ga y))\;.
$$

For every $s\in
\mathopen{]}\delta_{\Ga,\,F},\delta_{\Ga,\,F}+1\mathclose{]}$, we
hence have, for every $n$ large enough,
\begin{align*}
\mu_{z,\,s}(U_n)&\leq \sum_{i=n+1}^{+\infty}\;\;\sum_{\ov\ga\,\in\;\Ga_p\bs\Ga}\;
\frac{1}{\ov Q_{y,\,y}(s)}\;h(d(z,\alpha_i\ov\ga y))
\;e^{\int_z^{\alpha_i\overline{\ga}y}(\wt F-s)}\\ & \leq
e^{c'}\sum_{i=n+1}^{+\infty}\;
e^{\int_z^{\alpha_i z}(\wt F-(s-\epsilon))}
\;\;\sum_{\ov\ga\,\in\;\Ga_p\bs\Ga}\;
\frac{1}{\ov Q_{y,\,y}(s)}\;h(d(z,\ov\ga y))
\;e^{\int_z^{\overline{\ga}y}(\wt F-s)}\\ & \leq e^{c'}\;\|\mu_{z,s}\|
\sum_{i=n+1}^{+\infty}\;e^{\int_z^{\alpha_i z}(\wt F-(s-\epsilon))}\;.
\end{align*}
Taking $s=s_k$ for $k$ large enough, and letting $k$ tends to
$+\infty$, we get (note that the constant function $1$ on
$\partial_\infty \wt M$ is continuous with compact support and that
$\mu_z(\partial U_n)=0$ for every $n\in\NN$,  by the
choice of $\epsilon'$)
$$
\mu_{z}(U_n)\leq e^{c'}\;\|\mu_{z}\|
\sum_{i=n+1}^{+\infty}\;
e^{\int_z^{\alpha_i z}(\wt F-(\delta_{\Ga,\,F}-\epsilon))}\;.
$$
Since $\delta_{\Ga,\,F}-\epsilon>\delta_{\Ga_p,\,F}$ and since the
remainder of a convergent series tends to $0$, we have
$$
\mu_z(\{p\})=\lim_{n\ra +\infty} \mu_z(U_n)=0\;.
$$
Since $\mu_x$ and $\mu_z$ are absolutely continuous with respect
to each other, Lemma \ref{lem:nonatompara} follows.
\cqfd

\medskip Now, since the limit points of $\Ga$ that are not conical
ones are (countably many) bounded parabolic fixed points if $\Ga$ is
geometrically finite, this lemma implies that $\mu_x(\Lambda_c\Ga)=
\mu_x(\Lambda\Ga)>0$. Theorem \ref{theo:DOPA} hence follows from
Corollary \ref{coro:critdivergence}.  \cqfd

\bcoro Assume that $\delta_{\Ga,\,F}$ is finite, that $\Ga$
is geometrically finite, and that $\delta_{\Ga_p,\,F}<\delta_{\Ga,\,F}$
for every parabolic fixed point $p$, where $\Ga_p$ is the stabiliser
of $p$ in $\Ga$. Then the Gibbs measure $m_F$ is finite.  
\ecoro

In particular, if $\Ga$ is geometrically finite, if
$\delta_{\Ga_p}<\delta_{\Ga}$ for every $p\in\Par_\Ga$ (which is in
particular the case when $\wt M$ is a symmetric space), and if
$\|F\|_\infty$ is small enough (that is, if $\|F\|_\infty\leq
\frac{1}{2}(\delta_\Ga-\delta_{\Ga_p})$ for every $p\in\Par_\Ga$), then
$m_F$ is finite, by Lemma \ref{lem:elemproppressure} (iv).

\medskip
\dem This is immediate by Theorems
\ref{theo:DOPA} and \ref{theo:DOPB}. \cqfd

\medskip As in \cite{Coudene03}, this allows us to construct finite
(hence ergodic) Gibbs measures for any geometrically finite group
$\Ga$, by choosing an appropriate potential $F$ in the preimage in
$T^1M$ of the cuspidal parts of $M$.

\section{Growth and equidistribution of orbits and periods}
\label{sec:growth}

We will give in this chapter precise asymptotic results as $t$ goes to
$+\infty$ for the counting functions defined in Subsection
\ref{subsec:countfunct}, as corollaries of convergence results for
suitable measures and equidistribution results.

Let $(\wt M,\Ga,F)$ be as in the beginning of Chapter
\ref{sec:negacurvnot}: $\wt M$ is a complete simply connected
Riemannian manifold, with dimension at least $2$ and pinched sectional
curvature at most $-1$; $\Ga$ is a nonelementary discrete group of
isometries of $\wt M$; and $\wt F :T^1\wt M\ra \RR$ is a
H\"older-continuous $\Ga$-invariant map.  Fix $x,y$ in $\wt M$.
Assume that $\delta_{\Ga,\,F}<+\infty$. Let $\wt m_F$ be the Gibbs
measure on $T^1\wt M$ associated with a pair of Patterson densities
$(\mu^\iota_{x}) _{x\in\wt M}$ and $(\mu_{x})_{x\in\wt M}$ for
respectively $(\Ga,F\circ\iota)$ and $(\Ga,F)$ of (common) dimension
$\delta_{\Ga,\,F\circ \iota} = \delta_{\Ga,\,F}$.  Denote by $\|m_F\|$
the total mass of the measure $m_F$ on $T^1M=\Ga\bs T^1\wt M$ induced
by $\wt m_F$.

\subsection{Convergence of measures on the square product} 
\label{subsec:productconv}

The aim of this subsection is to prove the following weak-star
convergence result, whose consequences will be discussed in the next
subsections.

\btheo\label{maingrowthRoblin}%
\index{theorem@Theorem!Two point orbital \\equidistribution} Assume
that the critical exponent $\delta_{\Ga,\,F}$ of $(\Ga,F)$ is finite
and positive, and that the Gibbs measure $m_F$ is finite and mixing
under the geodesic flow on $T^1M$. As $t$ goes to $+\infty$, the
measures
$$
\nu_{x,\,y,\,F,\,t}=\delta_{\Ga,\,F}\;\|m_F\|\;e^{-\delta_{\Ga,\,F}\; t}
\;\sum_{\ga\in\Ga\;:\;d(x,\ga y)\leq t} \;
e^{\int_{x}^{\ga y} \wt F}\;\D_{\ga^{-1}x}\otimes\D_{\ga y}
$$
converge to the product measure $\mu^\iota_{y}\otimes\mu_{x}$ with
respect to the weak-star convergence of measures on $(\wt
M\cup\partial_\infty \wt M)^2$.  
\etheo

We will follow closely the proof for the case $F=0$ given in
\cite[Chap.~4]{Roblin03}. It relies on the next technical proposition,
improved in Proposition \ref{prop:steptwo}. To simplify the notation,
we set $\nu_T=\nu_{x,\,y,\,F,\,T}$ for every $T\geq 0$. We recall the
notation $\C^\pm_r(z,Z)$, $K^\pm(z,r,Z)$, $K(z,r)$, $v_{\xi,\,\eta,\,z}$,
$L_r(z,w)$ and $\OOO^\pm_r(z,w)$ of Subsection
\ref{subsec:geomnotation}, as well as the notation of the antipodal
map $\iota_z:\partial_\infty\wt M\ra \partial_\infty\wt M$ with
respect to any $z\in\wt M$ defined in Subsection
\ref{subsec:unittanbun}.

\bprop\label{prop:stepone} Assume that $\delta_{\Ga,\,F}$ is finite and
positive, and that $m_F$ is finite and mixing under the geodesic flow
on $T^1M$. For every $\epsilon>0$, for all $\xi_0$ and $\eta_0$ in
$\partial_\infty \wt M$ such that $\iota_x(\xi_0)\in \Lambda\Gamma$
and $\iota_y(\eta_0)\in \Lambda\Gamma$, there exist open
neighbourhoods $V$ and $W$ of respectively $\xi_0$ and $\eta_0$ in
$\partial_\infty \wt M$ such that for all Borel subsets $A\subset V$
and $B\subset W$,
$$
\limsup_{T\ra+\infty}\;\;\nu_T\big(\C^-_1(y,B)\times \C^-_1(x,A)\big)
\;\leq\; e^\epsilon\,\mu^\iota_{y}(B)\;\mu_{x}(A)
$$
$$
{\rm and}\;\;\;\liminf_{T\ra+\infty}\;\;\nu_T\big(\C^+_1(y,B)\times 
\C^+_1(x,A)\big)
\;\geq\; e^{-\epsilon}\;\mu^\iota_{y}(B)\;\mu_{x}(A)\;.
$$
\eprop

\dem To simplify the notation, let $\delta= \delta_{\Ga,\,F}>0$. We
will denote by $\epsilon_1,\epsilon_2,\epsilon_3, ...$ positive
functions of $\epsilon$ (depending only on $\delta$, on
$\max_{\pi^{-1}(B(z,\,2))}|\wt F|$ for $z=x,y$, and on the constants
$c_1,c_2, c_3,c_4>0$ appearing in Lemma \ref{lem:holderconseq} for
$\wt F$ and $\wt F\circ\iota$), that converge to $0$ as $\epsilon$
goes to $0$, and whose exact computation is, though possible,
unnecessary.

Let $\epsilon, \xi_0,\eta_0$ be as in the statement of Proposition
\ref{prop:stepone}. We first define the neighbourhoods $V$ and $W$
required in the statement of Proposition \ref{prop:stepone}.

Let $r$ be in $\mathopen{]}0,\min\{1,\epsilon\}\mathclose{[}$ possibly
outside a countable subset. By the finiteness of the Patterson
densities, we may assume that
\begin{equation}\label{eq:boundarynonvanishing}
\mu^\iota_{x}(\partial\OOO_{\xi_0}B(x,r))=
\mu_{y}(\partial\OOO_{\eta_0}B(y,r))=0\;.
\end{equation}
Since the support of the Patterson densities is $\Lambda\Ga$ (since
$m_F$ is finite, see Corollary \ref{coro:uniqpatdens} and
\ref{coro:uniqGibstate}), and since $\iota_x(\xi_0)$ and
$\iota_y(\eta_0)$ belong to $\Lambda\Ga$, we have
\begin{equation}\label{eq:defiCr}
C_r=\mu^\iota_{x}(\OOO_{\xi_0}B(x,r))\;
\mu_{y}(\OOO_{\eta_0}B(y,r))>0\;.
\end{equation}
By Equation \eqref{eq:boundarynonvanishing} and the convergence
property of $\OOO_r^\pm(z,r)$ to $\OOO_{\xi_0}B(x,r)$ as $z\ra\xi_0$
described in Subsection \ref{subsec:geomnotation}, there exist open
neighbourhoods $\wh V,\wh W$ of (respectively) $\xi_0,\eta_0$ in $\wt
M\cup\partial_\infty \wt M$ such that for every $r$ as above, for
all $z\in \wh V$ and $w\in \wh W$,
\begin{equation}\label{eq:defiwhV}
e^{-\epsilon}\mu^\iota_{x}(\OOO_{\xi_0}B(x,r))\leq 
\mu^\iota_{x}(\OOO^\pm_{r}(z,x))\leq 
e^{\epsilon}\mu^\iota_{x}(\OOO_{\xi_0}B(x,r))
\end{equation}
and
\begin{equation}\label{eq:defiwhW}
e^{-\epsilon}\mu_{y}(\OOO_{\eta_0}B(y,r))\leq 
\mu_{y}(\OOO^\pm_{r}(w,y))\leq 
e^{\epsilon}\mu_{y}(\OOO_{\eta_0}B(y,r))\;.
\end{equation}
Finally, let $V,W$ be open neighbourhoods of $\xi_0,\eta_0$ in
$\partial_\infty \wt M$ whose closures are contained in $\wh V,\wh W$
respectively.

\bigskip For all Borel subsets $A$ in $V$ and $B$ in $W$, we defined (see
Subsection \ref{subsec:geomnotation})
$$
K^+=K^+(x,r,A)\;\;\;{\rm and}\;\;\;K^-=K^-(y,r,B)\;.
$$
The heart of the proof is to give two pairs of upper and lower bounds
(respectively in Equation \eqref{eq:*1*2} and Equation
\eqref{eq:minomF}) of
$$
I_\pm(T)=\int_0^{T\pm 3r} \;e^{\delta\, t}\;\sum_{\ga\in \Ga} \;
\wt m_F(K^+\cap \phi_{-t}\ga K^-)\;dt\;.
$$
Both of them use the following estimate: For every $(\xi,\eta)
\in \partial^2_\infty\wt M$ such that the orthogonal projection $w$ of
$x$ on the geodesic line between $\xi$ and $\eta$ satisfies
$d(w,x)\leq r\leq\min\{1,\epsilon\}$, and in particular if
$(\xi,\eta)\in L_r(x,\ga y)$, we have, by Equation
\eqref{eq:ecartvisuel} and Lemma \ref{lem:holderconseq} (1),
\begin{equation}\label{eq:controlFgap}
e^{-\epsilon_1}\leq D_{F-\delta,\,x}(\xi,\eta)=
e^{-\frac{1}{2}(C_{F-\delta,\,\eta}(x,\,w)+C_{F\circ \iota-\delta,\,\xi}(x,\,w))}
\leq e^{\epsilon_1}\;.
\end{equation}

\medskip
\noindent {\bf First upper and lower bound. } We use the
definition of the Gibbs measure to estimate individually
$\wt m_F(K^+\cap \phi_{-t}\ga K^-)$.

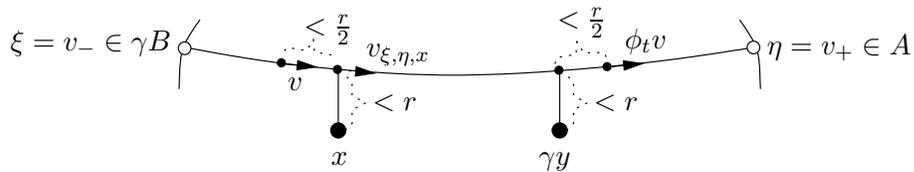
\begin{center}
\input{fig_kpcapkm.pstex_t}
\end{center}

\medskip For every $\ga\in\Ga$ such that $d(x,\ga y)>2r$ and for every
$t\geq 0$, we have, by the definition of the various geometric sets in
Subsection \ref{subsec:geomnotation}, that $K^+\cap \phi_{-t}\ga K^-$
is the set of $v\in T^1\wt M$ such that $v$ is $r$-close to $x$,
$v_+\in A$, $v_-\in \ga B$ and $\phi_t v$ is $r$-close to $\ga y$
(see the picture above).  Furthermore, by the definition of $\wt
m_{F}$,

\medskip
$\wt m_F(K^+\cap \phi_{-t}\ga K^-)=$
$$
\int_{(\xi,\,\eta)\in L_{r}(x,\,\ga y)\cap(\ga B\times A)}
\frac{d\mu^\iota_{x}(\xi)\;
d\mu_{y}(\eta)}{D_{F-\delta,\,x}(\xi,\eta)^2}\;
\int_{-\frac{r}{2}}^{\frac{r}{2}}
\mathbbm{1}_{K(\ga y,\,r)}(\phi_{t+s}v_{\xi,\,\eta,\,x})\;ds\;,
$$
where $\mathbbm{1}_Z$ is the characteristic function of a subset $Z$.

As in \cite[pages 59, 60]{Roblin03} which only uses an estimate as in
Equation \eqref{eq:controlFgap} and geometric arguments, there exists
a constant $c_1>0$ such that for every $T>3r$, if
$$
J_\pm(T)=\sum_{\mbox{\tiny$\begin{array}{c}
\ga\in\Ga\;:\;d(x,\ga y)\leq T,\\ 
\ga^{-1} x\in\C^\pm_1(y,B)\cap\wh W,\\
\ga y\in\C^\pm_1(x,A)\cap\wh V\end{array}$}}
\mu^\iota_{x}\big(\OOO^\pm_r(\ga y,
x)\big)\;\mu_{x}\big(\OOO^\pm_r(x,\ga y)\big)e^{\delta\, d(x,\ga y)}\;,
$$
then the following two inequalities hold
\begin{equation}\label{eq:*1*2}
I_-(T)\leq  \;e^{\epsilon_2}\;r^2\;J_+(T)\;+c_1\;\;\;{\rm and}\;\;\;
I_+(T)\geq  \;e^{-\epsilon_2}\;r^2\;J_-(T)\;\;-c_1\;.
\end{equation}

\medskip
\noindent {\bf Second upper and lower bound. } Now, we use the
mixing property of the geodesic flow to estimate the sum
$\sum_{\ga\in \Ga} \wt m_F(K^+\cap \phi_{-t}\ga K^-)$.

\medskip Since the geodesic flow is mixing (the definitions in Section
\ref{subsec:pushmeas} take care of the possible presence of torsion in
$\Ga$), we have, for $t$ large enough,
\begin{equation}\label{eq:mixingquoi}
e^{-\epsilon}\;\wt m_F(K^+)\;\wt m_F(K^-) \;\leq \; \|m_F\|\;
\sum_{\ga\in \Ga} \wt m_F(K^+\cap \phi_{-t}\ga K^-)\;\leq 
\;e^{\epsilon}\;\wt m_F(K^+)\;\wt m_F(K^-)\;.
\end{equation}
By the definition of $K^+=K^+(x,r,A)$, we have
$$
\wt m_F(K^+)=r\int_{\eta\in A}\;d\mu_{x}(\eta)
\int_{\zeta\in \OOO_{\eta}B(x,\,r)}
\;D_{F-\delta,\,x}(\zeta,\eta)^{-2}\;d\mu^\iota_{x}(\zeta)\;.
$$
By Equation \eqref{eq:controlFgap} and by Equation \eqref{eq:defiwhV}
applied to $z\in A$ since $A\subset \wh V$, we have
$$
e^{-\epsilon-2\epsilon_1}\;r\;
\mu_{x}(A)\;\mu^\iota_{x}(\OOO_{\xi_0}B(x,r))
\;\leq \wt m_F(K^+)\leq \;
e^{\epsilon+2\epsilon_1}\;r\;
\mu_{x}(A)\;\mu^\iota_{x}(\OOO_{\xi_0}B(x,r))\;.
$$
Similarly 
$$
e^{-\epsilon-2\epsilon_1}\;r\;
\mu^\iota_{y}(B)\;\mu_{y}(\OOO_{\eta_0}B(y,r))
\;\leq \wt m_F(K^-)\leq \;
e^{\epsilon+2\epsilon_1}\;r\;
\mu^\iota_{y}(B)\;\mu_{y}(\OOO_{\eta_0}B(y,r))\;.
$$
By taking the product of these inequalities, by multiplying Equation
\eqref{eq:mixingquoi} by $e^{\delta t}$ and integrating it over $t\in
\mathopen{[}0,T\pm3r\mathclose{]}$, we have since $\delta>0$ and
$r\leq \epsilon$, for some constant $c_2$ independent of $T$ (coming
from the fact that the estimations above are valid only for $t$ large
enough, hence we have to cut the integral over $t\in \mathopen{[}0,
T\pm3r\mathclose{]}$ in two), and by the definition of $C_r$ in
Equation \eqref{eq:defiCr}, for $T$ large enough,
\begin{equation}\label{eq:minomF}
\delta \; \|m_F\|\;I_-(T)\geq e^{-\epsilon_3}\;r^2\;C_r\;
\mu^\iota_{y}(B)\;\mu_{x}(A)\;e^{\delta T} \;\;- c_2
\end{equation}
$$
{\rm and}\;\;\;\delta \; \|m_F\|\;I_+(T)\leq e^{\epsilon_3}\;r^2\;C_r\;
\mu^\iota_{y}(B)\;\mu_{x}(A)\;e^{\delta T} \;\;+ c_2\;.
$$

\medskip
\noindent{\bf Estimate on $C_r$. }
Before showing how these two upper and lower bounds of $I_\pm(T)$
prove the result, let us give a lower and upper estimate on $C_r$.

\medskip By the defining properties of the Patterson densities (see
the equations \eqref{eq:equivdensity} and \eqref{eq:radonykodensity}),
we have
$$
\mu_{y}\big(\OOO^\pm_r(\ga^{-1}x,y)\big)
=\mu_{\ga y}\big(\OOO^\pm_r(x,\ga y)\big)= 
\int_{\eta\in \OOO^\pm_r(x,\ga y)}e^{C_{F-\delta,\,\eta}(x,\ga y)}
\;d\mu_{x}(\eta)\;.
$$
Note that $\OOO^\pm_r(x,\ga y)\subset\OOO_xB(\ga y,2r)$ by the
convexity of the distance. Hence, by Lemma \ref{lem:holderconseq} (2)
applied to $\wt F-\delta$ and $y'=\ga y$, since $\max_{\pi^{-1}(B(\ga
  y,\,2r))}|\wt F|\leq \max_{\pi^{-1}(B(y,\,2))}|\wt F|$ by
invariance, and since $r\leq \epsilon$, we have
$$ 
e^{-\epsilon_4}\;\mu_{y}\big(\OOO^\pm_r(\ga^{-1}x,y)\big)\;\leq\;
\mu_{x}\big(\OOO^\pm_r(x,\ga y)\big)\;e^{-\int_x^{\ga y} (\wt
  F-\delta)}\;\leq\;
e^{\epsilon_4}\;\mu_{y}\big(\OOO^\pm_r(\ga^{-1}x,y)\big)\;.
$$
If $(\ga y,\ga^{-1}x)\in\wh V\times \wh W$, we have, respectively by
the definition of $C_r$ (see Equation \eqref{eq:defiCr}), the
definitions of $\wh V$ and $\wh W$ (see the equations
\eqref{eq:defiwhV} and \eqref{eq:defiwhW}), and the previous
inequality,
\begin{align*}
C_r&=\mu^\iota_{x}(\OOO_{\xi_0}B(x,r))\;\mu_{y}(\OOO_{\eta_0}B(y,r))\\
& \geq e^{-2\epsilon}\;\mu^\iota_{x}(\OOO^+_{r}(\ga y,x))\;
\mu_{y}(\OOO^+_{r}(\ga^{-1}x,y))\\
& \geq e^{-2\epsilon-\epsilon_4}\;\mu^\iota_{x}(\OOO^+_{r}(\ga y,x))\;
\mu_{x}\big(\OOO^+_r(x,\ga y)\big)
\;e^{-\int_x^{\ga y} (\wt F-\delta)}\;.
\end{align*}
Similarly,
$$
C_r\leq e^{2\epsilon+\epsilon_4}\;
\mu^\iota_{x}(\OOO^-_{r}(\ga y,x))\;
\mu_{x}\big(\OOO^-_r(x,\ga y)\big)\;
e^{-\int_x^{\ga y} (\wt F-\delta)}\;.
$$

\medskip
\noindent{\bf Conclusion of the proof. }
Now, respectively by the definition of the measure $\nu_T=
\nu_{x,\,y,\,F,\,T}$, by using some constant $c_3>0$ independent of $T$, by
the previous lower bound on $C_r$ and the definition of $J_\pm(T)$, by
Equation \eqref{eq:*1*2}, and by Equation \eqref{eq:minomF} with some
constant $c_4$ independent of $T$, we have
\begin{align*}
C_r\;\nu_{T}\big(\C^+_1(y,B)\times \C^+_1(x,A)\big)&
=\delta \; \|m_F\|\;e^{-\delta T}\;\sum_{\mbox{\tiny$\begin{array}{c}
\ga\in\Ga\;:\;d(x,\ga y)\leq T,\\ 
(\ga^{-1} x,\ga y)\in\C^+_1(y,B)\times\C^+_1(x,A)\end{array}$}}
C_r\;e^{\int_x^{\ga y}\wt F}\\
&\geq \delta \; \|m_F\|\;e^{-\delta T}\;
\sum_{\mbox{\tiny$\begin{array}{c}
\ga\in\Ga\;:\;d(x,\ga y)\leq T,\\ 
\ga^{-1} x\in\C^+_1(y,B)\cap\wh W,\\
\ga y\in\C^+_1(x,A)\cap\wh V\end{array}$}}
C_r\;e^{\int_x^{\ga y}\wt F}\;\; -c_3\;e^{-\delta T}\\
& \geq \delta \; \|m_F\|\;e^{-\delta T}
\;e^{-2\epsilon -\epsilon_4}\;J_+(T) -c_3\;e^{-\delta T}\\ & \geq 
\delta \; \|m_F\|\;e^{-\delta T}\;e^{-\epsilon_2-2\epsilon -\epsilon_4}\;r^{-2}
\big( I_-(T) -c_1\big)-c_3\;e^{-\delta T}\\
& \geq e^{-\epsilon_5}\;C_r\;
\mu^\iota_{y}(B)\;\mu_{x}(A) -c_4\;e^{-\delta T}\;.
\end{align*}
Similarly, 
$$
C_r\;\nu_{T}\big(\C^-_1(y,B)\times \C^-_1(x,A)\big)
\leq e^{\epsilon_5}\;C_r\;
\mu^\iota_{y}(B)\;\mu_{x}(A) +c_4\;e^{-\delta T}\;.
$$
Dividing by $C_r>0$ (independent of $T$), the result follows.
\cqfd

\medskip The next result is a small improvement of Proposition
\ref{prop:stepone} which gets rid of the assumption that
$\iota_x(\xi_0),\iota_y(\eta_0)\in \Lambda\Gamma$, at the price of
replacing $1$-thickened/$1$-thinned cones by $R$-thickened/$R$-thinned
ones, for some $R>0$ large enough. Note that when
$\Lambda\Ga=\partial_\infty \wt M$, then this step is not necessary
(that is $R=1$ works). In particular the reader interested only in the
case when $\Ga$ is cocompact may skip this proposition.

\bprop\label{prop:steptwo} Assume that $\delta_{\Ga,\,F}$ is finite
and positive, and that $m_F$ is finite and mixing under the geodesic
flow on $T^1M$. For every $\epsilon'>0$, for all $\xi_0$ and
$\eta_0$ in $\partial_\infty \wt M$, there exist $R>0$ and open
neighbourhoods $V$ and $W$ of respectively $\xi_0$ and $\eta_0$ in
$\partial_\infty \wt M$ such that for all Borel subsets $A\subset V$
and $B\subset W$,
$$
\limsup_{T\ra+\infty}\;\;\nu_T\big(\C^-_R(y,B)\times \C^-_R(x,A)\big)
\;\leq\; e^{\epsilon'}\,\mu^\iota_{x}(B)\;\mu_{x}(A)
$$
$$
{\rm and}\;\;\;\liminf_{T\ra+\infty}\;\;\nu_T\big(\C^+_R(y,B)\times 
\C^+_R(x,A)\big)
\;\geq\; e^{-\epsilon'}\,\mu^\iota_{x}(B)\;\mu_{x}(A)\;.
$$
\eprop

\dem The idea of the proof is that two cones over the same subset of
$\partial_\infty \wt M$, with vertices two distinct points in $\wt M$,
are asymptotic near infinity, and that we have an appropriate change
of base point formula for the Patterson measures.

Let $\epsilon>0$ and $\xi_0,\eta_0\in\partial_\infty \wt M$. Fix
$\zeta_0\in\Lambda\Ga-\{\xi_0,\eta_0\}$ and $x_0$ (respectively $y_0$)
a point on the geodesic line between $\zeta_0$ and $\xi_0$
(respectively $\zeta_0$ and $\eta_0$). Let 
$$
R=1+\max\{d(x,x_0),\,d(y,y_0)\}
$$ 
and, to simplify the notation, let $\delta=
\delta_{\Ga,\,F}$.

Let us now define the neighbourhoods $V$ and $W$ required by the
statement.  Let $V_0,W_0$ be small enough open neighbourhoods of
$\xi_0,\eta_0$ respectively in $\partial_\infty \wt M$, so that
Proposition \ref{prop:stepone} holds with $x_0,y_0,V_0,W_0$ replacing
$x,y,V,W$ respectively.  Let $\wh{V_0},\wh{W_0}$ be open
neighbourhoods of $\xi_0,\eta_0$ in $\wt M\cup\partial_\infty \wt M$
such that $V_0$ contains $\wh{V_0}\cap \partial_\infty \wt M$ and
$W_0$ contains $\wh{W_0}\cap \partial_\infty \wt M$, respectively, and
such that for all $z\in \wh{V_0}\cap\wt M$, $w\in \wh{W_0}\cap\wt M$,
\begin{equation}\label{eq:approxbusemanncocycle}
|\;d(x_0,z)-d(x,z)-\beta_{\xi_0}(x_0,x)\;|\leq \epsilon\;,\;\;\;
|\;d(y_0,w)-d(y,w)-\beta_{\eta_0}(y_0,y)\;|\leq \epsilon\;,
\end{equation}
\begin{equation}\label{eq:approxgibbscocycle}
\Big|\;\int_{x}^{z}\wt F-\int_{x_0}^{z}\wt F+C_{F,\,\xi_0}(x,x_0)\;\Big|
\leq \epsilon\;\;\;{\rm and}\;\;\;
\Big|\;\int_{y}^{w}\wt F\circ\iota-\int_{y_0}^{w}\wt F\circ\iota+
C_{F\circ\iota,\,\eta_0}(y,y_0)\;\Big| \leq \epsilon\;.
\end{equation}
Let us consider neighbourhoods $V$ and $W$ of $\xi_0$ and $\eta_0$,
respectively, in $\partial_\infty \wt M$ whose closures are contained
in $\wh{V_0}$ and $\wh{W_0}$, respectively.

Now, we start the proof with some computations. If $\ga y_0,\ga
y\in\wh{V_0}$ and $\ga^{-1}x_0,\ga^{-1}x\in\wh{W_0} $, the formulae
\eqref{eq:approxbusemanncocycle} and \eqref{eq:approxgibbscocycle}
imply that
\begin{align*}
d(x_0,\ga y_0)\leq d(x,\ga y_0) +\beta_{\xi_0}(x_0,x) +\epsilon&
=d(y_0,\ga^{-1}x ) +\beta_{\xi_0}(x_0,x) +\epsilon\\ & \leq 
d(y,\ga^{-1}x) +\beta_{\eta_0}(y_0,y) +\beta_{\xi_0}(x_0,x)
+2\epsilon\\ & =
d(x,\ga y ) -\beta_{\xi_0}(x,x_0)-\beta_{\eta_0}(y,y_0) +2\epsilon\;,
\end{align*}
and that, by Equation \eqref{eq:timereversal},
\begin{align*}
  e^{\int_{x}^{\ga y}\wt F}&\leq e^{\int_{x_0}^{\ga y}\wt F -
    C_{F,\,\xi_0}(x,\,x_0)+\epsilon}= e^{\int_{y}^{\ga^{-1} x_0}\wt
    F\circ\iota-C_{F,\,\xi_0}(x,\,x_0)+\epsilon} \\ &
  \leq e^{\int_{y_0}^{\ga^{-1} x_0}\wt F\circ\iota -
    C_{F\circ\iota,\,\eta_0}(y,\,y_0)- C_{F,\,\xi_0}(x,\,x_0) +2\epsilon}
  =e^{\int_{x_0}^{\ga y_0}\wt F}\;e^{-C_{F\circ\iota,\,\eta_0}(y,\,y_0)-
    C_{F,\,\xi_0}(x,\,x_0)+2\epsilon}\;.
\end{align*}

Let $\wh{V_{-R}}=\{z\in\wt M\;:\; B(z,R)\subset \wh{V_0}\}$ and
$\wh{W_{-R}}=\{z\in\wt M\;:\; B(z,R)\subset \wh{W_0}\}$, and let $A$
and $B$ be Borel subsets of $V$ and $W$ respectively.  Note that if
$$
(\ga^{-1}x,\ga y)\in\big(\C^-_R(y,B)\times
\C^-_R(x,A)\big)\cap\big(\wh{W_{-R}}\times \wh{V_{-R}}\big)\;,
$$ 
then 
$$
(\ga^{-1}x_0,\ga y_0)\in\big(\C^-_1(y_0,B)\times
\C^-_1(x_0,A)\big)\cap\big(\wh{W_0}\times \wh{V_{0}}\big)
$$ 
by the choice of $R$ and the definition of the
$r$-thickened/$r$-thinned cones $\C^\pm_r(z,Z)$, and that
$\big(\C^-_R(y,B)\times \C^-_R(x,A)\big)\setminus \big(\wh{W_{-R}}
\times \wh{V_{-R}}\big)$ is a bounded subset of $\wt M\times\wt M$.

Now, by the definition of the measure $\nu_T=\nu_{x,\,y,\,F,\,T}$,
and by the previous computations, for some constant $c'_1>0$
independent of $T$, we have
\begin{align*}
  \mbox{}&\nu_{T}\big(\C^-_R(y,B)\times \C^-_R(x,A)\big)\\
  &\;\;\;\;\;\;\;\;\;\;\;\;\;\;\;\;\;\;=\delta \; 
\|m_F\|\;e^{-\delta T}\;
  \sum_{\mbox{\tiny$\begin{array}{c}
        \ga\in\Ga\;:\;d(x,\ga y)\leq T,\\
        (\ga^{-1} x,\ga y) \in \C^-_R(y,B) \times\C^-_R(x,A)
      \end{array}$}} \;e^{\int_x^{\ga y}\wt F} \\ 
&\;\;\;\;\;\;\;\;\;\;\;\;\;\;\;\;\;\; \leq
  \nu_{x_0,\,y_0,\,F,\,T-\beta_{\xi_0}(x,\,x_0)-\beta_{\eta_0}(y,\,y_0)
    +2\epsilon}\big(\C^-_1(y_0,B)\times \C^-_1(x_0,A)\big) \\
  &\;\;\;\;\;\;\;\;\;\;\;\;\;\;\;\;\;\;\;\;\;
\times\;e^{-C_{F\circ\iota,\,\eta_0}(y,\,y_0)-
    C_{F,\,\xi_0}(x,\,x_0)+2\epsilon}\;
  e^{-\delta\beta_{\xi_0}(x,\,x_0)-\delta\beta_{\eta_0}(y,\,y_0)
    +2\epsilon}\;\; +c'_1e^{-\delta T}\;.
\end{align*}
By Equation \eqref{eq:translatgibbscocycle} and Proposition
\ref{prop:stepone}, we then have
\begin{align*}\mbox{}&\limsup_{T\ra+\infty}\;
\nu_{T}\big(\C^-_R(y,B)\times \C^-_R(x,A)\big)\\ &
\;\;\;\;\;\;\;\;\;\;\;\;\;\;\;\;\;\;\;\;\;\;\;\;\;\;\;
\leq e^{5\epsilon}\;\mu^\iota_{y_0}(B)\;\mu_{x_0}(A)\;
e^{-C_{F\circ\iota-\delta,\,\eta_0}(y,\,y_0)}\;e^{
    -C_{F-\delta,\,\xi_0}(x,\,x_0)}\;.
\end{align*}
By continuity in the equations \eqref{eq:approxbusemanncocycle} and
\eqref{eq:approxgibbscocycle}, since $A\subset V\subset \wh{V_0}$ and
$B\subset W\subset \wh{W_0}$, we have, for all $\xi\in A$ and
$\eta\in B$,
$$
|C_{F-\delta,\,\xi}(x,x_0)-C_{F-\delta,\,\xi_0}(x,x_0)|\leq 
(1+\delta)\epsilon
$$
and
$$
|C_{F\circ\iota-\delta,\,\eta}(y,y_0)-
C_{F\circ\iota-\delta,\,\eta_0}(y,y_0)|
\leq (1+\delta)\epsilon\;.
$$
By Equation \eqref{eq:radonykodensity} for the Patterson densities
$(\mu_{x})_{x\in \wt M}$ and $(\mu^\iota_{x})_{x\in \wt M}$, we hence have
$$
\mu_{x_0}(A)\;e^{-C_{F-\delta,\,\xi_0}(x,x_0)}\leq
\mu_{x}(A)\;e^{(1+\delta)\epsilon}
$$
and
$$
\mu^\iota_{y_0}(B)\;
e^{-C_{F\circ\iota-\delta,\,\eta_0}(y,y_0)}\leq \mu^\iota_{y}(B)\;
e^{(1+\delta)\epsilon}\;.
$$
Therefore 
$$
\limsup_{T\ra+\infty}\;
\nu_{T}\big(\C^-_R(y,B)\times \C^-_R(x,A)\big)\;\leq 
\;e^{(7+2\delta)\epsilon}\;\mu^\iota_{y}(B)\;\mu_{x}(A)\;.
$$
The analogous lower bound is proved similarly, and Proposition 
\ref{prop:steptwo} follows.
\cqfd

\medskip
\noindent{\bf Proof of Theorem \ref{maingrowthRoblin}.}
The end of the proof of Theorem \ref{maingrowthRoblin} follows from
Proposition \ref{prop:steptwo} exactly as in the ``Troisi\`eme \'etape
: conclusion'' in \cite[page 62-63]{Roblin03}, by replacing $\mu_x$
there by $\mu_{x}$, as well as $\mu_y$ by $\mu^\iota_{y}$, 
$\nu^t_{x,y}$ by $\nu_{x,\,y,\,F,\,t}$, and $r$ by $R$.  \cqfd

\medskip As mentioned by the referee (who provided its proof), we have
the following version of Theorem \ref{maingrowthRoblin} when the
critical exponent is possibly nonpositive.  We use the convention that
$\frac{1-e^{-s}}{s}=1$ if $s=0$, and denote by $\weakstar$ the
weak-star convergence of measures.

\btheo\label{maingrowthRoblinbis}%
\index{theorem@Theorem!Two point orbital \\equidistribution} Assume that
$\delta_{\Ga,\,F}<+\infty$ and that the Gibbs measure $m_F$ is finite
and mixing under the geodesic flow on $T^1M$. For every $c>0$, as $t$
goes to $+\infty$, we have
$$
\frac{\|m_F\|\,\delta_{\Ga,\,F}}{1-e^{-c\delta_{\Ga,\,F}}}
\;e^{-\delta_{\Ga,\,F}\; t}
\;\sum_{\ga\in\Ga\;:\;t-c<d(x,\ga y)\leq t} \;
e^{\int_{x}^{\ga y} \wt F}\;\D_{\ga^{-1}x}\otimes\D_{\ga y}
\;\;\weakstar\;\;\mu^\iota_{y}\otimes\mu_{x}\,.
$$
\etheo

Note that Theorem \ref{maingrowthRoblin} follows from this result by
an easy geometric series argument.

\medskip
\dem The proof of Theorem \ref{maingrowthRoblin} could be adapted to
prove this version. But by applying Theorem \ref{maingrowthRoblin} by
replacing $F$ by $F+\kappa$ for $\kappa$ large enough so that
$\delta_{\Ga,\,F+\kappa}=\delta_{\Ga,\,F}+\kappa>0$, since proving the
weak-star convergence of measures can be checked by evaluation on
nonnegative continuous maps with compact support having positive
integral for the limit measure, Theorem \ref{maingrowthRoblinbis} is
also a consequence of the following classical lemma.

\blemm \label{lem:annulaireferee} Let $I$ be a discrete set and
$f,g:I\ra\mathopen{[}0,+\infty \mathclose{[}$ be maps with $f$
proper. Assume that there exist $\delta,\kappa\in\RR$ with
$\delta+\kappa>0$ such that, as $t\ra+\infty$,
$$
\sum_{i\in I,\;f(i) \leq t} e^{\kappa f(i)} 
\,g(i)\;\;\sim\;\; \frac{e^{(\delta+\kappa)t}}{\delta+\kappa}\,.
$$
Then for every $c>0$, we have, as $t\ra+\infty$,
$$
\sum_{i\in I,\;t-c<f(i) \leq t}
\,g(i) \;\;\sim\;\; \frac{1-e^{-c\,\delta}}{\delta}\;e^{\delta \,t}\,.
$$ 
\elemm

\dem For every $\epsilon>0$ small enough, if $s$ is large enough, we have
$$
e^{-\epsilon^2}\;\frac{e^{(\delta+\kappa)s}}{\delta+\kappa}\leq
\sum_{i\in I,\;f(i) \leq s} e^{\kappa f(i)} \,g(i)\leq 
e^{\epsilon^2}\;\frac{e^{(\delta+\kappa)s}}{\delta+\kappa}\;.
$$
Applying this at $s$ and $s-\epsilon$ and subtracting, we have, for
some $A>0$ large enough independent of $s$ and $\epsilon$,
$$
e^{-A\epsilon}\;\epsilon\;e^{\delta\,s}\leq
\sum_{i\in I,\;s-\epsilon< f(i) \leq s} g(i)\leq 
e^{A\epsilon}\;\epsilon\;e^{\delta\,s}\;.
$$
Let $n\in\NN$ and $t$ be large enough. Applying this with
$\epsilon=\frac{c}{n}$ and $s= t -\frac{ck}{n}$ for $k=0,\dots,n-1$
and summing over $k$, a geometric series argument gives
$$
e^{-A\frac{c}{n}}\;\frac{c}{n}\;\frac{1-e^{-c\,\delta}}{1-e^{-\frac{c\delta}{n}}}
\;e^{\delta\,t}\leq \sum_{i\in I,\;t-c< f(i) \leq t} g(i)\leq 
e^{A\frac{c}{n}}\;\frac{c}{n}\;\frac{1-e^{-c\,\delta}}{1-e^{-\frac{c\delta}{n}}}
\;e^{\delta\,t}\;,
$$
say when $\delta\neq 0$, but with the conventions obvious by
continuity, this also holds when $\delta= 0$. Since
$\frac{c(1-e^{-c\,\delta})}{n(1-e^{-\frac{c\delta}{n}})}\sim
\frac{1-e^{-c\,\delta}}{\delta}$ as $n$ goes to $+\infty$, the
result follows.  \cqfd \cqfd

\subsection{Counting orbit points of discrete groups} 
\label{subsec:growth}

We deduce in this subsection some corollaries of Theorem
\ref{maingrowthRoblin}, keeping the notation of the beginning of
Chapter \ref{sec:growth}.

\bcoro\label{coro:coroun} Assume that the critical exponent
$\delta_{\Ga,\,F}$ of $(\Ga,F)$ is finite and positive, and that the
Gibbs measure $m_F$ is finite and mixing under the geodesic
flow in $T^1M$. As $t$ goes to $+\infty$, the measures
$$
\delta_{\Ga,\,F}\;\|m_F\|\;e^{- \delta_{\Ga,\,F}\;t}\;
\sum_{\ga\in\Ga\;:\;d(x,\ga y)\leq t} \;
e^{\int_{x}^{\ga y} \wt F}\;\D_{\ga y}
$$
converge to the measure $\|\mu^\iota_{y}\|\;\mu_{x}$ for the weak-star
convergence of measures on $\wt M\cup\partial_\infty \wt M$, and the
measures
$$
\delta_{\Ga,\,F}\;\|m_F\|\;e^{- \delta_{\Ga,\,F}\;t}\;
\sum_{\ga\in\Ga\;:\;d(x,\ga y)\leq t} \;
e^{\int_{x}^{\ga y} \wt F}\;\D_{\ga^{-1} x}
$$
similarly converge to the measure $\|\mu_{x}\|\;\mu^\iota_{y}$.
\ecoro

\dem Since the push-forward of measures by a continuous map is linear
and weak-star continuous, the result follows from Theorem
\ref{maingrowthRoblin} using the projections $(\wt
M\cup\partial_\infty \wt M)^2\ra (\wt M\cup\partial_\infty \wt M)$ on
the first factor and on the second factor.  \cqfd

\medskip The second assertion of Corollary 1.4 may also be obtained by
exchanging $x$ and $y$ and $F$ and $F\circ\iota$ (see Equation
\eqref{eq:timereversal} and the last remark of Subsection
\ref{subsec:uniqueness}).

\bcoro\label{coro:corodeux} \index{theorem@Theorem!Orbital counting}
Assume that the critical exponent $\delta_{\Ga,\,F}$ of $(\Ga,F)$ is
finite and positive, and that the Gibbs measure $m_F$ is finite and
mixing under the geodesic flow on $T^1M$. As $t$ goes to $+\infty$,
$$
\sum_{\ga\in\Ga\;:\;d(x,\ga y)\leq t} \; e^{\int_{x}^{\ga y} \wt
  F}\;\sim\;\frac{\|\mu^\iota_{y}\|\;\|\mu_{x}\|}
{\delta_{\Ga,\,F}\;\|m_F\|}\;e^{\delta_{\Ga,\,F}\;t} 
\;. 
$$ 
\ecoro

\dem This follows by taking the total mass of the above measures on
the compact space $\wt M\cup\partial_\infty \wt M$. 
\cqfd

\medskip As another immediate corollary of Corollary
\ref{coro:coroun}, we obtain the following sharp asymptotic property
on the sectorial orbital counting function $G_{\Ga,\,F,\,x,\,y,\,U}$
(introduced in Subsection \ref{subsec:countfunct}), improving
Corollary \ref{coro:loggrowth} (1).

\bcoro\label{coro:corotrois} Assume that the critical exponent
$\delta_{\Ga,\,F}$ of $(\Ga,F)$ is finite and positive, and that the
Gibbs measure $m_F$ is finite and mixing under the geodesic flow on
$T^1M$. Let $U$ be an open subset of $\partial_\infty\wt M$ meeting
$\Lambda\Ga$ such that $\mu_x(\partial U)=0$. Then as $t$ goes to
$+\infty$,
$$
G_{\Ga,\,F,\,x,\,y,\,U}(t)\;\sim\;\frac{\|\mu^\iota_{y}\|\;\mu_{x}(U)}
{\delta_{\Ga,\,F}\;\|m_F\|}\;e^{\delta_{\Ga,\,F}\;t} 
\;. 
$$ 
\ecoro

\dem This follows by considering the characteristic functions of cones
$\C_xU$ on open subsets $U$ of $\partial_\infty\wt M$ (or by taking
$V=\partial_\infty\wt M$ in the next corollary). 
\cqfd

\medskip 
Similarly, we obtain the following sharp asymptotic property on the
bisectorial orbital counting function $G_{\Ga,\,F,\,x,\,y,\,U,\,V}$
(introduced in Subsection \ref{subsec:countfunct}), improving
Corollary \ref{coro:loggrowth} (2).

\bcoro\label{coro:coroquatre} Assume that the critical exponent
$\delta_{\Ga,\,F}$ of $(\Ga,F)$ is finite and positive, and that the
Gibbs measure $m_F$ is finite and mixing under the geodesic flow on
$T^1M$. Let $U$ and $V$ be two open subsets of $\partial_\infty\wt M$
meeting $\Lambda\Ga$ such that $\mu_x(\partial U)=\mu^\iota_y(\partial
V)=0$. Then as $t$ goes to $+\infty$,
$$
G_{\Ga,\,F,\,x,\,y,\,U,\,V}(t)\;\sim\;\frac{\mu_{x}(U)\;\mu^\iota_{y}(V)}
{\delta_{\Ga,\,F}\;\|m_F\|}\;e^{\delta_{\Ga,\,F}\;t} 
\;.
$$ 
\ecoro

\dem Use Theorem \ref{maingrowthRoblin} and consider the
product of the characteristic functions of the cones $\C_xU$ and $\C_yV$ on
the open subsets $U$ and $V$ of $\partial_\infty\wt M$. \cqfd

\medskip In the same way Theorem \ref{maingrowthRoblinbis} has been
deduced from Theorem \ref{maingrowthRoblin}, the following result,
which no longer assumes $\delta_{\Ga,\,F}>0$, can be deduced from
Corollaries \ref{coro:coroun}, \ref{coro:corodeux},
\ref{coro:corotrois} and \ref{coro:coroquatre}.

\bcoro\label{coro:corocinq} Assume that $\delta_{\Ga,\,F}<+\infty$,
and that the Gibbs measure $m_F$ is finite and mixing under the
geodesic flow on $T^1M$. Let $U$ and $V$ be two open subsets of
$\partial_\infty\wt M$ meeting $\Lambda\Ga$ such that $\mu_x(\partial
U)=\mu^\iota_x(\partial V)=0$. Then for every $c>0$, as $t$ goes to
$+\infty$,
$$
\frac{\delta_{\Ga,\,F}\;\|m_F\|}{1-e^{-c\,\delta_{\Ga,\,F}}}\;
e^{- \delta_{\Ga,\,F}\;t}\;\sum_{\ga\in\Ga\;:\;t-c<d(x,\ga y)\leq t} \;
e^{\int_{x}^{\ga y} \wt F}\;\D_{\ga y}
\;\;\weakstar\;\;\|\mu^\iota_{y}\|\;\mu_{x}\,;
$$
$$
\frac{\delta_{\Ga,\,F}\;\|m_F\|}{1-e^{-c\,\delta_{\Ga,\,F}}}\;
e^{- \delta_{\Ga,\,F}\;t}\;\sum_{\ga\in\Ga\;:\;t-c<d(x,\ga y)\leq t} \;
e^{\int_{x}^{\ga y} \wt F}\;\D_{\ga^{-1} x}
\;\;\weakstar\;\;\|\mu_{x}\|\;\mu^\iota_{y}\,;
$$
$$
\sum_{\ga\in\Ga\;:\;t-c<d(x,\ga y)\leq t} \; e^{\int_{x}^{\ga y} \wt
  F}\;\sim\;\frac{\|\mu^\iota_{y}\|\;\|\mu_{x}\|\;
(1-e^{-c\,\delta_{\Ga,\,F}})}{\delta_{\Ga,\,F}\;\|m_F\|}\;
e^{\delta_{\Ga,\,F}\;t} \;;
$$ 
$$
G_{\Ga,\,F,\,x,\,y,\,U,\,c}(t)\;\sim\;\frac{\|\mu^\iota_{y}\|\;\mu_{x}(U)\;
(1-e^{-c\,\delta_{\Ga,\,F}})}{\delta_{\Ga,\,F}\;\|m_F\|}\;
e^{\delta_{\Ga,\,F}\;t} \;;
$$ 
$$
G_{\Ga,\,F,\,x,\,y,\,U,\,V,\,c}(t)\;\sim\;\frac{\mu_{x}(U)\;\mu^\iota_{y}(V)\;
(1-e^{-c\,\delta_{\Ga,\,F}})}{\delta_{\Ga,\,F}\;\|m_F\|}\;
e^{\delta_{\Ga,\,F}\;t} 
\;.
$$ 
\ecoro

\subsection{Equidistribution and counting of periodic 
orbits of the  geodesic flow}
\label{subsec:equicountperiodorb}

The aim of this subsection is to use Theorem
\ref{maingrowthRoblin} to prove that the periodic orbits of the
geodesic flow in $T^1M$, appropriately weighted by their periods with
respect to the potential, equidistribute with respect to the Gibbs
measure. When $F=0$ and $M$ is convex-cocompact, the result is due to
Bowen \cite{Bowen72a,Bowen72b,Bowen73}, see also \cite{Parry84}. We
follow Roblin's proof when $F=0$ in \cite{Roblin03}.

Recall (see Subsection \ref{subsec:countfunct}) that for every
$t\in\RR$ and for every periodic (we will only consider primitive ones
in this subsection) orbit $g$ of length $\ell(g)$ of the geodesic flow
on $T^1M$, we denote by $\L_g$ the Lebesgue measure along $g$, by
$\int_g F=\L_g(F)$ the period of $g$ for the potential $F$, by
$\Per'(t)$ the set
of primitive periodic orbits of the geodesic flow on
$T^1M$ with length at most $t\in\RR$.

As in Bowen's two equidistribution results in the convex-cocompact
case, the first assertion of the theorem below claims the
equidistribution of the Lebesgue measures of the closed orbits when
weighted by the exponential of their periods for the potential, the
second one claims the equidistribution of their Lebesgue means.

\btheo\label{theo:equidisperorb}%
\index{theorem@Theorem!Equidistribution of periodic \\orbits} Assume
that the critical exponent $\delta_{\Ga,\,F}$ of $(\Ga,F)$ is
finite and positive, and that the Gibbs measure $m_F$ is finite and
mixing under the geodesic flow of $T^1M$.

(1) As $t$ goes to $+\infty$, the measures
$$
\delta_{\Ga,\,F}\;\;e^{-\delta_{\Ga,\,F}\; t}\;
\sum_{g\in\Per'(t)} \;e^{\L_g(F)}\;\L_{g}
$$
converge to the normalised Gibbs measure $\frac{m_F}{\|m_F\|}$ with
respect to the weak-star convergence of measures on $T^1 M$.

(2) As $t$ goes to $+\infty$, the measures
$$
\delta_{\Ga,\,F}\; t\;\;e^{-\delta_{\Ga,\,F}\; t}\;\sum_{g\in\Per'(t)} 
\;e^{\L_g(F)}\;\frac{\L_{g}}{\ell(g)}
$$
weak-star converge to $\frac{m_F}{\|m_F\|}$.
\etheo

\medskip \dem Let $\delta=\delta_{\Ga,\,F}$. Let us denote by $\H_\Ga$
the set of loxodromic elements of $\Ga$. For every $\ga\in \H_\Ga$,
let $\ell(\ga)$ be its translation length, let $\Axe_\ga$ be its
translation axis, oriented from the repulsive fixed point $\ga_-$ to
the attractive one $\ga_+$, and let $\L_\ga$ be the measure on $T^1\wt
M$ which is the lift to $T^1\wt M$ of the Lebesgue measure along the
translation axis of $\ga$: For every map $f:T^1\wt M\ra\RR$ which is
continuous with compact support, we have 
$$
\L_\ga(f)= \int_{t\in\RR} f(\phi_tv)\;dt
$$ 
where $v$ is any unit tangent vector to the oriented geodesic line
$\Axe_\ga$. Recall that the period of $\ga$ for $F$ is
$\operatorname{Per}_F(\ga)=\int_0^{\ell(\ga)}\wt F(\phi_tv)\;dt$ with
$v$ as above. Note that $\alpha_*\L_\ga= \L_{\alpha\ga\alpha^{-1}}$
for every $\alpha\in\Ga$.  The set ${\H} _{\Ga,\,t}$ of loxodromic
elements of $\Ga$ whose translation length is at most $t$, and its
subset ${\H}' _{\Ga,\,t}$ of primitive elements, are invariant under
conjugation by elements of $\Ga$.

\medskip (1) Let us prove the first statement of Theorem
\ref{theo:equidisperorb}. Every primitive periodic orbit $g$ of the
geodesic flow in $T^1M$ is the image by $T^1\wt M$ of the lift by
$\pi:T^1\wt M\ra \wt M$ of the (oriented) translation axis of a
primitive loxodromic element of $\Ga$, whose translation length is the
length of $g$, and the number of such loxodromic elements with given
translation axis is equal to the cardinality of the stabiliser of any
tangent vector to this translation axis. By the definition (and
continuity) of the induced measure in $T^1M$ of a $\Ga$-invariant
measure in $T^1\wt M$ by the branched cover $T^1\wt M\ra T^1M$, we
only have to show that, as $t$ goes to $+\infty$, the measures
$$
\nu'_{t}=\delta\;\;e^{-\delta\; t}\;\sum_{\ga\in{\H}'_{\Ga,\,t}}
\;e^{\operatorname{Per}_F(\ga)}\;\L_{\ga}
$$
converge to $\frac{\wt m_F}{\|m_F\|}$ with respect to the weak-star
convergence of measures on $T^1\wt M$.

\medskip
\noindent {\bf Step 1. } Let us first prove that, as $t$ goes to
$+\infty$, the measures
$$
\nu''_{t}=\delta\;\;e^{-\delta\; t}\;\sum_{\ga\in{\H}_{\Ga,\,t}}
\;e^{\operatorname{Per}_F(\ga)}\;\L_{\ga}
$$
converge to $\frac{\wt m_F}{\|m_F\|}$ with respect to the weak-star
convergence of measures on $T^1\wt M$.

By the definition of the weak-star convergence, for every fixed
compact subset $K$ of $\wt M$, we only need to prove this convergence
when the measures are restricted to the compact subset $\pi^{-1}(K)$
of $T^1\wt M$. Let $\epsilon\in\mathopen{]}0,
\frac{1}{2}\mathclose{]}$.  As in the proof of Proposition
\ref{prop:stepone}, we will denote by $\epsilon_1,\epsilon_2,
\epsilon_3$ positive functions of $\epsilon$ (depending only on
$\delta>0$, on $\max_{\pi^{-1}(\N_1(K))} |\wt F|$, and on the
constants $c_1,c_2, c_3,c_4>0$ appearing in Lemma
\ref{lem:holderconseq} for $\wt F$ as well as for $\wt F\circ\iota$),
that converge to $0$ as $\epsilon$ goes to $0$.

For every $x\in K$, let $V(x,\epsilon)$ be the set of pairs
$(\xi,\eta)$ of distinct elements in $\wt M\cup\partial_\infty \wt M$
such that the geodesic segment, ray or line between $\xi$ and $\eta$
meets $B(x,\epsilon)$.  As in Equation \eqref{eq:controlFgap}, by
Equation \eqref{eq:ecartvisuel} and Lemma \ref{lem:holderconseq} (1),
for every $(\xi,\eta)\in V(x,\epsilon)\cap \partial_\infty^2 \wt M$,
we have
\begin{equation}\label{eq:controlFgapper}
e^{-\epsilon_1}\leq D_{F-\delta,\,x}(\xi,\eta)\leq e^{\epsilon_1}\;.
\end{equation}
For all $x\in K$ and $t>0$, let $\nu$ and $\nu_t$ be the measures on
$(\wt M\cup\partial_\infty \wt M)^2$ defined by:
$$
d\nu(\xi,\eta)=
\frac{d\mu^\iota_{x}(\xi)\otimes d\mu_{x}(\eta)}
{D_{F-\delta,\,x}(\xi,\eta)^2}\;,\;\;\;\;\nu_t=
\delta\;\|m_F\|\; e^{-\delta t}\sum_{\ga\in\Ga,\;d(x,\,\ga x)\leq t}
e^{\int_x^{\ga x}\wt F}\;\D_{\ga^{-1}x}\otimes\D_{\ga x} \;.
$$
By Theorem \ref{maingrowthRoblin} (taking $y=x$), we know that $\nu_t$
weak-star converges to $\mu^\iota_{x}\otimes \mu_{x}$ as
$t\ra+\infty$.  Let $\psi$ be a continuous nonnegative map with
support in $V(x,\epsilon)$. By Equation \eqref{eq:controlFgapper}, we
hence have
$$
e^{-2\,\epsilon_1}\;\nu(\psi)\leq \liminf_{t\ra+\infty}\nu_t(\psi)\leq
\limsup_{t\ra+\infty}\nu_t(\psi)\leq e^{2\,\epsilon_1}\;\nu(\psi)\;.
$$
Let 
$$
\nu'''_{t}= \delta\;\|m_F\|\; e^{-\delta t}\sum_{\ga\in{\H}_{\Ga,\,t}} 
e^{{\rm Per}_F(\ga)}\;\;\D_{\ga_-}\otimes\D_{\ga_+}\;.
$$
Note that by Gromov's hyperbolicity criterion (see also Lemma
\ref{lem:crithyperboelem}), for every $\ga\in\Ga$, if $(\ga^{-1}x,\ga
x)\in V(x,\epsilon)$ and $d(x,\ga x)$ is large enough, then $\ga$ is
loxodromic and $x$ is at distance at most $2 \epsilon$ from the
translation axis $\Axe_\ga$ of $\ga$. In particular
$$
\ell(\ga)\leq d(x,\ga x)\leq \ell(\ga)+4\epsilon\;.
$$
Furthermore, $\ga^{\pm 1} x$ is close to $\ga_\pm$ uniformly in $\ga$,
hence so is $\D_{\ga^{\pm 1}x}$ to $\D_{\ga_\pm}$. By Lemma
\ref{lem:technicholder} and the remark following it, if $p$ is the
closest point on $\Axe_\ga$ to $x$, we hence have
\begin{align*}
\Big|\int_x^{\ga x}\wt F-{\rm Per}_F(\ga)\Big|&=
\Big|\int_x^{\ga x}\wt F-\int_p^{\ga p}\wt F\,\Big|
\\ & \leq 2\,c_3(2\epsilon)^{c_4}+
2\,\epsilon\,\Big(\max_{\pi^{-1}(B(x,\,2\epsilon))}|\wt F|+
\max_{\pi^{-1}(B(\ga x,\,2\epsilon))}|\wt F|\Big)\;.
\end{align*}
The right hand side of this inequality, which we will denote by
$\eta(2\epsilon)$ for future use, tends to $0$ as $\epsilon$ goes to
$0$ uniformly in $\ga$ by the $\Ga$-invariance of $\wt F$.  Hence if
$t$ is large enough, then
$$
e^{-\epsilon_2}
\nu_t(\psi) \leq \nu'''_{t}(\psi)\leq e^{\epsilon_2} \nu_{t+4\epsilon}(\psi)\;.
$$

Identifying $T^1\wt M$ with $\partial_\infty^2 \wt M\times\RR$ by the
Hopf parametrisation with respect to the base point $x$, note that
$$
\frac{1}{\|m_F\|}\;\nu'''_{t}\otimes ds=\nu''_{t}\,.
$$
Therefore, for every continuous map $\psi':T^1\wt M\ra\RR$ with
compact support in $V(x,\epsilon)\times\RR$, which is a product of two
continuous maps on each of the two variables of this product, we have
$$
e^{-\epsilon_3}\;\frac{m_F(\psi')}{\|m_F\|}\leq 
\liminf_{t\ra+\infty}\nu''_{t}(\psi')\leq
\limsup_{t\ra+\infty}\nu''_{t}(\psi')\leq 
e^{\epsilon_3}\;\frac{m_F(\psi')}{\|m_F\|}\;.
$$
Approximating uniformly continuous maps with compact support in
$V(x,\epsilon)\times\RR$ by finite linear combinations of product
maps, covering $\pi^{-1}(K)$ by finitely many sets
$V(x,\epsilon)\times\RR$, using a partition of unity and letting
$\epsilon$ goes to $0$, the first step  follows.

\medskip
\noindent {\bf Step 2. } Let us now prove that $\nu''_t-\nu'_t$
weak-star converges to $0$ as $t\ra+\infty$, which implies the first
assertion of Theorem \ref{theo:equidisperorb}. 

We start by proving the following lemma, provided by the referee,
which says that the contribution of the non-primitive periodic orbits
is negligible.

Let $x\in\wt M$, $\epsilon\in\mathopen{]}0,\frac{1}{2} \mathclose{]}$
and $U_-,U_+$ be small enough neighbourhoods in $\partial_\infty \wt
M$ of two distinct points in $\Lambda\Ga$ (with boundary of measure
$0$ for the Patterson measures $\mu^\iota_x, \mu_x$ respectively) so
that the geodesic line between any point in $U_-$ and any point in
$U_+$ passes at distance at most $\epsilon$ from $x$. In particular,
for every loxodromic element $\ga\in\Ga$ with $\ga_\pm\in U_\pm$, its
translation axis passes at distance at most $\epsilon$ from $x$, so
that, as seen above,
\begin{equation}\label{eq:controlengthpourcomptperiod}
\ell(\ga)\leq d(x,\ga x)\leq \ell(\ga)+2\epsilon\,,
\end{equation}
and, with the above notation $\eta(\cdot)$,
\begin{equation}\label{eq:controlperiodpourcomptperiod}
\Big|\int_x^{\ga x}\wt F-{\rm Per}_F(\ga)\Big|\leq \eta(\epsilon)\,.
\end{equation}
For every $k\in\NN-\{0\}$, let $H_k$ be the set of loxodromic elements
$\ga\in\Ga$ such that $\ga_{\pm}\in U_\pm$ and there exists $\alpha
\in \Ga$ with $\ga=\alpha^k$. Note that $H_k=\{\alpha^k\;:\; \alpha\in
H_1\}$ and $H_k\subset H_1$.

\blemm\label{lem:neglinonprimreferee} For every $\epsilon'>0$, there
exists $C>0$ and $\rho\in\mathopen{]}0,1\mathclose{[}$ such that, for
all $n,k\in\NN$ with $k\geq 2$, we have
$$
\sum_{\ga\in H_k,\;n-1< d(x,\,\ga x)\leq n} e^{\int_x^{\ga x}\wt F}\leq 
C\;\rho^n+\frac{\epsilon'}{2^k}\,.
$$
\elemm

\dem In order to obtain the two terms in the right hand side of this
inequality, we will separate the sum on its left hand side into the
elements of $H_k$ which are $k$-th powers of elements with not too
large translation length, and the others.

Let us fix $\epsilon''>0$. For every $\alpha\in H_1$, we have 
$$
e^{{\rm  Per}_F (\alpha)} < e^{\delta\,\ell(\alpha)}\,,
$$ 
otherwise for every $n\in\NN$, we have $e^{{\rm Per}_F(\alpha^n)}
\geq e^{\delta\,\ell(\alpha^n)}$, hence $e^{\int_x^{\alpha^n x} \wt F}
\geq e^{-\eta(\epsilon)-2\epsilon} \,e^{\delta\,d(x,_,\alpha^n x)}$ by
Equations \eqref{eq:controlengthpourcomptperiod} and
\eqref{eq:controlperiodpourcomptperiod}, a contradiction to
Proposition \ref{prop:atomlessreferee} (2). By discreteness, there
are only finitely many loxodromic elements in $\Ga$, whose translation
axis passes at distance at most $\epsilon$ from $x$, with translation
length as most some constant. Hence, again by Proposition
\ref{prop:atomlessreferee} (2), there exists $t_0=t_0(\epsilon'')>0$
such that if $\alpha\in H_1$ satisfies $\ell(\alpha)\geq t_0$, then
$$
e^{\int_x^{\alpha x}\wt F} \leq \epsilon''\, e^{\delta\,d(x,\,\alpha x)}\,.
$$

Since $\sum_{\alpha'\in H_1,\;\ell(\alpha')\leq t}e^{{\rm Per}_F(\alpha')} 
= \nu''_t(U_-\times U_+)$ and as seen above, there exists $t_1,c>0$
such that, for all $t\geq t_1$,
\begin{equation}\label{eq:minopournegnoprim}
\sum_{\alpha'\in H_1,\;t-1<\ell(\alpha')\leq t} e^{\int_x^{\alpha' x}\wt F} 
\geq e^{-\eta(\epsilon)} \sum_{\alpha'\in H_1,\;t-1<\ell(\alpha')\leq t} 
e^{{\rm Per}_F (\alpha')} \geq\frac{1}{c}\;e^{\delta\;t}\,.
\end{equation}

Let $t_2=\max\{t_0,t_1\}$. We define $\rho=\max_{\alpha\in
  H_1,\;\ell(\alpha)\leq t_2} e^{\frac{{\rm Per}_F
    (\alpha)}{\ell(\alpha)}-\delta} \in\mathopen{]}0,1\mathclose{[}$ and $C=
e^{\eta(\epsilon)}\;\rho^{-1-2\,\epsilon}\,\card\{\alpha\in H_1\;:\;
\ell(\alpha) \leq t_2\}$ (which both depend on $\epsilon''$).

Now, as indicated at the beginning, we write, for every $n,k\in\NN$
with $k\geq 2$, using Equations \eqref{eq:controlengthpourcomptperiod} and
\eqref{eq:controlperiodpourcomptperiod},
\begin{align*}
&\sum_{\substack{\ga\in H_k\\n-1< d(x,\,\ga x)\leq n}} e^{\int_x^{\ga x}\wt F} =
\sum_{\substack{\alpha\in H_1\\n-1< \,d(x,\,\alpha^k x)\leq n}} 
e^{\int_x^{\alpha^k x} \wt F}\leq 
\sum_{\substack{\alpha\in H_1\\n-1-2\epsilon\,< \ell(\alpha^k)\leq n}} 
e^{\eta(\epsilon)} e^{{\rm Per}_F(\alpha^k)} \\ \leq\; & 
\sum_{\substack{\alpha\in H_1,\;\ell(\alpha)\leq t_2\\n-1-2\epsilon\,< k\,\ell(\alpha)\leq n}} 
e^{\eta(\epsilon)}e^{k\,{\rm Per}_F(\alpha)}+
\sum_{\substack{\alpha\in H_1,\;\ell(\alpha)\geq t_2\\\frac{n}{k}-1< \ell(\alpha)\leq \frac{n}{k}}} 
e^{\eta(\epsilon)}e^{k\,{\rm Per}_F(\alpha)}\;.
\end{align*}
Let us denote by $I$ and $J$ the first and second sums in the last
line above. We have, by the definition of $\rho$ and $C$,
$$
I\leq e^{\eta(\epsilon)}
\sum_{\substack{\alpha\in H_1,\;\ell(\alpha)\leq t_2\\
n-1-2\epsilon\,< k\,\ell(\alpha)\leq n}} 
\rho^{k\,\ell(\alpha)}e^{k\,\delta\;\ell(\alpha)}
\leq e^{\eta(\epsilon)}
\sum_{\alpha\in H_1,\;\ell(\alpha)\leq t_2} 
\rho^{n-1-2\epsilon}e^{\delta\;n}\leq C\;\rho^n\;e^{\delta\;n}\,.
$$
By Equation \eqref{eq:controlperiodpourcomptperiod}, for every $\alpha
\in H_1$ such that $t_0\leq \ell(\alpha)\leq \frac{n}{k}$, we have
$$
e^{\int_x^{\alpha x} \wt F}\leq \epsilon''\, e^{\delta\,d(x,\,\alpha x)}\leq
\epsilon''\, e^{\delta\,(\ell(\alpha)+2\,\epsilon)}\leq
\epsilon''\, e^{\delta\,(\frac{n}{k}+2\,\epsilon)}\leq
c\,\epsilon''\, e^{2\,\epsilon\,\delta}
\sum_{\alpha'\in H_1,\;\frac{n}{k}-1<\ell(\alpha')\leq \frac{n}{k}} 
e^{\int_x^{\alpha' x}\wt F}\,.
$$
By Equation \eqref{eq:minopournegnoprim}, we hence have
\begin{align*}
J &\leq e^{(1+k)\eta(\epsilon)}
\sum_{\substack{\alpha\in H_1,\;\ell(\alpha)\geq t_0\\\frac{n}{k}-1< \ell(\alpha)\leq \frac{n}{k}}} 
e^{k\,\int_x^{\alpha x} \wt F}=e^{(1+k)\eta(\epsilon)}
\sum_{\substack{\alpha\in H_1,\;\ell(\alpha)\geq t_0\\\frac{n}{k}-1< \ell(\alpha)\leq \frac{n}{k}}} 
e^{\int_x^{\alpha x} \wt F}\big(e^{\int_x^{\alpha x} \wt F}\big)^{k-1}\\ & \leq
e^{(1+k)\eta(\epsilon)}\,(c\,\epsilon''\, e^{2\,\epsilon\,\delta})^{k-1}
\sum_{\substack{\alpha\in H_1,\;\ell(\alpha)\geq t_0\\\frac{n}{k}-1< \ell(\alpha)\leq \frac{n}{k}}}
e^{\int_x^{\alpha x} \wt F}
\Big(\sum_{\alpha'\in H_1,\;\frac{n}{k}-1< \ell(\alpha')\leq \frac{n}{k}} 
e^{\int_x^{\alpha' x} \wt F}\Big)^{k-1}\;.
\end{align*}
By the triangle inequality, for all $\alpha_1,\dots,
\alpha_k\in H_1$ with $\ell(\alpha_i)\leq \frac{n}{k}$, we have
$$
d(x,\alpha_1\dots\alpha_k x)\leq \sum_{i=1}^kd(x,\alpha_i x)\leq 
\sum_{i=1}^k\ell(\alpha_i) +2\,k\,\epsilon\leq n+2\,k\,\epsilon\,.
$$
By a quasi-geodesic argument as in Proposition
\ref{prop:schottkysemigroups} and by Lemma \ref{lem:technicholder},
there exists $c'>0$ such that if $\epsilon, U_\pm$ are small enough
and $t$ big enough, for every $k\geq 2$, the product elements
$\alpha_1\dots\alpha_k$ in $\Ga$, where $\alpha_1,\dots, \alpha_k\in
H_1$ with $t-1<d(x,\alpha_ix)\leq t+2\epsilon$ and
$d(\alpha_ix,\alpha_jx)\geq 6$ if $\alpha_i\neq \alpha_j$ are pairwise
distinct and
$$
\Big|\int_x^{\alpha_1\dots\alpha_k x} \wt F-
\sum_{i=1}^k\int_x^{\alpha_i x} \wt F\;\Big|\leq k\;c'\,.
$$
Hence, up to increasing $t_0$, by developing the above $(k-1)$-th
power and by a definite proportion argument to make the relevant orbit
points $6$-separated, there exists $c''>0$ (independent of
$\epsilon''$) such that
$$
J\leq (c''\epsilon'')^k\sum_{\ga\in\Ga,\;d(x,\,\ga x) \leq n+2\,k\,\epsilon} 
e^{\int_x^{\ga x} \wt F}\,.
$$ 
By Corollary \ref{coro:bigOloggrowth}, there exists $c'''>0$ such that
this sum is at most $c'''e^{\delta\,(n+2\,k\,\epsilon)}$. Hence, if
$\epsilon''>0$ is small enough, we have $J\leq\frac{ \epsilon'}{2^k}$,
as required.  \cqfd

\medskip Let us now conclude the proof of Step 2. Since the
translation lengths of the loxodromic elements of $\Ga$ whose
translation axis meets a given compact subset of $\wt M$ have a
positive lower bound by discreteness, there exists $a>0$ such that for
all $k\in\NN-\{0\}$ and $\ga\in H_k$ with $d(x,\ga x)\leq n$, we have
$k\leq a\,n$. By Equations \eqref{eq:controlengthpourcomptperiod} and
\eqref{eq:controlperiodpourcomptperiod}, and by a finite summation, we
hence have that
$$
\sum_{\ga\in\bigcup_{k\geq 2}H_k,\;n-1< \ell(\ga)\leq n} e^{{\rm Per}_F(\ga)}
$$
tends to $0$ as $n\ra+\infty$. By a standard covering argument by
small open sets (of the form $U_-\times U_+\times\mathopen{]}-\eta,
+\eta\mathclose{[}$ in Hopf's coordinates, with $U_\pm$ as above and
$\eta>0$ small enough) of the support of a continuous map $\psi'$ with
compact support on $T^1\wt M$, and since $|\L_\ga(\psi')|$ is
uniformly bounded (by the product of the maximum of $|\psi'|$ and the
diameter of its support) for every loxodromic element $\ga\in\Ga$,
this concludes the proof of Step 2.

\medskip (2) The deduction of the second assertion from the first one
is standard. Consider the measures 
$$
m'_t=\delta e^{-\delta t}
\sum_{g\in\Per'(t)} \;e^{\L_g(F)} \;\L_{g}\;\;\;{\rm and}\;\;\; 
m''_t=\delta t e^{-\delta t} \sum_{g\in\Per'(t)} \;e^{\L_g(F)}\;
\frac{\L_{g}}{\ell(g)}
$$
on $T^1M$.  Fix a continuous map $\psi:T^1M\ra
\mathopen{[}0,+\infty\mathclose{[}$ with compact support. For every
$\epsilon>0$, for every $t>0$, we have
\begin{align*}
  m''_t(\psi) & \geq m'_t(\psi)\geq \delta e^{-\delta t}
  \sum_{g\in\Per'(t)-\Per'(e^{-\epsilon}t)}
  \;e^{\L_g(F)} \;\L_{g}(\psi)\\
  & \geq e^{-\epsilon}\delta t e^{-\delta t}
  \sum_{g\in\Per'(t)-\Per'(e^{-\epsilon}t)} \;e^{\L_g(F)}
  \;\frac{\L_{g}(\psi)}{\ell(g)}\\
&= e^{-\epsilon} m''_t(\psi) - e^{-\epsilon}\delta t \;e^{-\delta t}
  \sum_{g\in\Per'(e^{-\epsilon}t)} \;e^{\L_g(F)}
  \;\frac{\L_{g}(\psi)}{\ell(g)}\;.
\end{align*}
Since the closed orbits meeting the support of $\psi$ have a positive
lower bound on their lengths, and by the first assertion of Theorem
\ref{theo:equidisperorb}, there exists a constant $c>0$ such that the
second term of the above difference is at most $c\,e^{-\delta t}
e^{\delta \,e^{-\epsilon}t}$, which tends to $0$ as $t$ tends to infinity.
Hence by applying twice the first assertion of Theorem
\ref{theo:equidisperorb}, we have
$$
\frac{m_F(\psi)}{\|m_F\|}= \lim_{t\ra+\infty}m'_t(\psi)\leq
\liminf_{t\ra+\infty}m''_t(\psi) \leq \limsup_{t\ra+\infty}
m''_t(\psi)\leq \lim_{t\ra+\infty} e^{\epsilon} m'_t(\psi)=
e^{\epsilon}\frac{m_F(\psi)}{\|m_F\|}\;,
$$ 
and the result follows by letting $\epsilon$ go to $0$ (and writing
any continuous map $\psi:T^1M\ra \RR$ with compact support into the
sum of its positive and negative parts).  \cqfd

\brema\label{rem:Gureprimit}
{\rm Step 2 of the above proof shows that the Gurevich pressure,
  if it is finite and if $m_F$ is finite and mixing, may be defined by
  considering only primitive periodic orbits: if $W$ is any relatively
  compact open subset of $T^1M$ meeting $\Omega\Ga$, for every $c>0$
  large enough, we have
$$
P_{Gur}(\Ga,F)=\lim_{t\ra+\infty}\;\frac{1}{t}\ln
\sum_{g\in \Per'(t)-\Per'(t-c),\;g\cap W\neq \emptyset}\;e^{\int_g F}\,,
$$
and if $P_{Gur}(\Ga,F)>0$, then
$$
P_{Gur}(\Ga,F)=\lim_{t\ra+\infty}\;\frac{1}{t}\ln
\sum_{g\in \Per'(t),\;g\cap W\neq \emptyset}\;e^{\int_g F}\,.
$$
}
\erema

\medskip In the same way Theorem \ref{maingrowthRoblinbis} has been
deduced from Theorem \ref{maingrowthRoblin}, the first claim of the
following result, which no longer assumes $\delta_{\Ga,\,F}>0$, can be
deduced from Theorem \ref{theo:equidisperorb}. The second claim
follows immediately.

\btheo\label{theo:equidisperorbbis}%
\index{theorem@Theorem!Equidistribution of periodic \\orbits}
Assume that $\delta_{\Ga,\,F}<+\infty$, and that the Gibbs measure
$m_F$ is finite and mixing under the geodesic flow on $T^1M$. Let $U$
and $V$ be two open subsets of $\partial_\infty\wt M$ meeting
$\Lambda\Ga$ such that $\mu_x(\partial U)=\mu^\iota_x(\partial
V)=0$. Then for every $c>0$, as $t$ goes to $+\infty$, 
$$
\frac{\delta_{\Ga,\,F}}{1-e^{-c\,\delta_{\Ga,\,F}}}
\,\;\;e^{-\delta_{\Ga,\,F}\; t}\;
\sum_{g\in\Per'(t)-\Per'(t-c)} \;e^{\L_g(F)}\;\L_{g}
\;\;\weakstar\;\;\frac{m_F}{\|m_F\|}\,,
$$
$$
\frac{\delta_{\Ga,\,F}}{1-e^{-c\,\delta_{\Ga,\,F}}}
\; t\;\;e^{-\delta_{\Ga,\,F}\; t}\;
\sum_{g\in\Per'(t)-\Per'(t-c)} 
\;e^{\L_g(F)}\;\frac{\L_{g}}{\ell(g)}
\;\;\weakstar\;\;\frac{m_F}{\|m_F\|}\,. 
$$
\etheo

\medskip Applying the equidistribution result \ref{theo:equidisperorb}
to characteristic functions, for every relatively compact open subset
$U$ of $T^1M$ whose boundary has Gibbs measure $0$, we have that, as
$t\ra+\infty$,
$$
\sum_{g\in\Per'(t)}  \;e^{{\rm Per}_F(g)}\frac{\ell(g\cap U)}{\ell(g)}
\;\sim \;
\frac{e^{\delta_{\Ga,\,F}t}}{\delta_{\Ga,\,F}t}\frac{m_F(U)}{\|m_F\|}\;.
$$
As there exists a compact subset of $T^1M$ containing all closed
orbits of the geodesic flow if $\Ga$ is convex-cocompact, the
following corollary holds.

\bcoro \label{coro:comptageperiode} 
Assume that $\Ga$ is convex-cocompact, that $(\phi_t)_{t\in\RR}$ is
topologically mixing and that $\delta_{\Ga,\,F}<+\infty$. Then for
every $c>0$, as $t\ra+\infty$,
$$
\sum_{g\in\Per'(t)-\Per'(t-c)}  \;e^{{\rm Per}_F(g)}
\sim\frac{1-e^{c\;\delta_{\Ga,\,F}}}{\delta_{\Ga,\,F}}\;
\frac{e^{\delta_{\Ga,\,F}\;t}}{t}\;.
$$
If furthermore $\delta_{\Ga,\,F}>0$, then as
$t\ra+\infty$,
$$
\sum_{g\in\Per'(t)}  \;e^{{\rm Per}_F(g)}
\sim\frac{e^{\delta_{\Ga,\,F}\;t}}{\delta_{\Ga,\,F}\;t}\;.
$$
\ecoro

Let $\C_b(T^1M)^*$ be the dual topological space of the Banach space
of real bounded continuous functions on $T^1M$. By the {\it narrow
  convergence}\index{convergence!narrow}%
\index{narrow convergence} of finite measures, we mean as usual the
convergence for the weak-star topology in $\C_b(T^1M)^*$. Corollary
\ref{coro:comptageperiode} can be improved: it holds when
``convex-cocompact'' is replaced by ``geometrically finite with finite
Gibbs measure'', though the extension requires more work. For this, we
improve, under the geometrically finiteness condition, Theorem
\ref{theo:equidisperorb} from weak-star convergence to narrow
convergence, as in \cite{Roblin03} for the case $F=0$. The main point
is to prove that there is no loss of mass in the cuspidal parts during
the convergence process.

\btheo\label{theo:equidisperorbetroite} Assume that the critical
exponent $\delta_{\Ga,\,F}$ of $(\Ga,F)$ is finite and positive,
that the Gibbs measure $m_F$ is finite and mixing under the geodesic
flow on $T^1M$, and that $\Ga$ is geometrically finite.

(1) As $t$ goes to $+\infty$, the measures
$$
\delta_{\Ga,\,F}\;\;e^{-\delta_{\Ga,\,F}\; t}\;\sum_{g\in\Per'(t)} 
\;e^{\L_g(F)}\;\L_{g}
$$
converge to the normalised Gibbs measure $\frac{m_F}{\|m_F\|}$ with
respect to the narrow convergence.

(2) As $t$ goes to $+\infty$,
$$
\delta_{\Ga,\,F}\; t\;\;e^{-\delta_{\Ga,\,F}\; t}\;\sum_{g\in\Per'(t)} 
\;e^{\L_g(F)}\;\frac{\L_{g}}{\ell(g)}
$$
converges to $\frac{m_F}{\|m_F\|}$ with respect to the narrow
convergence.
\etheo

\dem The deduction of Assertion (2) from Assertion (1) proceeds as in
the proof of Theorem \ref{theo:equidisperorb}, using the fact that
since $\Ga$ is geometrically finite, there exists a compact set of $M$
meeting (though not necessarily containing) all closed geodesics,
hence there exists $c'>0$ such that every closed geodesic of $M$ has
length at least $c'$.

Let us prove Assertion (1). Let $\delta= \delta_{\Ga,\,F}$.  For
$t\geq 0$, consider the measure
$$
m_t=\delta\,
e^{-\delta t} \sum_{g\in\Per'(t)} \;e^{\L_g(F)} \;\L_{g}
$$ 
on $T^1M$ and let $\pi_*m_t$ be its push-forward on $M$. We will use
the notation of the proof of Theorem \ref{theo:DOPB}, in particular
$\Par_\Ga,\H_p,\Ga_p,\F_p, x_0,\kappa$ for every parabolic fixed point
$p$ of $\Ga$.

For every $r\geq 0$ and for every parabolic fixed point $p$ of $\Ga$,
let $\H_p(r)$ be the horoball centred at $p$, contained in $\H_p$
such that the distance between the horospheres $\partial\H_p(r)$ and
$\partial\H_p$ is $r$.

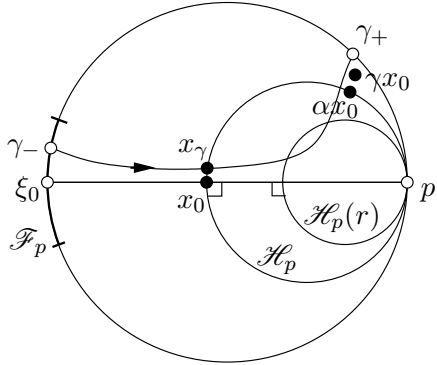
\begin{center}
\input{fig_paraborne3.pstex_t}
\end{center}

By the finiteness of the number of orbits of parabolic fixed points
under $\Ga$ and by the compactness of the quotient $\Ga\bs
\big(\C\Lambda\Ga -\bigcup_{p\in\Par_\Ga} \H_p(r)\big)$, we only have to
prove that for every parabolic fixed point $p$ of $\Ga$,
$$
\lim_{r\ra +\infty}\;\limsup_{t\ra+\infty} \pi_*m_t(\Ga_p\bs\H_p(r))=0\;.
$$
Fix such a point $p$. A closed geodesic in $M$ meeting $\Ga_p\bs\H_p$
has a unique lift in $\wt M$ meeting $\H_p$ starting from the
fundamental domain $\F_p$ for the action of $\Ga_p$ on
$\Lambda\Ga-\{p\}$; its other endpoint is in $\alpha\F_p$ for some
$\alpha\in\Ga_p$. For all $t\geq 0$ and $\alpha\in\Ga_p$, let
$\Ga(t,\alpha)$ be the set of loxodromic elements $\ga$ of $\Ga$ with
repulsive fixed point $\ga_-$ in $\F_p$, attractive fixed point
$\ga_+$ in $\alpha\F_p$ and translation length $\ell(\ga)$ at most
$t$.  Then (the inequality being an equality if the elements of
$\Ga(t,\alpha)$ were being assumed to be primitive)
\begin{equation}\label{eq:mesurcuspr}
\pi_*m_t(\Ga_p\bs\H_p(r))\leq \delta e^{-\delta t} \sum_{\alpha\in \Ga_p}
\sum_{\ga\in \Ga(t,\,\alpha)} \;e^{\operatorname{Per}_F(\ga)}
\;\Long(\Axe_\ga\cap \H_p(r))\;.
\end{equation}

Let $\ga\in \Ga(t,\alpha)$ be such that its translation axis $\Axe_\ga$
meets $\H_p(r)$ (note that for the others, we have
$\Long(\Axe_\ga\cap \H_p(r))=0$). We orient $\Axe_\ga$ from $\ga_-$ to
$\ga_+$, and we denote by $x_\ga$ the entrance point of $\Axe_\ga$ in
$\H_p$.  By Equation \eqref{eq:proxixx}, we have
$$
d(x_0,x_\ga)\leq \kappa\;.
$$
Since the distance between any point of $\partial \H_p$ and any point
of $\H_p(r)$ is at least $r$, we have (see the above picture),
\begin{equation}\label{eq:majolength}
\Long(\Axe_\ga\cap \H_p(r))\leq d(x_0,\alpha x_0)-2r +2\kappa\;.
\end{equation}
In particular,  
\begin{equation}\label{eq:bensifinalement}
d(x_0,\alpha x_0)\geq 2r -2\kappa\;.
\end{equation}

Consider the constants $R$ and $C$ given by Mohsen's shadow lemma
\ref{lem:shadowlemma} applied to $(\mu_x)_{x\in\wt M}$ as above and
$K=\{x_0\}$, and let $\OOO_\ga= \OOO_{x_0}B(\ga x_0,R)$. The point
$\alpha x_0$ is at distance uniformly bounded from $\Axe_\ga$. Since
$\H_p$ is precisely invariant (see the beginning of the proof of
Theorem \ref{theo:DOPB}) and $x_0\in\partial\H_p$, the point $\ga x_0$
does not belong to the interior of $\H_p$.  Hence $\alpha x_0$ is at
distance uniformly bounded from $\mathopen{[}x_0,\ga
x_0\mathclose{]}$. Therefore $x_0=\alpha^{-1}\alpha x_0$ is at
distance uniformly bounded from
$\mathopen{[}\alpha^{-1}x_0,\alpha^{-1}\ga x_0\mathclose{]}$. If $r$
is large enough, then $d(x_0,\alpha^{-1}x_0)$ is large, hence
$\alpha^{-1}x_0$ is close to $p$. Therefore, there exists a compact
subset $K$ (independent of $t,r,\alpha,\ga$) of $\partial_\infty \wt
M-\{p\}$ such that the set $\alpha^{-1}\OOO_\ga= \OOO_{\alpha^{-1}x_0}
B(\alpha^{-1}\ga x_0,R)$ is contained in $K$ for $r$ large enough.

Since $d(x_0,x_\ga)\leq \kappa$ and by (two applications of) Lemma
\ref{lem:technicholder}, there exists $c_1\geq 0$ such that
$$
\Big|\operatorname{Per}_F(\ga)-\int_{x_0}^{\ga x_0}\wt F\;\Big|
=\Big|\int_{x_\ga}^{\ga x_\ga}\wt F\;-\int_{x_0}^{\ga x_0}\wt F\;\Big|
\leq c_1\;.
$$
Since $x_0$ is at distance at most $\kappa$ from $\Axe_\ga$, we have
$$
d(x_0,\ga x_0)\leq \ell(\ga)+2\kappa\leq t+2\kappa\;,
$$
and if furthermore $\ga\notin\Ga(t-1,\alpha)$, then
$$
t-1\leq \ell(\ga)\leq d(x_0,\ga x_0)\;.
$$ 
Hence, by Mohsen's shadow lemma \ref{lem:shadowlemma},
$$
e^{\operatorname{Per}_F(\ga)}\leq 
e^{c_1}e^{\int_{x_0}^{\ga x_0}\wt F}\leq 
e^{c_1+2\delta \kappa}e^{\delta t}
e^{\int_{x_0}^{\ga x_0}(\wt F-\delta)}
\leq C\,e^{c_1+2\delta \kappa}e^{\delta t}\mu_{x_0}(\OOO_\ga)\;.
$$
By the discreteness of $\Ga x_0$, there exists $N>0$ (independent of
$t$) such that a point $\xi\in\partial_\infty \wt M$ belongs to at most
$N$ shadows $\OOO_\ga$ for $\ga\in \Ga$ with $t-1\leq d(x_0,\ga x_0)
\leq t+2\kappa$. Since the shadow $\OOO_\ga$ is contained in $\alpha
K$ if $\Axe_\ga$ meets $\H_p(r)$ for $r$ large enough, we hence have, for
$r$ large enough,
$$
\sum_{\ga\in \Ga(t,\,\alpha)-\Ga(t-1,\,\alpha)\,:\; \Axe_\ga\cap\H_p(r)\neq
  \emptyset}\;e^{\operatorname{Per}_F(\ga)}\leq 
CN\,e^{c_1+2\delta\kappa}\,e^{\delta t}\mu_{x_0}(\alpha K)\;.
$$
By a summation, there exists therefore a constant $c_2>0$ such that
\begin{equation}\label{eq:soixantedix}
\sum_{\ga\in \Ga(t,\,\alpha)\,:\; \Axe_\ga\cap\H_p(r)\neq
  \emptyset} \;e^{\operatorname{Per}_F(\ga)}\leq 
c_2\,e^{\delta t}\mu_{x_0}(\alpha K)\;.
\end{equation}
By the equations \eqref{eq:equivdensity}, 
\eqref{eq:radonykodensity} and \eqref{eq:reserviraDOPMoh}, we have,
for some constant $c_3>0$,
\begin{equation}\label{eq:ombrebis}
\mu_{x_0}(\alpha K)=\mu_{\alpha^{-1} x_0}(K)=
\int_{\xi\in K} e^{-C_{F-\delta,\,\xi}(\alpha^{-1} x_0,\,x_0)}\;d\mu_{x_0}(\xi)
\leq c_3\mu_{x_0}(K)\;e^{\int_{x_0}^{\alpha x_0} (\wt F-\delta)}\;.
\end{equation}
By the equations \eqref{eq:mesurcuspr}, \eqref{eq:majolength},
\eqref{eq:bensifinalement}, \eqref{eq:soixantedix} and
\eqref{eq:ombrebis}, for $r$ large enough and for every $t\geq 0$, we
hence have
\begin{align*}
&\pi_*m_t(\Ga_p\bs\H_p(r))\\ &
\leq \delta c_2 c_3 \mu_{x_0}(K)\sum_{\alpha\in
  \Ga_p\,:\;d(x_0,\,\alpha x_0)\geq 2r-2\kappa}
(d(x_0,\alpha x_0)-2r +2\kappa)\;e^{\int_{x_0}^{\alpha x_0} (\wt F-\delta)}\;.
\end{align*}
The result then follows from Theorem \ref{theo:DOPB}, since we assume
$m_F$ to be finite (which implies that $(\Ga,F)$ is of divergence
type by Corollary \ref{coro:uniqGibstate}).  
\cqfd

\medskip In a similar way to that where Theorem
\ref{maingrowthRoblinbis} has been deduced from Theorem
\ref{maingrowthRoblin}, the first claim of the following result, which
no longer assumes $\delta_{\Ga,\,F}>0$, can be deduced from Theorem
\ref{theo:equidisperorbetroite}. The second claim follows
immediately.

\medskip\btheo\label{theo:equidisperorbetroitebis} Assume that
$\delta_{\Ga,\,F}<+\infty$, that $m_F$ is finite and mixing under the
geodesic flow on $T^1M$, and that $\Ga$ is geometrically finite. Then
for every $c>0$, as $t\ra+\infty$, the measures
$$
\frac{\delta_{\Ga,\,F}}{1-e^{-c\,\delta_{\Ga,\,F}}}
\;\;e^{-\delta_{\Ga,\,F}\; t}\;\sum_{g\in\Per'(t)-\Per'(t-c)} 
\;e^{\L_g(F)}\;\L_{g}
$$
and
$$
\frac{\delta_{\Ga,\,F}}{1-e^{-c\,\delta_{\Ga,\,F}}}\; t\;\;
e^{-\delta_{\Ga,\,F}\; t}\;
\sum_{g\in\Per'(t)-\Per'(t-c)} \;e^{\L_g(F)}\;\frac{\L_{g}}{\ell(g)}
$$
converge to $\frac{m_F}{\|m_F\|}$ with respect to the narrow
convergence.
\etheo

\medskip The next result is immediate, applying the second assertions
of Theorems \ref{theo:equidisperorbetroite} and
\ref{theo:equidisperorbetroitebis} above to the (bounded) constant
function $1$. By Subsection \ref{subsec:finiteness}, the finiteness of
$m_F$ (which was automatic in the convex-cocompact case but no longer
is in the geometrically finite one) follows if
$\delta_{\Ga_p,\,F}<\delta_{\Ga,\,F}$ for every parabolic fixed point
$p$ of $\Ga$.

\bcoro \label{coro:comptageperiodegeomfini} 
\index{theorem@Theorem!Counting of periodic orbits}
Assume that $\Ga$ is geometrically finite, $(\phi_t)_{t\in\RR}$ is
topologically mixing, $\delta_{\Ga,\,F}<+\infty$ and $m_F$ is
finite. Then for every $c>0$, as $t\ra+\infty$, we have
$$
\sum_{g\in\Per'(t)-\Per'(t-c)}  \;e^{{\rm Per}_F(g)}
\;\sim\;\frac{1-e^{-c\,\delta_{\Ga,\,F}}}{\delta_{\Ga,\,F}}\;
\frac{e^{\delta_{\Ga,\,F}\;t}}{t}\;.
$$
If furthermore $\delta_{\Ga,\,F}>0$, then, as $t\ra+\infty$, we have
$$
\sum_{g\in\Per'(t)}  \;e^{{\rm Per}_F(g)}
\sim\frac{e^{\delta_{\Ga,\,F}\;t}}{\delta_{\Ga,\,F}\;t}\;.
$$
\ecoro

\subsection{The case of infinite Gibbs measure}
\label{subsec:infinitegibbs}

In this subsection, with the notation of the beginning of Chapter
\ref{sec:growth}, we still assume that $\delta_{\Ga,\,F}<+\infty$, but
we now assume that the Gibbs measure $m_F$ on $T^1M$ is infinite. We
state (weak) analogs of the previous asymptotic results in this
chapter, generalising the case $F=0$ in \cite{Roblin03}, leaving the
adaptation of the proofs to the reader.

Babillot's theorem \cite[Theo.~1]{Babillot02b}, saying that the
geodesic flow $(\phi_t)_{t\in\RR}$ on $T^1M$ is mixing for the measure
$m_F$, is also valid if $m_F$ is infinite, with the standard definition
of the mixing property of infinite measures.

\btheo [Babillot]
\label{theo:babmixinfini}\index{Babillot!mixing theorem}%
\index{theorem@Theorem!of Babillot!of mixing of the geodesic flow}
Assume that $\delta_{\Ga,\,F}<+\infty$, that the geodesic flow
$(\phi_t)_{t\in\RR}$ on $T^1M$ is topologically mixing on its
(topological) non-wandering set $\Omega\Ga$ (or equivalently that the
length spectrum of $M$ is not contained in a discrete subgroup of
$\RR$), and that $m_F$ is infinite.  For all bounded Borel subsets $A$
and $B$ in $T^1M$, we have
$$
\lim_{t\ra+\infty} m_F(A\cap \phi_t B)=0\;.
$$
\etheo

Equivalently, the mixing property may be written as follows: for all
bounded Borel subsets $A$ and $B$ in $T^1\wt M$,
$$
\lim_{t\ra+\infty} \;\sum_{\ga\in\Ga}\wt m_F(A\cap \phi_t \ga B)=0\;,
$$
which is the appropriate replacement for the use of the mixing
property in the proof of Proposition \ref{prop:stepone}.

Corollaries \ref{coro:corodeux} and \ref{coro:corocinq} on the
asymptotic of the orbital counting function become the following
result.

\btheo Assume that $\delta_{\Ga,\,F}<+\infty$, that the geodesic flow
on $T^1M$ is topologically mixing, and that $m_F$ is infinite.  Then,
as $t\ra+\infty$, for every $c>0$, 
$$
\sum_{\ga\in\Ga,\;t-c<d(x,\,\ga y)\leq t}e^{\int_x^{\ga y}\wt F}=
\;\operatorname{o}\,(e^{\delta_{\Ga,\,F}\;t})\,,
$$
and if $\delta_{\Ga,\,F}>0$, then 
$$
\sum_{\ga\in\Ga,\;d(x,\,\ga y)\leq t}e^{\int_x^{\ga y}\wt F}=
\;\operatorname{o}\,(e^{\delta_{\Ga,\,F}\;t})\,.
$$
\etheo

Theorems \ref{theo:equidisperorb} and \ref{theo:equidisperorbbis} on the
asymptotic of the period counting function become the following
result.

\btheo Assume that $\delta=\delta_{\Ga,\,F}<+\infty$, that the geodesic flow
on $T^1M$ is topologically mixing, and that $m_F$ is infinite.  Then,
as $t\ra+\infty$, for every $c>0$,
$$
e^{-\delta\; t}
\sum_{g\in\Per'(t)-\Per'(t-c)} e^{\L_g(F)}\;\L_{g}
\;\;\weakstar\;\;0\,,
\;\;\;\; t\;e^{-\delta\; t}
\sum_{g\in\Per'(t)-\Per'(t-c)} 
e^{\L_g(F)}\;\frac{\L_{g}}{\ell(g)}
\;\;\weakstar\;\;0\,,
$$
and if $\delta>0$, then 
$$
e^{-\delta\; t}
\sum_{g\in\Per'(t)} \;e^{\L_g(F)}\;\L_{g}
\;\;\weakstar\;\;0\,,
\;\;\;\; t\;e^{-\delta\; t}
\sum_{g\in\Per'(t)} 
e^{\L_g(F)}\;\frac{\L_{g}}{\ell(g)}
\;\;\weakstar\;\;0\,. 
$$
\etheo

\section{The ergodic theory of the strong 
unstable foliation}
\label{sec:ergtheounistabfolia}

Let $(\wt M,\Ga,\wt F)$ be as in the beginning of Chapter
\ref{sec:negacurvnot}: $\wt M$ is a complete simply connected
Riemannian manifold, with dimension at least $2$ and pinched sectional
curvature at most $-1$; $\Ga$ is a nonelementary discrete group of
isometries of $\wt M$; and $\wt F :T^1\wt M\ra \RR$ is a
H\"older-continuous $\Ga$-invariant map. Let $x_0\in\wt M$ and
$\delta=\delta_{\Ga,\,F}$. The aim of this chapter is to prove a
result of unique ergodicity for the strong unstable foliation of
$T^1M=\Ga\bs T^1\wt M$ endowed with the conditional measures of Gibbs
measures. We refer for instance to \cite{Walters82,BowMar77} for
generalities on unique ergodicity of measurable dynamical systems and
foliations, though in our case we will (and have to, in order to take
into account the potential) consider quasi-invariant measures instead
of invariant ones.

\subsection{Quasi-invariant transverse measures}
\label{subsec:defiquasinvtransv}

We first recall some definitions. Let $N$ be a smooth manifold of
dimension $n$ endowed with a smooth action of a discrete group $G$ and
a $\operatorname{C}^0$ foliation $\F$ of dimension $k$ invariant by
$G$. A {\em transversal}\index{transversal} to $\F$ is a
$\operatorname{C}^0$ submanifold $T$ of $N$ such that for every $x\in
T$, there exists a foliated chart $\varphi:U\ra\RR^k\times\RR^{n-k}$
at $x$, sending $x$ to $0$ and $T\cap U$ to $\{0\}\times
\RR^{n-k}$. Let $\T(\F)$ be the set of transversals to $\F$. A {\em
  holonomy map}\index{holonomy map} for $\F$ is a homeomorphism
$f:T\ra T'$ between two transversals to $\F$, such that $f(x)$ is in
the same leaf of $\F$ as $x$ for every $x\in T$.  A {\em $G$-invariant
  cocycle}\index{cocycle} for the foliation $\F$ is a continuous map
$c$ defined on the subspace of $N\times N$ of pairs of points in a
same leaf of $\F$, with real values, satisfying $c(u,v)+c(v,w)=c(u,w)$
and $c(u,v)=c(\ga u,\ga v)$ for all triples $(u,v,w)$ of points in a
same leaf of $\F$ and every $\ga\in G$.  Given a $G$-invariant cocycle
$c$ for $\F$, a {\em $G$-equivariant $c$-quasi-invariant transverse
  measure for $\F$}\index{transverse measure}%
\index{transverse measure!quasi-invariant} is a family
$\nu=(\nu_T)_{T\in\T(\F)}$, where $\nu_T$ is a (not necessarily
finite) locally finite (positive Borel) measure on the transversal
$T$, nonzero for at least one $T$, satisfying the following properties
for all $T,T'\in \T(\F)$:
\begin{itemize}
\item[(i)]
if $T'\subset T$, then  $(\nu_T)_{\mid T'}=\nu_{T'}$,
\item[(ii)] $\ga_*\nu_T=\nu_{\ga T}$ for every $\ga\in G$,
\item[(iii)] for all holonomy maps $f:T\to T'$, the measures
  $\nu_{T'}$ and $f_*\nu_T$ on $T'$ are equivalent, and their
  Radon-Nikodym derivative satisfies, for $\nu_T$-almost every $x\in T$,
$$
\frac{d\,\nu_{T'}}{d\,f_*\nu_T}(f(x))=e^{c(f(x),\,x)}\;.
$$
\end{itemize}
When $c=0$, such a family $\nu=(\nu_T)_{T\in\T(\F)}$ is said to be
{\it invariant under holonomy}\index{invariance under holonomy}%
\index{transverse measure!invariant under holonomy}: for all holonomy
maps $f:T\to T'$, we have $f_*\nu_T=\nu_{T'}$. Such a family
$\nu=(\nu_T)_{T\in\T(\F)}$
is {\em ergodic}\index{transverse measure!ergodic}\index{ergodic} if,
given a $G$-invariant subset $A$ of $N$, which is a union of leaves of
$\F$, whose intersection with any transversal is measurable, then
either for all transversals $T$ we have $\nu_T(A\cap T)=0$, or for
all transversals $T$ we have $\nu_T(\,^c\!A\cap T)=0$.

We will apply these definitions with $N=T^1\wt M$, $G=\Ga$, $\F$ the
strong unstable foliation $\wt\W^{\rm su}$ of $T^1\wt M$ (defined in
Subsection \ref{subsec:geodflowback}) and with the ($\Ga$-invariant)
cocycle $c=c_{\wt F}$ defined as follows: for all $v,w$ in the same
leaf of $\wt\W^{\rm su}$,
\begin{equation}\label{eq:deficsubwtF}
c_{\wt F}(v,w)=\lim_{t\to +\infty}\int_0^t \wt F(\phi_{-t}v)-\wt
F(\phi_{-t}w)\,dt=-C_{F\circ\iota-\delta,\,v_-}(\pi(v),\pi(w))\;.
\end{equation}
Note that when $F$ is constant, this cocycle is equal to $0$. Note
that $c_{\wt F}=c_{\wt F+\kappa}$ for every $\kappa\in\RR$, hence the
property of being a $c_{\wt F}$-quasi-invariant transverse measure is
unchanged by replacing $\wt F$ with $\wt F+\kappa$, for every
$\kappa\in\RR$.

We will say that a $\Ga$-equivariant $c_{\wt F}$-quasi-invariant
transverse measure $\nu$ for $\wt\W^{\rm su}$ {\it gives full measure
  to the negatively recurrent set}\index{giving full measure to the
  negatively recurrent set}%
\index{transverse measure!giving full measure to the negatively
  recurrent set} if for all transversals $T$ to $\wt\W^{\rm su}$, the
set of $v\in T$ such that $v_-$ belongs to the conical limit set
$\Lambda_c\Ga$ has full measure with respect to $\nu_T$, where
$\nu=(\nu_T)_{T\in\T(\wt\W^{\rm su})}$.

\medskip\noindent {\bf Remark (1). }  When $\Ga$ is torsion free, the
image of the strong unstable foliation $\wt\W^{\rm su}$ of $T^1\wt M$ by
the covering map $T^1\wt M\ra T^1M$ defines a foliation denoted by
$\W^{\rm su}$ on $T^1M$, the cocycle $c_{\wt F}$ defines by passing to the
quotient a cocycle $c_F$ for this foliation (already defined in
Subsection \ref{subsec:proofvaraprincip}, see Equation
\eqref{eq:defcocycleforunstab} and Equation \eqref
{eq:relatpetitcgrandC}), hence we may work directly with
$c_F$-quasi-invariant transverse measures for the foliation $\W^{\rm su}$
on $T^1M$ (defined as above by forgetting the equivariance property
(ii)). The case with torsion (useful when working with arithmetic
groups, for instance) makes it easier to work equivariantly on $T^1\wt
M$ than to work with singular foliations on $T^1M$.

\medskip\noindent {\bf Remark (2). } As seen in Subsection
\ref{subsec:condmesgibbs}, for every $v\in T^1M$, the stable
submanifold $T=W^{\rm s}(v)$ of $v$ in $T^1\wt M$ is a transversal to the
strong unstable foliation $\wt\W^{\rm su}$ of $T^1\wt M$, and every strong
unstable leaf $W^{\rm su}(w)$ with $w_-\neq v_+$ meets this transversal in
one and only one point $w'$.  If $U_T=\{w\in T^1\wt M\;:\;w_-\neq
v_+\}$, then the map $\psi_T:U_T\ra T$ defined by $w\mapsto w'$ is a
continuous fibration, whose fibre above $w'\in T$ is the strong
unstable leaf of $w'$.

In particular, any small enough transversal $T'$ to the strong
unstable foliation $\wt\W^{\rm su}$ of $T^1\wt M$ admits a holonomy map
$f:T'\ra T$ where $T$ is a submanifold of a stable submanifold. Hence
by the above properties (i) and (iii), a $\Ga$-equivariant $c_{\wt
  F}$-quasi-invariant transverse measure
$\nu=(\nu_T)_{T\in\T(\wt\W^{\rm su})}$ for $\wt\W^{\rm su}$ is uniquely
determined by the measures $\nu_T$ where $T$ is a stable submanifold of
$T^1\wt M$.  We will use this remark without further comment in what
follows.

\medskip\noindent {\bf Remark (3). } Let us give another justification
of our terminology of Gibbs measure, by further developing the
analogy with symbolic dynamics. We refer to Subsection
\ref{subsec:gibbsproperty} for the notation of a countable alphabet
$A$, the distance $d$ on $\Sigma=A^\ZZ$, the full shift
$\sigma:\Sigma\ra \Sigma$, the cylinders $\mathopen{[}a_{-n'},\dots,
a_n\mathclose{]}$ where $n,n'\in\NN$ and $a_{-n'},\dots, a_n\in A$, a
H\"older-continuous map $F:\Sigma\ra \RR$, and the definition of a
Gibbs measure for the potential $F$ on $\Sigma$ (see Equation
\eqref{eq:propGibbssymboldyn}).

For every $x=(x_i)_{i\in\ZZ}\in\Sigma$, we define the {\it strong unstable
  leaf}\index{strong!unstable leaf!in symbolic dynamics}%
\index{leaf!strong unstable!in symbolic dynamics} of $x$ by
$$
W^{\rm su}(x)=\{y=(y_i)_{i\in\ZZ}\in\Sigma\;:\; \exists\;
n\in\ZZ,\;\forall\; i\leq n,\;\;\; y_i=x_i\}\;.
$$
Note that these sets $W^{\rm su}(x)$ for all $x\in\Sigma$ are either
equal or disjoint. We define the {\it unstable Gibbs
  cocycle}\index{unstable!Gibbs cocycle}%
\index{Gibbs cocycle!unstable!in symbolic dynamics} of
$(\Sigma,\sigma, F)$ (compare with Equation \eqref{eq:deficsubwtF},
with \cite{Ruelle04,Zinsmeister96,Keller98}, and in particular with
\cite{HayRue92}) as the map $c_F$ which associates, to all $x,y$ in
the same strong unstable leaf, the real number
$$
c_F(x,y)=\;\sum_{i=1}^{+\infty}\;F(\sigma^{-i}x)-F(\sigma^{-i}y)\;,
$$
which exists by the H\"older-continuity of $F$, and satisfies the
cocycle properties $c_F(x,y)+c_F(y,z)=c_F(x,z)$ and $c_F(x,y) =
-c_F(y,x)$ if $x,y,z$ are in the same strong unstable leaf. Similarly,
let us define the {\it stable leaf}%
\index{stable leaf!in symbolic dynamics} of $x$ by
$$
W^{\rm s}(x)=\{y=(y_i)_{i\in\ZZ}\in\Sigma\;:\; \exists\;
n,k\in\ZZ,\;\forall\; i\geq n,\;\;\; y_i=x_{i+k}\}\;.
$$
We denote respectively by $\W^{\rm su}$ and $\W^{\rm s}$ the partitions of
$\Sigma$ in strong unstable leaves and in stable leaves.  

The projection map $p:\Sigma\ra \Sigma^+=A^\NN$ defined by
$(x_i)_{i\in\ZZ} \mapsto (x_i)_{i\in\NN}$ is a $1$-Lipschitz
fibration, whose fibres are contained in (strong) stable leaves. For
all $z=(z_i)_{i\in\NN}$ and $z'=(z'_i)_{i\in\NN}$ in $\Sigma^+$, the
map $f_{z',\,z}:p^{-1}(z)\ra p^{-1}(z')$, uniquely defined by
$(x_i)_{i\in\ZZ} \mapsto (y_i)_{i\in\ZZ}$ where $y_i=x_i$ for every
$i< 0$, is a Lipschitz homeomorphism. Note that $x$ and
$f_{z',\,z}(x)$ are in the same strong unstable leaf, for every $x\in
p^{-1}(z)$. Hence we will consider the family $(p^{-1}(z))_{z\in
  \Sigma^+}$ of fibres of $p$ as a family of transversals to the
partition into strong unstable leaves of $\Sigma$, and
$(f_{z',\,z}))_{z,\,z'\in \Sigma^+}$ as a family of holonomy maps
along the  strong unstable leaves of $\Sigma$.

We define a {\it $c_F$-quasi-invariant transverse
  measure}\index{transverse measure!in symbolic dynamics} for
$\W^{\rm su}$ to be a family $(\nu_z)_{z\in \Sigma^+}$, where $\nu_z$ is a
measure on $p^{-1}(z)$ for every $z\in \Sigma^+$, such that there exists
a constant $c>0$ such that for all $z,z'$ in $\Sigma^+$, with
$f=f_{z',\,z}$, the measures $\nu_{z'}$ and $f_*\nu_z$ have the same
measure class, and their Radon-Nikodym derivative satisfies, for
$\nu_z$-almost every $x\in p^{-1}(z)$,
$$
\frac{1}{c}\;e^{c_F(f(x),\,x)}\leq 
\frac{d\,\nu_{z'}}{d\,f_*\nu_z}(f(x))
\leq c\;e^{c_F(f(x),\,x)}\;.
$$
We will also consider $\nu_z$ as a measure on $\Sigma$ with support in
$p^{-1}(z)$.

By the cocycle property of the Radon-Nikodym derivatives, note that we
may take $c=1$, up to replacing the cocycle $c_F$ by a cohomologous
one.

\bprop \label{prop:quasiinvtransshift} Let $m_F$ be a probability
Gibbs measure for the potential $F$ on $\Sigma$. Then there exists a
$c_F$-quasi-invariant transverse measure $(\nu_z)_{z\in \Sigma^+}$ for
$\W^{\rm su}$ such that $m_F$ disintegrates, by the fibration
$p$, with family of conditional measures $(\nu_z)_{z\in \Sigma^+}$.
\eprop

When the alphabet $A$ is finite, there exists in fact a unique, up a
multiplicative constant, $c_F$-quasi-invariant transverse measure for
$\W^{\rm su}$, see \cite[Prop.~3.3]{BowMar77} when $F=0$ (compare with
\cite[Prop.~6 (i)]{PolSha94}) and \cite[Prop.~4.4]{BabLed96} for
general $F$. Analogous uniqueness results for geodesic flows are
the main goal of this chapter (see Remark (1) after Theorem
\ref{theo:uniergodtransv}).

\medskip \dem We refer for instance to Subsection
\ref{subsec:measupartentrop} and \cite{LeGall06} for basics on
conditional measures and martingales.

Let $(\nu_z)_{z\in \Sigma^+}$ be a family of conditional measures of
$m_F$ for the fibration $p$. By definition, $\nu_z$ is a probability
measure on $p^{-1}(z)$ such that for all Borel subsets $B$ in
$\Sigma$, the map $z\mapsto \nu_z(B)$ is measurable and we have
\begin{equation}\label{eq:desintshift}
\int_{x\in\Sigma} \mathbbm{1}_B(x) \,dm_F(x)=
\int_{z\in\Sigma^+}\nu_z(B)\,d (p_*m_F)(z)\;.
\end{equation}
Note that $(\nu_z)_{z\in \Sigma^+}$ is well defined up to a set of
$p_*m_F$-measure zero of $z\in\Sigma^+$.

\blemm \label{lem:martingalearg} 
For $p_*m_F$-almost every $z=(z_i)_{i\in\NN}\in\Sigma^+$, for
all $a_{-n'}, \dots, a_{-1}\in A$, we have
$$
\nu_z(\mathopen{[}a_{-n'}, \dots, a_{-1}\mathclose{]})= \lim_{k\ra+\infty}
\frac{m_F(\mathopen{[}a_{-n'}, \dots, a_{-1}\mathclose{]}\cap 
\mathopen{[}z_{0}, \dots, z_{k}\mathclose{]})}
{m_F(\mathopen{[}z_{0}, \dots, z_{k}\mathclose{]})}\;.
$$
\elemm

\dem For every $k\in\NN$, let $\F_k$ be the $\sigma$-algebra generated
by the cylinders $\mathopen{[}a_{0}, \dots, a_{k}\mathclose{]}$ in
$\Sigma$, where $a_0,\dots, a_k\in A$. Then $(\F_k)_{k\in\NN}$ is a
filtration in the Borel $\sigma$-algebra $\B$ of $\Sigma$, and a map
$f:\Sigma\ra \RR$ is $\F_\infty$-measurable for
$\F_\infty=\bigvee_{k\in\NN} \F_k$ if and only if $f(x)$ depends only
on $p(x)$ for all $x\in\Sigma$ and $z\mapsto f(p^{-1}(z)):\Sigma^+\ra
\RR$ is measurable. For every Borel subset $B$ in the probability
space $(\Sigma,m_F)$, the sequence $(X_k=E\mathopen{[}\mathbbm{1}_B
\mid \F_k\mathclose{]})_{k\in\NN}$ is a (basic example of a)
martingale adapted to the filtration $(\F_k)_{k\in\NN}$: we have
$E\mathopen{[}X_{k+1}\mid\F_{k}\mathclose{]}=X_k$ for every
$k\in\NN$. Since $\mathbbm{1}_B$ is $m_F$-integrable, by the $\LL^1$
martingale convergence theorem, the random variable $X_n$ converges
almost surely (and in $\LL^1$) to $X_\infty=
E\mathopen{[}\mathbbm{1}_B\mid\F_\infty \mathclose{]}$. By the
definition of the conditional expectation and by Equation
\eqref{eq:desintshift}, for $m_F$-almost every $x\in \Sigma$, we have
$E\mathopen{[}\mathbbm{1}_B \mid\F_\infty\mathclose{]}(x)=
\nu_{p(x)}(B)$. Since a Gibbs measure gives a nonzero mass to any
cylinder (see Equation \eqref{eq:propGibbssymboldyn}) and since a map
$f:\Sigma\ra\RR$ is $\F_k$-measurable if and only if $f(x)$ depends
only on the components $x_0,\dots, x_k$ of $x$, for $m_F$-almost every
$x=(x_i)_{i\in\ZZ} \in \Sigma$, we have $E\mathopen{[}\mathbbm{1}_B
\mid\F_k\mathclose{]}(x)= \frac{m_F(B\cap \mathopen{[}x_{0}, \, \dots,
  \, x_{k}\mathclose{]})}{m_F(\mathopen{[}x_{0}, \, \dots,
  \,x_{k}\mathclose{]})}$. The result follows.  \cqfd

\medskip For all sets $X$ and maps $g,g':X\ra \RR$, let us write $g
\asymp g'$ if there exists $c>0$ such that $\frac{1}{c}\;g(x)\leq
g'(x)\leq c\;g(x)$ for every $x\in X$. The above lemma and Equation
\eqref{eq:propGibbssymboldyn} imply that there exists
$\delta'_F\in\RR$ such that for $m_F$-almost every $x=(x_i)_{i\in\ZZ}
\in \Sigma$, we have
$$
\nu_{p(x)}(\mathopen{[}x_{-n'}, \dots, x_{-1}\mathclose{]})
\asymp \lim_{n\ra+\infty}
\frac{e^{\sum_{i=-n'}^n F(\sigma^i x)-(n+n'+1)\delta'_F}}
{e^{\sum_{i=0}^{n} F(\sigma^i x)-(n+1)\delta'_F}}=
e^{\sum_{i=1}^{n'} F(\sigma^{-i} x)-n'\,\delta'_F}\;.
$$
Hence for $p_*m_F$-almost all $z,z'\in\Sigma^+$ and for
$\nu_z$-almost every $x=(x_i)_{i\in\ZZ}\in p^{-1}(z)$, if $f=f_{z',\,z}$,
then
$$
\frac{\nu_{z'}(\mathopen{[}x_{-n'}, \dots, x_{-1}\mathclose{]})}
{\nu_{z}(\mathopen{[}x_{-n'}, \dots, x_{-1}\mathclose{]})}
\asymp e^{\sum_{i=1}^{n'} F(\sigma^{-i}f(x))-F(\sigma^{-i}x)}
\;.
$$
Proposition \ref{prop:quasiinvtransshift} now follows by taking limits
and modifying the family $(\nu_z)_{z\in \Sigma^+}$ on a set of
$p_*m_F$-measure zero.
 \cqfd

\subsection{Quasi-invariant measures on the space of 
horospheres}
\label{subsec:quasiinvhorospace}

Let us start this subsection by giving a few definitions. Given a
topological space $X$ endowed with a continuous action of a discrete
group $G$, a {\it real cocycle}\index{cocycle}\index{real cocycle} on
$X$ is a continuous map $\breve{c} :G\times X\ra \RR$ such that, for
all $x\in X$ and $\ga,\ga'\in G$, we have
$$
\breve c\,(\ga\ga',x)=\breve c\,(\ga,\ga' x)+\breve c\,(\ga',x)\;.
$$
A real cocycle $\breve c\,'$ on $X$ is {\it
  cohomologous}\index{cocycle!cohomologous}%
\index{real cocycle!cohomologous}\index{cohomologous} to $\breve c$
via a continuous map $f: X\ra \RR$ if $\breve c\,'(\ga,x)-\breve
c(\ga,x)=f(\ga x) -f(x)$ for all $x\in X$ and $\ga\in G$.  Given a
real cocycle $\breve c$ on $X$, a {\it $\breve c$-quasi-invariant
  measure}\index{measure!quasi-invariant}%
\index{quasi-invariant!measure} on $X$ is a (not necessarily finite)
locally finite nonzero (positive Borel) measure $\breve{\nu}$ on $X$
such that for every $\ga\in G$, the measures $\ga_*\breve{\nu}$ and
$\breve{\nu}$ are equivalent, with, for $\breve{\nu}$-almost every
$x\in X$,
$$
\frac{d\ga_*\breve{\nu}}{d\breve{\nu}}(x)=
e^{\breve c\,(\ga^{-1},\,x)}\;.
$$
We will be interested in the case when $X$ is the space $\H_{\wt M}$
of horospheres of $\wt M$, endowed with the topology of Hausdorff
convergence on compact subsets and the action of $G=\Ga$ on subsets of
$\wt M$, and when $\breve c$ is the real cocycle
$$
\breve c_{\wt F}=\breve c_{\wt F,\,x_0}:
(\ga,H)\mapsto C_{F\circ\iota,\,H_\infty}(x_0,\ga^{-1} x_0)\;,
$$
where $H$ is a horosphere centred at $H_\infty$ and $\ga\in\Ga$.

Note that $\breve c_{\wt F}=0$ if $\wt F=0$, that $\breve c_{\wt F}=
\breve c_{\wt F,\,x_0}$ does depend on $x_0$, but that $\breve c_{\wt
  F,\,x'_0}$ is cohomologous to $\breve c_{\wt F,\, x_0}$ via the map
$f:H\mapsto C_{F\circ\iota,\,H_\infty}(x'_0,x_0)$ for every
$x'_0\in\wt M$. If two potentials $\wt F^*$ and $\wt F$ are
cohomologous via the map $\wt G$, then $\wt F^*\circ\iota$ and $\wt
F\circ\iota$ are cohomologous via the map $- \wt G\circ \iota$, hence, by
Equation \eqref{eq:cohomologuecocycle}, the real cocycle $\breve
c_{\wt F^*}$ is cohomologous to $\breve c_{\wt F}$ via the map
$f:H\mapsto \wt G\circ\iota(v_{x_0\,H_\infty})$, with the notation of
Remark (1) just before Lemma \ref{lem:holderconseq}.

\medskip As explained for instance in \cite{Schapira03a, Schapira04a},
there exists a correspondence between $\Ga$-equivariant $c_{\wt
  F}$-quasi-invariant transverse measures $\nu=
(\nu_T)_{T\in\T(\wt\W^{\rm su}))}$ for $\wt\W^{\rm su}$ and $\breve c_{\wt
  F}$-quasi-invariant measures $\breve{\nu}$ on $\H_{\wt M}$, see
Proposition \ref{prop:corresqitrqimes} below.  This second viewpoint
is the one used in \cite{Roblin03} when $F=0$, in which case the
$\breve c_{\wt F,\,x_0}$-quasi-invariant measures on $\H_{\wt M}$ are
just the $\Ga$-equivariant measures, which explains why it is easier
to work with the second viewpoint then. The first viewpoint turns out
to be a bit more convenient when $F$ is not constant.

Let $\H^{\rm su}$ be the space of strong unstable leaves of $T^1\wt M$,
that is, the quotient topological space of $T^1\wt M$ by the
equivalence relation ``being in the same strong unstable leaf'', with
the quotient action of $\Ga$. We denote by $p^{\rm su}:T^1\wt M\ra
\H^{\rm su}$ the canonical projection.

Let $\H_{\wt M}$ be the space of horospheres in $\wt M$, endowed with
the Chabauty topology on closed subsets of $\wt M$ (induced by the
Hausdorff distance on compact subsets, see for instance \cite[\S
1.2]{Paulin07b}).

Note that the map which associates to a horosphere in $\wt M$ its outer
unit normal bundle is a $\Ga$-equivariant homeomorphism between
$\H_{\wt M}$ and $\H^{\rm su}$, whose inverse is the map (also denoted by
$\pi$) which associates to a strong unstable leaf $W$ the horosphere
$\pi(W)$. For every $v\in T^1\wt M$, the restriction of $\pi\circ
p^{\rm su}$ to the transversal $T=W^{\rm s}(v)$ to the foliation
$\wt\W^{\rm su}$ in $T^1\wt M$ is a homeomorphism onto its image in
$\H_{\wt M}$, and we will identify $T$ with its image by this map.

Let us endow $\partial_\infty\wt M\times\RR$ with the action of $\Ga$
defined as follows: for all $\ga\in \Ga$, $\xi\in \partial_\infty\wt
M$, and $t\in\RR$,
$$
\ga(\xi,t)=(\ga\xi, t+\beta_\xi(\ga^{-1}x_0,x_0))\;.
$$
Note that for all horospheres $H$ with centre $H_\infty$ and $u\in
H$, the real number $\beta_{H_\infty}(x_0,u)$ does not depend on
$u\in H$, it will be denoted by $\beta_{H_\infty}(x_0,H)$. The map
from $\H_{\wt M}$ to $\partial_\infty\wt M\times\RR$, defined by
$$
H\mapsto \big(H_\infty,\beta_{H_\infty}(x_0,H)\big)\;,
$$
is a $\Ga$-equivariant homeomorphism. It depends on $x_0$, but only up
to a translation (depending on the first factor) in the second factor
$\RR$.

\bprop\label{prop:corresqitrqimes} \hfill For every $\Ga$-equivariant
$c_{\wt F}$-quasi-invariant transverse measure $\nu=$ \\
$(\nu_T)_{T\in \T(\wt\W^{\rm su})}$ for $\wt\W^{\rm su}$, there exists one and
only one $\breve c_{\wt F}$-quasi-invariant measure
$\breve\nu=\breve\nu_{x_0}$ on $\H_{\wt M}$ such that for every $v\in
T^1\wt M$, if $T=W^{\rm s}(v)$, then for all $w\in T$ and $H=\pi(W^{\rm su}(w))$, we
have
\begin{equation}\label{eq:passqitransqimeas}
d\breve\nu(H)= e^{C_{F\circ\iota,\,w_-}(\pi(w),\,x_0)}\;d\nu_T(w)\;.
\end{equation}
Conversely, for every $\breve c_{\wt F}$-quasi-invariant measure
$\breve\nu$ on $\H_{\wt M}$, there exists one and only one
$\Ga$-equivariant $c_{\wt F}$-quasi-invariant transverse measure
$\nu=(\nu_T)_{T\in\T(\wt\W^{\rm su})}$ for $\wt\W^{\rm su}$ such that Equation
\eqref{eq:passqitransqimeas} holds for every stable leaf $T$.
\eprop

Note that the measure $\breve\nu=\breve\nu_{x_0}$ does depend on the
base point $x_0$, hence the aforementioned correspondence is not
canonical.

Note that the above measure $\breve\nu$ on $\H_{\wt M}$ gives full
measure to the set of horospheres centred at conical limit points if
and only if the transverse measure $\nu$ for $\wt\W^{\rm su}$ gives full
measure to the negatively recurrent set.

\medskip \dem The open subsets $T=W^{\rm s}(v)$ of $\H_{\wt M}$ cover
$\H_{\wt M}$ as $v$ ranges over $T^1\wt M$. Let us prove that the
measures on these open subsets defined by the member on the right of
Equation \eqref{eq:passqitransqimeas} glue together to define a
measure on $\H_{\wt M}$, that is, that two of them coincide on the
intersection of their domain.

\smallskip\noindent
\begin{minipage}{8.8cm} ~~~ Let $v,v'\in T^1\wt M$, and let $T=W^{\rm
    s}(v)$ and $T'=W^{\rm s}(v')$. The map from $\{w\in T\;:\;w_-\neq
  v'_+\}$ to $\{w'\in T'\;:\;w'_-\neq v_+\}$ which sends $w$ to the
  unique element $w'\in W^{\rm su}(w)$ such that $w'_+=v'_+$ is a
  holonomy map for the foliation $\wt\W^{\rm su}$. Note that
  $w'_-=w_-$ and that $\pi(w)$ and $\pi(w')$ lie on the same
  horosphere centred at $w'_-=w_-$. In particular,
  $\beta_{w_-}(\pi(w),\pi(w')=0$.
\end{minipage}\begin{minipage}{6cm}
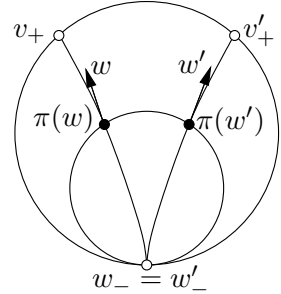
\begin{center}
\input{fig_changtransvuns.pstex_t}
\end{center}
\end{minipage}

\smallskip\noindent By the properties (i) and (iii) of $\nu$, we hence
have, for every $w\in \{w\in T\;:\;w_-\neq v'_+\}$,
\begin{align}
e^{C_{F\circ\iota,\,w'_-}(\pi(w'),\,x_0)}\;d\nu_{T'}(w')&=
e^{C_{F\circ\iota,\,w'_-}(\pi(w'),\,x_0)}\;e^{c_{\wt F}(w',\,w)}\;
d\nu_{T}(w)\nonumber \\& 
=e^{C_{F\circ\iota,\,w_-}(\pi(w'),\,x_0)-C_{F\circ\iota,\,w_-}(\pi(w'),\,\pi(w))}
\;d\nu_{T}(w)\nonumber\\& 
=e^{C_{F\circ\iota,\,w_-}(\pi(w),\,x_0)}\;d\nu_{T}(w)\;.\label{eq:holonqimeas}
\end{align}
This proves the existence and uniqueness of a nonzero measure
$\breve\nu$ on $\H_{\wt M}$ satisfying Equation
\eqref{eq:passqitransqimeas} for every stable leaf $T$.

Let us now prove that $\breve\nu$ is $\breve c_{\wt
  F}$-quasi-invariant.  Let $\ga\in\Ga$,  $v\in T^1\wt M$, and
$T=W^{\rm s}(v)$, so that $\ga T=W^{\rm s}(\ga v)$. Let $w\in \{u\in
\ga T\;:\; u_-\neq v_+\}$ and  $H=\pi(W^{\rm ss}(w))$. Using respectively

$\bullet$~ the definition of $\breve\nu$
for the second equality,

$\bullet$~ the invariance of the Gibbs cocycle (see
Equation \eqref{eq:equivcocprop}) and the property (ii) of $\nu$ for
the third equality,

$\bullet$~ the property (iii) of $\nu$ for the holonomy map
from $\{u\in T\;:\; u_-\neq \ga v_+\}$ to $\{u\in \ga T\;:\;
u_-\neq v_+\}$ for the foliation $\wt\W^{\rm su}$ sending $w'$ to the
unique element $w\in W^{\rm su}(w')$ such that $w_+=\ga v_+$ for the fourth
equality, 

$\bullet$~ the fact that $w'_-=w_-$ and the definition of the cocycle
$c_{\wt F}$ for the fifth equality,

$\bullet$~ the cocycle properties (see Equation
\eqref{eq:cocycleprop}) of the Gibbs cocycle for the sixth equality,

$\bullet$~ the definition of the cocycle $\breve c_{\wt F}$ and the
facts that $\pi(W^{\rm su}(w'))=\pi(W^{\rm su}(w))=H$ and $H_\infty=w'_-$ for
the seventh equality, we have
\begin{align}
  d\ga_*\breve\nu(H)&=d\breve\nu(\ga^{-1}H)= e^{C_{F\circ\iota,\,(\ga^{-1}w)_-}
    (\pi(\ga^{-1}w),\,x_0)}\;d\nu_T(\ga^{-1}w) \nonumber\\ &
  =e^{C_{F\circ\iota,\,w_-} (\pi(w),\,\ga x_0)}\;d\nu_{\ga T}(w) =
  e^{C_{F\circ\iota,\,w_-} (\pi(w),\,\ga x_0)}\;e^{c_{\wt F}(w,\,w')}\;
  d\nu_{T}(w') \nonumber\\ &= e^{C_{F\circ\iota,\,w'_-} (\pi(w),\,\ga
    x_0)-C_{F\circ\iota,\,w'_-}(\pi(w),\,\pi(w'))}\; d\nu_{T}(w')
  \nonumber\\ &= 
e^{C_{F\circ\iota,\,w'_-} (x_0,\,\ga x_0)}\;e^{C_{F\circ\iota,\,w'_-}
    (\pi(w'),\,x_0)}\; d\nu_{T}(w') \nonumber\\ & = e^{\breve c_{\wt
      F} (\ga^{-1},\, H)}\;d\breve\nu(H)\;.\label{eq:equivqimeas}
\end{align}
This proves that $\breve\nu$ is $\breve c_{\wt
  F}$-quasi-invariant.

\medskip Conversely, let $\breve\nu$ be a $\breve c_{\wt
  F}$-quasi-invariant nonzero measure on $\H_{\wt M}$.  For every
stable leaf $T=W^{\rm s}(v)$, let us define 
$$
d\nu_T(w) =e^{-C_{F\circ\iota,\,w_-}(\pi(w),\,x_0)}\;
d\breve\nu\big(\pi(W^{\rm su}(w))\big)\;.
$$ 
The computations \eqref{eq:holonqimeas} and
\eqref{eq:equivqimeas} show that it is possible to define, using
restrictions and holonomy maps, a measure $\nu_T$ for any transversal
$T$ to $\wt\W^{\rm su}$ such that the family
$\nu=(\nu_T)_{T\in\T(\wt\W^{\rm su})}$ satisfies the properties (i),
(ii) and (iii).  \cqfd

\subsection{Classification of quasi-invariant transverse 
measures for $\W^{\rm su}$}
\label{subsec:classifquasinvtransv}

Let us now state our unique ergodicity result of the strong unstable
foliation, starting by defining the transverse measure which will play
the central role.

By Subsection \ref{subsec:condmesgibbs}, if $\delta=
\delta_{\Ga,\,F}<+\infty$, to every Patterson density
$(\mu^\iota_x)_{x\in\wt M}$ of dimension $\delta$ for
$(\Ga,F\circ\iota)$ is associated, for any $v\in T^1\wt M$ with stable
leaf $T=W^{\rm s}(v)$, a nonzero measure $\mu^\iota_T$ on $T$, defined
(independently of $x_0$), using the homeomorphism from $T$ to
$(\partial_\infty \wt M - \{v_+\}) \times \RR$ defined by $w'\mapsto
\big(w'_-,t=\beta_{v_+} (\pi(v),\pi(w'))\big)$, by (see Equation
\eqref{eq:defims})
\begin{equation}\label{eq:definuniqmestrans}
d\mu^\iota_T(w')=e^{C_{F\circ\iota-\delta,\,w'_-}(x_0,\,\pi(w'))}
\;d\mu^\iota_{x_0}(w'_-)\,dt\;.
\end{equation}
This measure is finite on the compact subsets of $T$, and if
$\mu^\iota_{x_0}({}^c\Lambda_c\Ga)=0$, then the set of $w'\in T$ such
that $w'_-$ belongs to $\Lambda_c\Ga$ has full measure with respect to
$\mu^\iota_T$.  By Remark (2) in Subsection
\ref{subsec:defiquasinvtransv}, by Equation \eqref{eq:propunmesstabl}
and Equation \eqref{eq:proptrmesstabl}, there exists hence a unique
$\Ga$-equivariant $c_{\wt F}$-quasi-invariant transverse measure
$(\mu^\iota_T) _{T\in\T( \wt \W^{\rm su})}$ for $\wt\W^{\rm su}$, giving full
measure to the negatively recurrent set, such that if $T=W^{\rm s}(v)$ for
some $v\in T^1\wt M$, then $\mu^\iota_T$ is the above measure. Note
that $(\mu^\iota_T) _{T\in\T( \wt \W^{\rm su})}$ is unchanged by replacing
$F$ with $F+\kappa$, for any $\kappa\in\RR$.

\bigskip
The next result proves that such a transverse measure is unique, under
our usual assumptions (see Chapter \ref{sec:finimixGibbs}, in
particular Theorem \ref{theo:babmix}, for situations where they are
satisfied).

\btheo \label{theo:uniergodtransv}%
\index{theorem@Theorem!Unique ergodicity of the strong unstable
  foliation} Let $\wt M$ be a complete simply connected Riemannian
manifold, with dimension at least $2$ and pinched negative curvature
at most $-1$. Let $\Ga$ be a nonelementary discrete group of
isometries of $\wt M $.  Let $\wt F:\wt M\ra \RR$ be a
$\Gamma$-invariant H\"older-continuous map.  Assume that
$\delta_{\Gamma,\,F}<+\infty$ and that the Gibbs measure $m_F$ on
$T^1M$ is finite and mixing under the geodesic flow.

Then the family $(\mu^\iota_T)_{T\in\T(\wt\W^{\rm su})}$ is, up to a
multiplicative constant, the unique $\Ga$-equi\-va\-riant $c_{\wt
  F}$-quasi-invariant transverse measure, giving full measure to
the negatively recurrent set, for the strong unstable foliation on
$T^1\wt M$.  Furthermore, $(\mu^\iota_T)_{T\in\T(\wt\W^{\rm su})}$ is ergodic.
\etheo

\noindent{\bf Remark (1). }
If $F=0$, this theorem (with the formulation of Corollary
\ref{coro:horocyclicflow} in the next remark if $M$ has dimension $2$)
is due to Furstenberg \cite{Furstenberg73} when $M$ is a compact
hyperbolic surface, to Dani \cite{Dani81}, Dani-Smillie
\cite{DanSmi84}, Ratner \cite{Ratner92} when $M$ is a finite volume
hyperbolic surface, to Bowen-Marcus \cite{BowMar77} when $M$ is a
compact or convex-cocompact manifold, via a coding argument (see also
Coudène \cite{Coudene09} for a very short dynamical argument in the
case of finite volume surfaces), and to Roblin
\cite[Chap.~6]{Roblin03} in full generality.  When $F=0$, $\wt M$
is a higher rank symmetric space of noncompact type, $G$ is the
identity component of the isometry group of $\wt M$, and $\Ga$ is a
lattice in $G$, we refer to Ratner's classification theorem of the
measures on $\Ga\bs G$ invariant under a unipotent subgroup of $G$
(see \cite{Ratner90,Ratner91,MarTom94,Morris05}). It would be
interesting to study the case $F$ non-constant in higher rank.

For general potentials $F$ (with the slightly different definition of
Gibbs measure mentioned in the beginning of Chapter \ref{sec:GPS}),
this theorem was proved in Babillot-Ledrappier \cite{BabLed96} in a
symbolic framework or when $M$ is an abelian Riemannian cover of a
compact manifold, and in Schapira \cite{Schapira04a} when $M$ is
compact or convex-cocompact. This theorem (with $\Ga$ torsion free and
$\wt F$ bounded) is stated, and its proof is sketched, in
\cite[Chap.~8.2.2]{Schapira03a} following Roblin's arguments in
\cite[Chap.~6]{Roblin03}. In the proof below, we will also follow
closely Roblin's strategy of proof, though the presence of potentials
requires many adjustments and new lemmas.

\medskip
\noindent{\bf Remark (2). } When $\wt M$ is a surface, it is
well-known since Furstenberg's work (see for instance
\cite{Schapira03b}) that this theorem corresponds to a unique
ergodicity result for horocyclic flows. When $\wt M$ is a symmetric
space, this theorem corresponds to a unique ergodicity result for some
unipotent group actions.

More precisely, assume first that $\wt M$ is a surface, fix one of the
two orientations of $\wt M$ and assume that $\Ga$ preserves the
orientation.  Then each horosphere is naturally oriented, say
anticlockwise. A {\it horocyclic flow}\index{horocyclic flow} on
$T^1\wt M$ is a continuous one-parameter group of homeomorphisms
$(h_s)_{s\in\RR}$ of $T^1\wt M$, such that, for every $v\in T^1\wt M$,
the map $s\mapsto h_sv$ is an orientation preserving homeomorphism
from $\RR$ to $W^{\rm su}(v)$ and for every $s\in\RR$, the homeomorphism
$v\mapsto h_sv$ of $T^1\wt M$ is $\Ga$-equivariant. By passing to the
quotient, it induces a continuous one-parameter group of
homeomorphisms of $T^1M$, called the {\it horocyclic flow}%
\index{horocyclic flow} on $T^1M$ and also denoted by
$(h_s)_{s\in\RR}$.  One way to define $h_sv$ is as the element of
$W^{\rm su}(v)$ on the right of $v$ if $s\geq 0$, and on its left if
$s\leq 0$, such that the Riemannian length of the arc of the horocycle
$\pi(W^{\rm su}(v))$ from $\pi(v)$ to $\pi(h_s(v))$ is $s$. In constant
curvature $-1$, this horocycle flow satisfies furthermore that
$\phi_t\circ h_s= h_{s\,e^t}\circ \phi_t$ for all $s,t\in\RR$. But in
variable curvature, this is no longer always true, and the Riemannian
parametrisation is not always the appropriate one. See for instance
\cite{Marcus75} which constructs, when $M$ is compact, a horocycle
flow satisfying, besides more precise regularity properties, that
$\phi_t\circ h_s= h_{s\,e^{\delta_\Ga t}}\circ \phi_t$ for all
$s,t\in\RR$.


A locally finite (positive Borel) measure $\wt m$ on $T^1\wt M$ is {\it
  $c_{\wt F}$-quasi-invariant under the horocyclic flow}%
\index{measure!quasi-invariant!under the horocyclic flow}%
\index{quasi-invariant!under the horocyclic flow}
$(h_s)_{s\in\RR}$ if for every $s\in\RR$, the measures $\wt m$ and
$(h_s)_*\wt m$ are mutually absolutely continuous,
and satisfy, for $\wt m$-almost every $v\in T^1\wt M$,
$$
\frac{d\,(h_s)_*\wt m}{d\,\wt m}(h_s v)=e^{c_\wt F(v,\,h_sv)}\;.
$$

More generally, assume that there exists a continuous action (on the
left) of a unimodular real Lie group $\UU$ on $T^1\wt M$, commuting
with the action of $\Ga$, such that for every $v\in T^1\wt M$ the map
$u\mapsto u\cdot v$ is a homeomorphism from $\UU$ to $W^{\rm su}(v)$.
Such an action exists for instance when $\wt M$ is a symmetric space,
that is, a real, complex, quaternionic or octonionic hyperbolic space
(see for instance \cite{Parker07}), since a nilpotent Lie group is
unimodular. We say that a locally finite measure $\wt m$ on $T^1\wt M$
is {\it $c_{\wt F}$-quasi-invariant under the action of $\UU$}%
\index{measure!quasi-invariant!under a unipotent group}%
\index{quasi-invariant!under a unipotent group} if for every
$u\in\UU$, the measures $\wt m$ and $u_*\wt m$ are absolutely
continuous with respect to each other, and satisfy, for $\wt
m$-almost every $v\in T^1\wt M$,
$$
\frac{d\,u_*\wt m}{d\,\wt m}(u\cdot v)=e^{c_\wt F(v,\,u\cdot v)}\;.
$$

\bcoro \label{coro:horocyclicflow} Under the assumptions of Theorem
\ref{theo:uniergodtransv}, there exists, up to a
multiplicative constant, a unique locally finite $\Ga$-in\-va\-riant
measure on $T^1\wt M$, giving full measure to $\{w\in T^1\wt M \;:\;
w_-\in\Lambda_c\Ga\}$, which is $c_{\wt F}$-quasi-invariant under
the action of $\UU$.  
\ecoro

The special case when $\UU=\RR$ gives that if $\wt M$ has dimension
$2$ and if $(h_s)_{s\in\RR}$ is a horocyclic flow, then under the
assumptions of Theorem \ref{theo:uniergodtransv}, there exists, up to
a multiplicative constant, a unique locally finite $\Ga$-in\-va\-riant
measure on $T^1\wt M$, giving full measure to $\{w\in T^1\wt M \;:\;
w_-\in\Lambda_c\Ga\}$, which is $c_{\wt F}$-quasi-invariant under
$(h_s)_{s\in\RR}$.

\medskip \dem We start the proof by a construction which, when
$\UU=\RR$, is due to Burger \cite{Burger90} (see more generally
\cite[\S 1C]{Roblin03}) when $F=0$, and to the last author
\cite{Schapira03b,Schapira04a} (see also Lemma
\ref{lem:constructmesnumu} below) when $\Ga$ is torsion free, with a
local approach.

Recall (see Subsection \ref{subsec:condmesgibbs} before Proposition
\ref{prop:disintegrGibbs}, exchanging the stable and unstable
foliations) that, for every $v\in T^1\wt M$, with $T=W^{\rm s}(v)$, the
map $\psi_{T}$ from the open set $U_{T}=\{w\in T^1\wt M\;:\; w_-\neq
v_+\}$ to $T$, sending $w\in U_{T}$ to the unique element $w'$ in
$W^{\rm su}(w)\cap T$, is a continuous fibration over $T$, whose fibre
over $w'\in T$ is the strong unstable leaf of $w'$. (Note that $U_{T}$
is invariant under $\UU$, and this fibration is a principal fibration
under $\UU$.)

Let $\lambda$ be a Haar measure on $\UU$ (say the usual Lebesgue
measure when $\UU=\RR$) and, for every $v\in T^1\wt M$, let
$\lambda_{W^{\rm su}(v)}$ be the push-forward measure of $\lambda$ by the
map $u\mapsto u\cdot v$, whose support is $W^{\rm su}(v)$. Since $\UU$
acts transitively on $W^{\rm su}(v)$ and since $\lambda$ is invariant
under right translations, the measure $\lambda_{W^{\rm su}(v)}$ indeed
depends only on the strong unstable leaf of $v$. Since $\lambda$ is
invariant under left translations, we have $u_*\lambda_{W^{\rm su}(v)}=
\lambda_{W^{\rm su}(v)}$ for every $u\in\UU$. Since the action of $\Ga$ on
$T^1\wt M$ commutes with the action of $\UU$, for every $\ga\in\Ga$,
we have $\ga_*\lambda_{W^{\rm su}(v)}=\lambda_{W^{\rm su}(\ga v)}$. Note that
since the action is continuous, for every $v\in T^1\wt M$, with
$T=W^{\rm s}(v)$, for every $f\in \C_c(T^1\wt M;\RR)$ with support contained
in $U_T$, the map from $T$ to $\RR$ defined by
$$
w'\mapsto \int_{w\in W^{\rm su}(w')}f( w)\;e^{c_{\wt F} (w,\, w')} \;d
\lambda_{W^{\rm su}(w')}(w)= \int_{u\in\UU}f(u\cdot w')\;e^{c_{\wt F}
  (u\cdot w',\, w')} \;d\lambda(u)
$$ 
is continuous with compact support.

Let $\nu= (\nu_T)_{T\in\T(\wt \W^{\rm su})}$ be a $c_{\wt
  F}$-quasi-invariant transverse measure for $\wt\W^{\rm su}$. We claim
that there exists a unique (positive Borel) measure $\wt
m^{\nu,\,\lambda}$ on $T^1\wt M$ such that for every $v\in T^1\wt M$,
with $T=W^{\rm s}(v)$, the restriction to $U_{T}$ of the measure $\wt
m^{\nu,\,\lambda}$ disintegrates by the fibration $\psi_{T}$ over the
measure $\nu_{T}$, with conditional measure on the fibre $W^{\rm su}(w')$
of $w'\in T$ the measure $e^{c_{\wt F}(w,\, w')} \;d
\lambda_{W^{\rm su}(w')}(w)$: for every $f\in \C_c(T^1\wt M;\RR)$ with support
contained in $U_T$, we have
$$
\int_{w\in T^1\wt M}f(w)\;d\wt m^{\nu,\,\lambda}(w)=
\int_{w'\in T}\int_{u\in\UU}f(u\cdot w')\;e^{c_{\wt F}(u\cdot w',\, w')} 
\;d\lambda(u)\;d\nu_{T}(w')\;.
$$
We refer to the proof of Lemma \ref{lem:constructmesnumu} for a proof
of this claim.

\blemm The map $\nu\mapsto \wt m^{\nu,\,\lambda}$ is a bijection from the
set of $\Ga$-equi\-va\-riant $c_{\wt F}$-quasi-invariant transverse
measures for $\wt\W^{\rm su}$, which give full measure to the negatively
recurrent set, to the set of $\Ga$-in\-va\-riant locally finite
measures on $T^1\wt M$ which are $c_{\wt F}$-quasi-invariant under the
action of $\UU$ and give full measure to $\{w\in T^1\wt M \;:\;
w_-\in\Lambda_c\Ga\}$.  
\elemm

\dem By construction, $\nu$ is $\Ga$-equi\-va\-riant if and only if
$\wt m^{\nu,\,\lambda}$ is $\Ga$-invariant, and $\nu_T$ is locally
finite for every $T$ if and only if $\wt m^{\nu,\,\lambda}$ is locally
finite. For every $v\in T^1\wt M$, with $T=W^{\rm s}(v)$, since the
support of $\lambda_{W^{\rm su}(v)}$ is equal to $W^{\rm su}(v)$, the measure
$\nu_T$ gives full measure to $\{w'\in T\;:\; w'_-\in\Lambda_c\Ga\}$
if and only if the restriction of $\wt m^{\nu,\,\lambda}$ to $U_T$
gives full measure to $\{w\in U_T\;:\; w_-\in\Lambda_c\Ga\}$.

For all $v\in T^1\wt M$, $w\in T=W^{\rm s}(v)$ and $u\in\UU$, since
$\lambda_{W^{\rm su}(w')}$ is invariant under $\UU$ and by the cocycle
property of $c_{\wt F}$, we have
$$
d\wt m^{\nu,\,\lambda}(u\cdot w)=
\int_{w'\in T}e^{c_{\wt F}(u\cdot w,\, w')} 
\;d\lambda_{W^{\rm su}(w')}(u\cdot w)\;d\nu_{T}(w')=
e^{c_{\wt F}(u\cdot w,\, w)}\;d\wt m^{\nu,\,\lambda}(w) \;.
$$
Hence $\wt m^{\nu,\,\lambda}$ is $c_{\wt F}$-quasi-invariant under the
action of $\UU$.

Conversely, let $\wt m$ be a $\Ga$-invariant locally finite measure on
$T^1\wt M$ which is $c_{\wt F}$-quasi-invariant under the action of
$\UU$. Then given any $v\in T^1\wt M$, with $T=W^{\rm s}(v)$, for the
principal fibration $U_T\ra T$ with group $\UU$, the conditional
measures of $m_{\mid U_T}$ on the fibre $W^{\rm su}(w')=\UU\cdot w'$ of
almost every $w'\in T$ are invariant under $\UU$. This conditional
measure may be taken to be equal to $\lambda_{W^{\rm su}(w')}$, by the
uniqueness property of Haar measures, up to multiplying by a
measurable function the (locally finite) measure $\nu_T$ on $T$ over
which $\wt m$ disintegrates. Then the family $\nu= (\nu_T)_{T\in\T(\wt
  \W^{\rm su})}$ is clearly a $c_{\wt F}$-quasi-invariant transverse
measure for $\wt\W^{\rm su}$.  
\cqfd

\medskip Since for every $t>0$, we have $\wt m^{t\nu,\,\lambda}=t\,\wt
m^{\nu,\,\lambda}$, Corollary \ref{coro:horocyclicflow} indeed follows
from Theorem \ref{theo:uniergodtransv}.  
\cqfd

\bigskip 
\noindent{\bf Proof of Theorem \ref{theo:uniergodtransv}. }  Up to
adding a large enough constant to $F$, which does not change the
statement, we assume that $\delta=\delta_{\Ga,\,F}$ is positive. Let
$(\mu^\iota_{x}) _{x\in\wt M}$ and $(\mu_{x})_{x\in\wt M}$ be the
Patterson densities of the same dimension $\delta$ for $(\Ga, F\circ
\iota)$ and $(\Ga, F)$, such that $m_F$ is the Gibbs measure on $T^1M$
(induced by the Gibbs measure $\wt m_F$ on $T^1\wt M$) associated with
this pair of Patterson densities. Note that they are uniquely defined
up to a scalar multiple, and give full measure to the conical limit
set, by Corollary \ref{coro:uniqGibstate} and Corollary
\ref{coro:uniqpatdens}. Hence the family $\mu^\iota= (\mu^\iota_T)
_{T\in\T(\wt\W^{\rm su})}$ associated to $(\mu^\iota_{x}) _{x\in\wt M}$ by
Equation \eqref{eq:definuniqmestrans} is indeed a $\Ga$-equivariant
$c_{\wt F}$-quasi-invariant transverse measure for $\wt\W^{\rm su}$,
giving full measure to the negatively recurrent set.

The scheme of the eight steps proof of Theorem
\ref{theo:uniergodtransv} is the following one. Given a second such
family $\nu$, we first construct a measure $m^{\nu,\,\mu}$ which is a
quasi-product measure of the measures defined by $\nu$ on the stable
leaves and the measures defined by $m_F$ on the strong unstable
leaves, so that $m_F=m^{\mu^\iota,\,\mu}$. We then introduce diffusion
operators $(I_r)_{r>0}$ along the strong unstable leaves acting on the
locally integrable functions on $T^1M$. We use the mixing property of
the geodesic flow (see Subsection \ref{subsec:mixing}) combined with
the action of these diffusion operators to prove that $m_F=
m^{\mu^\iota,\,\mu}$ is absolutely continuous with respect to
$m^{\nu,\,\mu}$ (this is the most technical part, which will require
several steps). We then conclude that $\nu$ and $\mu^\iota$ are
homothetic by an easy argument of disintegration and ergodicity.

\bigskip\noindent{\bf Step 1. }  Let $\nu= (\nu_T)_{T\in\T(\wt
  \W^{\rm su})}$ be a $\Ga$-equivariant $c_{\wt F}$-quasi-invariant
transverse measure for $\wt\W^{\rm su}$, giving full measure to the
negatively recurrent set. In this first step, we construct a measure
$\wt m^{\nu,\,\mu}$ on $T^1\wt M$ which, for every stable leaf $T$,
disintegrates, with conditional measures in the class of the strong
unstable measures $\mu_{W^{\rm su}(w')}$ (as does the Gibbs measure), over
the measure $\nu_{T}$ (instead of over the measure $\mu^\iota_{T}$
as does the Gibbs measure).

Recall (see Subsection \ref{subsec:condmesgibbs} before Proposition
\ref{prop:disintegrGibbs}, exchanging the stable and unstable
foliations) that, for every $v\in T^1\wt M$, with $T=W^{\rm s}(v)$, the
map $\psi_{T}$ from the open set $U_{T}=\{w\in T^1\wt M\;:\; w_-\neq
v_+\}$ to $T$, sending $w$ to the unique element $w'$ in
$W^{\rm su}(w)\cap T$, is a continuous fibration over $T$, whose fibre
over $w'\in T$ is the strong unstable leaf of $w'$.

\blemm\label{lem:constructmesnumu} There exists a unique (positive
Borel) measure $\wt m^{\nu,\,\mu}$ on $T^1\wt M$ such that for every
$v\in T^1\wt M$, with $T=W^{\rm s}(v)$, the restriction to $U_{T}$ of
the measure $\wt m^{\nu,\,\mu}$ disintegrates by the fibration
$\psi_{T}$ over the measure $\nu_{T}$, with conditional measure
$e^{c_{\wt F}(w,\, w')} \;d\mu_{W^{\rm su}(w')} (w)$ on the fibre
$W^{\rm su}(w')$ of $w'\in T$: for every $w\in T^1\wt M$ such that
$w_-\neq v_+$, we have
\begin{equation}\label{eq:defimnumu}
d\wt m^{\nu,\,\mu}(w)=\int_{w'\in T}\;e^{c_{\wt F}(w,\, w')} 
\;d\mu_{W^{\rm su}(w')}(w)\;d\nu_{T}(w')\;.
\end{equation}
Furthermore, $\wt m^{\nu,\,\mu}$ is $\Ga$-invariant and locally
finite, and the set $\wt\Omega_c\Ga$ of elements $v\in T^1\wt M$ such
that $v_-,v_+\in \Lambda_c\Ga$ has full measure for $\wt m^{\nu,\,\mu}$.
\elemm

Hence the measure $\wt m^{\nu,\,\mu}$ defines a locally finite measure
$m^{\nu,\,\mu}$ on $T^1M$, which gives full measure to the set of
two-sided recurrent elements.

Note that by the disintegration of the Gibbs measure (see Proposition
\ref{prop:disintegrGibbs} (1), exchanging the stable and unstable
foliations), we have
$$
\wt m_F=\wt m^{\mu^\iota,\,\mu}\;\;\;{\rm and}\;\;\; m_F=m^{\mu^\iota,\,\mu}\;.
$$

\dem Let $v,v'\in T^1\wt M$, and $T=W^{\rm s}(v)$, $T'=W^{\rm s}(v')$. The map
from $U_{T'}\cap T=\{w'\in T\;:\;w'_-\neq v'_+\}$ to $U_{T}\cap
T'=\{w''\in T'\;:\;w''_-\neq v_+\}$ sending $w'$ to the unique element
$w''$ of $W^{\rm su}(w')\cap T'$ is a holonomy map for the strong unstable
foliation. Hence $d\nu_{T'}(w'')= e^{c_{\wt F}(w'',\,w')}d\nu_{T}(w')$
by the $c_{\wt F}$-quasi-invariance property of $\nu$. Since
$W^{\rm su}(w')= W^{\rm su}(w'')$ and by the cocycle property of $c_{\wt F}$,
for every $w\in U_T\cap U_{T'}$, we therefore have
$$
\int_{w'\in T}\;e^{c_{\wt F}(w,\, w')} 
\;d\mu_{W^{\rm su}(w')}(w)\;d\nu_{T}(w')=
\int_{w''\in T'}\;e^{c_{\wt F}(w,\, w'')} 
\;d\mu_{W^{\rm su}(w'')}(w)\;d\nu_{T'}(w'')\;.
$$
Hence the measures on the sets $U_T$ defined by the right hand side of
Equation \eqref {eq:defimnumu} glue together to define a measure $\wt
m^{\nu,\,\mu}$ on $T^1\wt M$. The uniqueness follows from the fact
that the open sets $U_T$ cover $T^1\wt M$.

The measure $\wt m^{\nu,\,\mu}$ is locally finite since $\nu_{T}$ and
$\mu_{W^{\rm su}(w')}$ are locally finite. It is $\Ga$-invariant since the
families $\nu$ and $(\mu_{W^{\rm su}(w')})_{w'\in T^1\wt M}$ are
$\Ga$-equivariant and the cocycle $c_{\wt F}$ is $\Ga$-invariant. The
last claim of the lemma follows from the fact that $\nu_T$ gives full
measure to the set of $w'\in T$ such that $w'_-\in \Lambda_c\Ga$ and
$\mu_{W^{\rm su}(w')}$ to the elements $w\in W^{\rm su}(w')$ such that $w_+\in
\Lambda_c\Ga$. \cqfd

\medskip Note that the measure $m^{\nu,\,\mu}$ is a priori not finite
(and could be infinite if $m_F$ was not assumed to be finite, see for
instance \cite{Burger90}). We will turn it into a probability measure
by multiplying it by a convenient bump function.

We fix a nonnegative (continuous) $\rho\in\C_c(T^1M;\RR)$ with large
enough (compact) support, so that
$$
\int_{T^1M}\rho \;dm_F>0\;\;\;{\rm and}\;\;\;
\int_{T^1M}\rho \;dm^{\nu,\,\mu}=1\;.
$$
With $\prT:T^1\wt M\ra T^1M$ the canonical projection, let $\wt
\rho=\rho\circ \prT$ be the lift of $\rho$ to $T^1\wt M$.  Consider the
measure $\wt \Pi=\wt \rho\; \wt m^{\nu,\,\mu}$ on $T^1\wt M$, which
induces a probability measure
$$
\Pi=\rho\; m^{\nu,\,\mu}
$$ 
on $T^1M$.

\bigskip\noindent{\bf Step 2. } In this step, we construct two
families $(I_r)_{r>0}$ and $(J_r)_{r>0}$ of diffusion operators along
the strong unstable leaves, acting on the locally integrable functions
on $T^1M$, and prove an adjointness property for them which will be
crucial for the proof of Theorem
\ref{theo:uniergodtransv}.

\medskip Recall (see Subsection \ref{subsec:geodflowback}) that
$d_{W^{\rm su} (v)}$ is the Hamenst\"adt distance on the strong
unstable leaf $W^{\rm su}(v)$ of $v\in T^1\wt M$.  Recall that we
denote by $B^{\rm su}(v,r)$ the open ball of centre $v$ and radius $r$
for the distance $d_{W^{\rm su}(v)}$.  

Recall that $\prT:T^1\wt M\ra T^1M$ is the canonical projection. For
every $r>0$ and for all measurable nonnegative maps $\psi:T^1M\ra
\RR$, let $\wt \psi=\psi\circ \prT$ be the {\it lift}\index{lift} of
$\psi$ to $ T^1\wt M$.  Let us define $\wt I_r\wt \psi: T^1\wt M\ra
\RR$ as an integral of $\wt \psi$ on the balls of radius $r$ of the
strong unstable leaves:
\begin{equation}\label{eq:horosphericaldiffusion}
\forall \;u\in T^1\wt M,\;\;\; \wt I_r\wt \psi\,(u)=
\int_{v\in B^{\rm su}(u,\,r)}\wt\psi(v)\;e^{c_{\wt F}(v,\,u)} 
\;d\mu_{W^{\rm su}(u)}(v)\;.
\end{equation}
We define similarly
\begin{equation}\label{eq:horosphericaldiffusionbis}
\forall \;u\in T^1\wt M,\;\;\; \wt J_r\wt \psi\,(u)=
\int_{v\in B^{\rm su}(u,\,r)}\wt\psi(v)\;d\mu_{W^{\rm su}(u)}(v)\;.
\end{equation}
The maps $\wt I_r\wt \psi$ and $\wt J_r\wt \psi$ are $\Ga$-invariant,
by Equation \eqref{eq:propunmesfortinstabl} and the $\Ga$-invariance
of $c_{\wt F}$ and $\wt\psi$. Hence they define measurable nonnegative
maps $I_r\psi,J_r\psi: T^1M\ra \mathopen{[}0,+\infty\mathclose{]}$
whose lifts are $\wt I_r\wt \psi,\wt J_r\wt \psi$. Note that
$I_r\psi=J_r\psi$ if $F=0$.

Since $\wt{\psi}$ is nonnegative, then for every $u\in T^1\wt M$, the
maps $r\mapsto \wt I_r\wt \psi\,(u)$ and $r\mapsto \wt J_r\wt
\psi\,(u)$ are nondecreasing. The operators $\wt I_r$, $\wt J_r$,
$I_r$ and $J_r$ are positive: For instance, if $\psi\leq\psi'$ then
$I_r\psi\leq I_r\psi'$.

If $\wt{\psi}$ is continuous and positive at an element $u\in T^1\wt
M$ belonging to the support of $\mu_{W^{\rm su}(u)}$ (which is the set
of vectors $v\in W^{\rm su}(u)$ such that $v_+\in\Lambda\Ga$), then
$\wt I_r\wt \psi\,(u)>0$. In particular, we have $\wt I_r\wt
\rho\,(u)>0$ if $\wt\rho(u)>0$ and $u_+\in\Lambda\Ga$. Hence the
$\Ga$-invariant map $\frac{\wt{\rho}} {\wt I_r\wt\rho}$, defined to
have value $0$ at $u\in T^1\wt M$ if $\wt\rho(u)=0$ or
$u_+\notin\Lambda\Ga$, induces a measurable nonnegative map
$\frac{{\rho}}{{I_r\rho}}$ on $T^1M$.

Let us denote by $1$ the constant map with value $1$ on $T^1M$ and on
$T^1\wt M$, which satisfies $\wt I_r 1(u)>0$ if $u_+\in\Lambda\Ga$. 
%
We will apply the mixing property of Theorem \ref{theo:mixingbab} when
$B$ is a Hamenst\"adt ball $B^{\rm su}(u,r)$, giving that, for every $u\in
T^1M$ such that $u_-\in\Lambda\Ga$ and every nonnegative
$\psi\in\C_c(T^1M;\RR)$,
\begin{equation}\label{eq:conseqmixingforuse}
\lim_{t\ra+\infty} I_{r}(\psi\circ\phi_t)(u)= 
\frac{I_r1\,(u)}{\|m_F\|}\int_{T^1M}\psi \;dm_F\;.
\end{equation}

\medskip 
The diffusion operators $I_r$ on the strong unstable balls have the
following commutation property with the geodesic flow.

\blemm \label{lem:comutflowoperators}
For all $t\in\RR$, $r>0$, $u\in T^1 M$ and $\psi:T^1M\ra\RR$
measurable nonnegative, we have
$$
I_{re^t}(\psi\circ\phi_{-t})(\phi_tu)=
e^{\delta\,t-\int_0^t F(\phi_s u)\,ds}\;I_r\psi\,(u)\;.
$$
\elemm


\dem By the equations \eqref{eq:propdemesfortinstabl},
\eqref{eq:changemoinsplus} and \eqref{eq:cocycleprop}, for all
$t\in\RR$, $u\in T^1\wt M$ and $w\in W^{\rm su}(u)$, we have
\begin{align*}
d\,\mu_{W^{\rm su}(\phi_t u)}(\phi_t w)&=
e^{-C_{F-\delta,\,w_+}(\pi(\phi_t w),\,\pi(w))}\;
d\,\mu_{W^{\rm su}(u)}(w)\\ & =
e^{-C_{F\circ \iota-\delta,\,w_-}(\pi(w),\,\pi(\phi_t w))}\;
d\,\mu_{W^{\rm su}(u)}(w)\;.
\end{align*}
Hence for all $t\in\RR$, $r>0$, $u\in T^1\wt M$ and
$\wt\psi:T^1M\ra\RR$ measurable nonnegative, by the definition of the
operators $\wt I_r$ and by Equation \eqref{eq:dilatboulhamen}, using
the change of variables $v=\phi_tw$ and the definition of the cocycle
$c_{\wt F}$ in Equation \eqref{eq:deficsubwtF}, we have
\begin{align}
\wt I_{re^t}(\wt\psi\circ\phi_{-t})(\phi_t u)& = 
\int_{v\in B^{\rm su}(\phi_t u,\,re^t)}\wt\psi(\phi_{-t}v)\;
e^{c_{\wt F}(v,\,\phi_t u)} 
\;d\mu_{W^{\rm su}(\phi_t u)}(v)\nonumber\\ & =
\int_{w\in B^{\rm su}(u,\,r)}\wt\psi(w)\;e^{c_{\wt F}(\phi_t w,\,\phi_t u)} 
e^{-C_{F\circ \iota-\delta,\,w_-}(\pi(w),\,\pi(\phi_t w))}\;d\,\mu_{W^{\rm su}(u)}(w)
\nonumber
\\ & =e^{-C_{F\circ \iota-\delta,\,u_-}(\pi(u),\,\pi(\phi_t u))}\;\wt I_r\wt\psi(u)
=e^{-\int_0^t (\wt F(\phi_s u)-\delta)\,ds}\;\wt I_r\wt\psi(u)\;.
\label{eq:precisioncomutIflow}
\end{align}
The result follows. \cqfd

\medskip The operators $I_r$ and $J_r$ satisfy the following crucial
adjointness property for the scalar product of $\LL^2(m^{\nu,\,\mu})$.

\blemm\label{lem:adjoint} For all $r>0$ and $\psi,\varphi:T^1M\ra\RR$
measurable nonnegative, we have
$$
\int_{T^1M} I_r\psi\;\varphi\;dm^{\nu,\,\mu}=
\int_{T^1M} \psi\;J_r\varphi\;dm^{\nu,\,\mu}\;.
$$
\elemm

\dem We assume that $\Ga$ is torsion free: the general case follows by
proving, for every $k\in\NN-\{0\}$, an analogous statement for the
restriction of $m^{\nu,\,\mu}$ to the image in $T^1M$ of the
$\Ga$-invariant Borel subset of points in $T^1\wt M$ whose stabiliser
in $\Ga$ has order $k$, and by summation (see Subsection
\ref{subsec:pushmeas}).

By \cite[\S 1C]{Roblin03}, there exists hence a weak
fundamental domain in $T^1\wt M$ with $\wt m^{\nu,\,\mu}$-negligible
boundary, that is an open subset $\D$ in $T^1\wt M$ such
that $\ga\D\cap \D$ is empty for every $\ga\in \Ga- \{\id\}$ and
$\bigcup_{\ga\in\Ga}\ga \,\overline{\D}= T^1\wt M$, and such that
$\wt m^{\nu,\,\mu}(\partial \D)=0$. 

Let $T=W^{\rm s}(v)$ be a stable leaf in $T^1\wt M$. Recall that $v_+$
is not an atom of $\mu_x$, since $m_F$ is finite, hence $(\Ga,F)$ is
of divergence type (see Lemma \ref{prop:atomlessreferee}). In
particular, $U_T=\{w\in T^1\wt M\;:\;w_-\neq v_+\}$ has full measure
for $\wt m^{\nu,\,\mu}$. Recall that $\mathbbm{1}_A$ is the
characteristic function of a subset $A$. Using

$\bullet$~ Lemma \ref{lem:constructmesnumu} and Equation
\eqref{eq:horosphericaldiffusion} for the second equality,

$\bullet$~ the equality $W^{\rm su}(w)=W^{\rm su}(w')$ and the cocycle
property of $c_{\wt F}$ for the third one, 

$\bullet$~ an elementary argument using $\Ga$ and $\D$, the fact that
$w\in B^{\rm su}(u,r)$ if and only if $u\in B^{\rm su}(w,r)$ and Fubini's
theorem for positive functions for the fourth one, we have
\begin{align*}
  \int_{T^1M} I_r\psi\;\varphi\;dm^{\nu,\,\mu} =&\int_{T^1\wt M}\;
  \wt I_r\wt \psi\;\;\wt\varphi\;\mathbbm{1}_\D\;d\wt m^{\nu,\,\mu}\\ 
=& \int_{w'\in
    T}\;\int_{w\in W^{\rm su}(w')}\;\Big(\int_{u\in B^{\rm su}(w,\,r)}
  \wt\psi(u)\;e^{c_{\wt F}(u,\,w) }\;d\mu_{W^{\rm su}(w)}(u)\Big)\\
  &\;\;\;\;\;\;\;\;\;\;\;\;\;\;\;
  \;\;\;\;\;\;\;\;\;\;\wt\varphi(w)\;\mathbbm{1}_\D(w)\; e^{c_{\wt F}(w,\,
    w')} \;d\mu_{W^{\rm su}(w')}(w)\;d\nu_{T}(w')\\ =& \int_{w'\in
    T}\;\int_{w\in W^{\rm su}(w')}\;\int_{u\in W^{\rm su}(w')}
  \wt\varphi(w)\;\wt\psi(u)\;\mathbbm{1}_\D(w)\;\mathbbm{1}_{B^{\rm su}(w,\,r)}(u) \\
  &\;\;\;\;\;\;\;\;\;\;\;\;\;\;\; \;\;\;\;\;\;\;\;\;\; e^{c_{\wt
      F}(u,\,w')}\;d\mu_{W^{\rm su}(w')}(u)\;d\mu_{W^{\rm su}(w')}(w)\;d\nu_{T}(w')
  \\ =& \int_{w'\in T}\;\int_{u\in W^{\rm su}(w')}\;\int_{w\in W^{\rm su}(w')}
  \wt\varphi(w)\;\wt\psi(u)\;\mathbbm{1}_\D(u)\;\mathbbm{1}_{B^{\rm su}(u,\,r)}(w) \\
  &\;\;\;\;\;\;\;\;\;\;\;\;\;\;\;\;\;\;\;\;\;\;\;\;\; e^{c_{\wt
      F}(u,\,w')}\;d\mu_{W^{\rm su}(w')}(w)\;d\mu_{W^{\rm su}(w')}(u)\;d\nu_{T}(w')
  \\ =&\int_{T^1M} \psi\;J_r\varphi\;dm^{\nu,\,\mu}\;,
\end{align*}
as required.
\cqfd

\medskip From now on, we consider a nonnegative $\psi\in\C_c(T^1M;\RR)$.
For every $r>0$, applying the above lemma with $\varphi=
\frac{\rho}{I_r(\rho)}$, we have
\begin{equation}\label{eq:adjonction}
\int_{T^1M}
\frac{I_{r}\psi\,}{I_{r}\rho}\;\rho\;dm^{\nu,\,\mu}
=\int_{T^1M}  \psi \;J_r\big(\frac{\rho}{I_r\rho}\big)\;dm^{\nu,\,\mu}\;.
\end{equation}

\medskip In the next three steps, we will prove that, up to  positive
multiplicative constants, the right hand side of this equality is
bounded from above by $\int_{T^1M} \psi\;dm^{\nu,\,\mu}$, and that
its left hand side is bounded from below by $\int_{T^1M}
\psi\;dm_F$ if $r$ is large enough. Since this is valid for any
nonnegative $\psi\in\C_c(T^1M;\RR)$, this will imply that $m_F$ is
absolutely continuous with respect to $m^{\nu,\,\mu}$.

\bigskip\noindent{\bf Step 3. } In this step, let us prove that the
map $J_r(\frac{\rho}{I_r\rho})$ is bounded from above on a subset of
$T^1M$ of full measure with respect to $m^{\nu,\,\mu}$, which implies
that the right hand side of Equation \eqref{eq:adjonction} is bounded
from above by a constant times $\int_{w\in T^1M}
\psi\;dm^{\nu,\,\mu}$.

Recalling that $\wt M$ has pinched negative curvature, the following
result of polynomial growth of horospheres is due to Roblin
\cite[Prop.~6.3 (a)]{Roblin03} (see also \cite{Bowditch95}). Recall
that $\wt\Omega_c\Ga$ is the set of $v\in T^1\wt M$ such that
$v_-,v_+\in\Lambda_c\Ga$, and that $\Omega_c\Ga$ is its image in
$T^1M$.

\blemm [Roblin]\label{lem:polygrowthhoro} There exists $N\in\NN-\{0\}$
such that for all $w\in\wt\Omega_c\Ga$ and $r>0$, there exist
$w_1,\dots, w_N\in W^{\rm su}(w)$ such that
$$
B^{\rm su}(w,r)\;\cap\,\overline{\wt\Omega_c\Ga}\;\subset \;
\bigcup_{i=1}^N \;B^{\rm su}\Big(w_i,\frac{r}{2}\Big)\;.
$$
\elemm

\blemm \label{lem:majoJrhoparN}
For all $r>0$ and $w\in \Omega_c\Ga$, we have
$$
J_r\Big(\frac{\rho}{I_r\rho}\Big)(w)\leq N\;.
$$
\elemm

\dem Let us fix $r>0$ and $w\in \wt\Omega_c\Ga$, and let $w_1,\dots,
w_N\in T^1\wt M$ be as in Lemma \ref{lem:polygrowthhoro}. For all
$i\in\{1,\dots,N\}$ and $u\in B^{\rm su}(w_i,\frac{r}{2})$, since
$B^{\rm su}(u,r)$ then contains $B^{\rm su}(w_i, \frac{r}{2})$, we have
\begin{align*}
\wt I_r\wt \rho\,(u)&
=\int_{v\in B^{\rm su}(u,\,r)}\wt\rho\,(v)\;e^{c_{\wt F}(v,\,u)} 
\;d\mu_{W^{\rm su}(u)}(v)\\&\geq e^{c_{\wt F}(w,\,u)} 
\;\int_{v\in B^{\rm su}(w_i, \,\frac{r}{2})}\wt\rho\,(v)\;e^{c_{\wt F}(v,\,w)} 
\;d\mu_{W^{\rm su}(w)}(v)\;.
\end{align*}
Hence
\begin{align*}
&\int_{u\in B^{\rm su}(w_i, \,\frac{r}{2})}\frac{\wt\rho\,(u)}{\wt 
I_r\wt \rho\,(u)}\;d\mu_{W^{\rm su}(w)}(u)\\\leq \;&
\int_{u\in B^{\rm su}(w_i, \,\frac{r}{2})}\frac{\wt\rho\,(u)\;e^{c_{\wt F}(u,\,w)}}
{\int_{v\in B^{\rm su}(w_i, \,\frac{r}{2})}\wt\rho\,(v)\;e^{c_{\wt F}(v,\,w)} 
\;d\mu_{W^{\rm su}(w)}(v)}
\;d\mu_{W^{\rm su}(w)}(u)=1\;.
\end{align*}
Using Lemma \ref{lem:polygrowthhoro}, since the support of the measure
$\mu_{W^{\rm su}(w)}$ is contained in the closed set $\{u\in T^1\wt
M\;:\; u_-=w_-, \;u_+ \in \Lambda \Ga\}$, which is contained in
$\overline{\wt\Omega_c\Ga}$ by the density of $\Lambda_c \Ga$ in
$\Lambda \Ga$, the result follows.  
\cqfd

\medskip As required, since $\Omega_c\Ga$ has full
$m^{\nu,\,\mu}$-measure by Lemma \ref{lem:constructmesnumu}, this
lemma implies that
\begin{equation}\label{eq:majorat}
\int_{T^1M}  \psi \;J_r\big(\frac{\rho}{I_r\rho}\big)
\;dm^{\nu,\,\mu}\leq N\;\int_{T^1M}  \psi\;dm^{\nu,\,\mu}\;.
\end{equation}

\bigskip\noindent{\bf Step 4. } In this most technical step, after
giving the useful notation, we prove some continuity properties of the
diffusion operators $I_r\psi\,(u)$ on the strong unstable balls, as
functions of $u$. For this, we will adapt the technical lemmas of
\cite[Sect.~1H]{Roblin03} to the presence of a possibly unbounded
potential, and add some more lemmas using the H\"older regularity of
$F$ needed to control the Gibbs cocycles.

\medskip\noindent\begin{minipage}{9.4cm} ~~~ 
For all $u,v\in T^1\wt M$ such that $u_-\neq v_-$, recall
(see Subsection \ref{subsec:condmesgibbs} before Equation
\eqref{eq:proptrmesfortinstabl}) that the map $\theta_{v,\,u}=
\theta^{su}_{v,\,u}: \{w\in W^{\rm su}(u)\;:\;w_+\neq v_-\} \ra \{w\in
W^{\rm su}(v)\;:\;w_+\neq u_- \}$ is the homeomorphism sending $w$ to the
unique element in $W^{\rm su}(v)\cap W^{\rm s}(w)$. Note that
$\theta_{v',\,u'}=\theta_{v,\,u}$ if $u'\in W^{\rm su}(u)$ and $v'\in
W^{\rm su}(v)$. We have $\theta_{u,\,v}= {\theta_{v,\,u}}^{-1}$ and,
whenever defined, $\theta_{w,\,u}
=\theta_{w,\,v}\circ\theta_{v,\,u}$. For every $t\in\RR$, we have
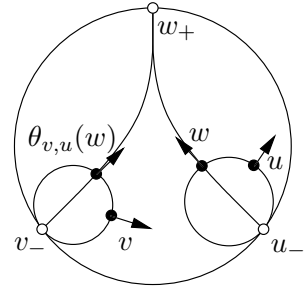
\begin{equation}\label{eq:comutthetaflow}
\phi_t\circ \theta_{v,\,u}=\theta_{\phi_tv,\,u}\;\;\;{\rm and}\;\;\;
\theta_{v,\,u}\circ \phi_t=\theta_{v,\,\phi_{-t}u}\;.
\end{equation}\end{minipage}
\begin{minipage}{5.5cm}
\begin{center}
\input{fig_familmesurebis.pstex_t}
\end{center}
\end{minipage}

\medskip 
We recall the definition of the canonical neighbourhoods (called {\it
  dynamical cells})\index{dynamical!cells} of the unit tangent
vectors, adapted to variable curvature. For all $u\in T^1\wt M$ and
$r_1,r_2,r_3>0$,
\begin{equation}\label{eq:cellulerie}
{\cal V}_{r_1,\,r_2,\,r_3}(u)=\bigcup_{|s|<r_3}\phi_s
\Big(\bigcup_{v\in B^{\rm ss}(u,\,r_1)} \theta_{v,\,u}(B^{\rm su}(u,r_2))\Big)\;.
\end{equation}
The sets ${\cal V}_{r_1,\,r_2,\,r_3}(u)$ are nondecreasing in
$r_1,r_2,r_3$ and form, as $r_1,r_2,r_3$ range over
$\mathopen{]}0,+\infty\mathclose{[}\,$, a neighbourhood basis of the
vector $u$ in $T^1\wt M$. Furthermore, for every $t\in\RR$, by the
equations \eqref{eq:comutthetaflow}, \eqref{eq:dilatboulhamen} and
\eqref{eq:contractboulhamen}, we have
\begin{equation}\label{eq:quasicommutflowcalV}
\phi_t\big({\cal V}_{r_1,\,r_2,\,r_3}(u)\big)=
{\cal V}_{e^{-t}r_1,\,e^{t}r_2,\,r_3}(\phi_tu)\;.
\end{equation}

\medskip For every $u\in T^1\wt M$, let
$$
u_\epsilon= \bigcup_{|s|<\epsilon} \phi_s B^{\rm ss}(u,\epsilon)\;.
$$
The sets $u_\epsilon$ are nondecreasing in $\epsilon$ and form, as
$\epsilon$ ranges over $\mathopen{]}0,1\mathclose{]}$, a neighbourhood
basis of the vector $u$ in its stable leaf $W^{\rm
  s}(u)$. Furthermore, for all $\ga\in\Ga$ and $t\geq 0$, we have
$\ga(u_\epsilon)=(\ga u)_\epsilon$ and, by Equation
\eqref{eq:contractboulhamen},
\begin{equation}\label{eq:relatflowusubeps}
\phi_t(u_\epsilon)\subset (\phi_tu)_\epsilon\;\;\;{\rm and}\;\;\;
(\phi_{-t}u)_\epsilon\subset \phi_{-t}(u_\epsilon)\;.
\end{equation}

The following result contains the elementary properties of the
notation introduced above. It says that if $v$ is in a small dynamical
cell around $u$, then the change of strong unstable leaf map
$\theta_{v,u}$ does not move points too much.
The important feature of this lemma is the uniformity of the
constants, which will be useful to control the variations of the
potential function and of the Gibbs cocycle.

\blemm\label{lem:voisdynlemrob} For every $\epsilon>0$, if
$r_1,r_2,r_3>0$ are small enough, then for every $u\in T^1\wt M$, we
have

(1) for all $v\in {\cal V}_{r_1,\,r_2,\,r_3}(u)$ and $w\in
B^{\rm su}(u,2)$, we have
$$
\theta_{v,u}(w)\in w_\epsilon\;,
$$

(2) for all $r\in \mathopen{[}\frac{1}{2},2\mathclose{]}$ and $v\in
{\cal V} _{r_1,\,r_2,\,r_3}(u)$, we have
$$
B^{\rm su}(v,r\,e^{-\epsilon})\subset \theta_{v,u}(B^{\rm su}(u,r))\subset 
B^{\rm su}(v,re^\epsilon)\;.
$$
%
\elemm

\dem (1) Let us recall the statement \cite[Lem.~1.13]{Roblin03}: for
every $\epsilon'>0$, there exists $r(\epsilon')>0$ such that for all
$u'\in T^1\wt M$, $v'\in B^{\rm ss}(u',r(\epsilon'))$ and $w'\in
B^{\rm su}(u',3)$ (the statement is given in loc.~cit.~with the number $2$
instead of this number $3$, but the proof is the same), we have
$$
\theta_{v',u'}(w')\in (w')_{\epsilon'}\;.
$$

Now, for every $\epsilon>0$, let $\epsilon'= \min\{r(\frac{\epsilon}
{3}), \frac{\epsilon}{3}\}>0$. Let
$r_1\in\mathopen{]}0,r(\epsilon')\mathclose{]}$, $r_2\in
\mathopen{]}0,1\mathclose{]}$ and
$r_3\in\mathopen{]}0,\frac{\epsilon}{3}\mathclose{]}$.

\medskip
\noindent\begin{minipage}{8.9cm}
  Let $u\in T^1\wt M$, $w\in B^{\rm su}(u,2)$ and $v\in {\cal V}_{r_1,
    \,r_2,\,r_3}(u)$.  By the definition of ${\cal
    V}_{r_1,\,r_2,\,r_3} (u)$, there exist
  $s_1\in\mathopen{]}-r_3,r_3\mathclose{[}\,$, $w_1\in B^{\rm ss}(u,r_1)$
  and $w_2\in B^{\rm su}(u, r_2)$ such that if
  $v_1=\theta_{w_1,\,u}(w_2)$, then $v= \phi_{s_1}v_1$.  By
  \cite[Lem.~1.13]{Roblin03} as recalled above, since $r_2\leq 3$ and
  $r_1\leq r(\epsilon')$, there exist $s_2\in \mathopen{]}-\epsilon',
  \epsilon'\mathclose{[}$ and $v_2\in B^{\rm ss}(w_2, \epsilon')$ such
  that $v_1=\phi_{s_2}v_2$, so that in particular
  $v=\phi_{s_1+s_2}v_2$.
\end{minipage}
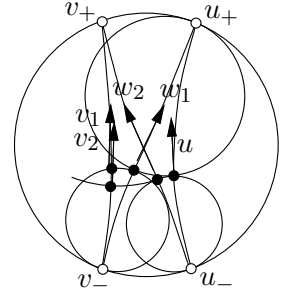
\begin{minipage}{6cm}
\begin{center}
\input{fig_unifconttheta.pstex_t}
\end{center}
\end{minipage}

\medskip By the triangle inequality, we have $w\in B^{\rm su}(u,2)\subset
B^{\rm su}(w_2, 2+r_2) \subset B^{\rm su}(w_2, 3)$. Since
$W^{\rm su}(w)=W^{\rm su}(u)= W^{\rm su}(w_2)$ and $v_2\in B^{\rm ss}(w_2,
r(\frac{\epsilon}{3}))$, applying again \cite[Lem.~1.13]{Roblin03} as
recalled above, there exist $s_3\in \mathopen{]}-\frac{\epsilon}{3},
\frac{\epsilon}{3}\mathclose{[}$ and $v_3\in B^{\rm ss}(w,
\frac{\epsilon}{3})$ such that $\theta_{v_2,u}(w)=\theta_{v_2,w_2}(w)=
\phi_{s_3}v_3$. Hence we have by Equation \eqref{eq:comutthetaflow}
that
$$
\theta_{v,\,u}(w)=\theta_{\phi_{s_1}v_1,\,u}(w)=
\theta_{\phi_{s_1+s_2}v_2,\,u}(w)=
\phi_{s_1+s_2}\theta_{v_2,\,u}(w)= 
\phi_{s_1+s_2+s_3}v_3\in  w_{\epsilon}\;.
$$

(2) By \cite[Coro.~1.14]{Roblin03}, for every $\epsilon>0$, there
exist $r'_1,r'_2,r'_3>0$ such that for all $u'\in T^1\wt M$ and $v'\in
{\cal V}_{r'_1,\,r'_2,\,r'_3} (u')$, we have
\begin{equation}\label{eq:Roblin114}
B^{\rm su}(v',e^{-\epsilon})\subset \theta_{v',u'}(B^{\rm su}(u',1))\subset 
B^{\rm su}(v',e^\epsilon)\;.
\end{equation}
Assume that $r_1\in\mathopen{]}0,\frac{r'_1}{2}\mathclose{]}$,
$r_2\in\mathopen{]}0,\frac{r'_2}{2}\mathclose{]}$ and $r_3\leq
r'_3$. For every $r\in\mathopen{[}\frac{1}{2},2\mathclose{]}$, let
$t=\ln r$. For every $v\in {\cal V}_{r_1, \,r_2,\,r_3}(u)$, let
$u'=\phi_{-t}u$ and $v'=\phi_{-t}v$. Then
$$
v'\in \phi_{-t}{\cal V}_{r_1,\,r_2,\,r_3}(u)={\cal V}_{e^tr_1,
\,e^{-t}r_2,\,r_3}(u')\subset {\cal V}_{r'_1,\,r'_2,\,r'_3} (u')\;,
$$
and the result follows by applying the map $\phi_t$ to the inclusions
\eqref{eq:Roblin114}, while using the equations
\eqref{eq:dilatboulhamen} and \eqref{eq:comutthetaflow}.
%
%
\cqfd

\medskip We will need the following consequence of the H\"older
regularity of $\wt F$ (which is empty if $\wt F=0$).

\blemm \label{lem:holderargumagain} For every $\epsilon>0$, if
$\epsilon'>0$ is small enough, then for all $w\in T^1\wt M$ and
$w'\in w_{\epsilon'}$, we have
$$
|\,C_{F-\delta,\,w_+}(\pi(w),\pi(w'))\,|\leq 
\epsilon\,\Big(1+\sup_{\pi^{-1}(B(\pi(w),\,1))}|\wt F|\Big)\;.
$$
\elemm

\dem 
First observe that 
\begin{equation}\label{eq:pointbaseproche}
\forall\;v\in u_\epsilon,\;\;\; d(\pi(v),\pi(u))< 2\epsilon\;.
\end{equation}
Indeed, if $v=\phi_sw$ where $|s|<\epsilon$ and $w\in
B^{\rm ss}(u,\epsilon)$, then $d(\pi(v),\pi(w))=|s|<\epsilon$ and by
Equation \eqref{eq:compardistHamen} (exchanging the stable and unstable
foliations), we have
$$
d(\pi(u),\pi(w))\leq d_{W^{\rm ss} (u)}(u,w)<\epsilon\;,
$$
so that Equation \eqref{eq:pointbaseproche} holds by the triangle
inequality.  

Now by Equation \eqref{eq:pointbaseproche} and  by Lemma
\ref{lem:holderconseq}, there exist two constants 
$c_1,c_2$ such that for every $w'\in w_{\epsilon'}$,
$$
|\,C_{F-\delta,\,w_+}(\pi(w),\pi(w'))\,|\leq c_1(2\epsilon')^{c_2}+
(2\epsilon')\Big(\max_{\pi^{-1}(B(\pi(w),\,2\epsilon'))}
|\wt F|+\delta\Big)
\;,
$$
which proves the result. 
\cqfd

\medskip A second consequence of the H\"older regularity of $\wt F$ that
we will need is the following result (again empty if $F=0$).

\blemm \label{lem:holderargumlast} For all $\epsilon>0$ and $T\geq
0$, if $r_1,\,r_2,\,r_3>0$ are small enough, then for all $u\in T^1\wt
M$ and $v\in {\cal V}_{r_1,\,r_2,\,r_3} (u)$, we have
$$
\Big|\,\int_0^T\wt F(\phi_tu)-\wt F(\phi_tv)\;dt\,\Big|
\leq \epsilon\;.
$$
\elemm

\dem This follows from the uniform continuity of the geodesic flow on
compact interval of times and the fact that the sets ${\cal
  V}_{r_1,\,r_2,\,r_3} (u)$ for $r_1,\,r_2,\,r_3>0$ form a
neighbourhood basis of $u\in T^1\wt M$. 
\cqfd
%
%

\medskip The following lemma (which is still empty if $F=0$),
giving a uniform continuity property of the cocycle $c_{\wt F}$, will
be useful to prove the continuity properties of $u\mapsto \wt I_r\wt
\varphi(u)$.

\blemm \label{lem:unifcontgibbs} For every $\epsilon>0$, if
$\epsilon',r'_1,r'_2,r'_3>0$ are small enough, then for all $u\in
T^1\wt M$, $w\in B^{\rm su}(u,2)$, $v\in {\cal V}_{r'_1,\,r'_2,\,r'_3} (u)$
and $w'\in w_{\epsilon'}\cap W^{\rm su}(v)$, we have
$$
|\,c_{\wt F}(w',\,v)-c_{\wt F} (w,\,u)\,|\leq \epsilon\;.
$$
\elemm

A less precise version of this result follows quite directly from the
H\"older-continuity of $c_{\wt F}$. The proof shows that this result
also holds if $w'\in {\cal V} _{r'_1,\,r'_2,\,r'_3} (u)\cap
W^{\rm su}(v)$. But we will use the explicit form given above (with
$\frac{\epsilon}{2}$ instead of $\epsilon$!).

\medskip \dem Let us fix $\epsilon>0$. Let $c>0$ and
$\alpha\in\mathopen{]}0,1\mathclose{]}$ be as in the beginning of the
proof of Lemma \ref{lem:holderargumagain}.  By Equation
\eqref{eq:compardistHamen} and Equation \eqref{eq:dilatHamdist}, there
exists a constant $c'>0$ such that, for all $t\geq 0$, $u'\in T^1\wt
M$ and $u''\in B^{\rm su}(u',3)$, we have
$$
d(\phi_{-t}u',\phi_{-t}u'')\leq 
c'\,d_{W^{\rm su}(\phi_{-t} u')}(\phi_{-t} u',\phi_{-t} u'')
\leq c'\,e^{-t}\,d_{W^{\rm su}(u')}(u',u'')\leq 3\,c'\,e^{-t}\;.
$$
Let $T=\max\{0,-\ln(\frac{1}{3c'}
\sqrt[\alpha]{\frac{\epsilon\alpha}{4c}}\,)\,\}$.  By the H\"older
regularity of $\wt F$, for every $w\in B^{\rm su}(u,2)\subset
B^{\rm su}(u,3)$, we have
$$
\Big|\,\int_T^{+\infty} \wt F(\phi_{-t}w)-\wt F(\phi_{-t}u)\;dt\,\Big|
\leq \int_T^{+\infty} c\,(3\,c'\,e^{-t})^\alpha=
\frac{c}{\alpha}\,(3\,c'\,e^{-T})^\alpha\leq \frac{\epsilon}{4}\;.
$$
If $\epsilon',r'_1,r'_2,r'_3$ are small enough, then for all $v\in
{\cal V}_{r'_1,\,r'_2,\,r'_3} (u)$ and $w'\in w_{\epsilon'}\cap W^{\rm
  su} (v)$, we have $w'\in B^{\rm su}(v,3)$ since $w\in B^{\rm su}
(u,2)$, hence similarly
$$
\Big|\,\int_T^{+\infty} \wt F(\phi_{-t}w')-\wt
F(\phi_{-t}v)\;dt\,\Big| \leq \frac{\epsilon}{4}\;.
$$

Since $w'\in w_{\epsilon'}$, there exist $s'\in\mathopen{]}-\epsilon',
\epsilon'\mathclose{[}$ and $w''\in B^{\rm ss}(w,\epsilon')$ such that
$$
w'=\phi_{s'}w''=\phi_{s'}\,\theta_{w'',\,w}(w)\;,
$$
this last equality holding since $w''\in W^{\rm s}(w)$. 
Since $w\in B^{\rm su}(w,\epsilon')$, we hence have
$$
w'\in {\cal V}_{\epsilon',\,\epsilon',\,\epsilon'}(w)\;.
$$

Let $r_1,r_2,r_3>0$ be the constants associated by Lemma
\ref{lem:holderargumlast} to $\frac{\epsilon}{4}$ and $T$. If
$\epsilon',r'_1,\,r'_2$, $r'_3>0$ are small enough, then
$\phi_{-T}v\in {\cal V} _{r_1,\,r_2,\,r_3} (\phi_{-T}u)$ and
$\phi_{-T}w'\in {\cal V} _{r_1,\,r_2,\,r_3} (\phi_{-T}w)$ by Equation
\eqref{eq:quasicommutflowcalV}.  Hence by Lemma
\ref{lem:holderargumlast}, using a change of variable $t'=T-t$, we
have
$$
\Big|\,\int_0^T\wt F(\phi_{-t}u)-\wt F(\phi_{-t}v)\;dt\,\Big|
=\Big|\,\int_0^T\wt F(\phi_{t'}(\phi_{-T}u))-
\wt F(\phi_{t'}(\phi_{-T}v))\;dt'\,\Big|
\leq \frac{\epsilon}{4}
$$
and similarly
$$
\Big|\,\int_0^T\wt F(\phi_{-t}w)-\wt F(\phi_{-t}w')\;dt\,\Big|
\leq \frac{\epsilon}{4}\;.
$$
Now, by the definition of the cocycle $c_{\wt F}$, we have
\begin{align*}
&|\,c_{\wt F}(w',\,v)-c_{\wt F} (w,\,u)\,|\\ =\; &
\Big|\int_0^{+\infty}\wt F(\phi_{-t}w')-\wt F(\phi_{-t}v)-
\wt F(\phi_{-t}w)+\wt F(\phi_{-t}u)\;dt\Big|
\\ \leq\;& \Big|\,\int_T^{+\infty} \wt F(\phi_{-t}w')-\wt
F(\phi_{-t}v)\;dt\,\Big| +
\Big|\,\int_T^{+\infty} \wt F(\phi_{-t}w)-\wt F(\phi_{-t}u)\;dt\,\Big|
\\ &+\Big|\,\int_0^T\wt F(\phi_{-t}u)-\wt F(\phi_{-t}v)\;dt\,\Big|
+\Big|\,\int_0^T\wt F(\phi_{-t}w)-\wt F(\phi_{-t}w')\;dt\,\Big|\\
\leq \;& \frac{\epsilon}{4}+\frac{\epsilon}{4}+\frac{\epsilon}{4}+
\frac{\epsilon}{4}=\epsilon\;.
\end{align*}
This proves the result.
\cqfd

\medskip Let us now introduce the functional version of the notation
$u_\epsilon$. For every $\epsilon>0$ and every nonnegative measurable
map $\varphi:T^1M\ra \RR$, whose lift to $T^1\wt M$ is denoted by $\wt
\varphi$, let us define
$$ 
\wt \varphi_{\,\overline\epsilon}:u\mapsto 
\sup_{v\in u_\epsilon}\wt \varphi(v)
\;\;\;{\rm and}\;\;\; \wt\varphi_ {\,\underline{\epsilon}}
:u\mapsto \inf_{v\in u_\epsilon}\wt \varphi(v)\;.
$$
Note that $\wt \varphi_{\,\overline{\epsilon'}}\leq\wt \varphi_{\,
  \overline \epsilon}$ and $\wt \varphi_{\, \underline{\epsilon'}}
\geq\wt \varphi_{\,\underline{\epsilon}}$ if $\epsilon'\leq
\epsilon$. The maps $\wt \varphi_{\,\overline\epsilon}$ and $\wt
\varphi_ {\,\underline{\epsilon}}$ are $\Ga$-invariant, hence define
nonnegative measurable maps $\varphi_{\,\overline\epsilon}$ and
$\varphi_ {\,\underline{\epsilon}}$ from $T^1M$ to $\RR$.  For every
$t\geq 0$, by Equation \eqref{eq:relatflowusubeps}, we have
\begin{equation}\label{eq:commutflowepsilondiff}
(\varphi\circ \phi_t)_{\,\overline\epsilon}\leq 
\varphi_{\,\overline\epsilon}\circ \phi_t \;\;\;{\rm and}\;\;\;
(\varphi\circ \phi_t)_{\,\underline{\epsilon}}\geq
\varphi_{\,\underline{\epsilon}}\circ\phi_t\;.
\end{equation}

\medskip
The following result is the regularity property in $u$ of
$\wt I_r\wt\varphi\,(u)$ that will be needed in the next steps.

\bprop\label{prop:continuitydiffusion} Let $K$ be a compact subset of
$T^1\wt M$. For all $\epsilon>0$ and $r\geq 1$, if $r_1,r_2,r_3>0$ are
small enough, then for all $u\in \Ga K$ and $v\in {\cal
  V}_{r_1,\,r_2,\,r_3}(u)$ and for every nonnegative measurable map
$\varphi:T^1M\ra \RR$, we have
\begin{equation}\label{eq:continuitydiffusion}
e^{-\epsilon}\wt I_{re^{-\epsilon}}(\wt\varphi_{\,\underline\epsilon})\;(u)
\leq \wt I_r\wt\varphi\;(v) \leq 
e^{\epsilon}\wt I_{re^{\epsilon}}(\wt\varphi_{\,\overline{\epsilon}})\;(u)\;.
\end{equation}
\eprop

\dem Let us first prove that we may assume without loss of generality
that $r=1$.

Assume that the result is true for $r=1$. Let $\epsilon>0$, $r\geq 1$
and $t=-\ln r\leq 0$. For all $r_1,r_2,r_3>0$, $u\in T^1\wt M$ and
$v\in {\cal V}_{e^tr_1,\,e^{-t}r_2,\,r_3}(u)$, we have $\phi_tv\in
{\cal V}_{r_1,\,r_2,\,r_3}(\phi_tu)$ by Equation
\eqref{eq:quasicommutflowcalV}. Hence by the case $r=1$ of this
proposition applied to $\phi_t K$, $\frac{\epsilon}{2}$ and
$\wt\varphi\circ \phi_{-t}$, we have, if $r_1,r_2,r_3$ are small
enough,
$$
e^{-\frac{\epsilon}{2}}\;\wt I_{e^{-\frac{\epsilon}{2}}}((\wt
\varphi\circ \phi_{-t})_{\,\underline{\frac{\epsilon}{2}}})\;
(\phi_{t}u)
\leq \wt I_1(\wt\varphi\circ \phi_{-t})\;(\phi_{t}v) \leq 
e^{\frac{\epsilon}{2}}\;\wt I_{e^{\frac{\epsilon}{2}}}(
(\wt\varphi\circ \phi_{-t})_{\,\overline{\frac{\epsilon}{2}}})\;
(\phi_{t}u)\;.
$$
By the monotonicity properties of the maps $s\mapsto \wt
I_s\wt{\varphi'}(u')$, $\epsilon'\mapsto
\wt{\varphi'}_{\,\overline{\epsilon'}}(u')$ and $\epsilon'\mapsto
\wt{\varphi'}_{\,\underline{\epsilon'}}(u')$ for every $u'\in T^1\wt M$ when
${\varphi'}:T^1\wt M\ra\RR$ is nonnegative, and by Equation
\eqref{eq:commutflowepsilondiff} since $-t\geq 0$, we have
$$
e^{-\frac{\epsilon}{2}}\;\wt I_{e^{-\epsilon}}(\wt
\varphi_{\,\underline{\epsilon}}\circ \phi_{-t})\;(\phi_{t}u)
\leq \wt I_1(\wt\varphi\circ \phi_{-t})\;(\phi_{t}v) \leq 
e^{\frac{\epsilon}{2}}\;\wt I_{e^{\epsilon}}(
\wt\varphi_{\,\overline{\epsilon}}\circ \phi_{-t})\;(\phi_{t}u)\;.
$$
By Lemma \ref{lem:comutflowoperators}, we have
\begin{align*}
e^{-\frac{\epsilon}{2}}\;e^{\delta t-\int_0^t\wt F(\phi_s u)\;ds}\;
&\wt I_{r\,e^{-\epsilon}}(\wt
\varphi_{\,\underline{\epsilon}})\;(u)
\\ & \leq
e^{\delta t-\int_0^t\wt F(\phi_s v)\;ds}\;
\wt I_r(\wt\varphi)\;(v)\leq
e^{\frac{\epsilon}{2}}\;e^{\delta t-\int_0^t\wt F(\phi_s u)\;ds}\;
\wt I_{r\,e^{\epsilon}}(
\wt\varphi_{\,\overline{\epsilon}})\;(u)\;.
\end{align*}
Dividing by $e^{\delta t}$ and applying Lemma
\ref{lem:holderargumlast} (with $\frac{\epsilon}{2}$ instead of
$\epsilon$), the general case of Proposition
\ref{prop:continuitydiffusion} follows.

\medskip Let us now prove Proposition \ref{prop:continuitydiffusion}
with $r=1$. Let $\epsilon,r_1,r_2,r_3>0$, $u\in T^1\wt M$, $v\in {\cal
  V}_{r_1,\,r_2,\,r_3}(u)$ and let $\varphi:T^1M\ra\RR$ be a
nonnegative measurable map. Since $e^{\epsilon}\wt I_{e^{\epsilon}}
(\wt\varphi_{\,\overline\epsilon})\;(u)$ is nondecreasing in
$\epsilon$ and since $e^{-\epsilon}\wt I_{e^{-\epsilon}}
(\wt\varphi_{\,\underline\epsilon})\;(u)$ is nonincreasing in
$\epsilon$, we may assume that $\epsilon\leq \ln 2$.

Let $\epsilon'\in\mathopen{]}0,\epsilon\mathclose{[}$ and $w\in
B^{\rm su}(u,e^{\pm\epsilon})$.  By Lemma \ref{lem:voisdynlemrob} (1), if
$r_1,r_2,r_3>0$ are small enough (depending only on $\epsilon'$), we
have
\begin{equation}\label{eq:thetabougepeu}
\theta_{v,\,u}(w)\in  w_{\epsilon'}\subset w_{\epsilon}\;.
\end{equation}

By Equation \eqref{eq:proptrmesfortinstabl}, using the
change of variables $w'=\theta_{v,\,u}(w)$, we have
\begin{align}
&\int_{w'\in\,\theta_{v,\,u}(B^{\rm su}(u,\,e^{\pm \epsilon}))}
\wt\varphi(w')\;e^{c_{\wt F}(w',\,v)}\;d\mu_{W^{\rm su}(v)}(w')
\nonumber\\= &\int_{w\in
B^{\rm su}(u,\,e^{\pm\epsilon})}\wt\varphi(\theta_{v,u}(w))
\;e^{c_{\wt F}(\theta_{v,u}(w),\,v)-c_{\wt F} (w,\,u)
-C_{F-\delta,\,w_+}(\pi(\theta_{v,u}(w)),\,\pi(w))}\nonumber
\\&
\;\;\;\;\;\;\;\;\;\;\;\;\;\;\;\;\;\;\;\;\;\;\;\;\;\;\;\;\;\;\;\;
\;\;\;\;\;\;\;\;\;\;\;\;\;\;\;\;\;\;\;\;\;\;\;\;\;\;\;\;\;\;\;\;
\;\;\;\;\;\;\;\;\;\;e^{c_{\wt F}(w,\,u)}\;d\mu_{W^{\rm su}(u)}(w)\;.
\label{eq:changvariabtheta}
\end{align}

By Lemma \ref{lem:voisdynlemrob} (2) with respectively $r=e^{-\epsilon}$
and $r=e^\epsilon$, we have 
$$
\theta_{v,\,u}(B^{\rm su}(u,e^{-\epsilon}))
\subset B^{\rm su}(v,1)\subset\theta_{v,\,u}(
B^{\rm su}(u,e^\epsilon))\;.
$$

Note that $\wt F$ is bounded on $\bigcup_{w\in\Ga K}
\pi^{-1}(B(\pi(w),1))$ if $K$ is compact, since $\wt F$ is
$\Ga$-invariant and continuous. By Equation \eqref{eq:thetabougepeu},
by Lemma \ref{lem:unifcontgibbs} and by Lemma
\ref{lem:holderargumagain}, the result follows from Equation
\eqref{eq:changvariabtheta}.  \cqfd

\bigskip\noindent{\bf Step 5. } In this step, in order to work in the
direction of giving a lower bound of the left hand side $\int_{T^1M}
\frac{I_{r}\psi\,}{I_{r}\rho}\;\rho\;dm^{\nu,\,\mu}$ of Equation
\eqref{eq:adjonction} by $\int_{T^1M} \psi\;dm_F$ if $r$ is large
enough, up to a positive multiplicative constant, we give a lower
bound of $\frac{I_{r}\psi\,}{I_{r}\rho}(u)$ for some $r$ depending on
$u$ in a full measure subset for $\Pi=\rho\;dm^{\nu,\,\mu}$. This
lower bound will be improved in the next step.

\medskip We fix a point $x_0$ in $\wt M$. For every
$R>0$, let $K_R$ be the closure of $\pi^{-1}(B(x_0,R))\cap
\wt\Omega_c\Ga$, which is a compact subset of $T^1\wt M$. Let
$\wt\Omega_R$ be the $\Ga$-invariant set of elements $v\in T^1\wt M$
such that there exists a sequence $(t_i)_{i\in\NN}$ in $\mathopen{[}
0,+\infty\mathclose{[}$ converging to $+\infty$ such that
$\phi_{-t_i}v\in \Ga K_R$ for every $i\in\NN$. Let $\Omega_R$ be the
image of $\wt\Omega_R$ in $T^1M$. Since the set of elements of $T^1\wt
M$ whose image in $T^1M$ is negatively recurrent under the geodesic
flow is the increasing union of the family $(\Omega_R)_{R>0}$, and
since $\nu$ gives full measure to the negatively recurrent vectors, we
fix $R>0$ large enough so that
$$
\Pi(\Omega_R)>1-\frac{1}{16N}\;,
$$
where $\Pi=\rho\;dm^{\nu,\,\mu}$ has been defined at the end of Step
1, and $N$ is given by Lemma \ref{lem:polygrowthhoro}.

\blemm\label{lem:step5} There exist $c_1,c_2>0$ and a $\Ga$-invariant
Borel map $r:\wt\Omega_R\ra\mathopen{[}0,+\infty\mathclose{[}$
such that for every nonnegative $\psi'\in\C_c(T^1M;\RR)$ and for 
every $u'\in\wt\Omega_R$,
$$
\wt I_{r(u')}\wt{\psi'}\,(u')\geq 
c_1\,\Big(\int_{T^1M}\psi'\;dm_F\;\Big)\;\wt I_{r(u')}\wt\rho\,(u')
$$
and 
$$
0<\wt I_{3\,r(u')}\wt\rho\,(u')\leq c_2\;\wt I_{r(u')}\wt\rho\,(u')\;.
$$
\elemm

\dem Let $c=\frac{1}{8}\,\inf_{v\in K_R}\;\wt I_{1}1\,(v)$ and
$c'=8\,\sup_{v\in K_R}\;\wt I_{12}1\,(v)$, which satisfy $c\leq
c'$. Since $K_R$ is contained in the support of $\wt\Pi= \wt\rho\,
\wt m^{\nu,\mu}$, for every $u\in K_R$, we have $\wt I_{\frac{1}{2}}1
\,(u)>0$ and $\wt I_{24}1\,(u)<+\infty$. By Proposition
\ref{prop:continuitydiffusion} applied with $K=K_R$, $\epsilon=\log
2$, $r=1$ or $r=12$, and $\wt\varphi=1$, the maps $v\mapsto \wt
I_{1}1\,(v)$ and $v\mapsto \wt I_{12}1\,(v)$ are hence locally with
positive lower bound and locally with finite upper bound on
$K_R$. Since $K_R$ is compact, we have $c>0$ and $c'<+\infty$. Define
$$
c_1=\frac{c}{c'\,\int_{T^1M}\rho\;dm_F}\;\;\;{\rm and}\;\;\;
c_2=\frac{c'}{c}\;.
$$

Let $\varphi = \rho$ (as fixed in the end of Step 1, in particular
$\int_{T^1M}\rho\;dm_F>0$) or $\varphi=\psi'$ (as in the statement). By
uniform continuity, let us fix $\epsilon\in\mathopen{]}0,\ln
2\mathclose{[}$ small enough so that, with the notation introduced
before Proposition \ref{prop:continuitydiffusion}, we have
\begin{equation}\label{eq:minoavecepsilon}
\int_{T^1M}\varphi_{\,\underline{\epsilon}}\;dm_F\geq 
\frac{1}{2}\int_{T^1M}\varphi\;dm_F\;\;\;{\rm and}\;\;\;
\int_{T^1M}\rho_{\,\overline{\epsilon}}\;dm_F\leq 
2\int_{T^1M}\rho\;dm_F\;.
\end{equation}
Proposition \ref{prop:continuitydiffusion}, applied with this
$\epsilon$, with $r=1$ or $r=3$, and with $K=K_R$, gives us
$r_1,r_2,r_3>0$ such that its conclusion holds.  

By compactness, let $S$ be a finite subset of $K_R$ such that $K_R$ is
covered by the dynamical cells ${\cal V}_{r_1,\,r_2,\,r_3}(u)$ as $u$
ranges over $S$. 

In the following sequence of inequalities, we use respectively

$\bullet$~ Proposition \ref{prop:continuitydiffusion} applied to the
function $\wt\varphi\circ \phi_t$;

$\bullet$~ Equation \eqref{eq:commutflowepsilondiff}
and monotonicity properties; 

$\bullet$~ Equation \eqref{eq:conseqmixingforuse} (which is the
consequence of the mixing property we are using) with $r=1$, since $S$
is finite and since for all $u\in\Ga S$, we have $u_-\in \Lambda\Ga$
as $\Ga S\subset\Ga K_R\subset \wt\Omega\Ga$;

$\bullet$~and Equation \eqref{eq:minoavecepsilon}:

\noindent There exists $t_0\geq 0$ such that for every $t\geq t_0$,
for all $u\in\Ga S$ and $v\in {\cal V}_{r_1,\,r_2,\,r_3}(u)$,
\begin{align}
\wt I_{2}(\wt\varphi\circ \phi_t)\,(v) & 
\geq e^{-\epsilon}\;\wt I_{2e^{-\epsilon}}
(\wt\varphi\circ \phi_t)_{\,\underline{\epsilon}}\,(u)\geq 
\frac{1}{2}\;\wt I_{1}
(\wt\varphi_{\,\underline{\epsilon}}\circ \phi_t)\,(u)\nonumber
\\ & \geq \frac{\wt I_{1}1\,(u)}{4\,\|m_F\|}\;
\int_{T^1M}\varphi_{\,\underline{\epsilon}}\;dm_F \geq 
\frac{c}{\|m_F\|}\;\int_{T^1M}\varphi\;dm_F\;.\label{eq:minoIwtde}
\end{align}
In particular, $\wt I_{2}(\wt\rho\circ \phi_t)\,(v)>0$. Similarly, an
upper bound is given by
\begin{align}
\wt I_{6}(\wt\rho\circ \phi_t)\,(v) 
&\leq e^{\epsilon}\;\wt I_{6e^{\epsilon}}
(\wt\rho\circ \phi_t)_{\,\overline{\epsilon}}\,(u) \leq 
2\;\wt I_{12}(\wt\rho_{\,\overline{\epsilon}}\circ \phi_t)\,(u)\nonumber
\\ & \leq \frac{4\,\wt I_{12}1\,(u)}{\|m_F\|}\;
\int_{T^1M}\rho_{\,\overline{\epsilon}}\;dm_F \leq 
\frac{c'}{\|m_F\|}\;\int_{T^1M}\rho\;dm_F\;.\label{eq:majoIwtsix}
\end{align}
Since $\wt I_{2}(\wt\rho\circ \phi_t)\,(v) \leq \wt I_{6}(\wt\rho\circ
\phi_t)\,(v) \leq \frac{c'}{\|m_F\|}\;\int_{T^1M}\rho\;dm_F$, we have,
by Equation \eqref{eq:minoIwtde} with $\varphi=\psi'$,
\begin{align*}
\wt I_{2}(\wt{\psi'}\circ \phi_t)\,(v) &
\geq\frac{c}{c'\,\int_{T^1M}\rho\;dm_F}\;\int_{T^1M}\psi'\;dm_F\;\;
\wt I_{2}(\wt\rho\circ \phi_t)\,(v)\\ &=
c_1\,\int_{T^1M}\psi'\;dm_F\;\;\wt I_{2}(\wt\rho\circ \phi_t)(v)\;.
\end{align*}

By Lemma \ref{lem:comutflowoperators} and a cancellation argument, we
hence have
$$
\wt I_{2e^t}\wt\psi\,(\phi_t v) 
\geq
c_1\,\int_{T^1M}\psi\;dm_F\;\;\wt I_{2e^t}\wt \rho\,(\phi_t v)\;.
$$
Similarly, since $0<\wt I_{6}(\wt\rho\circ \phi_t)\,(v) \leq
\frac{c'}{c}\;\wt I_{2}(\wt\rho\circ \phi_t)\,(v)$ by Equation
\eqref{eq:minoIwtde} and Equation \eqref{eq:majoIwtsix}, we have
$$
0<\wt I_{6\,e^t}\wt\rho\,(\phi_t v) \leq 
c_2\,\wt I_{2e^t}\wt \rho\,(\phi_t v)\;.
$$
Since for every $u'\in \wt \Omega_R$ there exists $t(u')\geq t_0\geq
0$ (which may be taken to be constant on orbits of $\Ga$ and
measurable) such that $u'\in\phi_{t(u')}\Ga K_R$, this concludes the
proof of Lemma \ref{lem:step5}, with $r(u')=2\,e^{t(u')}$ (note that
the map $r$ depends only on $R$).  \cqfd

\bigskip\noindent{\bf Step 6. } In this step, we start to improve the
previous one by giving a lower bound of $\frac{I_{r}\psi\,}
{I_{r}\rho}(u)$ for $u$ in some good subset $\wt G_r$ of $T^1\wt
M$, that will be shown in the next step to be large enough for our
purpose, if $r$ is large enough.

For every $r>0$, we first define the good set $\wt G_r$ (depending on
the nonnegative $\psi\in\C_c(T^1M;\RR)$ that has been fixed at the end
of Step 2) on which the lower bound will take place.

For every $s'>0$, let $\wt \E_{s'}$ be the $\Ga$-invariant Borel set
of $u\in\wt\Omega_R$ such that $r(u)\leq s'$, where the map $r$ is
given by Lemma \ref{lem:step5} in the previous step. The set
$\wt\Omega_R$ is the increasing union of the family $(\wt
\E_{s'})_{s'>0}$. Note that $\Pi(\Omega_R)> 1- \frac{1}{16N}$ by the
choice of $R$ in the previous step. Hence there exist $\sigma>0$ and a
$\Ga$-invariant closed subset $\wt \E$ of $\wt \E_{\sigma}$ whose
image $\E$ in $T^1M$ satisfies
$$
\Pi(\E)>1-\frac{1}{8N}\;.
$$
For every $r>0$, let $\wt \Delta_r$ be the $\Ga$-invariant Borel set
of $u\in \wt\Omega_R$ such that
\begin{equation}\label{eq:defiDeltasubr}
\wt I_r(\wt\rho\,\mathbbm{1}_{\wt \E})\,(u)\leq 
\frac{1}{2}\;\wt I_r\wt\rho\,(u)\;.
\end{equation}
For every $r>0$, let $\wt Z_r$ be the $\Ga$-invariant Borel set of
$u\in \wt\Omega_R$ such that 
\begin{equation}\label{eq:defigranZsubr}
\wt I_{r+\sigma}\wt\psi\,(u)>2\;\wt I_r\wt\psi\,(u)\;.
\end{equation}
Define 
$$
\wt G_r=\wt\Omega_c\Ga-(\wt \Delta_r\cup \wt Z_r)\;,
$$ 
which is a $\Ga$-invariant Borel subset of $T^1\wt M$ (depending on
$\psi$). Let $\Delta_r,Z_r,G_r$ be the images of $\wt \Delta_r,\wt
Z_r,\wt G_r$ in $T^1M$.

\medskip The next result now gives the important property of the set
$\wt G_r$ introduced above, saying that it is indeed a good set
concerning the problem of finding a lower bound of $\frac{\wt I_{r}\wt
  \psi\,} {\wt I_{r}\wt \rho}(u)$ by a constant times
$\int_{T^1M}\psi\;dm_F$.

\blemm\label{lem:step6} There exists $c_3>0$ such that for every
nonnegative $\psi\in\C_c(T^1M;\RR)$, for all $r>0$ and $u\in \wt G_r$,
we have
$$
\wt I_r\wt\psi\,(u)\geq 
c_3\;\int_{T^1M}\psi\;dm_F\;\;\wt I_r\wt\rho\,(u)
$$
\elemm

\dem Let $r>0$ and $u\in \widetilde{G}_r$. 

As pioneered in \cite{Rudolph83} to study the ergodic theory of the
strong unstable foliation, and as in \cite[page 89]{Roblin03}, we
start the proof by a Vitali covering argument. Consider the measure
$\mathbb{M}$ on $W^{\rm su}(u)$ defined by
$$
d\MM(w)=\wt\rho(w)\;e^{c_{\wt F}(w,\,u)}\;
d\mu_{W^{\rm su}(u)}(w)\;.
$$
For all $v\in B^{\rm su}(u,r)\cap\wt\E$, we have $v\in \wt\Omega_R$
and, by the second property of Lemma \ref{lem:step5},
$$
\MM\big(B^{\rm su}(v,2\,r(v))\big)\le 
c_2\; \MM\big(B^{\rm su}(v,r(v))\big)\,.
$$

Let us construct a finite or countable sequence of pairwise disjoint
balls $(B_i)_{1\leq i<i_*}$ where $i_*\in(\NN-\{0,1\})\cup\{+\infty\}$
amongst the open balls $B^{\rm su}(v,r(v))$ with $v\in B^{\rm su}(u,r)
\cap\wt\E$, as follows.  For the first one, choose a ball $B_1=B^{\rm
  su} (v_1,r(v_1))$ with maximal radius $r(v_1)$.  Assume that the
first $n$ balls $(B_i)_{1\le i\le n}$ are constructed. Let $i_*=n+1$
if for every $v\in B^{\rm su}(u,r) \cap \wt\E$, the ball $B^{\rm su}
(v,r(v))$ meets some ball $B_i$ with $1\le i\le n$. Otherwise, choose
a ball $B_{n+1}= B^{\rm su} (v_{n+1},r(v_{n+1}))$ disjoint from all
balls $B_i$ with $1\le i\le n$, with maximal radius $r(v_{n+1})$
amongst them.  If the process does not stop, let $i_*=+\infty$. Note
that $(r(v_i))_{1\leq i<i_*}$ is nonincreasing by the maximality
property at each step. Since $B^{\rm su}(u,r)\cap \wt\E$ is relatively
compact, for every $\rho>0$, there are only finitely many pairwise
disjoint balls centred at a point $v \in B^{\rm su}(u,r)\cap \wt\E$
of radius $r(v)$ at least $\rho$. Hence if $i_*=+\infty$, then
$r(v_i)$ tends to $0$ as $i\ra+\infty$.

Let us prove that $B^{\rm su}(u,r) \cap\wt\E$ is contained in
$\bigcup_{1\leq i<i_*} B^{\rm su}(v_i,2\,r(v_i))$.  Indeed, let $v\in
B^{\rm su} (u,r) \cap\wt\E$ and let $k=\max\{i\in[1,i_*[\,\cap\,\NN
\;:\; r(v_i)\geq r(v)\}$, which exists by the maximality of $r(v_1)$
and since $\lim_{i\ra+\infty} r(v_i)=0$ if $i_*=+\infty$. Then either
$k+1=i_*$ or $k+1<i_*$. In the first case, by construction, the open
ball $B^{\rm su}(v,r(v))$ meets $B_i$ for some $i\in \{1,\dots,k\}$,
say at a point $v'$. Hence, by the triangle inequality,
$$
d(v,v_i)\leq d(v,v')+ d(v',v_i) < r(v)+r(v_i)\leq 2\,r(v_i)\;.
$$  
Therefore $v\in B^{\rm su} (v_i,2\,r(v_i))$. If we have $k+1<i_*$,
then $r(v)> r(v_{k+1})$ by the maximality of $k$. Hence, by the
maximality property at the step $k+1$, the ball $B^{\rm su}(v, r(v))$
meets $B_i$ for some $i\in \{1,\dots,k\}$, and hence as above $v\in
B^{\rm su} (v_i, 2\,r(v_i))$.

This family of balls hence satisfies 
\begin{align}\label{argumentVitali}
\sum_{1\leq i<i_*}\MM(B_i) &\ge \frac{1}{c_2}\,
\sum_{1\leq i<i_*}\MM\big(B^{\rm su}(v_i,2\,r(v_i))\big)\ge 
\frac{1}{c_2}\;\MM\Big(\bigcup_{1\leq i<i_*}
B^{\rm su}(v_i,2\,r(v_i))\Big) \nonumber\\ &
\ge \frac{1}{c_2}\;\MM(B^{\rm su}(u,r)\cap \wt\E)\;.
\end{align}

Now, using respectively  

$\bullet$~ the fact that $u$ does not belong to $\wt Z_r$ (see
Equation \eqref{eq:defigranZsubr}),

$\bullet$~ the fact that the balls $B_i$ for $1\leq i<i_*$ are
pairwise disjoint and contained in $B^{\rm su}(u,r+\sigma)$ (since
$v_i\in \wt \E$ implies that $r(v_i)\leq \sigma$ by the definition of
$\wt \E$),

$\bullet$~ the equality $W^{\rm su}(v_i)=W^{\rm su}(u)$ and the cocycle
property of $c_{\wt F}$,

$\bullet$~ Lemma \ref{lem:step5} in Step 5 since $\wt \E\subset \wt
\Omega_R$,

$\bullet$~ the cocycle property of $c_{\wt F}$,

$\bullet$~ the conclusion (\ref{argumentVitali}) of the above Vitali
argument,

$\bullet$~ the fact that $u$ does not belong to $\wt \Delta_r$ (see
Equation \eqref{eq:defiDeltasubr}),

\noindent we have
\begin{align*}
  \wt I_r\wt\psi\,(u)&\geq\frac{1}{2}\;
\wt I_{r+\sigma}\wt\psi\,(u) \geq \frac{1}{2}\;\sum_{1\leq i <i_*}
\int_{w\in B_i}\wt\psi(w)\;e^{c_{\wt F}(w,\,u)}\;d\mu_{W^{\rm su}(u)}(w)\\ & =
\frac{1}{2}\;\sum_{1\leq i <i_*}\;
\wt I_{r(v_i)}\wt \psi\,(v_i)\;e^{c_{\wt F}(v_i,u)}\\ &
\geq \frac{c_1}{2}\;\int_{T^1M}\psi\;dm_F\;\sum_{1\leq i <i_*}\;
\wt I_{r(v_i)}\wt \rho\,(v_i)\;e^{c_{\wt F}(v_i,u)}\\ & 
=\frac{c_1}{2}\;\int_{T^1M}\psi\;dm_F\;\sum_{1\leq i <i_*}
\int_{w\in B_i}\wt\rho(w)\;e^{c_{\wt F}(w,\,u)}\;
d\mu_{W^{\rm su}(u)}(w)\\ & \geq 
\frac{c_1}{2\,c_2}\;\int_{T^1M}\psi\;dm_F\;
\int_{w\in B^{\rm su}(u,\,r)\,\cap \,\wt \E}
\wt\rho(w)\;e^{c_{\wt F}(w,\,u)}\;d\mu_{W^{\rm su}(u)}(w)\\ & \geq 
\frac{c_1}{4\,c_2}\;\int_{T^1M}\psi\;dm_F\;\;
\wt I_r\wt\rho\,(u)\;.
\end{align*}
This proves Lemma \ref{lem:step6} with $c_3=\frac{c_1}{4\,c_2}$.
\cqfd

\bigskip\noindent{\bf Step 7. } In this penultimate step, we conclude
the search of a lower bound of $\frac{\wt I_{r}\wt \psi\,} {\wt
  I_{r}\wt \rho}(u)$ for $r$ large enough, on a large enough subset of
$u\in T^1\wt M$, in order to obtain a lower bound of the left hand
side $\int_{T^1M} \frac{I_{r}\psi\,}
{I_{r}\rho}\;\rho\;dm^{\nu,\,\mu}$ of Equation \eqref{eq:adjonction}
by a positive constant times $\int_{T^1M} \psi\;dm_F$.

We start by proving that the good set $G_r$ has large measure for some
$r$ large enough, or equivalently that the bad sets $Z_r$ and
$\Delta_r$ have small measures.

\blemm \label{lem:majobadDeltar}
For every $r\geq 1$, we have $\Pi(\Delta_r)\leq \frac{1}{4}$.
\elemm

\dem
By the crucial adjointness property of the operators $I_r$ and $J_r$
(see Lemma \ref{lem:adjoint} in Step 2), we have
\begin{equation}\label{eq:majomesDeltar}
\int_{T^1M} I_r(\rho\;\mathbbm{1}_{\,^c\E})\;
\frac{\rho}{I_r\rho}\;\mathbbm{1}_{\,\Delta_r}\;dm^{\nu,\,\mu}=
\int_{T^1M} \rho\;\mathbbm{1}_{\,^c\E}\;
J_r\Big(\frac{\rho}{I_r\rho}\;\mathbbm{1}_{\,\Delta_r}\Big)
\;dm^{\nu,\,\mu}\;.
\end{equation}

Since $\mathbbm{1}_{\,^c\E}=1-\mathbbm{1}_{\E}$, the definition of
$\Delta_r$ (see Equation \eqref{eq:defiDeltasubr}) implies that for
every $u\in\Delta_r$, we have
$$
\frac{I_r(\rho\;\mathbbm{1}_{\,^c\E})}{I_r\rho}(u)=1-
\frac{I_r(\rho\;\mathbbm{1}_{\E})}{I_r\rho}(u)\geq\frac{1}{2}\;.
$$
Hence the left hand side of Equation \eqref{eq:majomesDeltar} is
bounded from below by $\frac{1}{2}\;\int_{T^1M}
\rho\;\mathbbm{1}_{\,\Delta_r}\;dm^{\nu,\,\mu}$
$=\frac{1}{2}\;\Pi(\Delta_r)$. Since $\mathbbm{1}_{\,\Delta_r}\leq 1$,
since $J_r$ is a positive operator, and by Lemma
\ref{lem:majoJrhoparN} in Step 3, the right hand side of Equation
\eqref{eq:majomesDeltar} is bounded from above by $N\,\Pi(^c\E)$,
which is at most $\frac{1}{8}$ by the definition of $\E$ in the
beginning of Step 6. This proves the result.  \cqfd

\medskip Before proving that the bad set $Z_r$ also has small measure
for some $r$ large enough, we start by showing that the map $r\mapsto
\wt I_r\wt \varphi\,(u)$, defined by integrating nonnegative Borel
maps $\wt \varphi$ (for the strong unstable measure weighted by the
cocycle $c_{\wt F}$) on strong unstable balls, has subexponential
growth in the radius of the ball. When $\wt F$ is bounded, the proof
below even proves its polynomial growth, recovering with a different
presentation the result of Roblin when $\wt F=0$.

\blemm \label{lem:polynomialgrowth} For every $u\in T^1M$ and for
every bounded nonnegative measurable map $\varphi:T^1M\ra \RR$,
there exists $c_u>0$ such that, as $r$ tends to $+\infty$,
$$
I_r\varphi\,(u)=\operatorname{O}(e^{c_u\,(\ln r)^2})\;.
$$
\elemm

\dem We may assume that $\varphi$ is constant with value $1$. We fix
$u\in T^1\wt M$ and we allow $v$ to vary in $W^{\rm su}(u)$.  Define
$x=\pi(u)$ and $y=\pi(v)$. By the definition of $\mu_{W^{\rm su}(u)}$ (see
Equation \eqref{eq:defmeasstrongunstable}), we have
$$
e^{c_{\wt F}(v,\,u)}\;d\mu_{W^{\rm su}(u)}(v)=
e^{c_{\wt F}(v,\,u)}e^{C_{F-\delta,\,v_+}(x,\,y)}\;d\mu_{x}(v_+)\;.
$$
Since the Patterson measure $\mu_{x}$ is finite, we hence only have to
prove that, as $v$ goes to infinity in
$W^{\rm su}(u)$, we have
$$
c_{\wt F}(v,u)+C_{F-\delta,\,v_+}(x,y)=\operatorname{O}\big((\ln
d_{W^{\rm su}(u)}(u,v))^2\big)\;.
$$

\medskip
\noindent\begin{minipage}{8.9cm}~~~ 
  By the convexity of horoballs, there exists $v'\in W^{\rm ss}(v)$ which
  is tangent to the geodesic ray from $x$ to $v_+$.  Let
  $y'=\pi(v')$. By the properties of geodesic triangles in
  $\operatorname{CAT}(-1)$ spaces (or the techniques of approximation by
  trees), there exists a universal constant $c''\geq 0$ such that if
  $T=\ln d_{W^{\rm su}(u)}(u,v)$ is large enough, then

  (i)~ $d_{W^{\rm ss}(\phi_{-T}v)}(\phi_{-T}v,\phi_{-T}v')\leq c''$,

  (ii)~ $d_{W^{\rm su}(\phi_{-T}v)}(\phi_{-T}v,\phi_{-T}u)\leq c''$,

  (iii)~ $|\,d(y',x)-2T\,|\leq c''$.
\end{minipage}
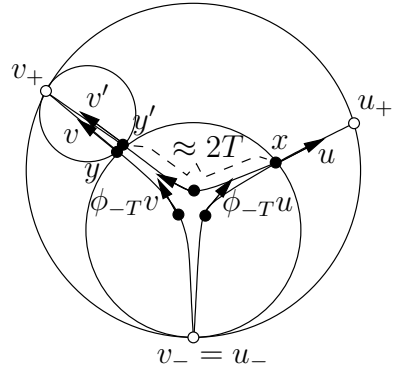
\begin{minipage}{6cm}
\begin{center}
\input{fig_subexpgrowth.pstex_t}
\end{center}
\end{minipage}

\medskip
By the definition of the cocycles, we have
\begin{align*}
c_{\wt F}(v,u)+C_{F-\delta,\,v_+}(x,y)&=
\int_0^{+\infty}\wt F(\phi_{-t}v)-\wt F(\phi_{-t}u)\;dt +
\int_0^{+\infty}\wt F(\phi_{t}v)-\wt F(\phi_{t}v')\;dt\\ &\;\;\;\; -
\int_{-d(y',\,x)}^0\big(\wt F(\phi_{t}v')-\delta\big)\;dt\\ &=
\int_T^{+\infty}\wt F(\phi_{-t}v)-\wt F(\phi_{-t}u)\;dt +
\int_{-T}^{+\infty}\wt F(\phi_{t}v)-\wt F(\phi_{t}v')\;dt\\ &\;\;\;\; 
+\delta\,d(y',x)-
\int_{-d(y',\,x)}^{-T}\wt F(\phi_{t}v')\;dt-
\int_{0}^{T}\wt F(\phi_{-t}u)\;dt\;.
\end{align*}
Let $c\geq 0$ and $\alpha\in\mathopen{]}0,1\mathclose{]}$ be the
H\"older constants of $\wt F$, as in the beginning of the proof of Lemma
\ref{lem:holderargumagain}.  By Equation \eqref{eq:compardistHamen}
and Equation \eqref{eq:dilatHamdist}, there exists $c'>0$ such that
$$
d(\phi_{-t}v,\phi_{-t}u)\leq
c'd_{W^{\rm su}(\phi_{-t}v)}(\phi_{-t}v,\phi_{-t}u)=
c'e^{T-t}d_{W^{\rm su}(\phi_{-T}v)}(\phi_{-T}v,\phi_{-T}u)\;.
$$
By (ii), we hence have
$$
\int_T^{+\infty}\wt F(\phi_{-t}v)-\wt F(\phi_{-t}u)\;dt\leq
\frac{c}{\alpha}(c'c'')^\alpha\;.
$$
Similarly by (i),  we
have 
$$
\int_{-T}^{+\infty}\wt F(\phi_{t}v)-\wt F(\phi_{t}v')\;dt\leq
\frac{c}{\alpha}(c'c'')^\alpha\;.
$$
As mentioned in the remark at the end of Subsection
\ref{subsec:holdercont}, by the H\"older-continuity of $\wt F$, we have
$$
|\,\wt F(w)-\wt F(w')\,|\leq 3\,c\, d(w,w')+c
$$ 
for all $w,w'\in T^1\wt M$. Hence by (iii) and since
$d(\phi_tw,w)=|t|$ for every $w\in T^1\wt M$, we have
\begin{align*}
&|\,c_{\wt F}(v,u)+C_{F-\delta,\,v_+}(x,y)\,|\\ \leq\;\; &
2\frac{c}{\alpha}(c'c'')^\alpha+
\delta\,(2T+c'')+2\,(T+c'')\big(\max_{\pi^{-1}(\{x\})}|\wt F|+
3c\,(T+c'')+c\big)\;.
\end{align*}
The result follows. \cqfd

\medskip Let us now prove that the bad set $Z_r$ also has small
measure for some $r$ large enough.

\blemm \label{lem:majobadZr} For every $r_0\geq 1$, there exists
$r\geq r_0$ such that $\Pi(Z_r)\leq \frac{1}{4}$.  
\elemm

\dem If the interior of the support of $\psi$ does not meet
$\Omega_c\Ga$, then $Z_r$ is empty, and the result is true. We
hence assume that the interior of the support of $\psi$ meets
$\Omega_c\Ga$.

Let us fix $u\in \wt\Omega_c\Ga$. For every $t\geq 0$, let $n(t)$ be
the number of integers $k$ in $\mathopen{[}0,\lfloor \frac{t}{\sigma}
\rfloor\mathclose{]}$ such that there exists $r\in \mathopen{[}
k\sigma, (k+1) \sigma\mathclose{[}$ with $u\in\wt Z_r$. Hence the
period of time spent by $u$ in the sets $\wt Z_r$ as $r$ ranges from
$0$ to $t$, that is $\int_0^t \mathbbm{1}_{\wt Z_r}(u)\;dr$, is at
most $n(t)\,\sigma$. To prove Lemma \ref{lem:majobadZr}, we will first
prove that $n(t)$ grows slowly in $t$.

Recall without proof the following result of \cite{DalBo00}
(generalising, to our (almost no) assumptions on $\wt M$ and $\Ga$,
results of Hedlund and others), which is valid since the geodesic flow
is mixing, hence topologically mixing on $\Omega\Ga$ by Babillot's
Theorem \ref{theo:babmix}.

\bprop[Dal'Bo]\label{prop:densitystrongunstable} For every
$v\in\wt\Omega_c\Ga$, the intersection with $\wt \Omega_c\Ga$ of the
orbit $\Ga W^{\rm su}(v)$ under $\Ga$ of the strong unstable leaf of $v$
is dense in $\wt \Omega_c\Ga$. 
\eprop

Hence $W^{\rm su}(u)$ meets the (nonempty) intersection with
$\wt\Omega_c\Ga$ of the interior of the support of $\wt \psi$. In
particular, there exists $t_0\geq 1$ such that $\wt I_{t_0}\wt\psi
\,(u)>0$. As the hinted proof in \cite[page 91]{Roblin03} is slightly
incorrect, let us prove the following claim.

\medskip\noindent{\bf Claim. } There exists $c>0$ (depending on $u$)
such that for every $t\geq t_0$,
$$
\wt I_{t+2\sigma}\wt\psi\,(u)\geq c\; 2^{\frac{n(t)}{2}}\;.
$$

\dem To simplify the notation, let $f$ be the nondecreasing map
$r\mapsto \wt I_r\wt\psi\,(u)$, and for every $t\geq t_0$, let
$n=n(t)$. Let $0\leq k_1<k_2<\dots<k_{n}\leq\lfloor
\frac{t}{\sigma}\rfloor$ and $r_i\in\mathopen{[}k_i\sigma,
(k_i+1)\sigma\mathclose{[}$ such that $u\in\wt Z_{r_i}$.  Note that
$r_{i+2}\geq r_i+\sigma$ (but $r_{i+1}-r_i$ could be strictly less
than $\sigma$). Let $p=\lfloor \frac{t_0}{\sigma}\rfloor+2$ and
$q=\lfloor \frac{n-p}{2}\rfloor$.  Note that $n-2q-1\geq p-1\geq 0$
and
$$
r_{n-2q}\geq k_{n-2q}\,\sigma\geq (n-2q - 1)\,\sigma\geq (p-1)\,
\sigma\geq t_0\;.
$$

By applying $q+1$ times the definition of the sets $\wt Z_r$ (see
Equation \eqref{eq:defigranZsubr}) and by monotonicity, we have
\begin{align*}
f(t+2\sigma)&\geq f(r_n+\sigma)\geq 2\,f(r_n)\geq 2\,f(r_{n-2}+\sigma)
\geq \dots\geq 2^qf(r_{n-2q}+\sigma)\\ & 
\geq 2^{q+1}f(r_{n-2q})\geq 2^{\frac{n-p}{2}}\,f(t_0)\;.
\end{align*}
Hence the result follows with $c=2^{-\frac{p}{2}}\,f(t_0)>0$.
\cqfd

\medskip This claim and Lemma \ref{lem:polynomialgrowth} imply the
following slow growth property of $n$: we have $n(t)=
\operatorname{O}((\ln t)^2)$ as $t\ra+\infty$.  Therefore, for every
$u\in\wt\Omega_c\Ga$, we have, by the definition of the map $t \mapsto
n(t)$,
$$
\lim_{t\ra+\infty}\frac{1}{t}\;\int_0^t\mathbbm{1}_{\wt Z_r}(u)\;dr=0\;.
$$
Since $\Omega_c\Ga$ has full measure with respect to the probability
measure $\Pi$, and by Lebesgue's dominated convergence theorem, we
have
$$
\lim_{t\ra+\infty}\frac{1}{t}\;\int_0^t\;\Pi(Z_r)\;dr=0\;.
$$
Lemma  \ref{lem:majobadZr} now follows immediately.
\cqfd

\medskip Let us now conclude the proof that the left hand side
$\int_{T^1M} \frac{I_{r}\psi\,} {I_{r}\rho}\;\rho\;dm^{\nu,\,\mu}$ of
Equation \eqref{eq:adjonction} is at least a positive constant times
$\int_{T^1M} \psi\;dm_F$.

By Lemma \ref{lem:majobadZr} and Lemma \ref{lem:majobadDeltar}, since
$G_r=\Omega_c\Ga-(\Delta_r\cup Z_r)$ and since $\Omega_c\Ga$ has full
measure with respect to $\Pi=\rho\, m^{\nu,\,\mu}$ by Lemma
\ref{lem:constructmesnumu}, there exists $r\geq 1$ large enough such
that
$$
\Pi(G_r)\geq 1-\Pi(Z_r)-\Pi(\Delta_r)\geq \frac{1}{2}\;.
$$
Hence, using Lemma \ref{lem:step6}, we have
\begin{align}
\int_{T^1M} \frac{I_{r}\psi\,} {I_{r}\rho}\;\rho\;dm^{\nu,\,\mu}&
=\int_{\Omega_c\Ga} \frac{I_{r}\psi\,} {I_{r}\rho}\;d\Pi\geq 
\int_{G_r} \frac{I_{r}\psi\,} {I_{r}\rho}\;d\Pi
\nonumber\\ & \geq 
c_3\;\Pi(G_r)\;\int_{T^1M}\psi\;dm_F\geq
\frac{c_3}{2}\;\int_{T^1M}\psi\;dm_F\;,
\label{eq:minmembregauche}
\end{align}
which is the lower bound we were looking for.

\bigskip\noindent{\bf Step 8. }  This step, which is an easy
ergodicity argument, is the final one to conclude the proof of Theorem
\ref{theo:uniergodtransv}: We prove now that $\nu$ (as given in Step
1) is proportional to $\mu^\iota$, which is ergodic.

By Equation \eqref{eq:adjonction}, Equation \eqref{eq:minmembregauche}
and Equation \eqref{eq:majorat}, with $c_4=\frac{2\,N}{c_3}$, we have,
for every nonnegative $\psi\in\C_c(T^1M;\RR)$,
$$
\int_{T^1M}\psi\;dm_F\leq c_4\;\int_{T^1M}\psi\;dm^{\nu,\,\mu}\;.
$$
Hence $m_F\leq c_4\;m^{\nu,\,\mu}$, and $m_F=m^{\mu^\iota,\,\mu}$ is
absolutely continuous with respect to $m^{\nu,\,\mu}$.  By
disintegration (see Lemma \ref{lem:constructmesnumu} in Step 1), for
every transversal $T$ to the strong unstable foliation, the measure
$\mu^\iota_T$ is hence absolutely continuous with respect to $\nu_T$. Note
that this is true for any $\Ga$-equivariant $c_{\wt
  F}$-quasi-invariant transverse measure $\nu$ which gives full
support to the negatively recurrent set.

This implies in particular that the $\Ga$-equivariant $c_{\wt
  F}$-quasi-invariant transverse measure $\mu^\iota=(\mu^\iota_T)
_{T\in\T( \wt \W^{\rm su})}$ is ergodic. Indeed, let $A$ be a
$\Ga$-invariant $\wt\W^{\rm su}$-saturated (that is, which is a union of
strong unstable leaves) subset of $T^1\wt M$, whose intersection with
any transversal is measurable, and assume that $\mu^\iota_T(A\cap
T)>0$ for at least one transversal $T$. Then $(\mathbbm{1}_{A\cap
  T}\;\mu^\iota_T) _{T\in\T( \wt \W^{\rm su})}$ is also a
$\Ga$-equivariant $c_{\wt F}$-quasi-invariant transverse measure for
the strong unstable foliation. By the above absolute continuity
claim applied to this family, the measure $\mu^\iota_T$ is absolutely
continuous with respect to $\mathbbm{1}_{A\cap T}\;\mu^\iota_T$, hence 
$\mu^\iota_T(\,^c\!A \cap T)=0$ for all transversals $T$.

Now, to prove Step 8, we may assume without loss of generality that
$\nu$ is ergodic.  Let $T$ be a transversal to the strong unstable
foliation, such that $\nu_T(T)>0$ (for the other transversals $T'$, we
have by absolute continuity $\mu^\iota_{T'}=\nu_{T'}=0$). Let $\G_T$
be the pseudo-group of holonomy maps between transversals to $\W^{\rm
  su}$ contained in $T$. Since $\mu^\iota_T$ is absolutely continuous
with respect to $\nu_T$, there exists a measurable map
$f_T=\frac{d\mu^\iota_T}{d\nu_T}:T\ra [0,+\infty[$ (well defined
$\nu_T$-almost everywhere) such that $\mu^\iota_T=f_T\,\nu_T$. Since
$\mu^\iota$ and $\nu$ are quasi-invariant under holonomy with respect
to the same cocycle $c_{\wt F}$, the map $f_T$ is $\nu_T$-almost
everywhere invariant under $\G_T$. By ergodicity, the map $f_T$ is
hence $\nu_T$-almost everywhere equal to a constant $c_T$. By the
density property of strong unstable leaves recalled in Proposition
\ref{prop:densitystrongunstable}, and since $T\cap\wt\Omega_c\Ga$ has
full measure with respect to both $\mu^\iota_T$ and $\nu_T$, the
constant $c_T$ is independent of $T$. Hence $\mu^\iota$ and $\nu$ are
proportional, which concludes the proof of Theorem
\ref{theo:uniergodtransv}.  \cqfd

\bigskip To conclude Chapter \ref{sec:ergtheounistabfolia}, we obtain
a complete classification of the $\Ga$-equivariant $c_{\wt
  F}$-quasi-invariant transverse measures for the strong unstable
foliation of $T^1\wt M$, under the assumption that $\Ga$ is
geometrically finite. We refer to the beginning of Subsection
\ref{subsec:finiteness} for a definition of a geometrically finite
group of isometries of $\wt M$.

The additional measures are the following ones. Let $\wt W$ be a
strong unstable leaf in $T^1\wt M$ such that the family $(\ga\wt
W)_{\ga\in\Ga}$ is locally finite (that is, for every compact subset
$K$ of $T^1\wt M$, the set of $\ga\in\Ga$ such that $\ga\wt W$ meets
$K$ is finite), or, equivalently, such that the image of $\wt W$ in
$T^1M$ is a closed subset of $T^1M$.  Fix $v\in \wt W$ and let
$\Ga_{\wt W}$ be the stabiliser of $\wt W$ in $\Ga$. We denote by
$\D_z$ the unit Dirac mass at a point $z\in T^1\wt M$. For every
transversal $T$ to $\wt\W^{\rm su}$, the measure
$$
\nu^{\wt W}_T=\sum_{\ga\in\Ga/\Ga_{\wt W}, \;w\in\, T\cap\ga \wt W} 
e^{c_{\wt F} (w,\,v)}\;\D_w
$$
is locally finite, and $(\nu^{\wt W}_T)_{T\in\T(\wt\W^{\rm su})}$ is a
nonzero $\Ga$-equivariant family of locally finite measures on
transversals to the strong unstable foliation, which is stable by
restrictions, and satisfies the property (iii) of $c_{\wt
  F}$-quasi-invariant transverse measures. Note that replacing $v$ by
another element of $\wt W$ changes this family only by a
multiplicative constant. Note that, since $\wt W=W^{\rm su}(v)$, the
(discrete, countable) support of $\nu^{\wt W}_T$ is the  set
$\{w\in T\;:\;w_-\in\Ga v_-\}$. When $\Ga$ is torsion free, if $T$ is
a transversal to the strong unstable foliation $\W^{\rm su}$ on $T^1 M$
and if $W$ is a leaf of $\W^{\rm su}$ which is a closed subset of $T^1M$,
the associated measure on $T$ is 
$$
\sum_{w\in\, T\cap W} e^{c_{F} (w,\,v)}\;\D_w
$$ 
where $c_{F}(w,v)=\int_0^{+\infty} F(\phi_{-t}w)-F(\phi_{-t}v)\;dt$
(see Equation \eqref{eq:defcocycleforunstab}).

\bcoro \label{coro:clasmesergodtransv} Under the assumptions of
Theorem \ref{theo:uniergodtransv}, assume furthermore that $\Ga$ is
geo\-metrically finite. Then, up to a multiplicative constant, the
only ergodic $\Ga$-equivariant $c_{\wt F}$-quasi-invariant transverse
measures for the strong unstable foliation on $T^1\wt M$ are

$\bullet$~ the family $(\mu^\iota_T)_{T\in\T(\wt\W^{\rm su})}$ or

$\bullet$ the families $(\nu^{\wt W}_T)_{T\in\T(\wt\W^{\rm su})}$ where $\wt
W$ is a strong unstable leaf of $T^1\wt M$ whose image in $T^1 M$ is
a closed subset. 
\ecoro

\dem Since $\Ga$ is geometrically finite, a point $\xi$ in
$\partial_\infty\wt M-\Lambda_c\Ga$ is either a point in
$\partial_\infty\wt M-\Lambda\Ga$ (in which case its stabiliser in
$\Ga$ is finite) or a (bounded) parabolic limit point (in which case
its stabiliser in $\Ga$ is infinite). In both cases, the image in $M$
of any horosphere centred at $\xi$ is a closed subset: This follows
from the fact that there exists, in both cases, a horoball centred at
$\xi$ which is precisely invariant under the stabiliser of $\xi$ in
$\Ga$ (see the beginning of Subsection \ref{subsec:finiteness} for the
second case).  In particular, for every $w\in T^1\wt M$ such that
$w_-\notin \Lambda_c\Ga$, the image in $T^1M$ of the strong unstable
leaf $W^{\rm su}(w)$ is a closed subset (see also \cite{DalBo00} for
further information on the topology of the strong unstable manifolds
in $M$).

Let $\nu=(\nu_T)_{T\in\T(\wt\W^{\rm su})}$ be an ergodic $\Ga$-equivariant
$c_{\wt F}$-quasi-invariant transverse measure for $\wt
\W^{\rm su}$. 

Assume first that there exists a transversal $T$ to $\wt \W^{\rm su}$
such that the measurable set $A=\{w\in T\;:\;w_-\notin\Lambda_c\Ga\}$
has positive measure with respect to $\nu_T$. Then by ergodicity, the
set $A$, which is saturated by (the intersections with $T$ of) the
leaves of $\wt \W^{\rm su}$, has full measure. Furthermore, there
exists $w\in A$ such that the support of $\nu_T$ is the closure in $T$
of the intersection with $T$ of $\Ga W^{\rm su}(w)$. Since
$w_-\notin\Lambda_c\Ga$, this intersection is already closed, hence by
quasi-invariance under holonomy of $\nu$, with $\wt W=W^{\rm su}(w)$,
we have that $(\nu_T)_{T\in\T(\wt\W^{\rm su})}$ is a multiple of
$(\nu^{\wt W}_T)_{T\in\T(\wt\W^{\rm su})}$.

Assume now on the contrary that $\nu$ gives full measure to the
negatively recurrent set. Then the result follows from Theorem
\ref{theo:uniergodtransv}.  \cqfd

\section{Gibbs states on Galois co\-vers} 
\label{sec:gibbscover}

Let $(\wt M,\Ga,\wt F)$ be as in the beginning of Chapter
\ref{sec:negacurvnot}: $\wt M$ is a complete simply connected
Riemannian manifold, with dimension at least $2$ and pinched sectional
curvature at most $-1$, $\Ga$ is a nonelementary discrete group of
isometries of $\wt M$, and $\wt F :T^1\wt M\ra \RR$ is a
H\"older-continuous $\Ga$-invariant map. Fix $x_0\in\wt M$.

A subgroup $\Ga_0$ of the isometry group $\Isom(\wt M)$ of $\wt M$ is
said {\it to normalise} $\Ga$ if it is contained in the {\it
  normaliser}\index{normaliser}
$$
\gls{normaliser}=\{\ga\in\Isom(\wt M)\;:\;\ga\Ga\ga^{-1}=\Ga\}
$$
of $\Ga$ in $\Isom(\wt M)$.

\medskip We study in this chapter the behaviour of the critical
exponents, the Patterson densities and the Gibbs measures associated
to normal subgroups of $\Ga$. The main point is that this gives a
natural framework in which one can study precisely these objects when
the Gibbs measure is not finite. We prove results analogous to those
in the chapters \ref{sec:GPS}, \ref{sec:critgur}, \ref{sec:GHTSappli},
\ref{sec:finimixGibbs}, \ref{sec:ergtheounistabfolia}, underlining the
fact that, contrarily to most of these chapters, the Gibbs measures
may be infinite. Most of the results are extensions of those when
$F=0$ in \cite{Roblin05}, with similar proofs, and we will concentrate
on the new features. This paper of Roblin furthermore mentions in its
introduction that its contents remain valid for Patterson densities
with potentials (referring to the multiplicative approach of
Ledrappier-Coudène, see the beginning of Chapter \ref{sec:GPS} for an
explanation of the differences with the additive convention considered
in this text).

\medskip The structure of this section is the following one. The first
two subsections give the tools used in the others, an extension of
Mohsen's shadow lemma and an extension of the Fatou-Roblin radial
convergence theorem, concerning a super-group $\Ga_0$ normalising $\Ga$
(Subsection \ref{subsec:shadowprinciple} is used (at least) in the
proofs of Subsections \ref{subsec:Fatou}, \ref{subsec:caracritexpo},
\ref{subsec:caracpattdens}, \ref{subsec:carachopftsuji},
\ref{subsec:amencovcritexp}, though Subsection \ref{subsec:Fatou} is
only used in the proofs of Subsections \ref{subsec:caracpattdens} and
\ref{subsec:classnilpotcov}). The next three subsections
\ref{subsec:caracritexpo}, \ref{subsec:caracpattdens},
\ref{subsec:carachopftsuji} concern extensions of Chapter
\ref{sec:GPS} and \ref{sec:GHTSappli}, where the effect of the
potential $F$ is perturbed by a character $\chi$ of $\Ga$. The last
two subsections give our results concerning the normal subgroups $\Ga'$
of $\Ga$, in particular the equality of its critical exponent with
potential to that of $\Ga$, and the ergodic theory of its Patterson
densities with potential.

\subsection{Improving Mohsen's shadow lemma} 
\label{subsec:shadowprinciple} 

Mohsen's shadow lemma \ref{lem:shadowlemma} for a Patterson density of
$(\Ga,F)$, as its original version by Sullivan \cite{Sullivan79} when
$F=0$ (see \cite[p.~10]{Roblin03}), requires the balls whose shadows
are considered to be centred at special points. The following result,
due to \cite[Théo.~1.1.1]{Roblin05} when $F=0$, improves the range of
validity of Mohsen's shadow lemma, in a uniform way. It will be useful
in order to study Patterson densities on Galois covers of $M$. Its proof
follows those of the above works of Sullivan, Mohsen and Roblin.

\bprop\label{prop:principombres} Let $\Ga_0$ be a subgroup of
$\Isom(\wt M)$, which contains and normalises $\Ga$, such that $\wt F$
is $\Ga_0$-invariant. For every compact subset $K$ of $\wt M$, there
exist nondecreasing maps $R=R_K$ and $C=C_K$ from
$\mathopen{[}\delta_{\Ga,\,F},+\infty\mathclose{[}$ to
$\mathopen{]}0,+\infty\mathclose{[}$ such that for every $\sigma\geq
\delta_{\Ga,\,F}$, for every Patterson density $(\mu_x)_{x\in \wt M}$
for $(\Ga,F)$ of dimension $\sigma$, for all $x,y\in \Ga_0K$, we have
\begin{equation}\label{eq:principombre}
\frac{1}{C(\sigma)}\;\|\mu_y\|\;e^{\int_x^y(\wt F-\sigma)}
\leq \mu_{x}\big(\OOO_xB(y,R(\sigma))\big)
\leq C(\sigma)\;\|\mu_y\|\;e^{\int_x^y(\wt F-\sigma)}\;.
\end{equation}
\eprop

\dem We start by the following key observation, which will be used
several times in this chapter.

\blemm \label{lem:normalpatdens}
For every $\alpha\in N(\Ga)$, if $(\mu_x)_{x\in \wt M}$ is a Patterson
density for $(\Ga,F)$ of dimension $\sigma$, then so is
$((\alpha^{-1})_* \mu_{\alpha x})_{x\in \wt M}$.  
\elemm

\dem For all $\ga\in\Ga$ and $x\in\wt M$, since $\alpha\ga\alpha^{-1}
\in \Ga$, we have
$$
\ga_*(\alpha^{-1})_*\mu_{\alpha x}=
(\alpha^{-1})_*(\alpha\ga\alpha^{-1})_*\mu_{\alpha x}
=(\alpha^{-1})_*\mu_{\alpha \ga x}\;,
$$
and for all $x,y\in\wt M$ and $\xi\in\partial_\infty \wt M$, we have
$$
\frac{d(\alpha^{-1})_*\mu_{\alpha x}}{d(\alpha^{-1})_*\mu_{\alpha y}}(\xi)
=\frac{d\mu_{\alpha x}}{d\mu_{\alpha y}}(\alpha \xi)
=e^{-C_{F-\sigma,\,\alpha \xi}(\alpha x,\,\alpha y)}=
e^{-C_{F-\sigma,\,\xi}(x,\,y)}\;,
$$
since $\wt F$ is $\Ga_0$-invariant. 
\cqfd

\medskip
Now, let $K$ be a fixed compact subset of $\wt M$. Since the statement
is empty if $\delta_{\Ga,\,F}=+\infty$, we assume that
$\delta_{\Ga,\,F}<+\infty$.  By the above observation, we only
have to prove Equation \eqref{eq:principombre} when $x\in K$, since
the validity of Equation \eqref{eq:principombre} when $x\in\alpha K$
for any $\alpha\in \Ga_0$ is obtained by replacing $(\mu_x)_{x\in \wt
  M}$ by $((\alpha^{-1})_* \mu_{\alpha x})_{x\in \wt M}$ and $y$ by
$\alpha^{-1}y$.

Let us prove that there exist non-decreasing maps $R',C':
\mathopen{[}\delta_{\Ga,\,F}, +\infty\mathclose{[}\ra
\mathopen{]}0,+\infty\mathclose{[}$ such that for every $\sigma\geq
\delta_{\Ga,\,F}$, for every Patterson density $(\mu_x)_{x\in \wt M}$
for $(\Ga,F)$ of dimension $\sigma$, for all $x\in K$ and $y\in
\Ga_0K$, we have
\begin{equation}\label{eq:principombrerac}
\frac{1}{C'(\sigma)}\;\|\mu_y\|
\leq \mu_{y}(\OOO_xB(y,R'(\sigma)))
\leq C'(\sigma)\;\|\mu_y\|\;.
\end{equation}
To prove that Equation \eqref{eq:principombrerac} implies Equation
\eqref{eq:principombre}, note that we have, by Equation
\eqref{eq:radonykodensity},
$$
\mu_{x}(\OOO_xB(y,R'(\sigma)))
=\int_{\xi\in\OOO_xB(y,\,R'(\sigma))} e^{-C_{F-\sigma,\,\xi}(x,\,y)}\;d\mu_{y}(\xi)\;.
$$
By Lemma \ref{lem:holderconseq} (2) applied with $r=r_0=R'(\sigma)$,
since $\wt F$ is $\Ga_0$-invariant, hence is uniformly bounded on
$\pi^{-1}(B(y,R'(\sigma)))$ for $y\in\Ga_0 K$, there exists
$C''(\sigma)>0$ (depending only on $r$, and that may be taken to be
nondecreasing in $\sigma$) such that for all $x\in K$ and
$\xi\in\OOO_xB(y,R'(\sigma))$,
$$
\big|\,C_{F-\sigma,\,\xi}(x,y)+
\int_x^y (\wt F-\sigma)\;\big|\leq C''(\sigma)\;.
$$
Hence 
\begin{align*}
e^{-C''(\sigma)}\;e^{\int_x^y (\wt F-\sigma)}\;\mu_{y}(\OOO_xB(y,R'(\sigma)))
&\leq \mu_{x}(\OOO_xB(y,R'(\sigma)))\\ & \leq 
e^{C''(\sigma)}\;e^{\int_x^y (\wt F-\sigma)}\;\mu_{y}(\OOO_xB(y,R'(\sigma)))\;.
\end{align*}
Therefore the result follows with $R=R'$ and $C=C'e^{C''}$.

\medskip Let us now prove the claim involving Equation
\eqref{eq:principombrerac}. 

The upper bound in this equation is clear with $C'(\sigma)=1$, whatever
$R'(\sigma)$ is. For future use, this proves that for every $R>0$, for
every compact subset $K$ of $\wt M$, there exists $c>0$ (depending only
on $K$ and $R$) such that, for all $x,y\in \Ga_0K$, we have
\begin{equation}\label{eq:principombreupper}
\mu_{x}(\OOO_xB(y,R))
\leq c\;\|\mu_y\|\;e^{\int_x^y(\wt F-\sigma)}\;.
\end{equation}

To prove the lower bound in Equation \eqref{eq:principombrerac}, we
assume for a contradiction that there exist sequences
$(x_i)_{i\in\NN}$ and $(x'_i)_{i\in\NN}$ in $K$, $(R_i)_{i\in\NN}$ in
$\mathopen{]}0,+\infty\mathclose{[}$ converging to $+\infty$,
$(\alpha_i)_{i\in\NN}$ in $\Ga_0$ and $(\mu^i)_{i\in\NN}$, where
$\mu^i$ is a Patterson density for $(\Ga,F)$ with bounded dimension
$\sigma_i$ such that
\begin{equation}\label{eq:hypoabsurprincipombr}
\lim_{i\ra+\infty} \;\frac{1}{\|\mu^i_{\alpha_i x_i}\|}\;
\mu^i_{\alpha_i x_i}(\OOO_{x'_i}B(\alpha_i x_i,R_i))\;=0\;.
\end{equation}
Up to extracting a subsequence, we have $\lim x_i=x\in K$, $\lim x'_i=
x'\in K$, $\lim\sigma_i= \sigma\geq \delta_{\Ga,\,F}$ and
$\lim\alpha_i^{-1} x'_i =\xi\in\wt M\cup\partial_\infty\wt M$. By
Equation \eqref{eq:hypoabsurprincipombr}, the shadow
$\OOO_{x'_i}B(\alpha_i x_i,R_i)$ is different from $\partial_\infty\wt
M$, hence $d(x'_i,\alpha_i x_i)\geq R_i$, which tends to $+\infty$ as
$i\ra+\infty$. Therefore $\xi\in\partial_\infty\wt M$. Let
$$
\nu^i= \;\frac{1}{\|\mu^i_{\alpha_i x_i}\|}\;
(\alpha_i^{-1})_*\mu^i_{\alpha_i x_i}\;,
$$
which is a probability measure on the compact metrisable space
$\partial_\infty\wt M$. By Banach-Alaoglu's theorem, up to extraction,
the sequence $(\nu^i)_{i\in\NN}$ weak-star converges to a probability
measure $\nu_x$. Define, for every $z\in\wt M$, a measure $\nu_z$ on
$\partial_\infty\wt M$ by 
$$
d\nu_z(\xi')=e^{-C_{F-\sigma,\,\xi}(z,\,x)}\; d\nu_x(\xi')\;.
$$
By the continuity of the Gibbs cocycle and by taking limits, the
family $(\nu_z)_{z\in\wt M}$ is a Patterson density for $(\Ga,F)$ of
dimension $\sigma$.

Let $V$ be a relatively compact open subset of $\partial_\infty\wt
M-\{\xi\}$.  Then $V$ is contained in $\OOO_{\alpha_i ^{-1}x'_i}
B(x_i,R_i)$ for $i$ large enough, since $\lim R_i=+\infty$ and
$(x_i)_{i\in\NN}$ stays in the compact subset $K$. Hence
$$
\nu_x(V)\leq \limsup_{i\ra+\infty}\nu^i(\OOO_{\alpha_i ^{-1}x'_i} B(x_i,R_i))=0
$$
by Equation \eqref{eq:hypoabsurprincipombr}. Therefore the measure
$\nu_x$ is supported in $\{\xi\}$. Since the support of a Patterson
density of $(\Ga,F)$ is $\Ga$-invariant, this implies that $\xi$ is
fixed by $\Ga$. This is a contradiction since $\Ga$ is nonelementary.
\cqfd

\medskip The following consequence is due to
\cite[Lem.~1.2.4]{Roblin05} when $F=0$.

\bcoro \label{coro:conseqprincipombr} Let $\Ga_0$ be a discrete
subgroup of $\Isom(\wt M)$, which contains and normalises $\Ga$, such
that $\wt F$ is $\Ga_0$-invariant.  Let $\sigma\in\RR$ and let
$(\mu_x)_{x\in\wt M}$ be a Patterson density for $(\Ga,F)$ of
dimension $\sigma$.

Then the critical exponent of the series $Q'(s)=\sum_{y\in\Ga_0 x_0}\;
\|\mu_y\|\;e^{\int_{x_0}^y(\wt F-s)}$ is at most $\sigma$. If
$\mu_{x_0}$ gives positive measure to $\Lambda_c\Ga_0$, then this
critical exponent is equal to $\sigma$, and this series diverges at
$s=\sigma$.  
\ecoro

\dem Note that $\sigma\geq \delta_{\Ga,\,F}$ by Corollary
\ref{coro:dimgeqpress} (2). Let $R=R_{\{x_0\}}(\sigma)$ and
$C=C_{\{x_0\}}(\sigma)$ be given by Proposition
\ref{prop:principombres} for $K=\{x_0\}$. For every $n\in\NN$, let
$$
E_n=\{y\in\Ga_0 x_0\;:\; n\leq d(x_0,y)<n+1\}\;.
$$
As seen in the proof of Corollary \ref{coro:dimgeqpress}, since
$\Ga_0$ is discrete, there exists $\kappa\in\NN$ such that for every
$n\in\NN$, any element of $\partial_\infty\wt M$ belongs to at most
$\kappa$ elements of the family $\big(\OOO_{x_0}B(y,R)\big)_{y\in E_n}$.

By the lower bound in Equation \eqref{eq:principombre} in Proposition
\ref{prop:principombres}, for all $n\in\NN$, we have
$$
\|\mu_{x_0}\|\geq \frac{1}{\kappa}\sum_{y\in E_n} 
\mu_{x_0}(\OOO_{x_0}B(y,R))
\geq \frac{1}{C\kappa}\sum_{y\in E_n} 
\|\mu_y\|\;e^{\int_{x_0}^y(\wt F-\sigma)}\;.
$$
Hence if $s>\sigma$, then the series $Q'(s)=\sum_{y\in\Ga_0 x_0}\;
\|\mu_y\|\;e^{\int_{x_0}^y(\wt F-s)}$ converges.

If the series $Q'(\sigma)$ converges, let us prove that
$\mu_{x_0}(\Lambda_c\Ga_0) =0$, which gives the last assertion by
contraposition. For every $n\in\NN$, let us define 
$$
A_n=\bigcap_{k\in\NN}\;
\bigcup_{y\in \Ga_0x_0-B(x_0,\,k)} \OOO_{x_0}B(y,n)\;.
$$
For every $n\in\NN$, by Equation \eqref{eq:principombreupper} in the
proof of Proposition \ref{prop:principombres} (with $K=\{x_0\}$) when
$n$ is large enough, there exists $c_n>0$ such that, for every
$k\in\NN$,
$$
\mu_{x_0}(A_n)\leq 
\sum_{y\in \Ga_0x_0-B(x_0,\,k)}\mu_{x_0}(\OOO_{x_0}B(y,n))
\leq c_n\;\sum_{y\in \Ga_0x_0-B(x_0,\,k)}
\|\mu_y\|\;e^{\int_{x_0}^y(\wt F-\sigma)}\;.
$$
Since the sum on the right hand side tends to $0$ as $k\ra+\infty$, we
have $\mu_{x_0}(A_n)=0$ for every $n\in\NN$. Since $\Lambda_c\Ga_0=
\bigcup_{n\in\NN} A_n$ by the definition of the conical limit set, we
have $\mu_{x_0}(\Lambda_c\Ga_0) =0$, as required.  
\cqfd

\medskip Here is another consequence of Proposition
\ref{prop:principombres}, due to \cite[Lem.~1.2.3]{Roblin05} when
$F=0$. For every $r>0$ and every nonelementary discrete group $\Ga_0$
of isometries of $\wt M$, let $\gls{lambdacongar}_0$ be the set of
elements $\xi\in \partial_\infty\wt M$ such that there exist $\rho<r$
and $(\ga_n)_{n\in\NN}$ in $\Ga_0$ such that $(\ga_nx_0)_{n\in\NN}$
converges to $\xi$ and $d(\ga_nx_0,\mathopen{[}x_0,\xi\mathclose{[})<
\rho$ (or equivalently $\xi\in\OOO_{x_0}B(\ga_nx_0,\rho)$). Note that
the sets $\Lambda_{c,\,r}\Ga_0$ are $\Ga_0$-invariant (by the
properties of asymptotic geodesic rays), are nondecreasing in $r$, and
satisfy $\Lambda_c\Ga_0=\bigcup_{r>0} \Lambda_{c,\,r}\Ga_0$.

\bcoro\label{coro:conseq2principombr} Let $\Ga_0$ be a discrete
subgroup of $\Isom(\wt M)$, which contains and normalises $\Ga$, such
that $\wt F$ is $\Ga_0$-invariant.  Let $\sigma\in\RR$ and let
$(\mu_x)_{x\in\wt M}$ be a Patterson density for $(\Ga,F)$ of
dimension $\sigma$. Then for every $r>R_{\{x_0\}}(\sigma)$, we have
$$
\mu_{x_0}(\Lambda_c\Ga_0-\Lambda_{c,\,r}\Ga_0)=0\;.
$$
\ecoro

\dem Let $\rho=R_{\{x_0\}}(\sigma)$ and $r>\rho$. First assume that
the quasi-invariant measure $\mu_{x_0}$ is ergodic for the action
of $\Ga$. If $\mu_{x_0}(\Lambda_c\Ga_0)=0$, the result is clear.
Otherwise, $\mu_{x_0}({}^c\Lambda_c\Ga_0)=0$ by ergodicity, and
$\mu_{x_0}(\!\Lambda_{c,\,R}\Ga_0)>0$ for $R>\rho$ large enough, so that
again by ergodicity $\mu_{x_0}({}^c\Lambda_{c,\,R}\Ga_0)=0$.

Recall without proof the following Vitali covering argument.

\blemm \label{lem:Vitalistrikesagain} (Roblin
\cite[Lem.~1.2.1]{Roblin05}) \label{lem:vitaliprincipombr} For every
$s>0$ and for every discrete subset $Z$ of $\wt M$, there exists a
subset $Z^*$ of $Z$ such that the shadows $\OOO_{x_0}B(z,s)$ for $z\in
Z^*$ are pairwise disjoint and $\bigcup_{z\in Z}\OOO_{x_0}B(z,s)
\subset \bigcup_{z\in Z^*}\OOO_{x_0}B(z,5s)$.  
\elemm

For every $n\in\NN$, applying this lemma with $s=R$ and $Z_n=\Ga_0
x_0-B(x_0,n)$ gives a subset $Z^*_n$ of $Z_n$. Since $R>\rho$, the
shadows $\OOO_{x_0}B(z,\rho)$ for $z\in Z_n^*$ are pairwise
disjoint. Define 
$$
A_n=\bigcup_{z\in Z_n^*}\OOO_{x_0}B(z,\rho)\;.
$$ 
We have $\sigma\geq \delta_{\Ga,\,F}$ by Corollary
\ref{coro:dimgeqpress} (2). By Proposition \ref{prop:principombres}
and the definition of $ \rho$, by Equation
\eqref{eq:principombreupper}, there exist $c,C>0$ such that, for every
$n\in\NN$, we have
\begin{align*}
\mu_{x_0}(A_n)&= \sum_{z\in Z_n^*}\mu_{x_0}\big(\OOO_{x_0}B(z,\rho)\big)
\geq\frac{1}{C}\sum_{z\in Z_n^*}\|\mu_{z}\| \;
e^{\int_{x_0}^z(\wt F-\sigma)}\\ & \geq\frac{1}{c\,C}
\sum_{z\in Z_n^*}\mu_{x_0}\big(\OOO_{x_0}B(z,5R)\big)\geq
\frac{1}{c\,C}\mu_{x_0}\big(\bigcup_{z\in Z_n}
\OOO_{x_0}B(z,R)\big)\\ & \geq \frac{1}{c\,C}
\mu_{x_0}(\Lambda_{c,\,R}\Ga_0)= \frac{\|\mu_{x_0}\|}{c\,C}\;.
\end{align*}
Since $r>\rho$, the set $\Lambda_{c,\,r}\Ga_0$ contains the
nonincreasing intersection of $\bigcup_{z\in Z_n} \OOO_{x_0}B(z,\rho)$
(which contains $A_n$). Hence $\mu_{x_0}(\Lambda_{c,\,r}\Ga_0)\geq
\frac{\|\mu_{x_0}\|}{c\,C(\sigma)}>0$, and again by ergodicity, we
have $\mu_{x_0}({}^c\Lambda_{c,\,r}\Ga_0)=0$.

By a Krein-Milman type of argument, since the map $(\mu_x)_{x\in\wt
  M} \mapsto \mu_{x_0}$ is a bijection from the set of Patterson
densities for $(\Ga,F)$ of dimension $\sigma$ to a convex cone of
finite quasi-invariant measures on $\partial_\infty\wt M$, whose
extremal points are the ergodic ones, the result follows.  
\cqfd

\subsection{The Fatou-Roblin radial convergence theorem} 
\label{subsec:Fatou} 

We state in this subsection the main tool for the classification of
Patterson densities on nilpotent covers of $M$, to be  given in
Subsection \ref{subsec:classnilpotcov}. It is (an immediate extension
of) one of the main results of \cite{Roblin05}, which we will call the
{\it Fatou-Roblin radial convergence}
 theorem. Recall that a sequence of points $(y_i)_{i\in \NN}$ in $\wt
M$ converges {\it radially}\index{convergence!radial}\index{radial
  convergence} to a point $\xi\in\partial_\infty \wt M$ if it
converges to $\xi$ while staying at bounded distance from a geodesic
ray.

Recall that Fatou's (ratio) radial convergence theorem says that given
$g$ and $h$ two positive harmonic functions on the open unit disc
$\DD^2$ in $\RR^2$, for almost every point $\zeta$ on the unit circle
$\SS^1$, the ratio $g(z)/h(z)$ has a limit when $z\in \DD^2$ radially
converges to $\zeta$. Consider a nonelementary Fuchsian group $\Ga'$
that is, when $\DD^2$ is endowed with Poincar\'e's metric (and hence is
a real hyperbolic plane), a non virtually cyclic discrete subgroup of
isometries of $\DD^2$. If $(\mu_x)_{x\in\DD^2}$ is a Patterson density
for $\Ga'$ (without potential), then Sullivan has proved that the map
$x\mapsto \|\mu_x\|$ is harmonic. Hence the ratio of the total masses
of two such Patterson densities has radial limits almost everywhere.

The next result extends this, and we refer to \cite{Roblin05} for more
motivation. 

\medskip Given two (Borel positive) measures $\mu$ and $\nu$ on
$\partial_\infty\wt M$, if $\mu=\mu'+\mu''$ is the unique
decomposition of $\mu$ as the sum of two (Borel positive) measures
$\mu',\mu''$ with $\mu'\ll \nu$ and $\mu''\perp \nu$, we denote by
$\frac{d\mu}{d\nu}$ the ($\nu$-almost everywhere well defined)
Radon-Nikodym derivative of the part $\mu'$ of $\mu$ which is
absolutely continuous with respect to $\nu$, so that
\begin{equation}\label{eq:minoderradnik}
\mu\geq \frac{d\mu}{d\nu}\;\nu\;.
\end{equation}

\btheo [Fatou-Roblin radial convergence
theorem]\label{theo:fatouroblin}
\index{theorem@Theorem!of Fatou-Roblin of radial convergence\\
  with potentials}%
\index{Fatou-Roblin radial convergence theorem \\with potentials} Let
$\Ga_0$ be a discrete \\ subgroup of $\Isom(\wt M)$, which contains
and normalises $\Ga$, such that $\wt F$ is $\Ga_0$-invariant. Let
$\sigma\in\RR$ and let $(\mu_{x})_{x\in\wt M}$ and $(\nu_{x})_{x\in\wt
  M}$ be two Patterson densities of dimension $\sigma$ for
$(\Ga,F)$. Then for $\nu_{x_0}$-almost every $\xi\in \Lambda_c\Ga_0$,
if $(y_i)_{i\in\NN}$ is a sequence in $\Ga_0x_0$ converging radially
to $\xi$, then
$$
\lim_{i\ra+\infty}\frac{\|\mu_{y_i}\|}{\|\nu_{y_i}\|}\;=\;
\frac{d\,\mu_{x_0}}{d\,\nu_{x_0}}(\xi)\;.
$$
\etheo

Note that this theorem is only interesting when $\nu_{x_0}
(\Lambda_c\Ga_0) >0$ and $\nu_{x_0}(\Lambda_c\Ga)=0$, that is, by the
Hopf-Tsuji-Sullivan-Roblin theorem \ref{theo:critnonergo}, when
$(\Ga_0,F)$ is of divergence type and when $(\Ga,F)$ is of convergence
type. In particular, it is trivial if $\Ga_0=\Ga$.

Indeed, the statement is empty if $\nu_{x_0} (\Lambda_c\Ga_0) =0$. And if
$\nu_{x_0}(\Lambda_c\Ga)>0$, then the Hopf-Tsuji-Sullivan-Roblin
theorem \ref{theo:critergo} and the uniqueness property in Corollary
\ref{coro:uniqpatdens} imply that $\Lambda_c\Ga$ has full measure both
for $\nu_{x_0}$ and $\mu_{x_0}$, and that $\nu_{x_0}$ and $\mu_{x_0}$
are proportional, so that the result is immediate.

\medskip
\dem The proof of \cite[Théo.~1.2.2]{Roblin05} extends almost
immediately, except the first of its five steps, that needs some
adaptation. 

We endow any infinite subset of $\Ga_0x_0$ with its {\it Fréchet
  filter}\index{Fréchet filter} of the complementary sets of its
finite subsets. Note that by the discreteness of the orbits and the
definition of the topology on $\wt M\cup \partial_\infty \wt M$, for
all $\rho\geq 0$ and $\xi\in\partial_\infty \wt M$, a sequence
$(y_i)_{i\in\NN}$ in $\Ga_0x_0$, which goes out of every finite subset
of $\Ga_0x_0$, and satisfies $d(y_i,\,\mathopen{[}x_0,\,\xi
\mathclose{[})\leq \rho$ for all $i\in\NN$, converges to $\xi$.

Let $R=R_{\{x_0\}}(\sigma)$ and $C=C_{\{x_0\}}(\sigma)$ be given by
Proposition \ref{prop:principombres} with $K=\{x_0\}$.  Let us fix
$r\in \mathopen{[}R,+\infty \mathclose{[}\,$. For every
$\xi\in\Lambda_{c,\,r}\Ga_0$, let
$$
D^+_r(\xi)=\lim_{\rho\ra r^-}
\;\limsup_{y\in\Ga_0x_0,\;d(y,\,\mathopen{[}x_0,\,\xi\mathclose{[})\leq
  \rho} \;\frac{\|\mu_{y}\|}{\|\nu_{y}\|}\;.
$$
This limit does exist, since the function of $\rho$ is nondecreasing,
and $D^+_r:\Lambda_{c,\,r}\Ga_0\ra \RR\cup\{+\infty\}$ is a Borel
$\Ga$-invariant map, by the properties of asymptotic geodesic rays and
the $\Ga$-equivariance property of $(\mu_{x})_{x\in\wt M}$ and
$(\nu_{x})_{x\in\wt M}$.

\medskip We claim that if $D^+_r(\xi)\geq 1$ for $\nu_{x_0}$-almost
every $\xi\in \Lambda_{c,\,r}\Ga_0$, then there exists $c'>0$ such
that $\mu_{x_0}\geq c'\nu_{x_0}$ on $\Lambda_{c,\,r}\Ga_0$.

Given this claim, the remainder of the proof of
\cite[Théo.~1.2.2]{Roblin05} yields without modification the proof of
Theorem \ref{theo:fatouroblin}.

To prove this claim, for every Borel subset $B$ of $\partial_\infty\wt
M$, and every open neighbourhood $V$ of $B$, let
$$
Z=\{z\in\Ga_0 x_0\;:\; \OOO_{x_0}B(z,r)\subset V\;\;{\rm and}
\;\;\|\mu_z\|\geq \frac{1}{2}\|\nu_z\|\}\;,
$$
and let $Z^*$ be the subset of $Z$ given by Lemma
\ref{lem:Vitalistrikesagain} with $s=r$. The assumption of the above
claim implies that, up to a set of zero $\nu_{x_0}$-measure, the set
$B\cap \Lambda_{c,\,r}\Ga_0$ is contained in $\bigcup_{z\in Z}
\OOO_{x_0}B(z,r)$. Hence, by the properties of $Z^*$, for the constant
$c>0$ given by Equation \eqref{eq:principombreupper}, which depends
only on $5r$ and $K=\{x_0\}$, and by Proposition
\ref{prop:principombres}, we have
\begin{align*}
\nu_{x_0}(B\cap\Lambda_{c,\,r}\Ga_0)& \leq 
\sum_{z\in Z^*}\nu_{x_0}\big(\OOO_{x_0}B(z,5r)\big)
\leq c\sum_{z\in Z^*}\|\nu_{z}\|\;e^{\int_{x_0}^z(\wt F-\sigma)}\\&
\leq 2c\sum_{z\in Z^*}\|\mu_{z}\|\;e^{\int_{x_0}^z(\wt F-\sigma)}
\leq 2c\,C\sum_{z\in Z^*}\mu_{x_0}\big(\OOO_{x_0}B(z,r)\big)\leq  
2c\,C\;\mu_{x_0}(V)\;.
\end{align*}
Since $V$ is an arbitrary open neighbourhood of $B$, we have
$\nu_{x_0}(B\cap\Lambda_{c,\,r}\Ga_0) \leq 2c\,C\;\mu_{x_0}(B)$, which
proves the claim.
\cqfd

\medskip The next result, especially interesting when $\Ga_0$ is
convex-cocompact, gives a criterion for a Patterson density of
$(\Ga,F)$ to be determined by its total masses at given points.

\bcoro \label{coro:uniqtwiPatdens} Let $\Ga_0$ be a discrete subgroup
of $\Isom(\wt M)$, which contains and normalises $\Ga$, such that $\wt
F$ is $\Ga_0$-invariant. Let $\sigma,\sigma'\in\RR$. Let
$(\mu_{x})_{x\in\wt M}$ and $(\nu_{x})_{x\in\wt M}$ be two Patterson
densities for $(\Ga,F)$ of dimension $\sigma$ and $\sigma'$
respectively, giving full measure to $\Lambda_c\Ga_0$. If for every
$y\in\Ga_0x_0$ we have $\|\mu_{y}\|=\|\nu_{y}\|$, then
$(\mu_{x})_{x\in\wt M}=(\nu_{x})_{x\in\wt M}$ (and in particular
$\sigma=\sigma'$).
\ecoro

\dem The proof is the same as when $F=0$ (see
\cite[Coro.~1.2.5]{Roblin05}).  
\cqfd

\subsection{Characters and critical exponents} 
\label{subsec:caracritexpo}

Recall that a (real) {\it character}\index{character} of a group $G$
is a group morphism from $G$ to the additive group $\RR$.

Note that the kernel $\Ker \chi$ of a character $\chi$ of $\Ga$ is a
nonelementary discrete group of isometries of $\wt M$, since otherwise
$\Ga$ is the extension of an abelian group by a virtually nilpotent
group, hence is virtually solvable, hence is elementary.

\medskip In this subsection, we fix a character $\chi$ of
$\Ga$. Define the {\it (twisted) Poincar\'e series of
  $(\Ga,F,\chi)$}\index{Poincar\'e series!twisted by a character} as
$$
Q_{\Ga,\,F,\,\chi}(s)=\gls{Poincarserietwi}(s)=\sum_{\ga\in\Ga} \;\;
e^{\chi(\ga)+\int_x^{\ga y} (\wt F-s)}\;.
$$
The {\it (twisted) critical exponent}%
\index{critical exponent!twisted by a character} of $(\Ga,F,\chi)$ is
the element $\delta_{\Ga,\,F,\,\chi}$ in $\mathopen{[}-\infty,
+\infty\mathclose{]}$ defined by
$$
\gls{criticalexponenttwi}=\limsup_{n\ra +\infty}\;\frac{1}{n}\ln
\;\sum_{\ga\in\Ga,\;n-1< d(x,\,\ga y)\leq n} 
\;e^{\chi(\ga)+\int_x^{\ga y} \wt F}\;.
$$
If $\delta_{\Ga,\,F,\,\chi}<+\infty$ (we will prove in the following
Proposition \ref{prop:proprieteelemexpocritchar} that
$\delta_{\Ga,\,F,\,\chi}>-\infty$), we say that $(\Ga,F,\chi)$ is {\it of
  divergence type}\index{divergence type!twisted} if the series
$Q_{\Ga,\,F,\,\chi}(\delta_{\Ga,\,F,\,\chi})$ diverges, and {\it of
  convergence type}\index{convergence type!twisted} otherwise. 

Since $\wt F$ is H\"older-continuous, the critical exponent
$\delta_{\Ga,\,F,\,\chi}$ of $(\Ga,F,\chi)$ does not depend on $x,y$,
and we will prove below that the above upper limit is a limit if
$\delta_{\Ga,\,F,\,\chi}>0$. When $F=0$, these objects were introduced
in \cite{Roblin05}, and we will extend in this chapter the results of
this reference, which itself extends results of \cite{BabLed96}.

For all $s\geq 0$ and $c>0$, for all $x,y\in\wt M$ and for all open
subsets $U$ and $V$ of $\partial_\infty\wt M$, let
$$
G_{\Ga,\,F,\,\chi,\,x,\,y,\,U,\,V}(t)= 
\sum_{\ga\in\Ga\;:\;d(x,\ga y)\leq t\,,\; 
\ga y\in \C_xU\,,\; \ga^{-1} x\in \C_yV} 
\;e^{\chi(\ga)+\int_x^{\ga y} \wt F}\;,
$$
and
$$
G_{\Ga,\,F,\,\chi,\,x,\,y,\,U,\,V,\,c}(t)= 
\sum_{\ga\in\Ga\;:\;t-c<d(x,\ga y)\leq t\,,\; 
\ga y\in \C_xU\,,\; \ga^{-1} x\in \C_yV} 
\;e^{\chi(\ga)+\int_x^{\ga y} \wt F}\;.
$$
The map $s\mapsto \gls{twibisectorbcountfunc}(s)$ will be called the
{\it (twisted) bisectorial orbital counting function}%
\index{counting function!twisted!bisectorial orbital}%
\index{twisted!bisectorial orbital counting \\function} of $(\Ga,
F,U,V)$ and $s\mapsto \gls{anntwibisectorbcountfunc}(s)$ the {\it
  (twisted) annular bisectorial orbital counting function}%
\index{counting function!twisted!annular bisectorial orbital}%
\index{twisted!annular!bisectorial orbital counting \\function} of
$(\Ga, F,U,V)$. When $V=\partial_\infty \wt M$, we denote them by
$s\mapsto \gls{twisectorbcountfunc}(s)$ and $s\mapsto
\gls{anntwisectorbcountfunc}(s)$, and call them the {\it (twisted)
  sectorial orbital counting function}%
\index{counting function!twisted!sectorial orbital}%
\index{twisted!sectorial orbital counting function} and {\it (twisted)
  annular sectorial orbital counting function}%
\index{counting function!twisted!annular sectorial orbital}%
\index{twisted!annular!sectorial orbital counting  \\function} of
$(\Ga,F,\chi,U)$. When $U=V=\partial_\infty \wt M$, we denote them by
$s\mapsto \gls{twiorbcountfunc}(s)$ and $s\mapsto
\gls{anntwiorbcountfunc}(s)$, and we call them the {\it (twisted) orbital
  counting function}%
\index{counting function!twisted!orbital}%
\index{twisted!orbital counting function} and {\it (twisted) annular
  orbital counting function}%
\index{counting function!twisted!annular orbital}%
\index{twisted!annular!orbital counting function} of $(\Ga,F,\chi)$.

Let $g$ be a periodic orbit of the geodesic flow on $T^1M$. Choose
$\ga_g$ one of the (finitely many up to conjugation) loxodromic
elements of $\Ga$ such that if $x$ is a point of the translation axis
of $\ga_g$, then $g$ is obtained by first lifting to $T^1\wt M$ (by
the unit tangent vectors) the geodesic segment $\mathopen{[}x,\ga_g
x\mathclose{]}$ oriented from $x$ to $\ga_g x$ and then by taking the
image by $T^1\wt M\ra T^1M$. Define $\chi(g)=\chi(\ga_g)$, which does
not depend on the choice of $\ga_g$, since for all $\ga,\ga'\in\Ga$,
we have $\chi(\ga)=\chi(\ga')$ if $\ga$ and $\ga'$ are conjugated or
if $\ga'\ga^{-1}$ has finite order. For every relatively compact open
subset $W$ of $T^1M$ meeting the (topological) non-wandering set
$\Omega\Ga$, we define the {\it (twisted) period counting series}%
\index{period!counting series!twisted by a character} of $(\Ga,F,\chi,
W)$ (see Subsection \ref{subsec:countfunct} when $\chi=0$) as
$$
\gls{periodcountingfucnttwi}(s)=
\sum_{g\in \Per(s),\;g\cap W\neq \emptyset}\;e^{\chi(g)+\int_g F}\;;
$$
We also define, for every $c>0$, the {\it (twisted) annular period
  counting series}%
\index{annular!period counting series!twisted by a character} of
$(\Ga,F,\chi, W)$ as
$$
\gls{annperiodcountingfucnttwi}(s)=
\sum_{g\in \Per(s)-\Per(s-c),\;g\cap W\neq \emptyset}\;e^{\chi(g)+\int_g F}\;;
$$
The {\it (twisted) Gurevich pressure}\index{Gurevich pressure!twisted
  by a character} of $(\Ga,F,\chi)$ is
\begin{equation}\label{eq:defintwigurepress}
\gls{Gurevichpressuretwi}=\limsup_{s\ra+\infty}\;\frac{1}{s} 
\ln Z_{\Ga,\,F,\,\chi,\,W,\,c}(s)\;.
\end{equation}
We will prove below that the twisted Gurevich pressure depends neither
on $W$ nor on $c>0$ and that the above upper limit is a limit if $c$
is large enough.

\medskip
\noindent{\bf Remark. } Assume in this remark that $\Ga$ is torsion
free and that every closed differential $1$-form on $M$ is
cohomologous to a (uniformly locally) H\"older-continuous one (this
second assumption is for instance satisfied if $M$ is compact). As
explained for instance in \cite{Babillot95} in a particular case, the
interpretation of the characters of $\Ga$ as closed $1$-forms on $M$
allows us to reduce the above objects twisted by a character to ones
associated to another potential, as we shall now see.

Given a differential $1$-form $\omega$, which is a smooth map from
$TM$ to $\RR$ (linear on each fibre of $\pi$) that we assume to be
H\"older-continuous, its restriction to $T^1M$, which will also be
denoted by $\omega$, is a (smooth) potential on $T^1M$. Note that
there are potentials which do not come from differential $1$-forms.

If two H\"older-continuous differential $1$-forms $\omega$ and
$\omega'$ are cohomologous (as differential $1$-forms), and if $f:M\ra
\RR$ is a smooth map such that $\omega'-\omega= df$, then their
associated potentials $\omega$ and $\omega'$ are cohomologous (as
potentials) via the map $G=f\circ \pi:T^1M\ra\RR$.

Let $H^1(M;\RR)$ be the space of cohomology classes of closed
differential $1$-forms on $M$. Consider the
map
$$
\mathopen{[}\omega\mathclose{]}\mapsto\Big(\ga\mapsto \int_\ga \omega=
\int_x^{\ga x}\;\wt \omega\Big)\;,
$$
where $x$ is any point of $\wt M$ and $\wt \omega$ is the potential
associated as explained above to the lift $\wt p^*\omega$ of $\omega$ to
$\wt M$ by the canonical projection $\wt p:\wt M\ra M$. By Hurewicz's
theorem, this map from $H^1(M;\RR)$ to $\Hom(\Ga,\RR)$ is a linear
isomorphism. For every $\chi\in \Hom(\Ga,\RR)$, let $\omega_\chi$ be a
H\"older-continuous differential $1$-form on $M$ such that
$\mathopen{[}\omega_\chi\mathclose{]}$ maps to $\chi$, and let $\wt
\omega_\chi$ be its lift to $\wt M$.

Under the hypotheses of this remark, we then immediately have
$$
Q_{\Ga,\,F,\,\chi,\, x,\,x}(s)=Q_{\Ga,\,F+\omega_\chi,\, x,\,x}(s)\;,
$$
$$
\delta_{\Ga,\,F,\,\chi}=\delta_{\Ga,\,F+\omega_\chi}\;,
$$
and
$$
P_{Gur}(\Ga,F,\chi)=P_{Gur}(\Ga,F+\omega_\chi)\;.
$$

When $F=0$, we recover the Poincar\'e series associated with a
cohomology class of $\Ga$ introduced in \cite{Babillot95}, and the
twisted critical exponent is then called the {\it cohomological
  pressure}.%
\index{cohomological pressure}\index{pressure!cohomological}

\medskip After this remark, let us give elementary properties
satisfied by the twisted critical exponents.

\bprop\label{prop:proprieteelemexpocritchar} Let $\chi$ be a character
of $\Ga$.

\smallskip
\noindent (i) The Poincar\'e series of $(\Ga,F,\chi)$ converges if
$s>\delta_{\Ga,\,F,\,\chi}$ and diverges if
$s<\delta_{\Ga,\,F,\,\chi}$. The Poincar\'e series of $(\Ga,F,\chi)$
diverges at $\delta_{\Ga,\,F,\,\chi}$ if and only if the Poincar\'e
series of $(\Ga,F+\kappa,\chi)$ diverges at
$\delta_{\Ga,\,F+\kappa,\,\chi}$, for every $\kappa\in\RR$, and in
particular
$$
\forall\; \kappa\in\RR,\;\;\;\delta_{\Ga,\,F+\kappa,\,\chi}=
\delta_{\Ga,\,F,\,\chi}+\kappa\;.
$$
(ii) We have
$$
\forall\;s\in\RR,\; Q_{\Ga,\,F\circ\iota,\,-\chi,\,x,\,y}(s)=
Q_{\Ga,\,F,\,\chi,\,y,\,x}(s)
\;\;\;{\rm and}\;\;\;
\delta_{\Ga,\,F\circ \iota,\,-\chi}=\delta_{\Ga,\,F,\,\chi}\;.
$$
In particular, $(\Ga,F\circ \iota,-\chi)$ is of divergence type if and
only if $(\Ga,F,\chi)$ is of divergence type.

\smallskip
\noindent (iii) If $\Ga'$ is a nonelementary subgroup of $\Ga$,
denoting by $F':\Ga'\backslash T^1\wt M\ra\RR$ the map induced by $\wt
F$, and by $\chi'$ the restriction to $\Ga'$ of $\chi$, we have
$$
\delta_{\Ga',\,F',\,\chi'}\leq \delta_{\Ga,\,F,\,\chi}\;.
$$
(iv) We have
\begin{equation}\label{eq:minocharpaschar}
\forall\;s\in\RR,\;
 Q_{\Ga,\,F,\,\chi,\,x,\,y}(s)+Q_{\Ga,\,F\circ\iota,\,\chi,\,y,\,x}(s)
\geq Q_{\Ga,\,F,\,x,\,y}(s)
\end{equation} 
and hence 
$$
\max\{\,\delta_{\Ga,\,F,\,\chi},\;\delta_{\Ga,\,F\circ
  \iota, \,\chi}\} \geq \delta_{\Ga,\,F}\;.
$$
(v) If there exists $c\geq 0$ such that $|\chi(\ga)|\leq
c\;d(x_0,\ga x_0)$ for every $\ga\in\Ga$, then
$$
\delta_{\Ga,\,F} - c\;\leq\; 
\delta_{\Ga,\,F,\,\chi}\;\leq \;
\delta_{\Ga,\,F} + c\;.
$$
(vi) We have  $\delta_{\Ga,\,F,\,\chi}>-\infty$.

\smallskip
\noindent (vii) The map $(F,\chi)\mapsto \delta_{\Ga,\,F,\,\chi}$ is
convex, that is, if $\chi^*$ is another character of $\Ga$, if
$\wt{F^*}:T^1\wt M\ra \RR$ is another H\"older-continuous
$\Ga$-invariant map, inducing $F^*:\Ga\backslash T^1\wt M\ra\RR$, if
$\delta_{\Ga,\,F,\,\chi}$ and $\delta_{\Ga,\,F^*,\,\chi^*}$ are
finite, then for every $t\in\mathopen{[}0,1\mathclose{]}$, we have
$$
\delta_{\Ga,\,tF+(1-t)F^*,\,t\chi+(1-t)\chi^*}\leq 
t\,\delta_{\Ga,\,F,\,\chi}+(1-t)\,\delta_{\Ga,\,F^*,\,\chi^*}\;.
$$
(viii) For every $c>0$, we have 
$$
\delta_{\Ga,\,F,\,\chi}=\limsup_{n\ra +\infty}\;\frac{1}{n}\ln
\;\sum_{\ga\in\Ga,\;n-c<d(x,\,\ga y)\leq n } \;e^{\chi(\ga)+\int_x^{\ga y} \wt F}\;,
$$
and if $\delta_{\Ga,\,F,\,\chi}\geq 0$, then 
$$
\delta_{\Ga,\,F,\,\chi}=\limsup_{n\ra +\infty}\;\frac{1}{n}\ln
\;\sum_{\ga\in\Ga,\;d(x,\,\ga y)\leq n } \;e^{\chi(\ga)+\int_x^{\ga y} \wt F}\;.
$$
(ix) If $\wt{F^*}:T^1\wt M\ra \RR$ is another H\"older-continuous
$\Ga$-invariant map, which is cohomologous to $\wt F$, and if
$F^*:T^1M\ra\RR$ is the induced map, then
$$
\delta_{\Ga,\,F,\,\chi}=\delta_{\Ga,\,F^*,\,\chi}\;,
$$
and $(\Ga,F^*,-\chi)$ is of divergence type if and
only if $(\Ga,F,\chi)$ is of divergence type.
\eprop

Before giving a proof, here are a few immediate
consequences. It follows from Assertion (iv) that if $F$ is
reversible, then
\begin{equation}\label{eq:reveregalchipaschi}
\delta_{\Ga,\,F,\,-\chi}=\delta_{\Ga,\,F,\,\chi}\geq \delta_{\Ga,\,F}\;.
\end{equation} 

We claim that the assumption of Assertion (v) is satisfied for
instance if $\Ga$ is convex-cocompact. More generally, assume that
$\Ga$ is finitely generated (which is the case if it is
convex-cocompact), and let $S$ be a finite generating set of
$\Ga$. For every $\ga\in\Ga$, let $\|\ga\|_S$ be the smallest length
of a word in $S\cup S^{-1}$ representing $\ga$. For every character
$\chi$ of $\Ga$, we have
$$
|\chi(\ga)|\leq \|\ga\|_S\;\max_{s\in S}|\chi(s)|\;.
$$
Assume also that the map $\ga\mapsto \ga x_0$ from $\Ga$ (endowed with
the word distance $(\ga,\ga')\mapsto \|\ga^{-1}\ga'\|_S$) to $\wt M$
is quasi-isometric, that is, assume that there exists $c>0$ such that
for every $\ga\in\Ga$, we have
$$
d(x_0,\ga x_0)\geq c\|\ga\|_S\;.
$$
This property does not depend on $x_0$ and is satisfied for instance
if $\Ga$ is convex-cocompact. Then the assumption of Assertion (v) is
indeed satisfied.

If $\Ga$ satisfies the hypothesis of Assertion (v) and if
$\delta_{\Ga,\,F}$ is finite (by Lemma \ref{lem:elemproppressure} (iv),
this is for instance the case if $\wt F$ is bounded on
$\pi^{-1}(\C\Lambda\Ga)$, which itself is the case if $\Ga$ is
convex-cocompact), then $\delta_{\Ga,\,F,\,\chi}$ is finite.

In particular, if $\Ga$ is convex-cocompact, then $\delta_{\Ga,\,F,
  \,\chi}$ is finite for every character $\chi$ of $\Ga$.

\medskip
\dem The proof is an immediate adaptation of the proof of Lemma
\ref{lem:elemproppressure}. For instance, to prove Equation
\eqref{eq:minocharpaschar}, we use
\begin{align*}
\sum_{\ga\in\Ga} \;\; e^{\chi(\ga)+\int_x^{\ga y} (\wt F-s)}+
\sum_{\ga\in\Ga} \;\;e^{\chi(\ga^{-1})+\int_y^{\ga^{-1} x} 
(\wt F\circ\iota-s)}
& =\sum_{\ga\in\Ga} \;\big(e^{\chi(\ga)}+e^{-\chi(\ga)}\big)\;
e^{\int_x^{\ga y} (\wt F-s)}
\\ & \geq \sum_{\ga\in\Ga} \;e^{\int_x^{\ga y} (\wt F-s)}
\;.\;\;\;\Box
\end{align*}

The following results concerning the twisted critical exponent and
twisted Gurevich pressure have  proofs similar to those when
$\chi=0$ seen in Chapter \ref{sec:critgur}.

\btheo\label{theo:twistedlogrsult} (1) For every $c>0$, we have
$G_{\Ga,\,F,\,\chi,\,x,\,y, \,U,\,V,\,c} (t)=
O(e^{\delta_{\Ga,\,F,\,\chi}\;t})$ as $t$ goes to $+\infty$.  If
$\delta_{\Ga,\,F,\,\chi}>0$, we have $G_{\Ga,\,F,\,\chi,\,x,\,y,
  \,U,\,V} (t)= O(e^{\delta_{\Ga,\,F,\,\chi}\;t})$ as $t$ goes to
$+\infty$.

(2) For every $c>$ large enough, we have 
$$
\delta_{\Ga,\,F,\,\chi}=\lim_{n\ra +\infty}\; \frac{1}{n} \ln
\;\sum_{\ga\in\Ga,\;n-c<d(x,\,\ga y)\leq n} \;e^{\chi(\ga)+\int_x^{\ga
    y} \wt F}=\lim_{s\ra +\infty}\; \frac{1}{s} \ln
\;G_{\Ga,\,F,\,\chi,\,x,\,y,\,c} (s)\;.
$$
If $\delta_{\Ga,\,F,\,\chi}>0$, then
$$
\delta_{\Ga,\,F,\,\chi}=\lim_{n\ra +\infty}\;
\frac{1}{n} \ln \;\sum_{\ga\in\Ga,\;d(x,\,\ga y)\leq n}
\;e^{\chi(\ga)+\int_x^{\ga y} \wt F}=\lim_{s\ra
  +\infty}\; \frac{1}{s} \ln \;G_{\Ga,\,F,\,\chi,\,x,\,y} (s)\;.
$$

(3) Let $U$ be an open subset of $\partial_\infty \wt M$ meeting
$\Lambda\Gamma$. For every $c>$ large enough, we have
$$
\lim_{t\ra+\infty} \frac{1}{t}\ln G_{\Ga,\,F,\,\chi,\,x,\,y,
  \,U,\,c}(t) = \delta_{\Ga,\,F,\,\chi}\,,
$$
and if $\delta_{\Ga,\,F,\,\chi}>0$, then
$$
\lim_{t\ra+\infty} \frac{1}{t}\ln G_{\Ga,\,F,\,\chi,\,x,\,y,
  \,U}(t) = \delta_{\Ga,\,F,\,\chi}\,.
$$

(4) Let $U$ and $V$ be any two open subsets of $\partial_\infty \wt M$
meeting the limit set $\Lambda\Gamma$. For every $c>$ large enough, we have
$$\lim_{t\ra+\infty}
\frac{1}{t}\ln G_{\Ga,\,F,\,\chi,\,x,\,y, \,U,\,V,\,c}(t) = 
\delta_{\Ga,\,F,\,\chi}\,,
$$
and if $\delta_{\Ga,\,F,\,\chi}>0$, then 
$$\lim_{t\ra+\infty}
\frac{1}{t}\ln G_{\Ga,\,F,\,\chi,\,x,\,y, \,U,\,V}(t) = 
\delta_{\Ga,\,F,\,\chi}\;.
$$

(5) The definition of $P_{Gur}(\Ga,F,\,\chi)$ in Equation
\eqref{eq:defintwigurepress} does not depend on $c>0$ and we have
$$
P_{Gur}(\Ga,F,\,\chi)=\delta_{\Ga,\,F,\,\chi}\;.
$$
Let $W$ be a relatively compact open subset of $T^1M$ meeting the
non-wandering set $\Omega\Ga$, if $c>0$ is large enough, then
$$
P_{Gur}(\Ga,F,\,\chi)=\lim_{s\ra+\infty}\;\frac{1}{s} 
\ln Z_{\Ga,\,F,\,\chi,\,W,\,c}(s)\;,
$$
and if $P_{Gur}(\Ga,F)>0$, then
$$
P_{Gur}(\Ga,F,\,\chi)=\lim_{s\ra+\infty}\;\frac{1}{s} 
\ln Z_{\Ga,\,F,\,\chi,\,W}(s)\;.
$$
\etheo

\dem (1) See the proof of Corollary \ref{coro:bigOloggrowth} (using
the twisted Patterson measure constructed in the coming Subsection
\ref{subsec:caracpattdens} and the extension of Mohsen's shadow lemma in
Proposition \ref{prop:principombres}).

\medskip (2) See the proof of Theorems \ref{theo:limsupenfaitlimbis}
and \ref{theo:limsupenfaitlim}.

\medskip (3) See the proof of Corollaries \ref{coro:loggrowthbis} (1)
and \ref{coro:loggrowth} (1), where Equation
\eqref{eq:equivmescontcon} should be replaced by
$$
a_{t,\,U',\,z}=e^{-\chi(\ga)}\;a_{t,\,\ga U',\,\ga z}\;,
$$
and the constant $c'_2$ in Equation \eqref{eq:findemloggrowthsect}
should be replaced by $c'_2+\max_{1\leq i\leq k} \;|\chi(\ga_i)|$.

\medskip (4) See the proof of Corollaries \ref{coro:loggrowthbis} (2)
and \ref{coro:loggrowth} (2), with adjustments similar to those
above, in particular since now
$$
b_{t,\,U',\,V',\,z,\,w}=e^{\chi(\alpha)}\;b_{t,\,U',\,\alpha V',\,z, \,
  \alpha w}\;.
$$

\medskip (5) The fact that the twisted Gurevich pressure does not
depend on $c>0$ is proved as the first claim of Assertion (vii) of
Lemma \ref{lem:elemproppressure}. For the other claims, see the proof
of Theorem \ref{theo:equalcritgur}, using the fact that for all
$\ga,\ga'\in\Ga$, we have $\chi(\ga)=\chi(\ga')$ if $\ga'\ga^{-1}$ has
finite order.  
\cqfd

\subsection{Characters and Patterson densities} 
\label{subsec:caracpattdens}

Let $\chi:\Ga\ra\RR$ be a character of $\Ga$.

For every $\sigma\in\RR$, a {\it twisted Patterson
  density}\index{Patterson density!twisted by a character} of
dimension $\sigma$ for $(\Ga,F,\chi)$ is a family of finite nonzero
(positive Borel) measures $(\mu_{x})_{x\in\wt M}$ on $\partial_\infty
\wt M$ such that, for every $\ga\in\Ga$, for all $x,y\in\wt M$, for
every $\xi\in\partial_\infty \wt M$, we have
\begin{equation}\label{eq:twistedequivpatdens}
\ga_*\mu_{x}=\;e^{-\chi(\ga)}\;\mu_{\ga x}\;,
\end{equation}
\begin{equation}\label{eq:twistedRNpatdens}
\frac{d\mu_{x}}{d\mu_{y}}(\xi)=e^{-C_{F-\sigma,\,\xi}(x,\,y)}\;.
\end{equation}
This definition is due to Babillot-Ledrappier \cite{BabLed96} and
Roblin \cite{Roblin05} when $F=0$.

Note that if $\Ga'$ is a nonelementary subgroup of $\Ker\chi$ (for
instance $\Ga'=\Ker\chi$), denoting by $F':\Ga'\backslash T^1\wt M
\ra\RR$ the map induced by $\wt F$, then a twisted Patterson density
of dimension $\sigma$ for $(\Ga,F, \chi)$ is a (standard) Patterson
density of dimension $\sigma$ for $(\Ga',F')$, and in particular
$\sigma\geq \delta_{\Ga',\,F'}$ by Corollary \ref{coro:dimgeqpress}
(2).

Similarly as in Subsection \ref{subsec:gibbspattersondens},
for every twisted Patterson density $(\mu_{x})_{x\in\wt M}$ of dimension
$\sigma$ for $(\Ga,F,\chi)$, we have:

$\bullet$~ $(\mu_{x})_{x\in\wt M}$ is a twisted Patterson density of
dimension $\sigma+s$ for $(\Ga,F+s,\chi)$, for every $s\in\RR$.

$\bullet$~ the support of $\mu_{x}$, which is independent of $x\in\wt
M$, contains $\Lambda\Ga$.

\medskip The twisted Patterson densities satisfy the following
elementary properties, due to \cite{Roblin05} when $F=0$.

\bprop\label{prop:proprielemtwistPattdens} Let $\chi$ be a character
of $\Ga$.

(1) If $\delta_{\Ga,\,F,\,\chi}<+\infty$, then there exists a twisted 
Patterson density of dimension $\delta_{\Ga,\,F,\,\chi}$ for
$(\Ga,F,\chi)$, whose support is exactly $\Lambda\Ga$.  

(2) For every $\sigma\in \RR$, if there exists a twisted Patterson
density of dimension $\sigma$ for $(\Ga,F,\chi)$, then $\sigma\geq
\delta_{\Ga,\,F,\,\chi}$.

(3) If $\Ga'$ is a nonelementary normal subgroup of $\Ga$, with
$F':\Ga'\backslash T^1\wt M\ra\RR$ the map induced by $\wt F$, if
$(\mu'_{x})_{x\in\wt M}$ is a Patterson density of dimension $\sigma
\in \RR$ for $(\Ga',F')$, which is ergodic for the action of $\Ga'$
and quasi-invariant for the action of $\Ga$, then there exists a
character $\chi'$ of $\Ga$, whose kernel contains $\Ga'$, such that
$(\mu'_{x})_{x\in\wt M}$ is a twisted Patterson density of dimension
$\sigma$ for $(\Ga,F,\chi')$.

(4) For every $\sigma\in \RR$ such that there exists a twisted
Patterson density of dimension $\sigma$ for $(\Ga,F,\chi)$, and for
all $x,y\in\wt M$, there exists $c>0$ such that for every $n\in\NN$,
we have
$$
\sum_{\ga\in\Ga\;:\;n-1< d(x,\,\ga y)\leq n}\;
e^{\chi(\ga)+\int_x^{\ga y} (\wt F-\sigma)}\leq c\;.
$$

(5) Assume that $\delta_{\Ga,\,F}<+\infty$, that $(\Ga,F)$ is of
divergence type and that $F$ is reversible. If $\delta_{\Ga,\,F}=
\delta_{\Ga,\,F,\,\chi}$, then $\chi=0$.

(6) Every twisted Patterson density $(\mu_x)_{x\in\wt M}$ of dimension
$\sigma\in \RR$ for $(\Ga,F,\chi)$ satisfies the following ``doubling
property of shadows'': For every compact subset $K$ of $\wt M$, for
every $R>0$ large enough, there exists $C=C(K,R)>0$ such that for all
$\ga \in \Ga$ and $x,y\in K$, we have
$$
\mu_x(\OOO_{x}B(\ga y,5R))\le C\,\mu_x(\OOO_xB(\ga y,R))\;.
$$
\eprop

Note that the reversibility assumption of Assertion (5) is necessary,
as pointed out by the referee: for instance if $\Ga$ is torsion free,
if $M$ is compact and if $\chi$ is a nontrivial character of $\Ga$,
let $\wt F=-\,\wt\omega_\chi$, with the notation of the remark above
Proposition \ref{prop:proprieteelemexpocritchar}, which is not
reversible since $F\circ \iota = -F$. Then
$\delta_{\Ga,\,F,\,0}=\delta_{\Ga,\,0,\,-\chi}$ and
$\delta_{\Ga,\,F,\,2\chi}=\delta_{\Ga,\,0,\,\chi}$ as seen in that
remark, and $\delta_{\Ga,\,0,\,-\chi}=\delta_{\Ga,\,0,\,\chi}$ by
Proposition \ref{prop:proprieteelemexpocritchar} (ii). Hence
$\delta_{\Ga,\,F,\,0}=\delta_{\Ga,\,F,\,2\chi}$, though $2\chi\neq 0$.

\medskip
\dem (1) As in \cite{Roblin05} when $F=0$ and in \cite{Mohsen07} when
$\chi=0$, the construction is similar to Patterson's (see the
proof of Proposition \ref{prop:existPattdens}).  For every $z\in\wt
M$, let $\D_z$ be the unit Dirac mass at $z$. Let $h:\mathopen{[}0,
+\infty \mathclose{[}\ra \mathopen{]}0,+\infty\mathclose{[}$ be a
nondecreasing map such that

$\bullet$~ for every $\epsilon>0$, there exists $r_\epsilon\geq 0$ such
that $h(t+r)\leq e^{\epsilon t}h(r)$ for all $t\geq 0$ and $r\geq
r_\epsilon$;

$\bullet$~ for all $x,y\in\wt M$, if $\overline{Q}_{x,\,y}(s)=
\sum_{\ga\in\Ga} \, e^{\chi(\ga)+\int_x^{\ga y} (\wt F-s)}\;
h(d(x,\ga y))$, then $\overline{Q}_{x,\,y}(s)$ diverges if and only if
$s\leq \delta_{\Ga,\,F,\,\chi}$.

Such a map is easy to construct, by taking $\log h$ to be piecewise
affine on $\mathopen{[}0,+\infty\mathclose{[}$, with positive slope
tending to $0$ near $+\infty$ (when the series
$Q_{\Ga,\,F,\,\chi,\,x,\,y}(s)$ diverges at
$s=\delta_{\Ga,\,F,\,\chi}$, we may take $h=1$ constant).

For $s>\delta=\delta_{\Ga,\,F,\,\chi}$, let us define the measure
$$
\mu_{x,\,s}= \frac{1}{\overline{Q}_{x_0,\,x_0}(s)}\sum_{\ga\in\Ga} \;\;
e^{\chi(\ga)+\int_x^{\ga x_0} (\wt F-s)}\;h(d(x,x_0))\;\D_{\ga x_0}
$$ 
on $\wt M$. By the weak-star compactness of the space of probability
measures on the compact space $\wt M\cup\partial_\infty \wt M$, there
exists a sequence $(s_k)_{k\in\NN}$ in $\mathopen{]}\delta,
+\infty\mathclose{[}$ converging to $\delta$ such that the sequence of
measures $(\mu_{x_0,\,s_k})_{k\in\NN}$ weak-star converges to a
measure $\mu_{x_0}$ on $\wt M\cup\partial_\infty \wt M$, with support
contained in $\Lambda\Ga$, hence equal to $\Lambda\Ga$ by
minimality. For every $x\in\wt M$, defining
$$
d\mu_{x}(\xi)= e^{-C_{F-\delta, \,\xi}(x,\,x_0)}\;d\mu_{x_0}(\xi)\;,
$$ 
it is easy to see that $(\mu_{x})_{x\in\wt M}= \lim_{k\ra+\infty}
(\mu_{x,\,s_k})_{x\in\wt M}$ is a twisted Patterson density of
dimension $\delta$ for $(\Ga,F,\chi)$, with support $\Lambda\Ga$.

\medskip (2) Let $(\mu_{x})_{x\in\wt M}$ be a twisted Patterson
density of dimension $\sigma$ for $(\Ga,F,\chi)$. Note that
$\|\mu_{\ga x_0}\| =e^{\chi(\ga)}\|\mu_{ x_0}\|$ by Equation
\eqref{eq:twistedequivpatdens}.  The result hence follows from the
first claim of Corollary \ref{coro:conseqprincipombr} applied by
taking $(\Ga_0,\Ga)$ therein to be $(\Ga,\Ker\chi)$.

\medskip (3) The proof is the same as the one when $F=0$ in
\cite[Lem.~2.1.1~(c)]{Roblin05}. Let $(\mu'_{x})_{x\in\wt M}$ be as in
the statement of Assertion (3). By Lemma \ref{lem:normalpatdens}
(applied with $\Ga$ therein equal to $\Ga'$), since $\Ga$ normalises
$\Ga'$, for every $\alpha\in \Ga$, the family $((\alpha^{-1})_*
\mu'_{\alpha x})_{x\in \wt M}$ is also a Patterson density for
$(\Ga',F')$ of dimension $\sigma$. By assumption, for every $x\in \wt
M$, the measure $(\alpha^{-1})_* \mu'_{\alpha x}$ is absolutely
continuous with respect to $\mu'_{x}$. By ergodicity, the
($\mu'_{x}$-almost everywhere well defined) Radon-Nikodym derivative
$\frac{d\,(\alpha^{-1})_* \mu'_{\alpha x}} {d\,\mu'_{x}}$ (which is
$\Ga'$-invariant and measurable) is $\mu'_{x}$-almost everywhere equal
to a constant $c(\alpha)>0$. We have $c(\ga')=1$ if $\ga'\in\Ga'$ and
for all $\alpha,\beta\in\Ga$, we have, for $\mu'_x$-almost every
$\xi\in \partial_\infty \wt M$,
\begin{align*}
c(\alpha\beta)&=
\frac{d((\alpha\beta)^{-1})_* \mu'_{\alpha\beta x}}{d\,\mu'_{x}}(\xi)=
\frac{d(\alpha^{-1})_* \mu'_{\alpha\beta x}}
{d\,\mu'_{\beta x}}(\beta \xi)\,
\frac{d(\beta^{-1})_* \mu'_{\beta x}}{d\,\mu'_{x}}(\xi)
\\ & =\frac{d(\alpha^{-1})_* \mu'_{\alpha\beta x}}
{d(\alpha^{-1})_* \mu'_{\alpha x}}(\beta \xi)\,
\frac{d(\alpha^{-1})_* \mu'_{\alpha x}}
{d\,\mu'_{x}}(\beta \xi)\,\frac{d\,\mu'_{x}}
{d\,\mu'_{\beta x}}(\beta \xi)\;c(\beta)
\\ &=c(\alpha)c(\beta)\;,
\end{align*}
this last equality holding since the action of $\beta$ on
$\partial_\infty\wt M$ preserves the measure class of $\mu'_x$ and
since both $((\alpha^{-1})_* \mu'_{\alpha x})_{x\in \wt M}$ and
$(\mu'_{x})_{x\in \wt M}$ are Patterson densities for $(\Ga',F')$ of
dimension $\sigma$ (see Equation \eqref{eq:twistedRNpatdens}).

Hence $\chi':\ga\mapsto\log c(\ga)$ is a character of $\Ga$, trivial on
$\Ga'$, and $(\mu'_{x})_{x\in\wt M}$ is a twisted Patterson density of
dimension $\sigma$ for $(\Ga,F,\chi')$.

\medskip 
(4) The proof is similar to that when $\chi=0$ in
Corollary \ref{coro:dimgeqpress} (1), using Proposition
\ref{prop:principombres} instead of Mohsen's shadow lemma
\ref{lem:shadowlemma}.  

\medskip (5) The proof is similar to that when $F=0$ (see
\cite[Lem.~2.1.2]{Roblin05}). Let $\delta=\delta_{\Ga,\,F}=
\delta_{\Ga,\,F,\,\chi}<+\infty$. Let $(\mu_{x})_{x\in\wt M}$ be a
twisted Patterson density of dimension $\delta=
\delta_{\Ga,\,F,\,\chi}$ for $(\Ga, F,\chi)$, which exists by
Assertion (1). Let $(\mu_{x}^*)_{x\in\wt M}$ be a twisted Patterson
density of dimension $\delta$ for $(\Ga,F,-\chi)$, which exists by
Assertion (1), since
$$
\delta=\delta_{\Ga,\,F,\,\chi}=\delta_{\Ga,\,F\circ\iota,\,-\chi}=
\delta_{\Ga,\,F,\,-\chi}
$$ 
using Proposition \ref{prop:proprieteelemexpocritchar} (ii) for the
second equality, and the reversibility of $F$ and Proposition
\ref{prop:proprieteelemexpocritchar} (ix) for the last equality.  Let
$(\mu_{x}^0)_{x\in\wt M}$ be a Patterson density of dimension
$\delta=\delta_{\Ga,\,F}$ for $(\Ga,F)$, which gives full measure to
$\Lambda_c\Ga$ by Corollary \ref{coro:uniqpatdens}, since $(\Ga,F)$ is
of divergence type. Let $x_0\in\wt M$, and let us normalise the above
Patterson densities so that
$$
\|\mu_{x_0}\|=\|\mu^*_{x_0}\|=\|\mu^0_{x_0}\|\;.
$$
Then $\big(\nu_x=\frac{1}{2}(\mu_{x}+\mu_{x}^*)\big)_{x\in\wt M}$ is a
Patterson density of dimension $\delta$ for $(\Ker \chi ,F)$, such
that, for every $\ga\in\Ga$,
$$
\|\nu_{\ga x_0}\|= \big\|\frac{1}{2}\big(e^{\chi(\ga)}\mu_{x_0}+
e^{-\chi(\ga)}\mu^*_{x_0}\big)\big\|=
\cosh\chi(\ga) \;\|\mu^0_{ x_0}\|\geq
\|\mu^0_{ x_0}\|=\|\mu^0_{\ga x_0}\|\,.
$$ 
By the Fatou-Roblin radial convergence theorem
\ref{theo:fatouroblin} applied with $(\Ga_0,\Ga)$ therein equal to
$(\Ga,\Ker\chi)$ and $\sigma=\delta$, and by Equation
\eqref{eq:minoderradnik}, we hence have $\nu_{x_0}\geq
\frac{d\nu_{x_0}}{d\mu^0_{x_0}} \;\mu^0_{x_0}\geq \;\mu^0_{x_0}$.
Since $\|\nu_{x_0}\|= \|\mu^0_{ x_0}\|$, we therefore have $\nu_{x_0}=
\mu^0_{ x_0}$. Thus $\chi=0$ by Equation
\eqref{eq:twistedequivpatdens} applied to $(\nu_{x})_{x\in\wt M}$ and
$x=x_0$, since $\ga_*\mu^0_{ x_0}=\mu^0_{\ga x_0}$ for every
$\ga\in\Ga$. 

\medskip (6) The proof is similar to the proof of Proposition
\ref{prop:doubling} (1), using Proposition \ref{prop:principombres}
with $(\Ga_0,\Ga)$ therein replaced by $(\Ga,\Ker \chi)$ (or more
precisely its proof, which shows that $R(\sigma)$ may be taken as
large as wanted, and that the constant $C(\sigma)$ appearing for the
upper and lower bounds can be made independent, though depending on
which $R(\sigma)$ is taken), instead of Mohsen's shadow lemma
\ref{lem:shadowlemma}. \cqfd

\medskip
\noindent{\bf Remark. }  Let $\wt F^* :T^1\wt M\ra \RR$ be a
H\"older-continuous $\Ga$-invariant map, which is cohomologous to $\wt
F$ via the $\Ga$-invariant map $\wt G :T^1\wt M\ra \RR$. Let
$(\mu_{x})_{x\in\wt M}$ be a twisted Patterson density of dimension
$\sigma\in\RR$ for $(\Ga,F,\chi)$. Then the family of measures
$(\mu^*_{x})_{x\in\wt M}$ defined by, for every
$\xi\in\partial_\infty\wt M$,
$$
d\mu^*_x(\xi)=e^{\wt G(v_{x\xi})}\;d\mu_x(\xi)\;,
$$
where $v_{x\xi}$ is the tangent vector at the origin of the geodesic
ray from a point $x\in\wt M$ to $\xi$, is  a twisted Patterson density of
dimension also $\sigma$ for $(\Ga,F^*,\chi)$.

\subsection{Characters and the Hopf-Tsuji-Sullivan-Roblin 
theorem} 
\label{subsec:carachopftsuji}

Let $\chi$ be a character of $\Ga$ and $\sigma\in\RR$.  Let
$(\mu^\iota_{x})_{x\in\wt M}$ and $(\mu_{x})_{x\in\wt M}$ be twisted
Patterson densities of the same dimension $\sigma$ for $(\Ga, F\circ
\iota,-\chi)$ and $(\Ga,F,\chi)$ respectively.

\medskip Using the Hopf parametrisation $v\mapsto (v_-,v_+,t)$ with
respect to the base point $x_0$ and Equation \eqref{eq:ecartvisuel},
we define the {\it (twisted) Gibbs measure on $T^1\wt M$ associated
  with the pair of twisted Patterson densities
  $\big((\mu^\iota_{x})_{x\in\wt M}, (\mu_{x})_{x\in\wt M}\big)$}%
\index{twisted!Gibbs measure!on $T^1\wt M$}%
\index{Gibbs measure!twisted!on $T^1\wt M$} as the measure $\wt m$ on
$T^1\wt M$ given by
\begin{align*}
d\,\wt m(v) &= \frac{d\mu^\iota_{x}(v_-)\,d\mu_{x}(v_+)\,dt}
{D_{F-\sigma,\,x}(v_-,v_+)^2} \\ &
= e^{C_{F\circ \iota-\sigma,\,v_-}(x,\,\pi(v))+C_{F-\sigma,\,v_+}(x,\,\pi(v))}
\;d\mu^\iota_{x}(v_-)\,d\mu_{x}(v_+)\,dt\;.
\end{align*}

The measure $\wt m$ on $T^1\wt M$ is independent of $x\in\wt M$ by the
equations \eqref{eq:twistedRNpatdens} and
\eqref{eq:cocycleprop}, hence is invariant under the action of $\Ga$
by the equations \eqref{eq:equivgap} and
\eqref{eq:twistedequivpatdens} (the cancellation
$e^{-\chi(\ga)}e^{\chi(\ga)}=1$ is the reason why it is important that
the Patterson density $(\mu^\iota_{x})_{x\in\wt M}$ is twisted by the
opposite character $-\chi$). It is invariant under the geodesic
flow. Hence it defines a measure $m$ on $T^1 M$ which is invariant
under the quotient geodesic flow, called the {\it (twisted) Gibbs
  measure on $T^1M$ associated with the pair of twisted Patterson
  densities $\big((\mu^\iota_{x})_{x\in\wt M}, (\mu_{x})_{x\in\wt
    M}\big)$}.\index{twisted!Gibbs measure!on $T^1 M$}%
\index{Gibbs measure!twisted!on $T^1 M$}

The measure $\iota_*\wt m$ is the twisted Gibbs measure on $T^1\wt
M$ associated with the switched pair of twisted Patterson densities
$\big((\mu_{x})_{x\in\wt M},(\mu^\iota_{x})_{x\in\wt M}\big)$ (and
similarly on $T^1 M$).

\medskip The following results have proofs similar to those of
Theorem \ref{theo:critnonergo} and Theorem \ref{theo:critergo}.

\btheo \label{theo:twicritnonergo}%
\index{theorem@Theorem!of Hopf-Tsuji-Sullivan-Roblin with potentials and characters!non-ergodic}\index{Hopf-Tsuji-Sullivan-Roblin theorem with \\ potentials and characters!non-ergodic} The
following assertions are equivalent.
\begin{enumerate}
\item[(1)] The twisted Poincar\'e series of $(\Ga,F,\chi)$ at the point
  $\sigma$ converges: $Q_{\Ga,\, F,\,\chi}(\sigma) < + \infty$.
\item[(2)] There exists $x\in \wt M$ such that $\mu_x(\Lambda_c
  \Ga)=0$.
\item[(2')] For every $x\in \wt M$, we have $\mu_x(\Lambda_c
  \Ga)=0$.
\item[(3)] There exists $x\in \wt M$ such that $\mu^\iota_x(\Lambda_c
  \Ga)=0$.
\item[(3')] For every $x\in \wt M$, we have $\mu^\iota_x(\Lambda_c
  \Ga)=0$.
\item[(4)] There exists $x\in \wt M$ such that the dynamical system
  $(\partial_\infty^2 \wt M, \Ga, (\mu^\iota_x\otimes
  \mu_x)_{\mid \partial_\infty^2 \wt M})$ is non-ergodic and
  completely dissipative.
\item[(4')] For every $x\in \wt M$, the dynamical system
  $(\partial_\infty^2 \wt M, \Ga, (\mu^\iota_x\otimes
  \mu_x)_{\mid \partial_\infty^2 \wt M})$ is non-ergodic and
  completely dissipative.
\item[(5)] The dynamical system $(T^1M, (\phi_t)_{t\in\RR}, m)$ is
  non-ergodic and completely dissipative.
\end{enumerate}
\etheo

\btheo\label{theo:twicritergo}%
\index{theorem@Theorem!of Hopf-Tsuji-Sullivan-Roblin with potentials and characters!ergodic}\index{Hopf-Tsuji-Sullivan-Roblin theorem with \\ potentials and characters!ergodic} The following assertions are equivalent.
\begin{enumerate}
\item[(i)]  The twisted Poincar\'e series of $(\Ga,F,\chi)$ at the point
  $\sigma$ diverges: $Q_{\Ga,\, F,\,\chi}(\sigma) = + \infty$.
\item[(ii)] There exists $x\in \wt M$ such that
  $\mu_x(\,^c\Lambda_c \Ga)=0$.
\item[(ii')] For every $x\in \wt M$, we have
  $\mu_x(\,^c\Lambda_c \Ga)=0$.
\item[(iii)] There exists $x\in \wt M$ such that
  $\mu^\iota_x(\,^c\Lambda_c \Ga)=0$.
\item[(iii')] For every $x\in \wt M$, we have
  $\mu^\iota_x(\,^c\Lambda_c \Ga)=0$.
\item[(iv)] There exists $x\in \wt M$ such that the dynamical system
  $(\partial_\infty^2 \wt M, \Ga, (\mu^\iota_x\otimes
  \mu_x)_{\mid \partial_\infty^2 \wt M})$ is ergodic and conservative,
  and $\mu^\iota_x\otimes \mu_x$ is atomless.
\item[(iv)] For every $x\in \wt M$, the dynamical system
  $(\partial_\infty^2 \wt M, \Ga, (\mu^\iota_x\otimes
  \mu_x)_{\mid \partial_\infty^2 \wt M})$ is ergodic and conservative,
  and $\mu^\iota_x\otimes \mu_x$ is atomless.
\item[(v)] The dynamical system $(T^1M, (\phi_t)_{t\in\RR}, m)$ is
  ergodic and conservative. 
\end{enumerate}  
\etheo

\dem Let us give a few words on the proofs of these two results,
besides referring to Subsection \ref{subsec:GHTSR}. We use Proposition
\ref{prop:principombres} with $(\Ga_0,\Ga)$ therein replaced by
$(\Ga,\Ker \chi)$, coupled with Equation \eqref{eq:principombreupper},
instead of Mohsen's shadow lemma \ref{lem:shadowlemma}. In particular,
using Proposition \ref{prop:proprielemtwistPattdens} (4) instead of
Lemma \ref{coro:dimgeqpress} (1), Lemma \ref{lem:lemtechunpourf} now
becomes the following statement:

Let $x\in\wt M$, if $Q_{\Gamma,\, F,\,\chi}(\sigma) = + \infty$ and if $R$ is
large enough, with $K=K(x,R)$, there exists $c'_8> 0$ such that
for every sufficiently large $T > 0$
$$
\int_0^T \int_0^T 
\Big(
\sum_{\alpha, \,\beta \in \Ga} 
\wt m(K \cap \phi_{-t} \alpha K \cap \phi_{-s-t} \beta K)
\Big)
\;ds\; dt \leq c'_8
\Big(
\;\sum_{\ga \in \Ga,\, d(x,\,\ga x) \leq T} 
e^{\chi(\ga)+\int_x^{\ga x}(\wt F -\sigma)}\;
\Big)^2\;.
$$

Similarly, Lemma \ref{lem:lemtechdeuxpourf} becomes:

Let $x\in\wt M$, if $Q_{\Gamma,\, F,\,\chi}(\sigma) = + \infty$ and if $R$
is large enough, with $K=K(x,R)$, there exists $c'_9 > 0$ such that for
every sufficiently large $T > 0$
$$
\int_0^T 
\Big(
\sum_{\ga \in \Ga} \wt m(K \cap \phi_{-t} \ga K)
\Big)
\; dt \geq c'_9 \;\sum_{\ga \in \Ga,\, d(x,\,\ga x) \leq T} 
e^{\chi(\ga)+\int_x^{\ga x}(\wt F -\sigma)} \;.
$$

To prove the analog of Assertion (g) in Proposition
\ref{prop:listequivHTSRG}, we first use Proposition
\ref{prop:proprielemtwistPattdens} (6) instead of Proposition
\ref{prop:doubling} (1). Then, since in the twisted case we have
$\|\mu_{\ga x}\|=e^{\chi(\ga)} \|\mu_x\|$ for every $\ga\in\Ga$,
Equation \eqref{eq:poinclefdansg} becomes
\begin{align*}
\frac{e^{-\chi(\ga_n)}}{c}\;\mu_{\ga_n x}(B\cap \OOO_xB(\ga_n x,R))
&\leq \frac{\mu_x(B\cap \OOO_xB(\ga_n x,R))}{\mu_x(\OOO_xB(\ga_n x,R))}
\\ &\leq c\;e^{-\chi(\ga_n)}\;\mu_{\ga_n x}(B\cap \OOO_xB(\ga_n x,R)) \;.
\end{align*}
and one concludes as in the proof of Assertion (g) using this time
the property \eqref{eq:twistedequivpatdens} of the twisted Patterson
densities.

Now the deduction of Theorem \ref{theo:critnonergo} and Theorem
\ref{theo:critergo} from Proposition \ref{prop:listequivHTSRG} remains
valid.  \cqfd

\medskip As when $\chi=0$ in Subsection \ref{subsec:uniqueness}, the
following results are extensions of respectively Corollary
\ref{coro:critdivergence}, Corollary \ref{coro:uniqpatdens},
Proposition \ref{prop:atomlessreferee} and Corollary
\ref{coro:uniqGibstate}, the last claim being an immediate consequence
of its assertions (1) and  (2) a).

\bcoro\label{coro:twiconseqHTSR}
  (1) Assume that $\delta_{\Ga,\,F,\,\chi}<+\infty$ and that there
  exists a twisted Patterson density $(\mu_x)_{x\in X}$ for
  $(\Ga,F,\chi)$ of dimension $\sigma\geq \delta_{\Ga,\,F,\,\chi}$ such that
  $\mu_{x_0}(\Lambda_c\Ga)>0$, then $\sigma=\delta_{\Ga,\,F,\,\chi}$ and
  $(\Ga,F,\chi)$ is of divergence type.

  (2) Assume that $\delta_{\Ga,\,F,\,\chi}<+\infty$ and that
  $(\Ga,F,\chi)$ is of divergence type.

  ~~~ a) There exists a twisted Patterson density $(\mu_x)_{x\in X}$
  for $(\Ga,F,\chi)$ of dimension $\delta_{\Ga,\,F,\,\chi}$, which is
  unique up to a scalar multiple. The measures $\mu_{x}$ are ergodic
  with respect to $\Ga$, give full measure to $\Lambda_c\Ga$, and are
  atomless. In particular their support is $\Lambda\Ga$.

~~~ b) If $\D_z$ is the unit Dirac mass at a point $z\in \wt M$, then
for every $x\in\wt M$, we have
\begin{equation}\label{eq:twiconstructPattOK}
  \frac{\mu_{x}}{\|\mu_{x}\|}=
  \lim_{s\ra \,{\delta_{\Ga,\,F,\,\chi}}^+} \;
\frac{1}{Q_{\Ga,\,F,\,\chi,\,x,\,x}(s)}\;
  \sum_{\ga \in\Ga} e^{\chi(\ga)+\int_x^{\ga x}(\wt F-s)}\;\D_{\ga x}\;.
\end{equation}

~~~ c)  We have $\mu_{x}(\Lambda\Ga-\Lambda_{\rm Myr} \Ga) =0$.

\medskip (3) Let $\sigma\in\RR$. If there exists a twisted Gibbs
measure of dimension $\sigma$ for $(\Ga,F,\chi)$ which is
finite on $T^1M$, then $\sigma= \delta_{\Ga,\,F,\,\chi}$, the triple
$(\Ga,F,\chi)$ is of divergence type, there exists a unique (up to a
scalar multiple) twisted Gibbs measure with potential $F$ on $T^1M$,
its support is $\Omega\Ga$, and the geodesic flow $(\phi_t)_{t\in\RR}$
on $T^1M$ is ergodic with respect to $m_F$.

\medskip (4) If $\chi$ is a character of $\Ga$, two twisted Patterson
densities for $(\Ga,F,\chi)$, giving full measure to $\Lambda_c\Ga$,
are proportional by a positive constant, and in particular have the
same dimension. 
\ecoro

\dem We only prove in Assertion (2) a) the claim that $\mu_{x}$ is
atomless, to indicate how the extensions proceed.

Let $\sigma=\delta_{\Ga,\,F,\,\chi}<+\infty$, which is equal to
$\delta_{\Ga,\,F\circ \iota,\,-\chi}$ by Proposition
\ref{prop:proprieteelemexpocritchar} (ii). As seen in the beginning of
Subsection \ref{subsec:caracpattdens}, if $\Ga'=\Ker\chi$, denoting by
$F':\Ga'\backslash T^1\wt M \ra\RR$ the map induced by $\wt F$, we
have $\sigma\geq \delta_{\Ga',\,F'}$.

Assume for a contradiction that $\xi\in\partial_\infty \wt M$ is an
atom of $\mu_{x}$. Since $\mu_{x}$ gives full measure to the
conical limit set by previous claim of Assertion (2) a), we have
$\xi\in\Lambda_c\Ga$, that is, there exists a sequence of elements
$(\ga_i)_{i\in\NN}$ in $\Ga$ such that $(\ga_i \,x)_{i\in\NN}$
converges to $\xi$ while staying at bounded distance from the geodesic
ray $\mathopen{[}x,\xi\mathclose{[}\,$, which implies in particular
that $\xi\in \OOO_xB(\ga_i x,R)$ for any $R>0$ large enough. Up to
increasing $R$, by Proposition \ref{prop:principombres}  applied by
taking $(\Ga_0,\Ga)$ therein to be $(\Ga,\Ker\chi)$ and with
$K=\{x\}$, there exists $C>0$ such that
$$
\mu_{x}(\OOO_xB (\ga_i x,R)) 
\leq C\;\|\mu_{\ga_i x}\|\;e^{\int_x^{\ga_i x} (\wt F-  \sigma)}=
C\;e^{\chi(\ga_i)}\;\|\mu_{x}\| \;e^{\int_x^{\ga_i x} (\wt F-  \sigma)}
$$ 
for all $i\in\NN$, using Equation \eqref{eq:twistedequivpatdens}.
Hence there exists $\epsilon>0$ such that for all $i\in\NN$, using
Equation \eqref{eq:timereversal},
$$
e^{\int_x^{\ga_i^{-1} x} (\wt F\circ\iota- \sigma)}=
e^{\int_x^{\ga_i x} (\wt F- \sigma)}\geq\epsilon\;e^{-\chi(\ga_i)}\,.
$$
Note that $(\Ga, F\circ \iota,-\chi)$ is of divergence type by
Proposition \ref{prop:proprieteelemexpocritchar} (ii). Let
$(\mu^\iota_{x})_{x\in\wt M}$ be a Patterson density of dimension
$\sigma=\delta_{\Ga,\,F\circ\iota,\,-\chi}$ for $(\Ga,F\circ\iota,
-\chi)$, which exists by the previous claims of Assertion (2)
a). Again by Proposition \ref{prop:principombres}, we have, for all
$i\in\NN$,
\begin{align*}
\mu^\iota_{x}(\OOO_xB (\ga_i^{-1} x,R)) &\geq
\frac{1}{C}\;\|\mu^\iota_{\ga_i^{-1} x}\|\;e^{\int_x^{\ga_i^{-1} x}
  (\wt F\circ\iota - \sigma)} = \frac{1}{C}\;e^{-\chi(\ga_i^{-1} )}
\;\|\mu^\iota_{x}\|\;e^{\int_x^{\ga_i^{-1}  x} (\wt F\circ\iota- \sigma)}
\\ & \geq \frac{\epsilon}{C}\;\|\mu^\iota_{x}\|>0\;.
\end{align*}
Up to extracting a subsequence, the sequence $(\ga_i^{-1}
x)_{i\in\NN}$ converges to $\eta \in \partial_\infty \wt M$, with
$\mu^\iota_{x}(\{\eta\}) >0$.  In particular, the measure
$\mu^\iota_{x}\otimes\mu_{x}$ has an atom at $(\eta,\xi)$, a
contradiction to Theorem \ref{theo:twicritergo} (iv').  \cqfd

\subsection{Galois covers and critical exponents} 
\label{subsec:amencovcritexp}

Let $\Ga'$ be a nonelementary subgroup of $\Ga$, and let
$F':\Ga'\backslash T^1\wt M\ra\RR$ be the map induced by $\wt F$.

The aim of this subsection is to compare the critical exponents of
$(\Ga,F)$ and $(\Ga',F')$. We saw in Lemma \ref{lem:elemproppressure}
that 
\begin{equation}\label{eq:croissanceexpocrit}
\delta_{\Ga',\,F'}\leq \delta_{\Ga,\,F}\;.
\end{equation}

The case of equality in Equation \eqref{eq:croissanceexpocrit} is
related to the amenability properties of the left quotient $\Ga'\bs\Ga$,
in the following sense.

We endow the left quotient $\Ga'\bs\Ga$ with the discrete topology,
and with the (right) action  of $\Ga$ by translations on the
right.  We denote by $R_\alpha:\mathopen{[}\ga\mathclose{]}
\mapsto\mathopen{[}\ga\alpha\mathclose{]}$ the action of
$\alpha\in\Ga$ on $\Ga'\bs\Ga$. A {\it right-invariant
  mean}\index{right-invariant mean}\index{mean (right-invariant)}
$\mmm$ on the left quotient $\Ga'\bs\Ga$ is a linear map
$\mmm:\ell^\infty(\Ga'\bs\Ga)\ra\RR$ such that $\mmm(1)=1$,
$\mmm(f)\geq 0$ if $f\geq 0$ and $\mmm(f\circ R_\alpha)=\mmm(f)$ for
every $\alpha\in\Ga$. The left quotient $\Ga'\bs\Ga$ is {\it
  amenable}\index{amenable} if it has a right-invariant mean.  When
$\Ga'$ is normal in $\Ga$, there are many equivalent definitions, see
for instance \cite{Greanleaf69,Paterson88}. For instance, if $\Ga'$ is
normal in $\Ga$ and if the group $\Ga'\bs\Ga$ is virtually solvable,
then the left quotient $\Ga'\bs\Ga$ is amenable.

When $F=0$, the following result is due to
\cite[Théo.~2.2.2]{Roblin05}.  We suspect that the normality
assumption on $\Ga'$ might not be necessary.

\btheo\label{theo:amenableegalexpocrit} If $F$ is reversible, if
$\Ga'$ is a nonelementary normal subgroup of $\Ga$ and if $\Ga'\bs\Ga$
is amenable, then $\delta_{\Ga',\,F'}= \delta_{\Ga,\,F}$.  
\etheo

\dem Let $\mmm:\ell^\infty(\Ga'\bs\Ga)$ be a right-invariant mean on
the left quotient $\Ga'\bs\Ga$. We may assume that $\delta'=
\delta_{\Ga',\,F'} <+\infty$, otherwise the result follows from
Equation \eqref{eq:croissanceexpocrit}.

Let $(\mu_{x})_{x\in\wt M}$ be a Patterson density of dimension
$\delta'$ for $(\Ga',F')$, which exists as seen in Proposition
\ref{prop:existPattdens}. For every $\alpha\in\Ga$, we claim that
there exists a map $f_\alpha\in\ell^\infty(\Ga'\bs\Ga)$ such that, for
every $\mathopen{[}\ga\mathclose{]}\in \Ga'\bs\Ga$,
$$
f_\alpha(\mathopen{[}\ga\mathclose{]})=
\ln\|\mu_{\ga\alpha x_0}\|-\ln\|\mu_{\ga x_0}\|\;.
$$
Since $(\mu_{x})_{x\in\wt M}$ is $\Ga'$-equivariant, for every
$\ga\in\Ga$, the value of $f_\alpha(\mathopen{[}\ga\mathclose{]})$
indeed depends only on the class of $\ga$ in $\Ga'\bs\Ga$. By Lemma
\ref{lem:holderconseq} (1), there exists $c_1>0$ such that
\begin{align*}
  |\,f_\alpha(\mathopen{[}\ga\mathclose{]})\,|&
\leq \max_{\xi\in\partial_\infty \wt M}
  \big|\,C_{F-\delta',\,\xi}(\ga\alpha x_0,\ga x_0)\,\big|\\& \leq
  c_1\,e^{d(\alpha x_0,\,x_0)} \;+\; d(\alpha x_0, x_0)\;
  \max_{\pi^{-1}(B(\alpha x_0,\,d(\alpha x_0, \,x_0)))}|\wt F-\delta'|\;.
\end{align*}
Hence $f_\alpha$ indeed belongs to $\ell^\infty(\Ga'\bs\Ga)$.

For all $\alpha,\beta,\ga\in\Ga$, we have
$$
f_{\alpha\beta}(\mathopen{[}\ga\mathclose{]})=
\ln\|\mu_{\ga\alpha\beta x_0}\|-\ln\|\mu_{\ga\alpha x_0}\|+
\ln\|\mu_{\ga\alpha x_0}\|-\ln\|\mu_{\ga x_0}\|\;.
$$
Therefore $f_{\alpha\beta}=f_{\beta}\circ R_\alpha+f_{\alpha}$, and
the map $\chi:\Ga\ra\RR$ defined by $\chi(\alpha)=\mmm(f_\alpha)$ is a
character of $\Ga$.

For every $x\in\wt M$ and every continuous map
$\varphi: \partial_\infty \wt M\ra\RR$, let $f_{x,\,\varphi}\in
\ell^\infty(\Ga'\bs\Ga)$ be the map defined by
$$
f_{x,\,\varphi}(\mathopen{[}\ga\mathclose{]})=
\frac{1}{\|\mu_{\ga x_0}\|}\;\int\varphi\;d\,(\ga^{-1})_*\mu_{\ga x}\;.
$$
This map is well defined, and bounded by an argument similar to the
above one. By Riesz's theorem,  
there exists a (positive Borel) measure $\nu_x$ on $\partial_\infty
\wt M$ such that 
$$
\int\varphi\; d\nu_x=\mmm(f_{x,\,\varphi})
$$ 
for every continuous map $\varphi: \partial_\infty \wt M\ra\RR$.  For
all $x\in\wt M$, $\ga\in\Ga$, $\ga'\in\Ga'$ and
$\varphi:\partial_\infty \wt M\ra \RR$ continuous, using the fact that
$\Ga'$ is normal in $\Ga$, we have
$$
\|\mu_{\ga\ga' x_0}\| = \|\mu_{(\ga\ga'\ga^{-1}) \ga x_0}\|
=\|(\ga\ga'\ga^{-1})_*\mu_{ \ga x_0}\| =\|\mu_{\ga x_0}\|
\;\;\;{\rm and}\;\;\;f_{\ga' x, \,\varphi}= 
f_{x,\,\varphi\circ\ga' }\circ R_{\ga'} \;.
$$
Using this and Lemma \ref{lem:normalpatdens}, it is then easy to check
that $(\nu_{x})_{x\in\wt M}$ is also a Patterson density of dimension
$\delta'$ for $(\Ga',F')$.

Since the proof of Jensen's inequality does not require the measure to
be $\sigma$-additive, by the convexity of the exponential map, for
every $f\in \ell^\infty(\Ga'\bs\Ga)$, we have $e^{\mmm(f)}\leq
\mmm(e^f)$. In particular, for every $\alpha\in\Ga$, we have 
$$
e^{\chi(\alpha)}=
e^{\mmm(f_\alpha)}\leq \mmm(e^{f_\alpha})=\mmm(f_{\alpha x_0,\,1})=
\|\nu_{\alpha x_0}\|\;.
$$
By the first assertion of Corollary \ref{coro:conseqprincipombr}
applied by taking $(\Ga_0,\Ga,(\mu_x)_{\in\wt M})$ therein equal to
$(\Ga,\Ga',(\nu_x)_{\in\wt M})$, we have $\delta'\geq
\delta_{\Ga,\,F,\,\chi}$. Since $F$ is reversible, we have
$\delta_{\Ga,\,F,\,\chi}\geq\delta_{\Ga,\,F}$ by Equation
\eqref{eq:reveregalchipaschi}. Hence
$\delta'\geq\delta_{\Ga,\,F}$, and the result follows from
Equation \eqref{eq:croissanceexpocrit}.  
\cqfd

\brema \label{rem:refereeComment19}
{\rm The reversibility assumption on $F$ in Theorem 
\ref{theo:amenableegalexpocrit} is necessary.
}
\erema

\dem Let us construct examples which prove that Theorem
\ref{theo:amenableegalexpocrit} fails without the reversibility
assumption on $F$. Assume that $\Ga$ is torsion free and cocompact,
and that there exists a nontrivial character $\chi$ of $\Ga$. Let
$\omega_\chi$ be the $1$-form on $M$ corresponding to $\chi$, see the
remark above Proposition \ref{prop:proprieteelemexpocritchar}. Let
$\ga\in\Ga$ be such that $\chi(\ga)>0$: In particular, $\ga$ is
loxodromic. Let $x$ be a point on the translation axis of $\ga$.

For every $N\in\NN$, let $\wt F=N\;\wt\omega_\chi$. The Poincar\'e
series of $\Ga$ satisfies, for all $s\geq 0$,
$$
Q_{\Ga,\,F,\,x,\,x}(s)\geq \sum_{k\in\NN}e^{\int_x^{\ga^kx}(\wt F-s)}
=\sum_{k\in\NN}e^{k\int_x^{\ga x}(N\;\wt\omega_\chi-s)}
=\sum_{k\in\NN}e^{k(N\,\chi(\ga)-s\,\ell(\ga))}\;,
$$ 
which diverges if $s\leq \frac{N\chi(\ga)}{\ell(\ga)}$. Hence
$\delta_{\Ga,\,F}\geq N\,\frac{\chi(\ga)}{\ell(\ga)}$.

If $\Ga'$ is the kernel of $\chi$ and $F':\Ga'\backslash T^1\wt
M\ra\RR$ the map induced by $\wt F$, then (using for instance the
equality between the critical exponent and the Gurevich pressure of
$(\Ga',F')$ in Theorem \ref{theo:equalcritgur}), we have
$\delta_{\Ga',\,F'}=\delta_{\Ga',\,0}$, which is independent of
$N$. Taking $N$ big enough, we then have $\delta_{\Ga,\,F}>
\delta_{\Ga',\,F'}$. Since the finitely generated group $\Ga'\bs\Ga$
is abelian, hence amenable, this proves the result.  \cqfd

\bigskip A nonelementary normal subgroup $\Ga'$ of $\Ga$ cannot be too
small. For instance, we have seen that the inclusion
$\Lambda\Ga'\subset \Lambda\Ga$ is an equality. Similarly, the default
of the inequality \eqref{eq:croissanceexpocrit} to be an equality
cannot be too large, even if $\Ga'\bs\Ga$ is non-amenable. If $\wt M$
is a real hyperbolic space $\HH^n_\RR$ with $n\geq 3$ or a complex
hyperbolic space $\HH^n_\CC$ with $n\geq 2$, and if $\Ga$ is a lattice
in $\Isom(\wt M)$, then Shalom \cite[Theo.~1.4]{Shalom00} has proved
some results of this type. When $F=0$, the following result is due to
\cite[Théo.~2.2.1]{Roblin05}.

\btheo \label{theo:petitrou} (1) If $\Ga'$ is any nonelementary normal
subgroup of $\Ga$, then 
$$
\delta_{\Ga',\,F'}\geq
\frac{1}{2}\;\delta_{\Ga,\,F+F\circ\iota}\;.
$$

(2) If $F$ is reversible, if $\delta_{\Ga,\,2F} <2\;
\delta_{\Ga,\,F}$ and if $(\Ga,2F)$ is of divergence type, then
$\delta_{\Ga',\,F'}> \frac{1}{2}\;\delta_{\Ga,\,2F}$.
\etheo

Note that the condition that $\delta_{\Ga,\,F+F\circ\iota} <2\;
\delta_{\Ga,\,F}$ is satisfied for instance if the upper bound
$\sup_{\pi^{-1}(\C\Lambda\Ga)}|F|$ is small enough, by Lemma
\ref{lem:elemproppressure} (iv).

\medskip \dem Let $\delta=\delta_{\Ga,\,F+F\circ\iota}$ and $\delta'=
\delta_{\Ga',\,F'}$. If $F$ is reversible, then $F+F\circ\iota$ is
cohomologous to $2F$, hence $\delta_{\Ga,\,F+F\circ\iota}=
\delta_{\Ga,\,2F}$ and $(\Ga,2F)$ is of divergence type if and only if
$(\Ga,F+F\circ\iota)$ is of divergence type (see the remark at the end
of Subsection \ref{subsec:GibbsPoincareseries}).

We may assume that $\delta'<+\infty$.
Let $(\mu'_{x})_{x\in\wt M}$ be a Patterson density of dimension
$\delta'$ for $(\Ga',F')$, such that the dynamical system
$(\partial_\infty\wt M,\Ga',\mu'_{x_0})$ is ergodic, which exists as
seen in Subsection \ref{subsec:gibbspattersondens}.  Let
$R'=R_{\{x_0\}}(\delta')$ and $C'=C_{\{x_0\}}(\delta')$ be given by
Proposition \ref{prop:principombres} applied by taking $(\Ga_0,\Ga,
(\mu_{x})_{x\in\wt M})$ therein equal to $(\Ga,\Ga',
(\mu'_{x})_{x\in\wt M})$.

\medskip (1) For every $y\in\Ga x_0$, Proposition
\ref{prop:principombres} gives that
\begin{equation}\label{eq:minomasstotdens}
\|\mu'_y\|\geq \mu'_y\big(\OOO_yB(x_0,R')\big)\geq
\frac{\|\mu'_{x_0}\|}{C'}\;e^{\int_y^{x_0}(\wt F-\delta')}\;.
\end{equation}
Hence for every $s>\delta'$, we have
$$
\sum_{\ga\in\Ga}\;
e^{\int_{x_0}^{\ga x_0}(\wt F+\wt F\circ\iota-\delta'-s)}
\leq \frac{C'}{\|\mu'_{x_0}\|}\;
\sum_{\ga\in\Ga}\;\|\mu'_{\ga x_0}\|\;
e^{\int_{x_0}^{\ga x_0}(\wt F-s)}\;.
$$
By the first claim of Corollary \ref{coro:conseqprincipombr}, the sum
on the right hand side converges. This proves that $\delta\leq
2\;\delta'$.

\medskip (2) Assume for a contradiction that $\delta'=
\frac{\delta}{2}$. Since $(\Ga,F+F\circ\iota)$ is of divergence
type, let $(\mu_{x})_{x\in\wt M}$ be a Patterson density of dimension
$\delta$ for $(\Ga,F+F\circ\iota)$, which gives full measure to
$\Lambda_c\Ga$ (see Corollary \ref{coro:uniqpatdens}). By Corollary
\ref{coro:conseq2principombr} applied with $\Ga_0=\Ga$, let $r\geq R'$
be large enough so that $\mu_{x_0}$ gives full measure to
$\Lambda_{c,\,r}\Ga$.

Let us first prove that $\mu_{x_0}$ is absolutely continuous with
respect to $\mu'_{x_0}$. Let $B$ be any Borel subset of
$\partial_\infty\wt M$, and let $V$ be any open neighbourhood of
$B$. Let
$$
Z=\{z\in \Ga x_0\;:\; \OOO_{x_0}B(z,r)\subset V\}\;,
$$
and let $Z^*$ be the subset of $Z$ constructed in Lemma
\ref{lem:Vitalistrikesagain} with $s=r$. Using respectively 

$\bullet$~ the fact that $\Lambda_{c,\,r}\Ga$ has full measure with
respect to $\mu_{x_0}$,

$\bullet$~ the fact that $B\cap\Lambda_{c,\,r}\Ga$ is contained in
$\bigcup_{z\in Z} \OOO_{x_0}B(z,r)$ and Lemma
\ref{lem:Vitalistrikesagain},

$\bullet$~ Equation \eqref{eq:principombreupper}, which provides a
constant $c>0$ depending only on $R=5r$ and $K=\{x_0\}$,

$\bullet$~ the $\Ga$-equivariance of $(\mu_{x})_{x\in\wt M}$,
which implies that $\|\mu_{\ga x_0}\|=\|\mu_{x_0}\|$ for every
$\ga\in\Ga$,

$\bullet$~ Equation \eqref{eq:minomasstotdens},
 
$\bullet$~ the equality $\delta=2\delta'$,
 
$\bullet$~ Proposition \ref{prop:principombres} applied by taking
$(\Ga_0,\Ga,(\mu_x)_{\in\wt M})$ therein equal to
$(\Ga,\Ga',(\mu'_x)_{\in\wt M})$, since $r\geq R'$, and setting
$c'=\frac{c\,(C')^2\;\|\mu_{x_0}\|}{\|\mu'_{x_0}\|}$,
 
$\bullet$~ the properties of $Z^*$ seen in Lemma
\ref{lem:Vitalistrikesagain} and the definition of $Z$,
 
\noindent we have
\begin{align*}
\mu_{x_0}(B)& = \mu_{x_0}(B\cap\Lambda_{c,\,r}\Ga)\leq 
\sum_{z\in Z^*}\mu_{x_0}\big(\OOO_{x_0}B(z,5r)\big)
\\ &
\leq \sum_{z\in Z^*}c\;\|\mu_{z}\|\;
e^{\int_{x_0}^z(\wt F+\wt F\circ\iota-\delta)}= c\;\|\mu_{x_0}\|
\sum_{z\in Z^*}\;e^{\int_{x_0}^z(\wt F+\wt F\circ\iota-\delta)}
\\ &
\leq c\;\|\mu_{x_0}\|
\sum_{z\in Z^*}\;\frac{C'\|\mu'_z\|}{\|\mu'_{x_0}\|}\;
e^{-\int_z^{x_0}(\wt F-\delta')}\;e^{\int_{x_0}^z(\wt F+\wt F\circ\iota-\delta)}
\\ &=
\frac{c\,C'\;\|\mu_{x_0}\|}{\|\mu'_{x_0}\|}\;
\sum_{z\in Z^*}\;\|\mu'_z\|\;e^{\int_{x_0}^z(\wt F-\delta')}
\\ &
\leq c'\sum_{z\in Z^*}\;\mu'_{x_0}\big(\OOO_{x_0}B(z,r)\big)\leq c'\;
\mu'_{x_0}(V)\;.
\end{align*}
Since $V$ is an arbitrary neighbourhood of $B$, we hence have
$\mu_{x_0}\leq c'\; \mu'_{x_0}$, as required.

Now, the Radon-Nikodym derivative $\frac{d\,\mu_{x_0}}{d\,\mu'_{x_0}}$
is quasi-invariant under $\Ga'$ and the measure $\mu_{x_0}$ is not the
zero measure. By ergodicity of $(\partial_\infty\wt M,\Ga',
\mu'_{x_0})$, the map $\frac{d\,\mu_{x_0}}{d\,\mu'_{x_0}}$ is hence
positive $\mu'_{x_0}$-almost everywhere, so that $\mu_{x_0}$ and
$\mu'_{x_0}$ are equivalent. Therefore $\mu'_{x_0}$ is quasi-invariant
under $\Ga$, as is $\mu_{x_0}$.  By Proposition
\ref{prop:proprielemtwistPattdens} (3), there exists a character
$\chi$ of $\Ga$ such that $(\mu'_{x})_{x\in\wt M}$ is a twisted
Patterson density of dimension $\delta'$ for $(\Ga,F,\chi)$.  Hence
$\delta'\geq \delta_{\Ga,\,F,\,\chi}$ by Proposition
\ref{prop:proprielemtwistPattdens} (2). Since $F$ is reversible, we
have $\delta_{\Ga,\,F,\,\chi}\geq\delta_{\Ga,\,F}$ by Equation
\eqref{eq:reveregalchipaschi}. Since $\delta_{\Ga,\,F}>
\frac{\delta}{2}$ by assumption, we have $\delta'> \frac{\delta}{2}$,
a contradiction.  \cqfd

\subsection{Classification of ergodic Patterson densities 
on nilpotent co\-vers}
\label{subsec:classnilpotcov}

Let $\Ga'$ be a normal subgroup of $\Ga$, and let
$F':\Ga'\backslash T^1\wt M\ra\RR$ be the map induced by $\wt F$.

\medskip We start by recalling the definition of a nilpotent
group. Let $G$ be a group. If $H$ and $H'$ are subgroups of $G$, let
$\mathopen{[}H,H'\mathclose{]}$ be the subgroup of $G$ generated by
the elements $\mathopen{[}h,h'\mathclose{]}=hh'h^{-1}{h'}^{-1}$ with
$h\in H$ and $h'\in H'$. The {\it lower central series}\index{lower
  central series} of a group $G$ is the sequence of normal subgroups
$(G_n)_{n\in\NN}$ defined by induction by
$$ 
G_0=G\;\;\;{\rm and}\;\;\;G_{n+1}=
\mathopen{[}G,G_n\mathclose{]}\;\;{\rm for}\;\;n\in\NN\;.
$$
The group $G$ is {\it nilpotent}\index{nilpotent group} if there
exists $n\in\NN$ such that $G_{n}=\{e\}$. It is well known that a
nilpotent group is in particular amenable.

Note that some fundamental groups $\Ga$ of compact connected
hyperbolic orbifolds $M$, as hyperbolic triangle groups, have no
nontrivial character, but contain many (non convex-cocompact, so that
Proposition \ref{prop:classifdimdens} (2) applies) normal subgroups
$\Ga'$, hence some restriction on the quotient $\Ga/\Ga'$ in the
following result is needed to be able to construct characters.

\btheo \label{theo:classdensnilpocov}%
\index{theorem@Theorem!Classification of Patterson densities on
  nilpotent covers} Let $\Ga'$ be a normal subgroup of $\Ga$ such that
$\Ga/\Ga'$ is nilpotent.  Assume that $\delta_{\Ga,\,F}<+\infty$. For
every $\sigma\in\mathopen{[}\delta_{\Ga,\,F},+\infty\mathclose{[}\,$,
the set of Patterson densities of $(\Ga',F')$ of dimension $\sigma$,
giving full measure to $\Lambda_c\Ga$ and ergodic with respect to the
action of $\Ga'$, is equal to the set of twisted Patterson densities
of $(\Ga,F,\chi)$ of dimension $\sigma$ for the characters $\chi$ of
$\Ga$ vanishing on $\Ga'$ such that $\delta_{\Ga,\,F,\,\chi}=\sigma$,
giving full measure to $\Lambda_c\Ga$.  \etheo

\dem As usual with nilpotent groups, the proof proceeds by induction
on the nilpotency order. By taking the preimages of the terms of the
lower central series of $\Ga/\Ga'$, there exist $n\in\NN$ and
subgroups $\Ga_k$ of $\Ga$ for $k=0,\dots,n$ such that
$\Ga_0=\Ga,\Ga_n=\Ga'$, $\Ga_{k+1}$ is a normal subgroup of $\Ga_k$
and $\Ga_{k+1}/\Ga'= \mathopen{[}\Ga/\Ga',\Ga_{k}/\Ga'\mathclose{]}$
for $k=0,\dots, n-1$.  We denote by $F_k:\Ga_k\bs T^1\wt M\ra\RR$ the
map induced by $\wt F$. Let us fix $\sigma\geq \delta_{\Ga,\,F}$.

If $n=0$, then $\Ga=\Ga'$ and the only character of $\Ga$ vanishing on
$\Ga'$ is $\chi=0$. The result is then immediate, as we have seen in
Corollary \ref{coro:critdivergence} that if there exists a Patterson
density $(\mu_x)_{x\in X}$ for $(\Ga,F)$ of dimension $\sigma\geq
\delta_{\Ga,\,F}$ such that $\mu_x(\Lambda_c\Ga)>0$, then
$\sigma=\delta_{\Ga,\,F}$ and $(\Ga,F)$ is of divergence type, hence
the measures $\mu_x$ are ergodic for the action of $\Ga$ by Theorem
\ref{theo:critergo}. Therefore we assume that $n\geq 1$.

\medskip
\noindent {\bf Step 1. } There exists $c>0$ such that for every $k\geq
1$, for every Patterson density $(\mu_{x})_{x\in\wt M}$ of dimension
$\sigma$ for $(\Ga_k,F_k)$, giving full measure to $\Lambda_c\Ga$, and
for all $\alpha\in\Ga_{k-1}$, we have
$$
(\alpha^{-1})_*\mu_{\alpha x_0}\geq c\;
e^{\int_{x_0}^{\alpha x_0}(\wt F\circ \iota-\sigma)}\;\mu_{x_0}\;.
$$

\medskip \dem Let $c=\frac{1}{C_{\{x_0\}}(\sigma)}>0$, with the
notation of Proposition \ref{prop:principombres}.  For all $\ga\in\Ga$
and $\alpha\in\Ga_{k-1}$, since $\mathopen{[}\alpha,\ga\mathclose{]}
\in \mathopen{[}\Ga_{k-1},\Ga\mathclose{]}$ belongs to $\Ga_k$, we
have, by this proposition,
\begin{align*}
\|(\alpha^{-1})_*\mu_{\alpha \ga x_0}\| & = \|\mu_{\alpha \ga x_0} \|= 
\|\mu_{\ga \alpha x_0} \|
\geq\frac{1}{C_{\{x_0\}}(\sigma)}\;\|\mu_{\ga x_0} \|\;
e^{\int_{\ga \alpha x_0}^{\ga x_0} (\wt F-\sigma)}\\ &  
=c\;\|\mu_{\ga x_0} \|\;e^{\int_{x_0}^{\alpha x_0}(\wt F\circ \iota-\sigma)}\;.
\end{align*}
By Lemma \ref{lem:normalpatdens} since $\Ga_{k-1}$ normalises $\Ga_k$,
and by the Fatou-Roblin radial convergence theorem
\ref{theo:fatouroblin} applied by replacing $(\Ga_0,\Ga)$ therein by
$(\Ga,\Ga_{k-1})$, this proves the first step. \cqfd

\medskip
\noindent {\bf Step 2. } Let $\mu=(\mu_{x})_{x\in\wt M}$ be a
Patterson density of dimension $\sigma$ for $(\Ga',F')$, giving full
measure to $\Lambda_c\Ga$, and ergodic for the action of $\Ga'$. Then
there exists a character $\chi$ of $\Ga$, vanishing on $\Ga'$, such that
$\mu$ is a twisted Patterson density of dimension $\sigma$ for
$(\Ga,F,\chi)$.

\medskip \dem Let us first prove that, for $k\in\{2,\dots,n\}$, if
$\mu$ is a Patterson density for $(\Ga_k,F_k)$, then $\mu$ is a
Patterson density for $(\Ga_{k-1},F_{k-1})$.

Since $\mu_{x_0}$ is ergodic under $\Ga'$ hence under $\Ga_{k}$, and
since $\mu_{x_0}$ is absolutely continuous with respect to
$(\alpha^{-1})_*\mu_{\alpha x_0}$ for every $\alpha\in\Ga_{k-1}$ by
Step 1, there exists (by Proposition
\ref{prop:proprielemtwistPattdens} (3) with $(\Ga,\Ga')$ therein
replaced by $(\Ga_{k-1},\Ga_k)$) a character $\chi_{k-1}$ of
$\Ga_{k-1}$, trivial on $\Ga_k$, such that $\mu_{x_0} =
e^{-\chi_{k-1}(\alpha)}\;(\alpha^{-1})_*\mu_{\alpha x_0}$ for every
$\alpha\in\Ga_{k-1}$. For every $\alpha\in\Ga_{k-1}$, we have
$$
\chi_{k-1}(\alpha)=\log\|\mu_{\alpha x}\|-\log\|\mu_{x}\|
$$ 
for every $x\in\wt M$ (by the definition of $\chi_{k-1}$ for $x=x_0$
and by Lemma \ref{lem:normalpatdens} and Equation
\eqref{eq:radonykodensity} to extend to any $x$).

Since $\Ga/\Ga'$ is amenable, let $\mmm$ be a right-invariant mean on
$\Ga'\bs\Ga=\Ga/\Ga'$. As seen in the proof of Theorem
\ref{theo:amenableegalexpocrit}, for every $\alpha\in\Ga$, the map
$f_\alpha$ which associates to $\mathopen{[}\ga\mathclose{]}
\in\Ga/\Ga'$ the real number $\log\|\mu_{\ga\alpha x_0}\|
-\log\|\mu_{\ga x_0}\|$ is well defined and bounded, and the map
$\chi:\Ga\ra\RR$ defined by $\chi(\alpha)=\mmm(f_\alpha)$ is a
character of $\Ga$.  This character hence coincides with $\chi_{k-1}$
on $\Ga_{k-1}$, by applying the above centred formula to $x=\ga x_0$
for every $ \ga\in \Ga$ and since, as seen in the proof of Step 1,
$\|\mu_{\ga\alpha x_0}\|= \|\mu_{\alpha\ga x_0}\|$ for every $\ga\in
\Ga$.  Since $\Ga_{k-1}/\Ga'\subset \mathopen{[}\Ga/\Ga',
\Ga_{k-2}/\Ga'\mathclose{]}$, the additive character $\chi$ vanishes
on $\Ga_{k-1}$, hence $\chi_{k-1}=0$. This proves the preliminary
claim.

\medskip Now, by induction, we hence have that $\mu$ is a Patterson
density for $(\Ga_1,F_1)$. Applying again Step 1 gives that $\mu$ is
quasi-invariant under $\Ga$. Applying Proposition
\ref{prop:proprielemtwistPattdens} (3) (using for this the ergodicity
of $\mu_{x_0}$ under $\Ga'$ hence under $\Ga_1$), the second step follows.
\cqfd

\medskip
\noindent {\bf Step 3. } If $\chi$ is a character of $\Ga$ vanishing
on $\Ga'$, if $(\mu_{x})_{x\in\wt M}$ is a twisted Patterson density
of dimension $\sigma$ for $(\Ga,F,\chi)$, giving full measure to
$\Lambda_c\Ga$, then $\sigma=\delta_{\Ga,\,F,\,\chi}$ and $\mu_{x_0}$
is ergodic for the action of $\Ga'$.

\medskip \dem Let $A$ be a Borel subset of $\Lambda_c\Ga$ invariant
under $\Ga'$, with nonzero $\mu_{x_0}$-measure, and let
$\mathbbm{1}_A$ be its characteristic function. Then
$(\mathbbm{1}_A\,\mu_{x})_{x\in\wt M}$ is a Patterson density of
dimension $\sigma$ for $(\Ga',F')$, giving full measure to
$\Lambda_c\Ga$. The preliminary  claim in the proof of Step 2 shows that
$(\mathbbm{1}_A\,\mu_{x})_{x\in\wt M}$ is also a Patterson density of
dimension $\sigma$ for $(\Ga_1,F_1)$. By Step 1,
$\mathbbm{1}_A\,\mu_{x}$ is quasi-invariant under $\Ga$, hence $A$ is
invariant up to a set of measure $0$ under $\Ga$. In particular,
$(\mathbbm{1}_A\,\mu_{x})_{x\in\wt M}$ is also a twisted Patterson
density of dimension $\sigma$ for $(\Ga,F,\chi)$, giving full measure
to $\Lambda_c\Ga$. By their uniqueness property seen in Corollary
\ref{coro:twiconseqHTSR} (1) and (2), we have
$\sigma=\delta_{\Ga,\,F,\,\chi}$ and there exists $c>0$ such that
$\mathbbm{1}_A\,\mu_{x_0}=c\,\mu_{x_0}$. In particular
$\mu_{x_0}({\,}^c\!A)=0$, which proves the ergodicity of $\mu_{x_0}$
under $\Ga'$.  \cqfd 

\medskip Step 2 and Step 3 combine together to prove Theorem
\ref{theo:classdensnilpocov}. \cqfd

\medskip The following result is then immediate. The assertions (2)
and (3) of Theorem \ref{theo:claspatdensnilpcovintro} in the
introduction follow from it and from Corollary
\ref{coro:twiconseqHTSR}.

\bcoro Let $\Ga'$ be a normal subgroup of $\Ga$ such that
$\Ga/\Ga'$ is nilpotent.  

(1) Assume that $\delta_{\Ga,\,F}<+\infty$, that $(\Ga,F)$ is of
divergence type and that $F$ is reversible. Then the unique (up to a
scalar multiple) Patterson density $(\mu_{\Ga,\,F,\,x})_{x\in\wt M}$
of dimension $\delta_{\Ga,\,F}$ for $(\Ga,F)$ is also the unique (up
to a scalar multiple) Patterson density of dimension
$\delta_{\Ga,\,F}$ for $(\Ga',F')$ giving full measure to
$\Lambda_c\Ga$.

(2) If $\Ga$ is convex-cocompact, then the set of ergodic Patterson
densities for $(\Ga',F')$ with support $\Lambda\Ga$ is the set of
twisted Patterson densities of $(\Ga,F,\chi)$ of dimension
$\delta_{\Ga,\,F,\,\chi}$ with support $\Lambda\Ga$ for the characters
$\chi$ of $\Ga$ vanishing on $\Ga'$.  
\ecoro

\dem (1) By Corollary \ref{coro:uniqpatdens}, there exists a unique (up to
a scalar multiple) Patterson density $\mu=(\mu_{\Ga,\,F,\,x})_{x\in\wt
  M}$ of dimension $\sigma=\delta_{\Ga,\,F}$ for $(\Ga,F)$, and it
gives full measure to $\Lambda_c\Ga$. By restriction, $\mu$ is also a
Patterson density of dimension $\sigma$ for $(\Ga',F')$ giving full
measure to $\Lambda_c\Ga$.

Let $\nu=(\nu_{x})_{x\in\wt M}$ be another Patterson density of
dimension $\sigma$ for $(\Ga',F')$ giving full measure to
$\Lambda_c\Ga$. Assume furthermore that $\nu$ is ergodic with respect
to the action of $\Ga'$. By Theorem \ref{theo:classdensnilpocov},
there exists a character $\chi$ of $\Ga$ such that $\nu$ is a twisted
Patterson density for $(\Ga,F,\chi)$ of dimension
$\delta_{\Ga,\,F,\,\chi} =\sigma$.  By Proposition
\ref{prop:proprielemtwistPattdens} (5) whose assumptions are
satisfied, we have $\chi=0$.  Hence $\nu$ is a scalar multiple of
$\mu$. By Krein-Milman's theorem, the ergodicity assumption on $\nu$
may be dropped. This proves Assertion (1).

\medskip (2) Since $\Lambda_c\Ga=\Lambda\Ga$ when $\Ga$ is
convex-cocompact, Assertion (2) follows immediately from Theorem
\ref{theo:classdensnilpocov}.
\cqfd

\medskip
We refer to \cite{Roblin05} for possible improvements of the second
assertion when $\Ga$ is only assumed to be geometrically finite (with
$\delta_{\Ga,\,F}<+\infty$).

It follows from the first assertion that not only the dynamical system
$(\partial_\infty\wt M,\Ga,\mu_{\Ga,\,F,\,x})$ is ergodic, but that
the dynamical system $(\partial_\infty\wt M,\Ga',\mu_{\Ga,\,F,\,x})$
is also ergodic. This extends results of Babillot-Ledrappier
\cite{BabLed96} when $F=0$ and $\Ga'/\Ga$ is isomorphic to $\ZZ^d$ for
$d\in\NN$, of Hamenst\"adt \cite{Hamenstadt02} when $M$ is compact,
and of Roblin \cite{Roblin05} when $F=0$ .

\medskip \rem Let $\Ga'$ be a subgroup of $\Ga$. Assume that
$\delta_{\Ga,\,F} <+\infty$ and that $(\Ga,F)$ is of divergence
type. Let $\wt m_F$ be the Gibbs measure of $(\Ga,\,F)$ on $T^1\wt M$
(which is unique up to a scalar multiple). Let $\phi=
(\phi_t)_{t\in\RR}$ be the geodesic flow on the (orbifold) cover
$M'=\Ga'\bs\wt M$ of $M$ defined by the group $\Ga'$. Let $m'_F$ be
the measure induced on $T^1 M'=\Ga'\bs T^1\wt M$ by $\wt m_F$.  Note
that if $\Ga'$ has infinite index in $\Ga$, the Gibbs measure $m'_F$
is infinite.

Even if $\Ga'$ is normal, if $\Ga'\bs\Ga$ is nilpotent and if $m_F$ is
finite, the dynamical system
$$
(T^1 M',\phi,m'_F)
$$
is not necessarily ergodic: the ergodicity of $(\partial_\infty\wt
M,\Ga',\mu_{\Ga,\,F,\,x})$ does not imply the ergodicity of
$(T^1 M',\phi,m'_F)$. See \cite{Rees81} for examples.

A necessary and sufficient condition on the left quotient $\Ga'\bs\Ga$
for $(T^1 M',\phi,m'_F)$ to be ergodic is still unknown (see the end
of the introduction for more details on open problems).


%




\newpage
\printglossaries

\newpage
\addcontentsline{toc}{section}{Index}
\input{countgibbs.ind}

\newpage
\addcontentsline{toc}{section}{References}
\bibliography{../biblio}

\bigskip
\noindent {\small 
\begin{tabular}{l} 
D\'epartement de Math\'ematique, Bât.~425, UMR 8628 CNRS\\ 
Université Paris-Sud\\ 
91505 ORSAY Cedex, FRANCE\\ 
{\it e-mail: frederic.paulin@math.u-psud.fr} 
\end{tabular} 
\\ 
 \medskip
\\ 
\begin{tabular}{l}  
Mathematics Institute,
\\  University of Warwick,
\\ Coventry, CV4 7AL, UK
\\{\it e-mail: mpollic@maths.warwick.ac.uk} 
\end{tabular} 
\\ 
 \medskip
\\ 
\begin{tabular}{l}  
LAMFA, UMR 7352 CNRS, UFR des Sciences,
\\  Université Picardie Jules Verne
\\ 33 rue Saint Leu
\\ 80 000 AMIENS, FRANCE
\\{\it e-mail: barbara.schapira@u-picardie.fr} 
\end{tabular}
}

\end{document}